\documentclass[10pt]{article}

\usepackage[a4paper, margin=1in]{geometry}

\usepackage{amsmath,amssymb,amsthm}
\usepackage{array}
\usepackage[shortlabels]{enumitem}
\usepackage{mleftright}

\usepackage{tikz}
\usetikzlibrary{arrows.meta}

\usepackage{cite}

\usepackage[colorlinks=true,allcolors=blue!50!black]{hyperref}

\allowdisplaybreaks
\numberwithin{equation}{section}

\setlist[enumerate,1]{leftmargin=8mm,label=(\roman*)}
\setlist[itemize,1]{leftmargin=6mm}

\newtheorem{theorem}{Theorem}[section]
\newtheorem{corollary}{Corollary}[theorem]
\newtheorem{proposition}[theorem]{Proposition}
\newtheorem{lemma}[theorem]{Lemma}

\newtheorem{definition}[theorem]{Definition}
\newtheorem{system}[theorem]{System}

\newtheorem{remark}[theorem]{Remark}

\DeclareMathOperator{\sgn}{sgn}
\DeclareMathOperator{\Res}{Res}

\renewcommand{\Re}{\operatorname{Re}}
\renewcommand{\Im}{\operatorname{Im}}

\newcommand{\lpc}{\Delta}
\newcommand{\grad}{\nabla}
\newcommand{\dr}{{\mkern-2mu\not\mkern2mu}\partial}

\newcommand{\<}{\langle}
\renewcommand{\>}{\rangle}

\newcommand{\norm}[1]{\mleft\lVert #1\mright\rVert}
\newcommand{\abs}[1]{\mleft\lvert #1\mright\rvert}
\newcommand{\brac}[1]{\mleft( #1\mright)}
\newcommand{\sbrac}[1]{\mleft[ #1\mright]}
\newcommand{\inner}[1]{\mleft\< #1\mright\>}
\newcommand{\set}[1]{\left\{#1\right\}}
\newcommand{\eva}[1]{\mleft.#1\mright|}

\newcommand{\nn}{|\mkern-2mu|\mkern-2mu|}

\newcommand{\Z} {\mathbb{Z}}
\newcommand{\N} {\mathbb{N}}
\newcommand{\C} {\mathbb{C}}
\newcommand{\R} {\mathbb{R}}

\renewcommand{\P}{\mathbb{P}}
\newcommand{\E}{\mathbb{E}}

\renewcommand{\S}{\mathcal{S}}

\newcommand{\mb}{\mathbf}
\newcommand{\bs}[1]{\boldsymbol{\mathbf{#1}}}

\renewcommand{\d}{\textup{d}}
\newcommand{\ph}{\,\cdot\,}

\makeatletter
\def\env@dmatrix{\hskip -\arraycolsep
  \let\@ifnextchar\new@ifnextchar
  \extrarowheight=2ex
  \array{*\c@MaxMatrixCols{>{\displaystyle}c}}}

\newenvironment{bdmatrix}
  {\left(\env@dmatrix}
  {\endmatrix\right)}
\makeatother

\newcommand{\mat}[1]{\begin{pmatrix}#1\end{pmatrix}}
\newcommand{\bmat}[1]{\begin{bdmatrix}#1\end{bdmatrix}}

\makeatletter
\newcommand{\pushright}[1]{\ifmeasuring@#1\else\omit\hfill$\displaystyle#1$\fi\ignorespaces}
\makeatother

\newcommand{\ENEQ}{\leavevmode\vspace*{-\dimexpr\abovedisplayskip + \baselineskip}}

\makeatletter
\newcommand{\trimtop}[2]{\mathpalette{\trimtop@{#1}}{#2}}
\newcommand{\trimtop@}[3]{\setbox\z@\hbox{$\m@th#2{#3}$}\ht\z@\dimexpr\ht\z@-#1\relax\box\z@}
\makeatother

\makeatletter
\newcommand{\trimbot}[2]{\mathpalette{\trimbot@{#1}}{#2}}
\newcommand{\trimbot@}[3]{\setbox\z@\hbox{$\m@th#2{#3}$}\dp\z@\dimexpr\dp\z@-#1\relax\box\z@}
\makeatother

\title{Stokes flow in a wedge with a moving contact point:\\ well-posedness with a fixed contact angle and Navier slip}

\author{Manuel V. Gnann\thanks{Delft Institute of Applied Mathematics, Delft University of Technology, 2628 CD Delft, Netherlands. Email: \href{mailto:M.V.Gnann@tudelft.nl}{M.V.Gnann@tudelft.nl}.}\and Hans Knüpfer\thanks{Institute for Mathematics and Interdisciplinary Center for Scientific Computing (IWR), Heidelberg University, 69120 Heidelberg, Germany. Email: \href{mailto:knuepfer@uni-heidelberg.de}{knuepfer@uni-heidelberg.de}.}\and King Ming Lam\thanks{Delft Institute of Applied Mathematics, Delft University of Technology, 2628 CD Delft, Netherlands. Email: \href{mailto:K.M.Lam@tudelft.nl}{K.M.Lam@tudelft.nl}.}\and Nader Masmoudi\thanks{NYUAD Research Institute, New York University Abu Dhabi, P.O.\ Box 129188, Abu Dhabi, United Arab Emirates, and Courant Institute of Mathematical Sciences, New York University, New York, NY 10012, USA. Email: \href{mailto:nm30@nyu.edu}{nm30@nyu.edu}.}\and Jonas Sauer\thanks{Institute of Mathematics, Friedrich Schiller University Jena, 07737 Jena, Germany. Email: \href{mailto:jonas.sauer@uni-jena.de}{jonas.sauer@uni-jena.de}.}}

\date{}

\begin{document}

\maketitle

\begin{abstract}
We study 2D Stokes flow in a wedge-like region bounded by a solid substrate with Navier slip and by a free surface driven by surface tension. The two meet at a moving contact point with a fixed contact angle. This is the model of a moving contact line at scales below the slip length, in its simplest setting. The slip length is the only length scale, and nothing is prescribed at the contact point beyond the angle: no relaxation law for the angle, and no contact-point dissipation. Slip removes the paradox found by Huh and Scriven, the infinite dissipation at a contact line that moves without slip, but in this model a weaker singularity remains. We give, to our knowledge, the first rigorous characterisation of this singularity for solutions of the free-boundary problem near the stationary state. For small contact angles and small initial data satisfying two compatibility conditions, we prove local existence, uniqueness and Lipschitz dependence on the data. The compatibility conditions are the price of the fixed angle for solutions regular in time up to the initial time. Rather than hiding the singularity, the solution class carries it explicitly in corner expansions. The velocity is continuous up to the contact point. The pressure has a logarithmic singularity, the slip-regularised version of that of Huh and Scriven. Its coefficient is proportional to the contact-point velocity. The slope of the free boundary has a corresponding logarithmic correction, the slip-scale counterpart of the Cox--Voinov relation. The linearised problem combines terms of different scaling, which prevents a direct treatment by the Mellin transform. We solve it by a Lax--Milgram argument, and an auxiliary test function problem shows that the variational solution satisfies the equations. Explicit Mellin representations then yield the corner expansions, and the nonlinear problem is solved by a contraction argument.
\end{abstract}

\tableofcontents

\section{Introduction}\label{introduction}

\subsection{The problem}

A viscous fluid resting on a solid substrate and bounded by a free surface meets the substrate along a contact line. At rest the configuration has been understood since Young \cite{Young}: the interface meets the substrate at the angle at which the surface tensions balance. In motion --- a drop spreading, a plate drawn out of a bath, a meniscus advancing along a capillary --- the neighbourhood of the contact line resists the classical description, and the difficulty is old and well documented. Huh and Scriven \cite{HuhScriven} computed the local flow under the no-slip condition, a similarity solution of the kind introduced for viscous corner flows by Dean and Montagnon \cite{DeanMontagnon}, Moffatt \cite{Moffatt} and Taylor \cite{Taylor}. Huh and Scriven found a stress of the order $\mu U/r$ at distance $r$ from a contact line moving with speed $U$, the Huh--Scriven singularity. The rate of dissipation and the force on the substrate therefore diverge logarithmically; in their phrase, ``not even Herakles could sink a solid if the physical model were entirely valid, which it is not''. Dussan and Davis \cite{DussanDavis} traced the difficulty to the kinematics rather than to any constitutive assumption: no slip and a moving contact line force the interface to roll onto the substrate and the velocity field to be multivalued at the contact line. This is the paradox of Huh and Scriven --- the classical model admits no moving contact line with finite dissipation --- and it has to be resolved below some small scale.

What resolves it is the subject of a large physical literature, for which we refer to the reviews of Dussan \cite{Dussan}, de Gennes \cite{deGennes}, Bonn, Eggers, Indekeu, Meunier and Rolley \cite{Bonn} and Snoeijer and Andreotti \cite{SnoeijerAndreotti}. The mechanisms proposed include precursor films, diffuse interfaces, evaporation and condensation near the contact line, and a non-Newtonian rheology of the fluid; for the last we refer to Weidner and Schwartz \cite{WeidnerSchwartz}, to the work of the first author with Lienstromberg and Nik \cite{GnannLienstrombergNik}, and to the references therein. The oldest and the most common is slip: following Navier \cite{Navier}, the no-slip condition is replaced by proportionality of the tangential velocity at the substrate to the shear rate, and the proportionality constant is called the slip length. The Navier slip condition has since been derived rigorously as the effective boundary condition of viscous flow over a rough surface \cite{JagerMikelic}, in a fluid domain without a contact line.

With slip the force on the substrate is finite \cite{Hocking}, and the flow organises itself into two regions, which carry two different notions of contact angle. In the inner region, at scales comparable with the slip length, the full free-boundary Stokes problem has to be solved, and the interface meets the substrate at the microscopic angle. This angle is a property of the materials in contact: it is either held at its equilibrium value or relaxed towards it by a constitutive law imposed at the contact line. In the outer region, at scales large compared with the slip length, the slope read off the interface is the apparent, or macroscopic, angle, which is not prescribed but produced by the flow. It is the one that experiments see. For small angles its cube differs from the cube of the microscopic angle by a multiple of the speed of the contact line, with a factor that depends logarithmically on the distance to the contact line divided by the slip length. This is the Cox--Voinov relation, derived by matched asymptotics in the works of Voinov \cite{Voinov}, Hocking \cite{Hocking} and Cox \cite{Cox}; see also \cite{HockingRival}. For a completely wetting fluid it leads to Tanner's law \cite{Tanner} for the rate at which a drop spreads.

The inner region is where the present paper is set. With Navier slip and a fixed contact angle the flow at a moving contact line is still not smooth. In this setting the contact line does not move if the velocity is continuously differentiable up to it, and asymptotic analysis gives a logarithmic singularity of the pressure; see \cite{FrickeKohneBothe} and the references therein. Unlike the Huh--Scriven stress of the order $\mu U/r$, the logarithm is integrable up to the contact line, so that the pressure exerts a finite force on the substrate. The question we address is how this residual singularity is organised: its exact form in the pressure, how the free surface bends above the contact point, and which data are admissible. We address this for two-dimensional Stokes flow with Navier slip, surface tension and a fixed contact angle. With a dynamic contact angle that relaxes to its equilibrium value the flow is not smooth either \cite{FrickeKohneBothe}; see Section \ref{comments and related work}.

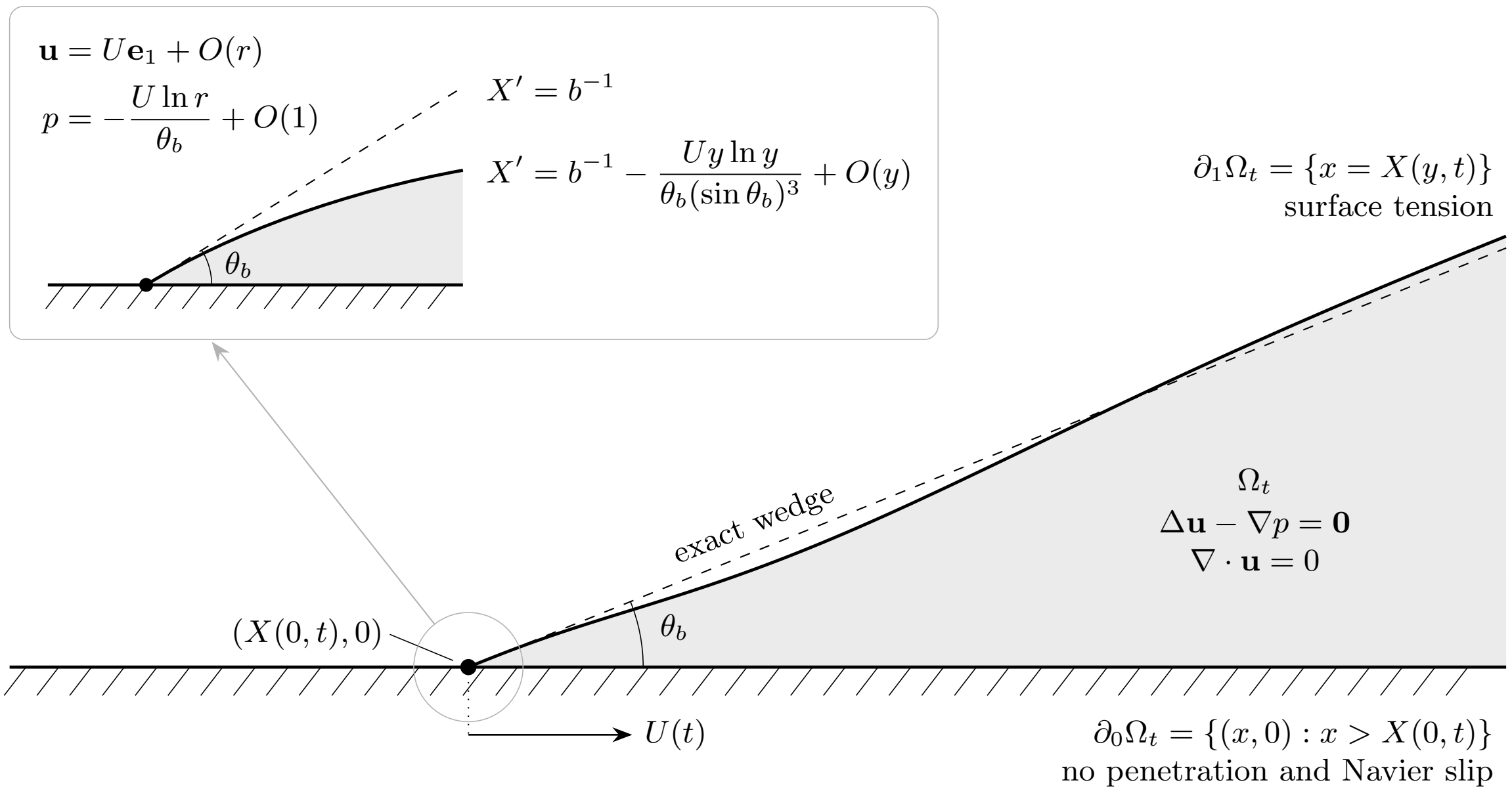
\begin{figure}[ht]
\centering
\begin{tikzpicture}[scale=1.1,>={Stealth[length=2mm]},line join=round]

\fill[black!8] plot[domain=0:9.5,samples=120,smooth]
  ({\x},{0.4040*\x+(\x*\x*\x*(\x-6)/73)*exp(-\x/1.6)}) -- (9.5,0) -- (0,0) -- cycle;

\draw[line width=0.9pt] (-4.2,0) -- (9.5,0);
\foreach \i in {-4.05,-3.75,...,9.4} \draw[line width=0.3pt] (\i,0) -- ++(-0.2,-0.26);

\draw[dashed] (0,0) -- (9.5,3.838);
\node[rotate=22,anchor=south,inner sep=3pt,font=\small] at (2.7,1.091) {exact wedge};

\draw[line width=0.9pt] plot[domain=0:9.5,samples=120,smooth]
  ({\x},{0.4040*\x+(\x*\x*\x*(\x-6)/73)*exp(-\x/1.6)});

\draw[line width=0.3pt] (1.6,0) arc[start angle=0,end angle=22,radius=1.6];
\node[font=\small] at (11:1.92) {$\theta_b$};

\fill (0,0) circle (2.1pt);
\draw[line width=0.3pt] (-0.72,0.30) -- (-0.14,0.06);
\node[anchor=east,font=\small,inner sep=1pt] at (-0.74,0.30) {$(X(0,t),0)$};
\draw[dotted] (0,-0.06) -- (0,-0.62);
\draw[->,line width=0.5pt] (0,-0.62) -- (1.5,-0.62) node[right,font=\small] {$U(t)$};

\node[font=\small,align=center] at (7.2,1.35)
  {$\Omega_t$\\[1pt]$\lpc\mb u-\grad p=\mb 0$\\$\grad\cdot\mb u=0$};
\node[anchor=south east,align=right,font=\small] at (9.5,4.02)
  {$\partial_1\Omega_t=\{x=X(y,t)\}$\\ surface tension};
\node[anchor=north east,align=right,font=\small] at (9.5,-0.36)
  {$\partial_0\Omega_t=\{(x,0):x>X(0,t)\}$\\ no penetration and Navier slip};

\draw[->,line width=0.4pt,black!30] (128.25:0.5) -- (-2.35,2.98);
\draw[line width=0.3pt,black!30] (0,0) circle (0.5);

\begin{scope}[shift={(-4.2,3.0)}]
\fill[white,rounded corners=4pt] (0,0) rectangle (8.5,3.05);
\draw[line width=0.3pt,black!30,rounded corners=4pt] (0,0) rectangle (8.5,3.05);
\fill[black!8] (1.25,0.5) .. controls (2.013,0.977) and (3.1,1.36) .. (4.15,1.55)
  -- (4.15,0.5) -- cycle;
\draw[line width=0.9pt] (0.35,0.5) -- (4.15,0.5);
\foreach \i in {0.5,0.75,...,4.15} \draw[line width=0.3pt] (\i,0.5) -- ++(-0.17,-0.22);
\draw[dashed] (1.25,0.5) -- (4.15,2.312);
\draw[line width=0.9pt] (1.25,0.5) .. controls (2.013,0.977) and (3.1,1.36) .. (4.15,1.55);
\draw[line width=0.3pt] (1.85,0.5) arc[start angle=0,end angle=32,radius=0.6];
\node[font=\small] at (2.10,0.685) {$\theta_b$};
\fill (1.25,0.5) circle (1.8pt);
\node[anchor=west,font=\small] at (4.25,2.31) {$X'=b^{-1}$};
\node[anchor=west,font=\small] at (4.25,1.50) {$\displaystyle X'=b^{-1}-{Uy\ln y\over\theta_b(\sin\theta_b)^3}+O(y)$};
\node[font=\small] at (1.55,2.26)
  {$\begin{aligned}\mb u&=U\mb e_1+O(r)\\p&=-{U\ln r\over\theta_b}+O(1)\end{aligned}$};
\end{scope}

\end{tikzpicture}
\caption{The setting of System \ref{dimensionless system}, drawn for a receding contact point, which moves to the right with the velocity $U>0$. The inset magnifies the corner: the pressure carries the logarithmic singularity of Theorem \ref{main theorem}(iv), and the slope of the free surface leaves $b^{-1}$ through a term $y\ln y$. This term lowers the apparent angle above a receding contact point and raises it above an advancing one. The resulting bend of the free surface below the exact wedge is drawn exaggerated, as is the angle in the inset.}
\label{setting figure}
\end{figure}

Precisely, we consider the fluid domain $\Omega_t=\{(x,y):y>0,\ x>X(y,t)\}$, bounded by the substrate $\{y=0\}$ and by the free surface $\{x=X(y,t)\}$, which meet at the contact point $(X(0,t),0)$. In the dimensionless variables introduced in Section \ref{setup}, in which lengths are measured in units of the slip length, the velocity $\mb u$, the pressure $p$ and the free boundary $X$ solve
\begin{alignat}{3}
\lpc\mb u-\grad p=\mb 0,\qquad\grad\cdot\mb u&=0\qquad&&\text{in}\quad\Omega_t,\nonumber\\
u_y=0,\qquad u_x-\partial_yu_x&=0\qquad&&\text{on}\quad\{y=0\},\nonumber\\
\mathbb T[\mb u,p]\mb n&=-X''(1+X'^2)^{-{3\over 2}}\mb n\qquad&&\text{on}\quad\{x=X(y,t)\},\label{intro system}\\
\partial_tX+u_y\partial_yX&=u_x\qquad&&\text{on}\quad\{x=X(y,t)\},\nonumber\\
X'(0,t)=b^{-1},\qquad\lim_{y\to\infty}X'(y,t)&=b^{-1},\qquad X(\ph,0)=X_0,\nonumber
\end{alignat}
where $\mathbb T[\mb u,p]=2\S(\mb u)-pI$ is the stress tensor, $\mb n$ the outward unit normal, and $b=\tan\theta_b$ with $\theta_b$ the prescribed contact angle. The domain is thus a perturbation of the wedge $\{0<y<b(x-X(0,t))\}$ of opening $\theta_b$ with its corner at the contact point, and the free boundary is a perturbation of the side $x=X(0,t)+b^{-1}y$ of this wedge. The contact point is free to move, its position being determined by the flow, while the angle at which the free surface meets the substrate is fixed at $\theta_b$; see Figure \ref{setting figure}. This is the model of a partially wetting fluid at scales below the slip length, where the angle is the microscopic one and not the apparent angle of the Cox--Voinov relation \cite{Voinov,Cox}.

\subsection{The main result}

To state the result we transform the free boundary into a profile on the reference wedge $\Omega=\{0<y<bx\}$, whose free side $\partial_1\Omega=\{(r\cos\theta_b,r\sin\theta_b):r>0\}$ is parametrised by the distance $r$ from the corner. The velocity is transformed to $\Omega$ by the Piola transform of Lemma \ref{Piola transform}, and the pressure by composition with the flattening $\bs\Psi_t$ of Lemma \ref{regularised flattening}. We call them the transformed velocity and the transformed pressure. The initial profile is
\begin{align}\label{intro profile}
\phi_0(r):=b^{-1}-X_0'(r\sin\theta_b),
\end{align}
the deviation of the initial slope from that of the wedge, which vanishes at the contact point by the contact angle condition and at infinity by the condition on the slope there. We ask $\phi_0$ to have a corner expansion
\begin{align}\label{intro expansion}
\phi_0(r)=c_{\phi_0,\dr_{-1}}r+c_{\phi_0,\dr_{-1}^2}r\ln r+\phi_0^{[2]}(r),
\end{align}
and we measure the size of $\phi_0$ by the data norm
\begin{align}\label{intro data norm}
\mathcal D:=\mathcal D(\phi_0):=(\sin\theta_b)^2\brac{\|\phi_0\|_{\partial_1\Omega,{1\over 2}-\kappa,l+{7\over 2}}+\|\phi_0\|_{\partial_1\Omega,{3\over 2}-\kappa,l+{7\over 2}}+\|\phi_0^{[2]}\|_{\partial_1\Omega,{5\over 2}-\kappa,l+{7\over 2}}},
\end{align}
where $\kappa\in(1/4,1/2)$ and the integer $l\geq 1$ are parameters of the theory. Here, for a weight $a$ and an integer $j\geq 0$,
\begin{align}\label{intro norm}
\|f\|_{\partial_1\Omega,a,j}^2:=\sum_{k=0}^j\int_0^\infty|r^{-a}(r\partial_r)^kf|^2\d r
\end{align}
is the weighted Sobolev norm of Definition \ref{weighted Sobolev norms}; for a half-integer $j$ the sum runs to $j-1/2$ and the term with $k=j$ is added, which Definition \ref{fractional counts} defines through the Mellin transform. The subscript $\dr_{-d}^k$ of a coefficient names the term that the coefficient multiplies: with $\dr_{-d}:=r\partial_r-d$, it is the lowest power of $\dr_{-d}$ that annihilates $r^d(\ln r)^{k-1}$. The remainder $\phi_0^{\trimtop{0.3pt}{[2]}}$ vanishes to the order $r^{2-\kappa}$ at the contact point. The reason for the term $r\ln r$ is the main structural fact about the problem, and it is best stated through the two compatibility conditions that the data must satisfy.

The first condition, \eqref{C1}, relates the coefficient of $r\ln r$ in \eqref{intro expansion} to the initial velocity $U(0)$ of the contact point,
\begin{align}\label{intro C1}
c_{\phi_0,\dr_{-1}^2}={U(0)\over\theta_b(\sin\theta_b)^2},\tag{C1}
\end{align}
in which $U(0)$ is itself determined by $\phi_0$: it is the horizontal velocity at the corner of the stationary Stokes problem that the initial profile drives (Theorem \ref{nonlinear initial data}(i)). The condition says that a moving contact point forces a logarithmic singularity $-U\theta_b^{-1}\ln r$ of the pressure at the corner, and that the normal stress balance on the free surface converts this singularity into a curvature $\sim\ln r$, hence a slope $\sim r\ln r$, of the free boundary. An initial interface prepared without this bend cannot be continued in our solution class.

The second condition, \eqref{C2}, is a regularity condition of the same kind one time derivative up. The equations determine from $\phi_0$ the first time derivative $\dot\phi_0$ of the profile at $t=0$ and, under \eqref{intro C1}, the second, $\ddot\phi_0$ (Theorem \ref{nonlinear initial data}); we call the triple $(\phi_0,\dot\phi_0,\ddot\phi_0)$ the initial jet. The second time derivative must vanish at the contact point in the weighted sense
\begin{align}\label{intro C2}
\|\ddot\phi_0\|_{\partial_1\Omega,{3\over 2}-\kappa,l+{3\over 2}}<\infty.\tag{C2}
\end{align}
Under the second-order corner expansion \eqref{intro second-order expansion} of the data and the strengthened angle condition \eqref{intro strengthened angle} below, Theorem \ref{jet corner structure} resolves \eqref{intro C2} into two equations for the coefficients of the corner expansion
\begin{align}\label{intro dot expansion}
\dot\phi_0(r)=c_{\dot\phi_0,\dr_{-1}}r+c_{\dot\phi_0,\dr_{-1}^2}r\ln r+{1\over 2}c_{\dot\phi_0,\dr_{-1}^3}r\ln^2r+\dot\phi_0^{[2]}(r)
\end{align}
of $\dot\phi_0$, namely
\begin{align}\label{intro C2 explicit}
c_{\dot\phi_0,\dr_{-1}^3}=0\qquad\text{ and }\qquad c_{\dot\phi_0,\dr_{-1}^2}={\dot U(0)\over\theta_b(\sin\theta_b)^2},\tag{C2$'$}
\end{align}
where $\dot U(0)$ is determined by $\phi_0$ as $U(0)$ is: it is the horizontal velocity at the corner of the time-differentiated Stokes problem (Theorem \ref{nonlinear initial data}(ii)), and formally the acceleration of the contact point at $t=0$. The second equation is the time derivative of \eqref{intro C1}. The first states that the coefficient of the second-order logarithm $r\ln^2r$, which the nonlinearity creates, vanishes.

The angle $\theta_b$ has to be small: the linear theory of Section \ref{linear existence} requires this without quantifying how small (Remark \ref{smallness remark}). In addition the angle enters through the explicit condition
\begin{align}\label{intro angle}
2-\kappa<\Lambda^*(\theta_b)\qquad\text{ and }\qquad\brac{2-{\kappa\over 2}}\theta_b<{\pi\over 2},
\end{align}
where $\Lambda^*(\theta_b)$ is the half-width of the largest strip $|\Re\lambda|<\Lambda^*(\theta_b)$ on which the characteristic function $\lambda\sin(2\theta_b)+\sin(2\lambda\theta_b)$ of the Stokes wedge has no zero other than $\lambda=0$; see Definition \ref{zero-free strip} and Lemma \ref{H lemma}. The first condition holds up to about $73^\circ$ (Remark \ref{H numerics remark}). The second condition is the more restrictive: it limits the angle to $\theta_b<\pi/(4-\kappa)$, and it comes from our estimate of the boundary nonlinearity after the flattening of the domain (Remark \ref{sharpness of the second angle condition}). The finer statements, (iii)--(v) of the theorem below, read the corner expansions one order further. For them the angle condition is strengthened to
\begin{align}\label{intro strengthened angle}
3-\kappa<\Lambda^*(\theta_b)\qquad\text{ and }\qquad(3-\kappa)\theta_b<{\pi\over 2},
\end{align}
which limits the angle to $\theta_b<\pi/(6-2\kappa)$, about $33^\circ$--$36^\circ$. The data are then asked to have a second-order corner expansion,
\begin{align}\label{intro second-order expansion}
\phi_0^{[2]}(r)=c_{\phi_0,\dr_{-2}}r^2+c_{\phi_0,\dr_{-2}^2}r^2\ln r+{1\over 2}c_{\phi_0,\dr_{-2}^3}r^2\ln^2r+\phi_0^{[3]}(r),
\end{align}
measured by the norm
\begin{align}\label{intro data norm +}
\mathcal D_+(\phi_0):=\mathcal D(\phi_0)+(\sin\theta_b)^2\brac{\sum_{k=1}^3|c_{\phi_0,\dr_{-2}^k}|+\|\phi_0^{[3]}\|_{\partial_1\Omega,{7\over 2}-\kappa,l+{7\over 2}}}.
\end{align}

\begin{theorem}[Main theorem]\label{main theorem}
Let $\kappa\in(1/4,1/2)$ and $l\geq 1$, let $\theta_b\in(0,\pi/2)$ satisfy \eqref{intro angle} and be small enough for the linear theory of Section \ref{linear existence}, and let $X_0$ be an initial free boundary with $X_0'(0)=b^{-1}$ whose profile $\phi_0$ of \eqref{intro profile} satisfies the two compatibility conditions \eqref{intro C1} and \eqref{intro C2}. There are constants $c_*,C_*>0$, depending only on $\theta_b,\kappa$ and $l$, such that the following holds whenever the data norm \eqref{intro data norm} satisfies $\mathcal D\leq c_*$.
\begin{enumerate}
\item (Existence, lifespan) The system \eqref{intro system} has a solution $(X,\mb u,p)$ on the interval $[0,T_*]$, where $T_*=(c_*/\mathcal D)^{2/9}$ for $\mathcal D>0$, and $T_*:=\infty$ for the zero datum, whose solution is the resting wedge. The profile $\phi(r,t)=b^{-1}-X'(r\sin\theta_b,t)$, the transformed velocity and the transformed pressure lie in the space-time class \eqref{ET definition} of Theorem \ref{nonlinear positive-time existence}, with norm at most $C_*(1+T_*)\mathcal D$ when the norms are formed with the weight $e^{-t/T_*}$. Without the weight, the same bound holds up to a factor $e$. The contact angle condition $X'(0,t)=b^{-1}$ holds classically for every $t$, and the contact point moves with the velocity $U(t)=u_x(X(0,t),0,t)$ of the fluid there.
\item (Uniqueness, continuous dependence) The solution is unique in the iteration space $B_{T_*}$ of \eqref{iteration ball}, formed with the weight $e^{-t/T_*}$ and the radius of Theorem \ref{nonlinear positive-time existence}. The profile, the transformed velocity and the transformed pressure depend Lipschitz continuously on the initial profile, uniformly on the common interval $[0,T]$ of existence of two solutions, from the norm $\mathcal D$ to the norm of the class \eqref{ET definition} formed with the weight $e^{-t/T}$, with Lipschitz constant $C_*(1+T)$.
\end{enumerate}
Suppose moreover that $\theta_b$ satisfies \eqref{intro strengthened angle} and that $\phi_0$ has the second-order corner expansion \eqref{intro second-order expansion} with $\mathcal D_+(\phi_0)<\infty$ in \eqref{intro data norm +}. Then (v) holds, and (iii) and (iv) hold for any fixed $\tilde\kappa\in(\kappa/2,\kappa)$ after $c_*$ has been reduced by a factor that depends also on $\tilde\kappa$.
\begin{enumerate}[resume]
\item (Regularity, propagation of the first compatibility condition) The solution has the improved regularity of Proposition \ref{positive-time regularity upgrade} on $[0,T_*]$; the constants of those estimates depend on $T_*$ and $\tilde\kappa$. The profile $\phi(t)$ satisfies \eqref{intro C1} at almost every time, and at these times the contact angle condition differentiated once in time holds pointwise.
\item (The contact point) The contact-point velocity $U(t)$ is continuous in time, and at almost every time we have
\begin{align*}
\mb u&=U\mb e_1+{U\over 2}\mb P_0(r,\theta-\theta_b)+O(r^{2+\tilde\kappa-\kappa}),\\
p&=-c_3(t)-{U\over\theta_b}\ln r+O(r^{1+\tilde\kappa-\kappa}),\\
\phi&=c(t)r+{U\over\theta_b(\sin\theta_b)^2}r\ln r+O(r^{2+\tilde\kappa-\kappa}),
\end{align*}
where $c(t)$ is a function of time alone, $\mb e_r,\mb e_\theta$ is the polar basis, and
\begin{align*}
\mb P_0(r,\theta-\theta_b)&=r\brac{\brac{{\theta-\theta_b\over\theta_b}-{\sin(2(\theta-\theta_b))\over\sin(2\theta_b)}}\mb e_\theta-\brac{{1\over 2\theta_b}-{\cos(2(\theta-\theta_b))\over\sin(2\theta_b)}}\mb e_r},\\
c_3(t)&=(\sin\theta_b)^2c(t)+\brac{{1\over 2\theta_b}+{1\over\sin(2\theta_b)}}U(t).
\end{align*}
\item (The compatible set) Near $\phi_0=0$, the initial profiles satisfying \eqref{intro C1} form a Lipschitz graph of codimension one in the space of data with $\mathcal D<\infty$; for this, \eqref{intro angle} suffices. Those satisfying \eqref{intro C1} and \eqref{intro C2} form a Lipschitz graph of codimension three in the space of data with a second-order corner expansion. Both graphs contain $\phi_0=0$.
\end{enumerate}
\end{theorem}

The statements are proved in Section \ref{nonlinear existence}: (i) is Theorem \ref{nonlinear positive-time existence} with Corollary \ref{polynomial lifespan}, returned to the moving domain by Corollary \ref{moving domain} and Theorem \ref{the contact point}(iii); (ii) is the uniqueness clause of Theorem \ref{nonlinear positive-time existence} and Proposition \ref{Lipschitz positive}; (iii) is Proposition \ref{positive-time regularity upgrade} and Corollary \ref{C1 along the flow}; (iv) is Theorem \ref{the contact point} and Remark \ref{contact point physics}, which also gives the dimensional form of the expansion of the free boundary; and (v) is Theorem \ref{compatible set}. The initial jet is constructed in Theorem \ref{nonlinear initial data}, which also explains the two compatibility conditions. The parameter $\tilde\kappa$ in (iii) and (iv) is introduced before Proposition \ref{time for radial trade}. It measures how close the exponents of the remainders in (iv) come to $2$ and $1$, which they do not reach. Under the hypotheses of (iii) and (iv) the remainder $\phi_0^{[2]}$ of the data is itself of the order $r^2\ln^2r$, by \eqref{intro second-order expansion}.

By Theorem \ref{main theorem}(iv), the velocity is continuous at the contact point with no logarithm, the pressure has the logarithmic singularity of a moving contact line with slip, and the slope of the free surface bends logarithmically, an apparent lowering of the angle above the contact point for a receding contact line and a raising for an advancing one.

The expansions of Theorem \ref{main theorem}(iv) are explicit up to the two functions $U$ and $c$ of time. The equations near the contact point do not determine these two functions: apart from the initial value $c(0)=\smash[b]{c_{\phi_0,\dr_{-1}}}$ in \eqref{intro expansion}, both are determined by the problem on the whole wedge. To first order in $\phi_0$, for instance, $U(0)$ is given by \eqref{Lambda formula}, a combination of an integral of $\phi_0$ over the whole free surface and an integral of the velocity over the whole substrate.

In the matched asymptotics \cite{Voinov,Hocking,Cox} that lead to the Cox--Voinov relation, the constants that the expansion in one region leaves free are determined by matching with the expansion in the adjacent region. Here the role of the matching is played by coercivity. After a Laplace transform in time, the Lax--Milgram argument of Section \ref{linear existence} gives exactly one variational solution of the linearised problem in a weighted space of velocity fields that vanish at infinity but need not vanish at the contact point. At infinity the slope of the free surface returns to $b^{-1}$, by \eqref{intro system}. The Cox--Voinov relation, which travelling waves of the thin-film equation satisfy at infinity \cite{GiacomelliGnannOtto,GnannWisse}, can therefore hold here only at intermediate distances from the contact point, a regime that we do not characterise.

\subsection{Comments and related work}\label{comments and related work}

Three features of the result deserve comment. First, the compatibility conditions are imposed at $t=0$ by the class in which we work, whose norms control the time derivative of the profile up to $t=0$, and by our construction of solutions in it. Without \eqref{intro C1} the derived profile $\partial_t\phi|_{t=0}$ has a nonzero value at the contact point, so that the differentiated contact angle condition fails at $t=0$ (Remark \ref{initial data closing remark}). The condition \eqref{intro C2} is exactly what is needed for the second time derivative of the profile to be an admissible datum of that construction (Theorem \ref{linear solution full regularity}). The bend described by \eqref{intro C1} is produced by the flow: the solution satisfies \eqref{intro C1} at almost every time, by part (iii) of the theorem.

A different class carrying time weights in the norms, so that the estimates are allowed to degenerate as $t\to 0$, should admit every small datum in the norm $\mathcal D$, whether or not it satisfies the two conditions. The bend would then form only for $t>0$. For the linear problem such data are already admitted: Theorem \ref{linear solution in time}, whose class controls no time derivative, asks of the initial profile only $\phi_0(0)=0$ and a weighted norm. Every small profile with compact support away from the contact point would then be admissible; we have not pursued this.

For solutions regular in time up to $t=0$, such as ours, the price of the fixed angle is a restriction of the data, near $\phi_0=0$, to Lipschitz graphs, of codimension one for \eqref{intro C1} and of codimension three for both conditions. Both graphs are large. Adding to a small smooth profile with compact support away from the contact point a suitable multiple of a fixed profile, equal to $r\ln r$ near the contact point, gives a profile satisfying \eqref{intro C1}. Adding to such a small smooth profile a suitable linear combination of three fixed profiles gives one satisfying \eqref{intro C1} and \eqref{intro C2}. A profile with compact support away from the contact point satisfies \eqref{intro C1} exactly when it leaves the contact point initially at rest. Every small smooth profile with compact support away from the contact point can be made to satisfy \eqref{intro C1} by adding a suitable multiple of a fixed profile that is itself smooth with compact support away from the contact point (Remark \ref{compatible data examples}).

Second, the lifespan is a power of the inverse size of the data, because the exponential weight $e^{-\varsigma t}$ of the linear theory can be chosen in terms of the length of the time interval, and the constants of the theory depend polynomially on $\varsigma^{-1}$ (Remark \ref{lifespan remark}). This weight is a device of the $L^2$-in-time framework and not a feature of the problem: the linear evolution is dissipative.

Third, the angle. The explicit condition \eqref{intro angle} identifies the two mechanisms that limit the method: the resonance of the harmonic extension used to flatten the domain, at the exponent $\pi/(2\theta_b)$, and the first Moffatt-type eigenvalue \cite{Moffatt} of the Stokes wedge with a free-slip substrate and a stress-free surface, at $\Re\lambda\theta_b\approx 2.1$. But the linear theory closes its coercivity estimate by absorbing terms that carry positive powers of $\theta_b$, and so asks for a smallness of the angle which we do not quantify. The range of angles actually covered is therefore not determined by \eqref{intro angle}.

The mathematical literature on moving contact lines is small in comparison with the physical one, and most of it is recent. The problem entered mathematical fluid dynamics with Pukhnachov and Solonnikov \cite{PukhnachovSolonnikov}, who made precise the incompatibility of the no-slip condition with a moving contact line and formulated the slip models on which the later mathematical work rests. Free-boundary problems for the Navier--Stokes equations with contact points and lines were then studied by Solonnikov \cite{Solonnikov}. Stationary flows with prescribed contact angles were studied by Socolowsky \cite{Socolowsky}, in two-dimensional coating configurations, and by Jin \cite{Jin}, in three dimensions at the angle $\pi/2$. A well-posed model with a contact angle of $\pi/2$ was constructed by Schweizer \cite{Schweizer}, through a reflection across the substrate which that angle alone permits, and which recurs in recent well-posedness theories for other free-boundary problems meeting a fixed boundary at right angles \cite{AbelsRaucheckerWilke,GarckeRauchecker}. For the inviscid counterpart --- the water waves equations with surface tension in a corner domain, with a dynamic contact angle --- a priori estimates and local well-posedness for small angles are due to Ming and Wang \cite{MingWangI,MingWangII}. For a floating object, Lannes \cite{Lannes} formulated the water waves problem without surface tension, with conditions at the contact line between the free surface and the object, and derived the corresponding shallow water models. Bresch and Desjardins \cite{BreschDesjardins} constructed global weak solutions of a viscous shallow water model. The depth of these solutions may vanish, and no condition is imposed where it does.

The nearest antecedent of the present work is the analysis of Darcy's flow in a wedge with a prescribed contact angle by the second and the fourth author of this paper. In \cite{KnupferMasmoudiCMP} they proved well-posedness and uniform bounds for the linearised evolution, a nonlocal equation of third order on an infinite wedge. In \cite{KnupferMasmoudi} they proved well-posedness of the nonlinear problem, together with the lubrication approximation, through a corner analysis by the Mellin transform that the present paper follows in spirit. The Stokes problem in the same geometry is the natural next step, and a substantial one: the harmonic pressure of Darcy's law is replaced by the full Stokes system, with its own family of corner exponents. It is in this model that the logarithms of the corner expansion and the compatibility conditions first appear.

A parallel line of work makes the contact angle dynamic. Ren and E proposed conditions at the contact line in which the contact-line velocity is proportional to the deviation of the angle from its equilibrium value \cite{RenE}, and later derived them from thermodynamic principles \cite{RenECMS}. In the lubrication approximation, the corresponding conditions follow from a formal gradient-flow structure with a dissipation at the contact line \cite{GiacomelliGnannPeschka}. With a condition of this type, Zheng and Tice \cite{ZhengTice} proved local well-posedness for the two-dimensional Stokes flow in a bounded vessel whose free surface meets the walls at contact points, Guo and Tice \cite{GuoTice,GuoTiceNS} the stability of the equilibrium there for Stokes and Navier--Stokes flow, and Tice and Wu \cite{TiceWu} the dynamics and stability of sessile drops. These analyses exploit the dissipation that the relaxation of the angle contributes at the contact points.

Ren and E asked whether the relaxation of the angle removes the singularity that remains with slip when the angle is fixed \cite{RenE}. E described its removal as a natural expectation and the question as open \cite[Section~5.3]{EMultiscale}. The relaxation is, however, constrained kinematically: for velocity fields that are continuously differentiable up to the contact line, Fricke, K\"ohne and Bothe \cite{FrickeKohneBothe} derived the evolution equation that transport imposes on the angle. They conclude that with Navier slip and any law for the contact angle under which the contact line dissipates energy, the angle can relax to its equilibrium value only if the solution is singular at the contact line, which answers the question in the negative. This conclusion depends on the model: by the same analysis it need not hold if the surface tension has a gradient at the contact line. For Stokes flow, Guo and Tice \cite{GuoTice} measure the velocity and the pressure in Sobolev norms weighted by a positive power of the distance to the contact points. For Navier--Stokes flow \cite{GuoTiceNS} the norms are based instead on Lebesgue spaces of exponent less than $2$. In both cases the norms are weak enough near the contact points that a pressure with a logarithmic singularity there has finite norm. A bound in these norms therefore does not show whether the singularity is present.

The present model is the formal limit of those in which the relaxation is instantaneous. The angle is fixed at its equilibrium value, and the model carries one constitutive parameter fewer: nothing is prescribed at the contact point beyond the angle itself. Prescribing only the angle is the boundary condition of the classical hydrodynamic theory of wetting \cite{Voinov,Cox}, on the picture that the angle equilibrates at molecular rates, fast compared with the rates of the flow \cite{deGennes,Bonn}. The dissipation that the relaxation provided vanishes with it. The fixed angle is in this sense the more fundamental condition, and the analytically harder one: no damping is available at the contact point, and well-posedness has to come out of the corner structure itself. One conclusion of this paper is that it does --- the contact point dissipation is not needed for well-posedness. What stands in its place, for solutions regular in time up to $t=0$, is the pair of compatibility conditions, a damping mechanism traded for a constraint on the data. The compatibility conditions tie the geometry of the initial interface to the motion of its contact point, and they cut the admissible data down to the Lipschitz graphs of finite codimension of Theorem \ref{main theorem}(v).

The geometry also separates the two settings. The analyses \cite{ZhengTice,GuoTice,GuoTiceNS,TiceWu} are set in a bounded vessel or drop. Ours is set in the unbounded wedge, which has no length of its own: the slip length is the only length of the problem, and no estimate depends on a further one, such as the size of a vessel or of a drop. In return every estimate has to be uniform across scales, and the corner analysis is the centre of the argument rather than a boundary effect.

A third line treats the contact line after the lubrication approximation, which collapses the flow in the wedge onto a degenerate parabolic equation of fourth order for the height of the film, the thin-film equation. There the theory is most developed. Well-posedness and regularity at the moving contact line are known for zero and nonzero contact angles and for various mobilities, in two and three dimensions; see for instance \cite{Otto,GiacomelliKnupferOtto,Knupfer,GiacomelliGnannKnupferOtto,John,KnupferGeneral,GnannPetrache,GnannWisseJFA}. The Cox--Voinov relation itself has been proved for travelling waves, in the complete wetting regime \cite{GiacomelliGnannOtto} and in the partial wetting regime \cite{GnannWisse}. Tanner's law \cite{Tanner} for the spreading of a drop is studied rigorously in \cite{GiacomelliOttoDroplet}, as an intermediate scaling law in time, and in \cite{DelgadinoMellet}, in the limit of small slip.

The lubrication approximation has been justified for Darcy's flow, with zero and with prescribed contact angle \cite{GiacomelliOtto,KnupferMasmoudi}, and for layers of Stokes flow without contact points \cite{GuntherProkert}. In that setting, free-boundary problems for Stokes flow driven by surface tension are by now well understood; see Prokert \cite{Prokert}, G\"unther and Prokert \cite{GuntherProkertDegenerate} and Matioc and Prokert \cite{MatiocProkert}. For Stokes flow with slip and a moving contact line the rigorous justification is, so far as we know, open, and an analysis of the flow at the scales at which the approximation breaks down, such as the present one, is a natural step towards it.

The corner itself is handled with the tools of the elliptic theory in domains with corners: weighted Sobolev spaces with powers of the distance to the corner as weights, the Mellin transform, and corner expansions read off from the poles of the resolvent. This theory goes back to Kondrat'ev \cite{Kondratev} and is presented in the books of Grisvard \cite{Grisvard} and of Kozlov, Maz'ya and Rossmann \cite{KozlovMazyaRossmann}, who devote a second book \cite{KozlovMazyaRossmannPencils} to the operator pencils of such problems, those of the Stokes system among them. Two circumstances prevent a direct appeal to the classical theory. First, the wedge is unbounded, with no scale to localise at, so that the analysis cannot be reduced to a neighbourhood of the corner plus a smooth remainder. Second, the problem is an evolution. So the corner expansion is not a fixed list of singular functions but a structure that has to be constructed at the initial time and then propagated along the flow: coefficients coupled to the motion of the contact point, and the two compatibility conditions. 

The Stokes system on the fixed wedge itself has been studied in the $L^p$ setting, with perfect slip on both boundaries: an $H^\infty$-calculus and maximal regularity in weighted spaces of Kondrat'ev type are due to Maier and Saal \cite{MaierSaal}, and optimal Sobolev regularity to K\"ohne, Saal and Westermann \cite{KohneSaalWestermann}. For the stationary Stokes system on a fixed wedge with Navier-slip boundary conditions, with both boundaries prescribed and no free surface, well-posedness and regularity up to the corner in weighted Sobolev spaces have recently been obtained by Bravin, Roodenburg and the first, second, fourth and fifth authors of this paper \cite{Bravin}. The linear theory of Section \ref{linear existence} concerns instead the problem with the Laplace parameter of the time transform. That problem has to be controlled uniformly in the parameter and coupled to the free boundary through the stress balance.

Against this background the contribution of the paper can be stated plainly. Theorem \ref{main theorem} appears to be the first well-posedness result for the Stokes equations with surface tension, Navier slip, a moving contact point and a contact angle fixed at a value other than $\pi/2$. It is also, to our knowledge, the first well-posedness theory of a moving contact line in which the logarithmic pressure singularity of the moving contact point is characterised within the solutions themselves. This singularity is what survives of the Huh--Scriven singularity once the fluid slips. Our solution class carries it explicitly, in corner expansions. The norms of their remainders, bounded in part (iii) of the theorem, are strong enough near the contact point to determine the singularity. Part (iv) exhibits the logarithm with its exact coefficient, the speed of the contact point divided by the angle. So the logarithm is present exactly when the contact point moves. Well-posedness holds in the presence of the logarithm; the two compatibility conditions are precisely what its propagation costs for solutions regular in time up to $t=0$.

The expansions are quantitative: in dimensional variables the coefficient of the pressure logarithm is the viscous stress at the slip scale divided by the angle, and the logarithmic bending of the free surface carries the capillary number of the contact line. They are statements about the free-boundary problem itself, not about a wedge with prescribed boundaries, and they are the counterpart, below the slip length, of the Cox--Voinov asymptotics above it (Remark \ref{contact point physics}). Throughout, the fluid above the free surface is passive.

\subsection{Method and structure of the paper}

The domain is flattened onto the reference wedge $\Omega$ by the harmonic extension of the deviation of the free boundary. This is not a bare change of variables, because the free boundary is an unknown. Section \ref{setup} sets up the transformation, the weighted spaces and the Mellin calculus. The Mellin transform acts on the radial variable and leaves the angular one. So it turns an elliptic problem on the wedge whose terms all have the same scaling, in the sense of Definition \ref{scaling definition}, into a two-point boundary value problem in the angle, which for the Laplace and the Stokes operator has an explicit solution formula. The elliptic problems solved in this way in Section \ref{setup} are auxiliary: the harmonic problems behind the flattening, the Helmholtz projection and the pressure.

Section \ref{linear existence} solves the linearised problem. After a Laplace transform in time it is a stationary Stokes problem with the Laplace parameter $s$. In this problem the Navier slip condition and the stress balance on the free surface each combine two terms of different scaling: the velocity and its derivative in the first, the stress and the surface tension, which carries the factor $s^{-1}$, in the second. A term of lower scaling shifts the line of the Mellin transform, so that the transform no longer reduces the problem to one in the angle and gives no solution formula. We solve the problem instead by the Lax--Milgram theorem in a weighted space $H_s^l$ with two weights. The smallness of the angle enters through the coercivity of the form, and the contact angle condition is encoded in the form rather than imposed.

Section \ref{corner expansions} is the corner analysis. Without a solution formula for the full problem, the corner expansions cannot be read off in one step. We treat the terms of lower scaling as data. The remaining problem has one scaling and an explicit Mellin representation. Since the data are built from the solution itself, the representation raises by one the scaling at which the solution is controlled, and it has to be iterated. Shifting the line of integration of the Mellin representation across its poles then produces the corner expansions rung by rung, together with the exact form of the compatibility conditions. The same representation, applied to time derivatives, trades time derivatives for radial regularity at positive times.

Section \ref{nonlinear existence} solves the nonlinear problem by two contractions. The first, at the initial time, constructs the compatible initial jet and identifies the compatibility conditions; the second works at positive times, in the space-time norms of the linear theory. Section \ref{nonlinear existence} then recovers the full regularity, propagates the first compatibility condition along the flow, proves the Lipschitz dependence and the codimension statement, and returns to the moving domain. Three appendices collect the identities of polar coordinates and of the Mellin transform, the weighted inequalities used throughout, and the resonance factors.

\section{The problem and its reformulation}\label{setup}

We study a two-dimensional Stokes flow occupying a wedge-like region with one solid boundary, on which a Navier slip condition is imposed, and one free boundary driven by surface tension, the two meeting at a contact point that is free to move. The purpose of this section is to set up the problem, to record the function spaces and the Mellin calculus that we shall work with, and to transform the moving domain into a fixed reference wedge. The transformation is not a bare change of variables: because the free boundary is one of the unknowns, we flatten the domain using a harmonic extension of the deviation of the profile from the exact wedge. This requires the elliptic theory of Subsection \ref{wedge elliptic systems}, so the transformation is postponed to Subsection \ref{fixed domain transformation}.

\begin{definition}[Notation]\label{notation}\leavevmode
\begin{itemize}
\item We denote by $\mb e_1,\mb e_2$ the standard basis of $\R^2$, and by $\mb e_r,\mb e_\theta$ the polar basis of $\R^2$. Subscripts denote components (not derivatives), e.g. $\mb u=u_x\mb e_1+u_y\mb e_2$, and $\mb u^\perp=(-u_y,u_x)$ and $\grad^\perp=(-\partial_y,\partial_x)$ are the rotations by $\pi/2$.
\item For functions $X(y,t)$ of exactly one space variable and one time variable, we use the shorthand notation $X'=\partial_y X$ and $\dot X=\partial_t X$. On functions of the angle $\theta$ alone, such as those arising on the Mellin side, the prime denotes $\partial_\theta$ instead. In statements comparing two solutions or two data a prime marks the second one, as in $(\phi',\mb u',p')$ or $\phi_0'$.
\item In the systems below every equation holds for all times $t$ considered, respectively for all $s$ on the Laplace side, unless it names its time, as the initial profile does; the qualifier after an equation says where in space it holds. A function satisfies a system strongly if its equations hold almost everywhere on the sets named by their qualifiers, the boundary conditions in the sense of traces.
\item For sesquilinear forms $B[\mb u,\mb v]$ we write $B[\mb u]$ to mean $B[\mb u,\mb u]$.
\item Objects collected in braces $\{\ph\}$ inside an estimate are a choice: the estimate holds with any one of them in place of the brace. Braces are chosen independently of one another, and independently on the two sides of the inequality; a power $\{\ph\}^k$ counts as $k$ separate braces, so that each of its $k$ factors may be chosen on its own. Thus $\{a,b\}^2$ stands for any of $a^2$, $ab$ and $b^2$, and an estimate $\{a,b\}\lesssim\{c,d\}$ asserts all four of $a\lesssim c$, $a\lesssim d$, $b\lesssim c$ and $b\lesssim d$. In a chain of estimates the choices in the intermediate members may depend on those in the first and the last member.
\item We write LHS and RHS for the left and right hand sides of an equation or an inequality.
\end{itemize}
\end{definition}

We will study the following system, in which the flow is in $\R^2$, equivalent to a flow in $\R^3$ with no dependence on the $z$ direction.

\begin{system}[Wedge-like Stokes flow with a free surface]\label{dimensional system}
Denote the symmetric gradient $\S(\mb u)={1\over 2}((\grad\mb u)+(\grad\mb u)^T)$ and the stress tensor $\mathbb T[\mb u,p]=2\mu\S(\mb u)-pI$, where $\mb u$ is the velocity field, $p$ the pressure and $\mu$ the dynamic viscosity. We consider the incompressible Stokes equation
\begin{align*}
\underbrace{\mu\lpc\mb u-\grad p}_{=\grad\cdot\mathbb T}&=\mb 0\\
\grad\cdot\mb u&=0
\end{align*}
on a wedge-like domain with free upper boundary $\Omega_t=\{(x,y)\in\R\times(0,\infty):x>X(y,t)\}$ where $X(\ph,t):[0,\infty)\to\R$ is a continuous function representing the extent of the fluid. We denote
\begin{align*}
\partial_0\Omega_t&=\{(x,0)\in\R^2:x>X(0,t)\}\\
\partial_1\Omega_t&=\{(x,y)\in\R\times(0,\infty): x = X(y,t)\}
\end{align*}
to be the lower and upper boundary of $\Omega_t$ respectively. We have the following boundary conditions
\begin{alignat*}{3}
u_y&=0\qquad\qquad&\text{on}\quad\partial_0\Omega_t \tag{no penetration}\\
u_x-\beta\partial_yu_x&=0\qquad\qquad&\text{on}\quad\partial_0\Omega_t\tag{Navier slip}\\
\mathbb T\mb n&=-\gamma X''(1+X'^2)^{-{3\over 2}}\mb n\qquad\qquad&\text{on}\quad\partial_1\Omega_t\tag{surface tension}\\
(\partial_t+u_y\partial_y)X&=u_x\qquad\qquad&\text{on}\quad\partial_1\Omega_t\tag{profile evolution}\\
X'(0,t)&=b^{-1}\tag{contact angle}\\
\lim_{y\to\infty}X'(y,t)&=b^{-1}\tag{slope at infinity}\\
X(\ph,0)&=X_0\tag{initial profile},
\end{alignat*}
where $\mb n$ denotes the outward unit normal to the surface and $b$ is some given positive constant that is the $\tan$ of the contact angle $\theta_b=\arctan b$.
\end{system}

The parameters $\mu$, $\gamma$ and $\beta$ can be scaled out. Measuring lengths in units of the slip length $\beta$, velocities in units of $\gamma/\mu$, times in units of $\beta\mu/\gamma$ and pressures in units of $\gamma/\beta$, that is, writing
\begin{alignat*}{4}
x &= \beta \tilde{x},\qquad\qquad & y &= \beta \tilde{y},\qquad\qquad & X &= \beta \tilde{X}\\
t &= {\beta\mu\over\gamma} \tilde{t},\qquad\qquad & \mb u &= {\gamma\over\mu} \tilde{\mb u},\qquad\qquad & p &= {\gamma\over\beta} \tilde{p},
\end{alignat*}
every equation of System \ref{dimensional system} loses its constants, since $X'$ is already dimensionless while $X''=\beta^{-1}\tilde X''$ and $\mathbb T=(\gamma/\beta)\tilde{\mathbb T}$ with $\tilde{\mathbb T}=2\S(\tilde{\mb u})-\tilde pI$. The only physical parameter left is $b$. Dropping the tildes, we may therefore consider the following system.

\begin{system}[Wedge-like Stokes flow with a free surface, dimensionless form]\label{dimensionless system}
System \ref{dimensional system} with $\mu=\gamma=\beta=1$; explicitly, with $\Omega_t$, $\partial_0\Omega_t$ and $\partial_1\Omega_t$ as above,
\begin{alignat*}{3}
\lpc\mb u-\grad p=\mb 0,\qquad\grad\cdot\mb u&=0\qquad&&\text{on}\quad\Omega_t\\
u_y=0,\qquad u_x-\partial_yu_x&=0\qquad&&\text{on}\quad\partial_0\Omega_t\\
\mathbb T[\mb u,p]\mb n&=-X''(1+X'^2)^{-{3\over 2}}\mb n\qquad&&\text{on}\quad\partial_1\Omega_t\\
(\partial_t+u_y\partial_y)X&=u_x\qquad&&\text{on}\quad\partial_1\Omega_t\\
X'(0,t)=b^{-1},\qquad\lim_{y\to\infty}X'(y,t)&=b^{-1},\qquad X(\ph,0)=X_0,
\end{alignat*}
where $\mathbb T[\mb u,p]=2\S(\mb u)-pI$ and the five lines are, in order, the Stokes equation and incompressibility, no penetration and Navier slip, surface tension, the profile evolution equation, and the contact angle, the slope at infinity and the initial profile.
\end{system}

We will transform the evolving domain $\Omega_t$ into a fixed wedge reference domain $\Omega$:
\begin{subequations}\label{fixed wedge reference domain}
\begin{align}
\Omega&:=\{(x,y)\in\R^2:0<y<bx\}
=\{(r\cos\theta,r\sin\theta)\in\R^2:r\in(0,\infty),\theta\in(0,\theta_b)\}\\
\partial_0\Omega&:=\{(x,0)\in\R^2:x>0\}
=\{(r,0)\in\R^2:r\in(0,\infty)\}\\
\partial_1\Omega&:=\{(x,bx)\in\R^2:x>0\}
=\{(r\cos\theta_b,r\sin\theta_b)\in\R^2:r\in(0,\infty)\}
\end{align}
\end{subequations}
Here $\overline\Omega=\Omega\cup\partial_0\Omega\cup\partial_1\Omega\cup\{\mb 0\}$ denotes the closed wedge, and likewise $\overline\Omega_t$ the closure of $\Omega_t$. For $0<\epsilon<R$ we write $\Omega_{\epsilon,R}:=\Omega\cap\{\epsilon<r<R\}$. Along $\partial_1\Omega$ we have $y=r\sin\theta_b$, where $\sin\theta_b=b/\sqrt{1+b^2}=(1+b^{-2})^{-1/2}$ since $b=\tan\theta_b$. A function $f$ of the height $y$ is regarded along $\partial_1\Omega$ as a function of the distance $r$ to the corner through this relation, its prime being the derivative in $y$ as in Definition \ref{notation}, so that
\begin{align}\label{r and y on the free boundary}
y=r\sin\theta_b,\qquad\qquad\partial_rf=(\sin\theta_b)f'\qquad\text{on}\quad\partial_1\Omega.
\end{align}

We will transform the problem into this reference domain, but in order to do so we first have to define the norms and establish the elliptic results on which the transform is based.

\subsection{Norms and Mellin transform}\label{norms and Mellin transform}

Radial derivatives enter everywhere in the scaled form $r\partial_r$, and it is convenient to fix a symbol for them first.

\begin{definition}[Scaled radial derivatives]\label{scaled radial derivatives}
We introduce the notation
\begin{align*}
\dr&=r\partial_r\\
\dr_{a}&=r\partial_r+a.
\end{align*}
When we take Mellin transform $\dr$ is transformed into $\lambda$; we denote analogously
\begin{align*}
\lambda_a=\lambda+a,
\end{align*}
and, on a Mellin line avoiding $\lambda=-a$, we write $\dr_a^{-1}$ for the Mellin multiplier $\lambda_a^{-1}$. We also write $\<\lambda\>=(1+|\lambda|^2)^{1/2}$ for $\lambda\in\C$, and $\<\dr\>^c$ for the Mellin multiplier $\<\lambda\>^c$.
\end{definition}

Note that $\<\lambda\>\leq 1+|\lambda|\leq\sqrt 2\<\lambda\>$.

The norms we shall work with are the following weighted Sobolev norms, where $a$ is the weight, $j$ the number of radial derivatives $\dr$ and $m$ the number of full derivatives $r\grad$. We call such numbers of derivatives counts. The parameter $\epsilon>0$ carried by the norms of Definition \ref{weighted Sobolev norms} weighs the successive derivative counts against each other; it will be fixed small in Section \ref{linear existence} and suppressed from the notation whenever it equals $1$.

\begin{definition}[Weighted Sobolev norms]\label{weighted Sobolev norms}
We introduce the notation
\begin{align}
\|\mb u\|_{O}^2&:=\|\mb u\|_{L^2(O)}^2\label{plain norm}\\
\|\mb u\|_{O,a}^2&:=\|r^{-a}\mb u\|_{O}^2=\|r^{-a}\mb u\|_{L^2(O)}^2\label{weighted norm}\\
\|\mb u\|_{O,a,j;\epsilon}^2&:=\sum_{j'=0}^j\|(\epsilon\dr)^{j'}\mb u\|_{O,a}^2=\sum_{j'=0}^j\|r^{-a}(\epsilon\dr)^{j'}\mb u\|_{L^2(O)}^2\label{radial norm}\\
\|\mb u\|_{O,a,j}^2&:=\|\mb u\|_{O,a,j;1}^2\nonumber\\
\|\mb u\|_{O,a,j,m;\epsilon}^2&:=\sum_{m'=0}^m\|(\epsilon r\grad)^{m'}\mb u\|_{O,a,j;\epsilon}^2=\sum_{j'=0}^j\sum_{m'=0}^m\|r^{-a}(\epsilon\dr)^{j'}(\epsilon r\grad)^{m'}\mb u\|_{L^2(O)}^2\label{full norm}\\
\|\mb u\|_{O,a,j,m}^2&:=\|\mb u\|_{O,a,j,m;1}^2,\nonumber
\end{align}
where $O$ is a placeholder and $(r\grad)^{m'}$ is a tensor product containing all derivatives of order $m'$.
\end{definition}

\begin{remark}[Integration by parts for $\dr$]\label{dr by parts}
Note that with integration by parts we have
\begin{align*}
\<u,r^a\dr v\>_{L^2(\partial\Omega)}&=\<u,r^{a+1}\partial_r v\>_{L^2(\partial\Omega)}=-\<(r^{a+1}\partial_r+(a+1)r^a)u,v\>_{L^2(\partial\Omega)}=-\<\dr_{a+1}u,r^av\>_{L^2(\partial\Omega)}\\
\<u,r^a\dr v\>_{L^2(\Omega)}&=\<u,r^{a+1}\partial_r v\>_{L^2(\Omega)}=-\<(r^{a+1}\partial_r+(a+2)r^a)u,v\>_{L^2(\Omega)}=-\<\dr_{a+2}u,r^av\>_{L^2(\Omega)},
\end{align*}
where the extra 1 in the second line comes from the $r$ in the integration measure $r\d r\d\theta$ for $\Omega$. The computation is carried out for $u$ and $v$ smooth and compactly supported away from the corner and from infinity, so that no boundary terms survive; they extend to more general $u$ and $v$ by approximation, provided the boundary terms $\bar ur^{a+1}v$, respectively $\bar ur^{a+2}v$, vanish at $r=0$ and as $r\to\infty$. Without such a restriction they are only formal.
\end{remark}

Detailed properties of the Mellin transform can be found for instance in Kozlov, Maz'ya and Rossmann \cite[Chapter~1]{KozlovMazyaRossmann}. Here is a brief summary.

\begin{definition}[Mellin transform]\label{Mellin transform}
We define the Mellin transform of a measurable $f:(0,\infty)\to\C$ to be
\begin{align*}
\hat f(\lambda):={1\over\sqrt{2\pi}}\int_0^\infty r^{-\lambda}f(r){\d r\over r}
\end{align*}
for $\lambda$ on the strip $\Re\lambda\in[\gamma_1,\gamma_2]$ of $\C$  (if it exists) for which the integral converges absolutely.
\end{definition}

$f$ can be recovered via the inverse Mellin transform
\begin{align*}
f(r)={1\over\sqrt{2\pi}}\int_{\Re\lambda=\gamma}r^\lambda\hat f(\lambda)\d\Im\lambda,
\end{align*}
where $\gamma\in[\gamma_1,\gamma_2]$, provided the integral converges absolutely.

\begin{lemma}[Properties of Mellin transform]\label{Mellin properties}
Whenever the following quantities are defined, we have
\begin{enumerate}
\item $(r^{-a}f)^\wedge(\lambda)=\hat f(\lambda+a)$
\item $(r\partial_r f)^\wedge(\lambda)=\lambda\hat f(\lambda)$
\item (Plancherel’s identity)
\begin{align}\label{Plancherel's identity}
\int_0^\infty r^{-2a}\overline{f(r)}g(r){\d r\over r}=\int_{\Re\lambda=a}\overline{\hat f(\lambda)}\hat g(\lambda)\d\Im\lambda
\end{align}
\end{enumerate}
\end{lemma}

Like the Fourier transform, the Mellin transform can be extended into an $L^2$ theory. In particular, there exists an isomorphism $L^2((0,\infty),r^{-2a-1}\d r)\to L^2(\{\lambda\in\C:\Re\lambda=a\},\d\Im\lambda)$ with $f\mapsto\hat f$.

\begin{definition}[Fractional and negative radial derivatives]\label{fractional counts}
For $c\in\R\setminus\N_0$, we shall denote
\begin{align*}
\|r^{-a-{1\over 2}}(r\partial_r)^cu\|_{L^2_r([0,\infty))}^2:=\int_{\Re\lambda=a}|\lambda|^{2c}|\hat u(\lambda)|^2\d\Im\lambda.
\end{align*}
For our purpose, it is sufficient to consider the LHS as a notation for the RHS. For $j\in(-1,\infty)\setminus\N_0$ and $m\in\N_0$, we denote
\begin{align*}
\|\mb u\|_{O,a,j,m;\epsilon}^2&=\|\mb u\|_{O,a,\lfloor j\rfloor,m;\epsilon}^2+\|(\epsilon\dr)^j\mb u\|_{O,a,0,m;\epsilon}^2,
\end{align*}
with $\|\mb u\|_{O,a,-1,m;\epsilon}:=0$.
\end{definition}

By Definition \ref{fractional counts}, a norm with the aggregate count $-1$ is simply absent, and one with a count in $(-1,0)$ consists of its fractional term alone.

\begin{remark}[The weight and the Mellin line]\label{weight and Mellin line}
Write $\sigma=1$ for $O=\Omega$ and $\sigma={1\over 2}$ for $O\in\{\partial_0\Omega,\partial_1\Omega\}$, the difference coming from the factor $r$ in the measure $r\d r\d\theta$ of $\Omega$. Then the weight $a$ of Definition \ref{weighted Sobolev norms} corresponds to the Mellin line
\begin{align}\label{weight and line}
\Re\lambda=a-\sigma.
\end{align}
Indeed, by \eqref{Plancherel's identity} and Lemma \ref{Mellin properties}(ii),
\begin{align*}
\|\dr^{j'}u\|_{\partial_\bullet\Omega,a}^2=\int_{\Re\lambda=a-{1\over 2}}|\lambda|^{2j'}|\hat u|^2\d\Im\lambda,\qquad\qquad
\|\dr^{j'}u\|_{\Omega,a}^2=\int_{\Re\lambda=a-1}\int_0^{\theta_b}|\lambda|^{2j'}|\hat u|^2\d\theta\d\Im\lambda
\end{align*}
for $j'\in\N_0$, and Definition \ref{fractional counts} reads the LHSs by the RHSs for fractional $j'$. Every norm of Definition \ref{weighted Sobolev norms} is a sum of such terms, so all of them are computed on the line \eqref{weight and line}; we use this below without further comment.
\end{remark}

\begin{definition}[Scaling]\label{scaling definition}
The dilations $\mb x\mapsto\omega\mb x$, $\omega>0$, map $\Omega$, $\partial_0\Omega$ and $\partial_1\Omega$ onto themselves, and we say that a quantity built from an unknown $u$ has scaling $a$ if it transforms under them as $r^{-a}u$ does. Accordingly, the norm $\|\ph\|_{O,a,j,m}$, which controls $r^{-a}u$ together with derivatives that leave the scaling unchanged, measures $u$ at the scaling $a$, and control at a higher, respectively lower, scaling means control in a norm with a larger, respectively smaller, weight $a$.
\end{definition}

By Definition \ref{scaling definition}, $r^{-1}$, $\partial_r$ and $\grad$ raise the scaling by one, $r$ lowers it by one, and $\dr$, $r\grad$ and $\partial_\theta$ leave it unchanged; a term of an equation has a higher scaling than another if it carries more powers of $r^{-1}$, as $r^{-1}\partial_\theta u_r$ does relative to $u_r$ in the Navier slip condition. A higher scaling asks for faster vanishing at the corner and allows faster growth at infinity: a function homogeneous of degree $d$ has a finite norm near the corner exactly at the scalings $a<d+\sigma$, and near infinity exactly at the scalings $a>d+\sigma$, with $\sigma$ as in Remark \ref{weight and Mellin line}.

An equation whose terms all have the same scaling is invariant under the dilations, and its Mellin transform acts on each line $\Re\lambda=a-\sigma$ separately, while a term of a lower scaling shifts the line; in Section \ref{corner expansions} such terms are treated as data. There, too, a time derivative, which by the profile evolution equation acts on the profile as $\partial_r$ does, is worth one unit of scaling, and it is traded for one in Proposition \ref{time for radial trade}.

\begin{lemma}[Elementary properties of the norms]\label{elementary properties}
Let $O\in\{\Omega,\partial_0\Omega,\partial_1\Omega\}$ and $j\geq 0$. Whenever the RHSs are finite, we have the following.
\begin{enumerate}
\item (Ladder) If $a_0\leq a\leq a_1$ then $\|u\|_{O,a,j,m}\lesssim_{a_0,a_1,j}\|u\|_{O,a_0,j,m}+\|u\|_{O,a_1,j,m}$, with constant $1$ when $j$ is an integer. Moreover $\|u\|_{O,a,j,m}$ is non-decreasing in $m$ and in integer values of $j$, and for $j\in\N_0$ we have $\|u\|_{O,a,j,m}\leq\|u\|_{O,a,j+1/2,m}\leq\sqrt 2\|u\|_{O,a,j+1,m}$.
\item (Weight shift) For $k\in\R$ we have $\|r^{-k}u\|_{O,a,j,m}\sim_{k,j,m}\|u\|_{O,a+k,j,m}$.
\item (Exchange) We have $[\dr,r\grad]=\mb 0$, and $\dr=\mb e_r\cdot(r\grad)$ with $\dr\mb e_r=\mb 0$. Consequently, for integers $k\geq 0$ and $m\geq 0$, and with $j$ an integer in the second inequality of \eqref{exchange},
\begin{align}
\|u\|_{\Omega,a,j+k,m}&\lesssim_{j,k,m}\|u\|_{\Omega,a,j,m+k},\qquad\qquad \|u\|_{\Omega,a,j,m}\lesssim_{j,m}\|u\|_{\Omega,a,0,j+m},\label{exchange}\\
\|u\|_{\Omega,a,j,m}&\lesssim_{j,m}\sum_{k=0}^{m}\|\grad^{k}u\|_{\Omega,a-k,j}\lesssim_{j,m}\|u\|_{\Omega,a,j,m}.\label{scaling}
\end{align}
\end{enumerate}
\end{lemma}
\begin{proof}
Throughout we use the correspondence \eqref{weight and line} between the weight $a$ and the Mellin line $\Re\lambda=a-\sigma$, with $\sigma$ as in Remark \ref{weight and Mellin line}.

(i) The monotonicity in $m$ and in integer values of $j$ is immediate from the definition, as is the first of the two inequalities for half-integer counts. For the second, Definition \ref{fractional counts} and the pointwise bound $|\lambda|^{2j+1}\leq|\lambda|^{2j}+|\lambda|^{2j+2}$ on the line \eqref{weight and line} give $\|\dr^{j+1/2}u\|_{O,a,0,m}^2\leq\|\dr^{j}u\|_{O,a,0,m}^2+\|\dr^{j+1}u\|_{O,a,0,m}^2$, so that
\begin{align*}
\|u\|_{O,a,j+{1\over 2},m}^2&=\|u\|_{O,a,j,m}^2+\|\dr^{j+{1\over 2}}u\|_{O,a,0,m}^2
\leq\|u\|_{O,a,j,m}^2+\|\dr^{j}u\|_{O,a,0,m}^2+\|\dr^{j+1}u\|_{O,a,0,m}^2\\
&\leq 2\|u\|_{O,a,j,m}^2+\brac{\|u\|_{O,a,j+1,m}^2-\|u\|_{O,a,j,m}^2}
\leq 2\|u\|_{O,a,j+1,m}^2.
\end{align*}
The second inequality uses that $\|\dr^{j}u\|_{O,a,0,m}\leq\|u\|_{O,a,j,m}$, the LHS being one of the summands of the RHS, together with the identity $\|\dr^{j+1}u\|_{O,a,0,m}^2=\|u\|_{O,a,j+1,m}^2-\|u\|_{O,a,j,m}^2$; the two copies of $\|u\|_{O,a,j,m}^2$ are thus the one carried by the definition of the fractional norm and the one just used. For an integer count the ladder is the pointwise bound $r^{-a}\leq r^{-a_0}+r^{-a_1}$ applied to $\dr^{j'}(r\grad)^{m'}u$ under the integral. A non-integer count needs one preliminary step. The multiplier $|\lambda|^{j}$ through which Definition \ref{fractional counts} reads $\dr^{j}u$ is prescribed line by line and is not analytic, so there is no single function whose norms on the three lines $\Re\lambda=a_0-\sigma$, $a-\sigma$ and $a_1-\sigma$ the ladder would compare. We therefore work with $\dr_{-\beta}^{j}$ instead, whose Mellin multiplier is $(\lambda-\beta)^{j}$, where $\beta<0$ is chosen so small that $2\Re\lambda\geq\beta$ on the closed strip $\Re\lambda\in[a_0-\sigma,a_1-\sigma]$. On that strip $\lambda-\beta$ stays off the negative real axis, so $(\lambda-\beta)^{j}$ is analytic and non-vanishing there, $\dr_{-\beta}^{j}u$ is a single function, and the pointwise argument above applies to it. The choice of $\beta$ gives $|\Re\lambda-\beta|\geq|\Re\lambda|$, hence $|\lambda-\beta|\geq|\lambda|$, and therefore
\begin{align*}
|\lambda|^{2j}\leq|\lambda-\beta|^{2j}\lesssim_{j,\beta}|\lambda|^{2j}+1
\end{align*}
on the strip. The first inequality passes from $\dr^{j}u$ to $\dr_{-\beta}^{j}u$ on the left of the ladder, the second passes back on the right, and the term $1$ produced there gives $\|u\|_{O,a_0,0,m}$ and $\|u\|_{O,a_1,0,m}$, which are contained in the norms on the right.

(ii) For real $k$ and integer counts this follows by the binomial theorem from
\begin{align*}
\dr^{j'}(r^{-k}u)=r^{-k}(\dr-k)^{j'}u,\qquad\qquad
(r\grad)^{m'}(r^{-k}u)=r^{-k}\brac{(r\grad)-k\mb e_r}^{m'}u.
\end{align*}
For a non-integer count we use instead that $(r^{-k}u)^\wedge(\lambda)=\hat u(\lambda+k)$, so that, substituting $\mu=\lambda+k$ and using $|\mu-k|^{2j}\lesssim_{k,j}1+|\mu|^{2j}$,
\begin{align*}
\|\dr^{j}(r^{-k}u)\|_{O,a}^2&=\int_{\Re\lambda=a-\sigma}|\lambda|^{2j}|\hat u(\lambda+k)|^2\d\Im\lambda
=\int_{\Re\mu=a+k-\sigma}|\mu-k|^{2j}|\hat u|^2\d\Im\mu\lesssim_{k,j}\|u\|_{O,a+k,j}^2.
\end{align*}
The converse is the same statement with $k$ replaced by $-k$.

(iii) From $[\dr,\grad]=-\grad$ and $[\dr,r]=r$ we get $[\dr,r\grad]=r\grad-r\grad=\mb 0$. The field $\mb e_r$ is homogeneous of degree $0$, so $\dr\mb e_r=\mb 0$, while $(r\grad)^{n}\mb e_r$ is bounded on $\Omega$ for every $n$. Hence $\mb e_r$ passes through every $\dr$ and
\begin{align*}
\dr^{j+k}(r\grad)^{m'}u=\mb e_r^{\otimes k}:\dr^{j}(r\grad)^{m'+k}u,
\end{align*}
which gives the first inequality of \eqref{exchange} since $|\mb e_r|=1$, the constant coming from the fact that several pairs $(j',m')$ on the left are bounded by the same term on the right; the second follows by writing every $\dr$ in this way. For \eqref{scaling}, an induction on $m'$ using $(r\grad)(r^{n}\grad^{n}u)=nr^{n}\mb e_r\otimes\grad^{n}u+r^{n+1}\grad^{n+1}u$ shows that $(r\grad)^{m'}u$ and $r^{m'}\grad^{m'}u$ are, in either direction, finite combinations of the $r^{n}\grad^{n}u$, respectively of the $(r\grad)^{n}u$, with $n\leq m'$ and with coefficients built from $\mb e_r$ and bounded on $\Omega$; radial derivatives are then moved across by $\dr(r^{n}\grad^{n}u)=r^{n}(\dr+n)\grad^{n}u$.
\end{proof}

By Lemma \ref{elementary properties}(i), the finiteness of the norm $\|u\|_{O,a,j,m}$ is monotone in $j$ through half-integer values as well, although the norm itself need not be, since $|\lambda|^{2j+1}>|\lambda|^{2j+2}$ for $|\lambda|<1$.

Finally, we record a mollifier in the radial variable.

\begin{lemma}[Radial smoothing]\label{radial smoothing}
For $\sigma\in[0,1]$ let $S_\sigma$ denote the Mellin multiplier $e^{\sigma^2\lambda^2}$, that is, convolution in the variable $s=\ln r$, at each fixed $\theta$ where there is one, with the Gaussian of integral one and variance $2\sigma^2$, with $S_0$ the identity. Let $O\in\{\Omega,\partial_0\Omega,\partial_1\Omega\}$, let the weight $a$ range over a compact set, on which the constants below are allowed to depend, and let $j\geq 0$, $m\in\N_0$, $k\in\N_0$. Then the following hold whenever the RHSs are finite.
\begin{enumerate}
\item (Boundedness and gain) For $n\in\N_0$,
\begin{align}\label{smoothing gain}
\|(\sigma\dr)^{2n}S_\sigma u\|_{O,a,j+k,m}\lesssim_{n,k}\sigma^{-k}\|u\|_{O,a,j,m}.
\end{align}
\item (Structure) $S_\sigma$ commutes with $\dr$, with $r\grad$, with $\partial_\theta$ and with restriction to $\partial_0\Omega$ and to $\partial_1\Omega$, and it acts as the identity on constants. Applied componentwise in polar components it preserves the conditions $\grad\cdot\mb v=0$ and $v_\theta|_{\partial_0\Omega}=0$.
\item (Scale derivative) $\partial_\sigma S_\sigma=2\sigma\dr^2S_\sigma$ for $0<\sigma\leq 1$, and $\sigma\mapsto S_\sigma u$ is continuous $[0,1]\to\|\ph\|_{O,a,j,m}$ whenever $\|u\|_{O,a,j,m}<\infty$; in particular $S_\sigma u\to u$ in that norm as $\sigma\to 0$, and $(\sigma\dr)^{2n}S_\sigma u\to 0$ in the same norm for every $n\geq 1$. Moreover, $\sigma\mapsto S_\sigma u$ is continuous $(0,1]\to\|\ph\|_{O,a,j+k,m}$.
\end{enumerate}
\end{lemma}
\begin{proof}
Write $c$ for the real part of $\lambda$ on the line carried by the weight $a$, that is $c=a-1$ for $O=\Omega$ and $c=a-{1\over 2}$ otherwise, which ranges over a compact set with $a$. On that line
\begin{align*}
|(\sigma\lambda)^{2n}\lambda^ke^{\sigma^2\lambda^2}|
=|\sigma\lambda|^{2n}|\lambda|^ke^{\sigma^2(c^2-(\Im\lambda)^2)}
\lesssim_{n,k,c}\sigma^{-k},
\end{align*}
since $|\lambda|\leq|c|+|\Im\lambda|$, every power of $x$ times $e^{-x^2}$ is bounded on $[0,\infty)$, and $\sigma\leq 1$. By Plancherel’s identity \eqref{Plancherel's identity} on each line, together with Definition \ref{fractional counts} for the fractional counts and the commutations of (ii) for the full derivative counts, this gives (i). For (ii), the convolution converges, and agrees with the multiplier on every line of a strip at once, by the superexponential decay of the Gaussian; commutation with $\dr$, $r\grad$, $\partial_\theta$ and with restriction to rays is then immediate from the Mellin representation, and $S_\sigma$ acts as the identity on constants since the Gaussian has integral one. Since $S_\sigma$ commutes with $r\grad$, componentwise smoothing gives $r\grad\cdot(S_\sigma\mb v)=S_\sigma(r\grad\cdot\mb v)$, which vanishes with $\grad\cdot\mb v$; and $v_\theta|_{\partial_0\Omega}=0$ is preserved since $S_\sigma$ commutes with the restriction. Part (iii) is the identity $\partial_\sigma e^{\sigma^2\lambda^2}=2\sigma\lambda^2e^{\sigma^2\lambda^2}$ together with dominated convergence on each line.
\end{proof}

By \eqref{smoothing gain}, $S_\sigma$ is bounded uniformly in $\sigma$, and raises the radial count by $k$ at the cost $\sigma^{-k}$.

\subsection{Wedge elliptic systems and Helmholtz projection}\label{wedge elliptic systems}

The two elliptic systems of this subsection are the scalar problems that the pressure of our Stokes system solves: a Dirichlet condition on the free boundary and a Neumann condition on the solid boundary in the first case, and Neumann conditions on both in the second. They are solved explicitly by the Mellin transform, which turns them into two-point boundary value problems in $\theta$, and the resulting solution formulae are what produce the two projections $\P$ and $\P^\theta$ that encode the two orthogonality constraints of our system.

Problems of this kind are the subject of the corner theory of Kondrat'ev, Grisvard, and Kozlov, Maz'ya and Rossmann \cite{Kondratev,Grisvard,KozlovMazyaRossmann,KozlovMazyaRossmannPencils}. We nevertheless solve them by hand, for two reasons. First, the wedge is unbounded, so that the analysis cannot be localised near the corner. Second, and more to the point here, the opening angle $\theta_b$ is a small parameter of the whole paper: it is the smallness of the angle that makes the coercivity of Section \ref{linear existence} work. So every constant has to be tracked in $\theta_b$, and the explicit formulae below make this possible. We first need two lemmas: one for estimating the $\cos$ and $\sin$ terms that appear in the solution formulae, and one for the angular integral that turns the exponential decay across the wedge into a factor $\{\theta_b,|\lambda|^{-1}\}$.

\begin{lemma}[Bounds for $\cos(\lambda\theta)$ and $\sin(\lambda\theta)$]\label{cosh bound}
Let $\lambda\in\C$, $\theta>0$ and $C>0$, and write, for $a\in\R$,
\begin{align}\label{sine resonance factor}
\Gamma_{a,\theta}:=\min\{1,\pi-|a\theta|\}\qquad\text{ for }\qquad|a\theta|<\pi,
\end{align}
the distance of $a\theta$ from the zeros $\pm\pi$ of the sine, capped at $1$. For $|\Re\lambda\theta|<\pi/2$ we have
\begin{align*}
\cos(\Re\lambda\theta)e^{|\Im\lambda|\theta}\lesssim|\cos(\lambda\theta)|&\lesssim e^{|\Im\lambda|\theta},
\end{align*}
and for $|\Re\lambda\theta|<\pi$ we have
\begin{align*}
e^{|\Im\lambda|\theta}\mb 1_{\{|\Im\lambda\theta|\geq C\}}+\Gamma_{\Re\lambda,\theta}|\lambda\theta|\mb 1_{\{|\Im\lambda\theta|<C\}}\lesssim_C|\sin(\lambda\theta)|&\lesssim_C e^{|\Im\lambda|\theta}\mb 1_{\{|\Im\lambda\theta|\geq C\}}+|\lambda\theta|\mb 1_{\{|\Im\lambda\theta|<C\}}\\
&\lesssim_C e^{|\Im\lambda|\theta}.
\end{align*}
\end{lemma}
\begin{proof}
Writing $x=\Re\lambda\theta$ and $y=\Im\lambda\theta$, the addition formulae give
\begin{align*}
\cos(\lambda\theta)&=\cos x\cosh y-i\sin x\sinh y,\qquad\qquad
\sin(\lambda\theta)=\sin x\cosh y+i\cos x\sinh y,
\end{align*}
hence, by $\cosh^2y=1+\sinh^2y$ and $\cos^2x+\sin^2x=1$,
\begin{align*}
|\cos(\lambda\theta)|^2=\cos^2x+\sinh^2y,\qquad\qquad|\sin(\lambda\theta)|^2=\sin^2x+\sinh^2y.
\end{align*}
The upper bounds follow from $|\cos x|,|\sin x|\leq 1$, $|\sinh y|\leq\cosh y\leq e^{|y|}$, and, when $|y|<C$, from $|\sin x|\leq|x|$ and $|\sinh y|\leq C^{-1}\sinh(C)|y|$. For the lower bounds we take the square root of $a^2+b^2\geq{1\over 2}(a+b)^2$, valid for $a,b\geq 0$, which gives $\sqrt{a^2+b^2}\geq{1\over\sqrt 2}(a+b)$ and hence
\begin{align*}
|\cos(\lambda\theta)|\geq{1\over\sqrt 2}\brac{|\cos x|+|\sinh y|},\qquad\qquad
|\sin(\lambda\theta)|\geq{1\over\sqrt 2}\brac{|\sin x|+|\sinh y|}.
\end{align*} When $|y|\geq C$ we have $|\sinh y|\geq{1\over 2}(1-e^{-2C})e^{|y|}$, and when $|y|<C$ we have $|\sinh y|\geq|y|$ and $\cos x\geq e^{-C}\cos(x)e^{|y|}$; this gives the bound for the cosine, $\cos x$ being positive when $|x|<\pi/2$, and the bound for the sine when $|y|\geq C$. When $|y|<C$ and $|x|<\pi$, the concavity of the sine on $[0,\pi]$ gives $\sin|x|\geq{2\over\pi}\min\{|x|,\pi-|x|\}$, and $\min\{|x|,\pi-|x|\}\geq\min\{|x|,1\}\Gamma_{\Re\lambda,\theta}$, as one checks separately for $|x|\leq 1$, where $\pi-|x|>1$, and for $|x|>1$; since moreover $\min\{|x|,1\}\geq|x|/\pi$ and $\Gamma_{\Re\lambda,\theta}\leq 1$, we get
\begin{align*}
|\sin(\lambda\theta)|\geq{1\over\sqrt 2}\brac{{2\over\pi^2}\Gamma_{\Re\lambda,\theta}|x|+|y|}\geq{\sqrt 2\over\pi^2}\Gamma_{\Re\lambda,\theta}(|x|+|y|)\geq{\sqrt 2\over\pi^2}\Gamma_{\Re\lambda,\theta}|\lambda\theta|,
\end{align*}
which is the remaining bound.
\end{proof}

By \eqref{sine resonance factor}, the factor $\Gamma_{\Re\lambda,\theta}$ in Lemma \ref{cosh bound} is $1$ unless $|\Re\lambda\theta|>\pi-1$.

The next lemma is the one place in this subsection where the exponential decay of the solution formulae across the wedge is converted into a gain. It is used in the proofs of both Lemma \ref{mixed wedge problem} and Lemma \ref{Neumann wedge problem}, on strips of different width, which is why the dependence of the constant on that width is recorded.

\begin{lemma}[The angular integral]\label{angular integral}
Let $\theta_b>0$, let $c>0$ and let $\lambda\in\C$ with $|\Re\lambda|\theta_b\leq c$. Then
\begin{align}\label{angular integral bound}
\int_0^{\theta_b}e^{-2|\Im\lambda|\theta}\d\theta\leq(c+1)\set{\theta_b,|\lambda|^{-1}},
\end{align}
the brace being a choice as in Definition \ref{notation}. The same bound holds for $\int_0^{\theta_b}e^{-2|\Im\lambda|(\theta_b-\theta)}\d\theta$.
\end{lemma}
\begin{proof}
The bound for $\int_0^{\theta_b}e^{-2|\Im\lambda|(\theta_b-\theta)}\d\theta$ follows from \eqref{angular integral bound} by the substitution $\theta\mapsto\theta_b-\theta$. Both choices have to be proved. The choice $\theta_b$ is immediate, the integrand being at most $1$. For the choice $|\lambda|^{-1}$ we split according to the size of $|\Im\lambda\theta_b|$. If $|\Im\lambda\theta_b|\geq 1$, then integrating and using $|\lambda|\theta_b\leq|\Re\lambda|\theta_b+|\Im\lambda|\theta_b\leq(c+1)|\Im\lambda\theta_b|$,
\begin{align*}
\int_0^{\theta_b}e^{-2|\Im\lambda|\theta}\d\theta\leq{1\over 2|\Im\lambda|}={\theta_b\over 2|\Im\lambda\theta_b|}\leq{c+1\over 2}|\lambda|^{-1}.
\end{align*}
If $|\Im\lambda\theta_b|<1$, then $|\lambda|\theta_b\leq|\Re\lambda|\theta_b+|\Im\lambda|\theta_b<c+1$, so that bounding the integrand by $1$ gives
\begin{align*}
\int_0^{\theta_b}e^{-2|\Im\lambda|\theta}\d\theta&\leq\theta_b<(c+1)|\lambda|^{-1}.\qedhere
\end{align*}
\end{proof}

Now we can solve the elliptic systems. The second estimate of the next lemma is used in Lemma \ref{test problem pressure estimate}.

\begin{lemma}[Mixed Dirichlet--Neumann problem on the wedge]\label{mixed wedge problem}
Consider the following three elliptic systems
\begin{alignat*}{3}
\lpc(q_{\mb F},q_{f_0},q_{f_1})&=(\grad\cdot(r^{-1}\mb F),0,0)\qquad&&\text{on}\qquad\Omega\\
\partial_\theta(q_{\mb F},q_{f_0},q_{f_1})&=(F_\theta,f_0,0)\qquad&&\text{on}\qquad\partial_0\Omega\\
(q_{\mb F},q_{f_0},q_{f_1})&=(0,0,f_1)\qquad&&\text{on}\qquad\partial_1\Omega.
\end{alignat*}
Let the weight $a\in\R$ satisfy $|a-1|\theta_b<\pi/2$. Then these systems have solutions $q_{\mb F},q_{f_0},q_{f_1}$ with regularity dependent on the regularity of $\mb F$, $f_0$, $f_1$ respectively, in the following sense, for $j\geq 0$ and $l\in\{0,1\}$:
\begin{align*}
\|\dr^j\partial_\theta^l q_{\mb F}\|_{\Omega,a} &\lesssim{1\over|\cos((1-a)\theta_b)|} \norm{\set{\theta_b,\dr^{-1}}^{1-l}\dr^{j}\mb F}_{\Omega,a}\\
\|\dr^j\partial_\theta^l q_{\mb F}\|_{\Omega,a} &\lesssim{1\over|\cos((1-a)\theta_b)|} \norm{\set{\theta_b,\dr^{-1}}\dr^{j+l}F_r}_{\Omega,a}\qquad\text{when}\qquad F_\theta=0\\
\|(r\grad)^{j+1}q_{\mb F}\|_{\Omega,a}&\lesssim_j{1\over|\cos((1-a)\theta_b)|}\|\mb F\|_{\Omega,a,0,j}\\
\|(r\grad)^jq_{f_0}\|_{\Omega,a}&\lesssim_j{1\over|\cos((1-a)\theta_b)|}\sum_{j'=\min\{1,j\}}^j\norm{\set{\theta_b^{1\over 2},\dr^{-{1\over 2}}}\dr^{j'-1}f_0}_{\partial_0\Omega,a-{1\over 2}}\\
\|(r\grad)^jq_{f_1}\|_{\Omega,a}&\lesssim_j{1\over|\cos((1-a)\theta_b)|}\sum_{j'=\min\{1,j\}}^j\norm{\set{\theta_b^{1\over 2},\dr^{-{1\over 2}}}\dr^{j'}f_1}_{\partial_1\Omega,a-{1\over 2}}\\
\|(r\grad)^jq_{f_1}\|_{\partial_0\Omega,a-{1\over 2}}
&\lesssim_j{\theta_b^{-j}\over|\cos((1-a)\theta_b)|}\|f_1\|_{\partial_1\Omega,a-{1\over 2}},
\end{align*}
where $\{\ph\}$ is a choice as in Definition \ref{notation}.
\end{lemma}
\begin{proof}
Mellin transform $r\mapsto\lambda$ turns the system of the statement, in the form $r^2\lpc q=r^2\grad\cdot(r^{-1}\mb F)$ on $\Omega$ with $q=f_1$ on $\partial_1\Omega$ and $\partial_\theta q=f_0+F_\theta$ on $\partial_0\Omega$, into the system
\begin{align}
(\lambda^2+\partial_\theta^2)\hat q(\lambda,\theta)&=\lambda\hat F_r(\lambda,\theta)+\partial_\theta\hat F_\theta(\lambda,\theta)\nonumber\\
\hat q(\lambda,\theta_b)&=\hat f_1(\lambda)\label{mixed Mellin system}\\
\partial_\theta\hat q(\lambda,0)&=\hat f_0(\lambda)+\hat F_\theta(\lambda,0),\nonumber
\end{align}
a two-point boundary value problem in $\theta$ for every fixed $\lambda$. We solve it by superposition: the two boundary data are carried by explicit solutions of the homogeneous equation, and the bulk RHS is expanded in eigenfunctions. The latter part is the one that needs an argument, and we give it first.

Let
\begin{align*}
h_k=\sqrt{2\over\theta_b}\cos\brac{\brac{{\pi\over 2}+k\pi}{\theta\over\theta_b}},\qquad k\in\N_0,
\end{align*}
be the $L^2([0,\theta_b])$-orthonormal eigenfunctions of $\partial_\theta^2$ with boundary condition $h(\theta_b)=0=\partial_\theta h(0)$, with eigenvalues $-\Lambda_k^2=-(\pi/2+k\pi)^2\theta_b^{-2}$. The functions $h_k'/\Lambda_k$ are again orthonormal, being the eigenfunctions of $\partial_\theta^2$ with the complementary boundary condition $h(0)=0=\partial_\theta h(\theta_b)$, and $\partial_\theta$ carries them back onto the $h_k$, since $(h_k'/\Lambda_k)'=-\Lambda_kh_k$. Since the RHS of \eqref{mixed Mellin system} carries $\hat F_r$ undifferentiated but $\hat F_\theta$ differentiated, we expand the two components in the two bases,
\begin{align*}
\hat F_{r,k}=\<h_k,\hat F_{r}\>_{L^2([0,\theta_b])},\qquad\qquad
\hat F_{\theta,k}=\inner{{h_k'\over\Lambda_k},\hat F_{\theta}}_{L^2([0,\theta_b])},
\end{align*}
the inner product being conjugate linear in its first slot, as in \eqref{Plancherel's identity}, so that these are the expansion coefficients $\int_0^{\theta_b}\hat F_rh_k\d\theta$ and $\int_0^{\theta_b}\hat F_\theta h_k'/\Lambda_k\d\theta$ of the two components. It is essential to pair $\hat F_\theta$ with $h_k'/\Lambda_k$ and not with $h_k$: expanding $\hat F_\theta$ in the $h_k$ and differentiating term by term would put $\hat q_{\mb F}$ on the basis $h_k'/\Lambda_k$, which satisfies the complementary boundary conditions and so is not admissible.

Write $\hat q_{\mb F}=\sum_kc_kh_k$ for the solution of \eqref{mixed Mellin system} with $\hat f_0=\hat f_1=0$. Testing the equation against $h_k$ and integrating by parts twice gives
\begin{align*}
(\lambda^2-\Lambda_k^2)c_k+\eva{\brac{\partial_\theta\hat q_{\mb F}h_k-\hat q_{\mb F}h_k'}}_0^{\theta_b}
=\lambda\hat F_{r,k}-\Lambda_k\hat F_{\theta,k}+\eva{\hat F_\theta h_k}_0^{\theta_b}.
\end{align*}
At $\theta=\theta_b$ every contribution vanishes, since $\hat q_{\mb F}(\theta_b)=0=h_k(\theta_b)$, and at $\theta=0$ so does the term carrying $h_k'(0)$, since $h_k'(0)=0$; what is left is $-\partial_\theta\hat q_{\mb F}(0)h_k(0)$ on the left and $-\hat F_\theta(0)h_k(0)$ on the right, and these two cancel each other precisely under the boundary condition $\partial_\theta q_{\mb F}|_{\partial_0\Omega}=F_\theta|_{\partial_0\Omega}$ of the statement, the natural condition for the divergence form RHS. We are left with $(\lambda^2-\Lambda_k^2)c_k=\lambda\hat F_{r,k}-\Lambda_k\hat F_{\theta,k}$, and therefore with the solution
\begin{align*}
\hat q&=\underbrace{-{\sin(\lambda(\theta_b-\theta))\over\lambda\cos(\lambda\theta_b)}\hat f_0}_{\hat q_{f_0}}+\underbrace{{\cos(\lambda\theta)\over\cos(\lambda\theta_b)}\hat f_1}_{\hat q_{f_1}}
+\underbrace{\sum_{k=0}^\infty{(\lambda\hat F_{r,k}-\Lambda_k\hat F_{\theta,k})h_k\over\lambda^2-\Lambda_k^2}}_{\hat q_{\mb F}},
\end{align*}
the first two terms being the solutions of the homogeneous equation that carry $\hat f_0$ and $\hat f_1$, as one checks directly.

We start by proving the estimates for the derivatives of $q_{\mb F}$. In the formula for $\hat q_{\mb F}$ given by the eigenfunction expansion, a derivative $\partial_\theta$ corresponds on the Mellin side to multiplication by a factor of $\Lambda_k$, while a derivative $\dr$ corresponds to multiplication by $\lambda$. The denominator $(\lambda^2-\Lambda_k^2)^{-1}$ therefore buys us two such factors. Precisely, we claim that for $l\in\{0,1,2\}$ and $\Re\lambda=a-1$,
\begin{align}\label{mixed multiplier bound}
\sup_{k\geq 0}\abs{\{\lambda,\Lambda_k\}^l\over\lambda^2-\Lambda_k^2}&\lesssim{\{\theta_b,|\lambda|^{-1}\}^{2-l}\over|\cos(\Re\lambda\theta_b)|},
\end{align}
the two braces being choices in the sense of Definition \ref{notation}, made independently of each other.

To see \eqref{mixed multiplier bound}, note first that $\lambda^2$ is even in $\lambda$, so we may assume $\Re\lambda\geq 0$; the hypothesis $|a-1|\theta_b<\pi/2$ then reads $0\leq\Re\lambda<\Lambda_0$. Hence $|\lambda+\Lambda_k|\geq\max\{\Lambda_k,|\Im\lambda|\}$, while for the other factor we have $\Lambda_0-\Re\lambda=\theta_b^{-1}(\pi/2-\Re\lambda\theta_b)\geq\theta_b^{-1}\cos(\Re\lambda\theta_b)={2\over\pi}\Lambda_0\cos(\Re\lambda\theta_b)$ when $k=0$, and $\Lambda_k-\Re\lambda\geq\Lambda_k-\Lambda_0\geq{2\over 3}\Lambda_k$ when $k\geq 1$, so that in either case
\begin{align*}
|\lambda-\Lambda_k|\geq\max\{\Lambda_k-\Re\lambda,|\Im\lambda|\}\gtrsim|\cos(\Re\lambda\theta_b)|\max\{\Lambda_k,|\Im\lambda|\}.
\end{align*}
Multiplying the two, and using $|\lambda|\leq\Lambda_k+|\Im\lambda|$, which holds because $\Re\lambda<\Lambda_0\leq\Lambda_k$, we obtain
\begin{align}\label{mixed denominator bound}
|\lambda^2-\Lambda_k^2|\gtrsim|\cos(\Re\lambda\theta_b)|\brac{\Lambda_k^2+|\lambda|^2}.
\end{align}
The numerator of \eqref{mixed multiplier bound} is at most $(\Lambda_k^2+|\lambda|^2)^{l/2}$, each of its $l$ factors being $|\lambda|$ or $\Lambda_k$, so that by \eqref{mixed denominator bound} the quotient is at most $|\cos(\Re\lambda\theta_b)|^{-1}(\Lambda_k^2+|\lambda|^2)^{-(2-l)/2}$. This is where the RHS of \eqref{mixed multiplier bound} comes from: each of the remaining $2-l$ factors $(\Lambda_k^2+|\lambda|^2)^{-1/2}$ may be estimated, independently of the others, either by $|\lambda|^{-1}$ or by $\Lambda_k^{-1}\leq\Lambda_0^{-1}=2\theta_b/\pi$, and the two options are exactly the two entries of the brace $\{\theta_b,|\lambda|^{-1}\}$.

Thus \eqref{mixed multiplier bound} lets us control up to one $\partial_\theta$ derivative of $q_{\mb F}$; indeed, for $j\geq 0$ and $l\in\{0,1\}$,
\begin{align}
\|\dr^j\partial_\theta^l q_{\mb F}\|_{\Omega,a}^2&=\int_{\Re\lambda=a-1}\int_0^{\theta_b}|\lambda^j\partial_\theta^l\hat q_{\mb F}|^2\d\theta\d\Im\lambda\nonumber\\
&=\int_{\Re\lambda=a-1}\int_0^{\theta_b}\abs{\lambda^j\partial_\theta^l\sum_{k=0}^\infty{(\lambda\hat F_{r,k}-\Lambda_k\hat F_{\theta,k})h_k\over\lambda^2-\Lambda_k^2}}^2\d\theta\d\Im\lambda\nonumber\\
&\leq\int_{\Re\lambda=a-1}\int_0^{\theta_b}\abs{\lambda^j\sum_{k=0}^\infty{(\lambda\hat F_{r,k}-\Lambda_k\hat F_{\theta,k})\Lambda_k^{l}(h_k^{(l)}/\Lambda_k^{l})\over\lambda^2-\Lambda_k^2}}^2\d\theta\d\Im\lambda\nonumber\\
&\leq\int_{\Re\lambda=a-1}\abs{|\lambda|^{j}{\{\theta_b,|\lambda|^{-1}\}^{1-l}\over|\cos(\Re\lambda\theta_b)|}}^2\sum_{k=0}^\infty\brac{|\hat F_{r,k}|^2+|\hat F_{\theta,k}|^2}\d\Im\lambda\nonumber\\
&\leq\int_{\Re\lambda=a-1}|\lambda|^{2j}{\{\theta_b^2,|\lambda|^{-2}\}^{1-l}\over|\cos(\Re\lambda\theta_b)|^2}\int_0^{\theta_b}|\hat{\mb F}|^2\d\theta\d\Im\lambda\nonumber\\
&\leq{1\over|\cos((1-a)\theta_b)|^2}\norm{\set{\theta_b,\dr^{-1}}^{1-l}\dr^{j}\mb F}_{\Omega,a}^2.\label{mixed bulk estimate}
\end{align}
When $\hat F_\theta=0$ the coefficients are $\lambda\hat F_{r,k}$ alone, and the extra factor $\lambda$ they carry can be played against the brace of \eqref{mixed multiplier bound} rather than simply pulled out; this is what buys the refinement. What has to be shown is that, for $l\in\{0,1\}$,
\begin{align}\label{mixed multiplier bound theta free}
\sup_{k\geq 0}\abs{\lambda\Lambda_k^{l}\over\lambda^2-\Lambda_k^2}\lesssim|\lambda|^{l}{\{\theta_b,|\lambda|^{-1}\}\over|\cos(\Re\lambda\theta_b)|},
\end{align}
with either entry of the brace, which is four statements, each an instance of \eqref{mixed multiplier bound}. For the entry $|\lambda|^{-1}$ we use \eqref{mixed multiplier bound} at exponent $2$, whose brace on the right is then empty. Reading the two factors of $\{\lambda,\Lambda_k\}^2$ as $\lambda,\Lambda_k$ gives $\sup_k|\lambda\Lambda_k/(\lambda^2-\Lambda_k^2)|\lesssim|\cos(\Re\lambda\theta_b)|^{-1}$, which is the case $l=1$; reading them as $\lambda,\lambda$ gives $\sup_k|\lambda^2/(\lambda^2-\Lambda_k^2)|\lesssim|\cos(\Re\lambda\theta_b)|^{-1}$, which is the case $l=0$ after one factor $|\lambda|$ has been divided out. For the entry $\theta_b$ we use \eqref{mixed multiplier bound} at exponent $1$, with the entry $\theta_b$ of the brace on its right. Reading the single factor of $\{\lambda,\Lambda_k\}$ as $\lambda$ gives the case $l=0$ at once; reading it as $\Lambda_k$ gives $\sup_k|\Lambda_k/(\lambda^2-\Lambda_k^2)|\lesssim\theta_b|\cos(\Re\lambda\theta_b)|^{-1}$, which is the case $l=1$ once the factor $|\lambda|$, which does not depend on $k$, has been taken out of the supremum. With \eqref{mixed multiplier bound theta free} in hand, the third line above is bounded, for $l\in\{0,1\}$, by
\begin{align*}
\int_{\Re\lambda=a-1}\abs{|\lambda|^{j+l}{\{\theta_b,|\lambda|^{-1}\}\over|\cos(\Re\lambda\theta_b)|}}^2\sum_{k=0}^\infty|\hat F_{r,k}|^2\d\Im\lambda
\leq{1\over|\cos((1-a)\theta_b)|^2}\norm{\set{\theta_b,\dr^{-1}}\dr^{j+l}F_r}_{\Omega,a}^2,
\end{align*}
which is the second estimate of the statement.

For two or more $\partial_\theta$ derivatives of $q_{\mb F}$, we use the equation $(\lambda^2+\partial_\theta^2)\hat q_{\mb F}=\lambda\hat F_r+\partial_\theta\hat F_\theta$ repeatedly: writing $\hat G:=\lambda\hat F_r+\partial_\theta\hat F_\theta$, an induction on $l\in\N_0$ gives
\begin{align*}
\partial_\theta^{2l+2}\hat q_{\mb F}&=(-\lambda^2)^{l+1}\hat q_{\mb F}+\sum_{l'=0}^{l}(-\lambda^2)^{l-l'}\partial_\theta^{2l'}\hat G,\\
\partial_\theta^{2l+3}\hat q_{\mb F}&=(-\lambda^2)^{l+1}\partial_\theta\hat q_{\mb F}+\sum_{l'=0}^{l}(-\lambda^2)^{l-l'}\partial_\theta^{2l'+1}\hat G,
\end{align*}
in which every elimination of a $\partial_\theta^2$ costs a factor $-\lambda^2$. Hence, for $l\in\N_0$, using \eqref{mixed bulk estimate} with $l=0$ and $j+2l+2$ in place of $j$, with the choice $\dr^{-1}$, for the first term,
\begin{align}
\|\dr^j\partial_\theta^{2l+2}q_{\mb F}\|_{\Omega,a}^2&=\int_{\Re\lambda=a-1}\int_0^{\theta_b}|\lambda^j\partial_\theta^{2l+2}\hat q_{\mb F}|^2\d\theta\d\Im\lambda\nonumber\\
&\lesssim_l\int_{\Re\lambda=a-1}\int_0^{\theta_b}\brac{|\lambda^{j+2l+2}\hat q_{\mb F}|^2+\sum_{l'=0}^l|\lambda^{j+2(l-l')}\partial_\theta^{2l'}\hat G|^2}\d\theta\d\Im\lambda\nonumber\\
&\lesssim_l{1\over|\cos((1-a)\theta_b)|^2}\sum_{l'=0}^l\|\dr^{j+2(l-l')}\partial_\theta^{2l'}\mb F\|_{\Omega,a,0,1}^2,\label{higher bulk 1}
\end{align}
the first term being bounded by $\|\dr^{j+2l+1}\mb F\|_{\Omega,a}^2\leq\|\dr^{j+2l}\mb F\|_{\Omega,a,0,1}^2$; and likewise, using \eqref{mixed bulk estimate} with $l=1$ and $j+2l+2$ in place of $j$ for the first term,
\begin{align}
\|\dr^j\partial_\theta^{2l+3}q_{\mb F}\|_{\Omega,a}^2&=\int_{\Re\lambda=a-1}\int_0^{\theta_b}|\lambda^j\partial_\theta^{2l+3}\hat q_{\mb F}|^2\d\theta\d\Im\lambda\nonumber\\
&\lesssim_l\int_{\Re\lambda=a-1}\int_0^{\theta_b}\brac{|\lambda^{j+2l+2}\partial_\theta\hat q_{\mb F}|^2+\sum_{l'=0}^l|\lambda^{j+2(l-l')}\partial_\theta^{2l'+1}\hat G|^2}\d\theta\d\Im\lambda\nonumber\\
&\lesssim_l{1\over|\cos((1-a)\theta_b)|^2}\sum_{l'=0}^l\|\dr^{j+2(l-l')}\partial_\theta^{2l'}\mb F\|_{\Omega,a,0,2}^2.\label{higher bulk 2}
\end{align}
In both displays the total number of derivatives falling on $\mb F$ is $j+2l+1$, respectively $j+2l+2$, that is, one less than the number falling on $q_{\mb F}$. On polar components, $\partial_\theta$ acts as $r\mb e_\theta\cdot\grad$ up to the bounded frame terms $\partial_\theta\mb e_r=\mb e_\theta$, $\partial_\theta\mb e_\theta=-\mb e_r$. So, exactly as in the proof of Lemma \ref{elementary properties}(iii), $(r\grad)^{j+1}q_{\mb F}$ is a combination with bounded coefficients of the $\dr^{m}\partial_\theta^{n}q_{\mb F}$ with $m+n\leq j+1$, and conversely the $\dr^{m}\partial_\theta^{n}\mb F$ with $m+n\leq j$ are bounded by $\|\mb F\|_{\Omega,a,0,j}$. The terms with $n\leq 1$ are covered by the first estimate of the statement, with the choice $\dr^{-1}$ when $n=0$, and the terms with $n\geq 2$ by the two displays, applied with the radial count $m$ in place of $j$ and with $l$ chosen so that $n=2l+2$, respectively $n=2l+3$. This gives $\|(r\grad)^{j+1}q_{\mb F}\|_{\Omega,a}\lesssim_j|\cos((1-a)\theta_b)|^{-1}\|\mb F\|_{\Omega,a,0,j}$, which completes the estimates for $q_{\mb F}$.

Now we will prove the estimates for $q_{f_0}$ and $q_{f_1}$. As for $q_{\mb F}$, $(r\grad)^jq$ is a combination with bounded coefficients of the $\dr^m\partial_\theta^nq$ with $\min\{1,j\}\leq m+n\leq j$, and we bound each of them by the term $j'=m+n$ of the statement. On the line $\Re\lambda=a-1$ we have $|\Re\lambda|\theta_b=|a-1|\theta_b<\pi/2$, so Lemma \ref{angular integral} applies with $c=\pi/2$ and gives
\begin{align*}
\int_0^{\theta_b}e^{-2|\Im\lambda|\theta}\d\theta\leq\brac{{\pi\over 2}+1}\set{\theta_b,|\lambda|^{-1}}.
\end{align*}
Hence, in light of Lemma \ref{cosh bound} we have
\begin{align*}
\|\dr^m\partial_\theta^nq_{f_0}\|_{\Omega,a}^2&=\int_{\Re\lambda=a-1}\int_0^{\theta_b}|\lambda^m\partial_\theta^n\hat q_{f_0}|^2\d\theta\d\Im\lambda\\
&=\int_{\Re\lambda=a-1}\int_0^{\theta_b}\abs{\hat f_0\lambda^m\partial_\theta^n{\sin(\lambda(\theta_b-\theta))\over\lambda\cos(\lambda\theta_b)}}^2\d\theta\d\Im\lambda\\
&\lesssim{1\over|\cos((1-a)\theta_b)|^2}\int_{\Re\lambda=a-1}\int_0^{\theta_b}|\lambda|^{2(j'-1)}|\hat f_0|^2e^{-2|\Im\lambda|\theta}\d\theta\d\Im\lambda\\
&\lesssim{1\over|\cos((1-a)\theta_b)|^2}\int_{\Re\lambda=a-1}\set{\theta_b,|\lambda|^{-1}}|\lambda|^{2(j'-1)}|\hat f_0|^2\d\Im\lambda\\
&\leq{1\over|\cos((1-a)\theta_b)|^2} \norm{r^{{1\over 2}-a} \set{\theta_b^{1\over 2},\dr^{-{1\over 2}}}\dr^{j'-1}f_0}_{L^2([0,\infty))}^2\\
&={1\over|\cos((1-a)\theta_b)|^2} \norm{\set{\theta_b^{1\over 2},\dr^{-{1\over 2}}}\dr^{j'-1}f_0}_{\partial_0\Omega,a-{1\over 2}}^2.
\end{align*}
The estimates for $q_{f_1}$ are the same, with $\cos(\lambda\theta)$ in place of $\sin(\lambda(\theta_b-\theta))/\lambda$, hence $|\lambda|^{2j'}$ in place of $|\lambda|^{2(j'-1)}$, $e^{2|\Im\lambda|(\theta-\theta_b)}$ in place of $e^{-2|\Im\lambda|\theta}$, and $\dr^{j'}f_1$ on $\partial_1\Omega$ in place of $\dr^{j'-1}f_0$ on $\partial_0\Omega$. Finally, on $\theta=0$ we have $|\lambda^m\partial_\theta^n\hat q_{f_1}|\leq|\lambda|^{m+n}|\hat f_1|/|\cos(\lambda\theta_b)|\leq\<\lambda\>^j|\hat f_1|/|\cos(\lambda\theta_b)|$, so that
\begin{align*}
\|(r\grad)^jq_{f_1}\|_{\partial_0\Omega,a-{1\over 2}}^2&\lesssim_j\int_{\Re\lambda=a-1}\abs{\<\lambda\>^j\hat f_1\over\cos(\lambda\theta_b)}^2\d\Im\lambda\\
&\lesssim{1\over|\cos((1-a)\theta_b)|^2}\int_{\Re\lambda=a-1}\<\lambda\>^{2j}e^{-2|\Im\lambda|\theta_b}|\hat f_1|^2\d\Im\lambda\\
&\lesssim_j{\theta_b^{-2j}\over|\cos((1-a)\theta_b)|^2}\int_{\Re\lambda=a-1}|\hat f_1|^2\d\Im\lambda
={\theta_b^{-2j}\over|\cos((1-a)\theta_b)|^2}\|f_1\|_{\partial_1\Omega,a-{1\over 2}}^2,
\end{align*}
where we used $\<\lambda\>\leq 1+|\Re\lambda|+|\Im\lambda|$, $\theta_b<\pi/2$, $|\Re\lambda|\theta_b<\pi/2$ and $\sup_{x\geq 0}(1+x)^{2j}e^{-2x}<\infty$.
\end{proof}

\begin{lemma}[Neumann problem on the wedge]\label{Neumann wedge problem}
Consider the following three elliptic systems
\begin{alignat*}{3}
\lpc(q_{\mb F},q_{f_0},q_{f_1})&=(\grad\cdot(r^{-1}\mb F),0,0)\qquad&&\text{on}\qquad\Omega\\
\partial_\theta(q_{\mb F},q_{f_0},q_{f_1})&=(F_\theta,f_0,0)\qquad&&\text{on}\qquad\partial_0\Omega\\
\partial_\theta(q_{\mb F},q_{f_0},q_{f_1})&=(F_\theta,0,f_1)\qquad&&\text{on}\qquad\partial_1\Omega.
\end{alignat*}
Let the weight $a\in\R$ satisfy $0<|1-a|\theta_b<\pi$, and recall the sine resonance factor $\Gamma_{c,\theta}$ of \eqref{sine resonance factor}. Then these systems have solutions $q_{\mb F},q_{f_0},q_{f_1}$ with regularity dependent on the regularity of $\mb F$, $f_0$, $f_1$ respectively, in the following sense, for $j\geq 0$ and $l\in\{0,1\}$:
\begin{align*}
\|\dr^j\partial_\theta^l q_{\mb F}\|_{\Omega,a}&\lesssim\Gamma_{1-a,\theta_b}^{-1}\|\dr^{j+l-1}\mb F\|_{\Omega,a}\\
\|(r\grad)^{j+1}q_{\mb F}\|_{\Omega,a}&\lesssim_j\Gamma_{1-a,\theta_b}^{-1}\|\mb F\|_{\Omega,a,0,j}\\
\|(r\grad)^jq_{f_\bullet}\|_{\Omega,a}&\lesssim_{j,|1-a|^{-1}}\norm{\set{\theta_b^{1\over 2},\dr^{-{1\over 2}}}\dr^{j-1}f_\bullet}_{\partial_\bullet\Omega,a-{1\over 2}} +{1\over\theta_b^{1/2}\Gamma_{1-a,\theta_b}} \|\dr^{j-2}f_\bullet\|_{\partial_\bullet\Omega,a-{1\over 2}},
\end{align*}
where $\{\ph\}$ is a choice as in Definition \ref{notation}.
\end{lemma}
\begin{proof}
Mellin transform $r\mapsto\lambda$ turns the system of the statement, in the form $\lpc q=\grad\cdot(r^{-1}\mb F)$ on $\Omega$ with $\partial_\theta q=f_\bullet+F_\theta$ on $\partial_\bullet\Omega$ for $\bullet\in\{0,1\}$, into the system
\begin{align}
(\lambda^2+\partial_\theta^2)\hat q(\lambda,\theta)&=\lambda\hat F_r(\lambda,\theta)+\partial_\theta\hat F_\theta(\lambda,\theta)\nonumber\\
\partial_\theta\hat q(\lambda,\theta_b)&=\hat f_1(\lambda)+\hat F_\theta(\lambda,\theta_b)\label{Neumann Mellin system}\\
\partial_\theta\hat q(\lambda,0)&=\hat f_0(\lambda)+\hat F_\theta(\lambda,0),\nonumber
\end{align}
again a two-point boundary value problem in $\theta$ for every fixed $\lambda$, and we solve it as before, deriving the bulk part first.

Let
\begin{align*}
h_k=\sqrt{2-\delta_{0k}\over\theta_b}\cos\brac{k\pi{\theta\over\theta_b}},\qquad k\in\N_0,
\end{align*}
be the $L^2([0,\theta_b])$-orthonormal eigenfunctions of $\partial_\theta^2$ with boundary condition $\partial_\theta h(\theta_b)=0=\partial_\theta h(0)$, with eigenvalues $-\Lambda_k^2=-(k\pi)^2\theta_b^{-2}$; and, exactly as in Lemma \ref{mixed wedge problem}, let
\begin{align*}
\hat F_{r,k}=\<h_k,\hat F_{r}\>_{L^2([0,\theta_b])},\qquad\qquad
\hat F_{\theta,k}=\inner{{h_k'\over\Lambda_k},\hat F_{\theta}}_{L^2([0,\theta_b])}\quad(k\geq 1),\qquad\qquad\hat F_{\theta,0}:=0,
\end{align*}
the functions $h_k'/\Lambda_k$ for $k\geq 1$ being the orthonormal eigenfunctions of $\partial_\theta^2$ with the complementary boundary condition $h(0)=0=h(\theta_b)$, and $(h_k'/\Lambda_k)'=-\Lambda_kh_k$. The $k=0$ term of the $\hat F_\theta$ sum does not occur, being multiplied by $\Lambda_0=0$.

Write $\hat q_{\mb F}=\sum_kc_kh_k$ for the solution of \eqref{Neumann Mellin system} with $\hat f_0=\hat f_1=0$. Testing the equation against $h_k$ and integrating by parts twice, as in the proof of Lemma \ref{mixed wedge problem}, gives
\begin{align*}
(\lambda^2-\Lambda_k^2)c_k+\eva{\brac{\partial_\theta\hat q_{\mb F}h_k-\hat q_{\mb F}h_k'}}_0^{\theta_b}
=\lambda\hat F_{r,k}-\Lambda_k\hat F_{\theta,k}+\eva{\hat F_\theta h_k}_0^{\theta_b}.
\end{align*}
Here $h_k'$ vanishes at both endpoints, so the terms carrying it drop out, and what is left is $\partial_\theta\hat q_{\mb F}h_k|_0^{\theta_b}$ on the left and $\hat F_\theta h_k|_0^{\theta_b}$ on the right; these cancel each other precisely under the boundary condition $\partial_\theta q_{\mb F}=F_\theta$ of the statement, which is now imposed on both boundaries and is, as in the proof of Lemma \ref{mixed wedge problem}, the natural one for the divergence form RHS. We are left with $(\lambda^2-\Lambda_k^2)c_k=\lambda\hat F_{r,k}-\Lambda_k\hat F_{\theta,k}$, and therefore with the solution
\begin{align*}
\hat q&=\underbrace{{\hat f_0\over\lambda\sin(\lambda\theta_b)}\cos(\lambda(\theta_b-\theta))}_{\hat q_{f_0}}\underbrace{-{\hat f_1\over\lambda\sin(\lambda\theta_b)}\cos(\lambda\theta)}_{\hat q_{f_1}}
+\underbrace{\sum_{k=0}^\infty{(\lambda\hat F_{r,k}-\Lambda_k\hat F_{\theta,k})h_k\over\lambda^2-\Lambda_k^2}}_{\hat q_{\mb F}},
\end{align*}
the first two terms being the solutions of the homogeneous equation that carry $\hat f_0$ and $\hat f_1$. The hypothesis on $a$ keeps the line $\Re\lambda=a-1$ off $\lambda=0$, where $\lambda^2-\Lambda_0^2$ and $\sin(\lambda\theta_b)$ both vanish, and off $\pm\Lambda_1=\pm\pi/\theta_b$.

Note that for $l\in\{0,1,2\}$ and $\Re\lambda=a-1$,
\begin{align}\label{Neumann multiplier bound}
\sup_{k\geq 0}\abs{\{\lambda,\Lambda_k\}^l\over\lambda^2-\Lambda_k^2}&\lesssim\Gamma_{1-a,\theta_b}^{-1}|\lambda|^{l-2},
\end{align}
the term $k=0$ being exactly $|\lambda|^{l-2}$ with the choice $\lambda$ and $0$ with the choice $\Lambda_0=0$. For $k\geq 1$ we claim that
\begin{align}\label{Neumann denominator bound}
|\lambda^2-\Lambda_k^2|\geq{\Gamma_{1-a,\theta_b}\over 4\pi}\brac{\Lambda_k^2+|\lambda|^2},
\end{align}
which suffices, since $\Lambda_k^2+|\lambda|^2\geq\{|\lambda|,\Lambda_k\}^l|\lambda|^{2-l}$. Since $\lambda^2$ is even in $\lambda$ we may assume $\Re\lambda\geq 0$; the hypothesis $|1-a|\theta_b<\pi$ says $\Re\lambda<\Lambda_1$. Taking square roots in $\alpha^2+\beta^2\geq{1\over 2}(\alpha+\beta)^2$, as in the proof of Lemma \ref{cosh bound}, gives for every $k\geq 1$
\begin{align*}
|\lambda+\Lambda_k|=\sqrt{(\Re\lambda+\Lambda_k)^2+(\Im\lambda)^2}\geq\sqrt{\Lambda_k^2+(\Im\lambda)^2}\geq{1\over\sqrt 2}\brac{\Lambda_k+|\Im\lambda|}.
\end{align*}
For the other factor, $\Lambda_1-\Re\lambda=\Lambda_1(\pi-\Re\lambda\theta_b)/\pi\geq\pi^{-1}\Gamma_{1-a,\theta_b}\Lambda_1$, while $\Lambda_k-\Re\lambda\geq\Lambda_k-\Lambda_1\geq\Lambda_k/2\geq\pi^{-1}\Gamma_{1-a,\theta_b}\Lambda_k$ for $k\geq 2$, the last step because $\Gamma_{1-a,\theta_b}\leq 1$; hence, for every $k\geq 1$,
\begin{align*}
|\lambda-\Lambda_k|=\sqrt{(\Lambda_k-\Re\lambda)^2+(\Im\lambda)^2}\geq{1\over\sqrt 2}\brac{{\Gamma_{1-a,\theta_b}\over\pi}\Lambda_k+|\Im\lambda|}\geq{\Gamma_{1-a,\theta_b}\over\pi\sqrt 2}\brac{\Lambda_k+|\Im\lambda|}.
\end{align*}
Multiplying the two gives $|\lambda^2-\Lambda_k^2|\geq(2\pi)^{-1}\Gamma_{1-a,\theta_b}(\Lambda_k+|\Im\lambda|)^2$, and \eqref{Neumann denominator bound} follows, since $\Lambda_k^2+|\lambda|^2\leq 2(\Lambda_k^2+(\Im\lambda)^2)\leq 2(\Lambda_k+|\Im\lambda|)^2$, the first step using $\Re\lambda<\Lambda_1\leq\Lambda_k$. So we have, by \eqref{Neumann multiplier bound}, for $j\geq 0$ and $l\in\{0,1\}$,
\begin{align*}
\|\dr^j\partial_\theta^l q_{\mb F}\|_{\Omega,a}^2&=\int_{\Re\lambda=a-1}\int_0^{\theta_b}|\lambda^j\partial_\theta^l\hat q_{\mb F}|^2\d\theta\d\Im\lambda\\
&=\int_{\Re\lambda=a-1}\int_0^{\theta_b} \abs{\lambda^j\partial_\theta^l\sum_{k=0}^\infty{(\lambda\hat F_{r,k}-\Lambda_k\hat F_{\theta,k})h_k\over\lambda^2-\Lambda_k^2}}^2\d\theta\d\Im\lambda\\
&\leq\int_{\Re\lambda=a-1}\int_0^{\theta_b} \abs{\lambda^j\sum_{k=l}^\infty{(\lambda\hat F_{r,k}-\Lambda_k\hat F_{\theta,k})\Lambda_k^{l}(h_k^{(l)}/\Lambda_k^{l})\over\lambda^2-\Lambda_k^2}}^2\d\theta\d\Im\lambda\\
&\lesssim\Gamma_{1-a,\theta_b}^{-2}\int_{\Re\lambda=a-1}\abs{|\lambda|^{j}|\lambda|^{l-1}}^2 \sum_{k=0}^\infty\brac{|\hat F_{r,k}|^2+|\hat F_{\theta,k}|^2}\d\Im\lambda\\
&\leq\Gamma_{1-a,\theta_b}^{-2}\int_{\Re\lambda=a-1}|\lambda|^{2(j+l-1)}\int_0^{\theta_b}|\hat{\mb F}|^2\d\theta\d\Im\lambda
\leq\Gamma_{1-a,\theta_b}^{-2}\|\dr^{j+l-1}\mb F\|_{\Omega,a}^2.
\end{align*}
A similar computation as \eqref{higher bulk 1} and \eqref{higher bulk 2}, with the powers of $\lambda$ carried in the same way and with the factor $\Gamma_{1-a,\theta_b}^{-1}$ in place of $|\cos((1-a)\theta_b)|^{-1}$, gives us $\|(r\grad)^{j+1}q_{\mb F}\|_{\Omega,a}\lesssim_j\Gamma_{1-a,\theta_b}^{-1}\|\mb F\|_{\Omega,a,0,j}$.

For the boundary data we again need the angular integral, now on the wider strip $|\Re\lambda|\theta_b=|1-a|\theta_b<\pi$: Lemma \ref{angular integral} applies with $c=\pi$ and gives
\begin{align*}
\int_0^{\theta_b}e^{-2|\Im\lambda|\theta}\d\theta\leq(\pi+1)\set{\theta_b,|\lambda|^{-1}},
\end{align*}
the wider strip costing the larger constant. As in the proof of Lemma \ref{mixed wedge problem}, we bound the $\dr^m\partial_\theta^nq_{f_0}$ with $m+n\leq j$. Since $|\lambda|\geq|1-a|>0$ on the line $\Re\lambda=a-1$, \eqref{Plancherel's identity} gives $\|\dr^m\partial_\theta^nq_{f_0}\|_{\Omega,a}\leq|1-a|^{m+n-j}\|\dr^{j-n}\partial_\theta^nq_{f_0}\|_{\Omega,a}$, so it suffices to take $m+n=j$. We have
\begin{align*}
\|\dr^m\partial_\theta^nq_{f_0}\|_{\Omega,a}^2
&=\int_{\Re\lambda=a-1}\int_0^{\theta_b}|\lambda^m\partial_\theta^n\hat q_{f_0}|^2\d\theta\d\Im\lambda\\
&=\int_{\Re\lambda=a-1}\int_0^{\theta_b}\abs{\hat f_0{\lambda^m\partial_\theta^n\over\lambda\sin(\lambda\theta_b)}\cos(\lambda(\theta_b-\theta))}^2\d\theta\d\Im\lambda\\
\intertext{so that, by Lemma \ref{cosh bound},}
&\lesssim\int_{\Re\lambda=a-1}\mb 1_{\{|\Im\lambda\theta_b|\geq 1\}}\int_0^{\theta_b}|\lambda|^{2(j-1)}|\hat f_0|^2e^{-2|\Im\lambda|\theta}\d\theta\d\Im\lambda\\
&\quad+\Gamma_{1-a,\theta_b}^{-2}\int_{\Re\lambda=a-1}\mb 1_{\{|\Im\lambda\theta_b|<1\}}\int_0^{\theta_b}{|\lambda|^{2(j-2)}\over|\theta_b|^2}|\hat f_0|^2\d\theta\d\Im\lambda\\
&\lesssim\int_{\Re\lambda=a-1}\set{\theta_b,|\lambda|^{-1}}|\lambda|^{2(j-1)}|\hat f_0|^2\d\Im\lambda
+\Gamma_{1-a,\theta_b}^{-2}\int_{\Re\lambda=a-1}{|\lambda|^{2(j-2)}\over|\theta_b|}|\hat f_0|^2\d\Im\lambda\\
&\leq\norm{r^{{1\over 2}-a}\set{\theta_b^{1\over 2},\dr^{-{1\over 2}}}\dr^{j-1}f_0}_{L^2([0,\infty))}^2 +{1\over\theta_b\Gamma_{1-a,\theta_b}^{2}}\|r^{{1\over 2}-a}\dr^{j-2}f_0\|_{L^2([0,\infty))}^2\\
&=\norm{\set{\theta_b^{1\over 2},\dr^{-{1\over 2}}}\dr^{j-1}f_0}_{\partial_0\Omega,a-{1\over 2}}^2 +{1\over\theta_b\Gamma_{1-a,\theta_b}^{2}}\|\dr^{j-2}f_0\|_{\partial_0\Omega,a-{1\over 2}}^2.
\end{align*}
The estimates for $q_{f_1}$ are the same, with $\cos(\lambda\theta)$ in place of $\cos(\lambda(\theta_b-\theta))$ and $e^{2|\Im\lambda|(\theta-\theta_b)}$ in place of $e^{-2|\Im\lambda|\theta}$, that is, under the reflection $\theta\mapsto\theta_b-\theta$, which exchanges the two boundaries; the powers of $\lambda$ are unchanged.
\end{proof}

In the following we define two projections that correspond to the orthogonality of our system.

\begin{definition}[Helmholtz and auxiliary projection]\label{the two projections}
For any sufficiently smooth $\mb v$, we define
\begin{enumerate}
\item \ENEQ\begin{align*}
\P\mb v:=\mb v-\grad q[\mb v],
\end{align*}
where $q[\mb v]$ is the solution to the system $\lpc q[\mb v]=\grad\cdot\mb v$ with $q[\mb v]|_{\partial_1\Omega}=0$ and $\partial_\theta q[\mb v]|_{\partial_0\Omega}=rv_\theta|_{\partial_0\Omega}$ given by Lemma \ref{mixed wedge problem}.
\item \ENEQ\begin{align*}
\P^\theta\mb v:=\mb v-\grad q^\theta[\mb v],
\end{align*}
where $q^\theta[\mb v]$ is the solution to the system $\lpc q^\theta[\mb v]=\grad\cdot\mb v$ with $\partial_\theta q^\theta[\mb v]|_{\partial_1\Omega}=rv_\theta|_{\partial_1\Omega}$ and $\partial_\theta q^\theta[\mb v]|_{\partial_0\Omega}=rv_\theta|_{\partial_0\Omega}$ given by Lemma \ref{Neumann wedge problem}.
\end{enumerate}
\end{definition}

\begin{remark}\label{P remark}
Note that we have $\grad\cdot\P\mb v=0$ and $(\P\mb v)_\theta|_{\partial_0\Omega}=0$. Moreover, for any sufficiently smooth $\mb w$ such that $\grad\cdot\mb w=0$ and $w_\theta|_{\partial_0\Omega}=0$, we have $\<\mb w,\mb v\>_{L^2(\Omega)}=\<\mb w,\P\mb v\>_{L^2(\Omega)}$ because the divergence theorem gives $\<\mb w,\grad q[\mb v]\>_{L^2(\Omega)}=0$.

Note that we have $\grad\cdot\P^\theta\mb v=0$ and $(\P^\theta\mb v)_\theta|_{\partial\Omega}=0$. In particular, we have $\<\grad p,\P^\theta\mb v\>_{L^2(\Omega)}=0$ for $p$ sufficiently smooth.
\end{remark}

\begin{lemma}[Mixed harmonic extension]\label{mixed harmonic extension}
For any smooth enough function $\phi$ on $\partial_1\Omega$, we define its mixed harmonic extension $\E\phi$ on $\Omega$ to be the solution to the system $\lpc\E\phi=0$ with $\E\phi|_{\partial_1\Omega}=\phi$ and $\partial_\theta\E\phi|_{\partial_0\Omega}=0$ given by Lemma \ref{mixed wedge problem}, that is, $\E\phi=q_{f_1}$ with $f_1=\phi$. Let the weight $a\in\R$ satisfy $|a-1|\theta_b<\pi/2$. Then for $k,j,m\in\N_0$ we have
\begin{align*}
\|(r\grad)^k\E\phi\|_{\Omega,a}&\lesssim_k{1\over|\cos((1-a)\theta_b)|}\|\phi\|_{\partial_1\Omega,a-{1\over 2},k-{1\over 2}}\\
\|r\grad\E\phi\|_{\Omega,a,j,m}&\lesssim_{j,m}{1\over|\cos((1-a)\theta_b)|}\|\dr\phi\|_{\partial_1\Omega,a-{1\over 2},\max\{j+m-{1\over 2},0\}},
\end{align*}
whenever the RHSs are finite. If moreover $a>1$, then the finiteness of $\|\dr\phi\|_{\partial_1\Omega,a-1/2}$ gives $\phi$ a limit $\phi(0)$ at the corner, and we set $\E\phi:=\phi(0)+\E(\phi-\phi(0))$; the second estimate then holds for such $\phi$ as well.
\end{lemma}
\begin{proof}
The first estimate is the bulk estimate for $q_{f_1}$ in Lemma \ref{mixed wedge problem} with the choice $\dr^{-1/2}$, since, by Definition \ref{fractional counts}, $\|\phi\|_{\partial_1\Omega,a-1/2,k-1/2}$ contains the term $\|\dr^{k-1/2}\phi\|_{\partial_1\Omega,a-1/2}$, and its integer terms bound the $\|\dr^{j'-1/2}\phi\|_{\partial_1\Omega,a-1/2}$ with $j'<k$, by the pointwise bound $|\lambda|^{2c}\leq|\lambda|^{2\lfloor c\rfloor}+|\lambda|^{2\lceil c\rceil}$ on the line \eqref{weight and line}. For the second, the second inequality of \eqref{exchange} and the definition of the norm give
\begin{align*}
\|r\grad\E\phi\|_{\Omega,a,j,m}\lesssim_{j,m}\|r\grad\E\phi\|_{\Omega,a,0,j+m}\leq\sum_{k=1}^{j+m+1}\|(r\grad)^k\E\phi\|_{\Omega,a},
\end{align*}
and we estimate each term by the same estimate of Lemma \ref{mixed wedge problem}, now with the choice $\theta_b^{1/2}$ in its term $j'=1$ and $\dr^{-1/2}$ in the others, which gives
\begin{align*}
\|r\grad\E\phi\|_{\Omega,a,j,m}\lesssim_{j,m}{1\over|\cos((1-a)\theta_b)|}\brac{\theta_b^{1\over 2}\|\dr\phi\|_{\partial_1\Omega,a-{1\over 2}}+\sum_{k=2}^{j+m+1}\|\dr^{k-{3\over 2}}(\dr\phi)\|_{\partial_1\Omega,a-{1\over 2}}}.
\end{align*}
In the bracket, the first term is at most $(\pi/2)^{1/2}\|\dr\phi\|_{\partial_1\Omega,a-1/2}$, the count $0$ term of the norm on the right of the statement, since $\theta_b<\pi/2$. The exponents $c=k-3/2$ of the sum, which is empty when $j+m=0$, run through the half-integers from $1/2$ to $j+m-1/2$: the last one gives the fractional term of that norm, and the others are bounded by its integer terms, by the same pointwise bound and since $\lceil c\rceil\leq j+m-1$ for $c<j+m-1/2$.

For the last statement, Proposition \ref{Hardy's inequality}, applied with $k=a-1>0$, gives the limit $\phi(0)$, which need not vanish, together with $\|\phi-\phi(0)\|_{\partial_1\Omega,a-1/2}\leq(a-1)^{-1}\|\dr\phi\|_{\partial_1\Omega,a-1/2}$, so that $\E(\phi-\phi(0))$ is defined; a constant being harmonic with vanishing Neumann datum, $\E\phi$ solves the same system, and the second estimate for $\phi-\phi(0)$ is the one for $\phi$, since $r\grad\E\phi=r\grad\E(\phi-\phi(0))$ and $\dr\phi=\dr(\phi-\phi(0))$. This is the case of the deviation $\chi$ of the profile in Lemma \ref{regularised flattening}, whose value at the corner is the position of the contact point.
\end{proof}

\subsection{Fixed domain transformation}\label{fixed domain transformation}

We now have the tools to carry out the transformation announced at the beginning of this section. The map from the reference wedge $\Omega$ onto $\Omega_t$, whose inverse flattens $\Omega_t$, is built from the mixed harmonic extension of the deviation $\chi=X-b^{-1}y$ of the profile from the exact wedge. This is why the elliptic theory of the previous subsection had to come first: the regularity of this map is governed by Lemma \ref{mixed harmonic extension}. The smallness that makes it a diffeomorphism is a smallness condition on the profile alone.

\begin{lemma}[Regularised flattening]\label{regularised flattening}
Let $\theta_b\in(0,\pi/2)$, let $\Omega_t=\{(x,y):y>0,\ x>X(y,t)\}$ be the fluid domain of System \ref{dimensional system}, and define
\begin{align*}
\bs\Psi_t(x,y):=(x+\E\chi(x,y,t),y),\qquad\text{ where }\qquad\chi(y,t):=X(y,t)-b^{-1}y,
\end{align*}
in which $\E$ is the mixed harmonic extension of Lemma \ref{mixed harmonic extension}, $\chi(\ph,t)$ being regarded as a function on $\partial_1\Omega$ through \eqref{r and y on the free boundary}. There exists $\epsilon_{\bs\Psi}=\epsilon_{\bs\Psi}(\theta_b)>0$ such that whenever $\|\chi'(t)\|_{\partial_1\Omega,1/2,3/2}\leq\epsilon_{\bs\Psi}$, $\bs\Psi_t$ is a diffeomorphism $\Omega\to\Omega_t$ which extends to a homeomorphism of $\overline\Omega$ onto $\overline\Omega_t$ mapping $\partial_0\Omega$ onto $\partial_0\Omega_t$ and $\partial_1\Omega$ onto $\partial_1\Omega_t$. Its Jacobian matrix and determinant are
\begin{align}\label{flattening Jacobian}
D\bs\Psi_t=\mat{1+\nu_1&\nu_2\\0&1},\qquad\qquad\det D\bs\Psi_t=1+\nu_1>{1\over 2},\qquad\text{ where }\qquad\bs\nu:=\grad\E\chi.
\end{align}
\end{lemma}
\begin{proof}
By Lemma \ref{embeddings}(ii), the weight shift of Lemma \ref{elementary properties}(ii), the second estimate of Lemma \ref{mixed harmonic extension} with $\phi=\chi$, $a=2$ and $j=m=1$, and $\dr\chi=(\sin\theta_b)r\chi'$, which follows from \eqref{r and y on the free boundary}, we have
\begin{align*}
\|\grad\E\chi\|_{L^\infty(\Omega)}\lesssim_{\theta_b}\|\grad\E\chi\|_{\Omega,1,1,1}\lesssim\|r\grad\E\chi\|_{\Omega,2,1,1}\lesssim_{\theta_b}\|\dr\chi\|_{\partial_1\Omega,{3\over 2},{3\over 2}}\lesssim|\sin\theta_b|\|\chi'\|_{\partial_1\Omega,{1\over 2},{3\over 2}};
\end{align*}
the second step is not an identity, since multiplication by $r$ does not commute with $\dr$ and $r\grad$. Picking $\epsilon_{\bs\Psi}$ small enough, we have $\|\grad\E\chi\|_{L^\infty(\Omega)}<1/2$. In particular $\E\chi$ is Lipschitz continuous on $\Omega$, and so extends continuously to $\overline\Omega$, with $\E\chi(\mb 0,t)=\chi(0,t)$, which need not vanish.

For each fixed $y\geq 0$ the map $x\mapsto x+\E\chi(x,y,t)$ has derivative $1+\partial_x\E\chi>1/2$, so it is a strictly increasing bijection of $[b^{-1}y,\infty)$ onto $[X(y,t),\infty)$: its value at $x=b^{-1}y$ is $b^{-1}y+\chi(y,t)=X(y,t)$, since $\E\chi=\chi$ on $\partial_1\Omega$, and its limit as $x\to\infty$ is $\infty$. Hence $\bs\Psi_t$ is a bijection of $\overline\Omega$ onto $\overline\Omega_t$, mapping $\Omega$ onto $\Omega_t$, $\partial_0\Omega$ onto $\partial_0\Omega_t$ and, by the same identity, $\partial_1\Omega$ onto $\partial_1\Omega_t$. Since $\overline\Omega$ is convex, $|\E\chi(\mb x,t)-\E\chi(\mb x',t)|\leq|\mb x-\mb x'|/2$, and hence $|\bs\Psi_t(\mb x)-\bs\Psi_t(\mb x')|\geq|\mb x-\mb x'|/2$, for $\mb x,\mb x'\in\overline\Omega$; so the inverse of $\bs\Psi_t$ is Lipschitz, and $\bs\Psi_t$ is a homeomorphism of $\overline\Omega$ onto $\overline\Omega_t$. Finally, differentiating the definition of $\bs\Psi_t$ gives \eqref{flattening Jacobian}, with $\det D\bs\Psi_t=1+\partial_x\E\chi>1/2$, so that $\bs\Psi_t$ is a diffeomorphism of the open set $\Omega$ onto $\Omega_t$ by the inverse function theorem.
\end{proof}

\begin{lemma}[Piola transform]\label{Piola transform}
Take the setting of Lemma \ref{regularised flattening}. Suppose $(\partial_t+u_y\partial_y)X=u_x$ on $\partial_1\Omega_t$. Define a coordinate $(\tilde{x},\tilde{y},\tilde{t})(x,y,t)$ via
\begin{align*}
(\tilde{x},\tilde{y})=\bs\Psi_t^{-1}(x,y),\qquad\qquad
\tilde{t}=t,
\end{align*}
so that $(\tilde x,\tilde y)$ are the coordinates on the reference wedge $\Omega$, which were called $(x,y)$ in Lemma \ref{regularised flattening}, and $(x,y)$ those on $\Omega_t$. Write $\phi(\tilde{y},\tilde{t})=-\chi'(\tilde{y},\tilde{t})$, regarded along $\partial_1\Omega$ as a function of $r$ through $\tilde{y}=r\sin\theta_b$, and $\bs\nu[\phi]=-(\sin\theta_b)\tilde{\grad}\E\dr^{-1}(r\phi)=\tilde{\grad}\E\chi$, the field $\bs\nu$ of \eqref{flattening Jacobian}; and let
\begin{align*}
\tilde{\mb u}(\tilde{x},\tilde{y};\tilde{t})=|\det D\bs\Psi_t|(D\bs\Psi_t)^{-1}\mb u\brac{\bs\Psi_{\tilde{t}}(\tilde{x},\tilde{y}),\tilde{t}}
\end{align*}
be the Piola transform of $\mb u$ \cite[Chapter~1]{Ciarlet}. Let $(\tilde{x},\tilde{y})=(r\cos\theta,r\sin\theta)$. Then
\begin{align*}
\grad&=\tilde{\grad}+\tilde{\grad}_{\phi}\qquad\text{where}\qquad \tilde{\grad}_{\phi}:=\mat{(1+\nu_1)^{-1}-1\\-(1+\nu_1)^{-1}\nu_2}\partial_{\tilde{x}}=-(1+\nu_1)^{-1}\bs\nu\partial_{\tilde{x}}\\
\mb u&=\tilde{\mb u}+\tilde{\mb u}_\phi\qquad\text{where}\qquad\tilde{\mb u}_\phi:=\mat{(1+\nu_1)^{-1}\nu_2\\(1+\nu_1)^{-1}-1}\tilde{u}_{\tilde{y}}=-(1+\nu_1)^{-1}\bs\nu^\perp\tilde{u}_{\tilde{y}}\\
\partial_{\tilde{t}}\phi&=(\sin\theta_b)^{-2}\partial_r\tilde{u}_\theta(r,\theta_b),
\end{align*}
where $\theta_b=\arctan b$. Moreover, we have $\tilde{\grad}\cdot\tilde{\mb u}=(1+\nu_1)\grad\cdot\mb u$.
\end{lemma}
\begin{proof}
We have
\begin{subequations}\label{transformation matrix}
\begin{align}
{\partial(x,y,t)\over\partial(\tilde{x},\tilde{y},\tilde{t})}&=\bmat{{\partial x\over\partial\tilde x}&{\partial x\over\partial\tilde y}&{\partial x\over\partial\tilde t}\\{\partial y\over\partial\tilde x}&{\partial y\over\partial\tilde y}&{\partial y\over\partial\tilde t}\\{\partial t\over\partial\tilde x}&{\partial t\over\partial\tilde y}&{\partial t\over\partial\tilde t}}
=\bmat{1+\partial_{\tilde{x}}\E\chi&\partial_{\tilde{y}}\E\chi&\partial_{\tilde{t}}\E\chi\\0&1&0\\0&0&1}
=\bmat{1+\nu_1&\nu_2&\partial_{\tilde{t}}\E\chi\\0&1&0\\0&0&1}\\
{\partial(\tilde{x},\tilde{y},\tilde{t})\over\partial(x,y,t)}
&=\bmat{{\partial\tilde{x}\over\partial x}&{\partial\tilde{x}\over\partial y}&{\partial\tilde{x}\over\partial t}\\{\partial\tilde{y}\over\partial x}&{\partial\tilde{y}\over\partial y}&{\partial\tilde{y}\over\partial t}\\{\partial\tilde{t}\over\partial x}&{\partial\tilde{t}\over\partial y}&{\partial\tilde{t}\over\partial t}}
=\bmat{(1+\nu_1)^{-1}&-(1+\nu_1)^{-1}\nu_2&-(1+\nu_1)^{-1}\partial_{\tilde{t}}\E\chi\\0&1&0\\0&0&1},
\intertext{the upper left $2\times 2$ blocks being $D\bs\Psi_t$ and $(D\bs\Psi_t)^{-1}$. Since $|\det D\bs\Psi_t|=1+\nu_1$ by \eqref{flattening Jacobian}, the definition of $\tilde{\mb u}$ and its inverse read}
\tilde{\mb u}&=\mat{1&-\nu_2\\0&1+\nu_1}\mb u=\mat{u_x-\nu_2u_y\\(1+\nu_1)u_y}\\
\mb u&=\mat{1&(1+\nu_1)^{-1}\nu_2\\0&(1+\nu_1)^{-1}}\tilde{\mb u}=\mat{\tilde u_{\tilde x}+(1+\nu_1)^{-1}\nu_2\tilde u_{\tilde y}\\(1+\nu_1)^{-1}\tilde u_{\tilde y}}.
\end{align}
\end{subequations}
\begin{align*}
\tilde{u}_\theta(\tilde{x},\tilde{y})=(\tilde{\mb e}_\theta\cdot\tilde{\mb u})(\tilde{x},\tilde{y})
={1\over\sqrt{\tilde{x}^2+\tilde{y}^2}}\mat{-\tilde{y}\\\tilde{x}}\cdot\mat{u_x-\nu_2u_y\\(1+\nu_1)u_y}
\end{align*}
so that, with $\bs\nu=\tilde\grad\E\chi$,
\begin{align}
\tilde{u}_\theta(b^{-1}\tilde{y},\tilde{y})&={1\over\sqrt{1+b^2}}((1+\nu_1)u_y-b(u_x-\nu_2u_y))
={1\over\sqrt{1+b^2}}((1+\nu_1+b\nu_2)u_y-bu_x)\nonumber\\
&={1\over\sqrt{1+b^2}}((1+(\partial_{\tilde{x}}+b\partial_{\tilde{y}})\E\chi)u_y-bu_x)
={1\over\sqrt{1+b^2}}\brac{\brac{1+\sqrt{1+b^2}\partial_r\chi}u_y-bu_x}\nonumber\\
&={1\over\sqrt{1+b^2}}\brac{\brac{1+\sqrt{1+b^2}(X'-b^{-1})\sin\theta_b}u_y-bu_x}\nonumber\\
&={b\over\sqrt{1+b^2}}\brac{X'u_y-u_x}={-b\partial_tX\over\sqrt{1+b^2}},\label{transformed u_theta}
\end{align}
where the second line uses the definition of $\bs\nu$, together with $\partial_{\tilde x}+b\partial_{\tilde y}=\sqrt{1+b^2}\partial_r$ along $\partial_1\Omega$ and $\E\chi=\chi$ there, the third uses $\chi'=X'-b^{-1}$ and \eqref{r and y on the free boundary}, and the last uses $\sqrt{1+b^2}\sin\theta_b=b$ and the profile evolution equation $\partial_tX=u_x-u_yX'$. We have
\begin{align*}
b\partial_{\tilde{t}}\chi
&=\brac{{\partial x\over\partial{\tilde{t}}}\partial_x+{\partial y\over\partial{\tilde{t}}}\partial_y+{\partial t\over\partial{\tilde{t}}}\partial_t}(bX-y)=b\partial_tX=-\sqrt{1+b^2}\tilde{u}_\theta
\end{align*}
on $\partial_1\Omega$. Differentiating this identity in $r$ along $\partial_1\Omega$, with $\phi=-\chi'$, \eqref{r and y on the free boundary} and $\sin\theta_b=b/\sqrt{1+b^2}$, gives $\partial_{\tilde t}\phi=(\sin\theta_b)^{-2}\partial_r\tilde u_\theta(r,\theta_b)$, the evolution equation for $\phi$. Finally, by the chain rule $\partial_{\tilde x}=(1+\nu_1)\partial_x$ and $\partial_{\tilde y}=\nu_2\partial_x+\partial_y$, that is $\tilde\grad=(D\bs\Psi_t)^T\grad$, we have
\begin{align*}
\tilde{\grad}\cdot\tilde{\mb u}&=\mat{(1+\nu_1)\partial_x\\\nu_2\partial_x+\partial_y}\cdot\mat{u_x-\nu_2u_y\\(1+\nu_1)u_y}
=(1+\nu_1)\brac{\grad\cdot\mb u-\partial_x(\nu_2u_y)+\nu_2\partial_xu_y}+u_y(\nu_2\partial_x+\partial_y)\nu_1\\
&=(1+\nu_1)\brac{\grad\cdot\mb u-u_y\partial_x\nu_2}+u_y(\nu_2\partial_x+\partial_y)\nu_1
=(1+\nu_1)\grad\cdot\mb u,
\end{align*}
where we used that
\begin{align*}
0&=\tilde{\grad}^\perp\cdot\tilde{\grad}\E\chi=\tilde{\grad}^\perp\cdot\bs\nu=(1+\nu_1)\partial_x\nu_2-(\nu_2\partial_x+\partial_y)\nu_1.\qedhere
\end{align*}
\end{proof}

Under this transformation in Lemma \ref{Piola transform}, System \ref{dimensionless system} becomes System \ref{reference wedge system} below. All the variables are with tildes, but for simplicity of notation we suppress the tildes. From now on we regard $\phi$ as a function of $r$ along $\partial_1\Omega$, the two descriptions being related by \eqref{r and y on the free boundary}, so that $\partial_r\phi=(\sin\theta_b)\phi'$. And note that a contact angle of $\theta_b=\arctan b$ is the condition that $X'(0,t)=b^{-1}$, or equivalently $\phi(0,t)=0$.

\begin{system}[Wedge-like Stokes flow on the fixed wedge]\label{reference wedge system}
\begin{alignat*}{3}
\lpc\mb u-\grad p&=\mb N[\phi,\mb u,p]\qquad&&\text{on}\quad\Omega\tag{Stokes equation}\\
\grad\cdot\mb u&=0\qquad&&\text{on}\quad\Omega\tag{incompressibility}\\
u_\theta&=0\qquad&&\text{on}\quad\partial_0\Omega\tag{no penetration}\\
u_r-r^{-1}\partial_\theta u_r&=0\qquad&&\text{on}\quad\partial_0\Omega\tag{Navier slip}\\
\mathbb T[\mb u,p]\mb e_\theta-(\sin\theta_b)^2(\partial_r\phi)\mb e_\theta&=\mb N_\partial[\phi,\mb u,p]\qquad&&\text{on}\quad\partial_1\Omega\tag{surface tension}\\
\partial_{t}\phi-(\sin\theta_b)^{-2}\partial_r u_\theta &=0\qquad&&\text{on}\quad\partial_1\Omega\tag{profile evolution}\\
\phi(0,t)&=0\tag{contact angle}\\
\lim_{r\to\infty}\phi(r,t)&=0\tag{slope at infinity}\\
\phi(\ph,0)&=\phi_0\tag{initial profile},
\end{alignat*}
where $\Omega,\partial_0\Omega,\partial_1\Omega$ are as defined in \eqref{fixed wedge reference domain} and
\begin{align*}
\mb N&=-(\grad+\grad_\phi)\cdot(\grad+\grad_\phi)\mb u_\phi-(\grad\cdot\grad_\phi+\grad_\phi\cdot\grad+\grad_\phi\cdot\grad_\phi)\mb u+\grad_\phi p\\
\mb N_\partial&=(\sin\theta_b)^{-1}\brac{\brac{1+(b^{-1}-\phi)^2}^{-{3\over 2}}-(\sin\theta_b)^3}(\partial_r\phi)\mb e_\theta\\
&\quad+\brac{(\sin\theta_b)^{-1}\brac{1+(b^{-1}-\phi)^2}^{-{3\over 2}}\partial_r\phi I-(\grad+\grad_\phi)(\mb u+\mb u_\phi)-((\grad+\grad_\phi)(\mb u+\mb u_\phi))^T+pI}\mb n_\phi\\
&\quad-\brac{\grad_\phi(\mb u+\mb u_\phi)+\grad\mb u_\phi+(\grad_\phi(\mb u+\mb u_\phi))^T+(\grad\mb u_\phi)^T}\mb e_\theta\\
\mb n_\phi&=\brac{1+(b^{-1}-\phi)^2}^{-{1\over 2}}((b^{-1}-\phi)\mb e_2-\mb e_1)-(\sin\theta_b)(b^{-1}\mb e_2-\mb e_1)
\end{align*}
and $\grad_\phi$ and $\mb u_\phi$ are as in Lemma \ref{Piola transform}.
\end{system}

In System \ref{reference wedge system} the curvature of the free boundary is $-X''(1+X'^2)^{-3/2}=(\sin\theta_b)^{-1}(1+(b^{-1}-\phi)^2)^{-3/2}\partial_r\phi$, by $X'=b^{-1}-\phi$ and $\partial_r\phi=(\sin\theta_b)\phi'$. The value of this curvature at $\phi=0$ is the term $(\sin\theta_b)^2\partial_r\phi$ on the LHS of the surface tension equation, since $(1+b^{-2})^{-3/2}=\sin^3\theta_b$. The free boundary $\{x=X(y,t)\}$ has the outward unit normal $(1+X'^2)^{-1/2}(X'\mb e_2-\mb e_1)=\mb e_\theta+\mb n_\phi$, where $\mb e_\theta=(\sin\theta_b)(b^{-1}\mb e_2-\mb e_1)$ is the outward unit normal of $\partial_1\Omega$. Let $\kappa_\flat$ denote the curvature of the free boundary. Then $\mb N_\partial$ is obtained by subtracting the exact condition $(2\S(\mb u+\mb u_\phi)-pI)(\mb e_\theta+\mb n_\phi)=\kappa_\flat(\mb e_\theta+\mb n_\phi)$, with the symmetric gradient formed with $\grad+\grad_\phi$, from its linearisation. The no penetration and Navier slip conditions, on the other hand, transform without any nonlinear correction. On $\partial_0\Omega$ we have $\nu_2=\partial_{\tilde y}\E\chi=r^{-1}\partial_\theta\E\chi=0$ by the Neumann condition of the mixed harmonic extension. Geometrically, the flattening moves points only horizontally, and at the substrate by an amount that is stationary in $y$, which makes $D\bs\Psi_t$ diagonal there. So on $\partial_0\Omega$, by \eqref{transformation matrix}, $\tilde u_{\tilde y}=(1+\nu_1)u_y$ and $\tilde u_{\tilde x}=u_x$, and $\partial_{\tilde y}\tilde u_{\tilde x}=\partial_yu_x$ once $u_y=0$.

It is convenient to have a name for this system with its two forcings left free. We write $\mathfrak W_t[\mb F,\mb f;\phi_0]$ for the system above with $\mb N$ and $\mb N_\partial$ replaced by given forcings $\mb F$ on $\Omega$ and $\mb f$ on $\partial_1\Omega$ and with the initial profile $\phi_0$, so that System \ref{reference wedge system} is $\mathfrak W_t[\mb N,\mb N_\partial;\phi_0]$, and $\mathfrak W_s[\mb F,\mb f;\phi_0]$ for its two-sided Laplace transform in time, in which every unknown is replaced by its transform and the profile evolution equation becomes $s\hat\phi-(\sin\theta_b)^{-2}\partial_r\hat u_\theta=\phi_0$.

Extending the variables by $0$ to negative time adds $\phi_0\delta(t)$ to the RHS of the profile evolution equation, and taking the two-sided Laplace transform in time we get $\mathfrak W_s[\hat{\mb N},\hat{\mb N}_\partial;\phi_0]$, that is,

\begin{system}[Wedge-like Stokes flow on the fixed wedge, after the Laplace transform]\label{reference wedge system on the Laplace side}
\begin{alignat*}{3}
\lpc\hat{\mb u}-\grad\hat p&=\hat{\mb N}[\phi,\mb u,p]\qquad&&\text{on}\quad\Omega\tag{Stokes equation}\\
\grad\cdot\hat{\mb u}&=0\qquad&&\text{on}\quad\Omega\tag{incompressibility}\\
\hat u_\theta&=0\qquad&&\text{on}\quad\partial_0\Omega\tag{no penetration}\\
\hat u_r-r^{-1}\partial_\theta\hat u_r&=0\qquad&&\text{on}\quad\partial_0\Omega\tag{Navier slip}\\
\mathbb T[\hat{\mb u},\hat p]\mb e_\theta-(\sin\theta_b)^2(\partial_r\hat\phi)\mb e_\theta&=\hat{\mb N}_\partial[\phi,\mb u,p]\qquad&&\text{on}\quad\partial_1\Omega\tag{surface tension}\\
s\hat\phi-(\sin\theta_b)^{-2}\partial_r\hat u_\theta &=\phi_0\qquad&&\text{on}\quad\partial_1\Omega\tag{profile evolution}\\
\hat\phi(0;s)&=0\tag{contact angle}\\
\lim_{r\to\infty}\hat\phi(r;s)&=0\tag{slope at infinity}
\end{alignat*}
\end{system}

Solving the profile evolution equation for $\hat\phi$ and substituting the result into the surface tension equation, the unknown $\hat\phi$ is eliminated and we get
\begin{align}\label{surface tension without profile}
\mathbb T[\hat{\mb u},\hat p]\mb e_\theta-s^{-1}\partial_r^2\hat u_\theta\mb e_\theta =\hat{\mb N}_\partial[\phi,\mb u,p]+s^{-1}(\sin\theta_b)^2(\partial_r\phi_0)\mb e_\theta\qquad\text{on}\quad\partial_1\Omega.
\end{align}
Similarly, the contact angle condition and the slope at infinity become conditions on $\partial_r\hat u_\theta$ alone,
\begin{align}\label{contact angle without profile}
\begin{aligned}
\partial_r\hat u_\theta(0,\theta_b)&=(\sin\theta_b)^2(s\hat\phi(0;s)-\phi_0(0))=0,\\
\lim_{r\to\infty}\partial_r\hat u_\theta(r,\theta_b)&=(\sin\theta_b)^2\lim_{r\to\infty}(s\hat\phi(r;s)-\phi_0(r))=0,
\end{aligned}
\end{align}
by the profile evolution equation: the contact angle and slope conditions of System \ref{reference wedge system on the Laplace side} hold for every $s$, being the transforms of $\phi(0,t)=0$ and $\lim_{r\to\infty}\phi(r,t)=0$ at all times, and the initial profile satisfies the same two conditions. It is the system \eqref{surface tension without profile}--\eqref{contact angle without profile}, with $\hat{\mb N}$ and $\hat{\mb N}_\partial$ frozen as given forcing, that we linearise in the next section.

\section{Existence for the linearised system}\label{linear existence}

In this section we solve the linearisation of System \ref{reference wedge system on the Laplace side}, in which the nonlinearities are frozen as given forcing. The solution is variational: for each number $l\in\N_0$ of radial derivatives to be controlled, we construct a sesquilinear form $B_{(l)}$ and a linear functional $\mathcal F_{(l)}$ on a weighted space $H_s^l$ whose identity $B_{(l)}[\mb u,\ph]=\mathcal F_{(l)}$ is, for sufficiently regular test functions, the weak form of the system, and apply the Lax--Milgram theorem. Two features make this more delicate than a standard energy argument. First, the form comes from testing not against $\mb v$ but against test functions containing derivatives, as for the stationary Stokes problem in \cite{Bravin} and for the heat equation with Robin boundary conditions in \cite{BravinHeat}. These test functions are the two families $\smash[b]{\mb V^r_{a,j}}(\mb v)$ and $\smash[b]{\mb V^\theta_{a,j}}(\mb v)$ of Definition \ref{test functions}, which control all derivatives up to order two together with radial derivatives of any order in one go. Second, the norm of $H_s^l$ combines two weights: $a=0$, the plain energy, and $a=-\kappa$, which sits clear of the poles crossed in the corner analysis of Section \ref{corner expansions}. At the weight $a=-\kappa$ the surface tension term, which carries the factor $1/s$, leaves a remainder, and this is paid for out of the plain energy at the cost $|s|^{-2\kappa}$ (see the discussion after Definition \ref{H_s^l space}). Since these test functions are not arbitrary, recovering the equation pointwise from the variational identity requires the ``test function problem'' of Subsections \ref{the test function problem} and \ref{solvability of the test function problem}.

\subsection{The linear system}

To solve our wedge-like Stokes problem with a free surface, we first consider the following linearisation.

\begin{system}[Linear system]\label{linear system}
\begin{alignat*}{3}
\lpc\mb u-\grad p_\P &=\P\mb F\qquad&&\text{on}\qquad\Omega\\
\grad\cdot\mb u&=0\qquad&&\text{on}\qquad\Omega\\
u_\theta&=0\qquad&&\text{on}\qquad\partial_0\Omega\\
r^{-1}\partial_\theta u_r-u_r&=0\qquad&&\text{on}\qquad\partial_0\Omega\\
2\mb e_r\cdot\S(\mb u)\cdot\mb e_\theta&=f_r\qquad&&\text{on}\qquad\partial_1\Omega\\
\partial_ru_\theta(0,\theta_b)&=0\\
\lim_{r\to\infty}\partial_ru_\theta(r,\theta_b)&=0\\
p_\P &=(g_{1,0}+g_{1,1})[\mb u]+g_{1,\phi_0}[\phi_0]+g_{1,f_\theta}[f_\theta]\qquad&&\text{on}\qquad\partial_1\Omega\\
\partial_\theta p_\P &=g_0[\mb u]\qquad&&\text{on}\qquad\partial_0\Omega,
\end{alignat*}
where $\Omega$ is defined in \eqref{fixed wedge reference domain}; $g_\bullet$ are defined below in Definition \ref{g definition}; and $\mb F,\mb f$ and $\phi_0$ are given data.
\end{system}

\begin{definition}[$g_\bullet$ definition]\label{g definition}
We denote
\begin{align*}
g_0[\mb u](r;s)&:=-\partial_\theta(\partial_ru_r-r^{-1}u_r)(r,0;s)\\
g_{1,0}[\mb u](r;s)&:=2(\mb e_\theta\cdot\S(\mb u)\cdot\mb e_\theta)(r,\theta_b;s)\\
g_{1,1}[\mb u](r;s)&:=-{1\over s}\partial_r^2u_\theta(r,\theta_b;s)\\
g_{1,\phi_0}[\phi_0](r;s)&:=-{1\over s}(\sin\theta_b)^{2}\partial_r\phi_0(r)\\
g_{1,f_\theta}[f_\theta](r;s)&:=-f_\theta(r;s)
\end{align*}
\end{definition}

Let us explain where these equations come from. The first seven are the equations of System \ref{reference wedge system on the Laplace side}, with the profile already eliminated as in \eqref{surface tension without profile}--\eqref{contact angle without profile} and with the nonlinearities frozen into the given forcing $\mb F$ and $\mb f$: the fifth equation is the $\mb e_r$ component of \eqref{surface tension without profile}, in which the pressure does not appear since $\mb e_r\cdot(pI)\cdot\mb e_\theta=0$, and the sixth and seventh are exactly \eqref{contact angle without profile}. Note that we have written the bulk forcing as $\P\mb F$ rather than $\mb F$; by Remark \ref{P remark} this does not change the velocity we are solving for, since $\mb F$ and $\P\mb F$ pair identically against every divergence-free field with vanishing $\theta$ component on $\partial_0\Omega$, and it is convenient to have a forcing that is already compatible with these two constraints.

The last two equations are the ones that are not part of the original system: they are the conditions that determine the pressure from the velocity. The first is the remaining scalar component of the surface tension condition, and the second is a consequence of the Stokes equation on $\partial_0\Omega$. The first of them is the $\mb e_\theta$ component of \eqref{surface tension without profile}, rearranged for $p_\P$: since $\mb e_\theta\cdot\mathbb T[\mb u,p_\P]\cdot\mb e_\theta=2\mb e_\theta\cdot\S(\mb u)\cdot\mb e_\theta-p_\P$, that component reads
\begin{align*}
p_\P=2(\mb e_\theta\cdot\S(\mb u)\cdot\mb e_\theta)-{1\over s}\partial_r^2u_\theta-{1\over s}(\sin\theta_b)^{2}\partial_r\phi_0-f_\theta\qquad\text{on}\qquad\partial_1\Omega,
\end{align*}
which is precisely the sum of the four terms $g_{1,0}$, $g_{1,1}$, $g_{1,\phi_0}$ and $g_{1,f_\theta}$ of Definition \ref{g definition}. The second is obtained by taking the $\mb e_\theta$ component of the Stokes equation on $\partial_0\Omega$: since $(\P\mb F)_\theta|_{\partial_0\Omega}=0$ by Remark \ref{P remark} and $\partial_\theta p_\P=r(\mb e_\theta\cdot\grad p_\P)$, we get $\partial_\theta p_\P=r(\mb e_\theta\cdot\lpc\mb u)=\dr(\grad^\perp\cdot\mb u)$ on $\partial_0\Omega$, since $\lpc\mb u=\grad^\perp(\grad^\perp\cdot\mb u)$ for divergence-free $\mb u$ and $\mb e_\theta\cdot\grad^\perp=\partial_r$, and along $\partial_0\Omega$ the no penetration condition $u_\theta=0$ leaves $\grad^\perp\cdot\mb u=-r^{-1}\partial_\theta u_r$, so that
\begin{align*}
\partial_\theta p_\P=-\dr\partial_\theta(r^{-1}u_r)=-\partial_\theta(\partial_ru_r-r^{-1}u_r)=g_0[\mb u]\qquad\text{on}\qquad\partial_0\Omega.
\end{align*}
Thus both the profile $\phi$ and the pressure $p_\P$ are expressed through $\mb u$ and the data, and $\mb u$ is left as the only unknown.

Taking the divergence of the first equation of System \ref{linear system}, we get $\lpc p_\P=0$. Together with the two boundary conditions for $p_\P$, this determines $p_\P$ and lets us write it in terms of $\mb u$ and the data $\phi_0$, $f_\theta$.

\begin{definition}[Pressure functionals]\label{pressure functionals}
Let $p_{\partial_0}$ and $p_{\partial_1,\star}$ be the harmonic functions such that $\partial_\theta p_{\partial_0}|_{\partial_0\Omega}=g_0$, $p_{\partial_0}|_{\partial_1\Omega}=0$ and $\partial_\theta p_{\partial_1,\star}|_{\partial_0\Omega}=0$, $p_{\partial_1,\star}|_{\partial_1\Omega}=g_{1,\star}$ given by Lemma \ref{mixed wedge problem} (i.e. $p_{\partial_0}=q_{f_0}$ with $f_0=g_0$, and $p_{\partial_1,\star}=\E g_{1,\star}$). Let
\begin{align*}
p_{\P|\mb u}[\mb u]&=(p_{\partial_0}+p_{\partial_1,0}+p_{\partial_1,1})[\mb u]\\
p_\P[\mb u,\phi_0,f_\theta]&=p_{\P|\mb u}[\mb u]+p_{\partial_1,\phi_0}[\phi_0]+p_{\partial_1,f_\theta}[f_\theta]
\end{align*}
\end{definition}

We will be seeking a variational solution to our system in the following weighted Sobolev space.

\begin{definition}[The space $H^l_s$]\label{H_s^l space}
Let $l\in\N_0$, $s\in\C\setminus\{0\}$, $\kappa\in(1/4,1/2]$ and $\epsilon>0$. We define
\begin{align*}
H_s^l:=\{\mb v:\|\mb v\|_{H_s^l}<\infty,\grad\cdot\mb v=0,v_\theta|_{\partial_0\Omega}=0\},
\end{align*}
where
\begin{align*}
\|\mb v\|_{H_s^l}^2 &:=\brac{1+|s|^{-2\kappa}}\|\mb v\|_{H_s^l,0}^2+\|\mb v\|_{H_s^l,-\kappa}^2,\\
\|\mb v\|_{H_s^l,a;\epsilon}^2 &:=\|\grad\mb v\|_{\Omega,a,l+1;\epsilon}^2 +\epsilon^2\|r\grad^2\mb v\|_{\Omega,a,l;\epsilon}^2 +\|v_r\|_{\partial_0\Omega,a,l+1;\epsilon}^2 +{1\over|s|}\|\partial_rv_\theta\|_{\partial_1\Omega,a,l+1;\epsilon}^2.
\end{align*}
\end{definition}

The norm $\|\mb v\|_{H_s^l,a}=\|\mb v\|_{H_s^l,a;1}$, the index $\epsilon$ being suppressed when it equals $1$ as in Subsection \ref{norms and Mellin transform}, is the energy that we get from testing our equation with a test function defined later in Definition \ref{test functions}, at the weight $a$. The two weights are $a=0$, the plain energy, and $a=-\kappa$, which sits clear of the poles at $\lambda\in\{0,1\}$ crossed in Section \ref{corner expansions}. The factor $1+|s|^{-2\kappa}$ in front of the plain part is what pays for the surface tension remainder of Lemma \ref{surface tension weighted remainder} at the weighted one, and is chosen so that both are available with a constant independent of $s$. The parameter $l$ counts how many radial derivatives we control.

The range $\kappa\in(1/4,1/2]$ is the one in which the coercivity of Proposition \ref{energy form coercivity} holds. The lower bound is used in the second estimate of Lemma \ref{pressure remainder} and in Lemma \ref{surface tension weighted remainder}, the upper bound in the interpolation steps of both and in Lemma \ref{weighted based control}. Both bounds reach that proposition only through these lemmas.

The norm controls $\grad\mb v$ but not $\mb v$ itself. Nevertheless every $\mb v\in H_s^l$ vanishes at infinity, by Lemma \ref{limit at infinity}, so that Hardy's inequality (Proposition \ref{Hardy's inequality}) is available at every weight $1-\varpi$ with $\varpi\in(0,\kappa]$; we use this below without further comment.

The following are two types of test functions that we will use to obtain our solution; testing with $\mb V^{r}_{a,j}$ will give us control of $j$ radial derivatives and one full derivative, while testing with $\mb V^{\theta}_{a,j}$ will give us control for $j$ radial derivatives on top of the Laplacian, which together with $\mb V^{r}_{a,j}$ will give us two full derivatives.

\begin{definition}[Test functions]\label{test functions}
We define the notation
\begin{align*}
\mb V^{r}_{a,j}(\mb v)&:=-r^{-2a}(-1)^j\dr^{2j}\mb v,\\
\mb V^{\theta}_{a,j}(\mb v)&:=(-1)^j\dr^{2j}\grad^\perp(r^{2-2a}\grad^\perp\cdot\mb v).
\end{align*}
\end{definition}

\begin{definition}[The counts $j^\pm$]\label{j plus minus}
For $j\in\N_0$, we write
\begin{align*}
j^+&:=j+\mb 1[j>0]\\
j^-&:=j-\mb 1[j>0],
\end{align*}
where $\mb 1[P]$ equals $1$ if the statement $P$ holds and $0$ otherwise.
\end{definition}

Testing $\mb V^{r}_{a,j}$ against $\mb V^{\theta}_{a,l}$ with $l>0$ and balancing the derivatives by integration by parts gives $j+l^+$ on each side, which is what the notation above is for.

For the test function problem that isolates the bulk, we will impose boundary conditions $\mb v|_{\partial_0\Omega}=\mb 0$ and $v_r|_{\partial_1\Omega}=0$, so we introduce the following space of functions. These two conditions are not those of System \ref{linear system}: the class serves only to recover the equation in the bulk, and its conditions are the ones under which the boundary terms of the identity of Proposition \ref{energy forms} vanish, as in the proof of Theorem \ref{linear solution Laplace side}, and under which all but one of the boundary terms vanish in the integration by parts of Lemma \ref{Vr against Vtheta}. The Navier slip condition is, like the $\mb e_r$ component of the surface tension condition, natural for the form, and it is recovered at the end of the proof of Theorem \ref{linear solution Laplace side}.

\begin{definition}[The class $C^\infty_{c,\sigma}(\Omega)$]\label{smooth bulk test functions}
We define
\begin{align*}
C^\infty_{c,\sigma}(\Omega)=\{\mb v\in C^\infty_c(\R^2)|_{\overline\Omega}:\grad\cdot\mb v=0,\mb v|_{\partial_0\Omega}=\mb 0,v_r|_{\partial_1\Omega}=0\},
\end{align*}
the subscript $\sigma$ standing for solenoidal.
\end{definition}

\subsection{Basic lemmas}

We collect a few useful lemmas here.

\begin{lemma}[The projections commute with $\dr$]\label{projections commute with dr}
We have $\P\dr=\dr\P$ and $q[\dr\ph]=\dr_{-1}q[\ph]$, where $\P$ and $q$ are as defined in Definition \ref{the two projections}. The same statements hold for $\P^\theta$ and $q^\theta$.
\end{lemma}
\begin{proof}
We have
\begin{align*}
\lpc q[\dr\mb v]&=\grad\cdot\dr\mb v=\dr_1\grad\cdot\mb v=\dr_1\lpc q[\mb v]=\lpc\dr_{-1}q[\mb v]\\
(\partial_\theta q[\dr\mb v])|_{\partial_0\Omega}&=(r\dr v_\theta)|_{\partial_0\Omega}=(\dr_{-1}rv_\theta)|_{\partial_0\Omega}=(\dr_{-1}\partial_\theta q[\mb v])|_{\partial_0\Omega}=(\partial_\theta\dr_{-1} q[\mb v])|_{\partial_0\Omega}\\
(\dr_{-1}q[\mb v])|_{\partial_1\Omega}&=0=q[\dr\mb v]|_{\partial_1\Omega}
\end{align*}
This means $q[\dr\mb v]=\dr_{-1}q[\mb v]$, so
\begin{align*}
\P\dr\mb v=\dr\mb v-\grad q[\dr\mb v]=\dr\mb v-\grad\dr_{-1}q[\mb v]=\dr\mb v-\dr\grad q[\mb v]=\dr\P\mb v.
\end{align*}
The proof for $\P^\theta$ is the same, with the Neumann condition on $\partial_1\Omega$ in place of the Dirichlet one.
\end{proof}

\begin{lemma}[Smoothness and decay of functions in $\P C_c^\infty(\Omega)$]\label{decay of projected test functions}
Suppose $\mb v\in C_c^\infty(\Omega)$ and let $a'\geq 0$ with $a'\theta_b<\pi/2$. Then $\|(r\grad)^{j}\P\mb v\|_{\Omega,a}<\infty$ for all $|a|\leq a'$ and $j\in\N_0$. ($\P$ was defined in Definition \ref{the two projections})
\end{lemma}
\begin{proof}
By definition $\P\mb v=\mb v-\grad q$, where $q$ is defined by $\lpc q=\grad\cdot\mb v$ with $q|_{\partial_1\Omega}=0$ and $\partial_\theta q|_{\partial_0\Omega}=rv_\theta|_{\partial_0\Omega}=0$. Since $\mb v\in C_c^\infty(\Omega)$, by Lemma \ref{mixed wedge problem} we have
\begin{align*}
\|(r\grad)^j\grad q\|_{\Omega,a}\lesssim_j\sum_{k=1}^{j+1}\|(r\grad)^kq\|_{\Omega,a+1}\lesssim_j{1\over|\cos(a\theta_b)|}\|r\mb v\|_{\Omega,a+1,0,j}<\infty,
\end{align*}
the lemma being applied at the weight $a+1$ with $\mb F=r\mb v$, which is admissible since $|a|\theta_b\leq a'\theta_b<\pi/2$. Hence $\|(r\grad)^j\P\mb v\|_{\Omega,a}\leq\|(r\grad)^j\mb v\|_{\Omega,a}+\|(r\grad)^j\grad q\|_{\Omega,a}<\infty$.
\end{proof}

\begin{lemma}[Trace]\label{trace}
Let $m\in\{0,1\}$ and $a\in\R$. For any $u\in H^1_{\rm loc}(\overline\Omega\setminus\{\mb 0\})$ for which the RHSs are finite, we have
\begin{align*}
\|u\|_{\partial_m\Omega,{1\over 2}+a}^2&\leq\theta_b^{-1}\|u\|_{\Omega,1+a}^2+2\|u\|_{\Omega,1+a}\|\grad u\|_{\Omega,a}\\
\|\dr^{1\over 2}u\|_{\partial_m\Omega,{1\over 2}+a}^2&\leq\|\dr u\|_{\Omega,1+a}\brac{\theta_b^{-1}\|u\|_{\Omega,1+a}+2\|\grad u\|_{\Omega,a}},
\end{align*}
where the meaning of the fractional power is laid out in Definition \ref{fractional counts}.
\end{lemma}
\begin{proof}
We will prove the case of the lower boundary $m=0$. The proof for the upper boundary $m=1$ is similar. We have
\begin{align*}
\|u\|_{\partial_0\Omega,{1\over 2}+a}^2&=\|r^{-{1\over 2}-a}u\|_{L^2(\partial_0\Omega)}^2=\int_{\partial_0\Omega}|u(r,0)|^2r^{-1-2a}\d S\\
&=\theta_b^{-1}\int_0^\infty\int_0^{\theta_b}\partial_\theta((\theta-\theta_b)|u(r,\theta)|^2)r^{-1-2a}\d \theta\d r\\
&=\theta_b^{-1}\int_0^\infty\int_0^{\theta_b}\brac{|u|^2+2(\theta-\theta_b)\Re(\bar u\partial_\theta u)}r^{-1-2a}\d \theta\d r\\
&\leq\theta_b^{-1}\|r^{-1-a}u\|_{L^2(\Omega)}^2+2\|r^{-1-a}u\|_{L^2(\Omega)}\|r^{-a}\grad u\|_{L^2(\Omega)}
\end{align*}
using $|\theta-\theta_b|\leq\theta_b$ and $|\partial_\theta u|\leq r|\grad u|$, and
\begin{align*}
\|r^{-{1\over 2}-a}\dr^{1\over 2}u\|_{L^2(\partial_0\Omega)}^2&=\int_{\partial_0\Omega}|\dr^{1\over 2}u(r,0)|^2r^{-1-2a}\d S
=\int_{\Re\lambda=a}|\lambda||\hat u(\lambda,0)|^2\d\Im\lambda\\
&=\theta_b^{-1}\int_{\Re\lambda=a}\int_0^{\theta_b}\partial_\theta((\theta-\theta_b)|\lambda||\hat u(\lambda,\theta)|^2)\d \theta\d\Im\lambda\\
&\leq\|r^{-1-a}\dr u\|_{L^2(\Omega)}(\theta_b^{-1}\|r^{-1-a}u\|_{L^2(\Omega)}+2\|r^{-a}\grad u\|_{L^2(\Omega)}).\qedhere
\end{align*}
\end{proof}

\begin{lemma}[Second derivatives from radial derivatives and the Laplacian]\label{second derivative bound}
For any $\mb u\in H^2_{\rm loc}(\Omega)$ we have, almost everywhere on $\Omega$,
\begin{align*}
|r\grad^2\mb u|^2&\lesssim|\grad(r\partial_r\mb u)|^2+|\grad\mb u|^2+|r\lpc\mb u|^2.
\end{align*}
\end{lemma}
\begin{proof}
We have
\begin{align*}
|r\grad^2\mb u|^2&\lesssim|r\partial_r\grad\mb u|^2+|\partial_\theta\grad\mb u|^2
\lesssim|r\partial_r\grad\mb u|^2+|\partial_\theta(\mb e_r\partial_r+\mb e_\theta r^{-1}\partial_\theta)\mb u|^2\\
&\lesssim|r\partial_r\grad\mb u|^2+|\grad\mb u|^2
+|\partial_r\partial_\theta\mb u|^2+|\partial_r(r\partial_r)\mb u|^2+|r\lpc\mb u|^2\\
&\lesssim|\grad(r\partial_r\mb u)|^2+|\grad\mb u|^2+|r\lpc\mb u|^2,
\end{align*}
using $\grad=\mb e_r\partial_r+\mb e_\theta r^{-1}\partial_\theta$, the polar form $r\lpc=\partial_r(r\partial_r)+r^{-1}\partial_\theta^2$ of the Laplacian for the term $r^{-1}\partial_\theta^2\mb u$, and, in the last step, $\grad(r\partial_r\mb u)=r\partial_r\grad\mb u+\grad\mb u$, of which $\partial_r\partial_\theta\mb u=r^{-1}\partial_\theta(r\partial_r\mb u)$ and $\partial_r(r\partial_r\mb u)$ are components.
\end{proof}

By Lemma \ref{second derivative bound}, in the construction of our solution it will be sufficient to test against $\mb V^r$ and $\mb V^\theta$ as defined in Definition~\ref{test functions}, which control $\grad(r\partial_r\mb u)$, $\grad\mb u$ and $r\lpc\mb u$.

The next lemma is Korn's inequality in the form in which we shall need it: no condition is imposed on $\partial_1\Omega$, so this is not the first Korn inequality, and the constant must not degenerate as $\theta_b\to 0$, since $\theta_b$ is taken small in Proposition \ref{energy form coercivity}.

\begin{lemma}[Korn's inequality on the wedge]\label{Korn on the wedge}
Let $\theta_b\in(0,\pi/2)$ and $a\in\R$. Suppose $\grad\cdot\mb u=0$ on $\Omega$, $u_\theta|_{\partial_0\Omega}=0$ and $\|\grad\mb u\|_{\Omega,a}<\infty$. Then
\begin{align}\label{Korn inequality}
\|\grad\mb u\|_{\Omega,a}^2\lesssim_{|a|}\|\S(\mb u)\|_{\Omega,a}^2.
\end{align}
\end{lemma}
\begin{proof}
In two dimensions $|\grad\mb u|^2=|\S(\mb u)|^2+{1\over 2}|\grad^\perp\cdot\mb u|^2$ pointwise, so \eqref{Korn inequality} says exactly that the vorticity is controlled by the strain. Some hypothesis is needed for that, a rigid rotation having vorticity but no strain, and here it is the boundary condition that supplies one, no rotation of the plane being tangent to $\partial_0\Omega$. Since $\Omega$ is a cone and the two sides scale alike we argue on the Mellin side, transforming $\grad\mb u$ and not $\mb u$, so that no weighted integrability of $\mb u$ --- which is not required to vanish at the corner --- is needed.

By the polar-coordinate formula \eqref{grad polar} the entries of $r\grad\mb u$ are $\dr u_r$, $\dr u_\theta$, $\partial_\theta u_r-u_\theta$ and $\partial_\theta u_\theta+u_r$, and $\|r\grad\mb u\|_{\Omega,a+1}=\|\grad\mb u\|_{\Omega,a}$ is finite by hypothesis, so that by Remark \ref{weight and Mellin line} these four functions may be Mellin transformed in $r$ on the line $\Re\lambda=a$. Write, suppressing the arguments $(\lambda,\ph)$,
\begin{align*}
g:=(\dr u_r)^\wedge,\qquad h:=(\dr u_\theta)^\wedge,\qquad \tilde g:=(\partial_\theta u_r-u_\theta)^\wedge,\qquad \tilde h:=(\partial_\theta u_\theta+u_r)^\wedge,
\end{align*}
so that
\begin{align*}
(r\grad\mb u)^\wedge=\mat{g&h\\\tilde g&\tilde h}.
\end{align*}
Incompressibility is the vanishing of the trace of $\grad\mb u$, so becomes $\tilde h=-g$; and applying $\dr$ to the second row and using that $\dr$ commutes with $\partial_\theta$ gives
\begin{align}\label{Korn row relation}
\lambda\tilde g=g'-h,\qquad\qquad\lambda\tilde h=h'+g,
\end{align}
and hence $h'=-\lambda_1g$, where $\lambda_1=\lambda+1$ as in Definition \ref{scaled radial derivatives}. Here and below the prime denotes $\partial_\theta$, the only derivative left after the transform, $\dr$ having become multiplication by $\lambda$. The identities \eqref{Korn row relation} hold in the sense of distributions in $\theta$ for almost every $\lambda$ on the line, and since $\tilde g,\tilde h,g,h$ are square integrable in $\theta$ they show that $g$ and $h$ lie in $H^1(0,\theta_b)$; this is all the regularity used below, and it gives $g$ and $h$ traces at both ends. That of $h$ at $\theta=0$ vanishes: $\dr u_\theta$ is a tangential derivative of $u_\theta$, which vanishes on $\partial_0\Omega$, so $\dr u_\theta(\ph,\theta)\to 0$ as $\theta\to 0$ in the sense of distributions in $r$, and hence $h(\ph,0)=0$. Finally, by incompressibility ($\tilde h=-g$),
\begin{align*}
|(r\grad\mb u)^\wedge|^2=2|g|^2+|h|^2+|\tilde g|^2,\qquad\qquad
|(r\S(\mb u))^\wedge|^2=2|g|^2+{1\over 2}|\psi|^2,\qquad\qquad \psi:=\tilde g+h,
\end{align*}
so, since $\|r\S(\mb u)\|_{\Omega,a+1}=\|\S(\mb u)\|_{\Omega,a}$ as well, Plancherel's identity \eqref{Plancherel's identity} reduces \eqref{Korn inequality} to the integral over $(0,\theta_b)$ of the first of these quantities being $\lesssim_{|a|}$ that of the second, for almost every $\lambda$ on the line. The strain therefore sees exactly $g$ and $\psi$, and since $\tilde g=\psi-h$ the one thing it does not see is $h$ --- the vorticity again, $(r\grad^\perp\cdot\mb u)^\wedge=h-\tilde g=2h-\psi$. So the lemma reduces to
\begin{align}\label{Korn key}
\int_0^{\theta_b}|h|^2\d\theta\lesssim_{|a|}\int_0^{\theta_b}\brac{|g|^2+|\psi|^2}\d\theta.
\end{align}

In terms of $v:=(g,h)$ the relations \eqref{Korn row relation}, together with $\tilde h=-g$, in the form $g'=h+\lambda\tilde g=-\lambda_{-1}h+\lambda\psi$ and $h'=-\lambda_1g$, are the first-order system
\begin{align}\label{Korn system}
v'=Mv+(\lambda\psi,0),\qquad\qquad M=\mat{0&-\lambda_{-1}\\-\lambda_1&0},\qquad\qquad h(0)=0.
\end{align}
Since $\lambda_{-1}=-\bar\lambda_1+2\Re\lambda$, the matrix $M$ is anti-Hermitian but for the entry $-2a$ in its upper right corner, and with that $-2a$ removed it has square $-|\lambda_1|^2I$: up to a perturbation of norm $2|a|$, bounded independently of $\lambda$ and of $\theta_b$, the system is a rotation of $\C^2$ of period $2\pi/|\lambda_1|$, conserving $|v|^2$ and carrying each component into the other. Fix $K\geq 1$, to be chosen below in terms of $|a|$ alone; the proof divides according to whether the wedge falls short of $K$ periods or exceeds them.

Suppose first that $|\lambda_1|\theta_b\leq 2\pi K$. The system then has too little room to turn, and $h$, which starts at zero and grows at the rate $h'=-\lambda_1g$, is controlled by $g$ directly: $h(0)=0$ and the Cauchy--Schwarz inequality give $|h(\theta)|^2\leq\theta\int_0^{\theta_b}|h'(\theta')|^2\d\theta'$, so that
\begin{align*}
\int_0^{\theta_b}|h|^2\d\theta \leq{\theta_b^2\over 2}|\lambda_1|^2\int_0^{\theta_b}|g|^2\d\theta \leq 2\pi^2K^2\int_0^{\theta_b}|g|^2\d\theta.
\end{align*}

Suppose now that $|\lambda_1|\theta_b>2\pi K$, so that, $\theta_b$ being less than $\pi/2$, also $|\lambda_1|>4K\geq 4$ and hence $|\lambda|\leq|\lambda_1|+1\leq 2|\lambda_1|$. The rotation is now completed many times, and it is the oscillation it imposes that replaces the boundary condition. Set $S:=\Re(\bar\lambda_1\bar gh)$, which by $h'=-\lambda_1g$ is exactly $-{1\over 2}(|h|^2)'$; differentiating it and using \eqref{Korn system}, with $-\bar\lambda_{-1}=\lambda_1-2a$ and $\Re\bar\lambda_1=a+1$,
\begin{align}\label{Korn exchange}
S'=|\lambda_1|^2\brac{|h|^2-|g|^2}-2a(a+1)|h|^2+\Re(\bar\lambda_1\bar\lambda\bar\psi h).
\end{align}
Dropping the last two terms leaves $(|h|^2)''=2|\lambda_1|^2(|g|^2-|h|^2)$: the imbalance between the two components is restored rather than accumulated, so that integrating \eqref{Korn exchange} across the wedge leaves only the endpoint values of $S$, of which $S(0)$ vanishes with $h(0)$ and $S(\theta_b)$ is a single point value, negligible against an integral over a wedge many periods wide. Precisely, \eqref{Korn exchange} integrates to
\begin{align}\label{Korn integrated}
\int_0^{\theta_b}|h|^2\d\theta=\int_0^{\theta_b}|g|^2\d\theta+{S(\theta_b)\over|\lambda_1|^2}+{2a(a+1)\over|\lambda_1|^2}\int_0^{\theta_b}|h|^2\d\theta-{1\over|\lambda_1|^2}\int_0^{\theta_b}\Re(\bar\lambda_1\bar\lambda\bar\psi h)\d\theta.
\end{align}
For the endpoint term, $|S|\leq|\lambda_1||g||h|\leq{1\over 2}|\lambda_1||v|^2$, and $M$ with its $-2a$ removed contributing nothing to the derivative of $|v|^2$,
\begin{align*}
(|v|^2)'=-4a\Re(\bar gh)+2\Re(\lambda\bar g\psi),\qquad\text{so that}\qquad \abs{(|v|^2)'}\leq 2|a||v|^2+2|\lambda||\psi||v|.
\end{align*}
Integrating this from $\theta$ to $\theta_b$ and averaging in $\theta$ bounds $|v(\theta_b)|^2$, and hence
\begin{align*}
{|S(\theta_b)|\over|\lambda_1|^2}\leq{|v(\theta_b)|^2\over 2|\lambda_1|} \leq\brac{{1\over 2|\lambda_1|\theta_b}+{|a|\over|\lambda_1|}}\int_0^{\theta_b}|v|^2\d\theta+{|\lambda|\over|\lambda_1|}\int_0^{\theta_b}|\psi||v|\d\theta.
\end{align*}
Now $|\lambda_1|\theta_b>2\pi K$ and $|\lambda_1|>4K$ make each of $(|\lambda_1|\theta_b)^{-1}$, $|a||\lambda_1|^{-1}$ and $|a(a+1)||\lambda_1|^{-2}$ at most $(1+|a|)^2/K$, and we fix $K$ in terms of $|a|$ so that this is at most $1/16$. With $|\lambda|\leq 2|\lambda_1|$ and $2xy\leq x^2/8+8y^2$ for the two terms carrying $\psi$, the three terms of \eqref{Korn integrated} after the first then come to at most half of $\int_0^{\theta_b}(|g|^2+|h|^2)\d\theta$ together with a multiple of $\int_0^{\theta_b}|\psi|^2\d\theta$; absorbing $\int_0^{\theta_b}|h|^2\d\theta$ gives \eqref{Korn key}. In either case the constant depends only on $K$, hence only on an upper bound for $|a|$, and not on $\theta_b$.
\end{proof}

\subsection{Pressure estimates}

We have the following estimate for our pressure functionals.

\begin{lemma}[Boundedness of $\P$]\label{projection norm estimate}
Let $|a|\theta_b<\pi/2$ and $j\in\N_0$. We have
\begin{align*}
  \|\P\mb F\|_{\Omega,a,j} &\lesssim{1\over|\cos(a\theta_b)|} \|\mb F\|_{\Omega,a,j}.
\end{align*}
\end{lemma}
\begin{proof}
By Lemma \ref{projections commute with dr}, $\dr^{j'}\P\mb F=\P\dr^{j'}\mb F$. So the estimate for $j=0$, applied to $\dr^{j'}\mb F$ for each $j'\leq j$, gives it for every $j$ with the same constant, and it suffices to take $j=0$. With $\mb G:=r\mb F$, the function $q[\mb F]$ of Definition \ref{the two projections} is the solution $q_{\mb G}$ of Lemma \ref{mixed wedge problem} with $f_0=f_1=0$: indeed $\grad\cdot(r^{-1}\mb G)=\grad\cdot\mb F$, and the Neumann datum $G_\theta|_{\partial_0\Omega}=rF_\theta|_{\partial_0\Omega}$ carried by $q_{\mb G}$ is exactly the one of Definition \ref{the two projections}, so that no boundary datum is left over. Since $\|\grad q\|_{\Omega,a}=\|r\grad q\|_{\Omega,a+1}$ and $r\grad q=\mb e_r\dr q+\mb e_\theta\partial_\theta q$, we apply the first estimate of that lemma at the weight $a+1$, whose admissibility condition is $|a|\theta_b<\pi/2$, once with $l=1$ and once with $l=0$, $1$ in place of $j$ and the choice $\dr^{-1}$:
  \begin{align*}
    \|\partial_\theta q[\mb F]\|_{\Omega,a+1}+\|\dr q[\mb F]\|_{\Omega,a+1} \lesssim{1\over|\cos(a\theta_b)|} \|\mb G\|_{\Omega,a+1}={1\over|\cos(a\theta_b)|} \|\mb F\|_{\Omega,a}.
  \end{align*}
 Hence $\|\grad q[\mb F]\|_{\Omega,a}\lesssim|\cos(a\theta_b)|^{-1}\|\mb F\|_{\Omega,a}$, and $\P\mb F=\mb F-\grad q[\mb F]$ finishes the proof.
\end{proof}

\begin{lemma}[Pressure estimates]\label{pressure estimates}
Let $K\geq 1$ and $a-1\in(-K,K)$. Suppose $\theta_b\in(0,\pi/3K)$. Then for $p_\bullet$ as defined in Definition \ref{pressure functionals}, for $\E$ as defined in Lemma \ref{mixed harmonic extension}, and, in (i), for $\mb u$ with $u_\theta|_{\partial_0\Omega}=0$, we have
\begin{alignat*}{2}
(i)\quad\qquad && \|(r\grad)^j\{p_{\partial_0},p_{\partial_1,0}\}[\mb u]\|_{\Omega,a}^2&\lesssim_j\|\grad\mb u\|_{\Omega,a,j}\brac{\theta_b^{-1}\|\grad\mb u\|_{\Omega,a,j^-}+\|r\grad^2\mb u\|_{\Omega,a,j^-}}\\
(ii)\quad\qquad && \|(r\grad)^jp_{\partial_1,1}[\mb u]\|_{\Omega,a}^2&\lesssim_j{1\over|s|^2}\|\dr u_\theta\|_{\partial_1\Omega,a+{3\over 2},1+j}^2\\
(iii)\quad\qquad && \|(r\grad)^jp_{\partial_1,1}[\mb u]\|_{\partial_0\Omega,-{1\over 2}+a}^2&\lesssim_j{\theta_b^{-2j}\over|s|^2}\|\dr u_\theta\|_{\partial_1\Omega,a+{3\over 2},1}^2\\
(iv)\quad\qquad && \|(r\grad)^j\E\phi\|_{\Omega,a}^2&\lesssim_j\|\phi\|_{\partial_1\Omega,a-{1\over 2},j-{1\over 2}}^2,
\end{alignat*}
where $\{\ph\}$ is a choice as in Definition \ref{notation} and $\phi$ in (iv) is any function on $\partial_1\Omega$ for which the RHS is finite.
\end{lemma}
\begin{proof}
For (i), since by Definition \ref{g definition} it holds that $g_0(r;s)=-\partial_\theta(\partial_r u_r-r^{-1}u_r)(r,0;s)=-\partial_\theta\dr(r^{-1}u_r)(r,0;s)$, we have $\hat g_0(\lambda;s)=-\lambda\partial_\theta\hat u_r(\lambda+1,0;s)$. Hence by Lemma \ref{mixed wedge problem}, with the choice $\theta_b^{1/2}$ when $j'=0$ and $\dr^{-1/2}$ when $j'\geq 1$, so that no negative power of $|\lambda-1|$ occurs on the line $\Re\lambda=a$, we have
\begin{align*}
\|(r\grad)^j p_{\partial_0}\|_{\Omega,a}^2 &\lesssim_j{1\over|\cos((1-a)\theta_b)|^2}\sum_{j'=\min\{1,j\}}^j\norm{\set{\theta_b^{1\over 2},\dr^{-{1\over 2}}}\dr^{j'-1}g_0}_{\partial_0\Omega,a-{1\over 2}}^2\\
&={1\over|\cos((1-a)\theta_b)|^2}\sum_{j'=\min\{1,j\}}^j\int_{\Re\lambda=-1+a}\set{\theta_b,|\lambda|^{-1}}|\lambda|^{2j'-2} |\lambda\partial_\theta\hat u_r(\lambda+1,0;s)|^2\d\Im\lambda\\
&={1\over|\cos((1-a)\theta_b)|^2}\sum_{j'=\min\{1,j\}}^j\int_{\Re\lambda=a}\set{\theta_b,|\lambda-1|^{-1}}|\lambda-1|^{2j'}|\partial_\theta\hat u_r(\lambda,0;s)|^2\d\Im\lambda\\
&\lesssim\int_{\Re\lambda=a}\brac{1+|\lambda|^{\max\{0,2j-1\}}} |\partial_\theta\hat u_r(\lambda,0;s)|^2\d\Im\lambda\\
&\lesssim\|w\|_{\partial_0\Omega,{1\over 2}+a}^2 +\|\dr^{\max\{0,j-{1\over 2}\}}w\|_{\partial_0\Omega,{1\over 2}+a}^2\\
&\lesssim\|\grad\mb u\|_{\Omega,a}(\theta_b^{-1}\|\grad\mb u\|_{\Omega,a}+\|r\grad^2\mb u\|_{\Omega,a})
+\|\grad\mb u\|_{\Omega,a,j}(\theta_b^{-1}\|\grad\mb u\|_{\Omega,a,j^-}+\|r\grad^2\mb u\|_{\Omega,a,j^-})\\
&\lesssim\|\grad\mb u\|_{\Omega,a,j}\brac{\theta_b^{-1}\|\grad\mb u\|_{\Omega,a,j^-}+\|r\grad^2\mb u\|_{\Omega,a,j^-}},
\end{align*}
where $w:=\mb e_r\cdot\partial_\theta\mb u$, whose trace on $\partial_0\Omega$ is that of $\partial_\theta u_r$: in the bulk $\partial_\theta u_r=\mb e_r\cdot\partial_\theta\mb u+u_\theta$, and $u_\theta$ vanishes on $\partial_0\Omega$ by hypothesis. This is where the hypothesis $u_\theta|_{\partial_0\Omega}=0$ enters, and it cannot be dropped: for $\mb u\equiv\mb e_2$ we have $g_0=r^{-1}\neq 0$ while $\grad\mb u=\mb 0$. We used the trace Lemma \ref{trace} for the second to last inequality, applied to $w$ and to $\dr^{j-1}w$, together with $|w|\leq r|\grad\mb u|$ and $|\grad w|\lesssim|\grad\mb u|+r|\grad^2\mb u|$; the constant depends on $K$ through $|\cos((1-a)\theta_b)|\geq 1/2$ and on $j$, and on nothing else. The proofs for the other estimates are a similar application of Lemma \ref{mixed wedge problem}. The bound for $p_{\partial_1,0}$ in (i) is the bulk estimate for $q_{f_1}$ of that lemma with $f_1=g_{1,0}$, together with Lemma \ref{trace}, and (ii) is the same estimate with $f_1=g_{1,1}=-s^{-1}r^{-2}\dr_{-1}\dr u_\theta$. Part (iii) is the square of the last estimate of Lemma \ref{mixed wedge problem}, for the trace of $q_{f_1}$ on $\partial_0\Omega$, with the same $f_1$. For (iv) it is the first estimate of Lemma \ref{mixed harmonic extension}.
\end{proof}

Since $p_{\partial_1,\star}=\E g_{1,\star}$, Lemma \ref{pressure estimates}(iv) covers both of the boundary pressures carried by the data: taking for $\phi$ the data $g_{1,f_\theta}$ and $g_{1,\phi_0}$ of Definition \ref{g definition}, it gives
\begin{align*}
\|(r\grad)^jp_{\partial_1,f_\theta}[f_\theta]\|_{\Omega,a}&\lesssim_j\|f_\theta\|_{\partial_1\Omega,a-{1\over 2},j-{1\over 2}}\\
\|(r\grad)^jp_{\partial_1,\phi_0}[\phi_0]\|_{\Omega,a}&\lesssim_j(\sin\theta_b)^{2}\|s^{-1}\partial_r\phi_0\|_{\partial_1\Omega,a-{1\over 2},j-{1\over 2}}.
\end{align*}

In the test function problems, the pressure terms are of the form $q^\theta_{a},q_a$ and $q^r_a$ defined and estimated as follows. They are the potentials that the projections $\P$ and $\P^\theta$ remove from the test functions of Definition \ref{test functions}. For instance, $\P(r^{-2a}\dr^{2j}\mb v)=r^{-2a}\dr^{2j}\mb v-\dr_{2a}^{2j}\grad q^r_a[\mb v]$ by Lemma \ref{projections commute with dr}.

\begin{lemma}[Pressures of the test function problem]\label{test problem pressure estimate}
Let $q^r_{a}[\mb v]:=q[r^{-2a}\mb v]$, $q_a[\mb v]:=q[\grad^\perp(r^{2-2a}\grad^\perp\cdot\mb v)]$ and $q^\theta_{a}[\mb v]:=q^\theta[\grad^\perp(r^{2-2a}\grad^\perp\cdot\mb v)]$ (as in Definition \ref{the two projections}), where in (i) we assume $\grad\cdot\mb v=0$ and $v_\theta|_{\partial_0\Omega}=0$. Let $K\geq 1$ and suppose $\theta_b\in(0,\pi/3K)$. Then for $1+a'-a\in(-K,K)\setminus\{0\}$, $l\in\{0,1\}$ we have
\begin{alignat*}{2}
(i)\quad &&\|\dr^j(r\grad)^lq^r_{a}\|_{\Omega,2-a+a'}^2
&\lesssim_{K,j,|a|,|1+a'-a|^{-1}}|2a|^2\theta_b^2\|v_r\|_{\Omega,1+a+a',(j+l)^-}^2\\
(ii)\quad &&\|(r\grad)^jq^{[\theta]}_{a}\|_{\Omega,2-a+a'}^2&\lesssim_{K,j,|a|,|1+a'-a|^{-1}}\brac{\|\grad\mb v\|_{\Omega,a+a',j}+\theta_b^{-1}\|\grad\mb v\|_{\Omega,a+a',j^-}}\\
&&&\qquad\qquad\qquad\qquad\qquad\brac{(1+\theta_b^{-1})\|\grad\mb v\|_{\Omega,a+a',j^-}+\|r\grad^\perp\grad^\perp\cdot\mb v\|_{\Omega,a+a',j^-}},
\end{alignat*}
where the $\theta$ in $[\theta]$ is optional.
\end{lemma}
\begin{proof}
For (i), note that $q^r_{a}=q[r^{-2a}\mb v]$ is such that $\lpc q^r_{a}=\grad\cdot(r^{-2a}\mb v)=-2a r^{-2a-1}v_r$, by $\grad\cdot\mb v=0$, with $q^r_{a}|_{\partial_1\Omega}=0$ and $\partial_\theta q^r_{a}|_{\partial_0\Omega}=r^{1-2a}v_\theta|_{\partial_0\Omega}=0$. This RHS is not in the divergence form of Lemma \ref{mixed wedge problem}, but it is brought into it by a radial antiderivative: with $G:=-2ar^{1-2a}v_r$ and $\mb F:=(\dr^{-1}G)\mb e_r$ we have $\grad\cdot(r^{-1}\mb F)=r^{-2}\dr\dr^{-1}G=-2ar^{-2a-1}v_r$ and $F_\theta=0$, so that $q^r_a=q_{\mb F}$ is the solution of that lemma with $f_0=f_1=0$, its Neumann datum $F_\theta|_{\partial_0\Omega}=0$ being the one just found. Since $F_\theta=0$ we may use the second estimate of Lemma \ref{mixed wedge problem}, at the weight $2-a+a'$ and with the choice $\theta_b$, which pays for the one $\partial_\theta$ derivative contained in $(r\grad)^l$ with a factor $\theta_b$ and one radial derivative on $F_r=\dr^{-1}G$; the same choice applied to $\dr^{j+l}q^r_a$ gives the same bound, so that
\begin{align*}
\|\dr^j(r\grad)^lq^r_{a}\|_{\Omega,2-a+a'}^2
&\lesssim{\theta_b^2\over|\cos((a-a'-1)\theta_b)|^2}\|\dr^{j+l-1}G\|_{\Omega,2-a+a'}^2\\
&\lesssim\theta_b^2\|2ar^{-2a+1}v_r\|_{\Omega,2-a+a',(j+l)^-}^2
\lesssim|2a|^2\theta_b^2\|v_r\|_{\Omega,1+a+a',(j+l)^-}^2,
\end{align*}
where at $j=l=0$ the factor $\dr^{-1}$ left on $G$ is bounded on the line $\Re\lambda=1+a'-a$ by $|1+a'-a|^{-1}$.

For (ii), note that $q^\theta_{a}=q^\theta[\grad^\perp(r^{2-2a}\grad^\perp\cdot\mb v)]$ is such that $\lpc q^\theta_{a}=0$ with $\partial_\theta q^\theta_{a}|_{\partial\Omega}=\dr(r^{2-2a}\grad^\perp\cdot\mb v)|_{\partial\Omega}$. So by Lemma \ref{Neumann wedge problem} at the weight $2-a+a'$, whose hypothesis $0<|a-a'-1|\theta_b<\pi$ follows from the one stated, which also gives $\Gamma_{a-a'-1,\theta_b}=1$ there, we have
\begin{align*}
&\|(r\grad)^jq^{\theta}_{a}\|_{\Omega,2-a+a'}^2\\
&\lesssim\norm{\dr^{j-{3\over 2}}\brac{\dr(r^{2-2a}\grad^\perp\cdot\mb v)}}_{\partial\Omega,{3\over 2}-a+a'}^2+\theta_b^{-1}\norm{\dr^{j-2}\brac{\dr(r^{2-2a}\grad^\perp\cdot\mb v)}}_{\partial\Omega,{3\over 2}-a+a'}^2\\
&=\norm{\dr^{j-{1\over 2}}\brac{r^{2-2a}\grad^\perp\cdot\mb v}}_{\partial\Omega,{3\over 2}-a+a'}^2+\theta_b^{-1}\norm{\dr^{j-1}\brac{r^{2-2a}\grad^\perp\cdot\mb v}}_{\partial\Omega,{3\over 2}-a+a'}^2\\
&\lesssim\brac{\norm{\dr^j\brac{r^{2-2a}\grad^\perp\cdot\mb v}}_{\Omega,2-a+a'}+\theta_b^{-1}\norm{\dr^{j-1}\brac{r^{2-2a}\grad^\perp\cdot\mb v}}_{\Omega,2-a+a'}}\\
&\qquad\brac{\theta_b^{-1}\norm{\dr^{j-1}\brac{r^{2-2a}\grad^\perp\cdot\mb v}}_{\Omega,2-a+a'}+\norm{\grad^\perp\dr^{j-1}\brac{r^{2-2a}\grad^\perp\cdot\mb v}}_{\Omega,1-a+a'}}\\
&\lesssim\brac{\|\grad^\perp\cdot\mb v\|_{\Omega,a+a',j}+\theta_b^{-1}\|\grad^\perp\cdot\mb v\|_{\Omega,a+a',j^-}
}\brac{(1+\theta_b^{-1})\|\grad^\perp\cdot\mb v\|_{\Omega,a+a',j^-}+\|r\grad^\perp\grad^\perp\cdot\mb v\|_{\Omega,a+a',j^-}},
\end{align*}
where we have used Lemma \ref{trace} to get the second inequality. The proof for $q_a$ is similar.
\end{proof}

The implicit constants of Lemma \ref{test problem pressure estimate}(i) depend on $|1+a'-a|^{-1}$ only when $j=0$, through the factor $\dr^{-1}$ that the count $0^-=0$ absorbs on the Mellin line $\Re\lambda=1+a'-a$.

\subsection{Coercivity of the linear system}

In the following definition we construct the sesquilinear form under which we will obtain solutions to System \ref{linear system} via the Lax--Milgram theorem. We call it an energy form throughout: the word refers to the quantity produced by testing the equation against the test functions of Definition \ref{test functions}, not to a conserved quantity of the evolution. The forms carry a small parameter $\epsilon>0$, which lets us treat the higher derivative counts as perturbations of the lowest one. It is chosen as $\epsilon=c_\epsilon\theta_b$ in Proposition \ref{energy form coercivity}.

\begin{definition}[Energy forms and functionals]\label{energy forms and functionals}
Let $l\in\N_0$, $a\in\R$ and $\epsilon>0$. Recall $\P^\theta$ and $p_\bullet$ as defined in Definition \ref{the two projections} and \ref{pressure functionals} respectively. For $j\in\N_0$ we define
\begin{align*}
B^r_{a,j}[\mb u,\mb v]
&:=2\inner{\S(\dr_{-2a}^{j}\mb u),r^{-2a}\S(\dr^{j}\mb v)}_{L^2(\Omega)}-2a\inner{\dr_{1-2a}^{j^-}\mb e_r\cdot(2\S(\mb u)-p_{\P|\mb u}I),r^{-2a-1}(-1)^{j-j^-}\dr^{j^+}\mb v}_{L^2(\Omega)}\\
&\quad+\inner{\dr_{1-2a}^{j}u_r,r^{-2a}\dr^{j}v_r}_{L^2(\partial_0\Omega)}+{1\over\bar s}\inner{\dr_{-2a-1}^{j}\dr u_\theta,r^{-2a-2}\dr_{-2a}\dr^{j}v_\theta}_{L^2(\partial_1\Omega)}\\
B^{\theta}_{a,j}[\mb u,\mb v]
&:=\inner{\dr_{1}^{j}(r\lpc\mb u),r^{-1}\dr^{j}\P^\theta\grad^\perp(r^{2-2a}\grad^\perp\cdot\mb v)}_{L^2(\Omega)}\\
\mathcal F^r_{a,j}[\mb v]&:=\inner{\dr_{1-2a}^{j^-}r(\P\mb F+\grad(p_{\partial_1,\phi_0}+p_{\partial_1,f_\theta})),-r^{-2a-1}(-1)^{j-j^-}\dr^{j^+}\mb v}_{L^2(\Omega)}\\
&\quad-\inner{\dr_{1-2a}^{j}f_r,-r^{-2a}\dr^{j}v_r}_{L^2(\partial_1\Omega)}\\
\mathcal F^{\theta}_{a,j}[\mb v]&:=\inner{\dr_{1}^{j}(r\P\mb F),r^{-1}\dr^{j}\P^\theta\grad^\perp(r^{2-2a}\grad^\perp\cdot\mb v)}_{L^2(\Omega)}
\end{align*}
and
\begin{align}
B_{(l)}[\mb u,\mb v]&:=\brac{1+|s|^{-2\kappa}}B_{0,(l)}[\mb u,\mb v]+\epsilon^5B_{-\kappa,(l)}[\mb u,\mb v]\label{main energy form}\\
B_{a,(l)}[\mb u,\mb v]&:=\sum_{j=0}^{l+1}\epsilon^jB^r_{a,j}[\mb u,\mb v]+\epsilon^{2}\sum_{j=0}^l\epsilon^jB^{\theta}_{a,j}[\mb u,\mb v]\\
\mathcal F_{(l)}[\mb v]&:=\brac{1+|s|^{-2\kappa}}\mathcal F_{0,(l)}[\mb v]+\epsilon^5\mathcal F_{-\kappa,(l)}[\mb v]\\
\mathcal F_{a,(l)}[\mb v]&:=\sum_{j=0}^{l+1}\epsilon^j\mathcal F^r_{a,j}[\mb v]+\epsilon^{2}\sum_{j=0}^l\epsilon^j\mathcal F^{\theta}_{a,j}[\mb v]\\
\mb V^{(1)}_{(l)}(\mb v)&:=\mb V^{(1),r}_{(l)}(\mb v)+\mb V^{(1),\theta}_{(l)}(\mb v),\qquad\quad \mb V^{(1),\bullet}_{(l)}(\mb v):=\brac{1+|s|^{-2\kappa}}\mb V^{\bullet}_{0,(l)}(\mb v)+\epsilon^5\mb V^{\bullet}_{-\kappa,(l)}(\mb v)\label{V 1}\\
\mb V^{r}_{a,(l)}(\mb v)&:=\sum_{j=0}^{l+1}\epsilon^j\mb V^r_{a,j}(\mb v),\qquad\qquad\qquad\ \,\mb V^{\theta}_{a,(l)}(\mb v):=\epsilon^{2}\sum_{j=0}^l\epsilon^j\P^\theta\mb V^\theta_{a,j}(\mb v),
\end{align}
where $\bullet$ is a placeholder for $r$ and $\theta$.
\end{definition}

The shift of $\dr$ in $B^\theta_{a,j}$ and $\mathcal F^\theta_{a,j}$ is $1$ for both weights: it is fixed by the pairing weight $r^{-1}$ alone, the weight $r^{2-2a}$ sitting inside the test function. The next proposition expresses the forms as the equation tested against the test functions, and bounds the functionals on $H_s^l$.

\begin{proposition}[Identities and bounds for the energy forms]\label{energy forms}
Let $l\in\N_0$ and $\epsilon>0$. For $\mb u,\mb v\in H_s^l$ (as defined in Definition \ref{H_s^l space}) and $a\in\{0,-\kappa\}$ we have the following: (recall $\mb V^\bullet_{a,j}$ as defined in Definition \ref{test functions})
\begin{enumerate}
\item For $j\in\{0,\cdots,l+1\}$ we have
\begin{align}
B^r_{a,j}[\mb u,\mb v]&=\<\lpc\mb u-\grad p_{\P|\mb u},\mb V^r_{a,j}(\mb v)\>_{L^2(\Omega)}\nonumber\\
&\quad+\sum_{m=0}^1(-1)^m\<(2\mb e_\theta\cdot\S(\mb u)\cdot\mb e_r-\delta_{0m}u_r),\mb e_r\cdot\mb V^r_{a,j}(\mb v)\>_{L^2(\partial_m\Omega)}\label{energy form undo 1}
\end{align}
provided $\mb v$ has the extra regularity $\|\mb v\|_{\Omega,a+1,2(l+1),1}^2 +\|v_\theta\|_{\partial_1\Omega,a+1,2(l+1)}^2 <\infty$.
\item For $j\in\{0,\cdots,l\}$ we have
\begin{align*}
B^{\theta}_{a,j}[\mb u,\mb v]&=\<\lpc\mb u-\grad p_{\P|\mb u},\P^\theta\mb V^\theta_{a,j}(\mb v)\>_{L^2(\Omega)}
\end{align*}
provided $\mb v$ has the extra regularity $\|(r\grad)\grad\mb v\|_{\Omega,a,2l}^2<\infty$.
\item For every $\varpi\in(0,\kappa]$ we have
\begin{align*}
\brac{1+|s|^{-2\kappa}}&\|\mathcal F^r_{0,0}\|_{H_s^l\to\C}\\
&\lesssim_\varpi\brac{1+|s|^{-2\kappa}}\bigg(\|r(\P\mb F+\grad p_{\partial_1,f_\theta})\|_{\Omega,\varpi} +\|f_r\|_{\partial_1\Omega,\varpi-{1\over 2}}+|s|^{\kappa-1}|(\sin\theta_b)^{2}\phi_0(0)|\bigg)\\
&\quad+\brac{1+|s|^{-\kappa}}(\sin\theta_b)^2\|s^{-{1\over 2}(1+\alpha)}\phi_0\|_{\partial_1\Omega,{\alpha\over 2}}
\end{align*}
and when $j\in\{1,\cdots,l+1\}$,
\begin{align*}
&\brac{1+|s|^{-2\kappa}}\|\mathcal F^r_{0,j}\|_{H_s^l\to\C}\\
&\lesssim\brac{1+|s|^{-\kappa}}\bigg(\|r(\P\mb F+\grad p_{\partial_1,f_\theta})\|_{\Omega,0,j^-} +\|\dr^{-{1\over 2}}f_r\|_{\partial_1\Omega,-{1\over 2},j}+(\sin\theta_b)^2\|s^{-{1\over 2}(1+\alpha)}\phi_0\|_{\partial_1\Omega,{\alpha\over 2},j+{1\over 2}}\bigg)
\end{align*}
for $\alpha\in[0,1]$. When $j\in\{0,\cdots,l+1\}$,
\begin{align*}
\|\mathcal F^r_{-\kappa,j}\|_{H_s^l\to\C}
&\lesssim\|r(\P\mb F+\grad(p_{\partial_1,\phi_0}+p_{\partial_1,f_\theta}))\|_{\Omega,-\kappa,j^-} +\|\dr^{-{1\over 2}}f_r\|_{\partial_1\Omega,-\kappa-{1\over 2},j}
\end{align*}
When the RHS is bounded to make $\mathcal F^r_{a,j}:H_s^l\to\C$, $j\in\{0,\cdots,l+1\}$, bounded linear maps, we have
\begin{align*}
\mathcal F^r_{a,j}[\mb v]&=\<\P\mb F+\grad(p_{\partial_1,\phi_0}+p_{\partial_1,f_\theta}),\mb V^r_{a,j}(\mb v)\>_{L^2(\Omega)}
-\<f_r,\mb e_r\cdot\mb V^r_{a,j}(\mb v)\>_{L^2(\partial_1\Omega)}
\end{align*}
provided $\mb v$ has the extra regularity $\|\mb v\|_{\Omega,a+1,2(l+1),1}^2<\infty$.
\item For $k\in\{0,1\}$ and $j\in\{0,\cdots,l\}$,
\begin{align*}
\brac{1+\delta_{0k}|s|^{-2\kappa}}\|\mathcal F^\theta_{-\kappa k,j}\|_{H_s^l\to\C}
\lesssim\brac{1+\delta_{0k}|s|^{-\kappa}}\|r(\P\mb F+\grad(p_{\partial_1,\phi_0}+p_{\partial_1,f_\theta}))\|_{\Omega,-\kappa k,j}
\end{align*}
When the RHS is bounded to make $\mathcal F^{\theta}_{a,j}:H_s^l\to\C$, $j\in\{0,\cdots,l\}$, bounded linear maps, we have
\begin{align*}
\mathcal F^{\theta}_{a,j}[\mb v]&=\<\P\mb F+\grad(p_{\partial_1,\phi_0}+p_{\partial_1,f_\theta}),\P^\theta\mb V^{\theta}_{a,j}(\mb v)\>_{L^2(\Omega)}
\end{align*}
provided $\mb v$ has the extra regularity $\|(r\grad)\grad\mb v\|_{\Omega,a,2l}^2<\infty$.
\end{enumerate}
In particular, if the assumptions of (i)--(iv) are satisfied, then we have
\begin{align}
B_{(l)}[\mb u,\mb v]-\mathcal F_{(l)}[\mb v]&=\<\lpc\mb u-\grad p_\P-\P\mb F,\mb V^{(1)}_{(l)}(\mb v)\>_{L^2(\Omega)}\nonumber\\
&\quad+\sum_{m=0}^1(-1)^m\<(2\mb e_\theta\cdot\S(\mb u)\cdot\mb e_r-\delta_{0m}u_r-\delta_{1m}f_r),\mb e_r\cdot\mb V^{(1),r}_{(l)}(\mb v)\>_{L^2(\partial_m\Omega)}.\label{undo integration main}
\end{align}
\end{proposition}
\begin{proof}
For (i), starting from $\<\lpc\mb u-\grad p_{\P|\mb u},\mb V^r_{a,j}(\mb v)\>_{L^2(\Omega)}$, our aim is to distribute derivatives equally in the two arguments of the inner product using integration by parts, in a way so that the energy form $B^r_{a,j}[\mb u,\mb v]$ obtained is as if the boundary conditions that we expect hold. Write $\mathbb T_{\P|\mb u}[\mb u]=2\S(\mb u)-p_{\P|\mb u}I$. We have
\begin{align}
\<\lpc\mb u - \grad p_{\P|\mb u},\mb V^r_{a,j}(\mb v)\>_{L^2(\Omega)}
&=-\<\lpc\mb u-\grad p_{\P|\mb u},r^{-2a}(-1)^j\dr^{2j}\mb v\>_{L^2(\Omega)}
=-\<\grad\cdot \mathbb T_{\P|\mb u},r^{-2a}(-1)^j\dr^{2j}\mb v\>_{L^2(\Omega)}\nonumber\\
&=\underbrace{\<\mathbb T_{\P|\mb u},r^{-2a}(-1)^j\grad\dr^{2j}\mb v\>_{L^2(\Omega)}}_{E}\underbrace{-2a\<\mb e_r\cdot \mathbb T_{\P|\mb u},r^{-2a-1}(-1)^j\dr^{2j}\mb v\>_{L^2(\Omega)}}_{R}\nonumber\\
&\quad+\underbrace{\sum_{m=0}^1(-1)^m\<\mb e_\theta\cdot \mathbb T_{\P|\mb u},r^{-2a}(-1)^j\dr^{2j}\mb v\>_{L^2(\partial_m\Omega)}}_{E_\partial}\label{energy form proof start}
\end{align}
Now we will show a series of equalities for $E$, $R$ and $E_\partial$ respectively, so that substituting them back into \eqref{energy form proof start} gives us exactly \eqref{energy form undo 1}. Firstly, for $E$:
\begin{align*}
E&=\<\mathbb T_{\P|\mb u},r^{-2a}(-1)^j\grad\dr^{2j}\mb v\>_{L^2(\Omega)}
=\<\mathbb T_{\P|\mb u},r^{-2a}(-1)^j\dr_1^{j}\grad\dr^{j}\mb v\>_{L^2(\Omega)}\\
\intertext{noting that $(p_{\P|\mb u}I):\grad\dr^{j}\mb v=p_{\P|\mb u}\grad\cdot\dr^{j}\mb v$ and that $\grad\cdot\dr^{j}\mb v=0$ (since $\mb v\in H_s^l$) we have}
&=2\<\S(\mb u),r^{-2a}(-1)^j\dr_1^{j}\grad\dr^{j}\mb v\>_{L^2(\Omega)}
=2\<\dr_{1-2a}^{j}\S(\mb u),r^{-2a}\grad\dr^{j}\mb v\>_{L^2(\Omega)}\\
&=2\<\S(\dr_{-2a}^{j}\mb u),r^{-2a}\grad\dr^{j}\mb v\>_{L^2(\Omega)}
=2\<\S(\dr_{-2a}^{j}\mb u),r^{-2a}\S(\dr^{j}\mb v)\>_{L^2(\Omega)}.
\end{align*}
Secondly, for $E_\partial$:
\begin{align*}
E_\partial&=\sum_{m=0}^1(-1)^m\<\mb e_\theta\cdot \mathbb T_{\P|\mb u},r^{-2a}(-1)^j\dr^{2j}\mb v\>_{L^2(\partial_m\Omega)}\\
&=\sum_{m=0}^1(-1)^m\<\mb e_\theta\cdot \mathbb T_{\P|\mb u}\cdot\mb e_r,r^{-2a}(-1)^j\dr^{2j}v_r\>_{L^2(\partial_m\Omega)}\\
&\quad+\sum_{m=0}^1(-1)^m\<\mb e_\theta\cdot \mathbb T_{\P|\mb u}\cdot\mb e_\theta,r^{-2a}(-1)^j\dr^{2j}v_\theta\>_{L^2(\partial_m\Omega)}\\
\intertext{using the fact that $v_\theta|_{\partial_0\Omega}=0$ (since $\mb v\in H_s^l)$ and $\mb e_\theta\cdot(p_{\P|\mb u}I)\cdot\mb e_r=0$ we get}
&=2\sum_{m=0}^1(-1)^m\<\mb e_\theta\cdot\S(\mb u)\cdot\mb e_r,r^{-2a}(-1)^j\dr^{2j}v_r\>_{L^2(\partial_m\Omega)}
-\<\mb e_\theta\cdot \mathbb T_{\P|\mb u}\cdot\mb e_\theta,r^{-2a}(-1)^j\dr^{2j}v_\theta\>_{L^2(\partial_1\Omega)}\\
\intertext{now writing $2\mb e_\theta\cdot\S(\mb u)\cdot\mb e_r=u_r+(2\mb e_\theta\cdot\S(\mb u)\cdot\mb e_r-u_r)$ on $\partial_0\Omega$, and using that $\mb e_\theta\cdot \mathbb T_{\P|\mb u}\cdot\mb e_\theta=\mb e_\theta\cdot(2\S(\mb u)-p_{\P|\mb u}I)\cdot\mb e_\theta=s^{-1}\partial_r^2u_\theta$ on $\partial_1\Omega$ (see Definition \ref{pressure functionals}) we get}
&=\<u_r,r^{-2a}(-1)^j\dr^{2j}v_r\>_{L^2(\partial_0\Omega)}
+\sum_{m=0}^1(-1)^m\<(2\mb e_\theta\cdot\S(\mb u)\cdot\mb e_r-\delta_{0m}u_r),r^{-2a}(-1)^j\dr^{2j}v_r\>_{L^2(\partial_m\Omega)}\\
&\quad-\inner{{1\over s}\partial_r^2u_\theta,r^{-2a}(-1)^j\dr^{2j}v_\theta}_{L^2(\partial_1\Omega)}\\
&=\<\dr_{1-2a}^{j}u_r,r^{-2a}\dr^{j}v_r\>_{L^2(\partial_0\Omega)}
+{1\over\bar s}\inner{\dr_{-2a}^{j}\partial_ru_\theta,r^{-2a-1}\dr_{-2a}\dr^{j}v_\theta}_{L^2(\partial_1\Omega)}\\
&\quad-\sum_{m=0}^1(-1)^m\<(2\mb e_\theta\cdot\S(\mb u)\cdot\mb e_r-\delta_{0m}u_r),\mb e_r\cdot\mb V^r_{a,j}(\mb v)\>_{L^2(\partial_m\Omega)},
\end{align*}
where $\smash[b]{\mb V^r_{a,j}}$ is as defined in Definition \ref{test functions}. Finally, for $R$:
\begin{align*}
R&=-2a\<\mb e_r\cdot \mathbb T_{\P|\mb u},r^{-2a-1}(-1)^j\dr^{2j}\mb v\>_{L^2(\Omega)}
=-2a\<\dr_{1-2a}^{j^-}\mb e_r\cdot \mathbb T_{\P|\mb u},r^{-2a-1}(-1)^{j-j^-}\dr^{j^+}\mb v\>_{L^2(\Omega)}.
\end{align*}
This completes (i). Note that each radial integration by parts above also produces endpoint terms at $r=0$ and at $r=\infty$, which vanish under the stated extra regularity of $\mb v$; when $j=0$ what is needed is exactly $\|r^{-1-a}\mb v\|_{L^2(\Omega)}+\|r^{-1-a}v_\theta\|_{L^2(\partial_1\Omega)}<\infty$.

For part (ii), moving $\P^\theta$ through the radial derivatives by Lemma \ref{projections commute with dr}, we have
\begin{align*}
\<\lpc\mb u-\grad p_{\P|\mb u},(-1)^j\dr^{2j}\P^\theta\grad^\perp(r^{2-2a}\grad^\perp\cdot\mb v)\>_{L^2(\Omega)}
&=\<\dr_1^{j}(r\lpc\mb u-r\grad p_{\P|\mb u}),r^{-1}\dr^{j}\P^\theta\grad^\perp(r^{2-2a}\grad^\perp\cdot\mb v)\>_{L^2(\Omega)}\\
&=\<\dr_1^{j}(r\lpc\mb u),r^{-1}\dr^{j}\P^\theta\grad^\perp(r^{2-2a}\grad^\perp\cdot\mb v)\>_{L^2(\Omega)},
\end{align*}
where the pressure term $\grad p_{\P|\mb u}$ vanishes against $\P^\theta$ (see Remark \ref{P remark}):
\begin{align*}
\<\dr_1^{j}r\grad p_{\P|\mb u},r^{-1}\dr^{j}\P^\theta\grad^\perp(r^{2-2a}\grad^\perp\cdot\mb v)\>_{L^2(\Omega)}
&=\<\grad\dr_1^{j} p_{\P|\mb u},\dr^{j}\P^\theta\grad^\perp(r^{2-2a}\grad^\perp\cdot\mb v)\>_{L^2(\Omega)}\\
&=-\<\dr_1^{j} p_{\P|\mb u},\grad\cdot\dr^{j}\P^\theta\grad^\perp(r^{2-2a}\grad^\perp\cdot\mb v)\>_{L^2(\Omega)}\\
&\quad+\sum_{m=0}^1(-1)^{m+1}\<\dr_1^{j} p_{\P|\mb u},\dr^{j}\mb e_\theta\cdot\P^\theta\grad^\perp(r^{2-2a}\grad^\perp\cdot\mb v)\>_{L^2(\partial_m\Omega)}\\
&=0.
\end{align*}

Parts (iii) and (iv) are the estimates of the operator norms of $\mathcal F^\bullet_{a,j}$, together with the integration by parts that produces the displayed identities. The bulk term is estimated by moving all the radial derivatives off $\mb v$,
\begin{multline*}
\abs{\inner{\dr_{1-2a}^{j^-}r(\P\mb F+\grad(p_{\partial_1,\phi_0}+p_{\partial_1,f_\theta})),r^{-2a-1}(-1)^{j-j^-}\dr^{j^+}\mb v}_{L^2(\Omega)}}\\
\lesssim\|r(\P\mb F+\grad(p_{\partial_1,\phi_0}+p_{\partial_1,f_\theta}))\|_{\Omega,a,j^-}\|\dr^{j^+}\mb v\|_{\Omega,a+1},
\end{multline*}
and then, for $j\geq 1$, using $\dr^{j^+}\mb v=r\dr_1^{j}\partial_r\mb v$, so that $\|\dr^{j^+}\mb v\|_{\Omega,a+1}=\|\dr_1^{j}\partial_r\mb v\|_{\Omega,a}\lesssim_j\|\grad\mb v\|_{\Omega,a,j}\leq\|\mb v\|_{H_s^l,a}$. For $j=0$ and $a=-\kappa$ we use instead Hardy's inequality (Proposition \ref{Hardy's inequality}) at the weight $1-\kappa$, in the form $\|\mb v\|_{\Omega,1-\kappa}\lesssim_\kappa\|\grad\mb v\|_{\Omega,-\kappa}$ of Lemma \ref{limit at infinity}.

The boundary term is estimated by Plancherel's identity \eqref{Plancherel's identity} on $\partial_1\Omega$, splitting the $j$ radial derivatives as $j-1/2$ on $f_r$ and $1/2$ on $v_r$ and using Lemma \ref{trace} for the latter, which produces the half derivative in $\|\dr^{-1/2}f_r\|_{\partial_1\Omega,\ph,j}$. The factors $1+|s|^{-2\kappa}$ and $1+|s|^{-\kappa}$ come from $\|\ph\|_{H_s^l}$.

The case $a=0$, $j=0$ of $\mathcal F^r_{0,0}$ is the one exception: there is no radial derivative to spare and Hardy's inequality has no margin at the weight $0$, so we pair at the weight $\varpi$ instead,
\begin{align*}
\abs{\inner{\P\mb F+\grad p_{\partial_1,f_\theta},\mb v}_{L^2(\Omega)}}&\leq\|r(\P\mb F+\grad p_{\partial_1,f_\theta})\|_{\Omega,\varpi}\|\mb v\|_{\Omega,1-\varpi},\\
\|\mb v\|_{\Omega,1-\varpi}&\lesssim_\varpi\|\grad\mb v\|_{\Omega,-\varpi}\leq\|\grad\mb v\|_{\Omega,0}+\|\grad\mb v\|_{\Omega,-\kappa}\leq 2\|\mb v\|_{H_s^l},
\end{align*}
by Hardy's inequality at the weight $1-\varpi$, whose margin is $\varpi>0$, applicable by Lemma \ref{limit at infinity} since $\mb v$ vanishes at infinity, and using $0\leq\varpi\leq\kappa$. The boundary term $\<f_r,v_r\>_{L^2(\partial_1\Omega)}$ is paired likewise at the weights $\varpi-{1\over 2}$ and ${1\over 2}-\varpi$, with $\|v_r\|_{\partial_1\Omega,1/2-\varpi}\lesssim\|\mb v\|_{\Omega,1-\varpi}+\|\grad\mb v\|_{\Omega,-\varpi}$ by Lemma \ref{trace}. The choice $\varpi=\kappa$ is the one for which the weight $-\varpi$ is itself an energy weight; the smaller choices are used in Subsection \ref{positive-time construction}.

We also note the following:
\begin{align*}
-\inner{\grad p_{\partial_1,\phi_0},\mb v}_{L^2(\Omega)}
&=-\inner{p_{\partial_1,\phi_0},v_\theta}_{L^2(\partial_1\Omega)}
={1\over\bar s}(\sin\theta_b)^{2}\inner{\partial_r\phi_0,v_\theta}_{L^2(\partial_1\Omega)}\\
&={1\over\bar s}(\sin\theta_b)^{2}\brac{-\overline{\phi_0(0)}v_\theta(0,\theta_b) -\inner{\phi_0,\partial_rv_\theta}_{L^2(\partial_1\Omega)}}.
\end{align*}
The endpoint term $-\bar s^{-1}(\sin\theta_b)^2\overline{\phi_0(0)}v_\theta(0,\theta_b)$ is controlled by the corner value of $v_\theta$, for which we claim
\begin{align}\label{corner value bound}
|v_\theta(0,\theta_b)|\lesssim_{\theta_b,\kappa}|s|^{\kappa}\|\mb v\|_{H_s^l};
\end{align}
this gives $|s|^{-1}|\phi_0(0)||v_\theta(0,\theta_b)|\lesssim|s|^{\kappa-1}|\phi_0(0)|\|\mb v\|_{H_s^l}$, which is the corner term of the statement.

To prove \eqref{corner value bound} we localise to $\Omega_{R,2R}$ with $R:=|s|^{-1}$. Since $\partial_rv_\theta\in L^2(\partial_1\Omega)$, the trace $v_\theta(\ph,\theta_b)$ is continuous on $[0,\infty)$ with $|v_\theta(0,\theta_b)|\leq|v_\theta(r,\theta_b)|+r^{1/2}\|\partial_rv_\theta\|_{L^2(\partial_1\Omega)}$, and averaging over $r\in[R,2R]$ gives
\begin{align*}
|v_\theta(0,\theta_b)|\leq R^{-{1\over 2}}\|v_\theta(\ph,\theta_b)\|_{L^2([R,2R])}+(2R)^{1\over 2}\|\partial_rv_\theta\|_{L^2(\partial_1\Omega)}.
\end{align*}
Since $v_\theta$ vanishes on $\partial_0\Omega$, the fundamental theorem of calculus in $\theta$ and the identity $\partial_\theta v_\theta=r\mb e_\theta\cdot\S(\mb v)\cdot\mb e_\theta-v_r$ of \eqref{stress compo} give $|v_\theta(r,\theta_b)|^2\leq\theta_b\int_0^{\theta_b}(r|\grad\mb v|+|\mb v|)^2\d\theta$, hence, with $r\leq 2R$ on $\Omega_{R,2R}$,
\begin{align*}
R^{-1}\|v_\theta(\ph,\theta_b)\|_{L^2([R,2R])}^2&\leq 4\theta_b\|\grad\mb v\|_{L^2(\Omega_{R,2R})}^2+2\theta_bR^{-2}\|\mb v\|_{L^2(\Omega_{R,2R})}^2,\\
R^{-1}\|\mb v\|_{L^2(\Omega_{R,2R})}&\leq 2R^{-\kappa}\|\mb v\|_{\Omega,1-\kappa}\leq 2\kappa^{-1}R^{-\kappa}\|\grad\mb v\|_{\Omega,-\kappa},
\end{align*}
the last step by Lemma \ref{limit at infinity}. Now $R^{-\kappa}=|s|^\kappa$ and $\|\grad\mb v\|_{\Omega,-\kappa}\leq\|\mb v\|_{H_s^l}$, while $\|\grad\mb v\|_{\Omega,0}$ and $|s|^{-1/2}\|\partial_rv_\theta\|_{L^2(\partial_1\Omega)}$ are bounded by $(1+|s|^{-2\kappa})^{-1/2}\|\mb v\|_{H_s^l}\leq|s|^\kappa\|\mb v\|_{H_s^l}$; inserting these bounds gives \eqref{corner value bound}. The scale $R=|s|^{-1}$ is the one at which the two parts $\|\grad\mb v\|_{\Omega,-\kappa}$ and $|s|^{-1/2}\|\partial_rv_\theta\|_{L^2(\partial_1\Omega)}$ of $\|\mb v\|_{H_s^l}$ contribute equally. The bound of Lemma \ref{regular test functions}(i) below, obtained at the fixed scale $R=1$, would not do here, being weaker at both ends of the range of $|s|$.

For the integral term we have
\begin{align*}
  |s|^{-1}\abs{\inner{\phi_0,\partial_rv_\theta}_{L^2(\partial_1\Omega)}} &\lesssim\inner{r^{-{\alpha\over 2}}|s|^{-{1\over 2}(1+\alpha)}|\phi_0|,r^{{\alpha\over 2}}|s|^{-{1\over 2}(1-\alpha)}|\partial_rv_\theta|^{1-\alpha}|\partial_rv_\theta|^\alpha}_{L^2(\partial_1\Omega)}\\
  &\lesssim\|s^{-{1\over 2}(1+\alpha)}\phi_0\|_{\partial_1\Omega,{\alpha\over 2}}\|\mb v\|_{H_s^l,0}
\end{align*}
for $\alpha\in[0,1]$. The last step is Hölder's inequality applied to the splitting
\begin{align*}
r^{{\alpha\over 2}}|s|^{-{1\over 2}(1-\alpha)}|\partial_rv_\theta|
=\brac{|s|^{-{1\over 2}}|\partial_rv_\theta|}^{1-\alpha}\brac{r^{1\over 2}|\partial_rv_\theta|}^{\alpha},
\end{align*}
of which the first factor is bounded in $L^2(\partial_1\Omega)$ by $\|\mb v\|_{H_s^l,0}$ directly, through its $\partial_1\Omega$ term, and the second by the same quantity through the bulk, since $\|r^{1/2}\partial_rv_\theta\|_{L^2(\partial_1\Omega)}=\|\partial_rv_\theta\|_{\partial_1\Omega,-1/2}\lesssim\|\grad\mb v\|_{\Omega}+\|r\grad^2\mb v\|_{\Omega}$ by Lemma \ref{trace}. The constant of the bound on the integral term carries a factor $\theta_b^{-\alpha/2}$ from the trace lemma, which we absorb into $\lesssim$ since $\theta_b$ is fixed in this section.

The same argument, applied after $j$ radial derivatives have been moved onto $\phi_0$, gives, when $j\in\{1,\cdots,l+1\}$,
\begin{align*}
  \abs{\inner{\dr_{1}^{j^-}(r\grad p_{\partial_1,\phi_0}),r^{-1}\dr^{j^+}\mb v}_{L^2(\Omega)}} &= \abs{\inner{\grad\dr_{1}^{j^-} p_{\partial_1,\phi_0},\dr^{j^+}\mb v}_{L^2(\Omega)}}\\
  &\lesssim\abs{\inner{\dr_{1}^{j^-} p_{\partial_1,\phi_0},\dr^{j^+}v_\theta}_{L^2(\partial_1\Omega)}}\\
  &=\abs{\inner{{1\over s}(\sin\theta_b)^2\dr_{1}^{j^-}\partial_r \phi_0,\dr^{j^+}v_\theta}_{L^2(\partial_1\Omega)}}\\
  &=\abs{\inner{{1\over s}(\sin\theta_b)^2\dr^{j} \phi_0,r^{-1}\dr^{j^+}v_\theta}_{L^2(\partial_1\Omega)}}\\
  &\lesssim(\sin\theta_b)^2\|s^{-{1\over 2}(1+\alpha)}\phi_0\|_{\partial_1\Omega,{\alpha\over 2},j+{1\over 2}}\|\mb v\|_{H_s^l,0}.
\end{align*}

Finally, \eqref{undo integration main} follows by subtracting the identities of (iii) and (iv) from those of (i) and (ii) and summing with the weights of $B_{(l)}$, $\mathcal F_{(l)}$ and $\smash{\mb V^{(1)}_{(l)}}$, since $p_\P=p_{\P|\mb u}+p_{\partial_1,\phi_0}+p_{\partial_1,f_\theta}$ by Definition \ref{pressure functionals}.
\end{proof}

In the next lemma, we show that the terms involving the pressure in Definition \ref{energy forms and functionals} are indeed controllable in our norms. The hypotheses of the lemma hold for $\mb u,\mb v\in H_s^l$, that on $\mb v$ by Lemma \ref{limit at infinity}.

\begin{lemma}[The pressure remainder]\label{pressure remainder}
Let $\kappa\in(1/4,1/2]$, and suppose that $u_\theta|_{\partial_0\Omega}=0$ and that $\mb v(r,\theta)\to\mb 0$ as $r\to\infty$ for a.e.\ $\theta$. Then
\begin{multline*}
\abs{\<\dr_{1+2\kappa}^j\{p_{\partial_0},p_{\partial_1,0}\}[\mb u],r^{2\kappa-1}\dr^lv_r\>_{L^2(\Omega)}}\\
\lesssim\|\grad\mb u\|_{\Omega,-\kappa,j}^{1\over 2}\brac{\theta_b^{-1}\|\grad\mb u\|_{\Omega,-\kappa,j^-}+\|r\grad^2\mb u\|_{\Omega,-\kappa,j^-}}^{1\over 2}\|\dr^l v_r\|_{\Omega,-\kappa+1}
\end{multline*}
\begin{align*}
\abs{\<\dr_{1+2\kappa}^jp_{\partial_1,1}[\mb u],r^{2\kappa-1}\dr^lv_r\>_{L^2(\Omega)}}
&\lesssim|s|^{-{1\over 2}-\kappa}\|\partial_ru_\theta\|_{\partial_1\Omega,0,j+1}\brac{|s|^{-\kappa}\|\grad\mb v\|_{\Omega,0,l^-}}^{{1\over 2\kappa}-1}\|\grad\mb v\|_{\Omega,-\kappa,l^-}^{2-{1\over 2\kappa}}.
\end{align*}
\end{lemma}
\begin{proof}
For the first inequality, using Lemma \ref{pressure estimates} we get
\begin{align*}
\abs{\<\dr_{1+2\kappa}^j\{p_{\partial_0},p_{\partial_1,0}\}[\mb u],r^{2\kappa-1}\dr^lv_r\>_{L^2(\Omega)}}
&\leq\|r^{\kappa}\dr_{1+2\kappa}^j\{p_{\partial_0},p_{\partial_1,0}\}[\mb u]\|_{L^2(\Omega)}\|r^{\kappa-1}\dr^lv_r\|_{L^2(\Omega)}\\
&\lesssim\|\grad\mb u\|_{\Omega,-\kappa,j}^{1\over 2}\brac{\theta_b^{-1}\|\grad\mb u\|_{\Omega,-\kappa,j^-}+\|r\grad^2\mb u\|_{\Omega,-\kappa,j^-}}^{1\over 2}\|\dr^lv_r\|_{\Omega,-\kappa+1}.
\end{align*}

For the last inequality, using Lemma \ref{pressure estimates} we get
\begin{align*}
\abs{\<\dr_{1+2\kappa}^jp_{\partial_1,1}[\mb u],r^{2\kappa-1}\dr^lv_r\>_{L^2(\Omega)}}
&=\abs{\<r^{1\over 2}\dr_{1+2\kappa}^jp_{\partial_1,1}[\mb u],r^{2\kappa-{3\over 2}}\dr^lv_r\>_{L^2(\Omega)}}\\
&\leq\|r^{1\over 2}\dr_{1+2\kappa}^jp_{\partial_1,1}[\mb u]\|_{L^2(\Omega)}\|r^{2\kappa-{3\over 2}}\dr^lv_r\|_{L^2(\Omega)}\\
&\lesssim{1\over|s|}\|\partial_ru_\theta\|_{\partial_1\Omega,0,j+1}\|v_r\|_{\Omega,-2\kappa+{3\over 2},l}\\
\intertext{since $\kappa>1/4$, using Hardy's inequality Proposition \ref{Hardy's inequality} we have}
&\lesssim{1\over|s|}\|\partial_ru_\theta\|_{\partial_1\Omega,0,j+1}\|\dr v_r\|_{\Omega,-2\kappa+{3\over 2},l^-}\\
&\lesssim{1\over|s|}\|\partial_ru_\theta\|_{\partial_1\Omega,0,j+1}\|\grad\mb v\|_{\Omega,-2\kappa+{1\over 2},l^-}\\
&\leq{1\over|s|}\|\partial_ru_\theta\|_{\partial_1\Omega,0,j+1}\|\grad\mb v\|_{\Omega,0,l^-}^{{1\over 2\kappa}-1}\|\grad\mb v\|_{\Omega,-\kappa,l^-}^{2-{1\over 2\kappa}}\\
&\leq|s|^{-{1\over 2}-\kappa}\|\partial_ru_\theta\|_{\partial_1\Omega,0,j+1}\brac{|s|^{-\kappa}\|\grad\mb v\|_{\Omega,0,l^-}}^{{1\over 2\kappa}-1}\|\grad\mb v\|_{\Omega,-\kappa,l^-}^{2-{1\over 2\kappa}},
\end{align*}
where the second-to-last line is by Hölder's inequality Proposition \ref{Holder's inequality} using the fact that $0\leq 2\kappa-1/2\leq\kappa$.
\end{proof}

In the following lemma, we show that $B^{\theta}_{a,j}$ gives us control of the Laplacian.

\begin{lemma}[The Laplacian from $B^\theta_{a,j}$]\label{Laplacian from B}
Let $l\in\N_0$ and let $a$ and $\theta_b$ be such that Lemma \ref{test problem pressure estimate} applies with $a'=0$. For small enough $\epsilon>0$ and divergence-free $\mb u$ we have
\begin{multline}\label{Laplacian from B estimate}
\sum_{j=0}^l\epsilon^j\brac{{1\over 4}\|r\lpc\dr^j\mb u\|_{\Omega,a}^2-\brac{\Re+\sgn(\Im(\bar s^{-1}))\Im}B^{\theta}_{a,j}[\mb u,\mb u]}\\
\lesssim\sum_{j=0}^l\epsilon^j\brac{(1+\theta_b^{-1})\|\dr^j\grad\mb u\|_{\Omega,a,1}^2+\theta_b^{-2}\|\dr^j\grad\mb u\|_{\Omega,a}^2}.
\end{multline}
\end{lemma}
\begin{proof}
We have
\begin{align*}
B^{\theta}_{a,j}[\mb u,\mb u]
&=\<\dr_1^{j}r\lpc\mb u,r^{-1}\dr^{j}\P^\theta\grad^\perp(r^{2-2a}\grad^\perp\cdot\mb u)\>_{L^2(\Omega)}\\
&=\<\dr_1^{j}r\lpc\mb u,r^{-1}\dr^{j}(\grad^\perp(r^{2-2a}\grad^\perp\cdot\mb u)-\grad q^\theta_{a}[\mb u])\>_{L^2(\Omega)}\\
&=\<\dr_1^{j}r\lpc\mb u,r^{-1}\dr^{j}(r^{2-2a}\lpc\mb u+(2-2a)r^{1-2a}\mb e_\theta\grad^\perp\cdot\mb u-\grad q^\theta_{a}[\mb u])\>_{L^2(\Omega)},
\end{align*}
where in the last line we used \eqref{lpc form} in the form $\grad^\perp\grad^\perp\cdot\mb u=\lpc\mb u$ (since $\mb u$ is divergence-free). From $\<\dr_1^{j}r\lpc\mb u,r^{-1}\dr^j(r^{2-2a}\lpc\mb u)\>_{L^2(\Omega)}$ we can extract a good top order energy $\|\dr_1^jr\lpc\mb u\|_{\Omega,a}^2$ with a lower order remainder $\<\dr_1^jr\lpc\mb u,r^{-2a}(\dr_{1-2a}^j-\dr_1^j)r\lpc\mb u\>_{L^2(\Omega)}$:
\begin{align*}
&\sum_{j=0}^l\epsilon^j\brac{\|\dr_1^jr\lpc\mb u\|_{\Omega,a}^2-\brac{\Re+\sgn(\Im(\bar s^{-1}))\Im}B^{\theta}_{a,j}[\mb u,\mb u]}\\
&=-\sum_{j=0}^l\epsilon^j\brac{\Re+\sgn(\Im(\bar s^{-1}))\Im}\Big\<\dr_1^jr\lpc\mb u,r^{-2a}(\dr_{1-2a}^j-\dr_1^j)r\lpc\mb u\\
&\pushright{+(2-2a)r^{-2a}\dr_{1-2a}^j\mb e_\theta\grad^\perp\cdot\mb u-r^{-2}\dr_{-1}^jr\grad q^\theta_{a}[\mb u]\Big\>_{L^2(\Omega)}}\\
&\leq 2\sum_{j=0}^l\epsilon^j\bigg|\Big\<\dr_1^jr\lpc\mb u,r^{-2a}(\dr_{1-2a}^j-\dr_1^j)r\lpc\mb u
+(2-2a)r^{-2a}\dr_{1-2a}^j\mb e_\theta\grad^\perp\cdot\mb u-r^{-2}\dr_{-1}^jr\grad q^\theta_{a}[\mb u]\Big\>_{L^2(\Omega)}\bigg|.
\end{align*}
Using Hölder's inequality and Young's inequality we get
\begin{align*}
&\sum_{j=0}^l\epsilon^j\brac{{1\over 2}\|\dr_1^jr\lpc\mb u\|_{\Omega,a}^2-\brac{\Re+\sgn(\Im(\bar s^{-1}))\Im}B^{\theta}_{a,j}[\mb u,\mb u]}\\
&\lesssim\sum_{j=0}^l\epsilon^j\bigg\|r^{a}\bigg(r^{-2a}(\dr_{1-2a}^j-\dr_1^j)r\lpc\mb u
+(2-2a)r^{-2a}\dr_{1-2a}^j\mb e_\theta\grad^\perp\cdot\mb u-r^{-2}\dr_{-1}^jr\grad q^\theta_{a}[\mb u]\bigg)\bigg\|_{L^2(\Omega)}^2
\end{align*}
For small enough $\epsilon>0$ we get
\begin{align*}
\sum_{j=0}^l\epsilon^j&\brac{{1\over 3}\|\dr_1^jr\lpc\mb u\|_{\Omega,a}^2-\brac{\Re+\sgn(\Im(\bar s^{-1}))\Im}B^{\theta}_{a,j}[\mb u,\mb u]}\\
&\lesssim\sum_{j=0}^l\epsilon^j\norm{r^{a}\brac{(2-2a)r^{-2a}\dr_{1-2a}^j\mb e_\theta\grad^\perp\cdot\mb u-r^{-2}\dr_{-1}^jr\grad q^\theta_{a}[\mb u]}}_{L^2(\Omega)}^2\\
&\lesssim(2-2a)^2\sum_{j=0}^l\epsilon^j\|\dr^j\grad\mb u\|_{\Omega,a}^2+\sum_{j=0}^l\epsilon^j\|r^{a-2}\dr^jr\grad q^\theta_{a}[\mb u]\|_{L^2(\Omega)}^2\\
&\lesssim(2-2a)^2\sum_{j=0}^l\epsilon^j\|\dr^j\grad\mb u\|_{\Omega,a}^2\\
&\quad+\sum_{j=0}^l\epsilon^j\brac{\|\dr^j\grad\mb u\|_{\Omega,a,1}+\theta_b^{-1}\|\dr^j\grad\mb u\|_{\Omega,a}}\brac{\theta_b^{-1}\|\dr^j\grad\mb u\|_{\Omega,a}+\|\dr^jr\grad^\perp\grad^\perp\cdot\mb u\|_{\Omega,a}},
\end{align*}
where we used Lemma \ref{test problem pressure estimate} to estimate the $q^\theta_{a}$ term. Using Young's inequality we can then absorb the last term into the LHS and get \eqref{Laplacian from B estimate}, since $\dr_1^jr\lpc\mb u=r\lpc\dr^j\mb u$.
\end{proof}

In $B^r_{-\kappa,0}[\mb u,\mb v]$ from Definition \ref{energy forms and functionals} we have the remainder term $|\bar s^{-1}\<\partial_ru_\theta,r^{2\kappa-1}v_\theta\>_{L^2(\partial_1\Omega)}|$. The following lemma shows that it is controllable. The hypotheses of the lemma on $\mb v$ hold for $\mb v\in H_s^l$, the decay by Lemma \ref{limit at infinity}.

\begin{lemma}[The weighted surface tension remainder]\label{surface tension weighted remainder}
Suppose $\kappa\in(1/4,1/2]$, $v_\theta|_{\partial_0\Omega}=0$ and $\mb v(r,\theta)\to\mb 0$ as $r\to\infty$ for a.e.\ $\theta$. Then
\begin{align*}
\abs{{1\over\bar s}\inner{\partial_ru_\theta,r^{2\kappa-1}v_\theta}_{L^2(\partial_1\Omega)}}
&\lesssim\theta_b^{1\over 2}|s|^{-\kappa-{1\over 2}}\|\partial_ru_\theta\|_{L^2(\partial_1\Omega)}\brac{|s|^{-\kappa}\|\grad\mb v\|_{L^2(\Omega)}}^{{1\over 2\kappa}-1}\|r^{\kappa}\grad\mb v\|_{L^2(\Omega)}^{2-{1\over 2\kappa}}.
\end{align*}
\end{lemma}
\begin{proof}
By the Cauchy--Schwarz inequality we have
\begin{align*}
\abs{{1\over\bar s}\inner{\partial_ru_\theta,r^{2\kappa-1}v_\theta}_{L^2(\partial_1\Omega)}}
\leq|s|^{-\kappa-{1\over 2}}\|\partial_ru_\theta\|_{L^2(\partial_1\Omega)}\underbrace{|s|^{\kappa-{1\over 2}}\|r^{2\kappa-1}v_\theta\|_{L^2(\partial_1\Omega)}}_I.
\end{align*}
Now we do not estimate the trace by Lemma \ref{trace}, which would cost a factor $\theta_b^{-1/2}$ here, but across the wedge as in \eqref{weighted based control f2}: since $v_\theta$ vanishes on $\partial_0\Omega$, the fundamental theorem of calculus and the Cauchy--Schwarz inequality in $\theta$ give
\begin{align*}
\|r^{2\kappa-1}v_\theta\|_{L^2(\partial_1\Omega)}^2
\leq\theta_b\norm{r^{2\kappa-{1\over 2}}\brac{r^{-1}\partial_\theta v_\theta}}_{L^2(\Omega)}^2,
\end{align*}
which is \eqref{weighted based control f2} with the term on $\partial_0\Omega$ absent. So, using $r^{-1}\partial_\theta v_\theta=\mb e_\theta\cdot\S(\mb v)\cdot\mb e_\theta-r^{-1}v_r$,
\begin{align*}
I&\lesssim\theta_b^{1\over 2}|s|^{\kappa-{1\over 2}}\brac{\|r^{2\kappa-{1\over 2}}\grad\mb v\|_{L^2(\Omega)}+\|r^{2\kappa-{3\over 2}}\mb v\|_{L^2(\Omega)}}
\lesssim\theta_b^{1\over 2}|s|^{\kappa-{1\over 2}}\|r^{2\kappa-{1\over 2}}\grad\mb v\|_{L^2(\Omega)}\\
&\lesssim\theta_b^{1\over 2}|s|^{\kappa-{1\over 2}}\|\grad\mb v\|_{L^2(\Omega)}^{{1\over 2\kappa}-1}\|r^{\kappa}\grad\mb v\|_{L^2(\Omega)}^{2-{1\over 2\kappa}}
=\theta_b^{1\over 2}\brac{|s|^{-\kappa}\|\grad\mb v\|_{L^2(\Omega)}}^{{1\over 2\kappa}-1}\|r^{\kappa}\grad\mb v\|_{L^2(\Omega)}^{2-{1\over 2\kappa}},
\end{align*}
where the second inequality is by the first inequality of Proposition \ref{Hardy's inequality} at $k=2\kappa-1/2$, using the fact that $\kappa>1/4$ and the decay $\mb v\to\mb 0$ at infinity assumed in the statement, without which the step fails (for $\mb v\equiv\mb e_1$ the RHS vanishes while $v_\theta(\ph,\theta_b)=-\sin\theta_b$), and the third inequality is by Hölder's inequality Proposition \ref{Holder's inequality} using the fact that $0\leq 2\kappa-1/2\leq\kappa$. Note the direction of the gain: Lemma \ref{trace} would give $\theta_b^{-1/2}$ in place of $\theta_b^{1/2}$, and it is the vanishing of $v_\theta$ on $\partial_0\Omega$ that buys the whole power of $\theta_b$ between them back.
\end{proof}

In the weighted case $a=-\kappa$, we have remainder terms in $B^r_{-\kappa,0}[\mb u,\mb v]$ (Definition \ref{energy forms and functionals}) that require $\|r^{\kappa-1}\mb u\|_{L^2(\Omega)}$ to control. But $\|r^{\kappa-1}\mb u\|_{L^2(\Omega)}$ is not inherently part of the energy of the system. The following lemma shows that we can control it indeed by the energies coming from both the weighted and unweighted case, with better bounds for small $\theta_b$.

\begin{lemma}[Control of the weighted norms by the unweighted ones]\label{weighted based control}
For all $\mb u$ such that the following quantities are defined, we have
\begin{align}
\|r^{\kappa-{1\over 2}}\mb u\|_{L^2(\partial_1\Omega)}^2&\leq 2\theta_b\|r^{\kappa}\grad\mb u\|_{L^2(\Omega)}^2+2\|r^{\kappa-{1\over 2}}\mb u\|_{L^2(\partial_0\Omega)}^2\label{weighted based control f2}\\
\|r^{\kappa-1}\mb u\|_{L^2(\Omega)}^2&\leq2\theta_b^2\|r^{\kappa}\grad\mb u\|_{L^2(\Omega)}^2+2\theta_b\|r^{\kappa-{1\over 2}}\mb u\|_{L^2(\partial_0\Omega)}^2\label{weighted based control f1}
\end{align}
For $\kappa\in(0,1/2]$ and $\theta_b\in(0,\pi)$, we have
\begin{align*}
\|r^{\kappa-{1\over 2}}\mb u\|_{L^2(\partial_1\Omega)}^2
&\leq3\theta_b\|r^{\kappa}\grad\mb u\|_{L^2(\Omega)}^2
+\theta_b^{6\kappa-3}C_{\kappa}(\|\grad\mb u\|_{L^2(\Omega)}^2+\|\mb u\|_{L^2(\partial_0\Omega)}^2)\\
\|r^{\kappa-1}\mb u\|_{L^2(\Omega)}^2
&\leq3\theta_b^2\|r^{\kappa}\grad\mb u\|_{L^2(\Omega)}^2
+\theta_b^{6\kappa-2}C_{\kappa}(\|\grad\mb u\|_{L^2(\Omega)}^2+\|\mb u\|_{L^2(\partial_0\Omega)}^2).
\end{align*}
\end{lemma}
\begin{proof}
Using the fundamental theorem of calculus and Jensen's inequality we have
\begin{align*}
\int_0^\infty|\mb u(r,\theta_b)|^2r^{2\kappa-1}\d r
&=\int_0^\infty\abs{\mb u(r,0)+\int_0^{\theta_b}\partial_\theta\mb u(r,\theta)\d\theta}^2r^{2\kappa-1}\d r\\
&\leq 2\int_0^\infty|\mb u(r,0)|^2r^{2\kappa-1}\d r+2\int_0^\infty\abs{\int_0^{\theta_b}\partial_\theta\mb u(r,\theta)\d\theta}^2r^{2\kappa-1}\d r\\
&\leq 2\int_0^\infty|\mb u(r,0)|^2r^{2\kappa-1}\d r+2\theta_b\int_0^\infty\int_0^{\theta_b}|\partial_\theta\mb u(r,\theta)|^2 r^{2\kappa-1}\d\theta\d r
\end{align*}
The proof for \eqref{weighted based control f1} is similar.

When $\kappa\in(0,1/2)$, by Lemma \ref{trace}, we have
\begin{align*}
\|r^{{\kappa\over 2}-{1\over 2}}\mb u\|_{L^2(\partial_0\Omega)}^2&\leq\theta_b^{-1}\|r^{{\kappa\over 2}-1}\mb u\|_{L^2(\Omega)}^2+2\|r^{{\kappa\over 2}-1}\mb u\|_{L^2(\Omega)}\|r^{\kappa\over 2}\grad\mb u\|_{L^2(\Omega)}\\
\intertext{using the first inequality of Proposition \ref{Hardy's inequality} at $k=\kappa/2$, which applies because $\mb u$ vanishes at infinity by Lemma \ref{limit at infinity}, the RHSs $\|r^\kappa\grad\mb u\|_{L^2(\Omega)}=\|\grad\mb u\|_{\Omega,-\kappa}$ and $\|\mb u\|_{L^2(\partial_0\Omega)}$ being finite, we have}
&\leq{4\over\theta_b\kappa^2}\|r^{\kappa\over 2}\grad\mb u\|_{L^2(\Omega)}^2+{4\over\kappa}\|r^{\kappa\over 2}\grad\mb u\|_{L^2(\Omega)}^2
\leq\brac{{4\over\theta_b\kappa^2}+{4\over\kappa}}\|\grad\mb u\|_{L^2(\Omega)}\|r^{\kappa}\grad\mb u\|_{L^2(\Omega)}\\
&\leq\brac{{4\over\theta_b\kappa^2}+{4\over\kappa}}\brac{{1\over 2\delta}\|\grad\mb u\|_{L^2(\Omega)}^2+{\delta\over 2}\|r^{\kappa}\grad\mb u\|_{L^2(\Omega)}^2}.
\end{align*}
We have
\begin{align*}
\|r^{\kappa-{1\over 2}}\mb u\|_{L^2(\partial_0\Omega)}^2
&=\|r^{2\kappa-1}|\mb u|^2\|_{L^1(\partial_0\Omega)}
=\norm{r^{2\kappa-1}|\mb u|^{2\kappa-1\over\kappa/2-1/2}|\mb u|^{2-{2\kappa-1\over\kappa/2-1/2}}}_{L^1(\partial_0\Omega)}\\
&\leq\norm{r^{2\kappa-1}|\mb u|^{2\kappa-1\over\kappa/2-1/2}}_{L^{{\kappa-1\over 2\kappa-1}}(\partial_0\Omega)}\norm{|\mb u|^{2-{2\kappa-1\over\kappa/2-1/2}}}_{L^{2\brac{2-{2\kappa-1\over\kappa/2-1/2}}^{-1}}(\partial_0\Omega)}\\
&\leq\|r^{{\kappa\over 2}-{1\over 2}}\mb u\|_{L^2(\partial_0\Omega)}^{4\kappa-2\over\kappa-1}\|\mb u\|_{L^2(\partial_0\Omega)}^{2-{4\kappa-2\over\kappa-1}}\\
&\leq\brac{2\kappa-1\over\kappa-1}\delta'^{\brac{2\kappa-1\over\kappa-1}^{-1}}\|r^{{\kappa\over 2}-{1\over 2}}\mb u\|_{L^2(\partial_0\Omega)}^2+\brac{1-{2\kappa-1\over\kappa-1}}\delta'^{-\brac{1-{2\kappa-1\over\kappa-1}}^{-1}}\|\mb u\|_{L^2(\partial_0\Omega)}^2\\
&\leq\brac{2\kappa-1\over\kappa-1}\brac{{4\over\theta_b\kappa^2}+{4\over\kappa}}\delta'^{\brac{2\kappa-1\over\kappa-1}^{-1}}\brac{{1\over 2\delta}\|\grad\mb u\|_{L^2(\Omega)}^2+{\delta\over 2}\|r^{\kappa}\grad\mb u\|_{L^2(\Omega)}^2}\\
&\quad+\brac{1-{2\kappa-1\over\kappa-1}}\delta'^{-\brac{1-{2\kappa-1\over\kappa-1}}^{-1}}\|\mb u\|_{L^2(\partial_0\Omega)}^2\\
\intertext{we can choose $\delta,\delta'$ small enough and depending on $\theta_b$ so that}
&\leq\epsilon'(1+\theta_b\kappa)\theta_b\|r^{\kappa}\grad\mb u\|_{L^2(\Omega)}^2+\theta_b^{3(2\kappa-1)} C_{\epsilon',\kappa}((1+\theta_b\kappa)\|\grad\mb u\|_{L^2(\Omega)}^2+\|\mb u\|_{L^2(\partial_0\Omega)}^2)
\end{align*}
(When $\kappa=1/2$ the two statements are trivial, since then $r^{\kappa-1/2}=1$ and $\theta_b^{6\kappa-3}=1$; so we may assume $\kappa<1/2$ below.) Here the factors of $\theta_b$ come from taking $\delta\propto\theta_b^{\alpha}$ and $\delta'\propto\theta_b^{\alpha'}$: the two exponent conditions have the unique solution $\alpha=2-3\kappa$ and $\alpha'=3\kappa(2\kappa-1)/(\kappa-1)$, giving the power $3(2\kappa-1)=6\kappa-3$ on $\partial_1\Omega$ and one power more in the bulk.
So we have
\begin{align*}
\|r^{\kappa-1}\mb u\|_{L^2(\Omega)}^2&\leq2\theta_b^2\|r^{\kappa}\grad\mb u\|_{L^2(\Omega)}^2+2\theta_b\|r^{\kappa-{1\over 2}}\mb u\|_{L^2(\partial_0\Omega)}^2\\
&\leq(1+\epsilon'(1+\theta_b\kappa))2\theta_b^2\|r^{\kappa}\grad\mb u\|_{L^2(\Omega)}^2
+2\theta_b^{1+3(2\kappa-1)}C_{\epsilon',\kappa}((1+\theta_b\kappa)\|\grad\mb u\|_{L^2(\Omega)}^2+\|\mb u\|_{L^2(\partial_0\Omega)}^2)
\end{align*}
Since $\theta_b$ and $\kappa$ are bounded, we have
\begin{align*}
\|r^{\kappa-1}\mb u\|_{L^2(\Omega)}^2
&\leq3\theta_b^2\|r^{\kappa}\grad\mb u\|_{L^2(\Omega)}^2
+\theta_b^{6\kappa-2}C_{\kappa}(\|\grad\mb u\|_{L^2(\Omega)}^2+\|\mb u\|_{L^2(\partial_0\Omega)}^2)
\end{align*}
and similarly for $\|r^{\kappa-1/2}\mb u\|_{L^2(\partial_1\Omega)}$, starting from \eqref{weighted based control f2} instead of \eqref{weighted based control f1}.
\end{proof}

\begin{proposition}[Coercivity and boundedness of the energy form]\label{energy form coercivity}
Let $\kappa\in(1/4,1/2]$ and $l\in\N_0$. There is $c_{\epsilon,0}>0$ depending only on $\kappa$ and $l$ such that, for every $c_\epsilon\in(0,c_{\epsilon,0}]$, there is $\theta_0=\theta_0(c_\epsilon,\kappa,l)>0$ with the following property: if $\theta_b\in(0,\theta_0]$ and $\epsilon=c_\epsilon\theta_b$, then for every $s$ with $\Re s>0$ we have
\begin{align*}
\|\mb u\|_{H_s^l}^2&\lesssim_{\kappa,l,c_\epsilon,\theta_b}|B_{(l)}[\mb u,\mb u]|\\
|B_{(l)}[\mb u,\mb v]|&\lesssim_{\kappa,l,c_\epsilon,\theta_b}\|\mb u\|_{H_s^l}\|\mb v\|_{H_s^l},
\end{align*}
where $H_s^l$ is as defined in Definition \ref{H_s^l space} and the sesquilinear form $B_{(l)}:H_s^l\times H_s^l\to\C$ is defined by \eqref{main energy form} in Definition \ref{energy forms and functionals}.
\end{proposition}
\begin{proof}
We first prove coercivity at the weight $0$. We then prove it at the weight $-\kappa$, where the remainders from the pressure and the surface tension are bounded by a small multiple of the weighted energy plus an unweighted remainder. The factor $\epsilon^5$ in $B_{(l)}$ then allows this remainder to be absorbed by the unweighted estimate. Boundedness is proved at the end.

\emph{Step 1: coercivity at the weight $0$.} Write
\begin{align*}
\|\mb u\|_a^2=\|\grad\mb u\|_{\Omega,a}^2+\|\mb u\|_{\partial_0\Omega,a}^2+{1\over|s|}\|\partial_r u_\theta\|_{\partial_1\Omega,a}^2.
\end{align*}
Note that
\begin{align*}
B^r_{0,0}[\mb u]=2\|\S(\mb u)\|_{\Omega}^2+\|u_r\|_{\partial_0\Omega}^2+{1\over\bar s}\|\partial_ru_\theta\|_{\partial_1\Omega}^2
\end{align*}
Note that since $\Re s>0$, we have $\Re(\bar s^{-1})=\Re(s)/|s|^2>0$, so
\begin{align*}
\abs{1\over s}=\abs{1\over\bar s}\lesssim\abs{\Re\brac{1\over\bar s}}+\abs{\Im\brac{1\over\bar s}}=\Re\brac{1\over\bar s}+\sgn(\Im(\bar s^{-1}))\Im\brac{1\over\bar s}.
\end{align*}
Using Korn's inequality (Lemma \ref{Korn on the wedge}), which applies because $\mb u$ is divergence-free with $u_\theta|_{\partial_0\Omega}=0$, and noting that $\|\mb u\|_{\partial_0\Omega}=\|u_r\|_{\partial_0\Omega}$ for the same reason, we then have
\begin{align*}
\|\mb u\|_0^2\lesssim\brac{\Re+\sgn(\Im(\bar s^{-1}))\Im}B^r_{0,0}[\mb u].
\end{align*}
Note that
\begin{align*}
\|\dr^j\mb u\|_0^2\lesssim\brac{\Re+\sgn(\Im(\bar s^{-1}))\Im}B^r_{0,j}[\mb u]+C\|\dr^j\mb u\|_0\sum_{j'=0}^{j-1}\|\dr^{j'}\mb u\|_0
\end{align*}
So, using Young's inequality, for small enough $\epsilon>0$ we have
\begin{align*}
\sum_{j=0}^{l+1}\epsilon^j\|\dr^j\mb u\|_0^2\lesssim\brac{\Re+\sgn(\Im(\bar s^{-1}))\Im}\sum_{j=0}^{l+1}\epsilon^jB_{0,j}^r[\mb u].
\end{align*}
Using Lemma \ref{second derivative bound} and \ref{Laplacian from B}, we see that
\begin{align*}
\sum_{j=0}^{l+1}\epsilon^j\|\dr^j\mb u\|_0^2+\epsilon^{2}\sum_{j=0}^l\epsilon^j\|r\grad^2\dr^j\mb u\|_{\Omega}^2
&\lesssim\sum_{j=0}^{l+1}\epsilon^j\|\dr^j\mb u\|_0^2+\epsilon^{2}\sum_{j=0}^l\epsilon^j\|r\lpc\dr^j\mb u\|_{\Omega}^2\\
&\lesssim\brac{\Re+\sgn(\Im(\bar s^{-1}))\Im}\brac{\sum_{j=0}^{l+1}\epsilon^jB_{0,j}^r[\mb u]+\epsilon^{2}\sum_{j=0}^l\epsilon^jB_{0,j}^{\theta}[\mb u]}\\
&\quad+\epsilon^{2}\sum_{j=0}^l\epsilon^j\brac{(1+\theta_b^{-1})\|\dr^j\grad\mb u\|_{\Omega,0,1}^2+\theta_b^{-2}\|\dr^j\grad\mb u\|_{\Omega,0}^2}.
\end{align*}
Since $\|\dr^j\grad\mb u\|_{\Omega,0,1}\leq\|\dr^j\mb u\|_0+\|\dr^{j+1}\mb u\|_0$, the last line is absorbed into the LHS once $\epsilon^2\theta_b^{-2}$ and $\epsilon\theta_b^{-1}$ are small, that is, once the ratio $c_\epsilon=\epsilon/\theta_b$ is small. Throughout this proof the small parameters are chosen as in the statement: the Young parameters $\delta$ below are fixed first, then the ratio $c_\epsilon=\epsilon/\theta_b$, and finally $\theta_b$, with $\epsilon=c_\epsilon\theta_b$, so that every smallness condition is a condition on $c_\epsilon$ or, once $c_\epsilon$ is fixed, on $\theta_b$ alone; absorptions already performed remain valid when $\theta_b$ is shrunk further, and we shall not repeat this remark at every step. Writing
\begin{align}\label{epsilon weighted energies}
E_a[\mb u]:=\sum_{j=0}^{l+1}\epsilon^j\|\dr^j\mb u\|_a^2+\epsilon^{2}\sum_{j=0}^l\epsilon^j\|r\grad^2\dr^j\mb u\|_{\Omega,a}^2,\qquad\qquad \epsilon^{l+2}\|\mb u\|_{H_s^l,a}^2\lesssim E_a[\mb u]\lesssim\|\mb u\|_{H_s^l,a}^2,
\end{align}
we have obtained
\begin{align}\label{unweighted coercivity}
E_0[\mb u]\lesssim\brac{\Re+\sgn(\Im(\bar s^{-1}))\Im}\brac{\sum_{j=0}^{l+1}\epsilon^jB_{0,j}^r[\mb u]+\epsilon^{2}\sum_{j=0}^l\epsilon^jB_{0,j}^{\theta}[\mb u]}.
\end{align}
We keep the $\epsilon$-weighted energies $E_a[\mb u]$, rather than the norms $\|\mb u\|_{H_s^l,a}$, until the very end: the conversion in \eqref{epsilon weighted energies} costs a factor $\epsilon^{l+2}$, which is harmless once, but which must not be paid before the remainder of the weighted estimate below has been absorbed into the unweighted one.

\emph{Step 2: coercivity at the weight $-\kappa$.} Using Lemma \ref{Korn on the wedge} at $a=-\kappa$ and the fact that $u_\theta|_{\partial_0\Omega}=0$, we have
\begin{align*}
\|\mb u\|_{-\kappa}^2&\lesssim\brac{\Re+\sgn(\Im(\bar s^{-1}))\Im}B^r_{-\kappa,0}[\mb u]\\
&\quad+C\kappa\brac{\abs{{1\over\bar s}\inner{\partial_ru_\theta,r^{2\kappa-1}u_\theta}_{L^2(\partial_1\Omega)}}+\abs{\<\mb e_r\cdot(2\S(\mb u)-p_{\P|\mb u}I),r^{2\kappa-1}\mb u\>_{L^2(\Omega)}}}\\
\intertext{and using Lemma \ref{surface tension weighted remainder}, which applies since $u_\theta$ vanishes on $\partial_0\Omega$, and Lemma \ref{pressure remainder} we have (the surface tension remainder in fact comes with a factor $\theta_b^{2\kappa}$, after Young's inequality, by Lemma \ref{surface tension weighted remainder}, which we discard; the $p_{\partial_1,1}$ part of the pressure remainder is of the same form without it, except that by Lemma \ref{pressure remainder} it carries $\|\partial_ru_\theta\|_{\partial_1\Omega,0,1}$, one radial derivative more than $\|\mb u\|_0$ contains, which is why the unweighted remainder below involves $\|\dr\mb u\|_0$ as well)}
&\lesssim\brac{\Re+\sgn(\Im(\bar s^{-1}))\Im}B^r_{-\kappa,0}[\mb u]
+C\kappa\brac{\delta\|\mb u\|_{-\kappa}^2+C_\delta|s|^{-2\kappa}\sum_{j'=0}^1\|\dr^{j'}\mb u\|_{0}^2}\\
&\quad+C\kappa\brac{(1+\theta_b^{-{1\over 2}})\|\grad\mb u\|_{\Omega,-\kappa}+\theta_b^{1\over 2}\|r\grad^2\mb u\|_{\Omega,-\kappa}}\|r^{\kappa-1}\mb u\|_{L^2(\Omega)}\\
\intertext{and using Lemma \ref{weighted based control} (for $\|r^{\kappa-1}\mb u\|_{L^2(\Omega)}$) we have}
&\lesssim\brac{\Re+\sgn(\Im(\bar s^{-1}))\Im}B^r_{-\kappa,0}[\mb u]+C\kappa\brac{\delta\|\mb u\|_{-\kappa}^2+C_\delta|s|^{-2\kappa}\sum_{j'=0}^1\|\dr^{j'}\mb u\|_{0}^2}\\
&\quad+C\kappa\brac{(1+\theta_b^{-{1\over 2}})\|\grad\mb u\|_{\Omega,-\kappa}+\theta_b^{1\over 2}\|r\grad^2\mb u\|_{\Omega,-\kappa}}\\
&\pushright{\brac{\theta_b\|r^{\kappa}\grad\mb u\|_{L^2(\Omega)}+\theta_b^{3\kappa-1}C_{\kappa}(\|\grad\mb u\|_{L^2(\Omega)}+\|\mb u\|_{L^2(\partial_0\Omega)})}}\\
&\lesssim\brac{\Re+\sgn(\Im(\bar s^{-1}))\Im}B^r_{-\kappa,0}[\mb u]+C\kappa\brac{\delta\|\mb u\|_{-\kappa}^2+C_\delta|s|^{-2\kappa}\sum_{j'=0}^1\|\dr^{j'}\mb u\|_{0}^2}\\
&\quad+C\kappa\brac{(\theta_b+\theta_b^{1\over 2})\|\grad\mb u\|_{\Omega,-\kappa}^2+\theta_b^{5\over 2}\|r\grad^2\mb u\|_{\Omega,-\kappa}^2+\theta_b^{6\kappa-{7\over 2}}C_{\kappa}\|\mb u\|_{0}^2}
\end{align*}
Taking $\delta$ and $\theta_b$ small enough, in the order described above, we get
\begin{align*}
\|\mb u\|_{-\kappa}^2
&\lesssim\brac{\Re+\sgn(\Im(\bar s^{-1}))\Im}B^r_{-\kappa,0}[\mb u]+C\kappa\theta_b^{5\over 2}\|r\grad^2\mb u\|_{\Omega,-\kappa}^2\\
&\quad+C\kappa\brac{\theta_b^{6\kappa-{7\over 2}}C_{\kappa}+C_\delta|s|^{-2\kappa}}\sum_{j'=0}^1\|\dr^{j'}\mb u\|_{0}^2.
\end{align*}

For $j>0$, using Lemma \ref{Korn on the wedge} applied to $\dr^j\mb u$, which is again divergence-free with vanishing $\theta$ component on $\partial_0\Omega$ since $\dr$ is tangential there, we have
\begin{align*}
\|\dr^j\mb u\|_{-\kappa}^2&\lesssim\brac{\Re+\sgn(\Im(\bar s^{-1}))\Im}B^r_{-\kappa,j}[\mb u]+C\sum_{j'=0}^{j}\|\dr^{j'}\mb u\|_{-\kappa}\sum_{j'=0}^{j-1}\|\dr^{j'}\mb u\|_{-\kappa}\\
&\quad+C\kappa\abs{\<\dr_{1-2a}^{j-1}\mb e_r\cdot(2\S(\mb u)-p_{\P|\mb u}I),r^{-2a-1}(-1)\dr^{j+1}\mb u\>_{L^2(\Omega)}}\\
&\lesssim\brac{\Re+\sgn(\Im(\bar s^{-1}))\Im}B^r_{-\kappa,j}[\mb u]+C\sum_{j'=0}^{j}\|\dr^{j'}\mb u\|_{-\kappa}\sum_{j'=0}^{j-1}\|\dr^{j'}\mb u\|_{-\kappa}\\
&\quad+C\kappa\abs{\<\dr_{1-2a}^{j}\mb e_r\cdot(2\S(\mb u)-(p_{\partial_0}+p_{\partial_1,0})I),r^{-2a-1}\dr^{j}\mb u\>_{L^2(\Omega)}}\\
&\quad+C\kappa\abs{\<\dr_{1-2a}^{j-1}\mb e_r\cdot p_{\partial_1,1}I,r^{-2a-1}\dr^{j+1}\mb u\>_{L^2(\Omega)}}\\
\intertext{with $a=-\kappa$, and using Lemma \ref{pressure remainder} we have}
&\lesssim\brac{\Re+\sgn(\Im(\bar s^{-1}))\Im}B^r_{-\kappa,j}[\mb u]+C\sum_{j'=0}^{j}\|\dr^{j'}\mb u\|_{-\kappa}\sum_{j'=0}^{j-1}\|\dr^{j'}\mb u\|_{-\kappa}\\
&\quad+C\kappa\brac{(1+\theta_b^{-{1\over 2}})\|\grad\mb u\|_{\Omega,-\kappa,j}+\theta_b^{1\over 2}\|r\grad^2\mb u\|_{\Omega,-\kappa,j-1}}\sum_{j'=0}^{j-1}\|\dr^{j'}\mb u\|_{-\kappa}\\
&\quad+C\kappa|s|^{-{1\over 2}-\kappa}\|\partial_ru_\theta\|_{\partial_1\Omega,0,j}\brac{|s|^{-\kappa}\|\grad\mb u\|_{\Omega,0,j}}^{{1\over 2\kappa}-1}\|\grad\mb u\|_{\Omega,-\kappa,j}^{2-{1\over 2\kappa}}.
\end{align*}
The last term is a product of powers of $|s|^{-\kappa}(\sum_{j'\leq j}\|\dr^{j'}\mb u\|_0^2)^{1/2}$ and of $\|\grad\mb u\|_{\Omega,-\kappa,j}$ with exponents adding up to $2$, so Young's inequality bounds it by $\delta\|\grad\mb u\|_{\Omega,-\kappa,j}^2+C_\delta|s|^{-2\kappa}\sum_{j'\leq j}\|\dr^{j'}\mb u\|_0^2$; the cross terms of the second line are handled by Young's inequality with a parameter of the size $\theta_b$, which is where $c_\epsilon$ has to be small. Summing with the weights $\epsilon^j$ we get
\begin{align*}
\sum_{j=0}^{l+1}\epsilon^j\|\dr^j\mb u\|_{-\kappa}^2
&\lesssim\sum_{j=0}^{l+1}\epsilon^j\brac{\Re+\sgn(\Im(\bar s^{-1}))\Im}B^r_{-\kappa,j}[\mb u]
+C\kappa\theta_b^{5\over 2}\|r\grad^2\mb u\|_{\Omega,-\kappa}^2\\
&\quad+C\epsilon\sum_{j=1}^{l+1}\epsilon^j\brac{(1+\theta_b^{-1})\|\grad\mb u\|_{\Omega,-\kappa,j}^2+\theta_b\|r\grad^2\mb u\|_{\Omega,-\kappa,j-1}^2}
+C\kappa\delta\sum_{j=1}^{l+1}\epsilon^j\|\grad\mb u\|_{\Omega,-\kappa,j}^2\\
&\quad+C\kappa\brac{\theta_b^{6\kappa-{7\over 2}}C_{\kappa}+C_\delta|s|^{-2\kappa}}\sum_{j=0}^{l+1}\epsilon^j\sum_{j'=0}^{\max\{j,1\}}\|\dr^{j'}\mb u\|_{0}^2.
\end{align*}
Taking $\delta$ and then $c_\epsilon$ small enough, the second line is absorbed into the LHS, except for the terms carrying $r\grad^2\mb u$, which we carry along, and we get
\begin{align*}
&\sum_{j=0}^{l+1}\epsilon^j\|\dr^j\mb u\|_{-\kappa}^2\\
&\lesssim\sum_{j=0}^{l+1}\epsilon^j\brac{\Re+\sgn(\Im(\bar s^{-1}))\Im}B^r_{-\kappa,j}[\mb u]
+C\brac{\theta_b^{5\over 2}\|r\grad^2\mb u\|_{\Omega,-\kappa}^2+\epsilon^{2}\theta_b\sum_{j=0}^{l}\epsilon^j\|r\grad^2\mb u\|_{\Omega,-\kappa,j}^2}\\
&\quad+C\brac{\theta_b^{6\kappa-{7\over 2}}C_{\kappa}+C_\delta|s|^{-2\kappa}}\epsilon^{-1}E_0[\mb u],
\end{align*}
where we used that $\sum_{j=0}^{l+1}\epsilon^j\sum_{j'\leq\max\{j,1\}}\|\dr^{j'}\mb u\|_0^2\lesssim\epsilon^{-1}\sum_{j'=0}^{l+1}\epsilon^{j'}\|\dr^{j'}\mb u\|_0^2\leq\epsilon^{-1}E_0[\mb u]$ for $\epsilon\leq{1\over 2}$, the factor $\epsilon^{-1}$ being the price of the one extra radial derivative at $j=0$. Using Lemma \ref{second derivative bound} and \ref{Laplacian from B}, we see that
\begin{align*}
E_{-\kappa}[\mb u]&\lesssim\sum_{j=0}^{l+1}\epsilon^j\|\dr^j\mb u\|_{-\kappa}^2+\epsilon^{2}\sum_{j=0}^l\epsilon^j\|r\lpc\dr^j\mb u\|_{\Omega,-\kappa}^2\\
&\lesssim\brac{\Re+\sgn(\Im(\bar s^{-1}))\Im}\brac{\sum_{j=0}^{l+1}\epsilon^jB^r_{-\kappa,j}[\mb u]+\epsilon^{2}\sum_{j=0}^l\epsilon^jB_{-\kappa,j}^{\theta}[\mb u]}\\
&\quad+\epsilon^{2}\sum_{j=0}^l\epsilon^j\brac{(1+\theta_b^{-1})\|\dr^j\grad\mb u\|_{\Omega,-\kappa,1}^2+\theta_b^{-2}\|\dr^j\grad\mb u\|_{\Omega,-\kappa}^2}\\
&\quad+C\brac{\theta_b^{5\over 2}\|r\grad^2\mb u\|_{\Omega,-\kappa}^2+\epsilon^{2}\theta_b\sum_{j=0}^{l}\epsilon^j\|r\grad^2\mb u\|_{\Omega,-\kappa,j}^2}\\
&\quad+C\brac{\theta_b^{6\kappa-{7\over 2}}C_{\kappa}+C_\delta|s|^{-2\kappa}}\epsilon^{-1}E_0[\mb u].
\end{align*}
The third line is absorbed into the LHS once $c_\epsilon$ is small, exactly as in the unweighted case; and the two terms carrying $r\grad^2\mb u$ are absorbed into the term $\epsilon^2\|r\grad^2\mb u\|_{\Omega,-\kappa}^2$ and its higher counts in $E_{-\kappa}[\mb u]$ once $\theta_b$ is small: the first because $\theta_b^{5/2}\ll\epsilon^2=c_\epsilon^2\theta_b^2$ amounts to $\theta_b^{1/2}\ll c_\epsilon^2$, which holds for $\theta_b$ small depending on $c_\epsilon$, and the second because $\theta_b\ll 1$. Note that $c_\epsilon$ is thus required to be small in some places and $\theta_b^{1/4}$ times a constant to be smaller than it in another, which is why $c_\epsilon$ is fixed first and $\theta_b$ afterwards. We arrive at
\begin{align}\label{weighted coercivity}
E_{-\kappa}[\mb u]&\lesssim\brac{\Re+\sgn(\Im(\bar s^{-1}))\Im}\brac{\sum_{j=0}^{l+1}\epsilon^jB_{-\kappa,j}^r[\mb u]+\epsilon^{2}\sum_{j=0}^l\epsilon^jB_{-\kappa,j}^{\theta}[\mb u]}\nonumber\\
&\quad+C\brac{\theta_b^{6\kappa-{7\over 2}}C_{\kappa}+C_\delta|s|^{-2\kappa}}\epsilon^{-1}E_0[\mb u].
\end{align}

\emph{Step 3: absorption of the unweighted remainder.} Combining \eqref{weighted coercivity} with the unweighted estimate \eqref{unweighted coercivity}, we get
\begin{align*}
\brac{1+|s|^{-2\kappa}}E_0[\mb u]+\epsilon^5E_{-\kappa}[\mb u]
&\lesssim\brac{\Re+\sgn(\Im(\bar s^{-1}))\Im}B_{(l)}[\mb u]
+C\epsilon^4\brac{\theta_b^{6\kappa-{7\over 2}}C_{\kappa}+C_\delta|s|^{-2\kappa}}E_0[\mb u].
\end{align*}
The last term is absorbed into the first term on the left: the part carrying $|s|^{-2\kappa}$ once $C C_\delta\epsilon^4\leq{1\over 2}$, and the part carrying $\theta_b^{6\kappa-7/2}$ once $CC_\kappa\epsilon^4\theta_b^{6\kappa-7/2}=CC_\kappa c_\epsilon^4\theta_b^{6\kappa+1/2}\leq{1\over 2}$, both of which hold for $\theta_b$ small; here it matters that the absorption is performed on the $\epsilon$-weighted energies, since on the norms $\|\ph\|_{H_s^l,a}$ the unweighted estimate carries the factor $\epsilon^{l+2}$ of \eqref{epsilon weighted energies}, against which $\epsilon^4$ has no room for $l\geq 2$. Hence, by \eqref{epsilon weighted energies},
\begin{align*}
\epsilon^{l+7}\|\mb u\|_{H_s^l}^2&\lesssim\epsilon^{l+2}\brac{1+|s|^{-2\kappa}}\|\mb u\|_{H_s^l,0}^2+\epsilon^{l+7}\|\mb u\|_{H_s^l,-\kappa}^2
\lesssim\brac{1+|s|^{-2\kappa}}E_0[\mb u]+\epsilon^5E_{-\kappa}[\mb u]\\
&\lesssim\brac{\Re+\sgn(\Im(\bar s^{-1}))\Im}B_{(l)}[\mb u]\lesssim|B_{(l)}[\mb u]|.
\end{align*}
Thus we get the coercivity estimate $\|\mb u\|_{H_s^l}^2\lesssim|B_{(l)}[\mb u]|$, with a constant depending on $\epsilon$, that is, on $c_\epsilon$ and $\theta_b$.

\emph{Step 4: boundedness.} The bound $|B_{(l)}[\mb u,\mb v]|\lesssim\|\mb u\|_{H_s^l}\|\mb v\|_{H_s^l}$ is the same inspection with the two arguments separated: the symmetric gradient, Navier slip and surface tension terms of $B^r_{a,j}$ are pairings of summands of $\|\mb u\|_{H_s^l,a}$ and $\|\mb v\|_{H_s^l,a}$ by the Cauchy--Schwarz inequality, the weights $r^{-2a-2}$ of $B^r_{a,j}$ and $r^{-2a}$ of $\|\ph\|_{H_s^l,a}$ matching in the surface tension term because $\dr u_\theta=r\partial_ru_\theta$. Two terms are exceptions. The term of $B^r_{a,j}$ with the factor $-2a$, present at $a=-\kappa$ only, pairs $\dr^{j^-}_{1+2\kappa}\mb e_r\cdot(2\S(\mb u)-p_{\P|\mb u}I)$ against $r^{2\kappa-1}\dr^{j^+}\mb v$, whose second factor is bounded by $\|\grad\mb v\|_{\Omega,-\kappa,j}$ for $j\geq 1$ and by Hardy's inequality at the weight $1-\kappa$ for $j=0$, and whose pressure part is bounded by Lemma \ref{pressure remainder}, the second estimate being closed by the weight $1+|s|^{-2\kappa}$ of $\|\ph\|_{H_s^l}$. In $B^\theta_{a,j}$ the test function is bounded by Lemmas \ref{projection norm estimate} and \ref{test problem pressure estimate}(ii). Summing over $j$ with the weights $\epsilon^j$ and $\epsilon^5$ of $B_{(l)}$ gives the bound, with a constant depending on $\kappa$, $l$, $c_\epsilon$ and $\theta_b$ but not on $s$.
\end{proof}

With coercivity of $B_{(l)}$ and boundedness of $\mathcal F_{(l)}$, by the Lax--Milgram theorem we can obtain $\mb u$ such that $B_{(l)}[\mb u,\mb v]-\mathcal F_{(l)}[\mb v]=0$ for all $\mb v\in H_s^l$. For regular enough $\mb v$, we have \eqref{undo integration main}. But to show that our equation holds a.e., we will have to show that the test functions $\smash{\mb V^{(1)}_{(l)}(\mb v)}$ exhaust the bulk and the boundary. That is the test function problem we will need to solve.

\subsection{The test function problem}\label{the test function problem}

The test functions $\smash[t]{\mb V^{(1)}_{(l)}}(\mb v)$ of \eqref{V 1} are built from those of Definition \ref{test functions} at the two weights $a=0$ and $a=-\kappa$. They therefore have two scalings in the bulk, the source of the weighted and unweighted norms in $H_s^l$. Hence, for given $\tilde{\mb v}\in C^\infty_c(\Omega)$, the equation $\P\smash[t]{\trimbot{0.6pt}{\mb V^{(1)}_{(l)}}}(\mb v)=\P\tilde{\mb v}$ cannot be solved for $\mb v$ directly in Mellin variables. We test $\smash[t]{\trimbot{0.5pt}{\mb V^{(1)}_{(l)}}}(\mb v)$ instead against a second set of test functions $\smash[t]{\trimbot{0.5pt}{\mb V^{(2)}_{(l)}}}(\mb v)$ with one scaling, defined in Proposition \ref{auxiliary variational problem}. This reduces the problem to solving for $\smash[t]{\mb V^{(2)}_{(l)}}(\mb v)$, which is done in Proposition \ref{solvability for V2}. For this we have to prove coercivity of $\smash{\<\P\mb V^{(1)}_{(l)}(\mb v),\mb V^{(2)}_{(l)}(\mb v)\>_{L^2(\Omega)}}$. Since $\trimtop{3.1pt}{\trimbot{1.5pt}{\mb V^{(1)}_{(l)}}}$ and $\trimtop{3.1pt}{\trimbot{1.5pt}{\mb V^{(2)}_{(l)}}}$ are each a sum of $\mb V^r_{a,j}$ and $\mb V^\theta_{a',l}$, we consider the product of these component terms and show that they are coercive modulo some good remainder terms. This is done in the next three lemmas.

\begin{lemma}[$\mb V^r$ against $\mb V^r$ coercivity]\label{Vr against Vr}
Let $j,l\in\N_0$ and $a,a'\in\R$ with $a\neq a'$. For $\theta_b$ small enough and $\mb v\in C^\infty_{c,\sigma}(\Omega)$ (as defined in Definition \ref{smooth bulk test functions}), we have
\begin{align}\label{Vr against Vr estimate}
{1\over 3}\|\dr^{j+l}\mb v\|_{\Omega,a+a'}^2&-\Re\<\P\mb V^r_{a,j}(\mb v),\mb V^r_{a',l}(\mb v)\>_{L^2(\Omega)}
\lesssim\|\mb v\|_{\Omega,a+a',j+l-1}^2.
\end{align}
\end{lemma}
\begin{proof}
We have, by Lemma \ref{projections commute with dr}, which is what makes the shift $\dr_{2a}^{2j}$ on the potential term correct,
\begin{align*}
\<\P\mb V^r_{a,j}(\mb v),\mb V^r_{a',l}(\mb v)\>_{L^2(\Omega)}
&=\inner{\P(r^{-2a}(-1)^j\dr^{2j}\mb v),r^{-2a'}(-1)^l\dr^{2l}\mb v}_{L^2(\Omega)}\\
&=\inner{r^{-2a}(-1)^j\dr^{2j}\mb v-(-1)^j\dr_{2a}^{2j}\grad q^r_{a}[\mb v],r^{-2a'}(-1)^l\dr^{2l}\mb v}_{L^2(\Omega)}\\
&=\inner{r^{-2a}\dr^{l}_{2-2(a+a')}\dr^{j}\mb v,r^{-2a'}\dr^{j}_{2-2(a+a')}\dr^{l}\mb v}_{L^2(\Omega)}\\
&\quad-\inner{\dr^{l}_{2-2a'}\dr_{2a}^{j}\grad q^r_{a}[\mb v],r^{-2a'}\dr_{2-2(a'+a)}^{j}\dr^{l}\mb v}_{L^2(\Omega)},
\end{align*}
where $q^r_a$ is as defined in Lemma \ref{test problem pressure estimate}. So
\begin{align*}
\|\dr^{j+l}\mb v\|_{\Omega,a+a'}^2&-\Re\<\P\mb V^r_{a,j}(\mb v),\mb V^r_{a',l}(\mb v)\>_{L^2(\Omega)}\\
&\lesssim\|\mb v\|_{\Omega,a+a',j+l}\brac{\|\mb v\|_{\Omega,a+a',j+l-1}+\|\grad q^r_{a}[\mb v]\|_{\Omega,a'-a,j+l}}
\end{align*}
Using Young's inequality, we can absorb the $\|\dr^{j+l}\mb v\|_{\Omega,a+a'}$ in $\|\mb v\|_{\Omega,a+a',j+l}$ to the LHS to get
\begin{align*}
{1\over 2}\|\dr^{j+l}\mb v\|_{\Omega,a+a'}^2-\Re\<\P\mb V^r_{a,j}(\mb v),\mb V^r_{a',l}(\mb v)\>_{L^2(\Omega)}
&\lesssim\|\mb v\|_{\Omega,a+a',j+l-1}^2+\|\grad q^r_{a}[\mb v]\|_{\Omega,a'-a,j+l}^2\\
&\lesssim\|\mb v\|_{\Omega,a+a',j+l-1}^2+\theta_b^2\|\mb v\|_{\Omega,a+a',j+l}^2,
\end{align*}
where we used Lemma \ref{test problem pressure estimate}(i), with $a'-1$ in place of $a'$, for the last line. Its condition $1+a'-a\neq 0$ then becomes the hypothesis $a\neq a'$.

Taking $\theta_b$ small enough, we can absorb the $\|\dr^{j+l}\mb v\|_{\Omega,a+a'}$ in $\|\mb v\|_{\Omega,a+a',j+l}$ to the LHS to get \eqref{Vr against Vr estimate}.
\end{proof}

\begin{lemma}[$\mb V^r$ against $\mb V^\theta$ coercivity]\label{Vr against Vtheta}
Let $j,l'\in\N_0$, $a,a'\in\R$ with $a\neq a'$, and $\delta>0$. For $\theta_b$ and $\epsilon>0$ small enough and $\mb v\in C^\infty_{c,\sigma}(\Omega)$ (as defined in Definition \ref{smooth bulk test functions}), we have
\begin{align}
&\sum_{l=0}^{l'}\epsilon^l\brac{{1\over C}\|r\grad(\dr^{j+l}\mb v)\|_{\Omega,a+a'}^2-\Re\<[\P^\theta]\mb V^\theta_{a,j}(\mb v),[\P]\mb V^r_{a',l}(\mb v)\>_{L^2(\Omega)}}\nonumber\\
&\lesssim\sum_{l=0}^{l'}\epsilon^l(\theta_b^{1\over 2}+\theta_b^{-1})\|r\grad\mb v\|_{\Omega,a+a',j+l-1}^2+\sum_{l=0}^{l'}\epsilon^l\theta_b^{-2}\|\mb v\|_{\Omega,a+a',j+l}^2\nonumber\\
&\quad+\sum_{l=-1}^{l'}\epsilon^l(\theta_b^{1\over 2}+\theta_b^{-1}+\delta^{-1}\theta_b^{1\over 2})\|\mb v\|_{\Omega,a+a',j+l+1}^2+\delta\theta_b^{1\over 2}\sum_{l=1}^{l'}\epsilon^{l}\|r^2\grad^\perp\grad^\perp\cdot\mb v\|_{\Omega,a'+a,j+l-1}^2,\label{Vr against Vtheta estimate}
\end{align}
where the projections in $[\ph]$ are optional.
\end{lemma}
\begin{proof}
We integrate by parts until the pairing becomes the principal term $2\<r^{2-2a}\S(\ph),r^{-2a'}\S(\ph)\>_{L^2(\Omega)}$, which controls $\|r\grad(\dr^{j+l}\mb v)\|_{\Omega,a+a'}^2$ by Korn's inequality (Lemma \ref{Korn on the wedge}). Three kinds of remainders appear: the one boundary term on $\partial_1\Omega$ that the boundary conditions of $C^\infty_{c,\sigma}(\Omega)$ leave, the two bulk terms in which a derivative falls on a power of $r$, and, when the projections are present, the potential terms, which are bounded by Lemma \ref{test problem pressure estimate}. Young's inequality bounds all of them by a small multiple of the principal term plus the RHS of the statement.

\emph{Step 1: the pairing without projections.} We have
\begin{align*}
\<\mb V^\theta_{a,j}(\mb v)&,\mb V^r_{a',l}(\mb w)\>_{L^2(\Omega)}\\
&=\inner{(-1)^j\dr^{2j}\grad^\perp(r^{2-2a}\grad^\perp\cdot\mb v),-r^{-2a'}(-1)^l\dr^{2l}\mb w}_{L^2(\Omega)}\\
&=\inner{\dr^{l^-}_{2-2a'}\dr^{j}\grad^\perp(r^{2-2a}\grad^\perp\cdot\mb v),-r^{-2a'}(-1)^{l-l^-}\dr^{j}_{2-2a'}\dr^{l^+}\mb w}_{L^2(\Omega)}\\
&=\inner{\dr^{l^-}_{2-2a'}\dr^{j}\brac{(2-2a)(r^{2-2a-1}\mb e_\theta\grad^\perp\cdot\mb v)+r^{2-2a}\lpc\mb v},-r^{-2a'}(-1)^{l-l^-}\dr^{j}_{2-2a'}\dr^{l^+}\mb w}_{L^2(\Omega)}\\
&=2\inner{r^{2-2a}\grad\cdot\S(\dr^{l^-}_{2-2(a+a')}\dr^{j}_{-2a}\mb v),-r^{-2a'}(-1)^{l-l^-}\dr^{j}_{2-2a'}\dr^{l^+}\mb w}_{L^2(\Omega)}\\
&\quad+(2-2a)\inner{r^{1-2a}\mb e_\theta\grad^\perp\cdot\dr^{l^-}_{2-2(a+a')}\dr^{j}_{-2a}\mb v,-r^{-2a'}(-1)^{l-l^-}\dr^{j}_{2-2a'}\dr^{l^+}\mb w}_{L^2(\Omega)}\\
&=2\inner{r^{2-2a}\S(\dr^{l^-}_{2-2(a+a')}\dr^{j}_{-2a}\mb v),r^{-2a'}(-1)^{l-l^-}\grad\dr^{j}_{2-2a'}\dr^{l^+}\mb w}_{L^2(\Omega)}\\
&\quad+2(2-2(a+a'))\inner{r^{1-2a}\mb e_r\cdot\S(\dr^{l^-}_{2-2(a+a')}\dr^{j}_{-2a}\mb v),r^{-2a'}(-1)^{l-l^-}\dr^{j}_{2-2a'}\dr^{l^+}\mb w}_{L^2(\Omega)}\\
&\quad-2\sum_{m=0}^1(-1)^m\inner{r^{2-2a}\mb e_\theta\cdot\S(\dr^{l^-}_{2-2(a+a')}\dr^{j}_{-2a}\mb v),-r^{-2a'}(-1)^{l-l^-}\dr^{j}_{2-2a'}\dr^{l^+}\mb w}_{L^2(\partial_m\Omega)}\\
&\quad+(2-2a)\inner{r^{1-2a}\mb e_\theta\grad^\perp\cdot\dr^{l^-}_{2-2(a+a')}\dr^{j}_{-2a}\mb v,-r^{-2a'}(-1)^{l-l^-}\dr^{j}_{2-2a'}\dr^{l^+}\mb w}_{L^2(\Omega)}\\
&=2\inner{r^{2-2a}\S(\dr^{l}_{2-2(a+a')}\dr^{j}_{-2a}\mb v),r^{-2a'}\grad\dr^{j}_{2-2a'}\dr^{l}\mb w}_{L^2(\Omega)}\\
&\quad-2\sum_{m=0}^1(-1)^m\inner{r^{1-2a}\dr^{l^-}_{2-2(a+a')}\dr^{j}_{-2a}\mb e_\theta\cdot r\S(\mb v),-r^{-2a'}(-1)^{l-l^-}\dr^{j}_{2-2a'}\dr^{l^+}\mb w}_{L^2(\partial_m\Omega)}\\
&\quad+2(2-2(a+a'))\inner{r^{1-2a}\mb e_r\cdot\S(\dr^{l^-}_{2-2(a+a')}\dr^{j}_{-2a}\mb v),r^{-2a'}(-1)^{l-l^-}\dr^{j}_{2-2a'}\dr^{l^+}\mb w}_{L^2(\Omega)}\\
&\quad+(2-2a)\inner{r^{1-2a}\mb e_\theta\grad^\perp\cdot\dr^{l^-}_{2-2(a+a')}\dr^{j}_{-2a}\mb v,-r^{-2a'}(-1)^{l-l^-}\dr^{j}_{2-2a'}\dr^{l^+}\mb w}_{L^2(\Omega)}\\
&=2\inner{r^{2-2a}\S(\dr^{l}_{2-2(a+a')}\dr^{j}_{-2a}\mb v),r^{-2a'}\S(\dr^{j}_{2-2a'}\dr^{l}\mb w)}_{L^2(\Omega)}\\
&\quad-2\sum_{m=0}^1(-1)^m\inner{r^{1-2a}\dr^{l^-}_{2-2(a+a')}\dr^{j}_{-2a}\mb e_\theta\cdot r\S(\mb v)\cdot\mb e_r,-r^{-2a'}(-1)^{l-l^-}\dr^{j}_{2-2a'}\dr^{l^+}w_r}_{L^2(\partial_m\Omega)}\\
&\quad-2\sum_{m=0}^1(-1)^m\inner{r^{1-2a}\dr^{l^-}_{2-2(a+a')}\dr^{j}_{-2a}\mb e_\theta\cdot r\S(\mb v)\cdot\mb e_\theta,-r^{-2a'}(-1)^{l-l^-}\dr^{j}_{2-2a'}\dr^{l^+}w_\theta}_{L^2(\partial_m\Omega)}\\
&\quad+2(2-2(a+a'))\inner{r^{1-2a}\mb e_r\cdot\S(\dr^{l}_{2-2(a+a')}\dr^{j}_{-2a}\mb v),r^{-2a'}\dr^{j}_{2-2a'}\dr^{l}\mb w}_{L^2(\Omega)}\\
&\quad+(2-2a)\inner{r^{1-2a}\mb e_\theta\grad^\perp\cdot\dr^{l}_{2-2(a+a')}\dr^{j}_{-2a}\mb v,-r^{-2a'}\dr^{j}_{2-2a'}\dr^{l}\mb w}_{L^2(\Omega)}\\
&=2\inner{r^{2-2a}\S(\dr^{l}_{2-2(a+a')}\dr^{j}_{-2a}\mb v),r^{-2a'}\S(\dr^{j}_{2-2a'}\dr^{l}\mb w)}_{L^2(\Omega)}\\
&\quad+2\inner{r^{1-2a}\dr^{l^-}_{2-2(a+a')}\dr^{j}_{-2a}\mb e_\theta\cdot r\S(\mb v)\cdot\mb e_\theta,-r^{-2a'}(-1)^{l-l^-}\dr^{j}_{2-2a'}\dr^{l^+}w_\theta}_{L^2(\partial_1\Omega)}\\
&\quad+2(2-2(a+a'))\inner{r^{1-2a}\mb e_r\cdot\S(\dr^{l}_{2-2(a+a')}\dr^{j}_{-2a}\mb v),r^{-2a'}\dr^{j}_{2-2a'}\dr^{l}\mb w}_{L^2(\Omega)}\\
&\quad+(2-2a)\inner{r^{1-2a}\mb e_\theta\grad^\perp\cdot\dr^{l}_{2-2(a+a')}\dr^{j}_{-2a}\mb v,-r^{-2a'}\dr^{j}_{2-2a'}\dr^{l}\mb w}_{L^2(\Omega)},
\end{align*}
where because of the boundary conditions in $C^\infty_{c,\sigma}(\Omega)$ all but one of the boundary terms vanish, and where we used $\lpc\mb v=2\grad\cdot\S(\mb v)$ for divergence-free $\mb v$. Noting that $\mb e_\theta\cdot r\S(\mb v)\cdot\mb e_\theta=-\dr v_r$ \eqref{stress compo}, we can estimate the boundary term:
\begin{align*}
&\abs{\inner{r^{1-2a}\dr^{l^-}_{2-2(a+a')}\dr^{j}_{-2a}\mb e_\theta\cdot r\S(\mb v)\cdot\mb e_\theta,-r^{-2a'}(-1)^{l-l^-}\dr^{j}_{2-2a'}\dr^{l^+}v_\theta}_{L^2(\partial_1\Omega)}}\\
&=\abs{\inner{r^{1-2a}\dr^{l^-}_{2-2(a+a')}\dr^{j}_{-2a}\dr v_r,r^{-2a'}\dr^{j}_{2-2a'}\dr^{l^+}v_\theta}_{L^2(\partial_1\Omega)}}
\intertext{using Plancherel's identity \eqref{Plancherel's identity} we have}
&=\abs{\int_{\Re\lambda=a+a'-1}(\bar\lambda^{l^-}_{2-2(a+a')}\bar\lambda^{j}_{-2a}\bar\lambda\bar v_r(\lambda,\theta_b))(\lambda^{j}_{2-2a'}\lambda^{l^+}\hat v_\theta(\lambda,\theta_b))\d\Im\lambda}\\
&=\abs{\int_{\Re\lambda=a+a'-1}\theta_b^{-1}\int_0^{\theta_b}\partial_\theta\brac{\theta(\bar\lambda^{l^-}_{2-2(a+a')}\bar\lambda^{j}_{-2a}\bar\lambda\bar v_r(\lambda,\theta))(\lambda^{j}_{2-2a'}\lambda^{l^+}\hat v_\theta(\lambda,\theta))}\d\Im\lambda}\\
&\lesssim\|\mb v\|_{\Omega,a+a',j+l+1}\|\mb v\|_{\Omega,a+a',j+l,1}+\theta_b^{-1}\|\mb v\|_{\Omega,a+a',j+l+1}\|\mb v\|_{\Omega,a+a',j+l}
\end{align*}
So we have
\begin{align*}
\|r\S(\dr^{j+l}\mb v)\|_{\Omega,a+a'}^2&-\Re\<\mb V^\theta_{a,j}(\mb v),\mb V^r_{a',l}(\mb v)\>_{L^2(\Omega)}\\
&\lesssim\|r\S(\mb v)\|_{\Omega,a+a',j+l}\|r\S(\mb v)\|_{\Omega,a+a',j+l-1}+\|\mb v\|_{\Omega,a+a',j+l+1}\|\mb v\|_{\Omega,a+a',j+l,1}\\
&\quad+\theta_b^{-1}\|\mb v\|_{\Omega,a+a',j+l+1}\|\mb v\|_{\Omega,a+a',j+l}+\|r\grad\mb v\|_{\Omega,a+a',j+l}\|\mb v\|_{\Omega,a+a',j+l}
\end{align*}
Using Korn's inequality (Lemma \ref{Korn on the wedge}, applied to $\dr^{j+l}\mb v$ at the weight $a+a'-1$) and Young's inequality we get
\begin{align*}
{1\over C}\|r\grad(\dr^{j+l}\mb v)\|_{\Omega,a+a'}^2&-\Re\<\mb V^\theta_{a,j}(\mb v),\mb V^r_{a',l}(\mb v)\>_{L^2(\Omega)}\\
&\lesssim\|r\grad\mb v\|_{\Omega,a+a',j+l-1}^2+\theta_b^{-2}\|\mb v\|_{\Omega,a+a',j+l}^2+\|\mb v\|_{\Omega,a+a',j+l+1}^2.
\end{align*}

\emph{Step 2: the potential terms.} When both projections are present, the second one is redundant by Remark \ref{P remark}, because $\P^\theta\mb V^\theta_{a,j}(\mb v)$ is divergence-free with vanishing $\theta$ component on $\partial_0\Omega$; so at most one potential term occurs. Now, by Lemma \ref{projections commute with dr} in each potential term,
\begin{align*}
\<\P^\theta\mb V^\theta_{a,j}(\mb v),&\mb V^r_{a',l}(\mb w)\>_{L^2(\Omega)}\\
&=\<\mb V^\theta_{a,j}(\mb v),\mb V^r_{a',l}(\mb w)\>_{L^2(\Omega)}-\inner{(-1)^j\dr^{2j}\grad q^\theta_{a}[\mb v],-r^{-2a'}(-1)^l\dr^{2l}\mb w}_{L^2(\Omega)}\\
&=\<\mb V^\theta_{a,j}(\mb v),\mb V^r_{a',l}(\mb w)\>_{L^2(\Omega)}-\inner{\dr^{l^-}_{2-2a'}\dr^{j}\grad q^\theta_{a}[\mb v],-r^{-2a'}(-1)^{l-l^-}\dr^{j}_{2-2a'}\dr^{l^+}\mb w}_{L^2(\Omega)}.
\end{align*}
\begin{align*}
\<\mb V^\theta_{a,j}(\mb v),&\P\mb V^r_{a',l}(\mb w)\>_{L^2(\Omega)}\\
&=\<\mb V^\theta_{a,j}(\mb v),\mb V^r_{a',l}(\mb w)\>_{L^2(\Omega)}-\inner{(-1)^j\dr^{2j}\grad^\perp(r^{2-2a}\grad^\perp\cdot\mb v),(-1)^l\dr^{2l}_{2a'}\grad q^r_{a'}[\mb w]}_{L^2(\Omega)}\\
&=\<\mb V^\theta_{a,j}(\mb v),\mb V^r_{a',l}(\mb w)\>_{L^2(\Omega)}-\inner{\dr^{l^-}_{2+2a'}\dr^{j}\grad^\perp(r^{2-2a}\grad^\perp\cdot\mb v),(-1)^{l-l^-}\dr^{j}_2\dr^{l^+}_{2a'}\grad q^r_{a'}[\mb w]}_{L^2(\Omega)},
\end{align*}
and using Lemma \ref{test problem pressure estimate} the remainder terms can be estimated as
\begin{align*}
&\abs{\inner{\dr^{l^-}_{2-2a'}\dr^{j}\grad q^\theta_{a}[\mb v],-r^{-2a'}(-1)^{l-l^-}\dr^{j}_{2-2a'}\dr^{l^+}\mb w}_{L^2(\Omega)}}\\
&\lesssim\|\grad q^\theta_{a}[\mb v]\|_{\Omega,a'-a,j+l^-}\|\mb w\|_{\Omega,a+a',j+l^+}\\
&\lesssim\brac{(1+\theta_b^{-{1\over 2}})\|r\grad\mb v\|_{\Omega,a'+a,j+l^-+1}+\theta_b^{-1}\|r\grad\mb v\|_{\Omega,a'+a,j+l^-}+\theta_b^{1\over 2}\|r^2\grad^\perp\grad^\perp\cdot\mb v\|_{\Omega,a'+a,j+l^-}}\\
&\qquad\|\mb w\|_{\Omega,a+a',j+l^+}.
\end{align*}
\begin{align*}
&\abs{\inner{\dr^{l^-}_{2+2a'}\dr^{j}\grad^\perp(r^{2-2a}\grad^\perp\cdot\mb v),(-1)^{l-l^-}\dr^{j}_2\dr^{l^+}_{2a'}\grad q^r_{a'}[\mb w]}_{L^2(\Omega)}}\\
&\lesssim\brac{\|r\grad\mb v\|_{\Omega,a+a',j+l^-}+\|r^2\grad^\perp\grad^\perp\cdot\mb v\|_{\Omega,a+a',j+l^-}}\|\grad q^r_{a'}[\mb w]\|_{\Omega,a-a',j+l^+}\\
&\lesssim|a'|\theta_b^{1\over 2}\brac{\|r\grad\mb v\|_{\Omega,a+a',j+l^-}+\|r^2\grad^\perp\grad^\perp\cdot\mb v\|_{\Omega,a+a',j+l^-}}\|\mb w\|_{\Omega,a+a',j+l^+},
\end{align*}
where we used Lemma \ref{test problem pressure estimate} to bound the $q$ term. That lemma is applied with $a'-1$ in place of $a'$ to $q^\theta_a$, and with $a'$ and $a-1$ in place of $a$ and $a'$ to $q^r_{a'}$. Each time its condition $1+a'-a\neq 0$ becomes the hypothesis $a\neq a'$. In both cases, they can be bounded by
\begin{align*}
&\brac{(1+\theta_b^{-{1\over 2}})\|r\grad\mb v\|_{\Omega,a'+a,j+l^-+1}+(\theta_b^{1\over 2}+\theta_b^{-1})\|r\grad\mb v\|_{\Omega,a'+a,j+l^-}+\theta_b^{1\over 2}\|r^2\grad^\perp\grad^\perp\cdot\mb v\|_{\Omega,a'+a,j+l^-}}\\
&\qquad\|\mb w\|_{\Omega,a+a',j+l^+}.
\end{align*}

\emph{Step 3: conclusion.} Combining Steps 1 and 2 by Young's inequality, we have
\begin{align*}
&{1\over C}\|r\grad(\dr^{j}\mb v)\|_{\Omega,a+a'}^2-\Re\<[\P^\theta]\mb V^\theta_{a,j}(\mb v),[\P]\mb V^r_{a',0}(\mb v)\>_{L^2(\Omega)}\\
&\lesssim\|r\grad\mb v\|_{\Omega,a+a',j-1}^2+\theta_b^{-2}\|\mb v\|_{\Omega,a+a',j}^2+\|\mb v\|_{\Omega,a+a',j+1}^2+((\delta'^{-1}+\delta''^{-1})(\theta_b^{1\over 2}+\theta_b^{-1})+\delta^{-1}\theta_b^{1\over 2})\|\mb v\|_{\Omega,a+a',j}^2\\
&\quad+\delta'\|r\grad\mb v\|_{\Omega,a'+a,j+1}^2+\delta''(\theta_b^{1\over 2}+\theta_b^{-1})\|r\grad\mb v\|_{\Omega,a'+a,j}^2+\delta\theta_b^{1\over 2}\|r^2\grad^\perp\grad^\perp\cdot\mb v\|_{\Omega,a'+a,j}^2
\end{align*}
and when $l>0$,
\begin{align*}
&{1\over C}\|r\grad(\dr^{j+l}\mb v)\|_{\Omega,a+a'}^2-\Re\<[\P^\theta]\mb V^\theta_{a,j}(\mb v),[\P]\mb V^r_{a',l}(\mb v)\>_{L^2(\Omega)}\\
&\lesssim(\theta_b^{1\over 2}+\theta_b^{-1})\|r\grad\mb v\|_{\Omega,a+a',j+l-1}^2+\theta_b^{-2}\|\mb v\|_{\Omega,a+a',j+l}^2+(\theta_b^{1\over 2}+\theta_b^{-1}+\delta^{-1}\theta_b^{1\over 2})\|\mb v\|_{\Omega,a+a',j+l+1}^2\\
&\quad+\delta\theta_b^{1\over 2}\|r^2\grad^\perp\grad^\perp\cdot\mb v\|_{\Omega,a'+a,j+l-1}^2.
\end{align*}
So for $\epsilon>0$ small enough we have \eqref{Vr against Vtheta estimate}.
\end{proof}

\begin{lemma}[$\mb V^\theta$ against $\mb V^\theta$ coercivity]\label{Vtheta against Vtheta}
Let $j,l\in\N_0$ and $a,a'\in\R$ with $a\neq a'$. For $\theta_b$ small enough and $\mb v\in C^\infty_{c,\sigma}(\Omega)$ (as defined in Definition \ref{smooth bulk test functions}) we have
\begin{multline}\label{Vtheta against Vtheta estimate}
{1\over 3}\|\dr^{j+l}r^2\grad^\perp\grad^\perp\cdot\mb v\|_{\Omega,a+a'}^2-\Re\<[\P^\theta]\mb V^\theta_{a,j}(\mb v),[\P]\mb V^\theta_{a',l}(\mb v)\>_{L^2(\Omega)}\\
\lesssim\|r^2\grad^\perp\grad^\perp\cdot\mb v\|_{\Omega,a+a',j+l-1}^2+(1+\theta_b^{-1})\|r\grad\mb v\|_{\Omega,a+a',j+l+1}^2+\theta_b^{-2}\|r\grad\mb v\|_{\Omega,a+a',j+l}^2,
\end{multline}
where the projections in $[\ph]$ are optional.
\end{lemma}
\begin{proof}
When both projections are present, the second is redundant by Remark \ref{P remark}, exactly as in the proof of Lemma \ref{Vr against Vtheta}, so that at most one of the two potential terms below occurs. We have, moving the projections through the radial derivatives by Lemma \ref{projections commute with dr},
\begin{align*}
\<[\P^\theta]_1\mb V^\theta_{a,j}(\mb v),[\P]_2\mb V^\theta_{a',l}(\mb w)\>_{L^2(\Omega)}
&=\inner{(-1)^j\dr^{2j}[\P^\theta]_1\grad^\perp(r^{2-2a}\grad^\perp\cdot\mb v),(-1)^l\dr^{2l}[\P]_2\grad^\perp(r^{2-2a'}\grad^\perp\cdot\mb w)}_{L^2(\Omega)}\\
&=\inner{\dr^{l}_2\dr^{j}[\P^\theta]_1\grad^\perp(r^{2-2a}\grad^\perp\cdot\mb v),\dr^{j}_2\dr^{l}[\P]_2\grad^\perp(r^{2-2a'}\grad^\perp\cdot\mb w)}_{L^2(\Omega)}\\
&=\inner{\dr^{l}_2\dr^{j}\grad^\perp(r^{2-2a}\grad^\perp\cdot\mb v),\dr^{j}_2\dr^{l}\grad^\perp(r^{2-2a'}\grad^\perp\cdot\mb w)}_{L^2(\Omega)}\\
&\quad-\sbrac{\inner{\dr^{l}_2\dr^{j}\grad q^\theta_{a}[\mb v],\dr^{j}_2\dr^{l}\grad^\perp(r^{2-2a'}\grad^\perp\cdot\mb w)}_{L^2(\Omega)}}_1\\
&\quad-\sbrac{\inner{\dr^{l}_2\dr^{j}\grad^\perp(r^{2-2a}\grad^\perp\cdot\mb v),\dr^{j}_2\dr^{l}\grad q_{a'}[\mb w]}_{L^2(\Omega)}}_2.
\end{align*}
So
\begin{align*}
\|\dr^{j+l}r^2\grad^\perp\grad^\perp\cdot\mb v\|_{\Omega,a+a'}^2&-\Re\<[\P^\theta]_1\mb V^\theta_{a,j}(\mb v),[\P]_2\mb V^\theta_{a',l}(\mb v)\>_{L^2(\Omega)}\\
&\lesssim\|r^2\grad^\perp\grad^\perp\cdot\mb v\|_{\Omega,a+a',j+l-1}\|r^2\grad^\perp\grad^\perp\cdot\mb v\|_{\Omega,a+a',j+l}\\
&\quad+\brac{\|\grad\mb v\|_{\Omega,a+a'-1,j+l}+\|\grad q^\theta_{a}[\mb v]\|_{\Omega,a'-a,j+l}+\|\grad q_{a'}[\mb v]\|_{\Omega,a-a',j+l}}\\
&\qquad\brac{\|\grad\mb v\|_{\Omega,a+a'-1,j+l}+\|r\grad^\perp\grad^\perp\cdot\mb v\|_{\Omega,a+a'-1,j+l}}
\end{align*}
Using Young's inequality we can absorb the top part of $\|r^2\grad^\perp\grad^\perp\cdot\mb v\|_{\Omega,a+a',j+l}$ to the LHS to get
\begin{align*}
&{1\over 2}\|\dr^{j+l}r^2\grad^\perp\grad^\perp\cdot\mb v\|_{\Omega,a+a'}^2-\Re\<[\P^\theta]_1\mb V^\theta_{a,j}(\mb v),[\P]_2\mb V^\theta_{a',l}(\mb v)\>_{L^2(\Omega)}\\
&\lesssim\|r^2\grad^\perp\grad^\perp\cdot\mb v\|_{\Omega,a+a',j+l-1}^2+\|r\grad\mb v\|_{\Omega,a+a',j+l}^2+\|\grad q^\theta_{a}[\mb v]\|_{\Omega,a'-a,j+l}^2+\|\grad q_{a'}[\mb v]\|_{\Omega,a-a',j+l}^2\\
&\lesssim\|r^2\grad^\perp\grad^\perp\cdot\mb v\|_{\Omega,a+a',j+l-1}^2+(1+\theta_b^{-1}+\delta^{-1})\|r\grad\mb v\|_{\Omega,a+a',j+l+1}^2+\theta_b^{-2}\|r\grad\mb v\|_{\Omega,a+a',j+l}^2\\
&\quad+\delta\|r^2\grad^\perp\grad^\perp\cdot\mb v\|_{\Omega,a+a',j+l}^2,
\end{align*}
where we used Lemma \ref{test problem pressure estimate} to estimate the $q$ norms. That lemma is applied with $a'-1$ in place of $a'$ to $q^\theta_a$, and with $a'$ and $a-1$ in place of $a$ and $a'$ to $q_{a'}$. Each time its condition $1+a'-a\neq 0$ becomes the hypothesis $a\neq a'$. Using Young's inequality we can again absorb the top part of $\|r^2\grad^\perp\grad^\perp\cdot\mb v\|_{\Omega,a+a',j+l}$ to the LHS to get \eqref{Vtheta against Vtheta estimate}.
\end{proof}

The next proposition combines the three lemmas.

\begin{proposition}[The auxiliary variational problem]\label{auxiliary variational problem}
Let $\tilde{\mb v}\in C^\infty_c(\Omega)$ and
\begin{align*}
H^{(1),l}&:=\overline{C^\infty_{c,\sigma}(\Omega)}^{\|\ph\|_{H^{(1),l}}}\\
\|\mb v\|_{H^{(1),l}}^2&:=\sum_{k=0}^1\|\mb v\|_{\Omega,1-\kappa k,2l+1,2}^2,
\end{align*}
where $C^\infty_{c,\sigma}$ is as defined in Definition \ref{smooth bulk test functions}. Define $B^{(1),l}:H^{(1),l}\times H^{(1),l}\to\C$ and $\mathcal F^{(1),l}_{\tilde{\mb v}}:H^{(1),l}\to\C$ by
\begin{align*}
B^{(1),l}[\mb v,\mb w]
&:=\<\P\mb V^{(1)}_{(l)}(\mb v),\mb V^{(2)}_{(l)}(\mb w)\>_{L^2(\Omega)}\\
\mathcal F^{(1),l}_{\tilde{\mb v}}[\mb w]
&:=\<\P\tilde{\mb v},\mb V^{(2)}_{(l)}(\mb w)\>_{L^2(\Omega)},
\end{align*}
where $\smash{\mb V^{(1)}_{(l)}}$ is as defined in Definition \ref{energy forms and functionals} and
\begin{align*}
\mb V^{(2)}_{(l)}(\mb w)
:=\sum_{j=0}^{l+2}\epsilon^j\mb V^r_{1,j}(\mb w)+\sum_{j=0}^{l+1}\epsilon^j\epsilon^{2}\mb V^\theta_{1,j}(\mb w).
\end{align*}
Then for small enough $\theta_b$ and $c_\epsilon$, where $\epsilon=c_\epsilon\theta_b$ as in Proposition \ref{energy form coercivity}, we have
\begin{align*}
\|\mb v\|_{H^{(1),l}}^2&\lesssim|B^{(1),l}[\mb v,\mb v]|\\
|B^{(1),l}[\mb v,\mb w]|&\lesssim\|\mb v\|_{H^{(1),l}}\|\mb w\|_{H^{(1),l}}\\
|\mathcal F^{(1),l}_{\tilde{\mb v}}[\mb w]|&\lesssim\|\mb w\|_{H^{(1),l}}.
\end{align*}
\end{proposition}
\begin{proof}
By Lemma \ref{Vr against Vr}, \ref{Vr against Vtheta} and \ref{Vtheta against Vtheta} we have
\begin{align*}
&\sum_{k=0}^1\epsilon^{5k}\brac{1+\delta_{0k}|s|^{-2\kappa}}\Bigg(\sum_{j=0}^{l+1}\sum_{j'=0}^{l+2}\epsilon^{j+j'}{1\over 3}\|\dr^{j+j'}\mb v\|_{\Omega,1-k\kappa}^2\\
&\quad+\Bigg(\sum_{j=0}^l\sum_{j'=0}^{l+2}+\sum_{j=0}^{l+1}\sum_{j'=0}^{l+1}\Bigg)\epsilon^{j+j'+2}{2\over C}\|r\grad(\dr^{j+j'}\mb v)\|_{\Omega,1-k\kappa}^2\\
&\quad+\epsilon^4\sum_{j=0}^l\sum_{j'=0}^{l+1}\epsilon^{j+j'}{1\over 3}\|\dr^{j+j'}r^2\grad^\perp\grad^\perp\cdot\mb v\|_{\Omega,1-k\kappa}^2\Bigg)\\
&\lesssim\Re B^{(1),l}[\mb v]+\sum_{k=0}^1\epsilon^{5k}\brac{1+\delta_{0k}|s|^{-2\kappa}}\Bigg(\sum_{j=0}^{l+1}\sum_{j'=0}^{l+2}\epsilon^{j+j'}\|\mb v\|_{\Omega,1-k\kappa,j+j'-1}^2\\
&\quad+\epsilon^2\Bigg(\sum_{j=0}^l\sum_{j'=0}^{l+2}+\sum_{j=0}^{l+1}\sum_{j'=0}^{l+1}\Bigg)\epsilon^{j+j'}\brac{(\theta_b^{1\over 2}+\theta_b^{-1})\|r\grad\mb v\|_{\Omega,1-k\kappa,j+j'-1}^2+\theta_b^{-2}\|\mb v\|_{\Omega,1-k\kappa,j+j'}^2}\\
&\quad+\epsilon^2\sum_{j=0}^l\epsilon^j\brac{\sum_{j'=-1}^{l+2}\epsilon^{j'}(\theta_b^{1\over 2}+\theta_b^{-1}+\delta^{-1}\theta_b^{1\over 2})\|\mb v\|_{\Omega,1-k\kappa,j+j'+1}^2+\delta\theta_b^{1\over 2}\sum_{j'=1}^{l+2}\epsilon^{j'}\|r^2\grad^\perp\grad^\perp\cdot\mb v\|_{\Omega,1-k\kappa,j+j'-1}^2}\\
&\quad+\epsilon^2\sum_{j=0}^{l+1}\epsilon^j\brac{\sum_{j'=-1}^{l+1}\epsilon^{j'}(\theta_b^{1\over 2}+\theta_b^{-1}+\delta^{-1}\theta_b^{1\over 2})\|\mb v\|_{\Omega,1-k\kappa,j+j'+1}^2+\delta\theta_b^{1\over 2}\sum_{j'=1}^{l+1}\epsilon^{j'}\|r^2\grad^\perp\grad^\perp\cdot\mb v\|_{\Omega,1-k\kappa,j+j'-1}^2}\\
&\quad+\epsilon^4\sum_{j=0}^l\sum_{j'=0}^{l+1}\epsilon^{j+j'}\brac{\|r^2\grad^\perp\grad^\perp\cdot\mb v\|_{\Omega,1-k\kappa,j+j'-1}^2+(1+\theta_b^{-1})\|r\grad\mb v\|_{\Omega,1-k\kappa,j+j'+1}^2+\theta_b^{-2}\|r\grad\mb v\|_{\Omega,1-k\kappa,j+j'}^2}\Bigg).
\end{align*}
Choose $\delta=\epsilon$. Then for $\theta_b$ and $c_\epsilon$ small enough, all the terms after $\Re B^{(1),l}[\mb v]$ can be absorbed into the LHS, we can also get the full second derivative on the LHS by Lemma \ref{second derivative bound}, so we get the coercivity estimate
\begin{align*}
\epsilon^{2l+10}\|\mb v\|_{H^{(1),l}}^2\lesssim\Re B^{(1),l}[\mb v]\leq|B^{(1),l}[\mb v]|,
\end{align*}
the power of $\epsilon$ being the largest one occurring on the left, namely $\epsilon^{5}\epsilon^{4}\epsilon^{2l+1}$.
Boundedness $|B^{(1),l}[\mb v,\mb w]|\lesssim\|\mb v\|_{H^{(1),l}}\|\mb w\|_{H^{(1),l}}$ is the same term by term inspection as the boundedness of $B_{(l)}$ in the proof of Proposition \ref{energy form coercivity}, the remainder terms being bounded by the norm in the coercivity estimate above. Finally, we have boundedness of $\smash{\mathcal F^{(1),l}_{\tilde{\mb v}}}$ since by Lemma \ref{decay of projected test functions} $\P\tilde{\mb v}$ has good decay.
\end{proof}

\subsection{Solvability of the test function problem}\label{solvability of the test function problem}

It remains to show that the test function map is onto, first in the bulk and then on the boundary; this is what turns the variational identity into the equation holding almost everywhere. We refer to these statements as solvability statements: what is asserted is that every projected bump function is of the form $\P\smash[b]{\mb V^\bullet_{(l)}}(\mb v)$, and it is in this sense that the test functions are dense. We proceed in two steps. Since it only has one scaling overall, the test function $\smash[t]{\trimbot{0.6pt}{\mb V^{(2)}_{(l)}}}(\mb v)$ can be solved for directly in Mellin variables, and the coercivity of Proposition \ref{auxiliary variational problem} then transfers this solvability to $\smash{\mb V^{(1)}_{(l)}}(\mb v)$.

\begin{proposition}[Solvability for $\mb V^{(2)}$]\label{solvability for V2}
For small enough $\theta_b$ and $c_\epsilon$ we have the following. For any $\tilde{\mb w}\in C^\infty_c(\Omega)$, there exists $\mb w\in H^{(1),l}$ such that $\P\smash{\mb V^{(2)}_{(l)}}(\mb w)=\P\tilde{\mb w}$, i.e.
\begin{align*}
\mb V^{(2)}_{(l)}(\mb w)-\grad q[\mb V^{(2)}_{(l)}(\mb w)]=\P\tilde{\mb w},
\end{align*}
where $q[\smash{\mb V^{(2)}_{(l)}}(\mb w)]$ is defined in Definition \ref{the two projections}.
\end{proposition}
\begin{proof}
Note that $r\grad^\perp\cdot\mb w=\dr_1w_\theta-\partial_\theta w_r$. We have
\begin{align*}
&r^{2}\brac{\mb V^{(2)}_{(l)}(\mb w)-\grad q[\mb V^{(2)}_{(l)}(\mb w)]}\\
&=-\sum_{j=0}^{l+2}\epsilon^j(-1)^{j}\dr^{2j}\mb w+\epsilon^2\sum_{j=0}^{l+1}\epsilon^jr^{2}(-1)^j\dr^{2j}\grad^\perp\grad^\perp\cdot\mb w-r^{2}\grad q[\mb V^{(2)}_{(l)}(\mb w)]\\
&=\brac{-\sum_{j=0}^{l+2}\epsilon^j(-1)^{j}\dr^{2j}w_r-\epsilon^2\sum_{j=0}^{l+1}\epsilon^jr^{2}(-1)^j\dr^{2j}r^{-1}\partial_\theta(r^{-1}(\dr_1w_\theta-\partial_\theta w_r))-r\dr q[\mb V^{(2)}_{(l)}(\mb w)]}\mb e_r\\
&\quad+\brac{-\sum_{j=0}^{l+2}\epsilon^j(-1)^{j}\dr^{2j}w_\theta+\epsilon^2\sum_{j=0}^{l+1}\epsilon^jr^{2}(-1)^j\dr^{2j}r^{-1}\dr(r^{-1}(\dr_1w_\theta-\partial_\theta w_r))-r\partial_\theta q[\mb V^{(2)}_{(l)}(\mb w)]}\mb e_\theta\\
&=\brac{-\sum_{j=0}^{l+2}\epsilon^j(-1)^{j}\dr^{2j}w_r-\epsilon^2\sum_{j=0}^{l+1}\epsilon^j(-1)^j\dr^{2j}_{-2}\partial_\theta(\dr_1w_\theta-\partial_\theta w_r)-r\dr q[\mb V^{(2)}_{(l)}(\mb w)]}\mb e_r\\
&\quad+\brac{-\sum_{j=0}^{l+2}\epsilon^j(-1)^{j}\dr^{2j}w_\theta+\epsilon^2\sum_{j=0}^{l+1}\epsilon^j(-1)^j\dr^{2j}_{-2}\dr_{-1}(\dr_1w_\theta-\partial_\theta w_r)-r\partial_\theta q[\mb V^{(2)}_{(l)}(\mb w)]}\mb e_\theta
\end{align*}
So in Mellin variables we want to solve
\begin{align*}
-\sum_{j=0}^{l+2}\epsilon^j(-1)^{j}\lambda^{2j}\hat w_r-\epsilon^2\sum_{j=0}^{l+1}\epsilon^j(-1)^j\lambda_{-2}^{2j}\partial_\theta(\lambda_1\hat w_\theta-\partial_\theta\hat w_r)
-\lambda_{-1}\hat q[\mb V^{(2)}_{(l)}(\mb w)](\lambda_{-1},\theta)
&=(\P\tilde{\mb w})_r^\wedge(\lambda_{-2},\theta)\\
-\sum_{j=0}^{l+2}\epsilon^j(-1)^{j}\lambda^{2j}\hat w_\theta+\epsilon^2\sum_{j=0}^{l+1}\epsilon^j(-1)^j\lambda_{-2}^{2j}\lambda_{-1}(\lambda_1\hat w_\theta-\partial_\theta\hat w_r)
-\partial_\theta\hat q[\mb V^{(2)}_{(l)}(\mb w)](\lambda_{-1},\theta)
&=(\P\tilde{\mb w})_\theta^\wedge(\lambda_{-2},\theta)
\end{align*}
Eliminating the $q$ term we have
\begin{align*}
\lambda_{-1}(\P\tilde{\mb w})_\theta^\wedge(\lambda_{-2},\theta)&-\partial_\theta(\P\tilde{\mb w})_r^\wedge(\lambda_{-2},\theta)\\
&=\sum_{j=0}^{l+2}\epsilon^j(-1)^{j}\lambda^{2j}\partial_\theta\hat w_r+\epsilon^2\sum_{j=0}^{l+1}\epsilon^j(-1)^j\lambda_{-2}^{2j}\partial_\theta^2(\lambda_1\hat w_\theta-\partial_\theta\hat w_r)\\
&\quad-\sum_{j=0}^{l+2}\epsilon^j(-1)^{j}\lambda^{2j}\lambda_{-1}\hat w_\theta+\epsilon^2\sum_{j=0}^{l+1}\epsilon^j(-1)^j\lambda_{-2}^{2j}\lambda_{-1}^2(\lambda_1\hat w_\theta-\partial_\theta\hat w_r)\\
&=-\sum_{j=0}^{l+2}\epsilon^j(-1)^{j}\lambda^{2j}(\lambda_{-1}\hat w_\theta-\partial_\theta\hat w_r)+\epsilon^2\sum_{j=0}^{l+1}\epsilon^j(-1)^j\lambda_{-2}^{2j}(\lambda_{-1}^2+\partial_\theta^2)(\lambda_1\hat w_\theta-\partial_\theta\hat w_r)
\end{align*}
The divergence constraint is $\lambda_1\hat w_r(\lambda,\theta)+\partial_\theta\hat w_\theta(\lambda,\theta)=0$. So eliminating $\hat w_r$ we get the following 4th order ODE in $\theta$:
\begin{align}
&\lambda_1\lambda_{-1}(\P\tilde{\mb w})_\theta^\wedge(\lambda_{-2},\theta)-\lambda_1\partial_\theta(\P\tilde{\mb w})_r^\wedge(\lambda_{-2},\theta)\nonumber\\
&=-\sum_{j=0}^{l+2}\epsilon^j(-1)^{j}\lambda^{2j}(\lambda_1\lambda_{-1}+\partial_\theta^2)\hat w_\theta+\epsilon^2\sum_{j=0}^{l+1}\epsilon^j(-1)^j\lambda_{-2}^{2j}(\lambda_{-1}^2+\partial_\theta^2)(\lambda_1^2+\partial_\theta^2)\hat w_\theta\nonumber\\
&=-\sum_{j=0}^{l+2}\epsilon^j(-1)^{j}\lambda^{2j}(\lambda_1\lambda_{-1}+\partial_\theta^2)\hat w_\theta+\epsilon^2\sum_{j=0}^{l+1}\epsilon^j(-1)^j\lambda_{-2}^{2j}\brac{\lambda_{-1}^2\lambda_1^2+(\lambda_{-1}^2+\lambda_1^2)\partial_\theta^2+\partial_\theta^4}\hat w_\theta\nonumber\\
&=\underbrace{\brac{\epsilon^2\sum_{j=0}^{l+1}\epsilon^j(-1)^j\lambda_{-2}^{2j}}}_{\mathfrak c_0(\lambda)}\partial_\theta^4\hat w_\theta-\underbrace{\brac{\sum_{j=0}^{l+2}\epsilon^j(-1)^{j}\lambda^{2j}-\epsilon^2\sum_{j=0}^{l+1}\epsilon^j(-1)^j\lambda_{-2}^{2j}(\lambda_{-1}^2+\lambda_1^2)}}_{\mathfrak c_1(\lambda)}\partial_\theta^2\hat w_\theta\nonumber\\
&\quad+\underbrace{\brac{-\sum_{j=0}^{l+2}\epsilon^j(-1)^{j}\lambda^{2j}\lambda_1\lambda_{-1}+\epsilon^2\sum_{j=0}^{l+1}\epsilon^j(-1)^j\lambda_{-2}^{2j}\lambda_{-1}^2\lambda_1^2}}_{\mathfrak c_2(\lambda)}\hat w_\theta\label{test function equation in mellin}
\end{align}
For any $\alpha_0>-1$, there exists $\epsilon$ small enough such that, for every $\lambda$ with $\Re\lambda\in[-1,\alpha_0]$,
\begin{align*}
\Re \mathfrak c_0(\lambda)&\geq{1\over 2}\epsilon^2\sum_{j=0}^{l+1}\epsilon^j|\lambda|^{2j}\\
\Re \mathfrak c_1(\lambda)&\geq{1\over 2}\sum_{j=0}^{l+2}\epsilon^j|\lambda|^{2j}\\
\Re \mathfrak c_2(\lambda)&\geq-C_{\alpha_0}+{1\over 2}\sum_{j=0}^{l+2}\epsilon^j|\lambda|^{2+2j}
\end{align*}
We have boundary condition $\mb w|_{\partial_0\Omega}=w_r|_{\partial_1\Omega}=\mb 0$, or equivalently
\begin{align*}
\hat w_\theta(\lambda,0)=\partial_\theta\hat w_\theta(\lambda,0)=\partial_\theta\hat w_\theta(\lambda,\theta_b)=0.
\end{align*}

If we want $\smash{\mb V^{(2)}_{(l)}(\mb w)-\grad q=\P\tilde{\mb w}}$ with $q|_{\partial_1\Omega}=0$, since $(\mb e_r\cdot\P\tilde{\mb w})|_{\partial_1\Omega}=0$, we also need $\smash{(\mb e_r\cdot\mb V^{(2)}_{(l)}(\mb w))|_{\partial_1\Omega}=0}$. That means we need $0=\partial_\theta(\lambda_1\hat w_\theta-\partial_\theta\hat w_r)=\lambda_1^{-1}\partial_\theta^3\hat w_\theta$. So we also impose the boundary condition $\partial_\theta^3\hat w_\theta(\lambda,\theta_b)=0$.

Conversely, given a solution $\hat w_\theta$ of this boundary value problem, define $\hat w_r$ by the divergence constraint and let $\mb D:=\smash[b]{\mb V^{(2)}_{(l)}}(\mb w)-\P\tilde{\mb w}$. The equation \eqref{test function equation in mellin} says that $\grad^\perp\cdot\mb D=0$, so that $\mb D=\grad q$ on the simply connected wedge; $D_r$ vanishes on $\partial_1\Omega$, so $q$ is constant there and may be normalised to vanish; and $D_\theta=(\smash[t]{\trimbot{0.5pt}{\mb V^{(2)}_{(l)}}}(\mb w))_\theta$ on $\partial_0\Omega$ and $\grad\cdot\mb D=\grad\cdot\smash[t]{\trimbot{0.5pt}{\mb V^{(2)}_{(l)}}}(\mb w)$ by Remark \ref{P remark}. Thus $q$ solves the mixed problem that defines $q[\smash[t]{\trimbot{0.5pt}{\mb V^{(2)}_{(l)}}}(\mb w)]$ in Definition \ref{the two projections}, and coincides with it, the difference being harmonic with vanishing data and hence zero by its Mellin representation; this is the identity $\P\smash[t]{\trimbot{2.0pt}{\mb V^{(2)}_{(l)}}}(\mb w)=\P\tilde{\mb w}$.

Multiplying \eqref{test function equation in mellin} by $\overline{\hat w_\theta}$, for any solution $\hat w_\theta$ of this boundary value problem, integrating in $\theta$ and then integrating by parts, we get
\begin{align*}
&\<\hat w_\theta(\lambda),(\lambda_1\lambda_{-1}(\P\tilde{\mb w})_\theta^\wedge-\lambda_1\partial_\theta(\P\tilde{\mb w})_r^\wedge)(\lambda_{-2})\>_{L^2_\theta([0,\theta_b])}\\
&=\mathfrak c_0(\lambda)\|\partial_\theta^2\hat w_\theta\|_{L^2_\theta([0,\theta_b])}^2+\mathfrak c_1(\lambda)\|\partial_\theta\hat w_\theta\|_{L^2_\theta([0,\theta_b])}^2+\mathfrak c_2(\lambda)\|\hat w_\theta\|_{L^2_\theta([0,\theta_b])}^2.
\end{align*}
Taking the real part we get, for $\Re\lambda\in[-1,\alpha_0]$,
\begin{align*}
&\epsilon^2\brac{\sum_{j=0}^{l+1}\epsilon^j|\lambda|^{2j}}\|\partial_\theta^2\hat w_\theta\|_{L^2_\theta}^2+\brac{\sum_{j=0}^{l+2}\epsilon^j|\lambda|^{2j}}\|\partial_\theta\hat w_\theta\|_{L^2_\theta}^2+\brac{-2C_{\alpha_0}+\sum_{j=0}^{l+2}\epsilon^j|\lambda|^{2+2j}}\|\hat w_\theta\|_{L^2_\theta}^2\\
&\leq2\Re\<\hat w_\theta(\lambda),(\lambda_1\lambda_{-1}(\P\tilde{\mb w})_\theta^\wedge-\lambda_1\partial_\theta(\P\tilde{\mb w})_r^\wedge)(\lambda_{-2})\>_{L^2_\theta}
\end{align*}
Using Jensen's inequality we have
\begin{align*}
\|\hat w_\theta\|_{L^2_\theta}^2&=\int_0^{\theta_b}|\hat w_\theta(\lambda,\theta)|^2\d\theta=\int_0^{\theta_b}\theta^2\abs{{1\over\theta}\int_0^\theta\partial_\theta\hat w_\theta(\lambda,\theta')\d\theta'}^2\d\theta
\leq\int_0^{\theta_b}\theta\int_0^\theta|\partial_\theta\hat w_\theta(\lambda,\theta')|^2\d\theta'\d\theta\\
&\leq{1\over 2}\theta_b^2\|\partial_\theta\hat w_\theta\|_{L^2_\theta}^2
\end{align*}
So for small enough $\theta_b$, we have
\begin{align*}
&\epsilon^2\brac{\sum_{j=0}^{l+1}\epsilon^j|\lambda|^{2j}}\|\partial_\theta^2\hat w_\theta\|_{L^2_\theta}^2+{1\over 2}\brac{\sum_{j=0}^{l+2}\epsilon^j|\lambda|^{2j}}\|\partial_\theta\hat w_\theta\|_{L^2_\theta}^2+\brac{1+\sum_{j=0}^{l+2}\epsilon^j|\lambda|^{2+2j}}\|\hat w_\theta\|_{L^2_\theta}^2\\
&\leq2\Re\<\hat w_\theta(\lambda),(\lambda_1\lambda_{-1}(\P\tilde{\mb w})_\theta^\wedge-\lambda_1\partial_\theta(\P\tilde{\mb w})_r^\wedge)(\lambda_{-2})\>_{L^2_\theta}
\end{align*}
Using Young's inequality, in which the coefficient $1+\sum_{j=0}^{l+2}\epsilon^j|\lambda|^{2+2j}\gtrsim_\epsilon(1+|\lambda|^2)^{l+3}$ of $\|\hat w_\theta\|_{L^2_\theta}^2$ on the left is kept, we get
\begin{align*}
&\epsilon^2\brac{\sum_{j=0}^{l+1}\epsilon^j|\lambda|^{2j}}\|\partial_\theta^2\hat w_\theta\|_{L^2_\theta}^2+{1\over 2}\brac{\sum_{j=0}^{l+2}\epsilon^j|\lambda|^{2j}}\|\partial_\theta\hat w_\theta\|_{L^2_\theta}^2+{1\over 2}\brac{1+\sum_{j=0}^{l+2}\epsilon^j|\lambda|^{2+2j}}\|\hat w_\theta\|_{L^2_\theta}^2\\
&\lesssim_\epsilon(1+|\lambda|^2)^{-l-2}\brac{(1+|\lambda|^2)\|(\P\tilde{\mb w})_\theta^\wedge(\lambda_{-2})\|_{L^2_\theta}^2+\|\partial_\theta(\P\tilde{\mb w})_r^\wedge(\lambda_{-2})\|_{L^2_\theta}^2}
\end{align*}
In particular the homogeneous problem has only the trivial solution, so that, the solution space of the fourth order equation being four-dimensional, the boundary value problem has a unique solution $\hat w_\theta$ for every RHS; and since its coefficients and its datum depend holomorphically on $\lambda$, so does $\hat w_\theta$ on the strip $-1<\Re\lambda\leq\alpha_0$, so that the inverse Mellin transforms taken on the lines $\Re\lambda=\alpha$ below all define one and the same field $\mb w$. Now we have control of $\hat w_\theta$ up to two $\partial_\theta$ derivatives. To get control of $\|\partial_\theta^4\hat w_\theta\|_{L^2_\theta}$: in equation \eqref{test function equation in mellin}, we move the last two terms to the left, take the square of the modulus, and then integrate in $\theta$ to get
\begin{align*}
\epsilon^4\brac{\sum_{j=0}^{l+1}\epsilon^j|\lambda|^{2j}}^2\|\partial_\theta^4\hat w_\theta\|_{L^2_\theta}^2&\lesssim\|\lambda_1\lambda_{-1}(\P\tilde{\mb w})_\theta^\wedge(\lambda_{-2},\theta)-\lambda_1\partial_\theta(\P\tilde{\mb w})_r^\wedge(\lambda_{-2},\theta)\|_{L^2_\theta}^2\\
&\quad+\brac{\sum_{j=0}^{2(l+2)}\epsilon^j|\lambda|^{2j}}\|\partial_\theta^2\hat w_\theta\|_{L^2_\theta}^2+\brac{1+\sum_{j=1}^{2(l+3)}\epsilon^{j-1}|\lambda|^{2j}}\|\hat w_\theta\|_{L^2_\theta}^2\\
&\lesssim_\epsilon(1+|\lambda|^2)\brac{(1+|\lambda|^2)\|(\P\tilde{\mb w})_\theta^\wedge(\lambda_{-2})\|_{L^2_\theta}^2+\|\partial_\theta(\P\tilde{\mb w})_r^\wedge(\lambda_{-2})\|_{L^2_\theta}^2}
\end{align*}
And to control $\|\partial_\theta^3\hat w_\theta\|_{L^2_\theta}$: the Gagliardo–Nirenberg interpolation inequality gives
\begin{align*}
\brac{\sum_{j=0}^{2l+3}\epsilon^j|\lambda|^{2j}}\|\partial_\theta^3\hat w_\theta\|_{L^2_\theta}^2&\lesssim\brac{\sum_{j=0}^{2l+3}\epsilon^j|\lambda|^{2j}}\brac{\|\partial_\theta^4\hat w_\theta\|_{L^2_\theta}\|\partial_\theta^2\hat w_\theta\|_{L^2_\theta}+\|\partial_\theta^2\hat w_\theta\|_{L^2_\theta}^2}\\
&\lesssim_\epsilon(1+|\lambda|^2)\brac{(1+|\lambda|^2)\|(\P\tilde{\mb w})_\theta^\wedge(\lambda_{-2})\|_{L^2_\theta}^2+\|\partial_\theta(\P\tilde{\mb w})_r^\wedge(\lambda_{-2})\|_{L^2_\theta}^2}
\end{align*}
Therefore we have
\begin{align*}
\sum_{j=0}^{2(l+1)}\sum_{n+m\leq 4}\|\lambda^{n+j}\partial_\theta^m\hat w_\theta\|_{L^2_\theta}^2\lesssim_\epsilon\sum_{j=0}^{2(l+1)}\sum_{n+m\leq 1}\|\lambda^{n+j}\partial_\theta^m(\P\tilde{\mb w})^\wedge(\lambda_{-2})\|_{L^2_\theta}^2
\end{align*}
Next we want to recover $\hat w_r$ from $\hat w_\theta$ using the divergence constraint. We have
\begin{align*}
1+|\lambda|^2\sim_\delta|\lambda+1|^2\qquad\text{for}\qquad \Re\lambda\in\R\setminus(-1-\delta,-1+\delta)
\end{align*}
so the divergence constraint $(\lambda+1)\hat w_r(\lambda,\theta)+\partial_\theta\hat w_\theta(\lambda,\theta)=0$, that is $\hat w_r=-(\lambda+1)^{-1}\partial_\theta\hat w_\theta$, gives
\begin{align*}
\sum_{j=0}^{2l+3}\sum_{n+m\leq 3}\|\lambda^{n+j}\partial_\theta^m\hat w_r\|_{L^2_\theta}^2&\lesssim_{\delta,\epsilon}\sum_{j=0}^{2l+2}\sum_{n+m\leq 1}\|\lambda^{n+j}\partial_\theta^m(\P\tilde{\mb w})^\wedge(\lambda_{-2})\|_{L^2_\theta}^2\\
&\qquad\text{for}\qquad \Re\lambda\in\R\setminus(-1-\delta,-1+\delta),
\end{align*}
each summand of the LHS being $\|\lambda^{n+j}(\lambda+1)^{-1}\partial_\theta^{m+1}\hat w_\theta\|_{L^2_\theta}$, bounded through the displayed equivalence by the summands of the previous estimate at $(n',m')$ with $n'+m'\leq 4$ and $j'\leq 2(l+1)$: the factor $(\lambda+1)^{-1}$ either absorbs one power of $\lambda$, moving $(n,j)$ with $n+j\geq 1$ to a pair with total $n'+j'=n+j-1$, or is bounded outright, and $\partial_\theta^{m+1}$ raises $m$ by one.
Now we know the full $\mb w$ and its regularity: for any $\alpha\in(-1,\alpha_0]$, we have
\begin{align*}
\sum_{j=0}^{2l+3}\sum_{n+m\leq 3}\|r^{-1-\alpha}\dr^{n+j}\partial_\theta^m\mb w\|_{L^2(\Omega)}^2
&\lesssim\int_{\Re\lambda=\alpha}\sum_{j=0}^{2l+3}\sum_{n+m\leq 3}\brac{\|\lambda^{n+j}\partial_\theta^m\hat w_\theta\|_{L^2_\theta}^2+\|\lambda^{n+j}\partial_\theta^m\hat w_r\|_{L^2_\theta}^2}\d\Im\lambda\\
&\lesssim \int_{\Re\lambda=\alpha}\sum_{j=0}^{2l+2}\sum_{n+m\leq 1}\|\lambda^{n+j}\partial_\theta^m(\P\tilde{\mb w})^\wedge(\lambda_{-2})\|_{L^2_\theta}^2\d\Im\lambda\\
&\lesssim \int_{\Re\lambda=\alpha-2}\sum_{j=0}^{2l+2}\sum_{n+m\leq 1}\|\lambda^{n+j}\partial_\theta^m(\P\tilde{\mb w})^\wedge(\lambda)\|_{L^2_\theta}^2\d\Im\lambda\\
&\lesssim\sum_{j=0}^{2l+2}\sum_{n+m\leq 1}\|r^{1-\alpha}\dr^{n+j}\partial_\theta^m\P\tilde{\mb w}\|_{L^2(\Omega)}^2,
\end{align*}
where the substitution $\mu=\lambda-2$ moves the line $\Re\lambda=\alpha$ to $\Re\mu=\alpha-2$, on which the bulk weight is $1+(\alpha-2)=\alpha-1$, and $|\mu+2|^{2(n+j)}\lesssim 1+|\mu|^{2(n+j)}$ is absorbed by the sum over $n+j$. Note that the norms for $\P\tilde{\mb w}$ on the RHS are all finite by Lemma \ref{decay of projected test functions}. Taking $\alpha_0=0$, we get the finiteness of $\|\mb w\|_{H^{(1),l}}$, together with $\grad\cdot\mb w=0$, $\mb w|_{\partial_0\Omega}=\mb 0$ and $w_r|_{\partial_1\Omega}=0$.

It remains to see that $\mb w$ lies in the closure $H^{(1),l}$ of $C^\infty_{c,\sigma}(\Omega)$, that is, that it can be approximated in the norm $\|\ph\|_{H^{(1),l}}$ by smooth, compactly supported, divergence-free fields with the same two boundary conditions. Write $\mb w=\grad^\perp\psi$ with a stream function normalised by $\psi=0$ on $\partial_0\Omega$, so that the two boundary conditions read $\psi=\partial_\theta\psi=0$ on $\partial_0\Omega$ and $\partial_\theta\psi=0$ on $\partial_1\Omega$, and let $\zeta$ be smooth on $[0,\infty)$ with $\zeta=0$ on $[0,1]$ and $\zeta=1$ on $[2,\infty)$, so that $\zeta_n(r):=\zeta(2^nr)(1-\zeta(2^{-n}r))$ has $\dr^k\zeta_n$ bounded uniformly in $n$ and supported in $\{2^{-n}\leq r\leq 2^{1-n}\}\cup\{2^n\leq r\leq 2^{n+1}\}$. Passing from $\psi$ to $\zeta_n\psi$ cuts off in $r$ with an error controlled by the derivatives of $\psi$, and hence, on integrating $\partial_\theta\psi=-rw_r$ from $\partial_0\Omega$, by those of $\mb w$ itself over the support of $\grad\zeta_n$. No Hardy inequality is needed, since $\|\mb w\|_{H^{(1),l}}$ contains $\|\mb w\|_{\Omega,1}$ and $\|\mb w\|_{\Omega,1-\kappa}$. The truncated stream function is then smoothed in the strip coordinates $(\ln r,\theta)\in\R\times[0,\theta_b]$, in which the weights are bounded above and below. We extend it evenly across $\theta=0$ and $\theta=\theta_b$, which preserves the two Neumann conditions; we mollify in $\ln r$ and $\theta$ with a compactly supported mollifier that is even in $\theta$; and we subtract the trace on $\theta=0$, which tends to $0$ with all its derivatives since $\psi=O(\theta^2)$ there. The rotated gradient of the result lies in $C^\infty_{c,\sigma}(\Omega)$ and converges to $\mb w$ in $\|\ph\|_{H^{(1),l}}$.
\end{proof}

Having solved $\P\smash{\mb V^{(2)}_{(l)}}(\mb w)=\P\tilde{\mb w}$ in Proposition \ref{solvability for V2}, we can now solve $\P\smash{\mb V^{(1)}_{(l)}}(\mb v)=\P\tilde{\mb v}$.

\begin{proposition}[Solvability in the bulk]\label{solvability in the bulk}
For small enough $\theta_b$ and $c_\epsilon$, given any $\tilde{\mb v}\in C^\infty_c(\Omega)$, there exists $\mb v\in H^{(1),l}$ such that $\P\smash[b]{\trimtop{0.4pt}{\mb V^{(1)}_{(l)}}}(\mb v)=\P\tilde{\mb v}$.
\end{proposition}
\begin{proof}
Applying the Lax--Milgram theorem to the setting of Proposition \ref{auxiliary variational problem}, we get that there exists $\mb v\in H^{(1),l}$ such that
\begin{align*}
0=B^{(1),l}[\mb v,\mb w]-\mathcal F^{(1),l}_{\tilde{\mb v}}[\mb w]=\<\P\mb V^{(1)}_{(l)}(\mb v)-\P\tilde{\mb v},\mb V^{(2)}_{(l)}(\mb w)\>_{L^2(\Omega)}\qquad\text{for all}\qquad\mb w\in H^{(1),l}.
\end{align*}
Proposition \ref{solvability for V2} tells us that, for any given $\tilde{\mb w}\in C^\infty_c(\Omega)$, there exists $\mb w\in H^{(1),l}$ such that $\P\smash{\mb V^{(2)}_{(l)}}(\mb w)=\P\tilde{\mb w}$. So we have
\begin{align*}
0=\<\P\mb V^{(1)}_{(l)}(\mb v)-\P\tilde{\mb v},\P\tilde{\mb w}\>_{L^2(\Omega)}=\<\P\mb V^{(1)}_{(l)}(\mb v)-\P\tilde{\mb v},\tilde{\mb w}\>_{L^2(\Omega)}\qquad\text{for all}\qquad\tilde{\mb w}\in C^\infty_c(\Omega).
\end{align*}
This means that $\P\smash{\mb V^{(1)}_{(l)}}(\mb v)=\P\tilde{\mb v}$.
\end{proof}

\begin{lemma}[From $H^{(1),l}$ to $H^l_s$]\label{test function upgrade}
Let $\tilde{\mb v}\in C^\infty_c(\Omega)$. Suppose we have $\mb v\in H^{(1),l}$ such that $\P\smash{\mb V^{(1)}_{(l)}}(\mb v)=\P\tilde{\mb v}$. Then for small enough $\theta_b$ and $c_\epsilon$ we have in fact $\mb v\in H_s^l$.
\end{lemma}
\begin{proof}
$\P\smash{\mb V^{(1)}_{(l)}}(\mb v)=\P\tilde{\mb v}$ a.e. gives us
\begin{align*}
\brac{1+|s|^{-2\kappa}}\P\brac{\mb V^r_{0,(l)}(\mb v)+\mb V^\theta_{0,(l)}(\mb v)}&=\P\tilde{\mb v}-\epsilon^5\P\brac{\mb V^r_{-\kappa,(l)}(\mb v)+\mb V^\theta_{-\kappa,(l)}(\mb v)}\qquad\text{ on }\Omega.
\end{align*}
Let
\begin{align*}
\mb V^{\dagger}_{a,l}(\mb v):=\sum_{j=0}^{l+2}\epsilon^j\mb V^r_{a,j}(\mb v)+\sum_{j=0}^{l+1}\epsilon^j\epsilon^2\mb V^\theta_{a,j}(\mb v)
\end{align*}
Note that we have
\begin{align*}
B&:=\inner{\brac{1+|s|^{-2\kappa}}\P\brac{\mb V^r_{0,(l)}(\mb v)+\mb V^\theta_{0,(l)}(\mb v)},\mb V^{\dagger}_{1+\kappa,l}(\mb v)}_{L^2(\Omega)}\\
&\:=\inner{\P\tilde{\mb v}-\epsilon^5\P\brac{\mb V^r_{-\kappa,(l)}(\mb v)+\mb V^\theta_{-\kappa,(l)}(\mb v)},\mb V^{\dagger}_{1+\kappa,l}(\mb v)}_{L^2(\Omega)}\\
&\:\lesssim_\epsilon\|\mb v\|_{\Omega,1,2l+1,2}(\|\P\tilde{\mb v}\|_{\Omega,1+2\kappa,2l+1,2}+\|\mb v\|_{\Omega,1,2l+1,2})
\lesssim\|\P\tilde{\mb v}\|_{\Omega,1+2\kappa,2l+1,2}^2+\|\mb v\|_{H^{(1),l}}^2<\infty,
\end{align*}
where $\|\P\tilde{\mb v}\|_{\Omega,1+2\kappa,2l+1,2}$ is finite for small enough $\theta_b$ by Lemma \ref{decay of projected test functions}. Now by Lemma \ref{Vr against Vr}, \ref{Vr against Vtheta} and \ref{Vtheta against Vtheta} we have
\begin{align*}
&\brac{1+|s|^{-2\kappa}}\Bigg(\sum_{j=0}^{l+1}\sum_{j'=0}^{l+2}\epsilon^{j+j'}{1\over 3}\|\dr^{j+j'}\mb v\|_{\Omega,1+\kappa}^2
+\Bigg(\sum_{j=0}^l\sum_{j'=0}^{l+2}+\sum_{j=0}^{l+1}\sum_{j'=0}^{l+1}\Bigg)\epsilon^{j+j'+2}{2\over C}\|r\grad(\dr^{j+j'}\mb v)\|_{\Omega,1+\kappa}^2\\
&\pushright{+\epsilon^4\sum_{j=0}^l\sum_{j'=0}^{l+1}\epsilon^{j+j'}{1\over 3}\|\dr^{j+j'}r^2\grad^\perp\grad^\perp\cdot\mb v\|_{\Omega,1+\kappa}^2\Bigg)}\\
&\lesssim\Re B+\brac{1+|s|^{-2\kappa}}\Bigg(\sum_{j=0}^{l+1}\sum_{j'=0}^{l+2}\epsilon^{j+j'}\|\mb v\|_{\Omega,1+\kappa,j+j'-1}^2\\
&\quad+\epsilon^2\Bigg(\sum_{j=0}^l\sum_{j'=0}^{l+2}+\sum_{j=0}^{l+1}\sum_{j'=0}^{l+1}\Bigg)\epsilon^{j+j'}\brac{(\theta_b^{1\over 2}+\theta_b^{-1})\|r\grad\mb v\|_{\Omega,1+\kappa,j+j'-1}^2+\theta_b^{-2}\|\mb v\|_{\Omega,1+\kappa,j+j'}^2}\\
&\quad+\epsilon^2\sum_{j=0}^l\epsilon^j\brac{\sum_{j'=-1}^{l+2}\epsilon^{j'}(\theta_b^{1\over 2}+\theta_b^{-1}+\delta^{-1}\theta_b^{1\over 2})\|\mb v\|_{\Omega,1+\kappa,j+j'+1}^2+\delta\theta_b^{1\over 2}\sum_{j'=1}^{l+2}\epsilon^{j'}\|r^2\grad^\perp\grad^\perp\cdot\mb v\|_{\Omega,1+\kappa,j+j'-1}^2}\\
&\quad+\epsilon^2\sum_{j=0}^{l+1}\epsilon^j\brac{\sum_{j'=-1}^{l+1}\epsilon^{j'}(\theta_b^{1\over 2}+\theta_b^{-1}+\delta^{-1}\theta_b^{1\over 2})\|\mb v\|_{\Omega,1+\kappa,j+j'+1}^2+\delta\theta_b^{1\over 2}\sum_{j'=1}^{l+1}\epsilon^{j'}\|r^2\grad^\perp\grad^\perp\cdot\mb v\|_{\Omega,1+\kappa,j+j'-1}^2}\\
&\quad+\epsilon^4\sum_{j=0}^l\sum_{j'=0}^{l+1}\epsilon^{j+j'}\brac{\|r^2\grad^\perp\grad^\perp\cdot\mb v\|_{\Omega,1+\kappa,j+j'-1}^2+(1+\theta_b^{-1})\|r\grad\mb v\|_{\Omega,1+\kappa,j+j'+1}^2+\theta_b^{-2}\|r\grad\mb v\|_{\Omega,1+\kappa,j+j'}^2}\Bigg).
\end{align*}
Choose $\delta=\epsilon$. Then for $\theta_b$ and $c_\epsilon$ small enough, all the terms after $\Re B$ are absorbed into the LHS, and Lemma \ref{second derivative bound} supplies the full second derivative, so that
\begin{align}\label{second test function a priori}
\|\mb v\|_{\Omega,1+\kappa,2l+1,2}^2\lesssim_\epsilon\|\P\tilde{\mb v}\|_{\Omega,1+2\kappa,2l+1,2}^2+\|\mb v\|_{\Omega,1,2l+1,2}^2.
\end{align}

This is an a priori estimate: the three lemmas are stated for $C^\infty_{c,\sigma}(\Omega)$, and the absorption presupposes that the LHS is finite, which is not known for $\mb v$, whose norms are those of $H^{(1),l}$. We therefore run the chain on the approximants $\mb w_{n,\sigma}:=\grad^\perp(\zeta(2^nr)\psi_\sigma)$. Here $\zeta$ is the cut-off of the proof of Proposition \ref{solvability for V2}, and $\psi_\sigma$ is the stream function of the radial smoothing $S_\sigma\mb v$ of Lemma \ref{radial smoothing}, normalised as in the last paragraph of that proof. The approximants are divergence-free, satisfy the two boundary conditions, vanish near the corner and have all their radial derivatives. So all their norms $\|\ph\|_{\Omega,a,j,2}$ with $a\geq 1-\kappa$ are finite, and the argument of that paragraph, run at these weights, approximates them by elements of $C^\infty_{c,\sigma}(\Omega)$ and so transfers the three lemmas to them. The approximants converge to $\mb v$ in $H^{(1),l}$, first in $n$ and then in $\sigma$, with $\|\mb w_{n,\sigma}\|_{H^{(1),l}}\lesssim\|\mb v\|_{H^{(1),l}}$ throughout. They satisfy the equation for $\mb v$ up to two errors. The first is the cut-off error $\trimtop{0.3pt}{\trimbot{1.6pt}{\mb V^{(1)}_{(l)}}}(\mb w_{n,\sigma}-S_\sigma\mb v)$, supported in $\{r\leq 2^{1-n}\}$ and vanishing as $n\to\infty$ at fixed $\sigma$. The second is the commutator of $S_\sigma$ with the factor $r^{2\kappa}$ carried by $\trimbot{2.6pt}{\mb V^r_{-\kappa,(l)}}+\trimbot{2.6pt}{\mb V^\theta_{-\kappa,(l)}}$. The pairing of this commutator against $\mb V^{\dagger}_{1+\kappa,l}(\mb w_{n,\sigma})$ sits at the weight $1$ on both sides and is bounded uniformly in $n$ and $\sigma$, since $S_\sigma$ and its adjoint with respect to $r\d r\d\theta$ are bounded on the weighted spaces. So \eqref{second test function a priori} holds for the approximants uniformly, and Fatou's lemma along the two limits gives $\|\mb v\|_{\Omega,1+\kappa,2l+1,2}<\infty$ with the same bound.

Now note that
\begin{multline*}
\inner{\P\tilde{\mb v}-\epsilon^5\P\brac{\mb V^r_{-\kappa,(l)}(\mb v)+\mb V^\theta_{-\kappa,(l)}(\mb v)},\mb V^{\dagger}_{1+2\kappa,l}(\mb v)}_{L^2(\Omega)}\\
\lesssim_\epsilon\|\mb v\|_{\Omega,1+\kappa,2l+1,2}(\|\P\tilde{\mb v}\|_{\Omega,1+3\kappa,2l+1,2}+\|\mb v\|_{\Omega,1+\kappa,2l+1,2})
<\infty,
\end{multline*}
where $\|\P\tilde{\mb v}\|_{\Omega,1+3\kappa,2l+1,2}$ is finite for small enough $\theta_b$ by Lemma \ref{decay of projected test functions}. So we repeat the above process with $\smash{\mb V^{\dagger}_{1+2\kappa,l}}(\mb v)$. The approximation argument is the same, with $\|\mb v\|_{\Omega,1+\kappa,2l+1,2}$ now finite: the pairings that were at the weight $1$ are now at the weight $1+\kappa$, and the cut-off error is measured with the factor $r^{1+2\kappa}$. We get $\|\mb v\|_{\Omega,1+2\kappa,2l+1,2}<\infty$. This gives us enough regularity to get $\mb v\in H_s^l$.

Let $a\in\{0,-\kappa\}$. The bulk terms of $\|\mb v\|_{H_s^l,a}$ are bounded by $\|\grad\mb v\|_{\Omega,a,l+1,1}\lesssim\|\mb v\|_{\Omega,1+a,2l+1,2}\leq\|\mb v\|_{H^{(1),l}}$. Since $\mb v|_{\partial_0\Omega}=\mb 0$, we have $\|v_r\|_{\partial_0\Omega,a,l+1}=0$. By Lemma \ref{trace}, the $\partial_1\Omega$ term satisfies $\|\partial_rv_\theta\|_{\partial_1\Omega,a,l+1}\lesssim\|\grad\mb v\|_{\Omega,a+1/2,l+1,1}$, which on $\{r<1\}$ is bounded by $\|\grad\mb v\|_{\Omega,2\kappa,2l+1,1}\lesssim\|\mb v\|_{\Omega,1+2\kappa,2l+1,2}$, since $2\kappa>1/2$, which is where $\kappa>1/4$ is used, and on $\{r>1\}$ by $\|\grad\mb v\|_{\Omega,a,l+1,1}$.
\end{proof}

We can also show that the test functions can be prescribed on the boundary:

\begin{proposition}[Solvability on the boundary]\label{solvability on the boundary}
Suppose $g\in C^\infty_c((0,\infty))$.
\begin{enumerate}
\item For small enough $\epsilon>0$, there exists $\mb v\in H_s^l$ satisfying the regularity stated in (i) and (ii) of Proposition \ref{energy forms} such that
\begin{alignat*}{3}
v_\theta&=0\qquad&&\text{on}\qquad\partial\Omega\\
-\sum_{j=0}^{l+1}\epsilon^j(-1)^j\dr^{2j}v_r&=g\qquad&&\text{on}\qquad\partial_1\Omega\\
v_r&=0\qquad&&\text{on}\qquad\partial_0\Omega
\end{alignat*}
\item For small enough $\epsilon>0$, there exists $\mb v\in H_s^l$ satisfying the regularity stated in (i) and (ii) of Proposition \ref{energy forms} such that
\begin{alignat*}{3}
v_\theta&=0\qquad&&\text{on}\qquad\partial\Omega\\
v_r&=0\qquad&&\text{on}\qquad\partial_1\Omega\\
-\sum_{j=0}^{l+1}\epsilon^j(-1)^j\dr^{2j}v_r&=g\qquad&&\text{on}\qquad\partial_0\Omega
\end{alignat*}
\end{enumerate}
\end{proposition}
\begin{proof}
In Mellin variable $-\sum_{j=0}^{l+1}\epsilon^j(-1)^j\dr^{2j}v_\star=g$ is
\begin{align*}
-\sum_{j=0}^{l+1}\epsilon^j(-1)^j\lambda^{2j}\hat v_\star=\hat g.
\end{align*}
\begin{enumerate}
\item Let
\begin{align*}
\hat v_r(\lambda,\theta)&=f(\theta){\hat g(\lambda)\over -\sum_{j=0}^{l+1}\epsilon^j(-1)^j\lambda^{2j}}\\
\hat v_\theta(\lambda,\theta)&=-\int_0^\theta f(\vartheta){(\lambda+1)\hat g(\lambda)\over -\sum_{j=0}^{l+1}\epsilon^j(-1)^j\lambda^{2j}}\d\vartheta,
\end{align*}
where $f$ is a smooth function such that $f(0)=0$, $f(\theta_b)=1$ and $\int_0^{\theta_b}f(\theta)\d\theta=0$. Then we have $(\lambda+1)\hat v_r(\lambda,\theta)+\partial_\theta\hat v_\theta(\lambda,\theta)=0$, i.e. $\grad\cdot\mb v=0$. Also, this $\hat v_r$ and $\hat v_\theta$ satisfies all the required boundary conditions. For $\Re\lambda\in(-1,\alpha_0]$,  we have $|-\sum_{j=0}^{l+1}\epsilon^j(-1)^j\lambda^{2j}|^{-1}\lesssim_{\alpha_0} 1$ for small enough $\epsilon$. Indeed, with $z:=-\epsilon\lambda^2$ the sum is $-\sum_{j=0}^{l+1}z^j$. If $|z|\leq 1/3$, then $|\sum_{j=0}^{l+1}z^j|\geq 1/2$. If $|z|>1/3$, then $|\Im\lambda|\gtrsim\epsilon^{-1/2}$ on the strip, so $|\arg z|\lesssim_{\alpha_0}\epsilon^{1/2}$, which is at most $\pi/(2l+2)$ for small $\epsilon$. In that case every $z^j$ with $j\leq l+1$ has nonnegative real part, and the real part of $\sum_{j=0}^{l+1}z^j$ is at least $1$. Then
\begin{align*}
\sum_{\substack{n\leq 2l+2\\m\leq 3}}\|r^{-1-\alpha}\dr^n\partial_\theta^m\mb v\|_{L^2(\Omega)}^2&\lesssim\int_{\Re\lambda=\alpha}\sum_{\substack{n\leq 2l+2\\m\leq 3}}\brac{\|\lambda^{n}\partial_\theta^m\hat v_\theta\|_{L^2_\theta}^2+\|\lambda^{n}\partial_\theta^m\hat v_r\|_{L^2_\theta}^2}\d\Im\lambda\\
&\lesssim\sum_{n\leq 2l+3}\|r^{-{1\over 2}-\alpha}\dr^ng\|_{L^2((0,\infty))}^2<\infty
\end{align*}
Picking an appropriate $\alpha_0$, we see that $\mb v\in H_s^l$ and satisfies the regularity stated in (i) and (ii) of Proposition \ref{energy forms}.
\item Similar to above, with $f(0)=1$ and $f(\theta_b)=0$ in place of $f(0)=0$ and $f(\theta_b)=1$.\qedhere
\end{enumerate}
\end{proof}

\subsection{Solution of the linearised system}

Finally, we can prove our main theorems. First, the existence of solutions to System \ref{linear system}.

\begin{theorem}[Solution of the linear system on the Laplace side]\label{linear solution Laplace side}
Let $\kappa\in(1/4,1/2]$, $\varpi\in(0,\kappa]$, $l\in\N_0$, $\Re s>0$ and $\alpha\in[0,1]$. Suppose $\mb F$, $\mb f$ and $\phi_0$ have regularity that makes the RHS of \eqref{linear solution Laplace side estimate} finite, and that $\phi_0(0)=0$. Then, for small enough $\theta_b>0$, the variational problem $B_{(l)}[\mb u_s,\mb v]=\mathcal F_{(l)}[\mb v]$ for all $\mb v\in H_s^l$ has a unique solution $\mb u_s\in H_s^l$, with $H_s^l$ as in Definition \ref{H_s^l space} and $B_{(l)}$, $\mathcal F_{(l)}$ as in Definition \ref{energy forms and functionals}. The solution satisfies System \ref{linear system} strongly, except for the condition $\partial_ru_\theta(0,\theta_b)=0$, which it satisfies weakly, as a natural boundary condition of the variational problem, and
\begin{align}
\|\mb u_s\|_{H_s^l}
&\lesssim_{\theta_b,\varpi,\kappa,l}\sum_{k=0}^1\brac{1+\delta_{0k}|s|^{-\kappa}}\brac{\|r\P\mb F\|_{\Omega,-\kappa k,l}+\|\mb f\|_{\partial_1\Omega,-\kappa k-{1\over 2},l+{1\over 2}}}\nonumber\\
&\quad+\brac{1+|s|^{-2\kappa}}\brac{\|r\P\mb F\|_{\Omega,\varpi}+\|\mb f\|_{\partial_1\Omega,\varpi-{1\over 2},{1\over 2}}}\nonumber\\
&\quad+(\sin\theta_b)^2\brac{\|s^{-1}\phi_0\|_{\partial_1\Omega,{1\over 2}-\kappa,l+{3\over 2}}+\brac{1+|s|^{-\kappa}}\|s^{-{1\over 2}(1+\alpha)}\phi_0\|_{\partial_1\Omega,{\alpha\over 2},l+{3\over 2}}}.\label{linear solution Laplace side estimate}
\end{align}
We write $\mathfrak D^{(l)}_s[\mb F,\mb f;\phi_0]$ for the RHS of \eqref{linear solution Laplace side estimate}.
\end{theorem}
\begin{proof}
By Proposition \ref{energy forms}, we have
\begin{align*}
\|\mathcal F_{(l)}\|_{H_s^l\to\C}
&\lesssim_{\epsilon,\varpi}\sum_{k=0}^1\brac{1+\delta_{0k}|s|^{-\kappa}}\brac{\|r(\P\mb F+\grad p_{\partial_1,f_\theta})\|_{\Omega,-\kappa k,l}+\|\dr^{-{1\over 2}}f_r\|_{\partial_1\Omega,-\kappa k-{1\over 2},l+1}}\\
&\quad+\brac{1+|s|^{-2\kappa}}\brac{\|r(\P\mb F+\grad p_{\partial_1,f_\theta})\|_{\Omega,\varpi}+\|f_r\|_{\partial_1\Omega,\varpi-{1\over 2}}}\\
&\quad+\brac{1+|s|^{-\kappa}}(\sin\theta_b)^2\|s^{-{1\over 2}(1+\alpha)}\phi_0\|_{\partial_1\Omega,{\alpha\over 2},l+{3\over 2}}
+\|r\grad p_{\partial_1,\phi_0}\|_{\Omega,-\kappa,l},
\end{align*}
where the corner term of Proposition \ref{energy forms}(iii) is absent since $\phi_0(0)=0$.
By Lemma \ref{pressure estimates}(iv), applied to the two extensions $p_{\partial_1,\phi_0}=\E g_{1,\phi_0}$ and $p_{\partial_1,f_\theta}=\E g_{1,f_\theta}$, we have
\begin{align*}
\|(r\grad)^jp_{\partial_1,\phi_0}[\phi_0]\|_{\Omega,a} &\lesssim(\sin\theta_b)^{2}\|s^{-1}\partial_r\phi_0\|_{\partial_1\Omega,a-{1\over 2},j-{1\over 2}}\\
\|(r\grad)^jp_{\partial_1,f_\theta}[f_\theta]\|_{\Omega,a} &\lesssim\|f_\theta\|_{\partial_1\Omega,a-{1\over 2},j-{1\over 2}}
\end{align*}
So we have
\begin{align*}
\|\mathcal F_{(l)}\|_{H_s^l\to\C} &\lesssim\mathfrak D^{(l)}_s[\mb F,\mb f;\phi_0]
\end{align*}

By Proposition \ref{energy form coercivity} the form $B_{(l)}$ is bounded and coercive on $H_s^l$, and by the estimate above $\mathcal F_{(l)}$ is a bounded linear functional on $H_s^l$. Applying the Lax--Milgram theorem, we therefore get that there exists $\mb u_s\in H_s^l$ with $\|\mb u_s\|_{H_s^l}\lesssim\|\mathcal F_{(l)}\|_{H_s^l\to\C}$ such that
\begin{align*}
0=B_{(l)}[\mb u_s,\mb v]-\mathcal F_{(l)}[\mb v]\qquad\text{for all}\qquad\mb v\in H_s^l.
\end{align*}
By Proposition \ref{solvability in the bulk}, given any $\tilde{\mb v}\in C^\infty_c(\Omega)$, there exists $\mb v\in H^{(1),l}$ such that $\smash{\P\mb V^{(1)}_{(l)}(\mb v)=\P\tilde{\mb v}}$, and by Lemma \ref{test function upgrade} this $\mb v$ lies in fact in $H_s^l$. In particular it satisfies the assumptions of (i) and (ii) in Proposition \ref{energy forms}. So by Proposition \ref{energy forms},
\begin{align*}
B_{(l)}[\mb u_s,\mb v]-\mathcal F_{(l)}[\mb v]&=\<\lpc\mb u_s-\grad p_\P-\P\mb F,\mb V^{(1)}_{(l)}(\mb v)\>_{L^2(\Omega)}\\
&\quad+\sum_{m=0}^1(-1)^m\<(2\mb e_\theta\cdot\S(\mb u_s)\cdot\mb e_r-\delta_{0m}u_{s,r}-\delta_{1m}f_r),\mb e_r\cdot\mb V^{(1),r}_{(l)}(\mb v)\>_{L^2(\partial_m\Omega)}\\
&=\<\lpc\mb u_s-\grad p_\P-\P\mb F,\mb V^{(1)}_{(l)}(\mb v)\>_{L^2(\Omega)},
\end{align*}
where the boundary terms all vanish because $\mb v\in H^{(1),l}$. So we have
\begin{align*}
0=\<\lpc\mb u_s-\grad p_\P-\P\mb F,\P\tilde{\mb v}\>_{L^2(\Omega)}=\<\lpc\mb u_s-\grad p_\P-\P\mb F,\tilde{\mb v}\>_{L^2(\Omega)}\qquad\text{for all}\qquad\tilde{\mb v}\in C^\infty_c(\Omega).
\end{align*}
Therefore, $\lpc\mb u_s-\grad p_\P-\P\mb F=0$.

Now using Proposition \ref{solvability on the boundary}, for any $g\in C^\infty_c((0,\infty))$, we can find $\mb v_g\in H_s^l$ satisfying the regularity stated in (i) and (ii) of Proposition \ref{energy forms} such that
\begin{align*}
(\mb e_r\cdot\mb V^{(1),r}_{(l)}(\mb v_g))|_{\partial_1\Omega}=((1+|s|^{-2\kappa})+\epsilon^5r^{2\kappa})g\qquad\text{and}\qquad(\mb e_r\cdot\mb V^{(1),r}_{(l)}(\mb v_g))|_{\partial_0\Omega}=0,
\end{align*}
the two weights $a=0$ and $a=-\kappa$ of $\smash{\mb V^{(1)}_{(l)}}$ producing the two summands of the prefactor. So we have
\begin{align*}
0&=-\<(2\mb e_\theta\cdot\S(\mb u_s)\cdot\mb e_r-f_r),\mb e_r\cdot\mb V^{(1),r}_{(l)}(\mb v_g)\>_{L^2(\partial_1\Omega)}\\
&=-\<(2\mb e_\theta\cdot\S(\mb u_s)\cdot\mb e_r-f_r),((1+|s|^{-2\kappa})+\epsilon^5r^{2\kappa})g\>_{L^2(\partial_1\Omega)}
\end{align*}
Since $g\in C^\infty_c((0,\infty))$ is arbitrary and the prefactor is positive, we have that $2\mb e_\theta\cdot\S(\mb u_s)\cdot\mb e_r-f_r=0$ on $\partial_1\Omega$. Similarly, we can show that $2\mb e_\theta\cdot\S(\mb u_s)\cdot\mb e_r-u_{s,r}=0$ on $\partial_0\Omega$.
\end{proof}

\begin{remark}[The corner condition and uniqueness]\label{corner condition and uniqueness}
The condition $\partial_ru_\theta(0,\theta_b)=0$ of System \ref{linear system} is encoded in the form $B_{(l)}$ rather than imposed on $H_s^l$. In the proof of Proposition \ref{energy forms}(i), the term of $B^r_{0,0}$ carrying $\bar s^{-1}$ comes from an integration by parts in $r$ on $\partial_1\Omega$ whose endpoint term at $r=0$ pairs $\partial_ru_\theta(0,\theta_b)$ with $v_\theta(0,\theta_b)$. This endpoint term vanishes for the test functions carrying the extra regularity of that proposition, but $H_s^l$ also contains fields with $v_\theta(0,\theta_b)\neq 0$, by Lemma \ref{regular test functions}(iii) below. So the condition is a natural boundary condition of the variational problem, and it holds once the solution is shown to have enough regularity.

Uniqueness in Theorem \ref{linear solution Laplace side} is asserted for the variational problem only, which is what the Lax--Milgram theorem gives. A strong solution in $H_s^l$ satisfies the variational identity for the test functions carrying the extra regularity of Proposition \ref{energy forms}, and these are not dense in $H_s^l$: by Lemma \ref{regular test functions}(ii) below, their closure is precisely the kernel of the corner value $v_\theta(0,\theta_b)$. Variational uniqueness is used on the Laplace side in Lemmas \ref{lowering the count}, \ref{stationary as variational} and \ref{time derivative of the solution}, and uniqueness of the time-dependent problem is proved separately in Lemma \ref{linear evolution uniqueness}.
\end{remark}

Now from $(\mb u,p,\phi_0)$ we can recover the profile $\phi$ using the profile evolution equation.

\begin{remark}[Where the angle is taken small]\label{smallness remark}
The hypothesis ``for small enough $\theta_b$'' of Theorem \ref{linear solution Laplace side} is not quantified, and it is the only unquantified smallness of the paper: every other condition on the angle is explicit, and is collected in \eqref{standing angle condition} and \eqref{third rung angle condition}. The smallness enters the linear theory only through absorption, and not by perturbation from the degenerate wedge $\theta_b=0$. The trace inequality of Lemma \ref{trace} costs $\theta_b^{-1/2}$, the pressure functionals of Lemma \ref{pressure estimates} cost $\theta_b^{-1}$ and $\theta_b^{-2}$, and the test functions carry the parameter $\epsilon=c_\epsilon\theta_b$ of the norms. Integrations across the wedge, on the other hand, gain powers of $\theta_b$, as in Lemma \ref{weighted based control}. The coercivity of Proposition \ref{energy form coercivity}, and with it the test function problem, is then obtained by absorbing remainders whose constants carry positive powers of $\theta_b$. Two explicit but loose conditions of the same kind are $\theta_b<\pi/(3K)$ in Lemmas \ref{pressure estimates} and \ref{test problem pressure estimate} and $a'\theta_b<\pi/2$ in Lemma \ref{decay of projected test functions}. All of these could be quantified, but the number so obtained would not reflect an obstruction: the explicit conditions \eqref{standing angle condition} and \eqref{third rung angle condition} come from the first Moffatt-type eigenvalue of the Stokes wedge and from the resonance of the flattening, see Remarks \ref{angle remark}, \ref{sharpness of the second angle condition} and \ref{H numerics remark}.
\end{remark}

\begin{corollary}[Recovering the profile on the Laplace side]\label{profile on the Laplace side}
Take the setting of Theorem \ref{linear solution Laplace side}. If in addition we define $\phi$ by
\begin{align}
  s\phi-(\sin\theta_b)^{-2}r^{-1}\dr u_\theta &= \phi_0\qquad&&\text{on}\qquad\partial_1\Omega\label{phi new def}
\end{align}
then $(\phi,\mb u,p_\P)$ has the following regularity:
\begin{align*}
\|\mb u\|_{H_s^l}&+(\sin\theta_b)^2\|\phi\|_{\partial_1\Omega,{1\over 2}-\kappa,l+{3\over 2}}+\sum_{k=0}^1\|p_\P\|_{\Omega,-\kappa k,l,1}
\lesssim_{\theta_b,\varpi,\kappa,l}\mathfrak D^{(l)}_s[\mb F,\mb f;\phi_0].
\end{align*}
Moreover, $(\phi,\mb u,p_\P)$ is a strong solution to the system
\begin{subequations}\label{system: linear with profile}
\begin{alignat}{3}
\lpc\mb u-\grad p_\P&=\P\mb F\qquad&&\text{on}\qquad\Omega\\
\grad\cdot\mb u&=0\qquad&&\text{on}\qquad\Omega\\
u_\theta&=0\qquad&&\text{on}\qquad\partial_0\Omega\\
r^{-1}\partial_\theta u_r-u_r&=0\qquad&&\text{on}\qquad\partial_0\Omega\\
2\mb e_r\cdot\S(\mb u)\cdot\mb e_\theta&=f_r\qquad&&\text{on}\qquad\partial_1\Omega\\
\mb e_\theta\cdot \mathbb T[\mb u, p_\P]\cdot\mb e_\theta &= (\sin\theta_b)^2r^{-1}\dr\phi+f_\theta\qquad&&\text{on}\qquad\partial_1\Omega\\
s\phi-(\sin\theta_b)^{-2}r^{-1}\dr u_\theta &= \phi_0\qquad&&\text{on}\qquad\partial_1\Omega
\intertext{and satisfies weakly}
\partial_ru_\theta(0,\theta_b)&=0.
\end{alignat}
\end{subequations}
\end{corollary}
\begin{proof}
Since $(\mb u,p_\P)$ is a solution to System \ref{linear system}, it satisfies the first five equations of \eqref{system: linear with profile} together with $\mb e_\theta\cdot\mathbb T[\mb u,p_\P]\cdot\mb e_\theta=s^{-1}(r^{-2}\dr\dr_{-1}u_\theta+(\sin\theta_b)^{2}\partial_r\phi_0)+f_\theta$ on $\partial_1\Omega$, by Definition \ref{g definition}.
Taking $\partial_r$ of \eqref{phi new def}, we get $s(\sin\theta_b)^2r^{-1}\dr\phi =r^{-2}\dr\dr_{-1}u_\theta+(\sin\theta_b)^{2}\partial_r\phi_0$. Therefore we get the system \eqref{system: linear with profile}. And from \eqref{system: linear with profile} we have
\begin{align*}
\dr\phi  &= (\sin\theta_b)^{-2}r\brac{\mb e_\theta\cdot \mathbb T[\mb u, p_\P]\cdot\mb e_\theta -f_\theta}
= (\sin\theta_b)^{-2}\brac{-2\dr u_r -rp_\P-rf_\theta} \qquad\qquad\text{on}\quad\partial_1\Omega
\end{align*}
so we get, since $|\lambda|\geq\kappa$ there,
\begin{align*}
\|\phi\|_{\partial_1\Omega,{1\over 2}-\kappa,l+{3\over 2}}
&\lesssim_\kappa(\sin\theta_b)^{-2}\|-2\dr u_r -rp_\P-rf_\theta\|_{\partial_1\Omega,{1\over 2}-\kappa,l+{1\over 2}}\\
&\leq(\sin\theta_b)^{-2}\brac{2\|\dr\mb u\|_{\Omega,1-\kappa,l,1}+\|p_\P\|_{\Omega,-\kappa,l,1}+\|f_\theta\|_{\partial_1\Omega,-{1\over 2}-\kappa,l+{1\over 2}}}\\
&\lesssim(\sin\theta_b)^{-2}\brac{\|\mb u\|_{H_s^l}+\|p_\P\|_{\Omega,-\kappa,l,1}+\|f_\theta\|_{\partial_1\Omega,-{1\over 2}-\kappa,l+{1\over 2}}}.
\end{align*}
The same argument does not bound $\|\phi\|_{\partial_1\Omega,1/2}$; see Remark \ref{weight half remark} below.
And from $\lpc\mb u-\grad p_\P =\P\mb F$ we have
\begin{align*}
\sum_{k=0}^1\|p_\P\|_{\Omega,-\kappa k,l,1}
\lesssim\sum_{k=0}^1\|r\grad p_\P\|_{\Omega,-\kappa k,l}
\leq\sum_{k=0}^1\|r(\lpc\mb u-\P\mb F)\|_{\Omega,-\kappa k,l}
\lesssim\|\mb u\|_{H_s^l}+\sum_{k=0}^1\|r\P\mb F\|_{\Omega,-\kappa k,l}
\end{align*}
Adding this to the estimate in Theorem \ref{linear solution Laplace side} we get the desired inequality.
\end{proof}

\begin{remark}[The weight ${1\over 2}$ of the profile]\label{weight half remark}
Corollary \ref{profile on the Laplace side} controls the velocity and the pressure at both weights, but the profile only in $\|\phi\|_{\partial_1\Omega,1/2-\kappa,l+3/2}$. The missing norm $\|\phi\|_{\partial_1\Omega,1/2}$ is finite if and only if $\int_0^1r^{-1}|\phi|^2\d r<\infty$, its behaviour at infinity being covered by $\|\phi\|_{\partial_1\Omega,1/2-\kappa}$; and since, in the corner expansions of Section \ref{corner expansions}, $\partial_ru_\theta(\ph,\theta_b)$ is the constant $\partial_ru_\theta(0,\theta_b)$ plus a function with this integrability, by \eqref{phi new def} the missing norm is finite exactly when the contact angle condition $\partial_ru_\theta(0,\theta_b)=0$ holds. The variational solution satisfies this condition weakly, but at fixed $s$ the expansion that makes it quantitative needs data two weights above those of Theorem \ref{linear solution Laplace side}, so the norm is not available here. It is obtained in Theorem \ref{linear solution full regularity}, from the rung $M=1$ of \eqref{Improved solution estimate}, and it enters the class of Lemma \ref{linear evolution uniqueness}.
\end{remark}

\begin{lemma}[The corner value and the regular test functions]\label{regular test functions}
Let $l\in\N_0$ and $\Re s>0$, and let $R_l\subseteq H_s^l$ be the set of those $\mb v\in H_s^l$ that carry the extra regularity of Proposition \ref{energy forms}(i)--(iv) at both weights $a\in\{0,-\kappa\}$, namely
\begin{align}\label{regular class}
\|\mb v\|_{\Omega,a+1,2(l+1),1}+\|v_\theta\|_{\partial_1\Omega,a+1,2(l+1)}+\|(r\grad)\grad\mb v\|_{\Omega,a,2l}<\infty.
\end{align}
\begin{enumerate}
\item For $\mb v\in H_s^l$ the trace $v_\theta(\ph,\theta_b)$ is continuous on $[0,\infty)$, and its corner value satisfies $|v_\theta(0,\theta_b)|\lesssim_{\theta_b,\kappa}(1+|s|^{1/2})\|\mb v\|_{H_s^l}$.
\item Every $\mb v\in R_l$ has $v_\theta(0,\theta_b)=0$, and the closure of $R_l$ in $H_s^l$ is $\{\mb v\in H_s^l:v_\theta(0,\theta_b)=0\}$.
\item Let $\eta$ be smooth on $[0,\infty)$ with $\eta=1$ on $[0,1]$ and $\eta=0$ on $[2,\infty)$, and let $\mb v_*:=\grad^\perp(y\eta(r))$. Then $\mb v_*\in H_s^l$, $\mb v_*=-\mb e_1$ on $\Omega\cap\{r\leq 1\}$, $\mb v_*=\mb 0$ for $r\geq 2$, and $v_{*,\theta}(0,\theta_b)=\sin\theta_b\neq 0$. Consequently, two bounded linear functionals on $H_s^l$ that agree on $R_l$ and at $\mb v_*$ agree on $H_s^l$.
\end{enumerate}
\end{lemma}
\begin{proof}
(i) The norm of $H_s^l$ contains $|s|^{-1/2}\|\partial_rv_\theta\|_{L^2(\partial_1\Omega)}$, so $v_\theta(\ph,\theta_b)$ is locally absolutely continuous on $(0,\infty)$ with $|v_\theta(r,\theta_b)-v_\theta(r',\theta_b)|\leq|r-r'|^{1/2}\|\partial_rv_\theta\|_{L^2(\partial_1\Omega)}$; hence it extends continuously to $r=0$, with $|v_\theta(0,\theta_b)|\leq|v_\theta(r,\theta_b)|+r^{1/2}\|\partial_rv_\theta\|_{L^2(\partial_1\Omega)}$. Averaging over $r\in[1,2]$ and using the trace inequality on the annulus $\Omega_{1,2}$,
\begin{align*}
|v_\theta(0,\theta_b)|\lesssim_{\theta_b}\|\mb v\|_{L^2(\Omega_{1,2})}+\|\grad\mb v\|_{L^2(\Omega)}+\|\partial_rv_\theta\|_{L^2(\partial_1\Omega)}
\lesssim_\kappa(1+|s|^{1\over 2})\|\mb v\|_{H_s^l},
\end{align*}
since $\|\mb v\|_{L^2(\Omega_{1,2})}\lesssim\|\mb v\|_{\Omega,1-\kappa}\lesssim_\kappa\|\grad\mb v\|_{\Omega,-\kappa}$ by Hardy's inequality at the weight $1-\kappa$, as in the proof of Proposition \ref{energy forms}(iii).

(ii) For $\mb v\in R_l$ the norm $\|v_\theta\|_{\partial_1\Omega,1}$ of \eqref{regular class} at $a=0$ contains $\|r^{-1}v_\theta\|_{L^2(\partial_1\Omega)}$, which is infinite unless the continuous function $v_\theta(\ph,\theta_b)$ vanishes at $r=0$; so $R_l$ lies in the kernel of the corner value, and so does its closure, by (i). Conversely, let $\mb v\in H_s^l$ with $v_\theta(0,\theta_b)=0$. We first show
\begin{align}\label{corner smallness}
\epsilon^{-1}\|\mb v\|_{L^2(\Omega_{\epsilon,2\epsilon})}\to 0\qquad(\epsilon\to 0).
\end{align}
Let $\mb c_\epsilon\in\C^2$ be the mean of $\mb v$ over $\Omega_{\epsilon,2\epsilon}$. The Poincar\'e inequality on $\Omega_{\epsilon,2\epsilon}$ and the trace inequality on its two straight sides, both of which are invariant under scaling, give $\|\mb v-\mb c_\epsilon\|_{L^2(\Omega_{\epsilon,2\epsilon})}\lesssim_{\theta_b}\epsilon\|\grad\mb v\|_{L^2(\Omega_{\epsilon,2\epsilon})}$ and $\int_\epsilon^{2\epsilon}|(\mb v-\mb c_\epsilon)(r,\theta)|^2\d r\lesssim_{\theta_b}\epsilon\|\grad\mb v\|_{L^2(\Omega_{\epsilon,2\epsilon})}^2$ for $\theta\in\{0,\theta_b\}$. On $\theta=0$ we have $(\mb v-\mb c_\epsilon)\cdot\mb e_\theta=-\mb c_\epsilon\cdot\mb e_\theta(0)$, and on $\theta=\theta_b$ we have $(\mb v-\mb c_\epsilon)\cdot\mb e_\theta=v_\theta(r,\theta_b)-\mb c_\epsilon\cdot\mb e_\theta(\theta_b)$ with $|v_\theta(r,\theta_b)|\leq r^{1/2}\|\partial_rv_\theta\|_{L^2([0,2\epsilon])}$ by (i) and the vanishing corner value; hence
\begin{align*}
|\mb c_\epsilon\cdot\mb e_\theta(0)|+|\mb c_\epsilon\cdot\mb e_\theta(\theta_b)|\lesssim_{\theta_b}\|\grad\mb v\|_{L^2(\Omega_{\epsilon,2\epsilon})}+\epsilon^{1\over 2}\|\partial_rv_\theta\|_{L^2([0,2\epsilon])}\to 0,
\end{align*}
and since $\mb e_\theta(0)$ and $\mb e_\theta(\theta_b)$ are linearly independent, $|\mb c_\epsilon|\to 0$. This gives \eqref{corner smallness}, because
\begin{align*}
\epsilon^{-1}\|\mb v\|_{L^2(\Omega_{\epsilon,2\epsilon})}\leq\epsilon^{-1}\|\mb v-\mb c_\epsilon\|_{L^2(\Omega_{\epsilon,2\epsilon})}+\epsilon^{-1}|\mb c_\epsilon||\Omega_{\epsilon,2\epsilon}|^{1\over 2}\lesssim_{\theta_b}\|\grad\mb v\|_{L^2(\Omega_{\epsilon,2\epsilon})}+|\mb c_\epsilon|.
\end{align*}

Now write $\mb v=\grad^\perp\psi$ with a stream function normalised by $\psi=0$ on $\partial_0\Omega$, which is possible since $v_\theta=\partial_r\psi$ vanishes there, and let $\mb v_n:=\grad^\perp(\zeta_n\psi)$ with the cut-off $\zeta_n$ of the proof of Proposition \ref{solvability for V2}, equal to $1$ on $[2^{1-n},2^n]$ and supported in $[2^{-n},2^{n+1}]$; then $\mb v_n$ is divergence-free with $v_{n,\theta}=0$ on $\partial_0\Omega$, and $\mb v_n\to\mb v$ in $H_s^l$. Indeed $\partial_\theta\psi=-rv_r$ in the bulk and $\partial_r\psi=v_\theta$ on $\partial_1\Omega$ bound $\psi$ and its derivatives by those of $\mb v$. The terms of $\mb v_n-\mb v$ in which a derivative falls on $\zeta_n$ are supported in the two annuli $\Omega_{\epsilon,2\epsilon}$, $\epsilon\in\{2^{-n},2^n\}$. So they are dominated there by pieces of the norm of $\mb v$ over sets of vanishing measure, by $\|\partial_rv_\theta\|_{L^2([0,2\epsilon])}$, and by $\|r^{-1-a}\mb v\|_{L^2(\Omega_{\epsilon,2\epsilon})}$. The norm $\|r^{-1-a}\mb v\|_{L^2(\Omega_{\epsilon,2\epsilon})}$ tends to $0$ at the corner by \eqref{corner smallness}, and at infinity because it is at most $\|r^{\kappa-1}\mb v\|_{L^2(\Omega_{\epsilon,2\epsilon})}$ there, with $r^{\kappa-1}\mb v\in L^2(\Omega)$ by Hardy's inequality as in (i). The terms in which no derivative falls on $\zeta_n$ go by dominated convergence. Finally, each $\mb v_n$, having compact support in $r$, is the limit in $H_s^l$ of its radial smoothings $S_\sigma\mb v_n$ of Lemma \ref{radial smoothing}, applied componentwise in polar components. These smoothings preserve $\grad\cdot\mb v_n=0$ and $v_{n,\theta}|_{\partial_0\Omega}=0$, and they lie in $R_l$: the Mellin multiplier $e^{\sigma^2\lambda^2}$ decays like a Gaussian on every line and so supplies the radial derivatives that \eqref{regular class} asks for. Hence $\mb v$ lies in the closure of $R_l$.

(iii) $\mb v_*$ is smooth on $\overline\Omega$, divergence-free as the rotated gradient of a stream function, and equals $\grad^\perp y=-\mb e_1$ where $\eta=1$; on $\partial_0\Omega$ its stream function $y\eta(r)$ vanishes, so $v_{*,\theta}=0$ there, and $v_{*,\theta}(0,\theta_b)=-\mb e_1\cdot\mb e_\theta(\theta_b)=\sin\theta_b$. Every norm in $\|\ph\|_{H_s^l}$ is finite, the derivatives of $\mb v_*$ vanishing near the corner and $\mb v_*$ vanishing for $r\geq 2$. Given $\mb v\in H_s^l$, the field $\mb v-(\sin\theta_b)^{-1}v_\theta(0,\theta_b)\mb v_*$ has vanishing corner value, hence lies in the closure of $R_l$ by (ii), on which the two functionals agree by continuity; and they agree at $\mb v_*$.
\end{proof}

\begin{lemma}[Lowering the derivative count]\label{lowering the count}
Let $l\leq l'$ be in $\N_0$, and let $\mb u\in H_s^{l'}$ be the solution of Theorem \ref{linear solution Laplace side} at the derivative count $l'$, for data with $\mathfrak D^{(l')}_s[\mb F,\mb f;\phi_0]<\infty$. Then $\mb u\in H_s^l$ and $B_{(l)}[\mb u,\mb v]=\mathcal F_{(l)}[\mb v]$ for all $\mb v\in H_s^l$; that is, $\mb u$ is also the solution at the count $l$.
\end{lemma}
\begin{proof}
The forms at the two counts are taken with the same parameter $\epsilon=c_\epsilon\theta_b$, admissible for both once $c_\epsilon$ and $\theta_b$ are small enough for both counts in Proposition \ref{energy form coercivity}, which is how the parameters are chosen whenever several counts are in play. Since $l\leq l'$ we have $\|\ph\|_{H_s^l}\leq\|\ph\|_{H_s^{l'}}$, so $\mb u\in H_s^l$, and both sides of the identity are bounded linear functionals of $\mb v\in H_s^l$, the RHS by Proposition \ref{energy forms}(iii)--(iv). By Lemma \ref{regular test functions}(iii) it suffices to check the identity on $R_l$ and at $\mb v_*$. On $R_l$ it is Proposition \ref{energy forms}: by \eqref{undo integration main} the difference of the two sides is the pairing of $\lpc\mb u-\grad p_\P-\P\mb F$ and of the two boundary expressions with the test functions, and these vanish since $\mb u$ solves System \ref{linear system} strongly. 

At $\mb v_*$ we compare with the identity at the count $l'$, which holds since $\mb v_*\in H_s^{l'}$. The forms and functionals at the two counts differ only by the terms $j\geq l+2$ of $B^r_{a,j}$ and $\mathcal F^r_{a,j}$ and the terms $j\geq l+1$ of $B^\theta_{a,j}$ and $\mathcal F^\theta_{a,j}$. The identities of Proposition \ref{energy forms} hold at $\mb v=\mb v_*$ individually, for (ii) and (iv) at every $j$ and for (i) and (iii) at every $j\geq 1$. For (ii) and (iv) this is the proposition itself, $\mb v_*$ having the extra regularity $\|(r\grad)\grad\mb v_*\|_{\Omega,a,2l'}<\infty$ required there; the test functions $\P^\theta\smash[b]{\mb V^\theta_{a,j}}(\mb v_*)$ are not supported away from the corner, the projection $\P^\theta$ being nonlocal, but no corner condition is encoded in these two identities. For (i) and (iii) at $j\geq 1$ the test functions $\mb V^r_{a,j}(\mb v_*)$ vanish near the corner, where $\mb v_*$ is constant, and for $r\geq 2$. So every endpoint term of the radial integrations by parts in the proof of the proposition, each of which carries at least one radial derivative of $\mb v_*$, vanishes identically. For the same reason the integrations by parts in the bulk take place away from the corner, where $\mb u$ is a strong solution. Each of the differences $B^\bullet_{a,j}[\mb u,\mb v_*]-\mathcal F^\bullet_{a,j}[\mb v_*]$ in question is therefore the pairing of $\lpc\mb u-\grad p_\P-\P\mb F$, and of the boundary expressions of \eqref{undo integration main}, with the corresponding test function, and vanishes by the strong equations; subtracting them from the identity at the count $l'$ leaves the identity at the count $l$, evaluated at $\mb v_*$.
\end{proof}

In the next lemma the Navier slip term is assumed finite, as the solution of Proposition \ref{stationary system} has it.

\begin{lemma}[Stationary solutions as variational solutions]\label{stationary as variational}
Let $l\in\N_0$ and $\Re s>0$. Let $(\mb u_\star,p_\star)$ solve the stationary system \eqref{system: stationary linear} a.e.\ with data $(\phi_\star,\mb F_\star,\mb f_\star)$ whose norms on the right of \eqref{data estimate} are finite at the count $l+1$. Suppose that $(\mb u_\star,p_\star)$ lies in the class of the left of \eqref{data estimate} at that count, with the Navier slip term $\sum_{k=0}^1\|u_{\star,r}\|_{\partial_0\Omega,-\kappa k,l+2}$ finite. Suppose also that $\mb u_\star$ has a corner value $\mb u_\star(\mb 0)\in\C^2$ with $\|\mb u_\star-\mb u_\star(\mb 0)\|_{\Omega,2-\kappa,l+1,2}<\infty$, and that the contact angle condition holds in the sense
\begin{align}\label{stationary corner condition}
\partial_ru_{\star,\theta}(r,\theta_b)\to 0\qquad(r\to 0).
\end{align}
Let $\dot\phi_\star:=(\sin\theta_b)^{-2}\partial_ru_{\star,\theta}$ on $\partial_1\Omega$, and suppose that $\mathfrak D^{(l)}_s[\mb F_\star,\mb f_\star;s\phi_\star-\dot\phi_\star]<\infty$. Then $\mb u_\star\in H_s^l$, and $\mb u_\star$ is the solution of that theorem for these data:
\begin{align}\label{stationary variational identity}
B_{(l)}[\mb u_\star,\mb v]=\mathcal F_{(l)}[\mb v]\qquad\text{for all}\quad\mb v\in H_s^l,
\end{align}
with $\mathcal F_{(l)}$ formed with the bulk forcing $\mb F_\star$, the boundary forcing $\mb f_\star$ and the profile datum $\phi_0=s\phi_\star-\dot\phi_\star$.
\end{lemma}
\begin{proof}
First, $\mb u_\star\in H_s^l$: the bulk terms and the $\partial_0\Omega$ term of the norm are among the norms assumed, and the $\partial_1\Omega$ term $|s|^{-1/2}\|\partial_ru_{\star,\theta}\|_{\partial_1\Omega,a,l+1}$, $a\in\{0,-\kappa\}$, is finite since, by Lemma \ref{trace},
\begin{align*}
\|\partial_ru_{\star,\theta}\|_{\partial_1\Omega,a,l+1}\lesssim_{\theta_b}\|\grad\mb u_\star\|_{\Omega,a+{1\over 2},l+1,1}\leq\|\grad\mb u_\star\|_{\Omega,-\kappa,l+1,1}+\|\grad\mb u_\star\|_{\Omega,1-\kappa,l+1,1},
\end{align*}
the two norms on the right being finite by the class of \eqref{data estimate} and by the assumption on $\mb u_\star-\mb u_\star(\mb 0)$. Likewise $\|\dr u_{\star,\theta}\|_{\partial_1\Omega,3/2-\kappa,l+1}\lesssim_{\theta_b}\|\mb u_\star-\mb u_\star(\mb 0)\|_{\Omega,2-\kappa,l+1,2}<\infty$.

Both sides of \eqref{stationary variational identity} are bounded linear functionals of $\mb v\in H_s^l$, the RHS by Proposition \ref{energy forms}(iii)--(iv) and the assumption on the data, so by Lemma \ref{regular test functions}(iii) it suffices to prove the identity on $R_l$ and at $\mb v_*$. 

We first express the pressure. The function $p_\star+q[\mb F_\star]-p_{\P|\mb u_\star}$, with $q$ as in Definition \ref{the two projections} and $p_{\P|\mb u_\star}$ as in Definition \ref{pressure functionals}, is harmonic, since $-\lpc p_\star=\grad\cdot\mb F_\star=\lpc q[\mb F_\star]$ by \eqref{system: stationary linear}. On $\partial_0\Omega$ its $\partial_\theta$ derivative vanishes: the $\mb e_\theta$ component of the Stokes equation gives $\partial_\theta p_\star=g_0[\mb u_\star]-rF_{\star,\theta}$ there, as in the derivation of System \ref{linear system}, while $\partial_\theta q[\mb F_\star]=rF_{\star,\theta}$ and $\partial_\theta p_{\P|\mb u_\star}=g_0[\mb u_\star]$. On $\partial_1\Omega$, where $q[\mb F_\star]=0$, the $\mb e_\theta$ component of the surface tension condition of \eqref{system: stationary linear} gives
\begin{align*}
p_\star-p_{\P|\mb u_\star}=-(\sin\theta_b)^2\partial_r\phi_\star-f_{\star,\theta}+{1\over s}\partial_r^2u_{\star,\theta}=g_{1,f_{\star,\theta}}+sg_{1,\phi_\star}-g_{1,\dot\phi_\star}
\end{align*}
in the notation of Definition \ref{g definition}. So $p_\star+q[\mb F_\star]-p_{\P|\mb u_\star}$ and $p_{\partial_1,f_{\star,\theta}}+sp_{\partial_1,\phi_\star}-p_{\partial_1,\dot\phi_\star}$ solve the mixed problem of Lemma \ref{mixed wedge problem} with the same data. Both have finite norm $\|\ph\|_{\Omega,-\kappa}$. For the first this follows from the class of $(\mb u_\star,p_\star)$, from Lemma \ref{mixed wedge problem} for $q[\mb F_\star]$, and from Lemma \ref{pressure estimates} for the pressure functionals of $\mb u_\star$; Lemma \ref{pressure estimates}(ii) uses $\|\dr u_{\star,\theta}\|_{\partial_1\Omega,3/2-\kappa,1}<\infty$. For the second it follows from Lemma \ref{pressure estimates}(iv), applied to the extensions of $f_{\star,\theta}$, $\partial_r\phi_\star$ and $\partial_r\dot\phi_\star$, with $\|\partial_r\phi_\star\|_{\partial_1\Omega,-1/2-\kappa}\leq\|\phi_\star\|_{\partial_1\Omega,1/2-\kappa,1}$ and the same for $\dot\phi_\star$; the norms on the right are finite by the assumptions on the data. The difference of the two functions is harmonic with vanishing data, and its Mellin transform on the line $\Re\lambda=-1-\kappa$, on which $\cos(\lambda\theta_b)\neq 0$, therefore vanishes; so the two functions coincide. Writing $p_{\partial_1}:=p_{\partial_1,\phi_0}+p_{\partial_1,f_\theta}$ for the data $(f_{\star,\theta},\phi_0)=(f_{\star,\theta},s\phi_\star-\dot\phi_\star)$, for which $p_{\partial_1,\phi_0}=sp_{\partial_1,\phi_\star}-p_{\partial_1,\dot\phi_\star}$ by the linearity of $\phi_0\mapsto p_{\partial_1,\phi_0}$, we have shown
\begin{align}\label{stationary pressure identity}
\lpc\mb u_\star-\grad p_{\P|\mb u_\star}=\mb F_\star+\grad(p_\star-p_{\P|\mb u_\star})=\P\mb F_\star+\grad p_{\partial_1},
\end{align}
which is the first equation of System \ref{linear system} for the data in question. The two boundary conditions $2\mb e_\theta\cdot\S(\mb u_\star)\cdot\mb e_r=u_{\star,r}$ on $\partial_0\Omega$ and $2\mb e_\theta\cdot\S(\mb u_\star)\cdot\mb e_r=f_{\star,r}$ on $\partial_1\Omega$ are the Navier slip condition and the $\mb e_r$ component of the surface tension condition of \eqref{system: stationary linear}. Consequently, for $\mb v\in R_l$, Proposition \ref{energy forms} gives \eqref{undo integration main} with $\mb u=\mb u_\star$, whose RHS vanishes by \eqref{stationary pressure identity} and the boundary conditions; this is \eqref{stationary variational identity} on $R_l$.

At $\mb v_*$ we argue as in the proof of Lemma \ref{lowering the count}: the differences $B^r_{a,j}[\mb u_\star,\mb v_*]-\mathcal F^r_{a,j}[\mb v_*]$ for $j\geq 1$ and $B^\theta_{a,j}[\mb u_\star,\mb v_*]-\mathcal F^\theta_{a,j}[\mb v_*]$ for all $j$ vanish, by the identities of Proposition \ref{energy forms} at $\mb v_*$ established there, together with \eqref{stationary pressure identity} and the boundary conditions. For $j=0$ the test function $\mb V^r_{a,0}(\mb v_*)=-r^{-2a}\mb v_*$ does not vanish at the corner, and we show $B^r_{a,0}[\mb u_\star,\mb v_*]=\mathcal F^r_{a,0}[\mb v_*]$ by running the proofs of Proposition \ref{energy forms}(i) and (iii) together on $\Omega\cap\{r>\epsilon\}$ with $\mb v=\mb v_*$ and letting $\epsilon\to 0$. Here the $\phi_0$ part of $\mathcal F^r_{0,0}$ is understood through the boundary form of the proof of (iii), $\bar s^{-1}(\sin\theta_b)^2(-\overline{\phi_0(0)}v_\theta(0,\theta_b)-\<\phi_0,\partial_rv_\theta\>_{L^2(\partial_1\Omega)})$, which is the form in which it is bounded on $H_s^l$, and $\phi_0(0)=s\phi_\star(0)$, since $\dot\phi_\star(0)=0$ by \eqref{stationary corner condition}. On $\{r>\epsilon\}$ the bulk pairings of (i) and (iii) coincide by \eqref{stationary pressure identity}, the boundary terms on $\partial_0\Omega$ and $\partial_1\Omega$ produced by the integrations by parts cancel by the boundary conditions, and the remaining integrals are absolutely convergent on $\Omega$, so that their truncations converge; two terms are left over. The integration by parts along $\partial_1\Omega$, which turns $-\bar s^{-1}\<\partial_r^2u_{\star,\theta},r^{-2a}v_{*,\theta}\>_{L^2(\partial_1\Omega)}$ into $\bar s^{-1}\<\partial_ru_{\star,\theta},\partial_r(r^{-2a}v_{*,\theta})\>_{L^2(\partial_1\Omega)}$, produces the endpoint term $\bar s^{-1}\epsilon^{-2a}\partial_ru_{\star,\theta}(\epsilon,\theta_b)v_{*,\theta}(\epsilon,\theta_b)$, which tends to $0$ by \eqref{stationary corner condition}. And the divergence theorem produces the integral over the arc $\{r=\epsilon\}$ of $\epsilon^{1-2a}(\mathbb T\mb e_r)\cdot\mb v_*$, in which the pressure functionals of (i) and the pressures $p_{\partial_1}$ and $q[\mb F_\star]$ of (iii) combine, by \eqref{stationary pressure identity}, into the pressure of the stationary pair: $\mathbb T=2\S(\mb u_\star)-(p_{\P|\mb u_\star}+p_{\partial_1}-q[\mb F_\star])I=\mathbb T[\mb u_\star,p_\star]$. This is the reason for running the two proofs together: $p_{\P|\mb u_\star}$ and $p_{\partial_1}$ separately need not be square integrable near the corner, since $p_{\partial_1,1}[\mb u_\star]=p_{\partial_1,\dot\phi_\star}$, which enters them with opposite signs, is controlled only in $\|\ph\|_{\Omega,-\kappa}$, whereas $\mathbb T[\mb u_\star,p_\star]$ is square integrable on bounded sets by the class of $(\mb u_\star,p_\star)$. The mean of the arc integral over $\epsilon\in[\epsilon_0,2\epsilon_0]$ is therefore bounded by a constant times $\epsilon_0^{-2a}\|\mathbb T[\mb u_\star,p_\star]\|_{L^2(\Omega\cap\{r<2\epsilon_0\})}$, which tends to $0$ as $\epsilon_0\to 0$. Letting $\epsilon\to 0$ along a suitable sequence gives $B^r_{a,0}[\mb u_\star,\mb v_*]=\mathcal F^r_{a,0}[\mb v_*]$ for $a\in\{0,-\kappa\}$, hence \eqref{stationary variational identity} at $\mb v_*$.
\end{proof}

The triple $(\mb u_\star,p_\star,\phi_\star)$ of Lemma \ref{stationary as variational} solves \eqref{system: linear with profile} with the data $(\mb F_\star,\mb f_\star,s\phi_\star-\dot\phi_\star)$ strongly, its last equation being the definition of $\dot\phi_\star$; the content of the lemma is that $\mb u_\star$ is the variational solution, which is where the contact angle condition enters.

Next we will prove existence of solutions to the original system in space and time, the system $\mathfrak W_t[\mb F,\mb f;\phi_0]$ defined after System \ref{reference wedge system}.

\begin{definition}[Time norms]\label{time norms}
\begin{align*}
  \nn u\nn_{[t_1,t_2],O:(\varsigma,a)(n,j,m)}^2 &= \sum_{n'=0}^n\norm{\|e^{-\varsigma t}\partial_t^{n'}u\|_{O,a,j+n-n',m}}_{L^2[t_1,t_2]}^2
\end{align*}
Here $\varsigma>0$ is an exponential weight in time, $a$ is the spatial weight, $n$ is the number of time derivatives, and $(j,m)$ are the radial and full derivative counts of Definition \ref{weighted Sobolev norms}. Entries of $(n,j,m)$ that are omitted are taken to be $0$; thus $\nn u\nn_{[t_1,t_2],O:(\varsigma,a)}=\|\|e^{-\varsigma t}u\|_{O,a}\|_{L^2[t_1,t_2]}$.
\end{definition}

The term carrying $\partial_t^{n'}$ in Definition \ref{time norms} sits at the radial count $j+n-n'$, so that time derivatives are traded against radial ones.

We first record the embedding that converts control of one time derivative into control of the supremum in time. It gives the continuity in time used below, and in Subsection \ref{positive-time construction} the bounds on suprema in time.

\begin{lemma}[Sup in time]\label{supremum in time}
Let $\varsigma,T>0$, $O\in\{\Omega,\partial_0\Omega,\partial_1\Omega\}$, $a\in\R$, $j\geq 0$ and $m\in\N_0$, and let $u$ be a function on $O\times[0,T]$ with $\nn u\nn_{[0,T],O:(\varsigma,a)(0,j,m)}$ and $\nn\partial_tu\nn_{[0,T],O:(\varsigma,a)(0,j,m)}$ finite. Then $u$ has a representative continuous in $t$ with values in the weighted space, and
\begin{align}
\sup_{t\in[0,T]}\|u(t)\|_{O,a,j,m}&\leq\|u(0)\|_{O,a,j,m}+e^{\varsigma T}T^{1\over 2}\nn\partial_tu\nn_{[0,T],O:(\varsigma,a)(0,j,m)},\label{sup in time 1}\\
\nn u\nn_{[0,T],O:(\varsigma,a)(0,j,m)}&\leq T^{1\over 2}\sup_{t\in[0,T]}\|u(t)\|_{O,a,j,m},\label{sup in time 2}\\
\sup_{t\in[0,T]}\|u(t)\|_{O,a,j,m}^2&\leq e^{2\varsigma T}\brac{\brac{T^{-1}+2\varsigma}\nn u\nn^2+2\nn u\nn\nn\partial_tu\nn},\label{sup in time 3}\\
\sup_{t\in[0,T]}\|u(t)\|_{O,a,j+{1\over 2},m}^2&\lesssim_je^{2\varsigma T}\brac{T^{-1}+2\varsigma+1}\nn u\nn_{[0,T],O:(\varsigma,a)(1,j,m)}^2,\label{sup in time 4}
\end{align}
where in \eqref{sup in time 3} the norm $\nn\ph\nn$ abbreviates $\nn\ph\nn_{[0,T],O:(\varsigma,a)(0,j,m)}$.
\end{lemma}
\begin{proof}
Write $\|\ph\|$ for the spatial norm $\|\ph\|_{O,a,j,m}$ and $\nn\ph\nn$ for the corresponding time norm. Since $\nn\partial_tu\nn<\infty$, the function $t\mapsto u(t)$ is absolutely continuous with values in the weighted space, which gives the continuous representative and lets us apply the fundamental theorem of calculus. For \eqref{sup in time 1},
\begin{align*}
\|u(t)\|\leq\|u(0)\|+\int_0^t\|\partial_tu(t')\|\d t'
\leq\|u(0)\|+\|e^{\varsigma t'}\|_{L^2[0,T]}\nn\partial_tu\nn
\leq\|u(0)\|+e^{\varsigma T}T^{1\over 2}\nn\partial_tu\nn,
\end{align*}
by the Cauchy--Schwarz inequality after inserting $e^{-\varsigma t'}e^{\varsigma t'}$. The bound \eqref{sup in time 2} is immediate from $e^{-\varsigma t}\leq 1$. For \eqref{sup in time 3}, pick $t_*\in[0,T]$ with $e^{-2\varsigma t_*}\|u(t_*)\|^2\leq T^{-1}\nn u\nn^2$, which exists by the mean value property of the integral; then for any $t$,
\begin{align*}
e^{-2\varsigma t}\|u(t)\|^2
&\leq e^{-2\varsigma t_*}\|u(t_*)\|^2+\int_0^T\abs{\partial_{t'}\brac{e^{-2\varsigma t'}\|u(t')\|^2}}\d t'\\
&\leq T^{-1}\nn u\nn^2+2\varsigma\nn u\nn^2+2\nn u\nn\nn\partial_tu\nn,
\end{align*}
and multiplying by $e^{2\varsigma t}\leq e^{2\varsigma T}$ gives \eqref{sup in time 3}. For \eqref{sup in time 4} we use that, for every count $k\geq 0$, $\|u\|_{O,a,k,m}^2$ is comparable, with constants depending only on $k$, to the sum of the integrals of $(1+|\lambda|^2)^k|\hat Y|^2$ over the Mellin line of the weight $a$, where $Y$ runs over the components of $(r\grad)^{m'}u$, $m'\leq m$, and over $\theta$ when $O=\Omega$. This follows from Definition \ref{fractional counts} and the pointwise bounds $|\lambda|^{2k'}\leq(1+|\lambda|^2)^k$ for $0\leq k'\leq k$ and $(1+|\lambda|^2)^k\leq 2^k(1+|\lambda|^{2k})$. We run the argument just given on this sum at $k=j+{1\over 2}$. The time derivative of its integrand is bounded by $2\varsigma$ times the integrand plus
\begin{align*}
2e^{-2\varsigma t}(1+|\lambda|^2)^{j+{1\over 2}}|\hat Y||\partial_t\hat Y|\leq e^{-2\varsigma t}\brac{(1+|\lambda|^2)^{j+1}|\hat Y|^2+(1+|\lambda|^2)^{j}|\partial_t\hat Y|^2},
\end{align*}
so that, with $t_*$ chosen as before but for this quantity, and since $(1+|\lambda|^2)^{j+1/2}\leq(1+|\lambda|^2)^{j+1}$,
\begin{align*}
\sup_{t\in[0,T]}e^{-2\varsigma t}\|u(t)\|_{O,a,j+{1\over 2},m}^2\lesssim_j\brac{T^{-1}+2\varsigma+1}\nn u\nn_{[0,T],O:(\varsigma,a)(0,j+1,m)}^2+\nn\partial_tu\nn_{[0,T],O:(\varsigma,a)(0,j,m)}^2.
\end{align*}
The RHS is at most $T^{-1}+2\varsigma+1$ times $\nn u\nn_{[0,T],O:(\varsigma,a)(1,j,m)}^2$, whose two summands are the two norms on the right, by Definition \ref{time norms}.
\end{proof}

\begin{theorem}[Solution of the linear system in time]\label{linear solution in time}
Let $\kappa\in(1/4,1/2)$, $\varpi\in(0,\kappa]$, $l\in\N_0$ and $\varsigma,T>0$. Suppose $\mb F:\Omega\times[0,T]\to\R^2$, $\mb f:\partial_1\Omega\times[0,T]\to\R^2$ and $\phi_0:\partial_1\Omega\to\R$ have regularity that makes the RHS of \eqref{time solution estimate} finite, and that $\phi_0(0)=0$. Then for small enough $\theta_b>0$ there exist $\phi:\partial_1\Omega\times[0,T]\to\R$ and $(\mb u,p):\Omega\times[0,T]\to\R^3$ satisfying a.e.
\begin{subequations}\label{system: time side}
\begin{alignat}{3}
\lpc\mb u-\grad p&=\mb F\qquad&&\text{on}\quad\Omega\\
\grad\cdot\mb u&=0\qquad&&\text{on}\quad\Omega\\
u_\theta&=0\qquad&&\text{on}\quad\partial_0\Omega\\
u_r-r^{-1}\partial_\theta u_r&=0\qquad&&\text{on}\quad\partial_0\Omega\\
\mathbb T[\mb u,p]\mb e_\theta-(\sin\theta_b)^2(\partial_r\phi)\mb e_\theta&=\mb f\qquad&&\text{on}\quad\partial_1\Omega\\
\partial_{t}\phi-(\sin\theta_b)^{-2}\partial_r u_\theta &=0\qquad&&\text{on}\quad\partial_1\Omega\\
\lim_{r\to\infty}\phi(r,t)&=0\qquad&&\text{for a.e.}\quad t\\
\phi(\ph,0)&=\phi_0
\intertext{and satisfies weakly}
\phi(0,t)&=0\qquad&&\text{for a.e.}\quad t
\end{alignat}
\end{subequations}
and $\partial_ru_\theta(0,\theta_b)=0$ weakly for a.e.\ $t$, with regularity
\begin{align}
&(\sin\theta_b)^2\nn\phi\nn_{[0,T],\partial_1\Omega:(\varsigma,{1\over 2}-\kappa)(0,l+{3\over 2})} +\nn\mb u\nn_{[0,T],\Omega:(\varsigma,1-\kappa)(0,l,2)}+\nn p\nn_{[0,T],\Omega:(\varsigma,-\kappa)(0,l,1)}\nonumber\\
&+\nn r\grad\mb u\nn_{[0,T],\Omega:(\varsigma,1)(0,l,1)}+\nn p\nn_{[0,T],\Omega:(\varsigma,0)(0,l,1)}\nonumber\\
&\lesssim_{\theta_b,\kappa,\varpi,l}\brac{1+\varsigma^{-1}}\Bigg(\sum_{k=0}^1\bigg(\nn\mb F\nn_{[0,T],\Omega:(\varsigma,-\kappa k-1)(0,l)} +\nn\mb f\nn_{[0,T],\partial_1\Omega:(\varsigma,-\kappa k-{1\over 2})(0,l+{1\over 2})}\bigg)\nonumber\\
&\quad+\nn\mb F\nn_{[0,T],\Omega:(\varsigma,\varpi-1)}+\nn\mb f\nn_{[0,T],\partial_1\Omega:(\varsigma,\varpi-{1\over 2})(0,{1\over 2})}+(\sin\theta_b)^2\|\phi_0\|_{\partial_1\Omega,{1\over 2}-\kappa,l+{3\over 2}}\Bigg).\label{time solution estimate}
\end{align}
We write $\mathfrak E_{n,j}[\phi,\mb u,p]$ and $\mathfrak D_{n,j}[\mb F,\mb f;\phi_0]$ for the two sides of \eqref{time solution estimate} with the time count $0$ and the derivative count $l$ replaced throughout by $n$ and $j$, so that \eqref{time solution estimate} reads $\mathfrak E_{0,l}[\phi,\mb u,p]\lesssim_{\theta_b,\kappa,\varpi,l}\mathfrak D_{0,l}[\mb F,\mb f;\phi_0]$.
\end{theorem}
\begin{proof}
Step 1 makes the solution of Corollary \ref{profile on the Laplace side} holomorphic in $s$, Step 2 bounds it on the half plane $\Re s\geq\varsigma$, and Step 3 returns to the time side by the theorems of Plancherel and Paley--Wiener.

\emph{Step 1: holomorphy in $s$.} Let $\hat{\mb F}$ and $\hat{\mb f}$ be the two-sided Laplace transform in time of $\mb F$ and $\mb f$, extended by zero outside $[0,T]$; they are entire in $s$ with values in the weighted spaces in which $\mb F$ and $\mb f$ are Bochner integrable in time. We use the solution of Corollary \ref{profile on the Laplace side} on the half plane $\Re s\geq\varsigma/2$, and need it to depend holomorphically on $s$ there, in order to apply Plancherel's theorem on the line $\Re s=\varsigma$ and the Paley--Wiener theorem below. The weight $1+|s|^{-2\kappa}$ of $B_{(l)}$ is not holomorphic in $s$, so on this half plane we replace it by the constant $W:=1+(\varsigma/2)^{-2\kappa}\geq 1+|s|^{-2\kappa}$, in $B_{(l)}$, in $\mathcal F_{(l)}$, in $\smash[t]{\trimbot{0.6pt}{\mb V^{(1)}_{(l)}}}$ and in $\|\ph\|_{H_s^l}$. Section \ref{linear existence} goes through verbatim with this constant weight. The remainder $C\epsilon^4C_\delta|s|^{-2\kappa}E_0[\mb u]$ of Proposition \ref{energy form coercivity} is absorbed into $WE_0[\mb u]$ as before, so the coercivity constant does not depend on $W$. The bounds of Proposition \ref{energy forms}(iii)--(iv) hold with $W$ and $W^{1/2}$ in place of $1+|s|^{-2\kappa}$ and $1+|s|^{-\kappa}$. In the test function problem the weight enters only as a positive constant. The solution so obtained satisfies the same strong system and the same estimates, with $W^{1/2}$ and $W$ in those two places.

On a compact subset of the half plane the norms $\|\ph\|_{H_s^l}$, which depend on $s$ only through the factor $|s|^{-1}$ in front of the $\partial_1\Omega$ term, are equivalent to one another, so that $H_s^l$ is one fixed Hilbert space $H^l$. On it the conjugates $\overline{B_{(l)}[\mb u,\mb v]}$ and $\overline{\mathcal F_{(l)}[\mb v]}$ depend on $s$ only through $s^{-1}$, carried by the surface tension term of $B^r_{a,j}$ and by $p_{\partial_1,1}$ and $p_{\partial_1,\phi_0}$, and through $\hat{\mb F}(s)$ and $\hat{\mb f}(s)$, which enter linearly after conjugation. Both conjugates are therefore holomorphic in $s$ for fixed $\mb u,\mb v$. Define the operators $A(s)$ on $H^l$ and the vectors $f(s)\in H^l$ by the Riesz representation theorem in the second slot,
\begin{align*}
\<\mb v,A(s)\mb u\>_{H^l}=\overline{B_{(l)}[\mb u,\mb v]}\qquad\text{ and }\qquad\<\mb v,f(s)\>_{H^l}=\overline{\mathcal F_{(l)}[\mb v]}\qquad\text{ for all }\mb v\in H^l.
\end{align*}
The inner product is conjugate linear in its first slot, as in \eqref{Plancherel's identity}. So both sides of the first identity are linear in $\mb u$ and conjugate linear in $\mb v$, which makes $A(s)$ a linear operator, and both sides of the second identity are conjugate linear in $\mb v$. The operators $A(s)$ are then holomorphic in the operator norm and the vectors $f(s)$ holomorphic in $H^l$, both being weakly holomorphic. Moreover, $A(s)$ is invertible with inverses bounded uniformly on compacts: $\|A(s)\mb u\|_{H^l}\|\mb u\|_{H^l}\geq|\<\mb u,A(s)\mb u\>_{H^l}|=|B_{(l)}[\mb u,\mb u]|\gtrsim\|\mb u\|_{H^l}^2$ by the coercivity, and $A(s)$ is onto because the variational problem $B_{(l)}[\mb u,\mb v]=\<g,\mb v\>_{H^l}$ for all $\mb v$ is solvable for every $g\in H^l$, by the Lax--Milgram theorem for the modified form. Conjugating the variational identity $B_{(l)}[\mb u_s,\mb v]=\mathcal F_{(l)}[\mb v]$ over $\mb v\in H^l$ gives $A(s)\mb u_s=f(s)$, so the solution $\mb u_s=A(s)^{-1}f(s)$ is holomorphic in $s$ on $\Re s>\varsigma/2$ with values in $H^l$, and so are $\hat p_\P$ and $\hat\phi$, which are obtained from $\mb u_s$, $\phi_0$ and the data by operations that are holomorphic in $s$. In particular $s\mapsto(\hat\phi,\hat{\mb u},\hat p_\P)$ is continuous, hence measurable, on the line $\Re s=\varsigma$. 

\emph{Step 2: the estimate on the Laplace side.} With this understood, by Corollary \ref{profile on the Laplace side} there exist $(\hat\phi,\hat{\mb u},\hat p_\P)$, holomorphic in $s$ on $\Re s>\varsigma/2$, satisfying strongly
\begin{subequations}\label{system: Laplace side}
\begin{alignat}{3}
\lpc\hat{\mb u}-\grad\hat p_\P&=\P\hat{\mb F}\qquad&&\text{on}\quad\Omega\\
\grad\cdot\hat{\mb u}&=0\qquad&&\text{on}\quad\Omega\\
\hat u_\theta&=0\qquad&&\text{on}\quad\partial_0\Omega\\
\hat u_r-r^{-1}\partial_\theta\hat u_r&=0\qquad&&\text{on}\quad\partial_0\Omega\\
\mathbb T[\hat{\mb u},\hat p_\P]\mb e_\theta-(\sin\theta_b)^2(\partial_r\hat\phi)\mb e_\theta&=\hat{\mb f}\qquad&&\text{on}\quad\partial_1\Omega\\
s\hat\phi-(\sin\theta_b)^{-2}\partial_r\hat u_\theta &=\phi_0\qquad&&\text{on}\quad\partial_1\Omega\\
\intertext{and satisfying weakly}
\hat\phi(0;s)&=0\qquad&&\text{for all}\quad s\\
\partial_r\hat u_\theta(0,\theta_b)&=0.
\end{alignat}
\end{subequations}
The factors $W^{1/2}$ and $W$, which replace $1+|s|^{-\kappa}$ and $1+|s|^{-2\kappa}$ in the constants of Corollary \ref{profile on the Laplace side}, are bounded by
\begin{align}\label{W in terms of tau}
W^{1\over 2}\leq 1+\brac{\varsigma\over 2}^{-\kappa}\leq\sqrt 2\brac{1+\varsigma^{-\kappa}}\qquad\text{ and }\qquad W=1+\brac{\varsigma\over 2}^{-2\kappa}\leq 3\brac{1+\varsigma^{-1}},
\end{align}
the former because $\sqrt{1+x^2}\leq 1+x$ and $2^\kappa\leq\sqrt 2$, the latter because $2\kappa\leq 1$ gives $\varsigma^{-2\kappa}\leq 1+\varsigma^{-1}$. So, when $\Re s\geq\varsigma$,
\begin{align*}
&\|\hat{\mb u}\|_{\Omega,1-\kappa,l,2}+\|r\grad\hat{\mb u}\|_{\Omega,1,l,1}+(\sin\theta_b)^2\|\hat\phi\|_{\partial_1\Omega,{1\over 2}-\kappa,l+{3\over 2}}+\sum_{k=0}^1\|\hat p_\P\|_{\Omega,-\kappa k,l,1}\\
&\lesssim_{\theta_b}\brac{1+\varsigma^{-1}}\brac{\sum_{k=0}^1\brac{\|\P\hat{\mb F}\|_{\Omega,-\kappa k-1,l}+\|\hat{\mb f}\|_{\partial_1\Omega,-\kappa k-{1\over 2},l+{1\over 2}}}+\|\P\hat{\mb F}\|_{\Omega,\varpi-1}+\|\hat{\mb f}\|_{\partial_1\Omega,\varpi-{1\over 2},{1\over 2}}}\\ 
&\quad+(\sin\theta_b)^2 \brac{\|s^{-1}\phi_0\|_{\partial_1\Omega,{1\over 2}-\kappa,l+{3\over 2}}+\brac{1+\varsigma^{-\kappa}}\|s^{\kappa-1}\phi_0\|_{\partial_1\Omega,{1\over 2}-\kappa,l+{3\over 2}}}.
\end{align*}
Here the dependence on $\varpi$ is absorbed into $\lesssim_{\theta_b}$. The profile terms are those of Corollary \ref{profile on the Laplace side} at $\alpha=1-2\kappa$, for which $\alpha/2={1\over 2}-\kappa$ and $s^{-(1+\alpha)/2}=s^{\kappa-1}$, the hypothesis $\phi_0(0)=0$ of that corollary being the one assumed here. The factor $W^{1/2}\lesssim 1+\varsigma^{-\kappa}$ of the profile term with $s^{\kappa-1}$ is kept explicit, and the factors $s^{-1}$ and $s^{\kappa-1}$ are converted into functions of time by Plancherel's theorem below. Let $\hat p=\hat p_\P-q[\hat{\mb F}]$; then we have $\lpc\hat{\mb u}-\grad\hat p =\hat{\mb F}$. Using Lemma \ref{projection norm estimate} we have
\begin{align}
&\|\hat{\mb u}\|_{\Omega,1-\kappa,l,2}+\|r\grad\hat{\mb u}\|_{\Omega,1,l,1}+(\sin\theta_b)^2\|\hat\phi\|_{\partial_1\Omega,{1\over 2}-\kappa,l+{3\over 2}}+\sum_{k=0}^1\|\hat p\|_{\Omega,-\kappa k,l,1}\nonumber\\
&\lesssim_{\theta_b}\brac{1+\varsigma^{-1}}\brac{\sum_{k=0}^1\brac{\|\P\hat{\mb F}\|_{\Omega,-\kappa k-1,l}+\|\hat{\mb f}\|_{\partial_1\Omega,-\kappa k-{1\over 2},l+{1\over 2}}}+\|\P\hat{\mb F}\|_{\Omega,\varpi-1}+\|\hat{\mb f}\|_{\partial_1\Omega,\varpi-{1\over 2},{1\over 2}}}\nonumber\\ 
&\quad+(\sin\theta_b)^2 \brac{\|s^{-1}\phi_0\|_{\partial_1\Omega,{1\over 2}-\kappa,l+{3\over 2}}+\brac{1+\varsigma^{-\kappa}}\|s^{\kappa-1}\phi_0\|_{\partial_1\Omega,{1\over 2}-\kappa,l+{3\over 2}}}.\label{Resolvent estimate for system}
\end{align}

\emph{Step 3: back to the time side.} By Plancherel's theorem on the line $\Re s=\varsigma$, on which the transforms are square integrable by the estimate just displayed and the finiteness of $\mathfrak D_{0,l}[\mb F,\mb f;\phi_0]$, we get 
\begin{align*}
&\nn\mb u\nn_{[0,T],\Omega:(\varsigma,1-\kappa)(0,l,2)}+\nn r\grad\mb u\nn_{[0,T],\Omega:(\varsigma,1)(0,l,1)}\\
&\quad+(\sin\theta_b)^2\nn\phi\nn_{[0,T],\partial_1\Omega:(\varsigma,{1\over 2}-\kappa)(0,l+{3\over 2})}+\sum_{k=0}^1\nn p\nn_{[0,T],\Omega:(\varsigma,-\kappa k)(0,l,1)}\\
&\lesssim_{\theta_b}\brac{1+\varsigma^{-1}}\Bigg(\sum_{k=0}^1\bigg(\nn\mb F\nn_{[0,T],\Omega:(\varsigma,-\kappa k-1)(0,l)}+\nn\mb f\nn_{[0,T],\partial_1\Omega:(\varsigma,-\kappa k-{1\over 2})(0,l+{1\over 2})}\bigg)\\
&\quad+\nn\mb F\nn_{[0,T],\Omega:(\varsigma,\varpi-1)}+\nn\mb f\nn_{[0,T],\partial_1\Omega:(\varsigma,\varpi-{1\over 2})(0,{1\over 2})}\Bigg)\\ 
&\quad+(\sin\theta_b)^2\nn\phi_0\nn_{[0,\infty),\partial_1\Omega:(\varsigma,{1\over 2}-\kappa)(0,l+{3\over 2})}\\
&\quad+(\sin\theta_b)^2\brac{1+\varsigma^{-\kappa}}\nn\Gamma(1-\kappa)^{-1}t^{-\kappa}\phi_0\nn_{[0,\infty),\partial_1\Omega:(\varsigma,{1\over 2}-\kappa)(0,l+{3\over 2})}\\
&\lesssim_{\theta_b}\brac{1+\varsigma^{-1}}\Bigg(\sum_{k=0}^1\bigg(\nn\mb F\nn_{[0,T],\Omega:(\varsigma,-\kappa k-1)(0,l)}+\nn\mb f\nn_{[0,T],\partial_1\Omega:(\varsigma,-\kappa k-{1\over 2})(0,l+{1\over 2})}\bigg)\\
&\quad+\nn\mb F\nn_{[0,T],\Omega:(\varsigma,\varpi-1)}+\nn\mb f\nn_{[0,T],\partial_1\Omega:(\varsigma,\varpi-{1\over 2})(0,{1\over 2})}+(\sin\theta_b)^2\|\phi_0\|_{\partial_1\Omega,{1\over 2}-\kappa,l+{3\over 2}}\Bigg),
\end{align*}
where the transforms of the data extended by zero are taken over $[0,\infty)$. The two functions of time in the profile terms, the inverse transforms of $s^{-1}$ and $s^{\kappa-1}$ on $\{t>0\}$, have been integrated explicitly:
\begin{align*}
\|e^{-\varsigma t}\|_{L^2(0,\infty)}&=(2\varsigma)^{-{1\over 2}}\leq 1+\varsigma^{-1},\\
\brac{1+\varsigma^{-\kappa}}{\|t^{-\kappa}e^{-\varsigma t}\|_{L^2(0,\infty)}\over\Gamma(1-\kappa)}&={\Gamma(1-2\kappa)^{1\over 2}\over\Gamma(1-\kappa)}\brac{1+\varsigma^{-\kappa}}(2\varsigma)^{\kappa-{1\over 2}}\lesssim_\kappa 1+\varsigma^{-1},
\end{align*}
the last step because $(1+\varsigma^{-\kappa})\varsigma^{\kappa-1/2}$ is at most $2$ for $\varsigma\geq 1$ and at most $2\varsigma^{-1/2}\leq 2\varsigma^{-1}$ for $\varsigma<1$. The norm of $t^{-\kappa}e^{-\varsigma t}$ is computed from $\int_0^\infty t^{-2\kappa}e^{-2\varsigma t}\d t=\Gamma(1-2\kappa)(2\varsigma)^{2\kappa-1}$, which is finite since $\kappa<{1\over 2}$. Hence we have \eqref{time solution estimate}. Taking inverse Laplace transform to \eqref{system: Laplace side} we get that $(\phi,\mb u,p)$ solves $\mathfrak W_t[\mb F,\mb f;\ph]$ with $\phi_0\delta(t)$ on the RHS of the profile evolution equation, and satisfies $\partial_ru_\theta(0,\theta_b)=0$ weakly.
Since $\mb F,\mb f$ and $\phi_0\delta(t)$ are causal (i.e.\ they vanish for negative time), so are $\mb u$, $p$ and $\phi$, by the Paley--Wiener theorem for causal functions \cite[Theorem~19.2]{Rudin}. It applies to the transforms $(\hat\phi,\hat{\mb u},\hat p)$, which are holomorphic on $\Re s>\varsigma/2$ with the estimate \eqref{Resolvent estimate for system} uniform on the lines $\Re s=\varsigma'\geq\varsigma$. The theorem is stated for scalar functions, and we apply it to the pairings of the transforms with fixed elements of the weighted spaces. By Lemma \ref{trace}, $\|r^{1/2}\partial_ru_\theta\|_{L^2(\partial_1\Omega)}\lesssim_{\theta_b}\|r\grad\mb u\|_{\Omega,1,0,1}$, so \eqref{time solution estimate} makes $r^{1/2}\partial_ru_\theta$ square integrable on $\partial_1\Omega\times[0,T]$. Since $\phi$ vanishes for negative time, the profile evolution equation with $\phi_0\delta(t)$ on its RHS then gives $\phi(t)=\phi_0+(\sin\theta_b)^{-2}\int_0^t\partial_ru_\theta(t')\d t'$ on $\partial_1\Omega$. So $\phi(t)-\phi_0$ is continuous in $t$ in the norm $\|\ph\|_{\partial_1\Omega,-1/2}$ and vanishes at $t=0$, which is the initial condition $\phi(\ph,0)=\phi_0$.

It remains to check the condition at infinity. By \eqref{time solution estimate} the norm $\|\phi(t)\|_{\partial_1\Omega,1/2-\kappa,l+3/2}$ is finite for a.e.\ $t\in[0,T]$, so that Lemma \ref{embeddings}(i) on the ray $\partial_1\Omega$, at $a=-\kappa$, gives
\begin{align}\label{profile decay}
\|r^{\kappa}\phi(t)\|_{L^\infty(\partial_1\Omega)}\lesssim_\kappa\|\phi(t)\|_{\partial_1\Omega,{1\over 2}-\kappa,1}
\end{align}
at every such $t$; since $\kappa>0$, this gives $\lim_{r\to\infty}\phi(r,t)=0$ for a.e.\ $t$, which is the condition at infinity. In the classes of Theorem \ref{linear solution full regularity} and of Theorem \ref{nonlinear positive-time existence}, where $\nn\partial_t\phi\nn_{[0,T],\partial_1\Omega:(\varsigma,1/2-\kappa)(0,1)}$ is finite, $t\mapsto\phi(t)$ has a representative continuous in $t$ in the norm $\|\ph\|_{\partial_1\Omega,1/2-\kappa,1}$, by Lemma \ref{supremum in time}; as \eqref{profile decay} is a bounded embedding, $t\mapsto r^\kappa\phi(\ph,t)$ is then continuous with values in $L^\infty(\partial_1\Omega)$, and the condition holds at every $t\in[0,T]$. Hence we get \eqref{system: time side}.
\end{proof}

The next lemma shows that there is only one solution with $\mathfrak E_{0,0}[\phi,\mb u,p]$ and $\nn\phi\nn_{[0,T],\partial_1\Omega:(\varsigma,1/2)(0,3/2)}$ finite, the latter being the quantitative form of the contact angle condition of Remark \ref{weight half remark}. That norm is supplied by Theorem \ref{linear solution full regularity} from its rung $M=1$, and the lemma is used in Step 4 of the proof of Theorem \ref{nonlinear positive-time existence}.

\begin{lemma}[Uniqueness for the linear evolution]\label{linear evolution uniqueness}
Let $\kappa\in(1/4,1/2)$ and $\varsigma,T>0$. Suppose $(\phi,\mb u,p)$ and $(\phi',\mb u',p')$ both satisfy \eqref{system: time side} a.e.\ on $[0,T]$ with the same data $(\mb F,\mb f,\phi_0)$, and that both make the six norms
\begin{align}\label{uniqueness class}
&\nn\phi\nn_{[0,T],\partial_1\Omega:(\varsigma,{1\over 2}-\kappa)(0,{3\over 2})},\qquad
\nn\phi\nn_{[0,T],\partial_1\Omega:(\varsigma,{1\over 2})(0,{3\over 2})},\qquad
\nn\mb u\nn_{[0,T],\Omega:(\varsigma,1-\kappa)(0,0,2)},\nonumber\\
&\nn r\grad\mb u\nn_{[0,T],\Omega:(\varsigma,1)(0,0,1)},\qquad
\nn p\nn_{[0,T],\Omega:(\varsigma,-\kappa)(0,0,1)},\qquad
\nn p\nn_{[0,T],\Omega:(\varsigma,0)(0,0,1)}
\end{align}
finite. Then $(\phi,\mb u,p)=(\phi',\mb u',p')$ a.e.
\end{lemma}
\begin{proof}
The difference, still denoted $(\phi,\mb u,p)$, satisfies \eqref{system: time side} with zero data and makes \eqref{uniqueness class} finite; we show that it vanishes. The class gives, at almost every time, the functions
\begin{align}\label{class functions}
r^{\kappa-1}\mb u,\qquad r^{\kappa}\grad\mb u,\qquad r^{1+\kappa}\grad^2\mb u,\qquad \grad\mb u,\qquad r\grad^2\mb u,\qquad r^{\kappa}p,\qquad p
\end{align}
in $L^2(\Omega)$ and $r^{\kappa-1/2}\phi$, $r^{-1/2}\phi$ in $L^2(\partial_1\Omega)$, with norms bounded by a function $m(t)\in L^2[0,T]$. For $\rho>0$ write $\nn g\nn_{\Omega_{\rho,2\rho}}$ for the norm of $g$ in $L^2(\Omega_{\rho,2\rho}\times[0,T])$, with the measure $r\d r\d\theta\d t$, and likewise on the two boundary pieces with $\d r\d t$; whenever $g$ is built from \eqref{class functions} and the two boundary functions, $\nn g\nn_{\Omega_{\rho,2\rho}}$ is bounded uniformly in $\rho$ and tends to $0$ as $\rho\to 0$ and as $\rho\to\infty$, by dominated convergence.

\emph{Step 1: two corner bounds.} Split $\mb u=\mb A+\mb w$, where $\mb A(r,t):=\theta_b^{-1}\int_0^{\theta_b}\mb u(r,\theta,t)\d\theta$ is the angular average. The fundamental theorem of calculus in $\theta$ and $|\partial_\theta\mb u|\leq r|\grad\mb u|$ give
\begin{align}\label{oscillation bound}
|\mb w(r,\theta,t)|\leq\theta_b^{1\over 2}r\brac{\int_0^{\theta_b}|\grad\mb u(r,\theta',t)|^2\d\theta'}^{1\over 2},
\end{align}
while for the average, comparing with its mean over $\Omega_{1,2}$ and integrating $\partial_r\mb A$ radially,
\begin{align}\label{average bound}
\sup_{\rho\leq r\leq 2}|\mb A(r,t)|\lesssim_{\theta_b}\brac{1+|\ln\rho|^{1\over 2}}m(t)\qquad(\rho\leq 1),
\end{align}
since $\int_\rho^2\int_0^{\theta_b}|\grad\mb u|\d\theta\d r\leq(\theta_b\ln(2/\rho))^{1/2}\|\grad\mb u\|_{L^2(\Omega\cap\{r<2\})}$ by the Cauchy--Schwarz inequality with the measure $r\d r\d\theta$.

\emph{Step 2: the truncated energy identity.} Let $0<\epsilon<1<R$. At almost every time the pair $(\mb u,p)$ lies in $H^2\times H^1$ of $\Omega_{\epsilon,R}$, so from $\grad\cdot \mathbb T[\mb u,p]=\lpc\mb u-\grad p=\mb 0$ and $\grad\cdot\mb u=0$, Green's identity gives
\begin{align*}
0=\int_{\Omega_{\epsilon,R}}\brac{\grad\cdot\mathbb T}\cdot\mb u=-2\int_{\Omega_{\epsilon,R}}|\S(\mb u)|^2+\int_{\partial\Omega_{\epsilon,R}}\brac{\mathbb T\mb n}\cdot\mb u.
\end{align*}
On $\partial_0\Omega$ the outward normal is $-\mb e_\theta$; there $u_\theta=0$, hence also $\partial_ru_\theta=0$ along the ray, so $(\mathbb T\mb e_\theta)\cdot\mb u=2(\mb e_r\cdot\S(\mb u)\cdot\mb e_\theta)u_r=(r^{-1}\partial_\theta u_r)u_r=u_r^2$ by the Navier slip condition. On $\partial_1\Omega$ the normal is $\mb e_\theta$ and $\mathbb T\mb e_\theta=(\sin\theta_b)^2(\partial_r\phi)\mb e_\theta$, so $(\mathbb T\mb n)\cdot\mb u=(\sin\theta_b)^2(\partial_r\phi)u_\theta$ there, and integrating by parts in $r$ and inserting the profile evolution equation $\partial_ru_\theta=(\sin\theta_b)^2\partial_t\phi$,
\begin{align*}
\int_\epsilon^R(\partial_r\phi)u_\theta\d r=\eva{\phi u_\theta}_\epsilon^R-(\sin\theta_b)^2\int_\epsilon^R\phi\partial_t\phi\d r.
\end{align*}
Writing $\Phi(\rho,t):=\rho\int_0^{\theta_b}((\mathbb T\mb e_r)\cdot\mb u)(\rho,\theta,t)\d\theta$ for the flux through the arc $\{r=\rho\}$ and collecting terms,
\begin{align}\label{truncated identity}
2\int_{\Omega_{\epsilon,R}}|\S(\mb u)|^2+\int_\epsilon^Ru_r^2(r,0,t)\d r+{(\sin\theta_b)^4\over 2}\partial_t\int_\epsilon^R\phi^2\d r
=\Phi(R,t)-\Phi(\epsilon,t)+(\sin\theta_b)^2\eva{\phi u_\theta}_\epsilon^R.
\end{align}
The function $t\mapsto\int_\epsilon^R\phi^2\d r$ is absolutely continuous with the derivative displayed, $\partial_t\phi$ being the trace of $(\sin\theta_b)^{-2}\partial_ru_\theta$, which lies in $L^2([\epsilon,R]\times[0,T])$, and it vanishes at $t=0$. Integrating \eqref{truncated identity} over $[0,T']$ for $T'\leq T$ therefore gives
\begin{align}
&\int_0^{T'}\brac{2\int_{\Omega_{\epsilon,R}}|\S(\mb u)|^2+\int_\epsilon^Ru_r^2\d r}\d t+{(\sin\theta_b)^4\over 2}\int_\epsilon^R\phi^2(r,T')\d r\nonumber\\
&\qquad=\int_0^{T'}\brac{\Phi(R,t)-\Phi(\epsilon,t)+(\sin\theta_b)^2\eva{\phi u_\theta}_\epsilon^R}\d t.\label{integrated identity}
\end{align}

\emph{Step 3: the RHS vanishes along sequences.} We bound averages over dyadic annuli of the time integral of the RHS, with $|\mathbb T|\leq 2|\grad\mb u|+|p|$ throughout. At infinity, matched weights gain powers of $R$: from $|\mathbb T||\mb u|\d\theta\d\rho=(r^\kappa|\mathbb T|)(r^{\kappa-1}|\mb u|)r^{-2\kappa}r\d\rho\d\theta$,
\begin{align*}
\int_R^{2R}\int_0^{T'}|\Phi(\rho,t)|\d t{\d\rho\over\rho}
\leq\int_0^{T'}\int\!\!\int_{\Omega_{R,2R}}|\mathbb T||\mb u|\d\theta\d\rho\d t
\leq R^{-2\kappa}\nn r^\kappa\mathbb T\nn_{\Omega_{R,2R}}\nn r^{\kappa-1}\mb u\nn_{\Omega_{R,2R}}\to 0,
\end{align*}
and, by the fundamental theorem of calculus in $\theta$ as in Step 1 together with the exponent identities $r^{-1-2\kappa}=(r^{\kappa-1})^2r^{-4\kappa}r$ and $r^{-1-2\kappa}\ph r=(r^{\kappa-1})(r^{\kappa})r^{-4\kappa}r$,
\begin{align*}
\int_0^{T'}\int_R^{2R}r^{-1-2\kappa}u_\theta^2(r,\theta_b,t)\d r\d t
\lesssim_{\theta_b}R^{-4\kappa}\brac{\nn r^{\kappa-1}\mb u\nn_{\Omega_{R,2R}}^2+\nn r^{\kappa-1}\mb u\nn_{\Omega_{R,2R}}\nn r^{\kappa}\grad\mb u\nn_{\Omega_{R,2R}}},
\end{align*}
so that, pairing $|\phi u_\theta|\rho^{-1}=(r^{\kappa-1/2}|\phi|)(r^{-\kappa-1/2}|u_\theta|)$,
\begin{align*}
\int_R^{2R}\int_0^{T'}|\phi u_\theta|(\rho,t)\d t{\d\rho\over\rho}
\leq\nn r^{\kappa-{1\over 2}}\phi\nn_{\Omega_{R,2R}}\brac{\int_0^{T'}\int_R^{2R}r^{-1-2\kappa}u_\theta^2\d r\d t}^{1\over 2}\to 0.
\end{align*}
At the corner, split $\mb u=\mb A+\mb w$. The oscillation pairs with no loss: by \eqref{oscillation bound} and the Cauchy--Schwarz inequality on the arcs,
\begin{align*}
\int_0^{T'}\int\!\!\int_{\Omega_{\epsilon,2\epsilon}}|\mathbb T||\mb w|\d\theta\d\rho\d t&\leq\theta_b\nn\mathbb T\nn_{\Omega_{\epsilon,2\epsilon}}\nn\grad\mb u\nn_{\Omega_{\epsilon,2\epsilon}}\to 0,\\
\int_0^{T'}\int_\epsilon^{2\epsilon}|\phi||\mb w|(\rho,t)\d t{\d\rho\over\rho}&\leq\theta_b^{1\over 2}\nn r^{-{1\over 2}}\phi\nn_{\Omega_{\epsilon,2\epsilon}}\nn\grad\mb u\nn_{\Omega_{\epsilon,2\epsilon}}\to 0.
\end{align*}
The average carries the half logarithm of \eqref{average bound}, which square summability over dyadic scales beats. Let
\begin{align*}
a_j^2:=\nn\grad\mb u\nn_{\Omega_{2^{-j},2^{1-j}}}^2+\nn p\nn_{\Omega_{2^{-j},2^{1-j}}}^2+\nn r^{-{1\over 2}}\phi\nn_{\Omega_{2^{-j},2^{1-j}}}^2,
\end{align*}
so that $\sum_ja_j^2<\infty$ and hence $\liminf_{j\to\infty}ja_j^2=0$, since otherwise the sum would diverge; along scales $\epsilon_j=2^{-j}$ realising this liminf, using \eqref{average bound} and the Cauchy--Schwarz inequality in $t$,
\begin{align*}
\int_0^{T'}\int\!\!\int_{\Omega_{\epsilon_j,2\epsilon_j}}|\mathbb T||\mb A|\d\theta\d\rho\d t+\int_0^{T'}\int_{\epsilon_j}^{2\epsilon_j}|\phi||\mb A|(\rho,t)\d t{\d\rho\over\rho}
\lesssim_{\theta_b}\brac{1+j^{1\over 2}}\|m\|_{L^2[0,T]}a_j\to 0.
\end{align*}
A nonnegative function on $[\rho_0,2\rho_0]$ takes a value at most $(\ln 2)^{-1}$ times its integral against $\d\rho/\rho$. Hence the bounds above give radii $\rho_j^-\in[\epsilon_j,2\epsilon_j]$ and $\rho_j^+\in[2^j,2^{j+1}]$ at which the sum of the time integrals of $|\Phi(\rho,\ph)|$ and $|\phi u_\theta|(\rho,\ph)$ tends to $0$ as $j\to\infty$. So the RHS of \eqref{integrated identity} with $\epsilon=\rho_j^-$ and $R=\rho_j^+$ tends to $0$.

\emph{Step 4: conclusion.} The LHS of \eqref{integrated identity} is nonnegative and nondecreasing as $\epsilon$ decreases and $R$ increases, so by Step 3 it vanishes for every $\epsilon,R$ and almost every $T'$. Hence $\phi=0$ a.e., $u_r|_{\partial_0\Omega}=0$ a.e., and $\S(\mb u(t))=0$ for a.e.\ $t$, so that $\mb u(t)$ is a rigid field $\mb a(t)+b(t)(-x_2,x_1)$; but no nonzero rigid field has finite $\|\ph\|_{\Omega,1-\kappa}$ norm, the integral $\int^\infty r^{2\kappa-1}\d r$ diverging, so $\mb u=\mb 0$. Then $\grad p=\lpc\mb u=\mb 0$, so $p$ is constant in space, and a nonzero constant has infinite $\|\ph\|_{\Omega,0}$ norm; so $p=0$ as well.
\end{proof}

The next lemma shows that the time derivative of the solution of Theorem \ref{linear solution in time} is the solution with the time derivatives of the data and with the initial profile derived from the stationary solution at $t=0$. It is applied twice in the proof of Theorem \ref{linear solution full regularity}.

\begin{lemma}[The time derivative of the solution]\label{time derivative of the solution}
Take the setting of Theorem \ref{linear solution in time} with $l+1$ in place of $l$, and let $(\phi,\mb u,p)$ be its solution. Let $(\mb u_0,p_0)$ satisfy the hypotheses of Lemma \ref{stationary as variational} at the count $l$, for every $s$ with $\Re s=\varsigma$, with the data $(\phi_0,\mb F(0),\mb f(0))$, and let $\dot\phi_0:=(\sin\theta_b)^{-2}\partial_ru_{0,\theta}$ on $\partial_1\Omega$. Suppose that Theorem \ref{linear solution in time} applies at the count $l$ to the data $(\partial_t\mb F,\partial_t\mb f,\dot\phi_0)$, and let $(\dot\phi,\dot{\mb u},\dot p)$ be its solution. Then $(\partial_t\phi,\partial_t\mb u,\partial_tp)=(\dot\phi,\dot{\mb u},\dot p)$ and $(\mb u,p)(0)=(\mb u_0,p_0)$.
\end{lemma}
\begin{proof}
Step 1 integrates $(\dot\phi,\dot{\mb u},\dot p)$ in time from $(\phi_0,\mb u_0,p_0)$, and Step 2 identifies the result with $(\phi,\mb u,p)$ by comparing Laplace transforms.

\emph{Step 1: the integrated triple.} Let
\begin{align*}
\check{\mb u}(t) = \mb u_0+\int_0^t\dot{\mb u}\,\d t',\qquad\quad
\check p(t) = p_0+\int_0^t\dot p\,\d t',\qquad\quad
\check\phi(t) = \phi_0+\int_0^t\dot\phi\,\d t'.
\end{align*}
Integrating the system \eqref{system: time side} of $(\dot\phi,\dot{\mb u},\dot p)$ in time and adding the stationary system \eqref{system: stationary linear} of $(\mb u_0,p_0)$, we get that $(\check\phi,\check{\mb u},\check p)$ solves \eqref{system: time side} with the data $(\mb F,\mb f,\phi_0)$, its profile evolution equation reading $\dot\phi-(\sin\theta_b)^{-2}\partial_r\check u_\theta=\partial_t\check\phi-(\sin\theta_b)^{-2}\partial_r\check u_\theta=0$ on $\partial_1\Omega$.

\emph{Step 2: $(\check\phi,\check{\mb u},\check p)=(\phi,\mb u,p)$.} By Step 1 both triples solve \eqref{system: time side} with the data $(\mb F,\mb f,\phi_0)$, so the claim is one of uniqueness. Lemma \ref{linear evolution uniqueness} does not give it, since its class contains $\nn\phi\nn_{[0,T],\partial_1\Omega:(\varsigma,1/2)(0,3/2)}$, the quantitative form of the contact angle condition, which neither triple is known to have. By Remark \ref{weight half remark} this norm is obtained from the rung $M=1$ of \eqref{Improved solution estimate}, through Proposition \ref{time for radial trade}, which needs control of $\partial_t\mb u$ and $\partial_tp$, available only through the present lemma. Differentiating \eqref{system: time side} in time would ask for uniqueness in a class containing $\partial_t\mb u$ and $\partial_tp$, on which Theorem \ref{linear solution in time} gives no bound. We therefore use the uniqueness of the variational solution of Theorem \ref{linear solution Laplace side}, in which the contact angle condition is a natural boundary condition (Remark \ref{corner condition and uniqueness}), and compare the Laplace transforms of $\mb u$ and $\check{\mb u}$ on the line $\Re s=\varsigma$.

The construction of Theorem \ref{linear solution in time} extends the data by zero beyond $T$, so that the Laplace transform of the extended $\partial_t\mb F$ carries the boundary term $e^{-sT}\mb F(T)$. For the comparison we therefore extend $(\mb F,\mb f)$ beyond $T$ to functions vanishing for $t\geq 2T$, by an extension operator that keeps the data norms of $(\mb F,\mb f)$ and $(\partial_t\mb F,\partial_t\mb f)$ finite, and consider the three triples on $[0,\infty)$, the dotted one with the data $(\partial_t\mb F,\partial_t\mb f)$ of the extension. On $[0,T]$ this changes none of the three triples, since a solution of Theorem \ref{linear solution in time} with $\phi_0=0$ and data vanishing on $[0,T]$ vanishes there. Indeed, the Laplace transforms of such data are $e^{-sT}$ times those of causal data, hence so is that of the solution, by linearity at each $s$, and the Paley--Wiener argument of the proof of Theorem \ref{linear solution in time} applies after a translation by $T$.

For $\Re s=\varsigma$ let $\mb u_s$ be the Laplace transform of $\mb u$. By construction, $\mb u_s$ is the solution of Theorem \ref{linear solution Laplace side} at the count $l+1$ with data $(\hat{\mb F}(s),\hat{\mb f}(s),\phi_0)$, for the form with the constant weight $W$ of the proof of Theorem \ref{linear solution in time}. Lemmas \ref{lowering the count} and \ref{stationary as variational} apply verbatim to this form, since nothing in their proofs depends on the value of the weight. Let $\dot{\mb u}^\wedge(s)$ be the Laplace transform of $\dot{\mb u}$, the solution at the count $l$ with data
\begin{align*}
\brac{(\partial_t\mb F)^\wedge(s),(\partial_t\mb f)^\wedge(s),\dot\phi_0}=\brac{s\hat{\mb F}(s)-\mb F(0),s\hat{\mb f}(s)-\mb f(0),\dot\phi_0},
\end{align*}
the data being continuous in time with values in the weighted spaces by Lemma \ref{supremum in time}. The Laplace transform of $\check{\mb u}$ is $s^{-1}(\mb u_0+\dot{\mb u}^\wedge(s))$.

By Lemma \ref{stationary as variational}, $\mb u_0$ is the solution at the count $l$ with data $(\mb F(0),\mb f(0),s\phi_0-\dot\phi_0)$. By the linearity of the solution operator of Theorem \ref{linear solution Laplace side}, $\mb u_0+\dot{\mb u}^\wedge(s)$ is then the solution at the count $l$ with data $s(\hat{\mb F}(s),\hat{\mb f}(s),\phi_0)$, that is, $s$ times the solution at the count $l$ with data $(\hat{\mb F}(s),\hat{\mb f}(s),\phi_0)$, which is $\mb u_s$ by Lemma \ref{lowering the count}. Hence $s^{-1}(\mb u_0+\dot{\mb u}^\wedge(s))=\mb u_s$ on the line $\Re s=\varsigma$, and $\check{\mb u}=\mb u$.

Then $\grad(\check p-p)=\lpc(\check{\mb u}-\mb u)=\mb 0$ at a.e.\ time, so $\check p-p$ is constant in space. It lies in $L^2(\Omega)$ by the weight $0$ of \eqref{time solution estimate} and of \eqref{data estimate}, so it vanishes. Moreover $\partial_t(\check\phi-\phi)=(\sin\theta_b)^{-2}\partial_r(\check u_\theta-u_\theta)=0$ with $(\check\phi-\phi)(0)=0$, so $\check\phi=\phi$. Thus $(\check\phi,\check{\mb u},\check p)=(\phi,\mb u,p)$, so that $(\partial_t\phi,\partial_t\mb u,\partial_tp)=(\dot\phi,\dot{\mb u},\dot p)$ and $(\mb u,p)(0)=(\mb u_0,p_0)$.
\end{proof}

\section{Corner expansions and improved regularity}\label{corner expansions}

In this section we will improve the regularity of solutions from Theorem \ref{linear solution in time}. The strategy is as follows. Heuristically, we have the relation $u\sim\phi$ from the surface tension equation of \eqref{system: time side} on the upper boundary $\partial_1\Omega$; and we have $\partial_t\phi\sim\partial_ru$ from the profile evolution equation of the same system. This gives us $\partial_t\phi\sim\partial_r\phi$ and $\partial_r u\sim\partial_t u$. So if we have a time derivative of $u$, we can use it to trade for a higher spatial scaling, in the sense of Definition \ref{scaling definition}. But in order to have $\partial_t u$, we need to solve for the time differentiated system with initial data $\partial_t\phi(0)$. This datum is not given; we relate it back to the initial data $\phi_0$ by solving another Stokes problem: by the profile evolution equation, $\partial_t\phi(0)=(\sin\theta_b)^{-2}\partial_ru_\theta(0)$ on $\partial_1\Omega$, where $\mb u(0)$ solves the stationary system of Proposition \ref{stationary system} with profile datum $\phi_0$. But due to the relation $\partial_t\phi\sim\partial_r\phi$, we need to assume control of the initial profile $\phi_0$ in a higher scaling norm. The next subsection, ``Initial-time regularity'', deals with this last point about the initial data, and the subsection after that deals with the gain for $u$ using $\partial_tu$ at positive time.

\subsection{Initial-time regularity}\label{initial-time regularity}

Our strategy for gaining spatial regularity is to rewrite our system with the higher spatial scaling terms on the LHS and the lower scaling ones on the RHS, which we treat as data. With all the terms on the LHS having the same spatial scaling, we can then obtain solution representations in Mellin variables. The LHS has a gain of $1$ in scaling compared with the RHS. Then if originally we have control of our solution in scaling $a$, via this solution representation, we get control of the solution in scaling $a+1$.

\subsubsection{The Mellin representation at the initial time}\label{Mellin representation initial time}

The Mellin solution representations are made up of trigonometric terms; the following lemma gives us estimates for them.

\begin{lemma}[Trigonometric estimates for the Green's function]\label{Green trigonometric bounds}
Let $C_{\rm Re},C_{\rm Im}>0$. Suppose $|\Re\lambda_{\pm 1}\theta_b|<C_{\rm Re}$ and $0\leq\theta\leq\theta_b\in(0,\pi)$. 
When $|\Im\lambda\theta_b|<C_{\rm Im}$ we have, with the quotient in \eqref{kernel estimates 1} understood as its limit at $\lambda=\pm 1$,
\begin{align}
\abs{\sin(\lambda_1\theta_b)\sin(\lambda_{-1}\theta)-\sin(\lambda_{-1}\theta_b)\sin(\lambda_1\theta) \over\lambda_{-1}\lambda_1} &\lesssim_{C_{\rm Re},C_{\rm Im}}{\<\lambda\>\theta_b^4}\label{kernel estimates 1}\\
|\sin(\lambda_1\theta_b)\cos(\lambda_{-1}\theta) -\sin(\lambda_{-1}\theta_b)\cos(\lambda_1\theta)| &\lesssim_{C_{\rm Re},C_{\rm Im}}\theta_b\label{kernel estimates 2}\\
|\cos(\lambda_1\theta_b)\sin(\lambda_{-1}\theta) -\cos(\lambda_{-1}\theta_b)\sin(\lambda_1\theta)| &\lesssim_{C_{\rm Re},C_{\rm Im}}\theta_b\label{kernel estimates 3}\\
|\cos(\lambda_{-1}(\theta-\theta_b)) - \cos(\lambda_1(\theta-\theta_b))| &\lesssim_{C_{\rm Re},C_{\rm Im}}|\lambda|(\theta_b-\theta)^2\label{kernel estimates 4}\\
|\sin(\lambda_{-1}(\theta-\theta_b)) - \sin(\lambda_1(\theta-\theta_b))| &\lesssim_{C_{\rm Re},C_{\rm Im}}\theta_b-\theta.\label{kernel estimates 5}
\end{align}
\end{lemma}
\begin{proof}
First we will prove \eqref{kernel estimates 1}. Trigonometric identities give
\begin{align*}
\sin(\lambda_{-1}\theta)\sin(\lambda_1\theta_b) -\sin(\lambda_1\theta)\sin(\lambda_{-1}\theta_b)
&={1\over 2}(\cos(\lambda(\theta_b-\theta)+\theta_b+\theta)-\cos(\lambda(\theta_b+\theta)+\theta_b-\theta)\\
&\qquad+\cos(\lambda(\theta_b+\theta)-(\theta_b-\theta))-\cos(\lambda(\theta_b-\theta)-(\theta_b+\theta)))\\
&=\sin(\lambda(\theta_b+\theta))\sin(\theta_b-\theta)-\sin(\lambda(\theta_b-\theta))\sin(\theta_b+\theta)
\end{align*}
When $|\Im\lambda\theta_b|<C_{\rm Im}$ (together with the assumptions of the lemma, $|\Re\lambda_{\pm 1}\theta_b|<C_{\rm Re}$ and $0\leq\theta\leq\theta_b\in(0,\pi)$), we have
\begin{align}
&\sin(\lambda_{-1}\theta)\sin(\lambda_1\theta_b) -\sin(\lambda_1\theta)\sin(\lambda_{-1}\theta_b)\nonumber\\
&=\brac{\lambda_{-1}\theta+O(|\lambda_{-1}\theta_b|^3)}\brac{\lambda_1\theta_b+O(|\lambda_1\theta_b|^3)} -\brac{\lambda_1\theta+O(|\lambda_1\theta_b|^3)}\brac{\lambda_{-1}\theta_b+O(|\lambda_{-1}\theta_b|^3)}\nonumber\\
&=O(|\lambda_1\theta_b||\lambda_{-1}\theta_b|^3)+O(|\lambda_{-1}\theta_b||\lambda_1\theta_b|^3)\label{sine sine expansion}
\end{align}
and
\begin{align}
\sin&(\lambda(\theta_b+\theta))\sin(\theta_b-\theta) -\sin(\lambda(\theta_b-\theta))\sin(\theta_b+\theta)\nonumber\\
&=\brac{\lambda(\theta_b+\theta)+O(|\lambda\theta_b|^3)}\brac{\theta_b-\theta+O(\theta_b^3)}-\brac{\lambda(\theta_b-\theta)+O(|\lambda\theta_b|^3)}\brac{\theta_b+\theta+O(\theta_b^3)}\nonumber\\
&=O(|\lambda|\theta_b^4)+O(|\lambda|^3\theta_b^4)=O(\<\lambda\>^3\theta_b^4)\label{sine sine rewritten}
\end{align}
When $|\lambda|<1/2$ or $|\lambda|>3/2$ we have $|\lambda_{\pm 1}|\sim\<\lambda\>$ and so from \eqref{sine sine rewritten} we have \eqref{kernel estimates 1}. When $1/2\leq|\lambda|\leq 3/2$, we have $|\lambda_{\pm 1}|\lesssim\<\lambda\>\sim 1$ and so from \eqref{sine sine expansion} we also have \eqref{kernel estimates 1}.

Now we will prove \eqref{kernel estimates 2}. Trigonometric identities give
\begin{align*}
\sin(\lambda_1\theta_b)\cos(\lambda_{-1}\theta) -\sin(\lambda_{-1}\theta_b)\cos(\lambda_1\theta)
&={1\over 2}(\sin(\lambda(\theta_b+\theta)+\theta_b-\theta)+\sin(\lambda(\theta_b-\theta)+\theta_b+\theta)\\
&\qquad-\sin(\lambda(\theta_b+\theta)-(\theta_b-\theta))-\sin(\lambda(\theta_b-\theta)-(\theta_b+\theta)))\\
&=\cos(\lambda(\theta_b+\theta))\sin(\theta_b-\theta)+\cos(\lambda(\theta_b-\theta))\sin(\theta_b+\theta)
\end{align*}
When $|\Im\lambda\theta_b|<C_{\rm Im}$ (together with the assumptions of the lemma, $|\Re\lambda_{\pm 1}\theta_b|<C_{\rm Re}$ and $0\leq\theta\leq\theta_b\in(0,\pi)$), we have $|\Im\lambda(\theta_b\pm\theta)|<2C_{\rm Im}$, so that the two cosines in the last line are bounded by $\cosh(2C_{\rm Im})$, while $|\sin(\theta_b\pm\theta)|\leq 2\theta_b$; hence
\begin{align}
|\sin(\lambda_1\theta_b)\cos(\lambda_{-1}\theta) -\sin(\lambda_{-1}\theta_b)\cos(\lambda_1\theta)|\lesssim_{C_{\rm Im}}\theta_b.
\end{align}
This yields \eqref{kernel estimates 2}. (Note that a Taylor expansion of the LHS in the four factors, as for \eqref{kernel estimates 1}, would only give $O(1)$ here: the error terms $|\lambda_{\pm 1}\theta_b|^3$ of such an expansion are of size $1$ when $|\lambda\theta_b|\sim 1$, and it is the identity above that produces the factor $\theta_b$.)

The proof of \eqref{kernel estimates 3} 
is similar to that of \eqref{kernel estimates 2}, starting from the identity $\cos(\lambda_1\theta_b)\sin(\lambda_{-1}\theta)-\cos(\lambda_{-1}\theta_b)\sin(\lambda_1\theta)=\cos(\lambda(\theta_b+\theta))\sin(\theta_b-\theta)-\cos(\lambda(\theta_b-\theta))\sin(\theta_b+\theta)$. 

Now we will prove \eqref{kernel estimates 4}. Trigonometric identities give
\begin{align*}
\cos(\lambda_{-1}(\theta-\theta_b)) - \cos(\lambda_1(\theta-\theta_b))
&=\cos(\lambda(\theta-\theta_b)-(\theta-\theta_b))-\cos(\lambda(\theta-\theta_b)+(\theta-\theta_b))\\
&=2\sin(\lambda(\theta-\theta_b))\sin(\theta-\theta_b)
\end{align*}
When $|\Im\lambda\theta_b|<C_{\rm Im}$ (together with the assumptions of the lemma, $|\Re\lambda_{\pm 1}\theta_b|<C_{\rm Re}$ and $0\leq\theta\leq\theta_b\in(0,\pi)$), we have $|\lambda\theta_b|\lesssim_{C_{\rm Re},C_{\rm Im}} 1$ and
\begin{align}
\cos(\lambda_{-1}(\theta-\theta_b)) &- \cos(\lambda_1(\theta-\theta_b))\nonumber\\
&=2\brac{\lambda(\theta-\theta_b)+O(|\lambda(\theta-\theta_b)|^3)}\brac{(\theta-\theta_b)+O(|\theta-\theta_b|^3)}\nonumber\\
&=2\lambda(\theta-\theta_b)^2+O(|\lambda|(\theta_b-\theta)^4)+O(|\lambda|^3(\theta_b-\theta)^4)=O(|\lambda|(\theta_b-\theta)^2).
\end{align}
This yields \eqref{kernel estimates 4}.

Now we will prove \eqref{kernel estimates 5}. Trigonometric identities give
\begin{align*}
\sin(\lambda_{-1}(\theta-\theta_b)) - \sin(\lambda_1(\theta-\theta_b))
&=\sin(\lambda(\theta-\theta_b)-(\theta-\theta_b))-\sin(\lambda(\theta-\theta_b)+(\theta-\theta_b))\\
&=-2\cos(\lambda(\theta-\theta_b))\sin(\theta-\theta_b)
\end{align*}
When $|\Im\lambda\theta_b|<C_{\rm Im}$ (together with the assumptions of the lemma, $|\Re\lambda_{\pm 1}\theta_b|<C_{\rm Re}$ and $0\leq\theta\leq\theta_b\in(0,\pi)$), we have $|\lambda\theta_b|\lesssim_{C_{\rm Re},C_{\rm Im}} 1$ and
\begin{align}
\sin(\lambda_{-1}(\theta-\theta_b)) &- \sin(\lambda_1(\theta-\theta_b))\nonumber\\
&=-2\brac{1+O(|\lambda(\theta-\theta_b)|^2)}\brac{(\theta-\theta_b)+O(|\theta-\theta_b|^3)}\nonumber\\
&=-2(\theta-\theta_b)+O(|\lambda|^2(\theta_b-\theta)^3)+O((\theta_b-\theta)^3)=O(\theta_b-\theta).
\end{align}
This yields \eqref{kernel estimates 5}.
\end{proof}

In obtaining the Mellin solution representations later, we will use the incompressibility condition to reduce the PDE system to an ODE system. The solution representation of that ODE system is a quotient whose denominator is the characteristic function
\begin{align}\label{H definition}
H(\lambda):=\lambda\sin(2\theta_b)+\sin(2\lambda\theta_b),
\end{align}
and the Mellin lines and strips on which the representation is used have to avoid the zeros of $H$. The name is the usual one: for the ODE system of Lemma \ref{Green's function initial time} below, the four homogeneous boundary conditions, evaluated on the fundamental system $\sin(\lambda_{\pm 1}\theta)$, $\cos(\lambda_{\pm 1}\theta)$ of the homogeneous equation, form a matrix whose determinant is a nonzero constant multiple of $\lambda^3\lambda_1\lambda_{-1}H(\lambda)$, so that for $\lambda\notin\{0,\pm 1\}$ the homogeneous problem has a nontrivial solution exactly at the zeros of $H$. We first record what we need about these zeros. One of them is $\lambda=0$; the strip named in the following definition contains no other, by Lemma \ref{H lemma}(i).

\begin{definition}[Zero-free strip of the characteristic function]\label{zero-free strip}
For $\theta_b\in(0,\pi/2)$ let $H$ be given by \eqref{H definition}, let
\begin{align}\label{Lambda star definition}
\Lambda^*(\theta_b):=\inf\set{\Re\lambda:H(\lambda)=0,\ \Re\lambda>0},
\end{align}
and call $\{\lambda\in\C:|\Re\lambda|<\Lambda^*(\theta_b)\}$ the zero-free strip of $H$. For real $c$ with $|c|<\Lambda^*(\theta_b)$ let
\begin{align}\label{Gamma H definition}
\Gamma^H_{c,\theta_b}:=\inf\set{{|H(\lambda)|\over\min\{|\lambda|\theta_b,1\}e^{2|\Im\lambda|\theta_b}}:\Re\lambda=c,\ \lambda\neq 0}.
\end{align}
\end{definition}

Appendix \ref{resonance factors} lists the resonance factors of this paper, several of which are written with the letter $\Gamma$.

\begin{lemma}[The characteristic function]\label{H lemma}
Let $\theta_b\in(0,\pi/2)$.
\begin{enumerate}
\item $\lambda=0$ is a simple zero of $H$ and its only zero in the open strip $|\Re\lambda|<\Lambda^*(\theta_b)$, and
\begin{align}\label{Lambda star bounds}
{\pi\over 4\theta_b}\leq\Lambda^*(\theta_b)\leq{3\pi\over 4\theta_b}.
\end{align}
\item For every $\lambda\neq 0$ with $|\Re\lambda|<\Lambda^*(\theta_b)$,
\begin{align}\label{H lower bound}
|H(\lambda)|\geq\Gamma^H_{\Re\lambda,\theta_b}\min\{|\lambda|\theta_b,1\}e^{2|\Im\lambda|\theta_b},
\end{align}
and $\Gamma^H_{c,\theta_b}>0$ for every $|c|<\Lambda^*(\theta_b)$; more precisely, $\inf_{c\in K}\Gamma^H_{c,\theta_b}>0$ for every compact $K\subset(-\Lambda^*(\theta_b),\Lambda^*(\theta_b))$.
\item If $|c|\theta_b\leq\pi/4$ then $\Gamma^H_{c,\theta_b}\geq c_H$ for an absolute constant $c_H>0$.
\end{enumerate}
\end{lemma}
\begin{proof}
(i) By \eqref{H definition}, $H$ is odd and entire with real Taylor coefficients, so that its zero set is symmetric under $\lambda\mapsto-\lambda$ and $\lambda\mapsto\bar\lambda$. The zero $\lambda=0$ is simple, since $H'(0)=\sin(2\theta_b)+2\theta_b>0$. For real $y\neq 0$ we have $H(iy)=i(y\sin(2\theta_b)+\sinh(2y\theta_b))\neq 0$ by $\sin(iw)=i\sinh w$, the bracket having the sign of $y$ because $\sin(2\theta_b)>0$. By the definition of $\Lambda^*(\theta_b)$ and the symmetry $\lambda\mapsto-\lambda$ there is no zero with $0<|\Re\lambda|<\Lambda^*(\theta_b)$; this proves the first assertion. For the lower bound in \eqref{Lambda star bounds}, let $|\Re\lambda|\theta_b\leq\pi/4$ and $\lambda\neq 0$. Writing
\begin{align*}
\lambda\sin(2\theta_b) &= (\Re\lambda+i\Im\lambda)\sin(2\theta_b),\\
\sin(2\lambda\theta_b) &= \sin(2\Re\lambda\theta_b)\cosh(2\Im\lambda\theta_b) +i\cos(2\Re\lambda\theta_b)\sinh(2\Im\lambda\theta_b),
\end{align*}
and using $\sin(2\theta_b)>0$ and $\cos(2\Re\lambda\theta_b)\geq 0$, we see that the real and imaginary parts of the two summands of $H$ have the same signs, or vanish, so that
\begin{align}\label{H sign bound}
|H(\lambda)|\geq|\sin(2\lambda\theta_b)|\qquad\text{ when }\qquad|\Re\lambda|\theta_b\leq{\pi\over 4}.
\end{align}
Since $|\sin(2\lambda\theta_b)|^2=\sin^2(2\Re\lambda\theta_b)+\sinh^2(2\Im\lambda\theta_b)$ vanishes in this strip only at $\lambda=0$, we get $\Lambda^*(\theta_b)\geq\pi/(4\theta_b)$.

For the upper bound we pass to the variable $z=2\lambda\theta_b$, in which $H(\lambda)=0$ reads $h(z):=\mu z+\sin z=0$ with $\mu:=\sin(2\theta_b)/(2\theta_b)\in(0,1)$, and find a zero of $h$ with $\pi<\Re z<3\pi/2$. For $y\geq 0$ let $z(y):=\pi+\arccos(\mu y/\sinh y)+iy$, where $y/\sinh y:=1$ at $y=0$. Since $0<\mu y/\sinh y\leq\mu<1$, the curve $z(y)$ is continuous, with $\pi<\Re z(y)<3\pi/2$ and $\Re z(y)\to 3\pi/2$ as $y\to\infty$. As $\cos(\Re z(y))=-\mu y/\sinh y$ and $\sin(\Re z(y))\to-1$,
\begin{align*}
\Im h(z(y))&=\mu y+\cos(\Re z(y))\sinh y=0,\\
\Re h(z(y))&=\mu\Re z(y)+\sin(\Re z(y))\cosh y\to-\infty\qquad\text{ as }\qquad y\to\infty.
\end{align*}
So $h$ is real on the path that runs along the real axis from $\pi$ to $z(0)$ and then along the curve $z(y)$. It is positive at the start of the path, $h(\pi)=\mu\pi$, and negative at $z(y)$ for large $y$. By the intermediate value theorem it vanishes at a point of the path, where $\pi<\Re z<3\pi/2$. So $H$ has a zero with $0<\Re\lambda<3\pi/(4\theta_b)$, which is the upper bound.

(ii) The estimate \eqref{H lower bound} is the definition \eqref{Gamma H definition}. Fix a compact $K\subset(-\Lambda^*(\theta_b),\Lambda^*(\theta_b))$ and let $c_K:=\max_K|c|$. The function $R(\lambda):=|H(\lambda)|/(\min\{|\lambda|\theta_b,1\}e^{2|\Im\lambda|\theta_b})$ is continuous and positive on $\{\Re\lambda\in K\}\setminus\{0\}$ by (i), and it extends continuously to $\lambda=0$ with the positive value $H'(0)/\theta_b$, the zero there being simple. For $|\Im\lambda|\theta_b\geq 1$ we have $\min\{|\lambda|\theta_b,1\}=1$ and
\begin{align*}
|H(\lambda)|\geq|\sin(2\lambda\theta_b)|-|\lambda|\sin(2\theta_b)\geq\sinh(2|\Im\lambda|\theta_b)-(c_K+|\Im\lambda|)\sin(2\theta_b)\geq{1\over 4}e^{2|\Im\lambda|\theta_b}
\end{align*}
as soon as $|\Im\lambda|\theta_b\geq Y_0$, for some $Y_0=Y_0(K,\theta_b)\geq 1$, by $\sinh(2y)\geq{1\over 2}e^{2y}-{1\over 2}$ and since the exponential dominates the linear term. So $R\geq 1/4$ on $\{\Re\lambda\in K,\ |\Im\lambda|\theta_b\geq Y_0\}$, while on the compact set $\{\Re\lambda\in K,\ |\Im\lambda|\theta_b\leq Y_0\}$ the continuous positive function $R$ is bounded below. Hence $\inf_{c\in K}\Gamma^H_{c,\theta_b}=\inf\{R(\lambda):\Re\lambda\in K,\lambda\neq 0\}>0$.

(iii) Let $|c|\theta_b\leq\pi/4$ and let $\lambda\neq 0$ with $\Re\lambda=c$. By \eqref{H sign bound} we have $|H(\lambda)|\geq|\sin(2\lambda\theta_b)|$, which we bound below by Lemma \ref{cosh bound} with $\theta=2\theta_b$ and $C=2$, whose two cases are then $|\Im\lambda|\theta_b\geq 1$ and $|\Im\lambda|\theta_b<1$. In the first case it gives $|\sin(2\lambda\theta_b)|\gtrsim e^{2|\Im\lambda|\theta_b}$, and there $\min\{|\lambda|\theta_b,1\}=1$. In the second case it gives $|\sin(2\lambda\theta_b)|\gtrsim\Gamma_{c,2\theta_b}|\lambda|\theta_b=|\lambda|\theta_b$, since $\Gamma_{c,2\theta_b}=\min\{1,\pi-2|c|\theta_b\}=1$ for $2|c|\theta_b\leq\pi/2$, and $|\lambda|\theta_b\geq e^{-2}\min\{|\lambda|\theta_b,1\}e^{2|\Im\lambda|\theta_b}$, since $2|\Im\lambda|\theta_b<2$. In both cases $|H(\lambda)|\gtrsim\min\{|\lambda|\theta_b,1\}e^{2|\Im\lambda|\theta_b}$ with an absolute constant, which is the claim by \eqref{Gamma H definition}.
\end{proof}

\begin{remark}[The zeros of $H$ numerically]\label{H numerics remark}
In the variable $z=2\lambda\theta_b$ the zeros of $H$ are the roots of $\sin z=-\mu z$, $\mu=\sin(2\theta_b)/(2\theta_b)$, the classical eigenvalue equation of the Stokes wedge with a free-slip substrate and a stress-free surface. For $\mu=1$ the root with the smallest positive real part is $z_1=4.2124+2.2507i$, and following it in $\mu$ gives $\Lambda^*(\theta_b)\theta_b=\Re z_1/2$, which increases slowly from $2.106$ as $\theta_b\to 0$ to $2.246$ at $\theta_b=73^\circ$; beyond that the pair of complex roots has merged into two real ones, and $\Lambda^*(\theta_b)\theta_b$ decreases to $\pi/2$ as $\theta_b\to\pi/2$. For $\kappa\in(1/4,1/2)$ the condition $2-\kappa<\Lambda^*(\theta_b)$ therefore holds for $\theta_b$ up to about $73^\circ$--$75^\circ$, and the condition $3-\kappa<\Lambda^*(\theta_b)$ up to about $45^\circ$--$50^\circ$.
\end{remark}

In the proof of Proposition \ref{Mellin system initial time} below, eliminating the pressure from the Mellin transformed Stokes system and expressing $\hat u_r$ through the divergence constraint leaves the system \eqref{u theta system}: a fourth order equation for $\hat u_\theta$ on $(0,\theta_b)$ with two boundary conditions at each end. The following lemma gives its solution representation and its estimates.

\begin{lemma}[The Green's function of the Mellin system]\label{Green's function initial time}
Suppose $\theta_b\in(0,\pi/2)$, $\lambda\neq 0$ and $|\Re\lambda|<\Lambda^*(\theta_b)$, with $\Lambda^*(\theta_b)$ as in Definition \ref{zero-free strip}. The system
\begin{align*}
(\lambda_{-1}^2+\partial_\theta^2)(\lambda_1^2+\partial_\theta^2)u &= F\qquad\text{on}\quad(0,\theta_b) \\
(\lambda_1\lambda_{-1}-\partial_\theta^2)u(0) &= F_1 \\
u(0) &= F_2 \\
(\lambda_1\lambda_{-1}-\partial_\theta^2)u(\theta_b) &= F_3 \\
(3\lambda^2+1+\partial_\theta^2)\partial_\theta u(\theta_b) &= F_4
\end{align*}
has a unique solution 
\begin{align*}
u(\theta) &= \int_0^{\theta_b} G_\lambda(\theta,\xi) F(\xi)\d\xi 
+ \partial_\xi G_\lambda(\theta,0) F_1 
- (3\lambda^2+1+\partial_\xi^2)\partial_\xi G_\lambda(\theta,0) F_2 \\
&\quad - \partial_\xi G_\lambda(\theta, \theta_b) F_3 
- G_\lambda(\theta, \theta_b) F_4,
\end{align*}
where the Green's function $G_\lambda(\theta,\xi)$ is, with the formulae understood as their limits at $\lambda=\pm 1$,
\begin{align*}
  G_\lambda(\theta,\xi) &= -{1\over 4\lambda H}\begin{cases} 
w_{11}(\theta)w_{21}(\xi) + w_{12}(\theta)w_{22}(\xi) & \qquad\text{ for }\qquad 0 \le \theta \le \xi \le \theta_b \\ 
w_{11}(\xi)w_{21}(\theta) + w_{12}(\xi)w_{22}(\theta) & \qquad\text{ for }\qquad 0 \le \xi \le \theta \le \theta_b
\end{cases}\\
H
&=\lambda\sin(2\theta_b)+\sin(2\lambda\theta_b)\\
w_{11}(\theta) &=  \lambda_{-1}^{-1}\sin(\lambda_1 \theta_b) \sin(\lambda_{-1} \theta) - \lambda_1^{-1}\sin(\lambda_{-1} \theta_b) \sin(\lambda_1 \theta) \\
w_{12}(\theta) &= \cos(\lambda_1 \theta_b) \sin(\lambda_{-1} \theta) - \cos(\lambda_{-1} \theta_b) \sin(\lambda_1 \theta)\\
w_{21}(\theta) &= \lambda_1 \cos(\lambda_{-1}(\theta-\theta_b)) - \lambda_{-1} \cos(\lambda_1(\theta-\theta_b)) \\
w_{22}(\theta) &= \sin(\lambda_{-1}(\theta-\theta_b)) - \sin(\lambda_1(\theta-\theta_b))
\end{align*}
The Green's function is symmetric, $G_\lambda(\theta,\xi)=G_\lambda(\xi,\theta)$, and for $\theta\leq\xi$ we have
\begin{align*}
  \lim_{\lambda\to 0}(\lambda^2G_\lambda(\theta,\xi)) &= -{\sin\theta\sin\xi\over\sin(2\theta_b)+2\theta_b}\\
  \lim_{\lambda\to 0}\partial_\lambda(\lambda^2G_\lambda(\theta,\xi)) &= 0\\
  \lim_{\lambda\to 1}(G_\lambda(\theta,\xi)) &= -{2\sin(2\theta_b)\theta + \sin(2\theta)\sin(2(\xi-\theta_b))\over 8\sin(2\theta_b)}\\
  \lim_{\lambda\to 1}\partial_\lambda G_\lambda(\theta,\xi) &= {N(\theta,\xi)\over 8\sin(2\theta_b)},
\end{align*}
where
\begin{align*}
N(\theta,\xi) &= 2\theta\sin(2\theta_b) + \theta\sin(2\xi) + \xi\sin(2\theta) - \theta\cos(2\theta)\sin(2(\xi-\theta_b)) - (\xi-\theta_b)\sin(2\theta)\cos(2(\xi-\theta_b)) \\
&\quad + \sin(2\theta)\sin(2(\xi-\theta_b)) \brac{{3\over 2}+\theta_b\cot(2\theta_b)}.
\end{align*}
Moreover, when $\xi\not=\theta$ we have
\begin{align}\label{Green's function pointwise estimate}
&|\{\partial_\theta,\partial_\xi\}^n G_\lambda(\theta,\xi)|\nonumber\\
&\qquad\lesssim_n{\<\lambda\>^{n-1}+\theta_b^{-1}\<\lambda\>^{n-2}\over|\lambda|^2\Gamma^H_{\Re\lambda,\theta_b}}\begin{cases}\displaystyle 1, & |\Im\lambda\theta_b|\leq 1\\
\displaystyle
\sum_{l=0}^1\brac{e^{-{1\over 2}|\Im\lambda|(2\theta_b-|\xi+(-1)^l\theta|)}+e^{-{1\over 2}|\Im\lambda||\xi+(-1)^l\theta|}}, & |\Im\lambda\theta_b|>1,
\end{cases}
\end{align}
where $\Gamma^H_{\Re\lambda,\theta_b}$ is the constant \eqref{Gamma H definition} of the line $\Re\lambda$, and
\begin{align}
\|\{\partial_\theta,\partial_\xi\}^n G_\lambda\|_{L^\infty_\theta L^1_\xi}+\|\{\partial_\theta,\partial_\xi\}^n G_\lambda\|_{L^\infty_\xi L^1_\theta}&\lesssim_n{\<\lambda\>^{n-2}\over|\lambda|^2\Gamma^H_{\Re\lambda,\theta_b}}\label{Green's function L1 estimate}\\
\|\{\partial_\theta,\partial_\xi\}^n G_\lambda(\{0,\theta_b\},\ph)\|_{L^2}+\|\{\partial_\theta,\partial_\xi\}^n G_\lambda(\ph,\{0,\theta_b\})\|_{L^2}&\lesssim_n{\<\lambda\>^{n-{3\over 2}}+\theta_b^{-{1\over 2}}\<\lambda\>^{n-2}\over|\lambda|^2\Gamma^H_{\Re\lambda,\theta_b}}\label{Green's function trace estimate}
\end{align}
all the norms in $\theta$ and $\xi$ being taken over $[0,\theta_b]$.
\end{lemma}
\begin{proof}
For $\lambda\neq\pm 1$ the homogeneous equation has the fundamental system $\sin(\lambda_{\pm 1}\theta)$, $\cos(\lambda_{\pm 1}\theta)$; what follows holds at $\lambda=\pm 1$ by continuity in $\lambda$. The functions $w_{11},w_{12}$ are solutions of the homogeneous equation satisfying the two homogeneous boundary conditions at $\theta=0$, and $w_{21},w_{22}$ are solutions satisfying the two at $\theta=\theta_b$. Hence, for fixed $\xi$, $G_\lambda(\ph,\xi)$ solves the homogeneous equation on $[0,\xi)$ and on $(\xi,\theta_b]$ and satisfies all four homogeneous boundary conditions; and since the formula for $G_\lambda$ is unchanged when $\theta$ and $\xi$ are exchanged, which is the symmetry, the same holds for $G_\lambda(\theta,\ph)$ with $\theta$ fixed. A direct computation shows that the derivatives of $G_\lambda$ of order at most two are continuous across $\xi=\theta$, and that $\partial_\theta^3G_\lambda(\theta,\theta^-)-\partial_\theta^3G_\lambda(\theta,\theta^+)=1$: the third derivative increases by $1$ as $\theta$ increases through $\xi$, which is the jump that makes $(\lambda_{-1}^2+\partial_\theta^2)(\lambda_1^2+\partial_\theta^2)G_\lambda(\ph,\xi)$ the delta function at $\xi$, and which fixes the normalisation $-1/(4\lambda H)$. Now let $u$ solve the system. Multiplying its equation, written in the variable $\xi$, by $G_\lambda(\theta,\xi)$, integrating in $\xi$ over $(0,\theta)$ and over $(\theta,\theta_b)$, and integrating by parts until all derivatives fall on $G_\lambda(\theta,\ph)$, which the operator annihilates, we are left with $\int_0^{\theta_b}G_\lambda(\theta,\xi)F(\xi)\d\xi$ on one side and boundary terms on the other. Those at $\xi=\theta$ cancel except for the jump of $\partial_\xi^3G_\lambda(\theta,\ph)$, which returns $u(\theta)$; those at $\xi=0$ and $\xi=\theta_b$ are, by the homogeneous boundary conditions of $G_\lambda(\theta,\ph)$, the four terms carrying $F_1,\ldots,F_4$. This is the solution formula. In particular the solution is unique, and therefore it exists, the four boundary conditions being four linear equations for the four constants of the general solution of the equation.

It remains to prove the limits and the estimates. In the limits we take $\theta\leq\xi$, as in the statement. We have
\begin{align*}
\lim_{\lambda\to 0}(\lambda^2G_\lambda(\theta,\xi))
 &= -{2\sin(\theta_b)\sin(\theta)2\cos(\xi-\theta_b) + 2\cos(\theta_b)\sin(\theta)2\sin(\xi-\theta_b)\over 4(\sin(2\theta_b)+2\theta_b)}
= -{\sin\theta\sin\xi\over\sin(2\theta_b)+2\theta_b}
\end{align*}
Notice that replacing $\lambda$ with $-\lambda$ leaves our system of ODEs and its boundary conditions unchanged, hence $G_\lambda=G_{-\lambda}$, i.e.\ $G$ is an even function of $\lambda$, which means $\lim_{\lambda\to 0}\partial_\lambda(\lambda^2G_\lambda(\theta,\xi))=0$. We have
\begin{align*}
  \lim_{\lambda\to 1}(G_\lambda(\theta,\xi))
 &= -{2\sin(2\theta_b)\theta + \sin(2\theta)\sin(2(\xi-\theta_b))\over 8\sin(2\theta_b)},
\end{align*}
and differentiating the formula for $G_\lambda$ in $\lambda$ and letting $\lambda\to 1$ gives the limit of $\partial_\lambda G_\lambda$ in the statement.

The denominator $H$ is controlled by Lemma \ref{H lemma}(ii): since $|\Re\lambda|<\Lambda^*(\theta_b)$,
\begin{align}\label{H bound in proof}
  |H(\lambda)|\geq\Gamma^H_{\Re\lambda,\theta_b}\min\{|\lambda|\theta_b,1\}e^{2|\Im\lambda|\theta_b}\geq\Gamma^H_{\Re\lambda,\theta_b}\begin{cases}{1\over 4}|\lambda|\theta_b, & |\Im\lambda\theta_b|\leq 1,\\ e^{2|\Im\lambda|\theta_b}, & |\Im\lambda\theta_b|>1,\end{cases}
\end{align}
where in the first regime we used $|\lambda|\theta_b\leq|\Re\lambda|\theta_b+|\Im\lambda|\theta_b<3\pi/4+1<4$, by \eqref{Lambda star bounds}. The factor $\Gamma^H_{\Re\lambda,\theta_b}$ is positive, and it is bounded below by an absolute constant when $|\Re\lambda|\theta_b\leq\pi/4$ by Lemma \ref{H lemma}(iii). We have
\begin{align*}
(-1)^n\partial_\theta^{2n} w_{11}(\theta)&=\lambda_{-1}^{2n-1}\sin(\lambda_1\theta_b)\sin(\lambda_{-1}\theta) - \lambda_1^{2n-1}\sin(\lambda_{-1}\theta_b)\sin(\lambda_1\theta)\\
&={1\over\lambda_{-1}\lambda_1}\brac{\lambda_{-1}^{2n}\lambda_1\sin(\lambda_1\theta_b)\sin(\lambda_{-1}\theta) - \lambda_1^{2n}\lambda_{-1}\sin(\lambda_{-1}\theta_b)\sin(\lambda_1\theta)}\\
(-1)^n\partial_\theta^{2n+1} w_{11}(\theta)&=\lambda_{-1}^{2n}\sin(\lambda_1\theta_b)\cos(\lambda_{-1}\theta) - \lambda_1^{2n}\sin(\lambda_{-1}\theta_b)\cos(\lambda_1\theta)
\end{align*}
When $|\Im\lambda\theta_b|\leq 1$ we have $\<\lambda\>\theta_b\lesssim 1$, as $|\lambda|\theta_b<4$ and $\theta_b<\pi/2$, and Lemma \ref{Green trigonometric bounds} applies with $C_{\rm Re}=5\pi/4$ and $C_{\rm Im}=2$, as $|\Re\lambda_{\pm 1}|\theta_b<3\pi/4+\pi/2$, so that the implicit constants below depend only on $n$. In the formulae above we replace $\lambda_{-1}^{2n}\lambda_1$ and $\lambda_1^{2n}\lambda_{-1}$ by $\lambda^{2n+1}$, and $\lambda_{\pm 1}^{2n}$ by $\lambda^{2n}$, with errors $O(\<\lambda\>^{2n})$ and $O(\<\lambda\>^{2n-1})$, and bound the error terms by $|\sin(\lambda_{\pm 1}\theta)|\lesssim|\lambda_{\pm 1}|\theta$ and $|\cos(\lambda_{\pm 1}\theta)|\lesssim 1$, which hold since $|\Im\lambda|\theta\leq 1$. By \eqref{kernel estimates 1} and \eqref{kernel estimates 2} we get
\begin{align*}
(-1)^n\partial_\theta^{2n} w_{11}(\theta)
&={\lambda^{2n+1}\over\lambda_{-1}\lambda_1}\brac{\sin(\lambda_1\theta_b)\sin(\lambda_{-1}\theta)-\sin(\lambda_{-1}\theta_b)\sin(\lambda_1\theta)}
+O\brac{\<\lambda\>^{2n}\theta_b\theta}=O\brac{\<\lambda\>^{2n}\theta_b^2}\\
(-1)^n\partial_\theta^{2n+1} w_{11}(\theta)
&=\lambda^{2n}\brac{\sin(\lambda_1\theta_b)\cos(\lambda_{-1}\theta) - \sin(\lambda_{-1}\theta_b)\cos(\lambda_1\theta)}
+O\brac{\<\lambda\>^{2n}\theta_b}=O\brac{\<\lambda\>^{2n}\theta_b}
\end{align*}
Therefore we have
\begin{align*}
  |\partial_\theta^{n} w_{11}(\theta)|\lesssim\<\lambda\>^{n-1}\theta_b
\end{align*}
Similarly, we can show that
\begin{align*}
  |\partial_\theta^{n} w_{12}(\theta)|,\ |\partial_\theta^{n} w_{22}(\theta)| &\lesssim\begin{cases}\<\lambda\>^{n}\theta_b, & n\text{ even}\\ \<\lambda\>^{n-1}, & n\text{ odd}\end{cases}\\
  |\partial_\theta^{n} w_{21}(\theta)| &\lesssim\begin{cases}\<\lambda\>^{n}, & n\text{ even}\\ \<\lambda\>^{n+1}\theta_b, & n\text{ odd}.\end{cases}
\end{align*}
The parity matters here and cannot be dispensed with: for instance $w_{21}(\theta_b)=\lambda_1-\lambda_{-1}=2$ exactly, and $\partial_\theta w_{22}(\theta_b)=-2$ exactly, neither of which is $O(\theta_b)$. In the Green's function these bounds enter through the products of $w_{11}$ with $w_{21}$ and of $w_{12}$ with $w_{22}$; using $\theta_b\lesssim\<\lambda\>^{-1}$ they give, for $n_1,n_2\in\N_0$ and $n=n_1+n_2$,
\begin{align*}
|\partial_\theta^{n_1}w_{11}(\theta)\partial_\xi^{n_2}w_{21}(\xi)|&\lesssim\<\lambda\>^{n-1}\theta_b\\
|\partial_\theta^{n_1}w_{12}(\theta)\partial_\xi^{n_2}w_{22}(\xi)|&\lesssim\<\lambda\>^{n-1}\theta_b+\mb 1_{\{\text{$n_1,n_2$ odd}\}}\<\lambda\>^{n-2}
\end{align*}
and the same with $\theta$ and $\xi$ interchanged. So when $|\Im\lambda\theta_b|\leq 1$, by \eqref{H bound in proof}, we have for $\xi\neq\theta$
\begin{align*}
|\{\partial_\theta,\partial_\xi\}^n G_\lambda(\theta,\xi)|\lesssim{\<\lambda\>^{n-1}\theta_b+\<\lambda\>^{n-2}\over|\lambda||H(\lambda)|} \lesssim{\<\lambda\>^{n-1}+\theta_b^{-1}\<\lambda\>^{n-2}\over|\lambda|^2\Gamma^H_{\Re\lambda,\theta_b}},
\end{align*}
which is \eqref{Green's function pointwise estimate} for $|\Im\lambda\theta_b|\leq 1$. Integrating over $\xi$ or $\theta$ across an interval of length $\theta_b$ multiplies this bound by $\theta_b$, and since $\theta_b\<\lambda\>^{n-1}\lesssim\<\lambda\>^{n-2}$ we get \eqref{Green's function L1 estimate} for $|\Im\lambda\theta_b|\leq 1$. Taking instead the $L^2$ norm over such an interval multiplies it by $\theta_b^{1/2}$, and since $\theta_b^{1/2}\<\lambda\>^{n-1}\lesssim\<\lambda\>^{n-3/2}$ we get \eqref{Green's function trace estimate} for $|\Im\lambda\theta_b|\leq 1$.

Now let $|\Im\lambda\theta_b|>1$. We have
\begin{align*}
w_{11}(\theta)w_{21}(\xi)+w_{12}(\theta)w_{22}(\xi)
&={\lambda^2J_0+\lambda J_1+J_2\over\lambda_{-1}\lambda_1},
\end{align*}
where
\begin{align*}
J_0&=(\sin(\lambda_1\theta_b)\sin(\lambda_{-1}\theta)-\sin(\lambda_{-1}\theta_b)\sin(\lambda_1\theta))(\cos(\lambda_{-1}(\xi-\theta_b))-\cos(\lambda_1(\xi-\theta_b)))\\
&\quad+(\cos(\lambda_1\theta_b)\sin(\lambda_{-1}\theta)-\cos(\lambda_{-1}\theta_b)\sin(\lambda_1\theta))(\sin(\lambda_{-1}(\xi-\theta_b))-\sin(\lambda_1(\xi-\theta_b)))\\
J_1&=(\sin(\lambda_1\theta_b)\sin(\lambda_{-1}\theta)-\sin(\lambda_{-1}\theta_b)\sin(\lambda_1\theta))(\cos(\lambda_{-1}(\xi-\theta_b))+\cos(\lambda_1(\xi-\theta_b)))\\
&\quad+(\sin(\lambda_1\theta_b)\sin(\lambda_{-1}\theta)+\sin(\lambda_{-1}\theta_b)\sin(\lambda_1\theta))(\cos(\lambda_{-1}(\xi-\theta_b))-\cos(\lambda_1(\xi-\theta_b)))\\
J_2&=(\sin(\lambda_1\theta_b)\sin(\lambda_{-1}\theta)+\sin(\lambda_{-1}\theta_b)\sin(\lambda_1\theta))(\cos(\lambda_{-1}(\xi-\theta_b))+\cos(\lambda_1(\xi-\theta_b)))\\
&\quad-(\cos(\lambda_1\theta_b)\sin(\lambda_{-1}\theta)-\cos(\lambda_{-1}\theta_b)\sin(\lambda_1\theta))(\sin(\lambda_{-1}(\xi-\theta_b))-\sin(\lambda_1(\xi-\theta_b)))
\end{align*}
Using trigonometric identities we get
\begin{align*}
J_0&=2(\sin(\lambda(\theta-\theta_b))\sin(\theta+\theta_b)-\sin(\lambda(\theta+\theta_b))\sin(\theta-\theta_b))\sin(\lambda(\xi-\theta_b))\sin(\xi-\theta_b)\\
&\quad+2(\cos(\lambda(\theta+\theta_b))\sin(\theta-\theta_b)+\cos(\lambda(\theta-\theta_b))\sin(\theta+\theta_b))\cos(\lambda(\xi-\theta_b))\sin(\xi-\theta_b)\\
&=2\brac{\sin(\theta+\theta_b)\cos(\lambda(\theta-\xi))+\sin(\theta-\theta_b)\cos(\lambda(\theta+\xi))}\sin(\xi-\theta_b)\\
J_1&=2(\sin(\lambda(\theta-\theta_b))\sin(\theta+\theta_b)-\sin(\lambda(\theta+\theta_b))\sin(\theta-\theta_b))\cos(\lambda(\xi-\theta_b))\cos(\xi-\theta_b)\\
&\quad+2(\cos(\lambda(\theta-\theta_b))\cos(\theta+\theta_b)-\cos(\lambda(\theta+\theta_b))\cos(\theta-\theta_b))\sin(\lambda(\xi-\theta_b))\sin(\xi-\theta_b)\\
&=\sin(\theta+\xi)\sin(\lambda(\theta+\xi-2\theta_b))-\sin(\theta+\xi-2\theta_b)\sin(\lambda(\theta+\xi))\\
&\quad+\sin(\theta-\xi+2\theta_b)\sin(\lambda(\theta-\xi))-\sin(\theta-\xi)\sin(\lambda(\theta-\xi+2\theta_b))\\
J_2&=2(\cos(\lambda(\theta-\theta_b))\cos(\theta+\theta_b)-\cos(\lambda(\theta+\theta_b))\cos(\theta-\theta_b))\cos(\lambda(\xi-\theta_b))\cos(\xi-\theta_b)\\
&\quad-2(\cos(\lambda(\theta+\theta_b))\sin(\theta-\theta_b)+\cos(\lambda(\theta-\theta_b))\sin(\theta+\theta_b))\cos(\lambda(\xi-\theta_b))\sin(\xi-\theta_b)\\
&=2\brac{\cos(\theta+\xi)\cos(\lambda(\theta-\theta_b))-\cos(\theta-\xi)\cos(\lambda(\theta+\theta_b))}\cos(\lambda(\xi-\theta_b))
\end{align*}
So when $|\Im\lambda\theta_b|>1$ and $\xi\geq\theta$ we have
\begin{align*}
&|\{\partial_\theta,\partial_\xi\}^n(w_{11}(\theta)w_{21}(\xi)+w_{12}(\theta)w_{22}(\xi))|\\
&\lesssim\<\lambda\>^{n-2}|\lambda|\brac{|\lambda|(\theta_b-\xi)+1}\brac{(\theta_b+\theta)e^{|\Im\lambda|(\xi-\theta)}+(\theta_b-\theta)e^{|\Im\lambda|(\xi+\theta)}}\\
&\quad+\<\lambda\>^{n-3}|\lambda|^2(\theta_b-\xi)\brac{e^{|\Im\lambda|(\xi-\theta)}+e^{|\Im\lambda|(\xi+\theta)}}\\
&\quad+\<\lambda\>^{n-2}|\lambda|\sum_{l=0}^1\brac{(\xi+(-1)^l\theta)e^{|\Im\lambda|(2\theta_b-\xi-(-1)^l\theta)}+(2\theta_b+(-1)^l\xi-\theta)e^{|\Im\lambda|(\xi+(-1)^l\theta)}}\\
&\quad+\<\lambda\>^{n-3}|\lambda|\sum_{l=0}^1\brac{e^{|\Im\lambda|(2\theta_b-\xi-(-1)^l\theta)}+e^{|\Im\lambda|(\xi+(-1)^l\theta)}}\\
&\quad+\<\lambda\>^{n-2}\brac{e^{|\Im\lambda|(2\theta_b-\xi-\theta)}+e^{|\Im\lambda|(2\theta_b-\xi+\theta)}}\\
&\quad+\<\lambda\>^{n-3}\brac{(\xi+\theta)e^{|\Im\lambda|(2\theta_b-\xi-\theta)}+(\xi-\theta)e^{|\Im\lambda|(2\theta_b-\xi+\theta)}}
\end{align*}
Using that
\[|\lambda|\lesssim|\Re\lambda|+|\Im\lambda|\leq\theta_b^{-1}\brac{{3\pi\over 4}+|\Im\lambda\theta_b|}\leq\theta_b^{-1}\brac{{3\pi\over 4}+1}|\Im\lambda\theta_b|\leq\brac{{3\pi\over 4}+1}|\Im\lambda|\]
and that exponential growth beats polynomial growth we have
\begin{align*}
|\{\partial_\theta,\partial_\xi\}^n(w_{11}(\theta)w_{21}(\xi)&+w_{12}(\theta)w_{22}(\xi))|e^{-2|\Im\lambda|\theta_b}\\
&\lesssim\<\lambda\>^{n-2}\sum_{l=0}^1\brac{e^{-{1\over 2}|\Im\lambda|(\theta_b+(-1)^l\theta+\theta_b-\xi)}+e^{-{1\over 2}|\Im\lambda|(\xi+(-1)^l\theta)}}
\end{align*}
Similarly, when $|\Im\lambda\theta_b|>1$ and $\theta\geq\xi$ we have
\begin{align*}
|\{\partial_\theta,\partial_\xi\}^n(w_{11}(\xi)w_{21}(\theta)&+w_{12}(\xi)w_{22}(\theta))|e^{-2|\Im\lambda|\theta_b}\\
&\lesssim\<\lambda\>^{n-2}\sum_{l=0}^1\brac{e^{-{1\over 2}|\Im\lambda|(\theta_b+(-1)^l\xi+\theta_b-\theta)}+e^{-{1\over 2}|\Im\lambda|(\theta+(-1)^l\xi)}}
\end{align*}
So when $|\Im\lambda\theta_b|>1$, using $|H(\lambda)|\geq\Gamma^H_{\Re\lambda,\theta_b}e^{2|\Im\lambda|\theta_b}$ from \eqref{H bound in proof},
\begin{align*}
|\{\partial_\theta,\partial_\xi\}^n G_\lambda(\theta,\xi)|
&\lesssim{\<\lambda\>^{n-2}\over|\lambda|\Gamma^H_{\Re\lambda,\theta_b}}
\sum_{l=0}^1\begin{cases}
e^{-{1\over 2}|\Im\lambda|(\theta_b+(-1)^l\theta+\theta_b-\xi)}+e^{-{1\over 2}|\Im\lambda|(\xi+(-1)^l\theta)}, & \xi>\theta;\\
e^{-{1\over 2}|\Im\lambda|(\theta_b+(-1)^l\xi+\theta_b-\theta)}+e^{-{1\over 2}|\Im\lambda|(\theta+(-1)^l\xi)}, &\theta>\xi.
                               \end{cases}\\
                               &\lesssim{\<\lambda\>^{n-2}\over|\lambda|\Gamma^H_{\Re\lambda,\theta_b}}
\sum_{l=0}^1\brac{e^{-{1\over 2}|\Im\lambda|(2\theta_b-|\xi+(-1)^l\theta|)}+e^{-{1\over 2}|\Im\lambda||\xi+(-1)^l\theta|}}
\end{align*}
which is \eqref{Green's function pointwise estimate} in this regime, since $\<\lambda\>/|\lambda|\geq 1$ (note that we exclude the line $\xi=\theta$ here since from $n\geq 4$ onwards delta functions form due to the limited smoothness of the Green's function across $\xi=\theta$).

It remains to prove \eqref{Green's function L1 estimate} and \eqref{Green's function trace estimate} in the regime $|\Im\lambda\theta_b|>1$, where by \eqref{Green's function pointwise estimate} we have, for $\xi\not=\theta$,
\begin{align*}
|\{\partial_\theta,\partial_\xi\}^n G_\lambda(\theta,\xi)|\lesssim{\<\lambda\>^{n-1}\over|\lambda|^2\Gamma^H_{\Re\lambda,\theta_b}}
\sum_{l=0}^1\brac{e^{-{1\over 2}|\Im\lambda|(2\theta_b-|\xi+(-1)^l\theta|)}+e^{-{1\over 2}|\Im\lambda||\xi+(-1)^l\theta|}},
\end{align*}
the term $\theta_b^{-1}\<\lambda\>^{n-2}$ of \eqref{Green's function pointwise estimate} being absent here since $\<\lambda\>\theta_b\gtrsim 1$ in this regime. We have
\begin{align*}
\sup_{\theta\in[0,\theta_b]}\int_0^{\theta_b}e^{-{1\over 2}|\Im\lambda||\xi+\theta|}\d\xi
&=\sup_{\theta\in[0,\theta_b]}{e^{-{1\over 2}|\Im\lambda|\theta}(1-e^{-{1\over 2}|\Im\lambda|\theta_b})\over{1\over 2}|\Im\lambda|}\lesssim{1\over|\Im\lambda|}\\
\sup_{\theta\in[0,\theta_b]}\int_0^{\theta_b}e^{-{1\over 2}|\Im\lambda||\xi-\theta|}\d\xi
&=\sup_{\theta\in[0,\theta_b]}{2-e^{-{1\over 2}|\Im\lambda|\theta}-e^{-{1\over 2}|\Im\lambda|(\theta_b-\theta)}\over{1\over 2}|\Im\lambda|}\lesssim{1\over|\Im\lambda|}\\
\sup_{\theta\in[0,\theta_b]}\int_0^{\theta_b}e^{-{1\over 2}|\Im\lambda|(2\theta_b-|\xi+\theta|)}\d\xi
&=\sup_{\theta\in[0,\theta_b]}{e^{-{1\over 2}|\Im\lambda|(\theta_b-\theta)}(1-e^{-{1\over 2}|\Im\lambda|\theta_b})\over{1\over 2}|\Im\lambda|}\lesssim{1\over|\Im\lambda|}\\
\sup_{\theta\in[0,\theta_b]}\int_0^{\theta_b}e^{-{1\over 2}|\Im\lambda|(2\theta_b-|\xi-\theta|)}\d\xi
&=\sup_{\theta\in[0,\theta_b]}{e^{-{1\over 2}|\Im\lambda|(2\theta_b-\theta)}+e^{-{1\over 2}|\Im\lambda|(\theta_b+\theta)}-2e^{-|\Im\lambda|\theta_b}\over{1\over 2}|\Im\lambda|}\lesssim{1\over|\Im\lambda|}
\end{align*}
Therefore, when $|\Im\lambda\theta_b|>1$, integrating in $\xi$, or likewise in $\theta$, we get \eqref{Green's function L1 estimate} with $\<\lambda\>^{n-1}|\Im\lambda|^{-1}$ in place of $\<\lambda\>^{n-2}$, which is the same thing since $|\Im\lambda|\leq\<\lambda\>\lesssim|\Im\lambda|$ in this regime; and there is no $\theta_b^{-1}$ term to absorb, the $O(\theta_b)$ bookkeeping above playing no role once the exponentials are present. Together with the regime $|\Im\lambda\theta_b|\leq 1$ treated above, this proves \eqref{Green's function L1 estimate}.

We have, for $\xi\in(0,\theta_b)$ and $|\Im\lambda\theta_b|>1$,
\begin{align*}
|\{\partial_\theta,\partial_\xi\}^n G_\lambda(0,\xi)|&\lesssim{\<\lambda\>^{n-1}\over|\lambda|^2\Gamma^H_{\Re\lambda,\theta_b}}\brac{e^{-{1\over 2}|\Im\lambda|(2\theta_b-\xi)}+e^{-{1\over 2}|\Im\lambda|\xi}}\\
|\{\partial_\theta,\partial_\xi\}^n G_\lambda(\theta_b,\xi)|&\lesssim{\<\lambda\>^{n-1}\over|\lambda|^2\Gamma^H_{\Re\lambda,\theta_b}}\sum_{l=0}^1e^{-{1\over 2}|\Im\lambda|(\theta_b+(-1)^l\xi)}
\end{align*}
and so, taking the $L^2_\xi$ norm and using $\|e^{-|\Im\lambda|\ph/2}\|_{L^2}\lesssim|\Im\lambda|^{-1/2}\lesssim\<\lambda\>^{-1/2}$, we get \eqref{Green's function trace estimate} for $G_\lambda(\{0,\theta_b\},\ph)$ when $|\Im\lambda\theta_b|>1$. By the symmetry of $G_\lambda$ we get the same bound for $G_\lambda(\ph,\{0,\theta_b\})$. Together with the regime $|\Im\lambda\theta_b|\leq 1$ treated above, this proves \eqref{Green's function trace estimate}.
\end{proof}

Using the solution representation of the ODE system of Lemma \ref{Green's function initial time}, we can obtain the solution representation of the PDE system in the following proposition; Remark \ref{Mellin initial time remark} below comments on the weights it excludes and on the constants in its estimate.

\begin{proposition}[The Mellin system at the initial time]\label{Mellin system initial time}
Let $\theta_b\in(0,\pi/2)$ and let $a\notin\{0,2\}$ satisfy $0<|a-1|<\Lambda^*(\theta_b)$, with $\Lambda^*(\theta_b)$ as in Definition \ref{zero-free strip}. Provided the RHS of \eqref{Mellin estimate initial time} is finite, the system
\begin{alignat*}{3}
r^2(\lpc\mb u-\grad p)&=\mb F_\Omega\qquad&&\text{on}\qquad\Omega\\
\grad\cdot\mb u&=0\qquad&&\text{on}\qquad\Omega\\
\partial_\theta u_r+\dr_{-1}u_\theta&=F_{0,0}\qquad&&\text{on}\qquad\partial_0\Omega\\
u_\theta&=F_{0,1}\qquad&&\text{on}\qquad\partial_0\Omega\\
\partial_\theta u_r+\dr_{-1}u_\theta&=F_{1,0}\qquad&&\text{on}\qquad\partial_1\Omega\\
2(\partial_\theta u_\theta+u_r) -rp
&=F_{1,1}\qquad&&\text{on}\qquad\partial_1\Omega
\end{alignat*}
has a unique solution with regularity
\begin{align}
  &\|\mb u\|_{\Omega,a,j,2}+\Gamma_{a-2}^{-1}\|p\|_{\Omega,a-1,j,1}\nonumber\\
  &\lesssim_j{\Gamma_{a-1}^2\over\Gamma^H_{a-1,\theta_b}}\brac{\|\mb F_\Omega\|_{\Omega,a,j}+\|F_{0,0}\|_{\partial_0\Omega,a-{1\over 2},j+{1\over 2}}+\|(F_{1,0},F_{1,1})\|_{\partial_1\Omega,a-{1\over 2},j+{1\over 2}} +\Gamma_{a}\|F_{0,1}\|_{\partial_0\Omega,a-{1\over 2},j+{3\over 2}}}\nonumber\\
  &\quad+{\Gamma_{a-1}^2\over\theta_b^{1\over 2}\Gamma^H_{a-1,\theta_b}}\brac{\|F_{0,0}\|_{\partial_0\Omega,a-{1\over 2},j}+\|(F_{1,0},F_{1,1})\|_{\partial_1\Omega,a-{1\over 2},j} +\Gamma_{a}\|F_{0,1}\|_{\partial_0\Omega,a-{1\over 2},j+1}},\label{Mellin estimate initial time}
\end{align}
where $\Gamma_c:=\<c\>/|c|$ for real $c\neq 0$, and $\Gamma^H_{a-1,\theta_b}$ is the constant \eqref{Gamma H definition} of the Mellin line $\Re\lambda=a-1$. Moreover, $(\mb u,p)$ has a solution representation in Mellin variables as
\begin{align*}
  \hat u_\theta(\lambda,\theta) &= \int_0^{\theta_b} \lambda_1\brac{\lambda_{-1}\hat F_{\Omega,\theta}(\lambda,\xi)+\hat F_{\Omega,r}(\lambda,\xi)\partial_\xi}G_\lambda(\theta,\xi)\d\xi+ \partial_\xi G_\lambda(\theta,0) \lambda_1\hat F_{0,0}(\lambda)\\
&\quad- (3\lambda^2+1+\partial_\xi^2)\partial_\xi G_\lambda(\theta,0) \hat F_{0,1}(\lambda)  - \partial_\xi G_\lambda(\theta, \theta_b) \lambda_1\hat F_{1,0}(\lambda) - G_\lambda(\theta, \theta_b) \lambda_1\lambda_{-1}\hat F_{1,1}(\lambda)\\
\hat u_r(\lambda,\theta) &= -\int_0^{\theta_b} \brac{\lambda_{-1}\hat F_{\Omega,\theta}(\lambda,\xi)+\hat F_{\Omega,r}(\lambda,\xi)\partial_\xi}\partial_\theta G_\lambda(\theta,\xi)\d\xi-\partial_\theta\partial_\xi G_\lambda(\theta,0)\hat F_{0,0}(\lambda)\\
&\quad + \lambda_1^{-1}(3\lambda^2+1+\partial_\xi^2)\partial_\theta\partial_\xi G_\lambda(\theta,0) \hat F_{0,1}(\lambda)  + \partial_\theta\partial_\xi G_\lambda(\theta, \theta_b)\hat F_{1,0}(\lambda) + \partial_\theta G_\lambda(\theta, \theta_b)\lambda_{-1}\hat F_{1,1}(\lambda)\\
\hat p(\lambda,\theta) &= -\int_0^{\theta_b} \brac{\hat F_{\Omega,\theta}(\lambda_1,\xi)+{1\over\lambda}\hat F_{\Omega,r}(\lambda_1,\xi)\partial_\xi}\partial_\theta G_{\lambda,p}(\theta,\xi)\d\xi\\
&\quad -{1\over\lambda}\partial_\xi\partial_\theta G_{\lambda,p}(\theta,0)\hat F_{0,0}(\lambda_1) + {1\over\lambda\lambda_2}(3\lambda_1^2+1+\partial_\xi^2)\partial_\xi\partial_\theta G_{\lambda,p}(\theta,0) \hat F_{0,1}(\lambda_1)\\
&\quad + {1\over\lambda}\partial_\xi\partial_\theta G_{\lambda,p}(\theta, \theta_b)\hat F_{1,0}(\lambda_1) + \partial_\theta G_{\lambda,p}(\theta, \theta_b)\hat F_{1,1}(\lambda_1),
\end{align*}
where $G_{\lambda,p}=(\lambda_2^2+\partial_\theta^2)G_{\lambda_1}$ and $G_\lambda$ is the Green's function from Lemma \ref{Green's function initial time}.
\end{proposition}
\begin{proof}
Since $\mb u$ is divergence-free, \eqref{lpc form} gives $\lpc\mb u=\grad^\perp\grad^\perp\cdot\mb u$, so that
\begin{align*}
r^2(\lpc\mb u-\grad p)
&=r^2(\grad^\perp\grad^\perp\cdot\mb u-\grad p)\\
&=\brac{-r\partial_\theta(r^{-1}(\dr_1u_\theta-\partial_\theta u_r))-r\dr p}\mb e_r+\brac{r\dr(r^{-1}(\dr_1u_\theta-\partial_\theta u_r))-r\partial_\theta p}\mb e_\theta\\
&=\brac{-\partial_\theta(\dr_1u_\theta-\partial_\theta u_r)-r\dr p}\mb e_r+\brac{\dr_{-1}(\dr_1u_\theta-\partial_\theta u_r)-r\partial_\theta p}\mb e_\theta.
\end{align*}
Hence, in Mellin variables, the system of the statement reads
\begin{subequations}
\begin{align}
-\partial_\theta(\lambda_1\hat u_\theta-\partial_\theta\hat u_r)-\lambda_{-1}\hat p(\lambda_{-1},\theta)&=\hat F_{\Omega,r}\label{p from u}\\
\lambda_{-1}(\lambda_1\hat u_\theta-\partial_\theta\hat u_r)-\partial_\theta\hat p(\lambda_{-1},\theta)&=\hat F_{\Omega,\theta}\label{d p from u}\\
\lambda_1\hat u_r+\partial_\theta\hat u_\theta&=0\label{Mellin divergence}\\
(\partial_\theta\hat u_r+\lambda_{-1}\hat u_\theta)(\lambda,0)&=\hat F_{0,0}\\
\hat u_\theta(\lambda,0)&=\hat F_{0,1}\\
(\partial_\theta\hat u_r+\lambda_{-1}\hat u_\theta)(\lambda,\theta_b)&=\hat F_{1,0}\\
2(\partial_\theta\hat u_\theta+\hat u_r)(\lambda,\theta_b)-\hat p(\lambda_{-1},\theta_b) &=\hat F_{1,1}\label{p boundary}
\end{align}
\end{subequations}
Taking $\lambda_{-1}$ times \eqref{d p from u} minus $\partial_\theta$ of \eqref{p from u}, and $\lambda_{-1}$ times \eqref{p boundary} minus \eqref{p from u} at $\theta=\theta_b$, we eliminate $\hat p$:
\begin{align*}
(\lambda_{-1}^2+\partial_\theta^2)(\lambda_1\hat u_\theta-\partial_\theta\hat u_r) &=\lambda_{-1}\hat F_{\Omega,\theta}-\partial_\theta\hat F_{\Omega,r}\\
2\lambda_{-1}(\partial_\theta\hat u_\theta+\hat u_r)(\lambda,\theta_b)+\partial_\theta(\lambda_1\hat u_\theta-\partial_\theta\hat u_r)(\lambda,\theta_b)+\hat F_{\Omega,r}(\lambda,\theta_b) &=\lambda_{-1}\hat F_{1,1}
\end{align*}
Using the divergence constraint \eqref{Mellin divergence} we have
\begin{subequations}\label{u theta system}
\begin{align}
(\lambda_{-1}^2+\partial_\theta^2)(\lambda_1^2+\partial_\theta^2)\hat u_\theta &=\lambda_1\brac{\lambda_{-1}\hat F_{\Omega,\theta}-\partial_\theta\hat F_{\Omega,r}}\\
(\lambda_1\lambda_{-1}-\partial_\theta^2)\hat u_\theta(\lambda,0)&=\lambda_1\hat F_{0,0}\\
\hat u_\theta(\lambda,0)&=\hat F_{0,1}\\
(\lambda_1\lambda_{-1}-\partial_\theta^2)\hat u_\theta(\lambda,\theta_b)&=\lambda_1\hat F_{1,0}\\
\brac{3\lambda^2+1+\partial_\theta^2}\partial_\theta\hat u_\theta(\lambda,\theta_b)&=\lambda_1\brac{\lambda_{-1}\hat F_{1,1}-\hat F_{\Omega,r}(\lambda,\theta_b)}
\end{align}
\end{subequations}
This is the system of Lemma \ref{Green's function initial time} with $F=\lambda_1(\lambda_{-1}\hat F_{\Omega,\theta}-\partial_\theta\hat F_{\Omega,r})$, $F_1=\lambda_1\hat F_{0,0}$, $F_2=\hat F_{0,1}$, $F_3=\lambda_1\hat F_{1,0}$ and $F_4=\lambda_1(\lambda_{-1}\hat F_{1,1}-\hat F_{\Omega,r}(\lambda,\theta_b))$, so that
\begin{align*}
\hat u_\theta(\theta) &= \int_0^{\theta_b} G_\lambda(\theta,\xi) \lambda_1\brac{\lambda_{-1}\hat F_{\Omega,\theta}-\partial_\xi\hat F_{\Omega,r}}(\lambda,\xi) \d\xi 
+ \partial_\xi G_\lambda(\theta,0) \lambda_1\hat F_{0,0}(\lambda)\\
&\quad- (3\lambda^2+1+\partial_\xi^2)\partial_\xi G_\lambda(\theta,0) \hat F_{0,1}(\lambda)  - \partial_\xi G_\lambda(\theta, \theta_b) \lambda_1\hat F_{1,0}(\lambda)\\
&\quad - G_\lambda(\theta, \theta_b) \lambda_1\brac{\lambda_{-1}\hat F_{1,1}(\lambda)-\hat F_{\Omega,r}(\lambda,\theta_b)}
\end{align*}
We have
\begin{align*}
\int_0^{\theta_b} G_\lambda(\theta,\xi) \lambda_1\brac{\lambda_{-1}\hat F_{\Omega,\theta}-\partial_\xi\hat F_{\Omega,r}}&(\lambda,\xi) \d\xi\\
&= G_\lambda(\theta,0)\lambda_1\hat F_{\Omega,r}(\lambda,0)-G_\lambda(\theta,\theta_b)\lambda_1\hat F_{\Omega,r}(\lambda,\theta_b)\\
&\quad+\int_0^{\theta_b} \lambda_1\brac{\lambda_{-1}\hat F_{\Omega,\theta}(\lambda,\xi)+\hat F_{\Omega,r}(\lambda,\xi)\partial_\xi}G_\lambda(\theta,\xi)\d\xi,
\end{align*}
and since $G_\lambda(\theta,0)=0$, as $w_{11}(0)=w_{12}(0)=0$, while the boundary term at $\xi=\theta_b$ cancels the $\hat F_{\Omega,r}(\lambda,\theta_b)$ term of the formula for $\hat u_\theta$, we get
\begin{align*}
\hat u_\theta(\lambda,\theta) &= \int_0^{\theta_b} \lambda_1\brac{\lambda_{-1}\hat F_{\Omega,\theta}(\lambda,\xi)+\hat F_{\Omega,r}(\lambda,\xi)\partial_\xi}G_\lambda(\theta,\xi)\d\xi+ \partial_\xi G_\lambda(\theta,0) \lambda_1\hat F_{0,0}(\lambda)\\
&\quad- (3\lambda^2+1+\partial_\xi^2)\partial_\xi G_\lambda(\theta,0) \hat F_{0,1}(\lambda)  - \partial_\xi G_\lambda(\theta, \theta_b) \lambda_1\hat F_{1,0}(\lambda) - G_\lambda(\theta, \theta_b) \lambda_1\lambda_{-1}\hat F_{1,1}(\lambda)
\end{align*}
Using the divergence constraint \eqref{Mellin divergence} we can recover the formula for $\hat u_r$:
\begin{align*}
\hat u_r(\lambda,\theta) &= -\int_0^{\theta_b} \brac{\lambda_{-1}\hat F_{\Omega,\theta}(\lambda,\xi)+\hat F_{\Omega,r}(\lambda,\xi)\partial_\xi}\partial_\theta G_\lambda(\theta,\xi)\d\xi-\partial_\theta\partial_\xi G_\lambda(\theta,0)\hat F_{0,0}(\lambda)\\
&\quad + \lambda_1^{-1}(3\lambda^2+1+\partial_\xi^2)\partial_\theta\partial_\xi G_\lambda(\theta,0) \hat F_{0,1}(\lambda)  + \partial_\theta\partial_\xi G_\lambda(\theta, \theta_b)\hat F_{1,0}(\lambda) + \partial_\theta G_\lambda(\theta, \theta_b)\lambda_{-1}\hat F_{1,1}(\lambda)
\end{align*}
Taking the second derivative in $\theta$ we have
\begin{align*}
\partial_\theta^2\hat u_r(\lambda,\theta)& = -\int_0^{\theta_b} \brac{\lambda_{-1}\hat F_{\Omega,\theta}(\lambda,\xi)+\hat F_{\Omega,r}(\lambda,\xi)\partial_\xi}\partial_\theta^3 G_\lambda(\theta,\xi)\d\xi +\hat F_{\Omega,r}(\lambda,\theta)\\
&\quad -\partial_\theta^3\partial_\xi G_\lambda(\theta,0)\hat F_{0,0}(\lambda) + \lambda_1^{-1}(3\lambda^2+1+\partial_\xi^2)\partial_\theta^3\partial_\xi G_\lambda(\theta,0) \hat F_{0,1}(\lambda)\\
&\quad + \partial_\theta^3\partial_\xi G_\lambda(\theta, \theta_b)\hat F_{1,0}(\lambda) + \partial_\theta^3 G_\lambda(\theta, \theta_b)\lambda_{-1}\hat F_{1,1}(\lambda),
\end{align*}
where the second term $\hat F_{\Omega,r}(\lambda,\theta)$ on the RHS comes from the discontinuity of the third derivative of the Green's function and the fact that $\partial_\xi\partial_\theta^2 G_\lambda(\theta,\theta^+)-\partial_\xi\partial_\theta^2 G_\lambda(\theta,\theta^-)=1$. To see this we use the Green's function's jump condition $\partial_\theta^3G_\lambda(\theta,\theta^-)-\partial_\theta^3G_\lambda(\theta,\theta^+)=1$. Since the second derivative of the Green's function is continuous, we have $\partial_\theta^2G_\lambda(\theta,\theta^+)=\partial_\theta^2G_\lambda(\theta,\theta^-)$. Taking the total derivative of this with respect to $\theta$ using the multivariate chain rule gives $(\partial_\theta+\partial_\xi)\partial_\theta^2G_\lambda(\theta,\theta^+)=(\partial_\theta+\partial_\xi)\partial_\theta^2G_\lambda(\theta,\theta^-)$. Hence we get $\partial_\xi\partial_\theta^2 G_\lambda(\theta,\theta^+)-\partial_\xi\partial_\theta^2 G_\lambda(\theta,\theta^-)=1$.

Using \eqref{p from u} and \eqref{d p from u}, we get the formula for $\hat p$ and $\partial_\theta\hat p$:
\begin{align*}
\hat p(\lambda,\theta) &= -\int_0^{\theta_b} \brac{\hat F_{\Omega,\theta}(\lambda_1,\xi)+{1\over\lambda}\hat F_{\Omega,r}(\lambda_1,\xi)\partial_\xi}\partial_\theta G_{\lambda,p}(\theta,\xi)\d\xi\\
&\quad -{1\over\lambda}\partial_\xi\partial_\theta G_{\lambda,p}(\theta,0)\hat F_{0,0}(\lambda_1) + {1\over\lambda\lambda_2}(3\lambda_1^2+1+\partial_\xi^2)\partial_\xi\partial_\theta G_{\lambda,p}(\theta,0) \hat F_{0,1}(\lambda_1)\\
&\quad + {1\over\lambda}\partial_\xi\partial_\theta G_{\lambda,p}(\theta, \theta_b)\hat F_{1,0}(\lambda_1) + \partial_\theta G_{\lambda,p}(\theta, \theta_b)\hat F_{1,1}(\lambda_1)\\
\partial_\theta\hat p(\lambda,\theta) &= \int_0^{\theta_b} \lambda\brac{\lambda\hat F_{\Omega,\theta}(\lambda_1,\xi)+\hat F_{\Omega,r}(\lambda_1,\xi)\partial_\xi}G_{\lambda,p}(\theta,\xi)\d\xi - \hat F_{\Omega,\theta}(\lambda_1,\theta)\\
&\quad + \lambda\partial_\xi G_{\lambda,p}(\theta,0)\hat F_{0,0}(\lambda_1) - {\lambda\over\lambda_2}(3\lambda_1^2+1+\partial_\xi^2)\partial_\xi G_{\lambda,p}(\theta,0) \hat F_{0,1}(\lambda_1)\\
&\quad - \lambda\partial_\xi G_{\lambda,p}(\theta, \theta_b)\hat F_{1,0}(\lambda_1) - \lambda^2 G_{\lambda,p}(\theta, \theta_b)\hat F_{1,1}(\lambda_1),
\end{align*}
where $G_{\lambda,p}=(\lambda_2^2+\partial_\theta^2)G_{\lambda_1}$.

For generic $G$ and $F$ we have, by the argument of Schur's test,
\begin{align*}
\norm{\int_0^{\theta_b}G(\ph,\xi)F(\xi)\d\xi}_{L^2([0,\theta_b])}^2
&=\int_0^{\theta_b}\abs{\int_0^{\theta_b}G(\theta,\xi)F(\xi)\d\xi}^2\d\theta\\
&\leq\int_0^{\theta_b}\norm{|G(\theta,\xi)|^{1\over 2}}_{L^2_\xi([0,\theta_b])}^2\norm{|G(\theta,\xi)|^{1\over 2}|F(\xi)|}_{L^2_\xi([0,\theta_b])}^2\d\theta\\
&\leq\int_0^{\theta_b}\|G(\theta,\xi)\|_{L^1_\xi([0,\theta_b])}\norm{|G(\theta,\xi)|^{1\over 2}|F(\xi)|}_{L^2_\xi([0,\theta_b])}^2\d\theta\\
&\leq\|G(\theta,\xi)\|_{L^\infty_\theta L^1_\xi([0,\theta_b]^2)}\int_0^{\theta_b}\int_0^{\theta_b}|G(\theta,\xi)||F(\xi)|^2\d\xi\d\theta\\
&\leq\|G(\theta,\xi)\|_{L^\infty_\theta L^1_\xi([0,\theta_b]^2)}\|G(\theta,\xi)\|_{L^\infty_\xi L^1_\theta([0,\theta_b]^2)}\|F\|_{L^2([0,\theta_b])}^2
\end{align*}
Therefore, writing
\begin{align}\label{Upsilon shorthand}
\Upsilon_\lambda:=\<\lambda\>^{1\over 2}+\theta_b^{-{1\over 2}}
\end{align}
for the factor supplied by the trace estimate \eqref{Green's function trace estimate} --- the two summands of $\Upsilon_\lambda$ corresponding to the two summands of that estimate --- we have for $l=0,1,2$
\begin{align*}
\Gamma^H_{\Re\lambda,\theta_b}\|\partial_\theta^l\hat u_\theta\|_{L^2_\theta([0,\theta_b])}
&\lesssim{\<\lambda\>^{l}\over|\lambda|^2}\|\hat{\mb F}_\Omega\|_{L^2_\theta([0,\theta_b])} +{\<\lambda\>^{l}\Upsilon_\lambda\over|\lambda|^2}\brac{|\hat F_{0,0}|+|\hat F_{1,0}|+|\hat F_{1,1}|} +{\<\lambda\>^{l+1}\Upsilon_\lambda\over|\lambda|^2}|\hat F_{0,1}|\\
\Gamma^H_{\Re\lambda,\theta_b}\|\partial_\theta^l\hat u_r\|_{L^2_\theta([0,\theta_b])}
&\lesssim{\<\lambda\>^{l}\over|\lambda|^2}\|\hat{\mb F}_\Omega\|_{L^2_\theta([0,\theta_b])} +{\<\lambda\>^{l}\Upsilon_\lambda\over|\lambda|^2}\brac{|\hat F_{0,0}|+|\hat F_{1,0}|+|\hat F_{1,1}|} +{\<\lambda\>^{l+2}\Upsilon_\lambda\over|\lambda|^2|\lambda_1|}|\hat F_{0,1}|
\end{align*}
Note that the $\hat{\mb F}_\Omega$ terms come from \eqref{Green's function L1 estimate} through Schur's test above and are therefore unaffected; it is only the boundary data, which are paired against traces of the Green's function, that pick up $\Upsilon_\lambda$. For the pressure we have
\begin{align*}
\Gamma^H_{\Re\lambda_1,\theta_b}\|\hat p\|_{L^2_\theta([0,\theta_b])}
&\lesssim{\<\lambda_1\>^{2}\over|\lambda_1|^2|\lambda|}\|\hat{\mb F}_\Omega(\lambda_1)\|_{L^2_\theta([0,\theta_b])} +{\<\lambda_1\>^{2}\Upsilon_{\lambda_1}\over|\lambda_1|^2|\lambda|}\brac{|\hat F_{0,0}(\lambda_1)|+|\hat F_{1,0}(\lambda_1)|}\\
&\quad +{\<\lambda_1\>\Upsilon_{\lambda_1}\over|\lambda_1|^2}|\hat F_{1,1}(\lambda_1)| +{\<\lambda\>^{4}\Upsilon_{\lambda}\over|\lambda_1|^2|\lambda_2||\lambda|}|\hat F_{0,1}(\lambda_1)|\\
\Gamma^H_{\Re\lambda_1,\theta_b}\|\partial_\theta\hat p\|_{L^2_\theta([0,\theta_b])}
&\lesssim{\<\lambda_1\>^{2}\over|\lambda_1|^2}\|\hat{\mb F}_\Omega(\lambda_1)\|_{L^2_\theta([0,\theta_b])} +{\<\lambda_1\>^{2}\Upsilon_{\lambda_1}\over|\lambda_1|^2}\brac{|\hat F_{0,0}(\lambda_1)|+|\hat F_{1,0}(\lambda_1)|}\\
&\quad +{\<\lambda_1\>\Upsilon_{\lambda_1}|\lambda|\over|\lambda_1|^2}|\hat F_{1,1}(\lambda_1)| +{\<\lambda\>^{4}\Upsilon_{\lambda}\over|\lambda_1|^2|\lambda_2|}|\hat F_{0,1}(\lambda_1)|
\end{align*}
Using the inequality $\<\lambda\>/|\lambda|\leq\<\Re\lambda\>/|\Re\lambda|=\Gamma_{\Re\lambda}$ we have
\begin{align}
  &\sum_{l=0}^2\<\lambda\>^{2-l}\|\partial_\theta^l\hat{\mb u}(\lambda)\|_{L^2_\theta([0,\theta_b])} +\sum_{l=0}^1\<\lambda\>^{1-l}\Gamma_{\Re\lambda_{-1}}^{l-1}\|\partial_\theta^l\hat p(\lambda_{-1})\|_{L^2_\theta([0,\theta_b])}\nonumber\\ &\lesssim{\Gamma_{\Re\lambda}^2\over\Gamma^H_{\Re\lambda,\theta_b}} \brac{\|\hat{\mb F}_\Omega\|_{L^2_\theta([0,\theta_b])} +\Upsilon_\lambda\brac{|\hat F_{0,0}|+|\hat F_{1,0}|+|\hat F_{1,1}|} +\Gamma_{\Re\lambda_1}\<\lambda\>\Upsilon_\lambda|\hat F_{0,1}|}\label{Mellin bound fixed lambda}
\end{align}
So for $j\geq 0$ we have
\begin{align*}
  &\|\mb u\|_{\Omega,a,j,2}^2+\Gamma_{a-2}^{-2}\|p\|_{\Omega,a-1,j,1}^2\\
   &\lesssim \sum_{j'=0}^j\brac{\sum_{l=0}^2\int_0^\infty\|\dr^{j'}(\dr,\partial_\theta)^l\mb u\|_{L^2_\theta([0,\theta_b])}^2r^{1-2a}\d r +\Gamma_{a-2}^{-2}\sum_{l=0}^1\int_0^\infty\|\dr^{j'}(\dr,\partial_\theta)^lp\|_{L^2_\theta([0,\theta_b])}^2r^{3-2a}\d r}\\
  &\lesssim \sum_{l=0}^2\int_{\Re\lambda=a-1}\|\<\lambda\>^{2+j-l}\partial_\theta^l\hat{\mb u}\|_{L^2_\theta([0,\theta_b])}^2\d\Im\lambda +\sum_{l=0}^1\int_{\Re\lambda=a-2}\Gamma_{\Re\lambda}^{-2}\|\<\lambda\>^{1+j-l}\partial_\theta^l\hat p\|_{L^2_\theta([0,\theta_b])}^2\d\Im\lambda\\
  &\lesssim \sum_{l=0}^2\int_{\Re\lambda=a-1}\|\<\lambda\>^{2+j-l}\partial_\theta^l\hat{\mb u}(\lambda)\|_{L^2_\theta([0,\theta_b])}^2\d\Im\lambda\\ 
  &\quad+\sum_{l=0}^1\int_{\Re\lambda=a-1}\Gamma_{\Re\lambda_{-1}}^{-2}\|\<\lambda\>^{1+j-l}\partial_\theta^l\hat p(\lambda_{-1})\|_{L^2_\theta([0,\theta_b])}^2\d\Im\lambda\\
  &\lesssim\int_{\Re\lambda=a-1}{\Gamma_{\Re\lambda}^4\over(\Gamma^H_{\Re\lambda,\theta_b})^2}\Big(\|\<\lambda\>^j\hat{\mb F}_\Omega\|_{L^2_\theta([0,\theta_b])} +\<\lambda\>^{j}\Upsilon_\lambda\brac{|\hat F_{0,0}|+|\hat F_{1,0}|+|\hat F_{1,1}|}\nonumber\\
  &\pushright{+\Gamma_{\Re\lambda_1}\<\lambda\>^{j+1}\Upsilon_\lambda|\hat F_{0,1}|\Big)^2\d\Im\lambda}\\
  &\lesssim{\Gamma_{a-1}^4\over(\Gamma^H_{a-1,\theta_b})^2}\bigg(\|\mb F_\Omega\|_{\Omega,a,j}^2+\|F_{0,0}\|_{\partial_0\Omega,a-{1\over 2},j+{1\over 2}}^2+\|(F_{1,0},F_{1,1})\|_{\partial_1\Omega,a-{1\over 2},j+{1\over 2}}^2+\Gamma_{a}^2\|F_{0,1}\|_{\partial_0\Omega,a-{1\over 2},j+{3\over 2}}^2\bigg)\\
  &\quad+{\Gamma_{a-1}^4\over\theta_b(\Gamma^H_{a-1,\theta_b})^2}\bigg(\|F_{0,0}\|_{\partial_0\Omega,a-{1\over 2},j}^2+\|(F_{1,0},F_{1,1})\|_{\partial_1\Omega,a-{1\over 2},j}^2+\Gamma_{a}^2\|F_{0,1}\|_{\partial_0\Omega,a-{1\over 2},j+1}^2\bigg).\qedhere
\end{align*}
\end{proof}

\begin{remark}[The excluded weights and the constants in \eqref{Mellin estimate initial time}]\label{Mellin initial time remark}
The weights $a=0$ and $a=2$ are excluded in Proposition \ref{Mellin system initial time} because factors of its solution formulae are singular on the Mellin lines there: the factor $\lambda_1^{-1}$ of the velocity formulae on the line $\Re\lambda=a-1$ at $a=0$, and the factors $\lambda^{-1}$ and $(\lambda\lambda_2)^{-1}$ of the pressure formulae on the line $\Re\lambda=a-2$ at $a=2$ and at $a=0$. The estimate bounds these factors separately from the kernels they multiply, and accordingly $\Gamma_a$ blows up at $a=0$ and $\Gamma_{a-2}$ at $a=2$. Every application below has $a-1\notin\Z$. The factor $\Gamma^H_{a-1,\theta_b}$ is positive by Lemma \ref{H lemma}(ii), bounded below by an absolute constant when $|a-1|\theta_b\leq\pi/4$ by Lemma \ref{H lemma}(iii), and degenerates only as $|a-1|\uparrow\Lambda^*(\theta_b)$. The terms of \eqref{Mellin estimate initial time} carrying the factor $\theta_b^{-1/2}$ come from the trace estimate \eqref{Green's function trace estimate}; they have half a radial derivative fewer on the boundary data than the corresponding terms without that factor, as in Lemma \ref{Neumann wedge problem}. Each of these terms is bounded by $\theta_b^{-1/2}$ times the corresponding term without that factor, and so they require no further regularity of the data.
\end{remark}

\subsubsection{The stationary system and the corner ladder}\label{stationary system and corner ladder}

To apply Proposition \ref{Mellin system initial time}, we first need to already have a $\mb u$ that solves the original system at a fixed given time. Since the initial data we have is not in the trace space of the positive-time estimate, as explained in Remark \ref{linear theory closing remark}, we need an additional existence result for $(\mb u,p)$ with respect to $(\phi,\mb F,\mb f)$ at one fixed time.

\begin{proposition}[The stationary system]\label{stationary system}
Let $\kappa\in(0,1/2]$ and $l\in\N_0$, and let $\theta_b>0$ be small enough, depending on $\kappa$ and $l$. If $(\phi,\mb F,\mb f)$ have regularity that makes the RHS of \eqref{data estimate} finite, then there exist $(\mb u,p)$ satisfying a.e.
\begin{subequations}\label{system: stationary linear}
\begin{alignat}{3}
\lpc\mb u-\grad p&=\mb F\qquad&&\text{on}\quad\Omega\label{stationary bulk equation}\\
\grad\cdot\mb u&=0\qquad&&\text{on}\quad\Omega\label{stationary divergence}\\
u_\theta&=0\qquad&&\text{on}\quad\partial_0\Omega\label{stationary no penetration}\\
u_r-r^{-1}\partial_\theta u_r&=0\qquad&&\text{on}\quad\partial_0\Omega\label{stationary Navier slip}\\
\mathbb T[\mb u,p]\mb e_\theta-(\sin\theta_b)^2(\partial_r\phi)\mb e_\theta&=\mb f\qquad&&\text{on}\quad\partial_1\Omega\label{stationary surface tension}
\end{alignat}
\end{subequations}
with regularity
\begin{align}
&\|\mb u\|_{\Omega,1-\kappa,l,2}+\|r\grad\mb u\|_{\Omega,1,l,1}+\sum_{k=0}^1\|p\|_{\Omega,-\kappa k,l,1}\nonumber\\
&\lesssim_{\theta_b,\kappa}\sum_{k=0}^1\brac{\|\mb F\|_{\Omega,-1-\kappa k,l}+\|\mb f\|_{\partial_1\Omega,-{1\over 2}-\kappa k,l+{1\over 2}}+(\sin\theta_b)^2\|\phi\|_{\partial_1\Omega,{1\over 2}-\kappa k,l+{3\over 2}}}\label{data estimate}
\end{align}
and this solution is unique within this class of norms. Its Navier slip term $\sum_{k=0}^1\|u_r\|_{\partial_0\Omega,-\kappa k,l+1}$ is moreover bounded by the same RHS.
\end{proposition}
\begin{proof}
In \eqref{system: stationary linear} the profile is a datum. The system is therefore System \ref{linear system} with the surface tension of the unknown profile switched off, with the profile entering without the factor $s^{-1}$ that its elimination produces in \eqref{surface tension without profile}, and without the two conditions \eqref{contact angle without profile} on $\partial_ru_\theta$ at the contact point and at infinity. In addition, $\mb F$ stands in place of $\P\mb F$, which by Remark \ref{P remark} does not change the velocity. Section \ref{linear existence} applies with three changes. First, in $B^r_{a,j}$ of Definition \ref{energy forms and functionals} the term carrying $\bar s^{-1}$ is dropped and $p_{\P|\mb u}$ is replaced by $p_{\partial_0}+p_{\partial_1,0}$. Second, in $\mathcal F^r_{a,j}$ and $\mathcal F^\theta_{a,j}$ the pressure $p_{\partial_1,\phi_0}$ is replaced by $p_{\partial_1,\phi}:=\E(-(\sin\theta_b)^2\partial_r\phi)$. Third, $\|\ph\|_{H_s^l}$ loses its $\partial_1\Omega$ term, which is the energy of the dropped term, together with the weight $1+|s|^{-2\kappa}$. That weight served only to absorb the two remainders attached to the dropped term, those of Lemma \ref{surface tension weighted remainder} and of the second estimate of Lemma \ref{pressure remainder}, so $\kappa\in(0,1/2]$ suffices here. The identities and bounds of Proposition \ref{energy forms}, the coercivity of Proposition \ref{energy form coercivity} and the test function problem of Subsections \ref{the test function problem} and \ref{solvability of the test function problem} then go through, with the $\bar s^{-1}$ term deleted throughout and the profile terms of (iii)--(iv) replaced by $(\sin\theta_b)^2\|\phi\|_{\partial_1\Omega,1/2-\kappa k,\ph}$ at the two weights $k\in\{0,1\}$. The corner term of the profile pairing is absent, since $\phi(0)=0$, which the finiteness of $\|\phi\|_{\partial_1\Omega,1/2,l+3/2}$ in \eqref{data estimate} forces.

The Lax--Milgram theorem and the argument of Theorem \ref{linear solution Laplace side} then give \eqref{data estimate}, with $p:=p_{\partial_0}+p_{\partial_1,0}+p_{\partial_1,\phi}+p_{\partial_1,f_\theta}-q[\mb F]$ as in Corollary \ref{profile on the Laplace side}, the Navier slip term being a summand of the modified energy. For the uniqueness, a solution in the class of \eqref{data estimate} lies in the modified space, its $\partial_0\Omega$ term being that Navier slip term. Such a solution satisfies the variational identity by the argument of Lemma \ref{stationary as variational} with the surface tension term deleted; the contact angle condition entered there only through the endpoint term that term produced. Two solutions in the class therefore coincide with the unique variational solution.
\end{proof}

\begin{remark}[The norms of the stationary solution]\label{stationary solution norms}
The bound on the Navier slip term in Proposition \ref{stationary system} is used in Lemma \ref{stationary as variational}. The estimate \eqref{data estimate} controls $\|r\grad\mb u\|_{\Omega,1,l,1}$ and not $\|\mb u\|_{\Omega,1,l,2}$, since $\mb u(\mb 0)$ is in general nonzero.
\end{remark}

The next proposition details the gain of spatial scaling using the Mellin solution representations we saw earlier.

Before stating it we record the notation for the expansion of the profile at the contact point. If $\hat\phi(\ph,\theta_b)$ is meromorphic across $\Re\lambda=1$ with at worst a double pole there, we write $\smash[b]{c_{\phi,\dr_{-1}}}$ and $\smash[b]{c_{\phi,\dr_{-1}^2}}$ for the coefficients defined by
\begin{align}\label{profile corner coefficients}
c_{\phi,\dr_{-1}}r+c_{\phi,\dr_{-1}^2}r\ln r:=-\sqrt{2\pi}\Res_{\lambda=1}\brac{r^\lambda\hat\phi(\lambda,\theta_b)},
\end{align}
the sign being chosen so that, by the residue theorem, if $\phi$ is the inverse Mellin transform of $\hat\phi(\ph,\theta_b)$ on a line just to the left of $\Re\lambda=1$, then $\phi-\smash[b]{c_{\phi,\dr_{-1}}}r-\smash[b]{c_{\phi,\dr_{-1}^2}}r\ln r$\vspace{1.6pt} is its inverse Mellin transform on a line just to the right. Subtracting these two terms, as in $\phi^{[2]}$ below, therefore removes the pole and leaves a profile controlled at the higher scaling.

\begin{proposition}[The corner ladder]\label{corner ladder}
Let $\theta_b\in(0,\pi/2)$ and $\kappa\in(0,1)$ with $2-\kappa<\Lambda^*(\theta_b)$, where $\Lambda^*(\theta_b)$ is as in Definition \ref{zero-free strip}, and let $j\geq 1$. Suppose we have $(\mb u,p,\phi,\mb F,\mb f)$ satisfying \eqref{system: stationary linear} a.e. with regularity
\begin{align*}
\|\mb u\|_{\Omega,1-\kappa,j,2} +\|p\|_{\Omega,-\kappa,j,1} +\|\mb F\|_{\Omega,-1-\kappa,j} +\|\mb f\|_{\partial_1\Omega,-{1\over 2}-\kappa,j+{1\over 2}} +\|\dr\phi\|_{\partial_1\Omega,{1\over 2}-\kappa,j+{1\over 2}} &< \infty.
\end{align*}
Let $\mb u^{\smash[t]{[M]}}$, $p^{\smash[t]{[M]}}$ and $\phi^{\smash[t]{[M]}}$, $M=0,1,2$, and the degree one fields $\smash[b]{\mb P_{\mb u_{F_{0,0}},\dr_{-1}}}$ and $\smash[b]{\mb P_{\mb u_{F_{1,1}},\dr_{-1}}}$ be given by
\begin{subequations}\label{data step relations}
\begin{align}
  \mb u^{[0]} &= \mb u\\
  \mb u^{[1]} &= \mb u-\mb u(\mb 0)\\
  \mb u^{[2]} &=\mb u - \mb u(\mb 0) - {u_1(\mb 0)\over 2}\mb P_{\mb u_{F_{0,0}},\dr_{-1}}-{1\over 4}(\sin\theta_b)^2c_{\phi,\dr_{-1}^2}\mb P_{\mb u_{F_{1,1}},\dr_{-1}}\label{ladder u2}\\
  p^{[0]} &= p^{[1]} = p\\
  p^{[2]} &= p +{u_1(\mb 0)\over\sin(2\theta_b)} +(\sin\theta_b)^2c_{\phi,\dr_{-1}}+{1\over 2}(\sin\theta_b)^2c_{\phi,\dr_{-1}^2}(2\ln r+1)\label{ladder p2}\\
  \phi^{[0]} &= \phi^{[1]} = \phi\\
  \phi^{[2]} &= \phi - c_{\phi,\dr_{-1}}r -c_{\phi,\dr_{-1}^2}r\ln r\label{ladder phi2}\\
  \mb P_{\mb u_{F_{0,0}},\dr_{-1}} &= r\brac{\brac{{\sin(2(\theta_b-\theta))\over\sin(2\theta_b)}-1}\mb e_\theta+\brac{\cos(2(\theta_b-\theta))\over\sin(2\theta_b)}\mb e_r}\label{ladder P F00}\\
  \mb P_{\mb u_{F_{1,1}},\dr_{-1}} &= r(2\theta\mb e_\theta -\mb e_r).\label{ladder P F11}
\end{align}
\end{subequations}
If $(\phi,\mb F,\mb f)$ have the additional regularity that makes the RHS below finite, then $(\mb u,p)$ also have the regularity of the LHS:
\begin{align}
\|\mb u^{[M]}&\|_{\Omega,M+1-\kappa,j,2} +\Gamma_{M-1-\kappa}^{-1}\|p^{[M]}\|_{\Omega,M-\kappa,j,1}\nonumber\\
  &\lesssim{\Gamma_{M-\kappa}^2\over\theta_b^{1\over 2}\Gamma^H_{M-\kappa,\theta_b}} \bigg(\|\mb F\|_{\Omega,M-1-\kappa,j} +\|u_r^{[M^-]}\|_{\partial_0\Omega,M-{1\over 2}-\kappa,j+{1\over 2}}\nonumber\\
  &\qquad\qquad\qquad\qquad+\|\mb f\|_{\partial_1\Omega,M-{1\over 2}-\kappa,j+{1\over 2}} +(\sin\theta_b)^2\|\dr\phi^{[M]}\|_{\partial_1\Omega,M+{1\over 2}-\kappa,j+{1\over 2}}\bigg)\label{um data reg master}
\end{align}
for $M=0,1,2$, where $M^-=\max\{M-1,0\}$. If the RHS of \eqref{um data reg master} is finite at $M=1$, then $\mb u$ is continuous at $\mb 0$, with $\mb u(\mb 0)=u_1(\mb 0)\mb e_1$.
Moreover, if the compatibility condition
\begin{align}
  c_{\phi,\dr_{-1}^2} &= {u_1(\mb 0)\over\theta_b(\sin\theta_b)^2}\label{C1 in the linear theory}
\end{align}
is satisfied, then the contact angle condition $\partial_ru_\theta(0,\theta_b)=0$ holds, as a limit along $\partial_1\Omega$.
\end{proposition}
\begin{proof}
We have from \eqref{system: stationary linear}, after multiplying the bulk equation by $r^2$, the system of Proposition \ref{Mellin system initial time}, the Navier slip term $ru_r$ being a datum: it has a lower scaling than $\partial_\theta u_r$, and as noted after Definition \ref{scaling definition} such terms are treated as data. Through this datum the trace of $u_r$ passes from one rung to the next in \eqref{um data reg master}.

For $M\in\{0,1\}$ let
\begin{subequations}\label{M01 data def}
\begin{align}
  \mb F_\Omega^{[M]} &= r^2\mb F\\
  F_{0,0}^{[M]} &= ru_r\label{M01 slip datum}\\
  F_{0,1}^{[M]} &= 0\\
  F_{1,0}^{[M]} &= rf_r\\
  F_{1,1}^{[M]} &= rf_\theta+(\sin\theta_b)^2\dr\phi,\label{M01 stress datum}
\end{align}
\end{subequations}
or in Mellin variables
\begin{subequations}\label{M01 Mellin data}
\begin{align}
  \hat{\mb F}_\Omega^{[M]}(\lambda,\theta) &= \hat{\mb F}(\lambda_{-2},\theta)\\
  \hat F_{0,0}^{[M]}(\lambda) &= \hat u_r(\lambda_{-1},0)\\
  \hat F_{0,1}^{[M]}(\lambda) &= 0\\
  \hat F_{1,0}^{[M]}(\lambda) &= \hat f_r(\lambda_{-1})\\
  \hat F_{1,1}^{[M]}(\lambda) &= \hat f_\theta(\lambda_{-1})+(\sin\theta_b)^2\lambda\hat\phi(\lambda,\theta_b).
\end{align}
\end{subequations}
By Proposition \ref{Mellin system initial time} there exist $(\mb u^{\smash[t]{[M]}},p^{\smash[t]{[M]}})$ such that
\begin{subequations}\label{data step equations}
\begin{alignat}{3}
r^2(\lpc\mb u^{[M]}-\grad p^{[M]})&=\mb F_\Omega^{[M]}\qquad&&\text{on}\qquad\Omega\label{data step bulk equation}\\
\grad\cdot\mb u^{[M]}&=0\qquad&&\text{on}\qquad\Omega\\
\partial_\theta u^{[M]}_r+\dr_{-1}u^{[M]}_\theta&=F_{0,0}^{[M]}\qquad&&\text{on}\qquad\partial_0\Omega\label{data step slip condition}\\
u^{[M]}_\theta&=F_{0,1}^{[M]}\qquad&&\text{on}\qquad\partial_0\Omega\\
\partial_\theta u^{[M]}_r+\dr_{-1}u^{[M]}_\theta&=F_{1,0}^{[M]}\qquad&&\text{on}\qquad\partial_1\Omega\\
2(\partial_\theta u^{[M]}_\theta+u^{[M]}_r) -rp^{[M]}
&=F_{1,1}^{[M]}\qquad&&\text{on}\qquad\partial_1\Omega\label{data step stress condition}
\end{alignat}
\end{subequations}
with regularity
\begin{align}
  &\|\mb u^{[M]}\|_{\Omega,M+1-\kappa,j,2}+\Gamma_{M-1-\kappa}^{-1}\|p^{[M]}\|_{\Omega,M-\kappa,j,1}\nonumber\\
  &\lesssim{\Gamma_{M-\kappa}^2\over\theta_b^{1\over 2}\Gamma^H_{M-\kappa,\theta_b}}\brac{\|\mb F_\Omega^{[M]}\|_{\Omega,M+1-\kappa,j}+\|F_{0,0}^{[M]}\|_{\partial_0\Omega,M+{1\over 2}-\kappa,j+{1\over 2}}+\|(F_{1,0}^{[M]},F_{1,1}^{[M]})\|_{\partial_1\Omega,M+{1\over 2}-\kappa,j+{1\over 2}}}\nonumber\\
  &\lesssim{\Gamma_{M-\kappa}^2\over\theta_b^{1\over 2}\Gamma^H_{M-\kappa,\theta_b}}\bigg(\|\mb F\|_{\Omega,M-1-\kappa,j}+\|u_r\|_{\partial_0\Omega,M-{1\over 2}-\kappa,j+{1\over 2}}\nonumber\\
  &\pushright{+\|\mb f\|_{\partial_1\Omega,M-{1\over 2}-\kappa,j+{1\over 2}}+(\sin\theta_b)^2\|\dr\phi\|_{\partial_1\Omega,M+{1\over 2}-\kappa,j+{1\over 2}}\bigg).\hspace{1cm}}\label{data um regularity}
\end{align}
Here the first inequality is \eqref{Mellin estimate initial time} at $a=M+1-\kappa$, its terms carrying the factor $\theta_b^{-1/2}$ being bounded as in Remark \ref{Mellin initial time remark}, and the second follows from \eqref{M01 data def}. At $M=0$ the norm $\|u_r\|_{\partial_0\Omega,-1/2-\kappa,j+1/2}$ on the right is finite under the hypotheses: on $\partial_0\Omega$ the Navier slip condition of \eqref{system: stationary linear} reads $ru_r=\partial_\theta u_r$, so that $\|u_r\|_{\partial_0\Omega,-1/2-\kappa,j+1/2}\lesssim_j\|\partial_\theta u_r\|_{\partial_0\Omega,1/2-\kappa,j+1/2}\lesssim_{\theta_b}\|\mb u\|_{\Omega,1-\kappa,j,2}$ by Lemma \ref{trace}. This proves \eqref{um data reg master} for $M=0,1$. 

From Proposition \ref{Mellin system initial time} we have the representation formula in Mellin variables
\begin{align*}
  \hat u_\theta(\lambda,\theta) &= \int_0^{\theta_b} \lambda_1\brac{\lambda_{-1}\hat F_{\Omega,\theta}(\lambda,\xi)+\hat F_{\Omega,r}(\lambda,\xi)\partial_\xi}G_\lambda(\theta,\xi)\d\xi+ \partial_\xi G_\lambda(\theta,0) \lambda_1\hat F_{0,0}(\lambda)\\
&\quad  - \partial_\xi G_\lambda(\theta, \theta_b) \lambda_1\hat F_{1,0}(\lambda) - G_\lambda(\theta, \theta_b) \lambda_1\lambda_{-1}\hat F_{1,1}(\lambda)\\
\hat u_r(\lambda,\theta) &= -\int_0^{\theta_b} \brac{\lambda_{-1}\hat F_{\Omega,\theta}(\lambda,\xi)+\hat F_{\Omega,r}(\lambda,\xi)\partial_\xi}\partial_\theta G_\lambda(\theta,\xi)\d\xi-\partial_\theta\partial_\xi G_\lambda(\theta,0)\hat F_{0,0}(\lambda)\\
&\quad  + \partial_\theta\partial_\xi G_\lambda(\theta, \theta_b)\hat F_{1,0}(\lambda) + \partial_\theta G_\lambda(\theta, \theta_b)\lambda_{-1}\hat F_{1,1}(\lambda)\\
\hat p(\lambda,\theta) &= -\int_0^{\theta_b} \brac{\hat F_{\Omega,\theta}(\lambda_1,\xi)+{1\over\lambda}\hat F_{\Omega,r}(\lambda_1,\xi)\partial_\xi}\partial_\theta G_{\lambda,p}(\theta,\xi)\d\xi -{1\over\lambda}\partial_\xi\partial_\theta G_{\lambda,p}(\theta,0)\hat F_{0,0}(\lambda_1)\\
&\quad + {1\over\lambda}\partial_\xi\partial_\theta G_{\lambda,p}(\theta, \theta_b)\hat F_{1,0}(\lambda_1) + \partial_\theta G_{\lambda,p}(\theta, \theta_b)\hat F_{1,1}(\lambda_1),
\end{align*}
where $G_\lambda$ is the Green's function of Lemma \ref{Green's function initial time} and $G_{\lambda,p}=(\lambda_2^2+\partial_\theta^2)G_{\lambda_1}$ as in Proposition \ref{Mellin system initial time}. The representation is the same for $M=0$ and $M=1$, so we can skip writing the superscript. Since $(\mb u,p)$ solves \eqref{data step equations} at $M=0$, as noted at the beginning of the proof, and $\|\mb u\|_{\Omega,1-\kappa,j,2}+\|p\|_{\Omega,-\kappa,j,1}$ is finite by hypothesis, the uniqueness in Proposition \ref{Mellin system initial time} at $a=1-\kappa$ gives $\mb u^{[0]}=\mb u$ and $p^{[0]}=p$. 

To show that the solutions of \eqref{data step equations} are the fields \eqref{data step relations}, we move the line of integration for the Mellin inversion across the poles lying between the two levels. This is legitimate for three reasons. First, the hypothesis $2-\kappa<\Lambda^*(\theta_b)$ places the closed strip $-\kappa\leq\Re\lambda\leq 2-\kappa$ inside $|\Re\lambda|<\Lambda^*(\theta_b)$, so by Lemma \ref{H lemma}(i) the only singularity of $G_\lambda$ in this strip is the double pole at $\lambda=0$. Second, the decay of $G_\lambda$ in $\Im\lambda$ given by \eqref{Green's function pointwise estimate} is uniform on the strip, since the factor $1/\Gamma^H_{\Re\lambda,\theta_b}$ in \eqref{Green's function pointwise estimate} is bounded there by Lemma \ref{H lemma}(ii). The transforms of the data are square integrable on the lines uniformly in the line, so together with this decay the two horizontal sides of the rectangular contour at $\Im\lambda=\pm Y$ drop out along a sequence $Y\to\infty$. Third, the poles at $\lambda=1$ met below come from the data $\hat F_\bullet$ and not from $G_\lambda$, their holomorphy being checked case by case. Throughout, if $\hat{\mb u}$ has a pole of order $m$ at $\lambda=d$ and no other singularity between the two lines, then
\begin{align}\label{residue rule}
{1\over\sqrt{2\pi}}\brac{\int_{\Re\lambda=d+\epsilon}-\int_{\Re\lambda=d-\epsilon}}r^\lambda\hat{\mb u}(\lambda,\theta)\d\Im\lambda&=\sqrt{2\pi}\Res_{\lambda=d}\brac{r^\lambda\hat{\mb u}}\nonumber\\
&={\sqrt{2\pi}\over(m-1)!}\lim_{\lambda\to d}\partial_\lambda^{m-1}\brac{(\lambda-d)^me^{\lambda\ln r}\hat{\mb u}},
\end{align}
and likewise for $\hat p$; we evaluate the RHS from the multipliers of Proposition \ref{Mellin system initial time} and their $\partial_\lambda$ derivatives at $\lambda=d$, the powers of $\ln r$ coming from the expansion of $e^{\lambda\ln r}$. 

At $\lambda=0$, where $G_\lambda$ has a double pole, inserting the limits of $\lambda^2G_\lambda$ of Lemma \ref{Green's function initial time} in the formula for $\hat u_\theta$ gives
\begin{align*}
  \lim_{\lambda\to 0}(\lambda^2\hat u_\theta(\lambda,\theta))(\sin(2\theta_b)+2\theta_b)
  &= -\int_0^{\theta_b}\brac{-\hat F_{\Omega,\theta}(0,\xi)\sin\xi+\hat F_{\Omega,r}(0,\xi)\cos\xi}\sin\theta\d\xi\\
&\quad - \hat F_{0,0}(0)\sin\theta + \hat F_{1,0}(0)\cos\theta_b\sin\theta - \hat F_{1,1}(0)\sin\theta_b\sin\theta\\
&= c_2\sin\theta,\\
  \lim_{\lambda\to 0}\partial_\lambda(\lambda^2\hat u_\theta(\lambda,\theta))(\sin(2\theta_b)+2\theta_b)
  &= (c_2+c_1)\sin\theta,
\end{align*}
where the first limit defines the constant $c_2$ and the second the constant $c_1$: since $\lim_{\lambda\to 0}\partial_\lambda(\lambda^2G_\lambda)=0$, the derivative $\partial_\lambda$ falls only on the multipliers and on the data, so that the second limit is again a multiple of $\sin\theta$. We need not compute $c_1$. By \eqref{Mellin divergence}, $\lambda^2\hat u_r=-\lambda_1^{-1}\partial_\theta(\lambda^2\hat u_\theta)$, so that
\begin{align*}
\lim_{\lambda\to 0}(\lambda^2\hat u_r(\lambda,\theta))(\sin(2\theta_b)+2\theta_b)
 &= -c_2\cos\theta,\\
  \lim_{\lambda\to 0}\partial_\lambda(\lambda^2\hat u_r(\lambda,\theta))(\sin(2\theta_b)+2\theta_b)
  &= -c_1\cos\theta.
\end{align*}
Since $\hat{\mb u}$ has at most a double pole at $\lambda=0$, \eqref{residue rule} with $m=2$ gives the sum of $\sqrt{2\pi}\ln r\lim_{\lambda\to 0}(\lambda^2\hat{\mb u})$ and $\sqrt{2\pi}\lim_{\lambda\to 0}\partial_\lambda(\lambda^2\hat{\mb u})$, and we get
\begin{align*}
  (\mb u^{[1]}-\mb u^{[0]})(r,\theta) &=\sqrt{2\pi}\Res_{\lambda=0}\brac{r^\lambda\hat{\mb u}} =c_{\mb u,\dr}\mb P_{\mb u,\dr} +c_{\mb u,\dr^2}\mb P_{\mb u,\dr^2},
\end{align*}
where
\begin{align*}
  \mb P_{\mb u,\dr} &= \cos\theta\mb e_r-\sin\theta\mb e_\theta = \mb e_1,\\
  \mb P_{\mb u,\dr^2} &= \mb e_1\ln r-\mb e_\theta\sin\theta,
\end{align*}
with $\smash[b]{c_{\mb u,\dr}}=-\sqrt{2\pi}c_1/(\sin(2\theta_b)+2\theta_b)$, and, by the limits above and \eqref{M01 Mellin data}, the coefficient $\smash[b]{c_{\mb u,\dr^2}}$ of $\smash[b]{\mb P_{\mb u,\dr^2}}$ with the singular $\ln r$ is given by
\begin{align*}
&{(\sin(2\theta_b)+2\theta_b)c_{\mb u,\dr^2}\over\sqrt{2\pi}}=-c_2\\
&=\int_0^{\theta_b}\brac{-\hat F_{\Omega,\theta}(0,\xi)\sin\xi+\hat F_{\Omega,r}(0,\xi)\cos\xi}\d\xi + \hat F_{0,0}(0) - \hat F_{1,0}(0)\cos\theta_b + \hat F_{1,1}(0)\sin\theta_b\\
&=\int_0^{\theta_b}\brac{-\hat F_{\theta}(-2,\xi)\sin\xi+\hat F_{r}(-2,\xi)\cos\xi}\d\xi + \hat u_r(-1,0) - \hat f_r(-1)\cos\theta_b\\
&\quad + \brac{\hat f_\theta(-1)+(\sin\theta_b)^2(\dr\phi)^\wedge(0)}\sin\theta_b.
\end{align*}

Substituting in the definition of the Mellin transform as an integral of the untransformed functions, we can then use the equations \eqref{system: stationary linear} to show that the coefficient is in fact $0$, and hence that there is no blow-up of the contact-point velocity:
\begin{align*}
(\sin(2\theta_b)+2\theta_b)c_{\mb u,\dr^2}
&=\int_0^\infty\int_0^{\theta_b}\brac{-F_{\theta}(r,\xi)\sin\xi+F_{r}(r,\xi)\cos\xi}r\d\xi\d r + \int_0^\infty u_r(r,0)\d r\\
&\quad - \int_0^\infty f_r(r)\cos\theta_b\d r + \int_0^\infty\brac{f_\theta(r)+(\sin\theta_b)^2(\partial_r\phi)(r)}\sin\theta_b\d r\\
&=\int_0^\infty\int_0^{\theta_b}\brac{\mb F\cdot\mb e_1}r\d\xi\d r + \int_0^\infty u_r(r,0)\d r\\
&\quad- \int_0^\infty(\mb f(r)+(\sin\theta_b)^2\partial_r\phi(r)\mb e_\theta)\cdot\mb e_1\d r\\
&=\int_\Omega\brac{\mb F\cdot\mb e_1}\d\mb x + \int_{\partial_0\Omega}\mb u\cdot\mb e_1\d r - \int_{\partial_1\Omega}(\mb f+(\sin\theta_b)^2\partial_r\phi\mb e_\theta)\cdot\mb e_1\d r\\
&=\int_\Omega\brac{\grad\cdot \mathbb T[\mb u,p]\cdot\mb e_1}\d\mb x + \int_{\partial_0\Omega}\mb u\cdot\mb e_1\d r - \int_{\partial_1\Omega}(\mb f+(\sin\theta_b)^2\partial_r\phi\mb e_\theta)\cdot\mb e_1\d r\\
&=\int_{\partial_0\Omega}(\mb u-\mb e_\theta\cdot \mathbb T[\mb u,p])\cdot\mb e_1\d r + \int_{\partial_1\Omega}(\mb e_\theta\cdot \mathbb T[\mb u,p]-\mb f-(\sin\theta_b)^2\partial_r\phi\mb e_\theta)\cdot\mb e_1\d r\\
&=\int_{\partial_0\Omega}(\mb u-\mb e_\theta\cdot \mathbb T[\mb u,p])\cdot\mb e_r\d r=0.
\end{align*}
Here $\mb F=\grad\cdot\mathbb T[\mb u,p]$ by \eqref{stationary bulk equation} and \eqref{stationary divergence}. The divergence theorem is applied with the outward normals $-\mb e_\theta$ on $\partial_0\Omega$ and $\mb e_\theta$ on $\partial_1\Omega$. The integral over $\partial_1\Omega$ then vanishes by \eqref{stationary surface tension}, and $\mb e_1=\mb e_r$ on $\partial_0\Omega$. Finally, $(\mb u-\mb e_\theta\cdot\mathbb T[\mb u,p])\cdot\mb e_r=u_r-r^{-1}\partial_\theta u_r$ on $\partial_0\Omega$ by \eqref{stationary no penetration}, and this vanishes by \eqref{stationary Navier slip}.

Since we have from \eqref{data um regularity} $\|\mb u^{[1]}\|_{\Omega,2-\kappa,j,2}<\infty$, Lemma \ref{continuity at the corner} tells us that $\mb u^{[1]}(\mb 0)=\mb 0$. Hence $\mb u=\mb u^{[1]}-\smash[b]{c_{\mb u,\dr}}\mb e_1$ is continuous at $\mb 0$, with $\mb u(\mb 0)=-\smash[b]{c_{\mb u,\dr}}\mb e_1$, i.e. $\smash[b]{c_{\mb u,\dr}}=-u_1(\mb 0)$, $\mb u(\mb 0)=u_1(\mb 0)\mb e_1$ and $\mb u^{[1]}=\mb u-\mb u(\mb 0)$.

Now we do the same for the pressure. Its weights $-\kappa$ and $1-\kappa$ at the levels $M=0,1$ correspond to the Mellin lines $\Re\lambda=-1-\kappa$ and $\Re\lambda=-\kappa$ by \eqref{weight and line}, and between them $G_{\lambda,p}=(\lambda_2^2+\partial_\theta^2)G_{\lambda_1}$ has a pole at $\lambda=-1$, coming from the double pole of $G_{\lambda_1}$ at $\lambda_1=0$. Since $1+\partial_\theta^2$ annihilates $\sin\theta\sin\xi$, the limits of Lemma \ref{Green's function initial time} at $\lambda=0$ give
\begin{align*}
  \lim_{\lambda\to-1}(\lambda_1^2G_{\lambda,p}) &=0,\\
  \lim_{\lambda\to-1}\partial_\lambda(\lambda_1^2G_{\lambda,p})
  &=\lim_{\lambda\to-1}(2\lambda_2\lambda_1^2G_{\lambda_1})
  =-{2\sin\theta\sin\xi\over\sin(2\theta_b)+2\theta_b},
\end{align*}
so that the pole of $\hat p$ at $\lambda=-1$ is simple, and we get
\begin{align*}
  (p^{[1]}-p^{[0]})(r,\theta) &=\sqrt{2\pi}\Res_{\lambda=-1}\brac{r^\lambda\hat p} = c_{p,\dr_1}P_{p,\dr_1},
\end{align*}
where
\begin{align*}
  P_{p,\dr_1} &= r^{-1}\cos\theta.
\end{align*}
From \eqref{data step equations} and the fact that $\mb F_\Omega^{\smash[t]{[0]}}=\mb F_\Omega^{\smash[t]{[1]}}$ we have $\lpc(\mb u^{\smash[t]{[1]}}-\mb u^{\smash[t]{[0]}})-\grad(p^{\smash[t]{[1]}}-p^{\smash[t]{[0]}})=\mb 0$, hence
\begin{align*}
  \mb 0 &= c_{\mb u,\dr}\lpc\mb P_{\mb u,\dr}-c_{p,\dr_1}\grad P_{p,\dr_1}=-c_{p,\dr_1}\grad P_{p,\dr_1}=c_{p,\dr_1}r^{-2}(\cos\theta\mb e_r+\sin\theta\mb e_\theta).
\end{align*}
Therefore $\smash[b]{c_{p,\dr_1}}=0$, i.e. $p^{\smash[t]{[1]}}=p$.

Now we will move on to $M=2$. Let
\begin{subequations}\label{M2 data def}
\begin{align}
  \mb F_\Omega^{[2]} &= r^2\mb F\\
  F_{0,0}^{[2]} &= ru^{[1]}_r\label{M2 slip datum}\\
  F_{0,1}^{[2]} &= 0\\
  F_{1,0}^{[2]} &= rf_r\\
  F_{1,1}^{[2]} &= rf_\theta+(\sin\theta_b)^2\dr(\phi-c_{\phi,\dr_{-1}}r-c_{\phi,\dr_{-1}^2}r\ln r).\label{M2 stress datum}
\end{align}
\end{subequations}
In Mellin variables
\begin{subequations}\label{M2 Mellin data}
\begin{align}
  \hat{\mb F}_\Omega^{[2]}(\lambda,\theta) &= \hat{\mb F}(\lambda_{-2},\theta)\\
  \hat F_{0,0}^{[2]}(\lambda) &= {1\over\lambda_{-1}}(\dr u^{[1]}_r)^\wedge(\lambda_{-1},0)\label{M2 Mellin slip datum}\\
  \hat F_{0,1}^{[2]}(\lambda) &= 0\\
  \hat F_{1,0}^{[2]}(\lambda) &= \hat f_r(\lambda_{-1})\\
  \hat F_{1,1}^{[2]}(\lambda) &= \hat f_\theta(\lambda_{-1}) +(\sin\theta_b)^2{\lambda\over\lambda_{-1}^2}(\dr_{-1}^2\phi)^\wedge(\lambda,\theta_b).\label{M2 Mellin stress datum}
\end{align}
\end{subequations}
As for $M=0,1$, \eqref{Mellin estimate initial time} at $a=3-\kappa$ and \eqref{M2 data def} give
\begin{align*}
  &\|\mb u^{[2]}\|_{\Omega,3-\kappa,j,2}+\Gamma_{1-\kappa}^{-1}\|p^{[2]}\|_{\Omega,2-\kappa,j,1}\\
  &\lesssim{\Gamma_{2-\kappa}^2\over\theta_b^{1\over 2}\Gamma^H_{2-\kappa,\theta_b}}\brac{\|\mb F_\Omega^{[2]}\|_{\Omega,3-\kappa,j}+\|F_{0,0}^{[2]}\|_{\partial_0\Omega,{5\over 2}-\kappa,j+{1\over 2}}+\|(F_{1,0}^{[2]},F_{1,1}^{[2]})\|_{\partial_1\Omega,{5\over 2}-\kappa,j+{1\over 2}}}\\
  &\lesssim{\Gamma_{2-\kappa}^2\over\theta_b^{1\over 2}\Gamma^H_{2-\kappa,\theta_b}}\bigg(\|\mb F\|_{\Omega,1-\kappa,j}+\|u^{[1]}_r\|_{\partial_0\Omega,{3\over 2}-\kappa,j+{1\over 2}}+\|\mb f\|_{\partial_1\Omega,{3\over 2}-\kappa,j+{1\over 2}}+(\sin\theta_b)^2\|\dr\phi^{[2]}\|_{\partial_1\Omega,{5\over 2}-\kappa,j+{1\over 2}}\bigg).
\end{align*}

Write $\mb u_{F_\bullet}^{\smash[t]{[M]}}$ and $p_{F_\bullet}^{\smash[t]{[M]}}$ for the part of $\mb u^{\smash[t]{[M]}}$ and $p^{\smash[t]{[M]}}$ that arises from $F_\bullet$. By the limit of $G_\lambda$ at $\lambda=1$ in Lemma \ref{Green's function initial time} and the symmetry of $G_\lambda$,
\begin{align*}
\partial_\xi G_1(\theta,0)
&={\sin(2(\theta_b-\theta))\over 4\sin(2\theta_b)}-{1\over 4},\\
\partial_\theta\partial_\xi G_1(\theta,0)
&=-{\cos(2(\theta_b-\theta))\over 2\sin(2\theta_b)},
\end{align*}
and the only pole of $\hat{\mb u}_{F_{0,0}}$ between the lines $\Re\lambda=1-\kappa$ and $\Re\lambda=2-\kappa$ is that of $\hat F^{\smash[t]{[2]}}_{0,0}$ at $\lambda=1$, which is at most simple by \eqref{M2 Mellin slip datum}, so we get
\begin{align*}
  (\mb u_{F_{0,0}}^{[2]}-\mb u_{F_{0,0}}^{[1]})(r,\theta) &=\sqrt{2\pi}\Res_{\lambda=1}\brac{r^\lambda\hat{\mb u}_{F_{0,0}}}\\
&=\sqrt{2\pi}\lim_{\lambda\to 1}\brac{r^\lambda(\lambda_1\partial_\xi G_\lambda(\theta,0)\mb e_\theta-\partial_\theta\partial_\xi G_\lambda(\theta,0)\mb e_r)(\dr u^{[1]}_r)^\wedge(\lambda_{-1},0)}\\
&=\sqrt{2\pi}r(2\partial_\xi G_1(\theta,0)\mb e_\theta-\partial_\theta\partial_\xi G_1(\theta,0)\mb e_r)(\dr u^{[1]}_r)^\wedge(0,0)\\
  &= c_{\mb u_{F_{0,0}},\dr_{-1}}\mb P_{\mb u_{F_{0,0}},\dr_{-1}},
\end{align*}
where
\begin{align*}
  \mb P_{\mb u_{F_{0,0}},\dr_{-1}} &= r\brac{\brac{{\sin(2(\theta_b-\theta))\over\sin(2\theta_b)}-1}\mb e_\theta+\brac{\cos(2(\theta_b-\theta))\over\sin(2\theta_b)}\mb e_r},
\end{align*}
since the values above give $2\partial_\xi G_1(\theta,0)\mb e_\theta-\partial_\theta\partial_\xi G_1(\theta,0)\mb e_r=\smash[b]{\mb P_{\mb u_{F_{0,0}},\dr_{-1}}}/(2r)$, the factor $2$ being $\lambda_1$ at $\lambda=1$.

On $\partial_0\Omega$ we have $u^{[M]}_\theta=F^{[M]}_{0,1}=0$ for every $M$, so \eqref{data step slip condition} reduces there to $\partial_\theta u^{[M]}_r=F^{[M]}_{0,0}$; subtracting the levels $M=2$ and $M=1$, by \eqref{M2 slip datum} and \eqref{M01 slip datum}, gives $\partial_\theta(u_{F_{0,0},r}^{\trimtop{1.2pt}{[2]}}-u_{F_{0,0},r}^{\trimtop{1.2pt}{[1]}})=r(u^{\trimtop{1.2pt}{[1]}}_r-u^{\trimtop{1.2pt}{[0]}}_r)$ on $\partial_0\Omega$. Evaluating this at $\theta=0$, where $\mb e_1=\mb e_r$, with \eqref{ladder P F00} and $\mb u^{[1]}-\mb u^{[0]}=\smash[b]{c_{\mb u,\dr}}\mb e_1$, we get
\begin{align*}
2rc_{\mb u_{F_{0,0}},\dr_{-1}}=c_{\mb u_{F_{0,0}},\dr_{-1}}\eva{\partial_\theta(\mb e_r\cdot\mb P_{\mb u_{F_{0,0}},\dr_{-1}})}_{\theta=0}=\eva{r(\mb u^{[1]}-\mb u^{[0]})\cdot\mb e_r}_{\theta=0}=rc_{\mb u,\dr},
\end{align*}
hence
\begin{align}
2c_{\mb u_{F_{0,0}},\dr_{-1}}
&=c_{\mb u,\dr}.\label{F00 velocity coefficient}
\end{align}

By $G_{\lambda,p}=(\lambda_2^2+\partial_\theta^2)G_{\lambda_1}$, the limits of $G_\lambda$ and $\partial_\lambda G_\lambda$ at $\lambda=1$ in Lemma \ref{Green's function initial time} and the symmetry of $G_\lambda$, the kernel $\partial_\xi\partial_\theta G_{\lambda,p}(\theta,0)$ equals $\partial_\theta(4+\partial_\theta^2)\partial_\xi G_1(\theta,0)=0$ at $\lambda=0$, and
\begin{align*}
\lim_{\lambda\to 0}\brac{{1\over\lambda}\partial_\xi\partial_\theta G_{\lambda,p}(\theta,0)} &= \partial_\theta\brac{4\partial_\xi G_1+(4+\partial_\theta^2)\partial_\xi\eva{\partial_\lambda G_\lambda}_{\lambda=1}}(\theta,0)=\partial_\theta{\theta\over\sin(2\theta_b)}={1\over\sin(2\theta_b)}.
\end{align*}
So we get, for a constant $\smash[b]{c_{p_{F_{0,0}},\dr}}$,
\begin{align*}
  (p_{F_{0,0}}^{[2]}-p_{F_{0,0}}^{[1]})(r,\theta) &=\sqrt{2\pi}\Res_{\lambda=0}\brac{r^\lambda\hat p_{F_{0,0}}}
  =-\sqrt{2\pi}\lim_{\lambda\to 0}\brac{r^\lambda\brac{{1\over\lambda}\partial_\xi\partial_\theta G_{\lambda,p}(\theta,0)}(\dr u^{[1]}_r)^\wedge(\lambda,0)}\\
  &= c_{p_{F_{0,0}},\dr}P_{p_{F_{0,0}},\dr},
\end{align*}
where
\begin{align*}
  P_{p_{F_{0,0}},\dr} &= -{1\over\sin(2\theta_b)}.
\end{align*}
Since the $F_{0,0}$-parts have the datum $0$ in \eqref{data step stress condition}, we have $2(\partial_\theta(u^{\smash[t]{[2]}}_{F_{0,0},\theta}-u^{\smash[t]{[1]}}_{F_{0,0},\theta})+u^{\smash[t]{[2]}}_{F_{0,0},r}-u^{\smash[t]{[1]}}_{F_{0,0},r})-r(p^{\smash[t]{[2]}}_{F_{0,0}}-p^{\smash[t]{[1]}}_{F_{0,0}})=0$ on $\partial_1\Omega$, and we get
\begin{align*}
  2c_{\mb u_{F_{0,0}},\dr_{-1}} r \brac{-{2\over\sin(2\theta_b)} + {1\over\sin(2\theta_b)}} + c_{p_{F_{0,0}},\dr} r {1\over\sin(2\theta_b)} &= 0,
\end{align*}
which means
\begin{align}
  c_{p_{F_{0,0}},\dr} &= 2c_{\mb u_{F_{0,0}},\dr_{-1}}.\label{F00 pressure coefficient}
\end{align}

Since $\hat f_r(\lambda_{-1})$ is holomorphic on $\Re\lambda\in[1-\kappa,2-\kappa]$, we have $\mb u_{F_{1,0}}^{\smash[t]{[2]}}=\mb u_{F_{1,0}}^{\smash[t]{[1]}}$ and $p_{F_{1,0}}^{\smash[t]{[2]}}=p_{F_{1,0}}^{\smash[t]{[1]}}$.

The datum $F_{1,1}$ enters $\hat{\mb u}$ through $\lambda_{-1}\hat F_{1,1}(\lambda)$ and $\hat p$ through $\hat F_{1,1}(\lambda_1)$, by Proposition \ref{Mellin system initial time}. By \eqref{M2 Mellin stress datum}, at the level $M=2$ the first has at most a simple pole at $\lambda=1$ and the second at most a double pole at $\lambda=0$, so \eqref{residue rule} multiplies them by $\lambda-1$ and by $\lambda^2$, respectively. Their $f_\theta$-parts then become $\lambda_{-1}^2\hat f_\theta(\lambda_{-1})$ and $\lambda^2\hat f_\theta(\lambda)$, which vanish to second order at $\lambda=1$ and at $\lambda=0$, respectively, since $\hat f_\theta$, like $\hat f_r$, is holomorphic on $-\kappa\leq\Re\lambda\leq 1-\kappa$. Only the profile parts $(\sin\theta_b)^2\lambda(\dr_{-1}^2\phi)^\wedge(\lambda,\theta_b)$ and $(\sin\theta_b)^2\lambda_1(\dr_{-1}^2\phi)^\wedge(\lambda_1,\theta_b)$ therefore contribute.

By the limit of $G_\lambda$ at $\lambda=1$ in Lemma \ref{Green's function initial time},
\begin{align*}
G_1(\theta,\theta_b)
&=-{\theta\over 4},\\
\partial_\theta G_1(\theta,\theta_b)
&=-{1\over 4},
\end{align*}
and we get
\begin{align*}
  (\mb u_{F_{1,1}}^{[2]}-\mb u_{F_{1,1}}^{[1]})(r,\theta) &=\sqrt{2\pi}\Res_{\lambda=1}\brac{r^\lambda\hat{\mb u}_{F_{1,1}}}\\
  &=\sqrt{2\pi}\lim_{\lambda\to 1}\bigg(r^{\lambda}\big(\partial_\theta G_\lambda(\theta,\theta_b)\mb e_r-\lambda_1 G_\lambda(\theta,\theta_b)\mb e_\theta\big)(\sin\theta_b)^2\lambda(\dr_{-1}^2\phi)^\wedge(\lambda,\theta_b)\bigg)\\
  &= c_{\mb u_{F_{1,1}},\dr_{-1}}\mb P_{\mb u_{F_{1,1}},\dr_{-1}},
\end{align*}
where
\begin{align*}
  \mb P_{\mb u_{F_{1,1}},\dr_{-1}} &= r\brac{2\theta\mb e_\theta-\mb e_r}.
\end{align*}
By $G_{\lambda,p}=(\lambda_2^2+\partial_\theta^2)G_{\lambda_1}$ and the limits of $G_\lambda$ and $\partial_\lambda G_\lambda$ at $\lambda=1$ in Lemma \ref{Green's function initial time},
\begin{align*}
\lim_{\lambda\to 0}\brac{\partial_\theta G_{\lambda,p}(\theta,\theta_b)} &= \partial_\theta\brac{(4+\partial_\theta^2)G_1}(\theta,\theta_b)=-1,\\
\lim_{\lambda\to 0}\partial_\lambda\brac{\partial_\theta G_{\lambda,p}(\theta,\theta_b)} &= \partial_\theta\brac{4G_1+(4+\partial_\theta^2)\eva{\partial_\lambda G_\lambda}_{\lambda=1}}(\theta,\theta_b)=\partial_\theta{\theta\over 2}={1\over 2},
\end{align*}
and we get
\begin{align*}
  (p_{F_{1,1}}^{[2]}-p_{F_{1,1}}^{[1]})(r,\theta) &=\sqrt{2\pi}\Res_{\lambda=0}\brac{r^\lambda\hat p_{F_{1,1}}}\\
  &=\sqrt{2\pi}\lim_{\lambda\to 0}\partial_\lambda\bigg(e^{\lambda\ln r}\brac{\partial_\theta G_{\lambda,p}(\theta,\theta_b)}
  (\sin\theta_b)^2\lambda_1(\dr_{-1}^2\phi)^\wedge(\lambda_1,\theta_b)\bigg)\\
  &= c_{p_{F_{1,1}},\dr}P_{p_{F_{1,1}},\dr} +c_{p_{F_{1,1}},\dr^2}P_{p_{F_{1,1}},\dr^2},
\end{align*}
where
\begin{align*}
  P_{p_{F_{1,1}},\dr} &= -1,\\
  P_{p_{F_{1,1}},\dr^2} &= 1-2\ln r.
\end{align*}

Since the $F_{1,1}$-parts have no bulk datum, \eqref{data step bulk equation} gives $\lpc(\mb u^{\smash[t]{[2]}}_{F_{1,1}}-\mb u^{\smash[t]{[1]}}_{F_{1,1}})-\grad(p^{\smash[t]{[2]}}_{F_{1,1}}-p^{\smash[t]{[1]}}_{F_{1,1}})=\mb 0$ on $\Omega$, and we get
\begin{align*}
  \mb 0 &= c_{\mb u_{F_{1,1}},\dr_{-1}}\lpc\mb P_{\mb u_{F_{1,1}},\dr_{-1}} - \brac{c_{p_{F_{1,1}},\dr}\grad P_{p_{F_{1,1}},\dr} + c_{p_{F_{1,1}},\dr^2}\grad P_{p_{F_{1,1}},\dr^2}}\\
  &= c_{\mb u_{F_{1,1}},\dr_{-1}}\brac{-{4\over r}\mb e_r} - c_{p_{F_{1,1}},\dr^2}\brac{-{2\over r}\mb e_r}
  = {2\over r}\brac{-2c_{\mb u_{F_{1,1}},\dr_{-1}} + c_{p_{F_{1,1}},\dr^2}}\mb e_r,
\end{align*}
since $\smash[b]{\mb P_{\mb u_{F_{1,1}},\dr_{-1}}}=2\theta\mb x^\perp-\mb x$ with $\mb x=r\mb e_r$, and $\lpc(\theta\mb x^\perp)=2\grad^\perp\theta=-2r^{-1}\mb e_r$ as $\theta$ is harmonic and $\mb x^\perp$ is linear. This means
\begin{align}
  c_{p_{F_{1,1}},\dr^2} &= 2c_{\mb u_{F_{1,1}},\dr_{-1}}.\label{F11 bulk identity}
\end{align}

Since, by \eqref{data step stress condition}, \eqref{M2 stress datum} and \eqref{M01 stress datum},
\begin{align*}
2(\partial_\theta(u^{[2]}_{F_{1,1},\theta}-u^{[1]}_{F_{1,1},\theta})+u^{[2]}_{F_{1,1},r}-u^{[1]}_{F_{1,1},r})-r(p^{[2]}_{F_{1,1}}-p^{[1]}_{F_{1,1}})
=-(\sin\theta_b)^2\dr(c_{\phi,\dr_{-1}}r+c_{\phi,\dr_{-1}^2}r\ln r)
\end{align*}
on $\partial_1\Omega$ and
\begin{align*}
\eva{\brac{2\partial_\theta P_{\mb u_{F_{1,1}},\dr_{-1},\theta} + 2P_{\mb u_{F_{1,1}},\dr_{-1},r}}}_{\theta=\theta_b} &= 2r,
\end{align*}
we get
\begin{align*}
  &c_{\mb u_{F_{1,1}},\dr_{-1}} (2r) - c_{p_{F_{1,1}},\dr} (-r) - c_{p_{F_{1,1}},\dr^2} \brac{-2r\ln r + r}\\
  &= -(\sin\theta_b)^2(r c_{\phi,\dr_{-1}} + r c_{\phi,\dr_{-1}^2} + r\ln r c_{\phi,\dr_{-1}^2}),
\end{align*}
which means, matching the coefficients of $r\ln r$ and of $r$,
\begin{align}
  2c_{p_{F_{1,1}},\dr^2} &= -(\sin\theta_b)^2 c_{\phi,\dr_{-1}^2},\label{F11 logarithm identity}\\
  2c_{\mb u_{F_{1,1}},\dr_{-1}} + c_{p_{F_{1,1}},\dr} - c_{p_{F_{1,1}},\dr^2} &= -(\sin\theta_b)^2(c_{\phi,\dr_{-1}} + c_{\phi,\dr_{-1}^2}),\nonumber
\end{align}
and in the second line the terms with $\smash[b]{c_{\mb u_{F_{1,1}},\dr_{-1}}}$ and $\smash[b]{c_{p_{F_{1,1}},\dr^2}}$ cancel by \eqref{F11 bulk identity}, leaving
\begin{align}
  c_{p_{F_{1,1}},\dr} &= -(\sin\theta_b)^2(c_{\phi,\dr_{-1}} + c_{\phi,\dr_{-1}^2}).\label{F11 pressure coefficient}
\end{align}

Since $\hat{\mb F}(\lambda_{-2},\theta)$ is holomorphic on $\Re\lambda\in[1-\kappa,2-\kappa]$, we have $\mb u_{\mb F_\Omega}^{\smash[t]{[2]}}=\mb u_{\mb F_\Omega}^{\smash[t]{[1]}}$ and $p_{\mb F_\Omega}^{\smash[t]{[2]}}=p_{\mb F_\Omega}^{\smash[t]{[1]}}$.

In summary, \eqref{F00 velocity coefficient} and $\smash[b]{c_{\mb u,\dr}}=-u_1(\mb 0)$ give $\smash[b]{c_{\mb u_{F_{0,0}},\dr_{-1}}}=-u_1(\mb 0)/2$, and \eqref{F11 bulk identity} and \eqref{F11 logarithm identity} give $\smash[b]{c_{\mb u_{F_{1,1}},\dr_{-1}}}=-(\sin\theta_b)^2\smash[b]{c_{\phi,\dr_{-1}^2}}/4$. Adding the contributions of $F_{0,0}$ and $F_{1,1}$, with the pressure coefficients \eqref{F00 pressure coefficient} and \eqref{F11 pressure coefficient}, we have
\begin{align*}
\mb u^{[2]} &=\mb u^{[1]}- {u_1(\mb 0)\over 2}\mb P_{\mb u_{F_{0,0}},\dr_{-1}}-{1\over 4}(\sin\theta_b)^2c_{\phi,\dr_{-1}^2}\mb P_{\mb u_{F_{1,1}},\dr_{-1}}\\
&=\mb u - \mb u(\mb 0) - {u_1(\mb 0)\over 2}\mb P_{\mb u_{F_{0,0}},\dr_{-1}}-{1\over 4}(\sin\theta_b)^2c_{\phi,\dr_{-1}^2}\mb P_{\mb u_{F_{1,1}},\dr_{-1}}\\
p^{[2]} &= p^{[1]} +{u_1(\mb 0)\over\sin(2\theta_b)} +(\sin\theta_b)^2(c_{\phi,\dr_{-1}}+c_{\phi,\dr_{-1}^2})+{1\over 2}(\sin\theta_b)^2c_{\phi,\dr_{-1}^2}(2\ln r-1)\\
&= p +{u_1(\mb 0)\over\sin(2\theta_b)} +(\sin\theta_b)^2c_{\phi,\dr_{-1}}+{1\over 2}(\sin\theta_b)^2c_{\phi,\dr_{-1}^2}(2\ln r+1).
\end{align*}

Since we have from \eqref{data um regularity} $\|\mb u^{\smash[t]{[2]}}\|_{\Omega,3-\kappa,j,2}<\infty$, Lemma \ref{continuity at the corner} tells us that $\grad\mb u^{\smash[t]{[2]}}(\mb 0)=\mb 0$, so that $\partial_ru_\theta^{\smash[t]{[2]}}(0,\theta_b)=0$ by \eqref{radial shear cartesian}. By \eqref{ladder P F00} and \eqref{ladder P F11}, $\partial_r(\mb e_\theta\cdot\smash[b]{\mb P_{\mb u_{F_{0,0}},\dr_{-1}}})=-1$ and $\partial_r(\mb e_\theta\cdot\smash[b]{\mb P_{\mb u_{F_{1,1}},\dr_{-1}}})=2\theta_b$ at $\theta=\theta_b$, so that \eqref{ladder u2} gives
\begin{align*}
\partial_ru_\theta(0,\theta_b)={\theta_b\over 2}(\sin\theta_b)^2c_{\phi,\dr_{-1}^2}-{u_1(\mb 0)\over 2}.
\end{align*}
This means that having $\partial_ru_\theta(0,\theta_b)=0$ is equivalent to \eqref{C1 in the linear theory}.
\end{proof}

\subsubsection{The initial jet}\label{initial jet subsection}

The initial jet is the profile together with its first two time derivatives at $t=0$. As explained at the beginning of Section \ref{corner expansions}, by the profile evolution equation of \eqref{system: time side} the first time derivative, the derived profile, is $(\sin\theta_b)^{-2}\partial_ru_\theta$ on $\partial_1\Omega$, where $(\mb u,p)$ solves the stationary system with the data $(\phi,\mb F,\mb f)$. The second is obtained in the same way from the dotted pair $(\dot{\mb u},\dot p)$, which solves the stationary system with the time derivatives of these data. The next theorem bounds the two pairs and the derived profiles by the data, by Propositions \ref{stationary system} and \ref{corner ladder} applied to both pairs.

\begin{theorem}[The initial jet]\label{initial jet}
Let $\kappa\in(0,1/2]$ and $l\in\N_0$, and let $\theta_b\in(0,\pi/2)$ satisfy $2-\kappa<\Lambda^*(\theta_b)$ and be small enough for Proposition \ref{stationary system} to apply. Let $(\mb u,p)$ be the solution with respect to data $(\phi,\mb F,\mb f)$ according to Propositions \ref{stationary system} and \ref{corner ladder} that satisfies the compatibility condition \eqref{C1 in the linear theory}. Let $\dot\phi:=(\sin\theta_b)^{-2}\partial_ru_\theta(\ph,\theta_b)$, let $(\dot{\mb F},\dot{\mb f})$ be given, and let $(\dot{\mb u},\dot p)$ be the solution with respect to data $(\dot\phi,\dot{\mb F},\dot{\mb f})$ according to Proposition \ref{stationary system} and to Proposition \ref{corner ladder} at the rungs $M\leq 1$ only. Let $\ddot\phi:=(\sin\theta_b)^{-2}\partial_r\dot u_\theta(\ph,\theta_b)$. Then $(\mb u,p,\dot{\mb u},\dot p,\dot\phi,\ddot\phi)$ has the following regularity provided the data $(\phi,\mb F,\mb f,\dot{\mb F},\dot{\mb f})$ have regularity that makes the RHS below finite:
\begin{align*}
&\sum_{M=0}^2\brac{(\sin\theta_b)^2\|\dot\phi\|_{\partial_1\Omega,M-{1\over 2}-\kappa,l+{5\over 2}}+\|\mb u^{[M]}\|_{\Omega,M+1-\kappa,l+2,2}+\|p^{[M]}\|_{\Omega,M-\kappa,l+2,1}}\\
&\quad+\sum_{M=0}^1\brac{(\sin\theta_b)^2\|\ddot\phi\|_{\partial_1\Omega,M-{1\over 2}-\kappa,l+{3\over 2}}+\|\dot{\mb u}^{[M]}\|_{\Omega,M+1-\kappa,l+1,2}+\|\dot p^{[M]}\|_{\Omega,M-\kappa,l+1,1}}\\
&\lesssim_{\theta_b,\kappa}\sum_{M=0}^2\brac{\|\mb F\|_{\Omega,M-1-\kappa,l+2}+\|\mb f\|_{\partial_1\Omega,M-{1\over 2}-\kappa,l+{5\over 2}}+(\sin\theta_b)^2\|\phi^{[M]}\|_{\partial_1\Omega,M+{1\over 2}-\kappa,l+{7\over 2}}}\\
&\quad+\sum_{M=0}^1\brac{\|\dot{\mb F}\|_{\Omega,M-1-\kappa,l+1}+\|\dot{\mb f}\|_{\partial_1\Omega,M-{1\over 2}-\kappa,l+{3\over 2}}}.
\end{align*}
\end{theorem}
\begin{proof}
From Propositions \ref{stationary system} and \ref{corner ladder} we have, absorbing the factors $\Gamma_{M-\kappa}^2\theta_b^{-1/2}(\Gamma^H_{M-\kappa,\theta_b})^{-1}$ and $\Gamma_{M-1-\kappa}$ of \eqref{um data reg master} into the constant,
\begin{align*}
\sum_{M=0}^2&\brac{\|\mb u^{[M]}\|_{\Omega,M+1-\kappa,l+2,2}+\|p^{[M]}\|_{\Omega,M-\kappa,l+2,1}}\nonumber\\
&\lesssim_{\theta_b,\kappa}\sum_{M=0}^2\brac{\|\mb F\|_{\Omega,M-1-\kappa,l+2}+\|\mb f\|_{\partial_1\Omega,M-{1\over 2}-\kappa,l+{5\over 2}}+(\sin\theta_b)^2\|\phi^{[M]}\|_{\partial_1\Omega,M+{1\over 2}-\kappa,l+{7\over 2}}}.
\end{align*}
On $\partial_1\Omega$ we have $\partial_ru_\theta=\partial_ru^{\smash[t]{[M]}}_\theta$ for $M\leq 2$. For $M\leq 1$ this follows from \eqref{data step relations}. For $M=2$, \eqref{ladder u2}, \eqref{ladder P F00} and \eqref{ladder P F11} give $\partial_r(u_\theta-u^{\smash[t]{[2]}}_\theta)=\theta_b(\sin\theta_b)^2\smash[b]{c_{\phi,\dr_{-1}^2}}/2-u_1(\mb 0)/2$ on $\partial_1\Omega$, which vanishes by the compatibility condition \eqref{C1 in the linear theory}. Hence, by Lemma \ref{trace},
\begin{align*}
\sum_{M=0}^2(\sin\theta_b)^2\|\dot\phi\|_{\partial_1\Omega,M-{1\over 2}-\kappa,l+{5\over 2}} &=\sum_{M=0}^2\|\partial_ru^{[M]}_\theta(\ph,\theta_b)\|_{\partial_1\Omega,M-{1\over 2}-\kappa,l+{5\over 2}}\\
&\lesssim\sum_{M=0}^2\|\dr u^{[M]}_\theta(\ph,\theta_b)\|_{\partial_1\Omega,M+{1\over 2}-\kappa,l+{5\over 2}}
\lesssim_{\theta_b}\sum_{M=0}^2\|\mb u^{[M]}\|_{\Omega,M+1-\kappa,l+2,2}.
\end{align*}
From Propositions \ref{stationary system} and \ref{corner ladder} we have, for the dotted pair,
\begin{align*}
\sum_{M=0}^1&\brac{\|\dot{\mb u}^{[M]}\|_{\Omega,M+1-\kappa,l+1,2}+\|\dot p^{[M]}\|_{\Omega,M-\kappa,l+1,1}}\nonumber\\
&\lesssim_{\theta_b,\kappa}\sum_{M=0}^1\brac{\|\dot{\mb F}\|_{\Omega,M-1-\kappa,l+1}+\|\dot{\mb f}\|_{\partial_1\Omega,M-{1\over 2}-\kappa,l+{3\over 2}}+(\sin\theta_b)^2\|\dot\phi^{[M]}\|_{\partial_1\Omega,M+{1\over 2}-\kappa,l+{5\over 2}}}.
\end{align*}
In the same way as for $\dot\phi$, now for $M\leq 1$ only,
\begin{align*}
\sum_{M=0}^1(\sin\theta_b)^2\|\ddot\phi\|_{\partial_1\Omega,M-{1\over 2}-\kappa,l+{3\over 2}} &=\sum_{M=0}^1\|\partial_r\dot u^{[M]}_\theta(\ph,\theta_b)\|_{\partial_1\Omega,M-{1\over 2}-\kappa,l+{3\over 2}}\\
&\lesssim\sum_{M=0}^1\|\dr\dot u^{[M]}_\theta(\ph,\theta_b)\|_{\partial_1\Omega,M+{1\over 2}-\kappa,l+{3\over 2}}
\lesssim_{\theta_b}\sum_{M=0}^1\|\dot{\mb u}^{[M]}\|_{\Omega,M+1-\kappa,l+1,2}.
\end{align*}
Combining these four estimates, we get the estimate of the theorem.
\end{proof}

In Theorem \ref{initial jet} no compatibility condition is asked of the dotted pair, whose decomposition at $M=2$ is not used there.

\subsubsection{The third rung}\label{third rung subsection}

The scheme of Proposition \ref{corner ladder} can be iterated once more, and for the nonlinear theory of Section \ref{nonlinear existence} we shall need one further rung: it is what gives the \emph{dotted} pair of Theorem \ref{initial jet} its $M=2$ level, on which the second-order contact angle condition rests. Two new phenomena appear at this level, and both are already visible in the data that the nonlinear problem produces. First, the poles that the Mellin line crosses are no longer at worst double: second-order logarithms enter, because the quadratic nonlinearity multiplies the corner behaviours $r$ and $r\ln r$ of the profile. For instance, the leading part of the curvature term of $\mb N_\partial$ in System \ref{reference wedge system} is $3\cos\theta_b(\sin\theta_b)^3\phi(\partial_r\phi)\mb e_\theta$, and in $\phi\partial_r\phi$ the term $\smash[b]{c_{\phi,\dr_{-1}^2}}r\ln r$ of $\phi$ produces $\smash[b]{c_{\phi,\dr_{-1}^2}^2}r\ln^2r$. We therefore extend the notation \eqref{profile corner coefficients} to poles of higher order, and the next lemma shows how the coefficients are read off from the Laurent expansion of the Mellin transform. For $d=1$ and $m=2$ the extension reduces to \eqref{profile corner coefficients}, and in its notation $\dr_{-d}^k$ is the lowest power of $\dr_{-d}$ that annihilates $r^d(\ln r)^{k-1}$.

\begin{lemma}[Corner blocks and Laurent coefficients]\label{blocks and Laurent coefficients}
Let $g$ be a function on $(0,\infty)$ whose Mellin transform on a line to the left of $\Re\lambda=d$ continues to a function $\hat g$ meromorphic across $\Re\lambda=d$, with a pole of order at most $m$ at $\lambda=d$, and define the block of $g$ at $(d,m)$, with its coefficients $\smash[b]{c_{g,\dr_{-d}}},\ldots,\smash[b]{c_{g,\dr_{-d}^m}}$, by
\begin{align}\label{higher corner coefficients}
\sum_{k=1}^{m}c_{g,\dr_{-d}^k}{r^d(\ln r)^{k-1}\over(k-1)!}:=-\sqrt{2\pi}\Res_{\lambda=d}\brac{r^\lambda\hat g(\lambda)}.
\end{align}
Then, as $\lambda\to d$,
\begin{align}\label{principal part}
\hat g(\lambda)=-{1\over\sqrt{2\pi}}\sum_{k=1}^m{c_{g,\dr_{-d}^k}\over(\lambda-d)^k}+O(1).
\end{align}
\end{lemma}
\begin{proof}
Write $\hat g(\lambda)$ as $\sum_{k=1}^ma_k(\lambda-d)^{-k}$ plus a function holomorphic at $\lambda=d$. The latter remains holomorphic there after multiplication by $r^\lambda$, so it does not contribute to the residue. Since $r^\lambda=r^de^{(\lambda-d)\ln r}$, the Taylor series of the exponential gives
\begin{align*}
r^\lambda\sum_{k=1}^m{a_k\over(\lambda-d)^k}=r^d\sum_{n\geq 0}{(\ln r)^n\over n!}(\lambda-d)^n\sum_{k=1}^m{a_k\over(\lambda-d)^k}=r^d\sum_{n\geq 0}\sum_{k=1}^m{a_k(\ln r)^n\over n!}(\lambda-d)^{n-k}.
\end{align*}
The residue at $\lambda=d$ is the coefficient of $(\lambda-d)^{-1}$, which collects exactly the terms with $n=k-1$, so that
\begin{align*}
\Res_{\lambda=d}\brac{r^\lambda\hat g(\lambda)}=r^d\sum_{k=1}^ma_k{(\ln r)^{k-1}\over(k-1)!}.
\end{align*}
Comparing with \eqref{higher corner coefficients}, the functions $r^d(\ln r)^{k-1}$ being linearly independent, gives $a_k=-\smash[b]{c_{g,\dr_{-d}^k}}/\sqrt{2\pi}$, which is \eqref{principal part}.
\end{proof}

On lines to the right of $\Re\lambda=d$, $\hat g$ is not the Mellin transform of $g$: if $g$ minus its block is square integrable against $r^{-2a_1-1}\d r$ for some $a_1>d$, then $\hat g$ is the Mellin transform of $g$ minus its block on the line $\Re\lambda=a_1$, by Lemma \ref{block matching}, so that the inverse Mellin transforms of $\hat g$ on the two sides of $\Re\lambda=d$ are $g$ and $g$ minus its block.

Second, the residues are no longer spanned by a fixed pair of explicit fields: the bulk forcing now contributes as well, with Laurent coefficients that are functions of $\theta$, so the subtracted corner blocks are best described as elements of the classes
\begin{align}\label{corner classes}
\mathcal Q_d^m:=\set{r^d\sum_{k=0}^{m-1}(\ln r)^k\mb q_k(\theta):\mb q_k\in H^2([0,\theta_b])},\qquad
\norm{r^d\sum_{k}(\ln r)^k\mb q_k}_{\mathcal Q_d^m}:=\sum_{k=0}^{m-1}\|\mb q_k\|_{H^2([0,\theta_b])},
\end{align}
which are finite dimensional in $r$, and similarly for scalar functions.

\begin{proposition}[The third rung]\label{third rung}
Let $\theta_b\in(0,\pi/2)$ and $\kappa\in(0,1)$ with $3-\kappa<\Lambda^*(\theta_b)$, and let $j\geq 1$.
\begin{enumerate}
\item (Second rung with a triple pole) Suppose that in Proposition \ref{corner ladder} the Mellin transform of $\phi$ on the line $\Re\lambda=1-\kappa$ continues to a function $\hat\phi(\ph,\theta_b)$ meromorphic across $\Re\lambda=1$, with a pole of order at most three at $\lambda=1$. Let $\phi^{[2]}$, $\mb u^{[2]}$ and $p^{[2]}$ be given, in place of \eqref{ladder phi2}, \eqref{ladder u2} and \eqref{ladder p2}, by
\begin{align}\label{triple subtracted profile}
\phi^{[2]}:=\phi-c_{\phi,\dr_{-1}}r-c_{\phi,\dr_{-1}^2}r\ln r-{1\over 2}c_{\phi,\dr_{-1}^3}r\ln^2r,
\end{align}
with the coefficients \eqref{higher corner coefficients} of $\phi$ at $(d,m)=(1,3)$, and
\begin{align}
\mb u^{[2]} &=\mb u-\mb u(\mb 0)-{u_1(\mb 0)\over 2}\mb P_{\mb u_{F_{0,0}},\dr_{-1}}-{(\sin\theta_b)^2c_{\phi,\dr_{-1}^2}\over 4}\mb P_{\mb u_{F_{1,1}},\dr_{-1}}-{(\sin\theta_b)^2c_{\phi,\dr_{-1}^3}\over 4}\mb P_{\mb u_{F_{1,1}},\dr_{-1}^2},\nonumber\\
p^{[2]} &=p+{u_1(\mb 0)\over\sin(2\theta_b)}+(\sin\theta_b)^2c_{\phi,\dr_{-1}}+{1\over 2}(\sin\theta_b)^2c_{\phi,\dr_{-1}^2}(2\ln r+1)+{1\over 2}(\sin\theta_b)^2c_{\phi,\dr_{-1}^3}Q^p_\star,\label{log decomposition}
\end{align}
where
\begin{align}
\mb P_{\mb u_{F_{1,1}},\dr_{-1}^2}&:=\mb P_{\mb u_{F_{1,1}},\dr_{-1}}\ln r+r\mb Q_\star(\theta),\label{log field}\\
\mb Q_\star(\theta)&:=\brac{{1\over 2}+{\theta_b\cos(2\theta)\over\sin(2\theta_b)}}\mb e_r-{\theta_b\sin(2\theta)\over\sin(2\theta_b)}\mb e_\theta,\nonumber\\
Q^p_\star&:=\ln^2r-\theta^2+\ln r+\theta_b^2+\theta_b\cot(2\theta_b)-{1\over 2}.\label{log pressure}
\end{align}
Then \eqref{um data reg master} holds at $M=2$.
\item (Third rung) Suppose $(\mb u,p,\phi,\mb F,\mb f)$ satisfy the hypotheses of Proposition \ref{corner ladder}, with the RHS of \eqref{um data reg master} finite for $M\leq 2$, and suppose the data admit the corner decompositions
\begin{align}\label{third rung data classes}
\mb F=\mb F^{[3]}_{\rm c}+\mb F^{[3]},\qquad\quad
\mb f=\mb f^{[3]}_{\rm c}+\mb f^{[3]},\qquad\quad
\phi^{[2]}=\phi^{[3]}_{\rm c}+\phi^{[3]},
\end{align}
where $\mb F^{\smash[t]{[3]}}_{\rm c}\in\mathcal Q_0^2$, $\mb f^{\smash[t]{[3]}}_{\rm c}\in\mathcal Q_1^3$ with $\theta$-independent coefficients, and $\phi^{\smash[t]{[3]}}_{\rm c}$ is the block \eqref{higher corner coefficients} of $\phi$ at $(d,m)=(2,3)$, and
\begin{align}
\mathcal R_j(\mb F,\mb f)&:=\|\mb F^{[3]}\|_{\Omega,2-\kappa,j}+\|\mb F^{[3]}_{\rm c}\|_{\mathcal Q_0^2}+\|\mb f^{[3]}\|_{\partial_1\Omega,{5\over 2}-\kappa,j+{1\over 2}}+\|\mb f^{[3]}_{\rm c}\|_{\mathcal Q_1^3},\nonumber\\
\mathcal R_j&:=\mathcal R_j(\mb F,\mb f)+(\sin\theta_b)^2\brac{\|\dr\phi^{[3]}\|_{\partial_1\Omega,{7\over 2}-\kappa,j+{1\over 2}}+\sum_{k=1}^3|c_{\phi,\dr_{-2}^k}|}<\infty.\label{third rung data norm}
\end{align}
Then $(\mb u^{\smash[t]{[2]}},p^{\smash[t]{[2]}})$ splits in turn,
\begin{align*}
\mb u^{[2]}=\mb u^{[3]}_{\rm c}+\mb u^{[3]},\qquad\qquad p^{[2]}=p^{[3]}_{\rm c}+p^{[3]},
\end{align*}
with $\mb u^{\smash[t]{[3]}}_{\rm c}\in\mathcal Q_2^3$ and $p^{\smash[t]{[3]}}_{\rm c}\in\mathcal Q_1^3$, and
\begin{align}\label{third rung estimate}
&\|\mb u^{[3]}\|_{\Omega,4-\kappa,j,2}+\Gamma_{2-\kappa}^{-1}\|p^{[3]}\|_{\Omega,3-\kappa,j,1}+\|\mb u^{[3]}_{\rm c}\|_{\mathcal Q_2^3}+\|p^{[3]}_{\rm c}\|_{\mathcal Q_1^3}\nonumber\\
&\qquad\lesssim_{\theta_b,\kappa}\mathcal R_j+\|u_r^{[2]}\|_{\partial_0\Omega,{5\over 2}-\kappa,j+{1\over 2}}+|u_1(\mb 0)|+(\sin\theta_b)^2|c_{\phi,\dr_{-1}^2}|.
\end{align}
\item (The derived profile) In the setting of part (ii), suppose in addition that the compatibility condition \eqref{C1 in the linear theory} holds, and let $\dot\phi:=(\sin\theta_b)^{-2}\partial_ru_\theta(\ph,\theta_b)$. Then $\dot\phi^\wedge$ is meromorphic across $\Re\lambda=1$ with at worst a triple pole there, and, writing $\dot\phi^{[2]}$ for $\dot\phi$ minus its block \eqref{higher corner coefficients} at $(d,m)=(1,3)$,
\begin{align}\label{derived profile third rung}
(\sin\theta_b)^2&\brac{\sum_{k=1}^3|c_{\dot\phi,\dr_{-1}^k}|+\|\dot\phi^{[2]}\|_{\partial_1\Omega,{5\over 2}-\kappa,j+{1\over 2}}}\nonumber\\
&\qquad\lesssim_{\theta_b,\kappa}\mathcal R_j+\|u_r^{[2]}\|_{\partial_0\Omega,{5\over 2}-\kappa,j+{1\over 2}}+|u_1(\mb 0)|+(\sin\theta_b)^2|c_{\phi,\dr_{-1}^2}|.
\end{align}
\end{enumerate}
\end{proposition}
\begin{proof}
The three parts follow the scheme of the proof of Proposition \ref{corner ladder}: the solution at the higher weight is given by Proposition \ref{Mellin system initial time}, and its difference from the solution at the lower weight is read off \eqref{residue rule} at the poles between the two Mellin lines. We call the data $\mb F_\Omega$, $F_{0,0}$, $F_{0,1}$, $F_{1,0}$, $F_{1,1}$ of Proposition \ref{Mellin system initial time}, and those of Proposition \ref{Mellin system positive time} below, channels. The Navier slip condition enters through the slip channel $F_{0,0}$, and the profile through the surface tension channel $F_{1,1}$, which we also call the profile channel. In part (i) the poles of the profile channel have one order more, in part (ii) the poles at $\lambda=2$ come from the corner blocks of the data, and part (iii) restricts the result of part (ii) to $\partial_1\Omega$.

\emph{Part (i).} Against the proof of Proposition \ref{corner ladder} only the profile part of the channel $F_{1,1}$ changes, and only at the rung $M=2$, since the rungs $M\leq 1$ involve $\hat\phi$ only on lines to the left of $\lambda=1$. At $M=2$ the datum $F^{\smash[t]{[2]}}_{1,1}=rf_\theta+(\sin\theta_b)^2\dr\phi^{\smash[t]{[2]}}$ of \eqref{M2 data def} is formed with \eqref{triple subtracted profile}, so that in Mellin variables it carries one more annihilator,
\begin{align*}
\hat F_{1,1}^{[2]}(\lambda)=\hat f_\theta(\lambda_{-1})+(\sin\theta_b)^2{\lambda\over\lambda_{-1}^3}(\dr_{-1}^3\phi)^\wedge(\lambda,\theta_b),
\end{align*}
and Proposition \ref{Mellin system initial time} gives \eqref{um data reg master} at $M=2$ as in the proof of Proposition \ref{corner ladder}. In that proof the line of integration is moved across $\lambda=1$ using that the annihilated profile $\dr_{-1}^2\phi=\dr_{-1}^2\phi^{[2]}$ has finite norms at the weights of both lines, so that its Mellin transform is holomorphic on the strip $1-\kappa\leq\Re\lambda\leq 2-\kappa$. For $\dr_{-1}^3\phi$ such norms would need two radial derivatives of $\dr\phi$, while \eqref{um data reg master} provides $j+1/2$ of them, fewer than two when $j<3/2$. We sidestep this obstacle with Lemma \ref{block matching} of the appendix, which needs no derivatives: applied to the pair $(\dr\phi,\dr\phi^{[2]})$ with the norms $\|\dr\phi\|_{\partial_1\Omega,3/2-\kappa,j+1/2}$ and $\|\dr\phi^{[2]}\|_{\partial_1\Omega,5/2-\kappa,j+1/2}$ of \eqref{um data reg master}, it shows that $(\dr_{-1}^3\phi)^\wedge(\lambda,\theta_b)=\lambda_{-1}^3\hat\phi(\lambda,\theta_b)$ is holomorphic on that strip, and by \eqref{principal part}
\begin{align}\label{triple pole datum}
(\dr_{-1}^3\phi)^\wedge(\lambda,\theta_b)=-{1\over\sqrt{2\pi}}\brac{c_{\phi,\dr_{-1}^3}+c_{\phi,\dr_{-1}^2}\lambda_{-1}+c_{\phi,\dr_{-1}}\lambda_{-1}^2}+O(\lambda_{-1}^3).
\end{align}
It remains to recompute the residues of the profile part of the channel $F_{1,1}$, the term $\hat f_\theta(\lambda_{-1})$ being holomorphic across $\Re\lambda=1$.

For the velocity, the channel $F_{1,1}$ enters the representation of Proposition \ref{Mellin system initial time} as $(\partial_\theta G_\lambda(\theta,\theta_b)\mb e_r-\lambda_1G_\lambda(\theta,\theta_b)\mb e_\theta)\lambda_{-1}\hat F_{1,1}(\lambda)$, which now has a double pole at $\lambda=1$. Here $G_\lambda(\theta,\theta_b)=-w_{11}(\theta)/(2\lambda H(\lambda))$ by Lemma \ref{Green's function initial time}, since $w_{21}(\theta_b)=2$ and $w_{22}(\theta_b)=0$, and expanding at $\lambda=1$ we get
\begin{align*}
G_\lambda(\theta,\theta_b)&=-{\theta\over 4}+\brac{{3\theta\over 8}+{\theta_b\sin(2\theta)\over 8\sin(2\theta_b)}}\lambda_{-1}+O(\lambda_{-1}^2),\\
r^\lambda\lambda\brac{\partial_\theta G_\lambda(\theta,\theta_b)\mb e_r-\lambda_1G_\lambda(\theta,\theta_b)\mb e_\theta}&={1\over 4}\mb P_{\mb u_{F_{1,1}},\dr_{-1}}+{1\over 4}\mb P_{\mb u_{F_{1,1}},\dr_{-1}^2}\lambda_{-1}+O(\lambda_{-1}^2),
\end{align*}
with the fields of \eqref{data step relations} and \eqref{log field}. By \eqref{residue rule}, multiplying the last expansion by \eqref{triple pole datum}, we get
\begin{align*}
(\mb u_{F_{1,1}}^{[2]}-\mb u_{F_{1,1}}^{[1]})(r,\theta)
&=\sqrt{2\pi}\lim_{\lambda\to 1}\partial_\lambda\brac{\lambda_{-1}^2r^\lambda\hat{\mb u}_{F_{1,1}}}\\
&=\sqrt{2\pi}(\sin\theta_b)^2\lim_{\lambda\to 1}\partial_\lambda\bigg(r^\lambda\lambda\big(\partial_\theta G_\lambda(\theta,\theta_b)\mb e_r-\lambda_1G_\lambda(\theta,\theta_b)\mb e_\theta\big)(\dr_{-1}^3\phi)^\wedge(\lambda,\theta_b)\bigg)\\
&=-{(\sin\theta_b)^2\over 4}\brac{c_{\phi,\dr_{-1}^2}\mb P_{\mb u_{F_{1,1}},\dr_{-1}}+c_{\phi,\dr_{-1}^3}\mb P_{\mb u_{F_{1,1}},\dr_{-1}^2}}.
\end{align*}
The term with $\smash[b]{c_{\phi,\dr_{-1}^2}}$\vspace{1pt} is the one found in the proof of Proposition \ref{corner ladder}, and the term with $\smash[b]{c_{\phi,\dr_{-1}^3}}$ is the new term of $\mb u^{[2]}$ in \eqref{log decomposition}.

For the pressure, the channel $F_{1,1}$ enters the representation as $\partial_\theta G_{\lambda,p}(\theta,\theta_b)\hat F_{1,1}(\lambda_1)$, which now has a triple pole at $\lambda=0$. Here $G_{\lambda,p}=(\lambda_2^2+\partial_\theta^2)G_{\lambda_1}$, and in $G_{\lambda_1}(\theta,\theta_b)=-w_{11}(\theta)/(2\lambda_1H(\lambda_1))$, with $w_{11}$ taken at $\lambda_1$ in place of $\lambda$, the operator $\lambda_2^2+\partial_\theta^2$ annihilates the term with $\sin(\lambda_2\theta)$ and multiplies the term with $\sin(\lambda\theta)$ by $\lambda_2^2-\lambda^2=4\lambda_1$. Expanding at $\lambda=0$ we get
\begin{align*}
\partial_\theta G_{\lambda,p}(\theta,\theta_b)&=-{2\sin(\lambda_2\theta_b)\over H(\lambda_1)}\cos(\lambda\theta)=-1+{\lambda\over 2}-\brac{{1\over 4}+{\theta_b^2-\theta^2\over 2}+{\theta_b\cot(2\theta_b)\over 2}}\lambda^2+O(\lambda^3),\\
r^\lambda\lambda_1\partial_\theta G_{\lambda,p}(\theta,\theta_b)&=-1-\brac{\ln r+{1\over 2}}\lambda-{1\over 2}Q^p_\star\lambda^2+O(\lambda^3),
\end{align*}
the first two coefficients of the first line being the values used in the proof of Proposition \ref{corner ladder}. By \eqref{residue rule}, multiplying the last expansion by \eqref{triple pole datum} at $\lambda_1$ in place of $\lambda$, we get
\begin{align*}
(p_{F_{1,1}}^{[2]}-p_{F_{1,1}}^{[1]})(r,\theta)
&={\sqrt{2\pi}\over 2}\lim_{\lambda\to 0}\partial_\lambda^2\brac{\lambda^3r^\lambda\hat p_{F_{1,1}}}
={\sqrt{2\pi}\over 2}(\sin\theta_b)^2\lim_{\lambda\to 0}\partial_\lambda^2\brac{r^\lambda\lambda_1\partial_\theta G_{\lambda,p}(\theta,\theta_b)(\dr_{-1}^3\phi)^\wedge(\lambda_1,\theta_b)}\\
&=(\sin\theta_b)^2\brac{c_{\phi,\dr_{-1}}+{1\over 2}c_{\phi,\dr_{-1}^2}(2\ln r+1)+{1\over 2}c_{\phi,\dr_{-1}^3}Q^p_\star}.
\end{align*}
The terms with $\smash[b]{c_{\phi,\dr_{-1}}}$ and $\smash[b]{c_{\phi,\dr_{-1}^2}}$\vspace{1pt} are those found in the proof of Proposition \ref{corner ladder}, and the term with $\smash[b]{c_{\phi,\dr_{-1}^3}}$ is the new term of $p^{[2]}$. As a check, the two new terms solve the system of Proposition \ref{Mellin system initial time} with all data zero except $F_{1,1}=-{1\over 2}(\sin\theta_b)^2\smash[b]{c_{\phi,\dr_{-1}^3}}\dr(r\ln^2r)$, which is the part of $F^{[2]}_{1,1}-F^{[1]}_{1,1}$ carrying $\smash[b]{c_{\phi,\dr_{-1}^3}}$; for instance $\grad\cdot\smash[b]{\mb P_{\mb u_{F_{1,1}},\dr_{-1}^2}}=0$ and $\lpc\smash[b]{\mb P_{\mb u_{F_{1,1}},\dr_{-1}^2}}=-2\grad Q^p_\star$.

\emph{Part (ii), Step 1: the solution at the weight $4-\kappa$.} Let
\begin{align}\label{M3 data def}
\mb F^{[3]}_\Omega:=r^2\mb F^{[3]},\quad
F^{[3]}_{0,0}:=ru^{[2]}_r,\quad
F^{[3]}_{0,1}:=0,\quad
F^{[3]}_{1,0}:=rf^{[3]}_r,\quad
F^{[3]}_{1,1}:=rf^{[3]}_\theta+(\sin\theta_b)^2\dr\phi^{[3]}.
\end{align}
In the estimate below, the norm $\|u^{\smash[t]{[2]}}_r\|_{\partial_0\Omega,5/2-\kappa,j+1/2}$, through which the datum $F^{\smash[t]{[3]}}_{0,0}=ru^{\smash[t]{[2]}}_r$ enters, is finite by Lemma \ref{trace}, since $\|\mb u^{\smash[t]{[2]}}\|_{\Omega,3-\kappa,j,2}$ is finite; the norms of the other data there are finite since $\mathcal R_j<\infty$. Proposition \ref{Mellin system initial time} at the weight $a=4-\kappa$, which is admissible since $|a-1|=3-\kappa<\Lambda^*(\theta_b)$, therefore gives a solution $(\mb u^{[3]},p^{[3]})$ with the data \eqref{M3 data def}, and as in \eqref{data um regularity}
\begin{align*}
&\|\mb u^{[3]}\|_{\Omega,4-\kappa,j,2}+\Gamma_{2-\kappa}^{-1}\|p^{[3]}\|_{\Omega,3-\kappa,j,1}\\
&\lesssim{\Gamma_{3-\kappa}^2\over\theta_b^{1\over 2}\Gamma^H_{3-\kappa,\theta_b}}\bigg(\|\mb F^{[3]}\|_{\Omega,2-\kappa,j}+\|u^{[2]}_r\|_{\partial_0\Omega,{5\over 2}-\kappa,j+{1\over 2}}+\|\mb f^{[3]}\|_{\partial_1\Omega,{5\over 2}-\kappa,j+{1\over 2}}+(\sin\theta_b)^2\|\dr\phi^{[3]}\|_{\partial_1\Omega,{7\over 2}-\kappa,j+{1\over 2}}\bigg),
\end{align*}
which bounds the norms of $(\mb u^{\smash[t]{[3]}},p^{\smash[t]{[3]}})$ in \eqref{third rung estimate}.

\emph{Step 2: the data of the two rungs.} By \eqref{data step equations}, $(\mb u^{[2]},p^{[2]})$ is the solution of Proposition \ref{Mellin system initial time} at the weight $3-\kappa$ with the data \eqref{M2 data def}. These differ from the data \eqref{M3 data def} by corner blocks of degree two: in the channels $\mb F_\Omega$, $F_{1,0}$ and $F_{1,1}$ by those formed from the blocks of \eqref{third rung data classes}, and in the channel $F_{0,0}$ by
\begin{align}\label{slip block}
r(u^{[1]}_r-u^{[2]}_r)(\ph,0)=r^2\brac{{u_1(\mb 0)\over 2}\cot(2\theta_b)-{(\sin\theta_b)^2c_{\phi,\dr_{-1}^2}\over 4}},
\end{align}
read off \eqref{data step relations} with the $\mb e_r$ components $r\cot(2\theta_b)$ and $-r$ of $\smash[b]{\mb P_{\mb u_{F_{0,0}},\dr_{-1}}}$ and $\smash[b]{\mb P_{\mb u_{F_{1,1}},\dr_{-1}}}$ at $\theta=0$. Through this block, which comes from the Navier slip condition, the degree one fields of \eqref{data step relations} enter the third rung.

To move the line of integration in Step 3, the Mellin transforms of the data \eqref{M2 data def} on the line $\Re\lambda=2-\kappa$ and of \eqref{M3 data def} on the line $\Re\lambda=3-\kappa$ must be the restrictions of one function meromorphic on the strip between, with known principal parts. As in part (i), this cannot be obtained as in the proof of Proposition \ref{corner ladder}, from finite norms of the annihilated data at the weights of both lines. The blocks of $\mb F$, $\mb f$ and $\phi^{[2]}$ are annihilated by $\dr^2$, $\dr_{-1}^3$ and $\dr_{-2}^3$. Finite norms of $\dr_{-1}^3\mb f$, for instance, would need three radial derivatives of $\mb f$, while \eqref{um data reg master} and $\mathcal R_j$ provide $j+1/2$ of them, fewer than three when $j<5/2$. Lemma \ref{block matching} again sidesteps this obstacle.

We apply Lemma \ref{block matching} to the pairs $(\mb F,\mb F^{\smash[t]{[3]}})$ at each $\theta$, $(u^{\smash[t]{[1]}}_r,u^{\smash[t]{[2]}}_r)(\ph,0)$, $(\mb f,\mb f^{\smash[t]{[3]}})$ and $(\dr\phi^{\smash[t]{[2]}},\dr\phi^{\smash[t]{[3]}})$, with the norms of \eqref{um data reg master} at $M=2$, of $\mathcal R_j$ and of Step 1. The lemma shows that the Mellin transform of the first member of each pair continues meromorphically to that of the second, with its only pole at the degree of the block and with the coefficients of the block in its principal part \eqref{principal part}. By Lemma \ref{Mellin properties}(i)--(ii), the data of both rungs are therefore the restrictions of
\begin{align}\label{M3 Mellin data}
&\hat{\mb F}_\Omega(\lambda,\theta)=\hat{\mb F}(\lambda_{-2},\theta),\qquad
\hat F_{0,0}(\lambda)=\hat u^{[1]}_r(\lambda_{-1},0),\qquad
\hat F_{0,1}(\lambda)=0,\nonumber\\
&\hat F_{1,0}(\lambda)=\hat f_r(\lambda_{-1}),\qquad
\hat F_{1,1}(\lambda)=\hat f_\theta(\lambda_{-1})+(\sin\theta_b)^2\lambda\hat\phi^{[2]}(\lambda,\theta_b),
\end{align}
with these continuations, written without superscript as in the proof of Proposition \ref{corner ladder}. So the only pole of $\hat F_\bullet$ on the strip is at $\lambda=2$, of order at most two in the channel $\mb F_\Omega$, one in $F_{0,0}$ and three in $F_{1,0}$ and $F_{1,1}$, and the coefficients of its principal parts there are bounded by those of the blocks of \eqref{third rung data classes} and of \eqref{slip block}.

\emph{Step 3: the corner blocks of the solution.} The representation of Proposition \ref{Mellin system initial time} with the data \eqref{M3 Mellin data} defines $\hat{\mb u}$ and $\hat p$, meromorphic on the strip, and by the uniqueness clause of that proposition they are the Mellin transforms of $(\mb u^{[2]},p^{[2]})$ on the lines $\Re\lambda=2-\kappa$ and $\Re\lambda=1-\kappa$, and those of $(\mb u^{[3]},p^{[3]})$ on the lines $\Re\lambda=3-\kappa$ and $\Re\lambda=2-\kappa$. By Lemma \ref{H lemma}(i) and $3-\kappa<\Lambda^*(\theta_b)$, $G_\lambda$ is regular on the strip $2-\kappa\leq\Re\lambda\leq 3-\kappa$, and so are $G_{\lambda,p}=(\lambda_2^2+\partial_\theta^2)G_{\lambda_1}$ and the factors $\lambda^{-1}$ of the pressure formula on $1-\kappa\leq\Re\lambda\leq 2-\kappa$. So the only pole of $\hat{\mb u}$ is at $\lambda=2$, and that of $\hat p$, which carries the data at $\lambda_1$, at $\lambda=1$, both of order at most three; the line of integration is moved across them as in the proof of Proposition \ref{corner ladder}, since by Lemma \ref{block matching} the data are square integrable on the parts $|\Im\lambda|\geq 1$ of the lines uniformly in the line, and \eqref{residue rule} gives
\begin{align*}
\mb u^{[3]}_{\rm c}:=\mb u^{[2]}-\mb u^{[3]}=-\sqrt{2\pi}\Res_{\lambda=2}\brac{r^\lambda\hat{\mb u}},\qquad\qquad
p^{[3]}_{\rm c}:=p^{[2]}-p^{[3]}=-\sqrt{2\pi}\Res_{\lambda=1}\brac{r^\lambda\hat p}.
\end{align*}
Thus, at each $\theta$, $\mb u^{\smash[t]{[3]}}_{\rm c}$ is the block \eqref{higher corner coefficients} of $\mb u^{\smash[t]{[2]}}$ at $(d,m)=(2,3)$ and $p^{\smash[t]{[3]}}_{\rm c}$ that of $p^{\smash[t]{[2]}}$ at $(d,m)=(1,3)$, their coefficients being, by \eqref{principal part}, $-\sqrt{2\pi}$ times the Laurent coefficients of $\hat{\mb u}$ at $\lambda=2$, respectively of $\hat p$ at $\lambda=1$.

\emph{Step 4: the bound on the corner blocks.} The blocks depend only on the singular parts of the data at $\lambda=2$. Let $\hat F^{\rm s}_\bullet$ be obtained from \eqref{M3 Mellin data} by replacing each continued Mellin transform by its principal part \eqref{principal part}, with the coefficients of Step 2,
\begin{align}
\hat{\mb F}^{\rm s}_\Omega(\lambda,\theta)&=-{1\over\sqrt{2\pi}}\sum_{k=1}^2{c_{\mb F,\dr^k}(\theta)\over\lambda_{-2}^k},\nonumber\\
\hat F^{\rm s}_{0,0}(\lambda)&=-{1\over\sqrt{2\pi}\lambda_{-2}}\brac{{u_1(\mb 0)\over 2}\cot(2\theta_b)-{(\sin\theta_b)^2c_{\phi,\dr_{-1}^2}\over 4}},\nonumber\\
\hat F^{\rm s}_{1,0}(\lambda)&=-{1\over\sqrt{2\pi}}\sum_{k=1}^3{c_{f_r,\dr_{-1}^k}\over\lambda_{-2}^k},\nonumber\\
\hat F^{\rm s}_{1,1}(\lambda)&=-{1\over\sqrt{2\pi}}\sum_{k=1}^3{c_{f_\theta,\dr_{-1}^k}+(\sin\theta_b)^2\lambda c_{\phi,\dr_{-2}^k}\over\lambda_{-2}^k},\label{singular data}
\end{align}
and $\hat F^{\rm s}_{0,1}=0$, where $\smash[b]{c_{\mb F,\dr^k}}$ and $\smash[b]{c_{f_r,\dr_{-1}^k}}$, $\smash[b]{c_{f_\theta,\dr_{-1}^k}}$ are the coefficients \eqref{higher corner coefficients} of the blocks $\mb F^{\smash[t]{[3]}}_{\rm c}$ and $\mb f^{\smash[t]{[3]}}_{\rm c}$. Then $\hat F_\bullet-\hat F^{\rm s}_\bullet$ is holomorphic at $\lambda=2$. Let $(\hat{\mb u}^{\rm s},\hat p^{\rm s})$ be given by the representation of Proposition \ref{Mellin system initial time} with the data $\hat F^{\rm s}_\bullet$. Its multipliers being regular on the strip, by Step 3, $\hat{\mb u}-\hat{\mb u}^{\rm s}$ and $\hat p(\lambda_{-1})-\hat p^{\rm s}(\lambda_{-1})$ are holomorphic at $\lambda=2$, while $\hat{\mb u}^{\rm s}$ and $\hat p^{\rm s}(\lambda_{-1})$ are holomorphic on and inside the circle $\mathcal C:=\{|\lambda-2|=\min\{\kappa,1-\kappa\}/2\}$, which lies in the strip, except at its centre.

Write $\mb u^{\smash[t]{[3]}}_{\rm c}=r^2\sum_{k=0}^2(\ln r)^k\mb q_k$ and $p^{\smash[t]{[3]}}_{\rm c}=r\sum_{k=0}^2(\ln r)^kq_k$, as in \eqref{corner classes}. By Step 3 and \eqref{principal part}, $-(k-1)!\mb q_{k-1}/\sqrt{2\pi}$ and $-(k-1)!q_{k-1}/\sqrt{2\pi}$ are the Laurent coefficients at $\lambda_{-2}^{-k}$ of $\hat{\mb u}^{\rm s}(\lambda)$ and $\hat p^{\rm s}(\lambda_{-1})$, so that the Cauchy integral formula for Laurent coefficients gives, for $k=1,2,3$,
\begin{align}\label{Cauchy block coefficients}
\mb q_{k-1}(\theta)=-{\sqrt{2\pi}\over 2\pi i(k-1)!}\oint_{\mathcal C}\lambda_{-2}^{k-1}\hat{\mb u}^{\rm s}(\lambda,\theta)\d\lambda,\qquad
q_{k-1}(\theta)=-{\sqrt{2\pi}\over 2\pi i(k-1)!}\oint_{\mathcal C}\lambda_{-2}^{k-1}\hat p^{\rm s}(\lambda_{-1},\theta)\d\lambda.
\end{align}
From these formulas,
\begin{align*}
\|\mb u^{[3]}_{\rm c}\|_{\mathcal Q_2^3}+\|p^{[3]}_{\rm c}\|_{\mathcal Q_1^3}
&\lesssim\max_{\lambda\in\mathcal C}\brac{\|\hat{\mb u}^{\rm s}(\lambda)\|_{H^2([0,\theta_b])}+\|\hat p^{\rm s}(\lambda_{-1})\|_{H^2([0,\theta_b])}}\\
&\lesssim_{\theta_b,\kappa}\max_{\lambda\in\mathcal C}\brac{\|\hat{\mb F}^{\rm s}_\Omega(\lambda)\|_{H^1([0,\theta_b])}+|\hat F^{\rm s}_{0,0}(\lambda)|+|\hat F^{\rm s}_{1,0}(\lambda)|+|\hat F^{\rm s}_{1,1}(\lambda)|}\\
&\lesssim_{\theta_b,\kappa}\mathcal R_j+|u_1(\mb 0)|+(\sin\theta_b)^2|c_{\phi,\dr_{-1}^2}|,
\end{align*}
which with Step 1 proves \eqref{third rung estimate}. Here the first inequality follows from the Cauchy formulas, $\mathcal C$ having radius less than $1$. The second is \eqref{Mellin bound fixed lambda} with the data $\hat F^{\rm s}_\bullet$, together with the $\theta$ derivative of \eqref{d p from u},
\begin{align*}
\partial_\theta^2\hat p^{\rm s}(\lambda_{-1},\theta)=\lambda_{-1}(\lambda_1\partial_\theta\hat u^{\rm s}_\theta-\partial_\theta^2\hat u^{\rm s}_r)(\lambda,\theta)-\partial_\theta\hat F^{\rm s}_{\Omega,\theta}(\lambda,\theta),
\end{align*}
for the second angular derivative of the pressure. The constants of \eqref{Mellin bound fixed lambda} are bounded on $\mathcal C$: there $\Re\lambda$ lies in the compact subset $[2-\kappa,3-\kappa]$ of $(-\Lambda^*(\theta_b),\Lambda^*(\theta_b))$, so that $\Gamma^H_{\Re\lambda,\theta_b}\gtrsim_{\theta_b,\kappa}1$ by Lemma \ref{H lemma}(ii), while $\<\lambda\>$, $\Upsilon_\lambda$, $\Gamma_{\Re\lambda}$ and $\Gamma_{\Re\lambda_{-1}}$ are bounded in terms of $\theta_b$. The third inequality follows from \eqref{singular data}, since $|\lambda_{-2}|=\min\{\kappa,1-\kappa\}/2$ and $|\lambda|\leq 3$ on $\mathcal C$, the coefficients of the blocks $\mb F^{\trimtop{1.2pt}{[3]}}_{\rm c}$, $\mb f^{\trimtop{1.2pt}{[3]}}_{\rm c}$ and $\phi^{\trimtop{1.2pt}{[3]}}_{\rm c}$ being bounded by $\mathcal R_j$, by \eqref{corner classes} and \eqref{third rung data norm}.

\emph{Part (iii).} On $\theta=\theta_b$ the $\mb e_\theta$ components of $\smash[b]{\mb P_{\mb u_{F_{0,0}},\dr_{-1}}}$ and $\smash[b]{\mb P_{\mb u_{F_{1,1}},\dr_{-1}}}$ are $-r$ and $2\theta_br$, so that by \eqref{data step relations} and part (ii)
\begin{align}\label{u trace expansion}
u_\theta(r,\theta_b)=\mb u(\mb 0)\cdot\mb e_\theta(\theta_b)+r\brac{-{u_1(\mb 0)\over 2}+{\theta_b(\sin\theta_b)^2c_{\phi,\dr_{-1}^2}\over 2}}+u^{[3]}_{{\rm c},\theta}(r,\theta_b)+u^{[3]}_\theta(r,\theta_b),
\end{align}
in which the bracket vanishes by \eqref{C1 in the linear theory}. Applying $\partial_r$,
\begin{align*}
(\sin\theta_b)^2\dot\phi=\partial_ru^{[2]}_\theta(\ph,\theta_b)=\partial_ru^{[3]}_{{\rm c},\theta}(\ph,\theta_b)+\partial_ru^{[3]}_\theta(\ph,\theta_b).
\end{align*}
The remainder $\partial_ru^{\smash[t]{[3]}}_\theta(\ph,\theta_b)$ is bounded by Lemma \ref{trace},
\begin{align*}
\|\partial_ru^{[3]}_\theta(\ph,\theta_b)\|_{\partial_1\Omega,{5\over 2}-\kappa,j+{1\over 2}}\lesssim\|\dr u^{[3]}_\theta(\ph,\theta_b)\|_{\partial_1\Omega,{7\over 2}-\kappa,j+{1\over 2}}\lesssim_{\theta_b}\|\mb u^{[3]}\|_{\Omega,4-\kappa,j,2}.
\end{align*}

It remains to identify $(\sin\theta_b)^{-2}\partial_ru^{\smash[t]{[3]}}_{{\rm c},\theta}(\ph,\theta_b)$ as the block of $\dot\phi$. By Lemma \ref{trace}, applied to $\mb u^{[2]}$, $\|\dot\phi\|_{\partial_1\Omega,3/2-\kappa}$ is finite, and by Lemma \ref{Mellin properties}(i)--(ii) the Mellin transform of $\dot\phi$ on the line $\Re\lambda=1-\kappa$ is $\dot\phi^\wedge(\lambda)=(\sin\theta_b)^{-2}\lambda_1\hat u_\theta(\lambda_1,\theta_b)$, with $\hat{\mb u}$ as in Step 3 of part (ii). Hence $\dot\phi^\wedge$ continues meromorphically across $\Re\lambda=1$, with a pole of order at most three at $\lambda=1$, to the Mellin transform of $(\sin\theta_b)^{-2}\partial_ru^{\smash[t]{[3]}}_\theta(\ph,\theta_b)$ on the line $\Re\lambda=2-\kappa$. Since $\lambda_1r^\lambda=\partial_rr^{\lambda_1}$, Step 3 of part (ii) shows that the block \eqref{higher corner coefficients} of $\dot\phi$ at $(d,m)=(1,3)$ is $(\sin\theta_b)^{-2}\partial_ru^{\smash[t]{[3]}}_{{\rm c},\theta}(\ph,\theta_b)$, so that $\dot\phi^{[2]}=(\sin\theta_b)^{-2}\partial_ru^{\smash[t]{[3]}}_\theta(\ph,\theta_b)$. With the coefficients $\mb q_k$ of $\mb u^{\trimtop{2.5pt}{[3]}}_{\rm c}$ from Step 4 of part (ii), this block is
\begin{align}\label{derived profile block}
(\sin\theta_b)^{-2}\partial_ru^{[3]}_{{\rm c},\theta}(r,\theta_b)=(\sin\theta_b)^{-2}r\sum_{k=0}^2(\ln r)^k\brac{2\mb q_k(\theta_b)+(k+1)\mb q_{k+1}(\theta_b)}\cdot\mb e_\theta(\theta_b),\qquad\mb q_3:=\mb 0,
\end{align}
whose coefficients are bounded by $(\sin\theta_b)^{-2}\|\mb u^{\smash[t]{[3]}}_{\rm c}\|_{\mathcal Q_2^3}$ up to a constant depending on $\theta_b$, by the embedding $H^1([0,\theta_b])\subset C([0,\theta_b])$. With the bound on the remainder above and \eqref{third rung estimate}, this proves \eqref{derived profile third rung}.
\end{proof}

\begin{remark}\label{third rung remark}
Part (iii) is the statement that matters for the nonlinear problem: it shows that the profile derived from a solution with the third-rung structure carries, at the contact point, not only the terms $r$ and $r\ln r$ of \eqref{profile corner coefficients} but also the term ${1\over 2}\smash[b]{c_{\dot\phi,\dr_{-1}^3}}r\ln^2r$. By Lemma \ref{profile channel} below, $\smash[b]{c_{\dot\phi,\dr_{-1}^3}}=\gamma\smash[b]{c_{\phi,\dr_{-2}^3}}+\ell_3(\mb F^{\smash[t]{[3]}}_{\rm c},\mb f^{\smash[t]{[3]}}_{\rm c})$, and by its proof only the coefficients of $r\ln^2r$ in $\mb f^{\trimtop{1.3pt}{[3]}}_{\rm c}$ enter $\ell_3$. In the nonlinear problem these are created by the quadratic nonlinearity, see the proof of Proposition \ref{nonlinear corner structure}(i). It is $\smash[b]{c_{\dot\phi,\dr_{-1}^3}}$ that obstructs the second-order contact angle condition in Theorem \ref{nonlinear initial data}, see \eqref{C2 explicit}.

Parts (ii) and (iii) do not use part (i): their hypothesis is that of Proposition \ref{corner ladder}, whose profile has at worst a double pole at $\lambda=1$. Part (i) is applied in Section \ref{nonlinear existence} to the dotted pair of Theorem \ref{initial jet} at $M=2$, whose profile datum $\dot\phi$ has, by part (iii), a block of order three at $d=1$; see Proposition \ref{nonlinear corner structure}(ii) and Theorem \ref{jet corner structure}. There it makes Lemma \ref{second-order corner values} applicable, which expresses the second compatibility condition \eqref{C2} by the two equations \eqref{C2 explicit}.

In the setting of part (i), by \eqref{log decomposition}, the field \eqref{log field} contributes the term ${1\over 2}\theta_b(\sin\theta_b)^2\smash[b]{c_{\phi,\dr_{-1}^3}}\ln r$ plus a constant to $\partial_ru_\theta(\ph,\theta_b)$, and the other terms have a limit as $r\to 0$, as in the proof of Proposition \ref{corner ladder}. So $\partial_ru_\theta(r,\theta_b)$ has a limit as $r\to 0$ only if $\smash[b]{c_{\phi,\dr_{-1}^3}}=0$, and then the last statement of Proposition \ref{corner ladder} holds unchanged.
\end{remark}

The next lemma expresses the coefficients of the block of the derived profile at $d=1$ through those of the block of $\phi$ at $d=2$ and of the forcing blocks. The dependence is triangular, with the nonvanishing diagonal entry $\gamma(\theta_b)$, and the proof of Theorem \ref{compatible set}(ii) uses it, together with \eqref{nondegeneracy condition}, to solve \eqref{C2 explicit} for the coefficients of $r^2\ln r$ and $r^2\ln^2r$ in the initial profile.

\begin{lemma}[The block of the derived profile]\label{profile channel}
In the setting of Proposition \ref{third rung}(iii), for $m\in\{2,3\}$,
\begin{align}\label{triangular profile channel}
c_{\dot\phi,\dr_{-1}^m}=\gamma c_{\phi,\dr_{-2}^m}+\sum_{k=m+1}^3t_{mk}c_{\phi,\dr_{-2}^k}+\ell_m\brac{\mb F^{[3]}_{\rm c},\mb f^{[3]}_{\rm c}},
\end{align}
where the constants $t_{mk}$ and the linear functionals $\ell_m$ of the coefficients of the forcing blocks depend only on $\theta_b$, and where
\begin{align}\label{gamma formula}
\gamma=\gamma(\theta_b):=-12G_2(\theta_b,\theta_b)={\sin 3\theta_b\sin\theta_b\over\sin 2\theta_b(1+\cos 2\theta_b)}={\sin 3\theta_b\over 4\cos^3\theta_b},
\end{align}
$G_2$ being the Green's function of Lemma \ref{Green's function initial time} at $\lambda=2$. In particular $\gamma(\theta_b)\neq 0$ for every $\theta_b\in(0,\pi/3)$.
\end{lemma}
\begin{proof}
The block of $\dot\phi$ at $d=1$ is given by \eqref{derived profile block} in terms of the coefficients $\mb q_k$ of $\mb u^{\smash[t]{[3]}}_{\rm c}$, and \eqref{Cauchy block coefficients} expresses these through $\hat{\mb u}^{\rm s}$, which is given by the representation of Proposition \ref{Mellin system initial time} with the data \eqref{singular data}. We compute the part of $\smash[b]{c_{\dot\phi,\dr_{-1}^m}}$ that carries the block of $\phi$ at $d=2$, and then evaluate $\gamma$.

\emph{The block of $\dot\phi$.} Comparing the coefficients of $r(\ln r)^{m-1}$ in \eqref{derived profile block} and \eqref{higher corner coefficients} gives
\begin{align*}
c_{\dot\phi,\dr_{-1}^m}=(m-1)!(\sin\theta_b)^{-2}\brac{2\mb q_{m-1}(\theta_b)+m\mb q_m(\theta_b)}\cdot\mb e_\theta(\theta_b),\qquad m=1,2,3.
\end{align*}
By \eqref{singular data}, the coefficients $\smash[b]{c_{\phi,\dr_{-2}^k}}$ enter only $\hat F^{\rm s}_{1,1}$, and by Proposition \ref{Mellin system initial time} this datum enters $\hat u_\theta$ as $-G_\lambda(\theta,\theta_b)\lambda_1\lambda_{-1}\hat F_{1,1}(\lambda)$. So the part of $\hat u^{\rm s}_\theta$ carrying these coefficients is
\begin{align*}
{(\sin\theta_b)^2\over\sqrt{2\pi}}\lambda\lambda_1\lambda_{-1}G_\lambda(\theta,\theta_b)\sum_{k=1}^3{c_{\phi,\dr_{-2}^k}\over\lambda_{-2}^k},
\end{align*}
where $\lambda\lambda_1\lambda_{-1}G_\lambda(\theta,\theta_b)$ is holomorphic at $\lambda=2$, by Lemma \ref{H lemma}(i) and $3-\kappa<\Lambda^*(\theta_b)$, with the value $6G_2(\theta,\theta_b)$ there. Hence, by Taylor's formula at $\lambda=2$, the Laurent coefficient of this part at $\lambda_{-2}^{-k}$ is
\begin{align*}
{(\sin\theta_b)^2\over\sqrt{2\pi}}\sum_{k'=k}^3{1\over(k'-k)!}\eva{\partial_\lambda^{k'-k}\brac{\lambda\lambda_1\lambda_{-1}G_\lambda(\theta,\theta_b)}}_{\lambda=2}c_{\phi,\dr_{-2}^{k'}}.
\end{align*}
The terms of \eqref{singular data} with the coefficients of the forcing blocks enter $\hat u^{\rm s}_\theta$ through multipliers that are regular at $\lambda=2$, and so contribute to its Laurent coefficients linear functionals of the forcing coefficients that, at $\theta=\theta_b$, depend only on $\theta_b$. The datum $\hat F^{\rm s}_{0,0}$ has a simple pole, so it reaches only $\mb q_0$, which enters the display for $\smash[b]{c_{\dot\phi,\dr_{-1}^m}}$\vspace{1pt} only at $m=1$.

By \eqref{Cauchy block coefficients}, $\mb q_{k-1}\cdot\mb e_\theta$ is $-\sqrt{2\pi}/(k-1)!$ times the Laurent coefficient of $\hat u^{\rm s}_\theta$ at $\lambda_{-2}^{-k}$. In the display for $\smash[b]{c_{\dot\phi,\dr_{-1}^m}}$\vspace{1pt}, the term with $\mb q_m$ thus involves only the $\smash[b]{c_{\phi,\dr_{-2}^k}}$ with $k>m$, while the term with $\mb q_{m-1}$ contributes the diagonal entry
\begin{align*}
2(m-1)!(\sin\theta_b)^{-2}\cdot\brac{-{\sqrt{2\pi}\over(m-1)!}}\cdot{(\sin\theta_b)^2\over\sqrt{2\pi}}\cdot 6G_2(\theta_b,\theta_b)=-12G_2(\theta_b,\theta_b),
\end{align*}
and otherwise again only the $\smash[b]{c_{\phi,\dr_{-2}^k}}$ with $k>m$. This is \eqref{triangular profile channel}, with $\gamma$ as in \eqref{gamma formula}. At $m=3$ the term with $\mb q_3=\mb 0$ vanishes, and $\mb q_2$ is reached only by the data of \eqref{singular data} with a pole of order three, so that $\ell_3$ involves only $\smash[b]{c_{f_r,\dr_{-1}^3}}$ and $\smash[b]{c_{f_\theta,\dr_{-1}^3}}$\vspace{1pt}.

\emph{The value of $\gamma$.} At $\lambda=2$ we have $\lambda_1=3$ and $\lambda_{-1}=1$, and in the notation of Lemma \ref{Green's function initial time}
\begin{align*}
w_{22}(\theta_b)=0,\qquad
w_{21}(\theta_b)=2,\qquad
w_{11}(\theta_b)={2\over 3}\sin 3\theta_b\sin\theta_b,\qquad
H(2)=2\sin 2\theta_b\brac{1+\cos 2\theta_b},
\end{align*}
so that
\begin{align*}
G_2(\theta_b,\theta_b)=-{w_{11}(\theta_b)w_{21}(\theta_b)+w_{12}(\theta_b)w_{22}(\theta_b)\over 8H(2)}=-{\sin 3\theta_b\sin\theta_b\over 12\sin 2\theta_b(1+\cos 2\theta_b)},
\end{align*}
which is \eqref{gamma formula}, the last form there following from $\sin 2\theta_b=2\sin\theta_b\cos\theta_b$ and $1+\cos 2\theta_b=2\cos^2\theta_b$. Finally $\gamma=0$ if and only if $\sin 3\theta_b=0$, and the second condition of \eqref{third rung angle condition} with $\kappa\leq{1\over 2}$ gives $\theta_b<\pi/5<\pi/3$.
\end{proof}

The absorption arguments of Section \ref{nonlinear existence} bound quantities of the shape of \eqref{third rung estimate} by estimates in which the quantity itself appears on the right with a small prefactor, and are vacuous unless the quantity is known to be finite. The finiteness is obtained by propagating the decompositions along the Picard iterations and passing to the limit with the following lemma, whose content is that no nonzero corner block lies in a weighted remainder class.

\begin{lemma}[Corner blocks under limits]\label{corner blocks under limits}
Let $K\subset\N_0\times\N_0$ be a finite set of pairs, let
\begin{align*}
\mathcal Q_K:=\set{\sum_{(d,k)\in K}r^d(\ln r)^k\,\mb q_{d,k}(\theta):\mb q_{d,k}\in H^2([0,\theta_b])},
\end{align*}
normed by the sum of the $H^2$ norms of the coefficients, and let $a,j,m$ be such that $a>d+1$ for every $(d,k)\in K$.
\begin{enumerate}
\item Suppose $u_n=B_n+v_n$ with $B_n\in\mathcal Q_K$,
\begin{align*}
S:=\liminf_{n\to\infty}\brac{\|B_n\|_{\mathcal Q_K}+\|v_n\|_{\Omega,a,j,m}}<\infty,
\end{align*}
and $u_n\to u$ in the sense of distributions on $\Omega$. Then $u=B+v$ with $B\in\mathcal Q_K$ and $\|B\|_{\mathcal Q_K}+\|v\|_{\Omega,a,j,m}\leq S$, and this decomposition of $u$ is the unique one of its kind: no nonzero element of $\mathcal Q_K$ has finite $\|\ph\|_{\Omega,a}$ norm.
\item Part (i) holds for scalar functions on $\partial_1\Omega$, with real coefficients $q_{d,k}$, the norms $\|\ph\|_{\partial_1\Omega,a,j}$ and the condition $a>d+{1\over 2}$ in place of $a>d+1$.
\item Parts (i) and (ii) hold for families depending on $t\in[0,T]$, with $\|B_n(t)\|_{\mathcal Q_K}$ bounded in $L^2[0,T]$ against the weight $e^{-\varsigma t}$, the remainders $v_n$ bounded in a norm of Definition \ref{time norms}, and convergence in the sense of distributions in space and time; the uniqueness then holds at almost every time.
\end{enumerate}
\end{lemma}
\begin{proof}
(i) Pass to a subsequence realising the liminf. The coefficients of $B_n$ are bounded in $H^2([0,\theta_b])$ and $v_n$ is bounded in the Hilbert space with norm $\|\ph\|_{\Omega,a,j,m}$, so along a further subsequence the coefficients converge weakly in $H^2$, hence uniformly, and $v_n\rightharpoonup v$ weakly; both convergences imply convergence in the sense of distributions, so $u=B+v$ with $B$ the block with the limiting coefficients, and the norm bound holds by weak lower semicontinuity. For the uniqueness, if $B+v=B'+v'$ then $B-B'=v'-v$ has finite $\|\ph\|_{\Omega,a}$ norm; but for a nonzero element of $\mathcal Q_K$, if $(d_0,k_0)$ is the lexicographically extremal pair with nonvanishing coefficient among those minimising $d$ and then maximising $k$, the element is $\sim r^{d_0}(\ln r)^{k_0}\mb q_{d_0,k_0}(\theta)$ as $r\to 0$ on a set of $\theta$ of positive measure, and
\begin{align*}
\int_0^1r^{-2a}r^{2d_0}(\ln r)^{2k_0}r\d r=\infty
\end{align*}
since $-2a+2d_0+1<-1$. Hence $B=B'$ and $v=v'$.

(ii) The proof of (i) applies verbatim with the measure $\d r$ on $\partial_1\Omega$, the divergent integral being $\int_0^1r^{-2a}r^{2d_0}(\ln r)^{2k_0}\d r$, since $-2a+2d_0<-1$.

(iii) The coefficients and the remainders are now bounded in Hilbert spaces of functions of $t$, so the weak limits are taken in these spaces, and they again imply convergence in the sense of distributions. The norm bound follows by weak lower semicontinuity as before. For the uniqueness, the norms of the two decompositions of the limit are finite at almost every time, and there (i) or (ii) applies.
\end{proof}

\subsection{Positive-time regularity}\label{positive-time regularity}

Similar to the last subsection, here we will need to develop Mellin solution representations for the positive-time system.

We stay in the time variable $t$, although the solution was constructed on the Laplace side. The gain of regularity below uses two time derivatives of the solution, and on the Laplace side these are not obtained by multiplying by $s$. Indeed, the solution is extended by zero to negative times and so jumps at $t=0$, and differentiating the extension produces multiples of $\delta(t)$ and of its derivative:
\begin{align*}
s\hat{\mb u}(s)=(\partial_t\mb u)^\wedge(s)+\mb u(0),\qquad\qquad s^2\hat{\mb u}(s)=(\partial_t^2\mb u)^\wedge(s)+s\mb u(0)+\partial_t\mb u(0),
\end{align*}
where $\partial_t\mb u$ and $\partial_t^2\mb u$ are the derivatives on $t>0$. Unless the initial values vanish, the derivatives of the extension are therefore not square integrable in time, and $s\hat{\mb u}$ and $s^2\hat{\mb u}$ are not square integrable on the line $\Re s=\varsigma$. They are outside the $L^2$-in-time framework of Section \ref{linear existence}, which passes between the two sides by Plancherel's theorem on this line. In $t$ no such terms occur: the derivatives are taken on $t>0$, the initial values $\mb u(0)$ and $\partial_t\mb u(0)$ are the velocities of the two stationary solutions of Theorem \ref{initial jet}, and the trade of Proposition \ref{time for radial trade} is made at a fixed time, with $t$ as a parameter.

\subsubsection{The Mellin representation at positive times}\label{Mellin representation positive time}

The Mellin solution representations are made up of trigonometric terms, and the following lemmas are estimates for them.

\begin{lemma}[Elementary bounds on $\lambda_{\pm 1}$]\label{lambda shift bounds}
For $\lambda\in\C$, we have
\begin{enumerate}
  \item $|\lambda_{\pm 1}|\sim|\lambda|$ when $|\lambda|>2$
  \item $|\lambda_{\pm 1}|\sim 1$ when $|\lambda|\leq 1/2$
  \item $|\lambda_{\pm 1}|\sim|\lambda|$ when $|\Im\lambda|>1/\pi$.
\end{enumerate}
\end{lemma}
\begin{proof}
All three follow from $\big||\lambda|-1\big|\leq|\lambda_{\pm 1}|\leq|\lambda|+1$; for the third, $|\Im\lambda_{\pm 1}|=|\Im\lambda|>1/\pi$ gives $1<\pi\min\{|\lambda|,|\lambda_{\pm 1}|\}$, so the additive $1$ is absorbed on either side.
\end{proof}

In the proof of Proposition \ref{Mellin system positive time} below, the boundary data on $\partial_1\Omega$ enter the velocity through the following functions $A_{\pm 1}(\lambda,\theta)$, for instance $\hat u_{\theta,F_{1,0}}=-\lambda_1A_1(\lambda,\theta)\hat F_{1,0}/(4\lambda)$, and those on $\partial_0\Omega$ through $A_{\pm 1}(\lambda,\theta_b-\theta)$.

\begin{lemma}[Estimates for the coefficients $A_\bullet$]\label{A coefficient estimates}
Suppose $|\Re\lambda_{\pm 1}\theta_b|<\pi$ and $\theta_b\in(0,\pi)$, and write $\Gamma(\lambda,\theta_b):=\Gamma_{\Re\lambda_{-1},\theta_b}\Gamma_{\Re\lambda_1,\theta_b}$ with the sine resonance factor $\Gamma_{a,\theta_b}=\min\{1,\pi-|a\theta_b|\}$ of \eqref{sine resonance factor}. Let
\begin{align*}
A_{\pm 1}(\lambda,\theta) &={\pm\sin(\lambda_{-1}\theta)\sin(\lambda_1\theta_b) -\sin(\lambda_1\theta)\sin(\lambda_{-1}\theta_b) \over\sin(\lambda_{-1}\theta_b)\sin(\lambda_1\theta_b)}.
\end{align*}
Then for $\ell\in\{0,1\}$ and $l>1$ we have
\begin{align*}
  \int_0^{\theta_b}|\partial_\theta^{\ell} A_{\pm 1}|^2\d\theta
  &\lesssim \set{\theta_b^2|\lambda|, {\theta_b|\lambda|\over\<\lambda\>}, {|\lambda|\over\<\lambda\>^2}}^{1\pm 1} {|\lambda|^{2\mb 1_{\{\ell=1\}}}(\{\theta_b,\<\lambda\>^{-1}\}+\mb 1_{\{\ell=1\}}\theta_b^{-1}|\lambda|^{-2}) \over\Gamma(\lambda,\theta_b)^2}\\
  \int_0^{\theta_b}|\partial_\theta^{l} A_{\pm 1}|^2\d\theta
  &\lesssim \set{|\lambda|, {|\lambda|\over\<\lambda\>}, {|\lambda|\over\<\lambda\>^2}}^{1\pm 1} {(|\lambda|^{2\mb 1_{\{\text{$l$ odd}\}}}+|\lambda|^{2l})(\{\theta_b,\<\lambda\>^{-1}\}+\mb 1_{\{\text{$l$ odd}\}}\theta_b^{-1}|\lambda|^{-2}) \over\Gamma(\lambda,\theta_b)^2},
\end{align*}
where $\{\ph\}$ is a choice as in Definition \ref{notation}. In particular, we have
\begin{align*}
\int_0^{\theta_b}|A_{1}|^2\d\theta
  &\lesssim \set{{\theta_b|\lambda|\over\<\lambda\>}, {|\lambda|\over\<\lambda\>^2}}^2 {\{\theta_b,\<\lambda\>^{-1}\} \over\Gamma(\lambda,\theta_b)^2}\\
\int_0^{\theta_b}|\partial_\theta A_{1}|^2\d\theta
  &\lesssim \set{\theta_b, {1\over\<\lambda\>}}^2 {|\lambda|^{2}\{\theta_b,\<\lambda\>^{-1}\} \over\Gamma(\lambda,\theta_b)^2}.
\end{align*}
\end{lemma}
\begin{proof}
Note that
\begin{align*}
\sin(\lambda_{-1}\theta)\sin(\lambda_1\theta_b) &-\sin(\lambda_1\theta)\sin(\lambda_{-1}\theta_b)\\
&={1\over 2}(\cos(\lambda(\theta_b-\theta)+\theta_b+\theta)-\cos(\lambda(\theta_b+\theta)+\theta_b-\theta)\\
&\qquad+\cos(\lambda(\theta_b+\theta)-(\theta_b-\theta))-\cos(\lambda(\theta_b-\theta)-(\theta_b+\theta)))\\
&=\sin(\lambda(\theta_b+\theta))\sin(\theta_b-\theta)-\sin(\lambda(\theta_b-\theta))\sin(\theta_b+\theta),
\end{align*}
\begin{align*}
\cos(\lambda_{-1}\theta)\sin(\lambda_1\theta_b)
&-\cos(\lambda_1\theta)\sin(\lambda_{-1}\theta_b)\\
&={1\over 2}(\sin(\lambda(\theta_b+\theta)+\theta_b-\theta)+\sin(\lambda(\theta_b-\theta)+\theta_b+\theta)\\
&\qquad-\sin(\lambda(\theta_b-\theta)-(\theta_b+\theta))-\sin(\lambda(\theta_b+\theta)-(\theta_b-\theta)))\\
&=\cos(\lambda(\theta_b+\theta))\sin(\theta_b-\theta)+\cos(\lambda(\theta_b-\theta))\sin(\theta_b+\theta),
\end{align*}
\begin{align*}
\cos(\lambda_{-1}\theta)\sin(\lambda_1\theta_b)
&+\cos(\lambda_1\theta)\sin(\lambda_{-1}\theta_b)\\
&={1\over 2}(\sin(\lambda(\theta_b+\theta)+\theta_b-\theta)+\sin(\lambda(\theta_b-\theta)+\theta_b+\theta)\\
&\qquad+\sin(\lambda(\theta_b-\theta)-(\theta_b+\theta))+\sin(\lambda(\theta_b+\theta)-(\theta_b-\theta)))\\
&=\sin(\lambda(\theta_b+\theta))\cos(\theta_b-\theta)+\sin(\lambda(\theta_b-\theta))\cos(\theta_b+\theta),
\end{align*}
and so
\begin{align}
  A_1 &= {\sin(\lambda(\theta_b+\theta))\sin(\theta_b-\theta) -\sin(\lambda(\theta_b-\theta))\sin(\theta_b+\theta) \over\sin(\lambda_{-1}\theta_b)\sin(\lambda_1\theta_b)}\nonumber\\
  \partial_\theta A_1 &= {\lambda_{-1}\cos(\lambda_{-1}\theta)\sin(\lambda_1\theta_b) -\lambda_1\cos(\lambda_1\theta)\sin(\lambda_{-1}\theta_b) \over\sin(\lambda_{-1}\theta_b)\sin(\lambda_1\theta_b)}\nonumber\\
  &=\lambda{\cos(\lambda(\theta_b+\theta))\sin(\theta_b-\theta) +\cos(\lambda(\theta_b-\theta))\sin(\theta_b+\theta) \over\sin(\lambda_{-1}\theta_b)\sin(\lambda_1\theta_b)}\nonumber\\
  &\quad-{\sin(\lambda(\theta_b+\theta))\cos(\theta_b-\theta)+\sin(\lambda(\theta_b-\theta))\cos(\theta_b+\theta) \over\sin(\lambda_{-1}\theta_b)\sin(\lambda_1\theta_b)}.\label{A1 sine forms}
\end{align}
Each of these different but equal formulae for $A_1$ and $\partial_\theta A_1$ is suited to a different region. So for estimates in different regions in the following, we will use the appropriate formula for each case.

Note that when $|\Im\lambda\theta_b|<1$ (together with the assumptions of the lemma, $|\Re\lambda_{\pm 1}\theta_b|<\pi$ and $\theta_b\in(0,\pi)$, so that $|\lambda_{\pm 1}\theta_b|<\pi+1$), we have
\begin{align*}
&\sin(\lambda_{-1}\theta)\sin(\lambda_1\theta_b) -\sin(\lambda_1\theta)\sin(\lambda_{-1}\theta_b)\\
&=\brac{\lambda_{-1}\theta+O(|\lambda_{-1}\theta_b|^3)}\brac{\lambda_1\theta_b+O(|\lambda_1\theta_b|^3)} -\brac{\lambda_1\theta+O(|\lambda_1\theta_b|^3)}\brac{\lambda_{-1}\theta_b+O(|\lambda_{-1}\theta_b|^3)}\\
&=O(|\lambda_1\theta_b||\lambda_{-1}\theta_b|^3)+O(|\lambda_{-1}\theta_b||\lambda_1\theta_b|^3)\\
&\lambda_{-1}\cos(\lambda_{-1}\theta)\sin(\lambda_1\theta_b)-\lambda_1\cos(\lambda_1\theta)\sin(\lambda_{-1}\theta_b)\\
&=\lambda_{-1}\brac{1+O(|\lambda_{-1}\theta|^2)}\brac{\lambda_1\theta_b+O(|\lambda_1\theta_b|^3)} -\lambda_1\brac{1+O(|\lambda_1\theta|^2)}\brac{\lambda_{-1}\theta_b+O(|\lambda_{-1}\theta_b|^3)}\\
&=O(|\lambda_1||\lambda_{-1}\theta_b|^3)+O(|\lambda_{-1}||\lambda_1\theta_b|^3);
\end{align*}
i.e.
\begin{align}
A_1\sin(\lambda_{-1}\theta_b)\sin(\lambda_1\theta_b)
&=O(|\lambda_1\theta_b||\lambda_{-1}\theta_b|^3)+O(|\lambda_{-1}\theta_b||\lambda_1\theta_b|^3)\nonumber\\
\partial_\theta A_1\sin(\lambda_{-1}\theta_b)\sin(\lambda_1\theta_b)
&=O(|\lambda_1||\lambda_{-1}\theta_b|^3)+O(|\lambda_{-1}||\lambda_1\theta_b|^3);\label{A1 first expansion}
\end{align}
also we have
\begin{align*}
\sin&(\lambda(\theta_b+\theta))\sin(\theta_b-\theta) -\sin(\lambda(\theta_b-\theta))\sin(\theta_b+\theta)\\
&=\brac{\lambda u-{(\lambda u)^3\over 6}+O(|\lambda\theta_b|^5)}\brac{v-{v^3\over 6}+O(\theta_b^5)} -\brac{\lambda v-{(\lambda v)^3\over 6}+O(|\lambda\theta_b|^5)}\brac{u-{u^3\over 6}+O(\theta_b^5)}\\
&={\lambda(1-\lambda^2)\over 6}uv(u^2-v^2)+O(|\lambda|\<\lambda\>^4\theta_b^6)
={2\over 3}\lambda(1-\lambda^2)\theta_b\theta(\theta_b^2-\theta^2)+O(|\lambda|\<\lambda\>^4\theta_b^6)\\
&=O(|\lambda|\<\lambda\>^2\theta_b^4),
\end{align*}
where we have written $u=\theta_b+\theta$ and $v=\theta_b-\theta$, so that $uv=\theta_b^2-\theta^2$ and $u^2-v^2=4\theta_b\theta$. Here we are in the regime $|\Im\lambda\theta_b|<1$, in which $\<\lambda\>\theta_b\lesssim 1$, so that the remainder $O(|\lambda|\<\lambda\>^4\theta_b^6)$ is $O(|\lambda|\<\lambda\>^2\theta_b^4)$, as the leading term is. Every term of the expansion carries a factor $\lambda$, since the LHS vanishes at $\lambda=0$; this is what makes the remainder $O(|\lambda|\<\lambda\>^4\theta_b^6)$ rather than merely $O(\<\lambda\>^5\theta_b^6)$, and it is needed when $|\lambda|$ is small. The cancellation of the two leading terms $\lambda uv$ is what produces the fourth power of $\theta_b$, and the second factor $1-\lambda^2$ must be kept, being of size $\<\lambda\>^2$ and not $O(1)$. The remainder differentiates to $O(|\lambda|\<\lambda\>^4\theta_b^5)$, since $u,v\leq 2\theta_b$ and the derivative of the Taylor remainder of order $n$ of $\sin$ at $z$ is $O(|z|^{n-1})$ for bounded $z$, so that a $\theta$-derivative costs a factor $\theta_b^{-1}$. Differentiating in $\theta$, we have
\begin{align*}
\lambda\cos&(\lambda(\theta_b+\theta))\sin(\theta_b-\theta)-\sin(\lambda(\theta_b-\theta))\cos(\theta_b+\theta)\\
&\quad+\lambda\cos(\lambda(\theta_b-\theta))\sin(\theta_b+\theta)-\sin(\lambda(\theta_b+\theta))\cos(\theta_b-\theta)\\
&={2\over 3}\lambda(1-\lambda^2)\theta_b(\theta_b^2-3\theta^2)+O(|\lambda|\<\lambda\>^4\theta_b^5)
=O(|\lambda|\<\lambda\>^2\theta_b^3),
\end{align*}
i.e.
\begin{align}
A_1\sin(\lambda_{-1}\theta_b)\sin(\lambda_1\theta_b)
&=O(|\lambda|\<\lambda\>^2\theta_b^4)\nonumber\\
\partial_\theta A_1\sin(\lambda_{-1}\theta_b)\sin(\lambda_1\theta_b)
&=O(|\lambda|\<\lambda\>^2\theta_b^3).\label{A1 second expansion}
\end{align}

Recall that
\begin{align*}
  \Gamma(\lambda,\theta_b) &= \Gamma_{\Re\lambda_{-1},\theta_b}\Gamma_{\Re\lambda_1,\theta_b},
\end{align*}
which we abbreviate to $\Gamma$ in the displays below. When $|\Im\lambda\theta_b|<1$, Lemma \ref{cosh bound} gives $|\sin(\lambda_{\pm 1}\theta_b)|\gtrsim\Gamma_{\Re\lambda_{\pm 1},\theta_b}|\lambda_{\pm 1}\theta_b|$, and when $|\Im\lambda\theta_b|\geq 1$ it gives $|\sin(\lambda_{\pm 1}\theta_b)|\gtrsim e^{|\Im\lambda|\theta_b}$ with no resonance factor.

When $|\Im\lambda\theta_b|<1$ (which means that $\<\lambda\>\lesssim\theta_b^{-1}$), we bound the numerators of $A_{-1}$ and $\partial_\theta A_{-1}$ term by term, by $|\sin z|\lesssim|z|$ and $|\cos z|\lesssim 1$ for $|z|<\pi+1$, which gives
\begin{align*}
\int_0^{\theta_b}|A_{-1}|^2\d\theta
&\lesssim\int_0^{\theta_b} \brac{|\lambda_1\theta||\lambda_{-1}\theta_b| +|\lambda_{-1}\theta||\lambda_1\theta_b| \over\Gamma_{\Re\lambda_{-1},\theta_b}|\lambda_{-1}\theta_b| \Gamma_{\Re\lambda_{1},\theta_b}|\lambda_1\theta_b|}^2\d\theta
\lesssim\Gamma^{-2}\theta_b\\ &\lesssim\Gamma^{-2}\{\theta_b,\<\lambda\>^{-1}\},\\
\int_0^{\theta_b}|\partial_\theta A_{-1}|^2\d\theta
&\lesssim\int_0^{\theta_b} \brac{|\lambda_1||\lambda_{-1}\theta_b| +|\lambda_{-1}||\lambda_1\theta_b| \over\Gamma_{\Re\lambda_{-1},\theta_b}|\lambda_{-1}\theta_b| \Gamma_{\Re\lambda_{1},\theta_b}|\lambda_1\theta_b|}^2\d\theta\\
&\lesssim\Gamma^{-2}\theta_b^{-1}.
\end{align*}
In this regime we distinguish three cases for $A_1$ and $\partial_\theta A_1$.

When $|\Im\lambda\theta_b|<1$ and $|\lambda|\leq 1/2$ (which means $|\lambda_{\pm 1}|\sim 1\sim\<\lambda\>$ by Lemma \ref{lambda shift bounds}(ii)), the two lines of \eqref{A1 second expansion} read $O(|\lambda|\theta_b^4)$ and $O(|\lambda|\theta_b^3)$, respectively, and give
\begin{align*}
\int_0^{\theta_b}|A_1|^2\d\theta
&\lesssim\int_0^{\theta_b} \brac{O(|\lambda|\theta_b^4) \over\Gamma_{\Re\lambda_{-1},\theta_b}|\lambda_{-1}\theta_b| \Gamma_{\Re\lambda_{1},\theta_b}|\lambda_1\theta_b|}^2\d\theta
\lesssim\int_0^{\theta_b} {|\lambda|^2\theta_b^4 \over\Gamma^2}\d\theta\\
&=\Gamma^{-2}|\lambda|^2\theta_b^5
\lesssim\Gamma^{-2}\{\theta_b,\<\lambda\>^{-1}\} \set{\theta_b^2|\lambda|, {\theta_b|\lambda|\over\<\lambda\>}, {|\lambda|\over\<\lambda\>^2}}^2,\\
\int_0^{\theta_b}|\partial_\theta A_1|^2\d\theta
&\lesssim\int_0^{\theta_b} \brac{O(|\lambda|\theta_b^3) \over\Gamma_{\Re\lambda_{-1},\theta_b}|\lambda_{-1}\theta_b| \Gamma_{\Re\lambda_{1},\theta_b}|\lambda_1\theta_b|}^2\d\theta
\lesssim\int_0^{\theta_b} {|\lambda\theta_b|^2\over\Gamma^2}\d\theta\\
&=\Gamma^{-2}|\lambda|^2\theta_b^3
\lesssim\Gamma^{-2}\theta_b^{-1}\set{\theta_b^2|\lambda|, {\theta_b|\lambda|\over\<\lambda\>}, {|\lambda|\over\<\lambda\>^2}}^2.
\end{align*}
In the last steps we used $\<\lambda\>\theta_b\lesssim 1$, by which $\theta_b$ is $\lesssim$ every choice of $\{\theta_b,\<\lambda\>^{-1}\}$ and $\theta_b^2|\lambda|$ is $\lesssim$ every choice of the brace $\{\theta_b^2|\lambda|,\theta_b|\lambda|/\<\lambda\>,|\lambda|/\<\lambda\>^2\}$. The last steps of the next two cases are the same.

When $|\Im\lambda\theta_b|<1$ and $1/2<|\lambda|\leq 2$, the factors $\lambda_{\pm 1}$ may be close to $0$, and we use \eqref{A1 first expansion} instead, each term of which carries the factor $|\lambda_{-1}\theta_b||\lambda_1\theta_b|$ of the lower bound for the denominator. Since $|\lambda_{\pm 1}|\leq 3$ and $|\lambda|>1/2$, this gives
\begin{align*}
\int_0^{\theta_b}|A_{1}|^2\d\theta
&\lesssim\int_0^{\theta_b} \brac{O(|\lambda_1\theta_b||\lambda_{-1}\theta_b|^3)+O(|\lambda_{-1}\theta_b||\lambda_1\theta_b|^3) \over\Gamma_{\Re\lambda_{-1},\theta_b}|\lambda_{-1}\theta_b| \Gamma_{\Re\lambda_{1},\theta_b}|\lambda_1\theta_b|}^2\d\theta
\lesssim\Gamma^{-2}\theta_b^5\\
&\lesssim\Gamma^{-2}|\lambda|^2\theta_b^5
\lesssim\Gamma^{-2}\{\theta_b,\<\lambda\>^{-1}\} \set{\theta_b^2|\lambda|, {\theta_b|\lambda|\over\<\lambda\>}, {|\lambda|\over\<\lambda\>^2}}^2,\\
\int_0^{\theta_b}|\partial_\theta A_{1}|^2\d\theta
&\lesssim\int_0^{\theta_b} \brac{O(|\lambda_1||\lambda_{-1}\theta_b|^3)+O(|\lambda_{-1}||\lambda_1\theta_b|^3) \over\Gamma_{\Re\lambda_{-1},\theta_b}|\lambda_{-1}\theta_b| \Gamma_{\Re\lambda_{1},\theta_b}|\lambda_1\theta_b|}^2\d\theta
\lesssim\Gamma^{-2}\theta_b^3\\
&\lesssim\Gamma^{-2}|\lambda|^2\theta_b^3
\lesssim\Gamma^{-2}\theta_b^{-1} \set{\theta_b^2|\lambda|, {\theta_b|\lambda|\over\<\lambda\>}, {|\lambda|\over\<\lambda\>^2}}^2.
\end{align*}
When $|\Im\lambda\theta_b|<1$ and $|\lambda|>2$ (which means that $|\lambda_{\pm 1}|\sim|\lambda|\sim\<\lambda\>$ by Lemma \ref{lambda shift bounds}(i), so that the lower bound for the denominator is $\sim\Gamma|\lambda|^2\theta_b^2$), \eqref{A1 second expansion} gives
\begin{align*}
\int_0^{\theta_b}|A_1|^2\d\theta
&\lesssim\int_0^{\theta_b} \brac{O(|\lambda|\<\lambda\>^2\theta_b^4) \over\Gamma_{\Re\lambda_{-1},\theta_b}|\lambda_{-1}\theta_b| \Gamma_{\Re\lambda_{1},\theta_b}|\lambda_1\theta_b|}^2\d\theta
\lesssim\int_0^{\theta_b} {|\lambda|^2\theta_b^4\over\Gamma^2}\d\theta\\
&=\Gamma^{-2}|\lambda|^2\theta_b^5
\lesssim\Gamma^{-2}\{\theta_b,\<\lambda\>^{-1}\} \set{\theta_b^2|\lambda|, {\theta_b|\lambda|\over\<\lambda\>}, {|\lambda|\over\<\lambda\>^2}}^2,\\
\int_0^{\theta_b}|\partial_\theta A_1|^2\d\theta
&\lesssim\int_0^{\theta_b} \brac{O(|\lambda|\<\lambda\>^2\theta_b^3) \over\Gamma_{\Re\lambda_{-1},\theta_b}|\lambda_{-1}\theta_b| \Gamma_{\Re\lambda_{1},\theta_b}|\lambda_1\theta_b|}^2\d\theta
\lesssim\int_0^{\theta_b} {|\lambda|^2\theta_b^2\over\Gamma^2}\d\theta\\
&=\Gamma^{-2}|\lambda|^2\theta_b^3
\lesssim\Gamma^{-2}\theta_b^{-1} \set{\theta_b^2|\lambda|, {\theta_b|\lambda|\over\<\lambda\>}, {|\lambda|\over\<\lambda\>^2}}^2.
\end{align*}

When $|\Im\lambda\theta_b|\geq 1$ (which means that $|\lambda_{\pm 1}|\sim|\lambda|$ by Lemma \ref{lambda shift bounds}(iii), since $|\Im\lambda|\geq\theta_b^{-1}>1/\pi$), also note that
\begin{align*}
\<\lambda\>&\lesssim|\lambda_{-1}|+|\lambda_1|\lesssim{|\Re\lambda_{\pm 1}\theta_b|+|\Im\lambda\theta_b|\over\theta_b}\leq{\pi+|\Im\lambda\theta_b|\over\theta_b}\leq{(\pi+1)|\Im\lambda\theta_b|\over\theta_b}\\
&\lesssim|\Im\lambda|.
\end{align*}
We bound $|\sin z|$ and $|\cos z|$ by $e^{|\Im z|}$, and for real $z$ by $|z|$ and $1$, in the form of $A_1$ and the second form of $\partial_\theta A_1$ in \eqref{A1 sine forms} and in the definition of $A_{-1}$. Then we have
\begin{align*}
\int_0^{\theta_b}|A_1|^2\d\theta
&\lesssim\int_0^{\theta_b}\brac{e^{-2|\Im\lambda|(\theta_b+\theta)}(\theta_b+\theta)^2+e^{-2|\Im\lambda|(\theta_b-\theta)}(\theta_b-\theta)^2}\d\theta
\leq 2\int_0^\infty e^{-2|\Im\lambda|\xi}\xi^2\d\xi\\
&={1\over 2|\Im\lambda|^3}
\lesssim\<\lambda\>^{-3}
\lesssim\{\theta_b,\<\lambda\>^{-1}\} \set{\theta_b^2|\lambda|, {\theta_b|\lambda|\over\<\lambda\>}, {|\lambda|\over\<\lambda\>^2}}^2,\\
\int_0^{\theta_b}|A_{-1}|^2\d\theta
&\lesssim\int_0^{\theta_b}e^{-2|\Im\lambda|(\theta_b-\theta)}\d\theta
\leq{1\over 2|\Im\lambda|}
\lesssim\<\lambda\>^{-1}
\lesssim\{\theta_b,\<\lambda\>^{-1}\},\\
\int_0^{\theta_b}|\partial_\theta A_1|^2\d\theta
&\lesssim\int_0^{\theta_b}\brac{(|\lambda|^2(\theta_b+\theta)^2+1)e^{-2|\Im\lambda|(\theta_b+\theta)}+(|\lambda|^2(\theta_b-\theta)^2+1)e^{-2|\Im\lambda|(\theta_b-\theta)}}\d\theta\\
&\lesssim{|\lambda|^2\over|\Im\lambda|^3}+{1\over|\Im\lambda|}
\lesssim\<\lambda\>^{-1}
\lesssim\{\theta_b,\<\lambda\>^{-1}\} |\lambda|^2\set{\theta_b^2|\lambda|, {\theta_b|\lambda|\over\<\lambda\>}, {|\lambda|\over\<\lambda\>^2}}^2,\\
\int_0^{\theta_b}|\partial_\theta A_{-1}|^2\d\theta
&\lesssim\int_0^{\theta_b} (|\lambda_1|^2+|\lambda_{-1}|^2) e^{-2|\Im\lambda|(\theta_b-\theta)}\d\theta
\lesssim{|\lambda|^2\over|\Im\lambda|}
\lesssim\<\lambda\>^{-1}|\lambda|^2
\lesssim\{\theta_b,\<\lambda\>^{-1}\}|\lambda|^2.
\end{align*}
In the last two steps of each bound we used $\<\lambda\>\lesssim|\Im\lambda|$, and then $|\lambda\theta_b|\geq 1$ and $|\lambda|\sim\<\lambda\>$, by which every choice of $\{\theta_b,\<\lambda\>^{-1}\}$ and of the brace $\{\theta_b^2|\lambda|,\theta_b|\lambda|/\<\lambda\>,|\lambda|/\<\lambda\>^2\}$ is $\gtrsim\<\lambda\>^{-1}$.

Since $\Gamma\leq 1$, adding the bounds of the cases $|\Im\lambda\theta_b|<1$ and $|\Im\lambda\theta_b|\geq 1$ gives
\begin{align*}
  \int_0^{\theta_b}|A_1|^2\d\theta &\lesssim \Gamma^{-2}\{\theta_b,\<\lambda\>^{-1}\} \set{\theta_b^2|\lambda|, {\theta_b|\lambda|\over\<\lambda\>}, {|\lambda|\over\<\lambda\>^2}}^2,\\
  \int_0^{\theta_b}|A_{-1}|^2\d\theta &\lesssim \Gamma^{-2}\{\theta_b,\<\lambda\>^{-1}\},\\
  \int_0^{\theta_b}|\partial_\theta A_1|^2\d\theta &\lesssim \Gamma^{-2}(\{\theta_b,\<\lambda\>^{-1}\}|\lambda|^2 +\theta_b^{-1}) \set{\theta_b^2|\lambda|, {\theta_b|\lambda|\over\<\lambda\>}, {|\lambda|\over\<\lambda\>^2}}^2,\\
  \int_0^{\theta_b}|\partial_\theta A_{-1}|^2\d\theta &\lesssim \Gamma^{-2}(\{\theta_b,\<\lambda\>^{-1}\}|\lambda|^2+\theta_b^{-1}).
\end{align*}
Now
\begin{align*}
\partial_\theta^{2l} A_{(-1)^k} &=(-1)^l{(-1)^k\lambda_{-1}^{2l}\sin(\lambda_{-1}\theta)\sin(\lambda_1\theta_b) -\lambda_1^{2l}\sin(\lambda_1\theta)\sin(\lambda_{-1}\theta_b) \over\sin(\lambda_{-1}\theta_b)\sin(\lambda_1\theta_b)}\\
&=(-1)^l\sum_{j=0}^{2l}\binom{2l}{j}\lambda^jA_{(-1)^{k+j}}\\
\partial_\theta^{2l+1} A_{(-1)^k} &=(-1)^l{(-1)^k\lambda_{-1}^{2l+1}\cos(\lambda_{-1}\theta)\sin(\lambda_1\theta_b)-\lambda_1^{2l+1}\cos(\lambda_1\theta)\sin(\lambda_{-1}\theta_b) \over\sin(\lambda_{-1}\theta_b)\sin(\lambda_1\theta_b)}\\
&=(-1)^l\sum_{j=0}^{2l}\binom{2l}{j}\lambda^j\partial_\theta A_{(-1)^{k+j}},
\end{align*}
so that
\begin{align*}
\int_0^{\theta_b}|\partial_\theta^{2l}A_1|^2\d\theta &\lesssim \Gamma^{-2}\{\theta_b,\<\lambda\>^{-1}\} \set{|\lambda|, {|\lambda|\over\<\lambda\>}, {|\lambda|\over\<\lambda\>^2}}^2 \brac{1+|\lambda|^{4l}}\\
\int_0^{\theta_b}|\partial_\theta^{2l}A_{-1}|^2\d\theta &\lesssim \Gamma^{-2}\{\theta_b,\<\lambda\>^{-1}\} \brac{1+|\lambda|^{4l}}\\
\int_0^{\theta_b}|\partial_\theta^{2l+1} A_1|^2\d\theta &\lesssim \Gamma^{-2}\brac{\{\theta_b,\<\lambda\>^{-1}\}|\lambda|^2 +\theta_b^{-1}}\set{|\lambda|, {|\lambda|\over\<\lambda\>}, {|\lambda|\over\<\lambda\>^2}}^2 \brac{1+|\lambda|^{4l}}\\
\int_0^{\theta_b}|\partial_\theta^{2l+1} A_{-1}|^2\d\theta &\lesssim \Gamma^{-2}\brac{\{\theta_b,\<\lambda\>^{-1}\}|\lambda|^2+\theta_b^{-1}} \brac{1+|\lambda|^{4l}}.
\end{align*}
The last two bounds of the statement follow from the first one with the sign $+$. The bound for $A_1$ is the first one at $\ell=0$ with fewer choices in its braces. For the bound for $\partial_\theta A_1$ we take, in the first one at $\ell=1$, $\theta_b^2|\lambda|$ in both factors of the squared brace and $\theta_b$ in $\{\theta_b,\<\lambda\>^{-1}\}$ when $\theta_b\<\lambda\>\leq 1$, and $|\lambda|/\<\lambda\>^2$ and $\<\lambda\>^{-1}$ when $\theta_b\<\lambda\>>1$. Using $\theta_b|\lambda|\leq 1$, respectively $\theta_b^{-1}<\<\lambda\>$ and $|\lambda|\leq\<\lambda\>$, the RHS of the first bound is then at most $2\Gamma^{-2}|\lambda|^2\min\{\theta_b,\<\lambda\>^{-1}\}^3$, which is at most twice every choice in the last bound of the statement.
\end{proof}

Similar to the initial-time problem, we rewrite our system with the higher spatial scaling terms on the LHS and the lower scaling ones on the RHS, which we treat as data $F_\bullet$. The boundary operator $\partial_\theta u_r+\dr_{-1}u_\theta=r\mb e_r\cdot\mathbb T[\mb u,p]\cdot\mb e_\theta$ below is $r$ times the tangential stress; on $\partial_0\Omega$, where $u_\theta=0$, the Navier slip condition of \eqref{system: time side} reads $\partial_\theta u_r+\dr_{-1}u_\theta=ru_r$. With all the terms on the LHS having the same spatial scaling, we can then obtain solution representations and estimates in Mellin variables:

\begin{proposition}[The Mellin system at positive times]\label{Mellin system positive time}
Let $\theta_b\in(0,\pi)$, and let $a\notin\{0,1,2\}$ with $|a|\theta_b<\pi$ and $|a-2|\theta_b<\pi$. Consider the system
\begin{subequations}\label{system: Mellin positive time}
\begin{alignat}{3}
r^2(\lpc\mb u-\grad p)&=\mb F_\Omega\qquad&&\text{on}\qquad\Omega\\
\grad\cdot\mb u&=0\qquad&&\text{on}\qquad\Omega\\
\partial_\theta u_r+\dr_{-1}u_\theta&=F_{0,0}\qquad&&\text{on}\qquad\partial_0\Omega\\
u_\theta&=F_{0,1}\qquad&&\text{on}\qquad\partial_0\Omega\\
\partial_\theta u_r+\dr_{-1}u_\theta&=F_{1,0}\qquad&&\text{on}\qquad\partial_1\Omega\\
u_\theta&=F_{1,1}\qquad&&\text{on}\qquad\partial_1\Omega.
\end{alignat}
\end{subequations}
Solutions $(\mb u_{F_\bullet},p_{F_\bullet})$ to \eqref{system: Mellin positive time} with all forcing zero except $F_\bullet$ exist depending on the regularity of $F_\bullet$ in the following way (here $j\geq 0$, $l\in\{0,1,2\}$ and $\ell\in\{0,1\}$):
\begin{align*}
  \|\dr^j\partial_\theta^l\mb u_{\mb F_\Omega}\|_{\Omega,a} &\lesssim{1\over\Gamma_{a,\theta_b}\Gamma_{a-2,\theta_b}}\|\{\theta_b,\dr^{-1}\}^{2-l}\dr^j\mb F_\Omega\|_{\Omega,a} \\
  \|\dr^j\partial_\theta^\ell p_{\mb F_\Omega}\|_{\Omega,a-1} &\lesssim{1\over\Gamma_{a-2,\theta_b}}\|\{\theta_b,\dr^{-1}\}^{1-\ell}\dr^{j}\mb F_\Omega\|_{\Omega,a} \\
  \|(r\grad)^{j+2}\mb u_{\mb F_\Omega}\|_{\Omega,a}  &\lesssim{1\over\Gamma_{a,\theta_b}\Gamma_{a-2,\theta_b}}\|\mb F_\Omega\|_{\Omega,a,0,j}\\
  \|(r\grad)^{j+1}p_{\mb F_\Omega}\|_{\Omega,a-1}  &\lesssim{1\over\Gamma_{a-2,\theta_b}}\|\mb F_\Omega\|_{\Omega,a,0,j}\\
  \|\dr^j\partial_\theta^\ell\mb u_{F_{\bullet,0}}\|_{\Omega,a} 
  &\lesssim {1\over\Gamma_{a,\theta_b}\Gamma_{a-2,\theta_b}}\norm{\set{\theta_b^{1\over 2},\<\dr\>^{-{1\over 2}}}\set{\theta_b^{1-\ell},\<\dr\>^{-1}}\dr^jF_{\bullet,0}}_{\partial_\bullet\Omega,a-{1\over 2},\ell}\\
  \|(r\grad)^{j+1}\mb u_{F_{\bullet,0}}\|_{\Omega,a} 
  &\lesssim {1\over\Gamma_{a,\theta_b}\Gamma_{a-2,\theta_b}}\Bigg(\norm{\set{\theta_b^{1\over 2},\<\dr\>^{-{1\over 2}}}F_{\bullet,0}}_{\partial_\bullet\Omega,a-{1\over 2},j} +\theta_b^{-{1\over 2}}\|F_{\bullet,0}\|_{\partial_\bullet\Omega,a-{1\over 2},j^-}\Bigg)\\
  \|(r\grad)^jp_{F_{\bullet,0}}\|_{\Omega,a-1}
  &\lesssim\norm{\set{\theta_b^{1\over 2},\<\dr\>^{-{1\over 2}}} \dr_{-1}^jF_{\bullet,0}}_{\partial_\bullet\Omega,a-{1\over 2}}+{\theta_b^{-{1\over 2}}\over\Gamma_{a-2,\theta_b}}\|\dr_{-1}^{j-1} F_{\bullet,0}\|_{\partial_\bullet\Omega,a-{1\over 2}}\\
  \|(r\grad)^{j}\mb u_{F_{\bullet,1}}\|_{\Omega,a} 
  &\lesssim{1\over\Gamma_{a,\theta_b}\Gamma_{a-2,\theta_b}}\Bigg(\norm{\set{\theta_b^{1\over 2},\<\dr\>^{-{1\over 2}}}\brac{1,\dr_1^{-1}\dr}F_{\bullet,1}}_{\partial_\bullet\Omega,a-{1\over 2},j}\\
  &\pushright{+\theta_b^{-{1\over 2}}\brac{\|F_{\bullet,1}\|_{\partial_\bullet\Omega,a-{1\over 2},j^-}+\|\dr_1^{-1}F_{\bullet,1}\|_{\partial_\bullet\Omega,a-{1\over 2},j}}\Bigg).}
\end{align*}
The pressure estimate at $F_{\bullet,1}$ is the one at $F_{\bullet,0}$ with $\dr F_{\bullet,1}$ in place of $F_{\bullet,0}$. Here $\Gamma_{a,\theta_b}$ is the sine resonance factor \eqref{sine resonance factor}, $\{\ph\}$ is a choice as in Definition \ref{notation}, and $j^-$ is as in Definition \ref{j plus minus}. The solution is moreover unique in the class of $(\mb u,p)$ for which the LHSs above with $j=0$ are finite.
\end{proposition}
\begin{proof}
Since $\mb u$ is divergence-free, \eqref{lpc form} gives $\lpc\mb u=\grad^\perp\grad^\perp\cdot\mb u$, so that
\begin{align*}
r^2(\lpc\mb u-\grad p)
&=r^2(\grad^\perp\grad^\perp\cdot\mb u-\grad p)\\
&=\brac{-r\partial_\theta(r^{-1}(\dr_1u_\theta-\partial_\theta u_r))-r\dr p}\mb e_r+\brac{r\dr(r^{-1}(\dr_1u_\theta-\partial_\theta u_r))-r\partial_\theta p}\mb e_\theta\\
&=\brac{-\partial_\theta(\dr_1u_\theta-\partial_\theta u_r)-r\dr p}\mb e_r+\brac{\dr_{-1}(\dr_1u_\theta-\partial_\theta u_r)-r\partial_\theta p}\mb e_\theta.
\end{align*}
Hence, in Mellin variables, \eqref{system: Mellin positive time} reads
\begin{subequations}
\begin{align}
-\partial_\theta(\lambda_1\hat u_\theta-\partial_\theta\hat u_r)-\lambda_{-1}\hat p(\lambda_{-1},\theta)&=\hat F_{\Omega,r}\label{p from u positive time}\\
\lambda_{-1}(\lambda_1\hat u_\theta-\partial_\theta\hat u_r)-\partial_\theta\hat p(\lambda_{-1},\theta)&=\hat F_{\Omega,\theta}\label{d p from u positive time}\\
\lambda_1\hat u_r+\partial_\theta\hat u_\theta&=0\label{Mellin divergence positive time}\\
(\partial_\theta\hat u_r+\lambda_{-1}\hat u_\theta)(\lambda,0)&=\hat F_{0,0}\\
\hat u_\theta(\lambda,0)&=\hat F_{0,1}\\
(\partial_\theta\hat u_r+\lambda_{-1}\hat u_\theta)(\lambda,\theta_b)&=\hat F_{1,0}\\
\hat u_\theta(\lambda,\theta_b)&=\hat F_{1,1}.
\end{align}
\end{subequations}
Taking $\lambda_{-1}$ times \eqref{d p from u positive time} minus $\partial_\theta$ of \eqref{p from u positive time}, we have
\begin{align*}
(\lambda_{-1}^2+\partial_\theta^2)(\lambda_1\hat u_\theta-\partial_\theta\hat u_r)=\lambda_{-1}\hat F_{\Omega,\theta}-\partial_\theta\hat F_{\Omega,r}.
\end{align*}
Using \eqref{Mellin divergence positive time}, we have
\begin{align*}
(\lambda_{-1}^2+\partial_\theta^2)(\lambda_1^2+\partial_\theta^2)\hat u_\theta &=\lambda_1(\lambda_{-1}\hat F_{\Omega,\theta}-\partial_\theta\hat F_{\Omega,r})\\
(\lambda_1\lambda_{-1}-\partial_\theta^2)\hat u_\theta(\lambda,0)&=\lambda_1\hat F_{0,0}\\
\hat u_\theta(\lambda,0)&=\hat F_{0,1}\\
(\lambda_1\lambda_{-1}-\partial_\theta^2)\hat u_\theta(\lambda,\theta_b)&=\lambda_1\hat F_{1,0}\\
\hat u_\theta(\lambda,\theta_b)&=\hat F_{1,1}.
\end{align*}
Let $f_k$ be the orthonormal eigenfunctions of $\partial_\theta^2$ on $[0,\theta_b]$ with $f(0)=0=f(\theta_b)$, and $-\Lambda_k^2$ their eigenvalues:
\begin{align*}
f_k=\sqrt{2\over\theta_b}\sin\brac{k\pi{\theta\over\theta_b}},\qquad\Lambda_k={k\pi\over\theta_b}\qquad\text{for}\qquad k\in\N.
\end{align*}
They are also the eigenfunctions of $(\lambda_{-1}^2+\partial_\theta^2)(\lambda_1^2+\partial_\theta^2)$ with the boundary condition $f(0)=\partial_\theta^2 f(0)=f(\theta_b)=\partial_\theta^2 f(\theta_b)=0$, with eigenvalue $(\lambda_{-1}^2-\Lambda_k^2)(\lambda_1^2-\Lambda_k^2)$. The functions $f_k'/\Lambda_k$ are again orthonormal, being the eigenfunctions of $\partial_\theta^2$ with the complementary boundary condition $\partial_\theta f(0)=0=\partial_\theta f(\theta_b)$ and the same eigenvalues, and $\partial_\theta$ carries them back onto the $f_k$, since $(f_k'/\Lambda_k)'=-\Lambda_kf_k$. Since the RHS above carries $\hat F_{\Omega,\theta}$ undifferentiated but $\hat F_{\Omega,r}$ differentiated, we expand the two components in the two bases,
\begin{align*}
\hat F_{\Omega,\theta,k}=\<f_k,\hat F_{\Omega,\theta}\>_{L^2([0,\theta_b])},\qquad\qquad
\hat F_{\Omega,r,k}=\inner{{f_k'\over\Lambda_k},\hat F_{\Omega,r}}_{L^2([0,\theta_b])},
\end{align*}
the inner product being conjugate linear in its first slot as in \eqref{Plancherel's identity}, so that $\hat F_{\Omega,\theta}=\sum_k\hat F_{\Omega,\theta,k}f_k$ and $\partial_\theta\hat F_{\Omega,r}=-\sum_k\Lambda_k\hat F_{\Omega,r,k}f_k$, and the whole RHS is expanded in the basis $f_k$, which is the one carrying the boundary conditions of $\hat u_\theta$. Here, unlike in Lemmas \ref{mixed wedge problem} and \ref{Neumann wedge problem}, no boundary condition on $\mb F_\Omega$ is needed for this: $f_k$ vanishes at both ends, so the boundary term of the integration by parts $\<\partial_\theta\hat F_{\Omega,r},f_k\>=-\<\hat F_{\Omega,r},f_k'\>$ is absent. Note that $\{f_k'/\Lambda_k\}_{k\geq 1}$ omits the constant mode $\theta_b^{-1/2}$, which completes it to an orthonormal basis of $L^2([0,\theta_b])$, so that
\begin{align}
\sum_{k=1}^\infty\hat F_{\Omega,r,k}{f_k'\over\Lambda_k}=\hat F_{\Omega,r}-{1\over\theta_b}\int_0^{\theta_b}\hat F_{\Omega,r}(\ph,\xi)\d\xi.\label{mean of F r}
\end{align}
The omission of the mean is harmless for $\mb u_{\mb F_\Omega}$, which sees $\hat F_{\Omega,r}$ only through $\partial_\theta\hat F_{\Omega,r}$, and $\partial_\theta$ annihilates the mean. The mean is carried instead by the standalone term $\hat F_{\Omega,r}(\lambda_1,\theta)$ of the formula for $\hat p_{\mb F_\Omega}$ below, and this term is what makes the residue \eqref{pressure residue Omega} of $p_{\mb F_\Omega}$ nonzero.

Under the reflection $\theta\mapsto\theta_b-\theta$, the system \eqref{system: Mellin positive time} is unchanged if $u_r$, $p$ and $F_{\Omega,r}$ are kept, $u_\theta$ and $F_{\Omega,\theta}$ change sign, and $(F_{0,0},F_{0,1})$ and $(-F_{1,0},-F_{1,1})$ are exchanged. This carries the formulae for $F_{1,\bullet}$ below into those for $F_{0,\bullet}$. We have
\begin{align}
\hat u_{\theta,\mb F_\Omega}&=\sum_{k=1}^\infty{\lambda_1(\lambda_{-1}\hat F_{\Omega,\theta,k}+\Lambda_k\hat F_{\Omega,r,k})f_k \over(\lambda_{-1}^2-\Lambda_k^2)(\lambda_1^2-\Lambda_k^2)} \nonumber\\
\hat u_{\theta,F_{1,1}} &={\lambda_{1}\sin(\lambda_{-1}\theta)\sin(\lambda_1\theta_b) -\lambda_{-1}\sin(\lambda_1\theta)\sin(\lambda_{-1}\theta_b) \over 2\sin(\lambda_{-1}\theta_b)\sin(\lambda_1\theta_b)}\hat F_{1,1} \nonumber\\
&={\lambda A_1(\lambda,\theta)-A_{-1}(\lambda,\theta)\over 2}\hat F_{1,1} \nonumber\\
\hat u_{\theta,F_{0,1}} &={\lambda_{1}\sin(\lambda_{-1}(\theta_b-\theta))\sin(\lambda_1\theta_b) -\lambda_{-1}\sin(\lambda_1(\theta_b-\theta))\sin(\lambda_{-1}\theta_b) \over 2\sin(\lambda_{-1}\theta_b)\sin(\lambda_1\theta_b)}\hat F_{0,1} \nonumber\\
&={\lambda A_1(\lambda,\theta_b-\theta)- A_{-1}(\lambda,\theta_b-\theta)\over 2}\hat F_{0,1} \nonumber\\
\hat u_{\theta,F_{1,0}} &=-\lambda_1{\sin(\lambda_{-1}\theta)\sin(\lambda_1\theta_b) -\sin(\lambda_1\theta)\sin(\lambda_{-1}\theta_b) \over 4\lambda\sin(\lambda_{-1}\theta_b)\sin(\lambda_1\theta_b)}\hat F_{1,0} \nonumber\\
&=-{\lambda_1\over 4\lambda}A_1(\lambda,\theta)\hat F_{1,0} \nonumber\\
\hat u_{\theta,F_{0,0}} &=-\lambda_1{\sin(\lambda_{-1}(\theta_b-\theta))\sin(\lambda_1\theta_b) -\sin(\lambda_1(\theta_b-\theta))\sin(\lambda_{-1}\theta_b) \over 4\lambda\sin(\lambda_{-1}\theta_b)\sin(\lambda_1\theta_b)}\hat F_{0,0} \nonumber\\
&=-{\lambda_1\over 4\lambda}A_1(\lambda,\theta_b-\theta)\hat F_{0,0},\label{u theta formulae positive time}
\end{align}
where $A_{\pm 1}$ is as defined in Lemma \ref{A coefficient estimates}. Using \eqref{Mellin divergence positive time}, we get the version for the $r$ component:
\begin{align}
\hat u_{r,\mb F_\Omega}&=-\sum_{k=1}^\infty{\Lambda_k(\lambda_{-1}\hat F_{\Omega,\theta,k}+\Lambda_k\hat F_{\Omega,r,k})(f_k'/\Lambda_k) \over(\lambda_{-1}^2-\Lambda_k^2)(\lambda_1^2-\Lambda_k^2)} \nonumber\\
\hat u_{r,F_{1,1}} &=-{1\over 2\lambda_1}\brac{\lambda\partial_\theta A_1(\lambda,\theta)-\partial_\theta A_{-1}(\lambda,\theta)}\hat F_{1,1} \nonumber\\
\hat u_{r,F_{0,1}} &={1\over 2\lambda_1}\brac{\lambda\partial_\theta A_1(\lambda,\theta_b-\theta)-\partial_\theta A_{-1}(\lambda,\theta_b-\theta)}\hat F_{0,1} \nonumber\\
\hat u_{r,F_{1,0}} &={1\over 4\lambda}\partial_\theta A_1(\lambda,\theta)\hat F_{1,0} \nonumber\\
\hat u_{r,F_{0,0}} &=-{1\over 4\lambda}\partial_\theta A_1(\lambda,\theta_b-\theta)\hat F_{0,0}.\label{u r formulae positive time}
\end{align}
Using \eqref{p from u positive time} and \eqref{Mellin divergence positive time} with $\lambda_1$ in place of $\lambda$, we have
\begin{align*}
\hat p(\lambda,\theta)=-{1\over\lambda}\brac{{1\over\lambda_2}\partial_\theta(\lambda_2^2+\partial_\theta^2)\hat u_\theta(\lambda_1,\theta)+\hat F_{\Omega,r}(\lambda_1,\theta)},
\end{align*}
in which $\lambda_2^2+\partial_\theta^2$ multiplies $f_k$ by $\lambda_2^2-\Lambda_k^2$, annihilates $\sin(\lambda_2\theta)$ and multiplies $\sin(\lambda\theta)$ by $\lambda_2^2-\lambda^2=4\lambda_1$. This gives the formulae for $\hat p$:
\begin{align}
  \hat p_{\mb F_\Omega} 
  &=-{1\over\lambda}\brac{\hat F_{\Omega,r}(\lambda_1,\theta)+\sum_{k=1}^\infty {\Lambda_k(\lambda_2^2-\Lambda_k^2)(\lambda\hat F_{\Omega,\theta,k}(\lambda_1)+\Lambda_k\hat F_{\Omega,r,k}(\lambda_1))(f_k'/\Lambda_k) \over(\lambda^2-\Lambda_k^2)(\lambda_2^2-\Lambda_k^2)}} \nonumber\\
  &=-{1\over\lambda}\brac{\hat F_{\Omega,r}(\lambda_1,\theta)+\sum_{k=1}^\infty {\Lambda_k(\lambda\hat F_{\Omega,\theta,k}(\lambda_1)+\Lambda_k\hat F_{\Omega,r,k}(\lambda_1))(f_k'/\Lambda_k) \over\lambda^2-\Lambda_k^2}} \nonumber\\
  \hat p_{F_{1,1}}
  &=-{1\over\lambda} {\lambda(\lambda_2^2-\lambda^2)\cos(\lambda\theta)\sin(\lambda_2\theta_b) \over 2\sin(\lambda\theta_b)\sin(\lambda_2\theta_b)}\hat F_{1,1}(\lambda_1)
  =-2\lambda_1{\cos(\lambda\theta)\over \sin(\lambda\theta_b)}\hat F_{1,1}(\lambda_1) \nonumber\\
  \hat p_{F_{0,1}}
  &=2\lambda_1{\cos(\lambda(\theta_b-\theta)) \over \sin(\lambda\theta_b)}\hat F_{0,1}(\lambda_1) \nonumber\\
  \hat p_{F_{1,0}}
  &=-{1\over\lambda}{(\lambda^3-\lambda_2^2\lambda)\cos(\lambda\theta)\sin(\lambda_2\theta_b)\over 4\lambda_1\sin(\lambda\theta_b)\sin(\lambda_2\theta_b)}\hat F_{1,0}(\lambda_1)
  ={\cos(\lambda\theta)\over \sin(\lambda\theta_b)}\hat F_{1,0}(\lambda_1) \nonumber\\
  \hat p_{F_{0,0}}
  &=-{\cos(\lambda(\theta_b-\theta))\over \sin(\lambda\theta_b)}\hat F_{0,0}(\lambda_1).\label{p formulae positive time}
\end{align}
Since $a\notin\{0,1,2\}$, $|a|\theta_b<\pi$ and $|a-2|\theta_b<\pi$, the factors $1/\lambda$, $1/\lambda_1$ and $1/(\sin(\lambda_{-1}\theta_b)\sin(\lambda_1\theta_b))$ of the velocity formulae are finite on the line $\Re\lambda=a-1$, and the factors $1/\lambda$ and $1/\sin(\lambda\theta_b)$ of the pressure formulae on the line $\Re\lambda=a-2$.
The estimates for $\mb u_{\mb F_\Omega}$ are proved similarly to those in Lemmas \ref{mixed wedge problem} and \ref{Neumann wedge problem}. On the line $\Re\lambda=a-1$, on which $\Re\lambda_1=a$ and $\Re\lambda_{-1}=a-2$, the conditions $|a|\theta_b<\pi$ and $|a-2|\theta_b<\pi$ give, exactly as in the proof of Lemma \ref{Neumann wedge problem}, $|\lambda_{\pm 1}^2-\Lambda_k^2|\gtrsim\Gamma_{\Re\lambda_{\pm 1},\theta_b}(\Lambda_k^2+|\lambda_{\pm 1}|^2)$ for every $k\geq 1$. This is where the factor $(\Gamma_{a,\theta_b}\Gamma_{a-2,\theta_b})^{-1}$ comes from. The factor $\Gamma_{a-2,\theta_b}^{-1}$ of $p_{\mb F_\Omega}$ arises in the same way from the denominators $\lambda^2-\Lambda_k^2$ on the line $\Re\lambda=a-2$ of the pressure.

For $\mb u_{F_{1,\bullet}}$ we use Lemma \ref{A coefficient estimates}. On the line $\Re\lambda=a-1$ its factor $\Gamma(\lambda,\theta_b)$ is $\Gamma_{a,\theta_b}\Gamma_{a-2,\theta_b}$, whose inverse square is the prefactor of the estimates below. There also $|\lambda_{-1}|\geq|a-2|>0$, $|\lambda|\geq|a-1|>0$ and $|\lambda_1|\geq|a|>0$, since $a\notin\{0,1,2\}$, so that these three are $\sim_a\<\lambda\>$; below, $|\lambda_{\{-1,0,1\}}|$ is a choice among them as in Definition \ref{notation}. In the bounds for $(r\grad)^j$, $(\lambda,\partial_\theta)^j$ stands for the $\lambda^m\partial_\theta^n$ with $m+n\leq j$, and $\lesssim_j$ accounts for the frame terms $\partial_\theta\mb e_r=\mb e_\theta$ and $\partial_\theta\mb e_\theta=-\mb e_r$, as in the proof of Lemma \ref{mixed wedge problem}. We have
\begin{align*}
\|\dr^j\partial_\theta^\ell u_{\theta,F_{1,0}}\|_{\Omega,a}^2
&=\int_{\Re\lambda=a-1}\int_0^{\theta_b}|\lambda^j\partial_\theta^\ell\hat u_{\theta,F_{1,0}}|^2\d\theta\d\Im\lambda\\
&=\int_{\Re\lambda=a-1} \abs{\lambda_1\over4\lambda}^2 \int_0^{\theta_b}|\lambda^j\partial_\theta^\ell A_1|^2\d\theta|\hat F_{1,0}|^2\d\Im\lambda\\
&\lesssim\int_{\Re\lambda=a-1}|\lambda_1|^2|\lambda|^{2j}\{\theta_b,\<\lambda\>^{-1}\}^4 {|\lambda|^{2\ell}\{\theta_b,\<\lambda\>^{-1}\}+\ell\theta_b^{-1}\over\Gamma(\lambda,\theta_b)^2}|\hat F_{1,0}|^2\d\Im\lambda\\
&\lesssim\int_{\Re\lambda=a-1}|\lambda|^{2j}\{\theta_b,|\lambda_{\{-1,0,1\}}|^{-1}\}^2 {\{\theta_b,|\lambda_{\{-1,0,1\}}|^{-1}\}(1+|\lambda|^{2\ell})\over\Gamma(\lambda,\theta_b)^2}|\hat F_{1,0}|^2\d\Im\lambda\\
&\lesssim\Gamma_{a,\theta_b}^{-2}\Gamma_{a-2,\theta_b}^{-2}\norm{\set{\theta_b^{1\over 2},\<\dr\>^{-{1\over 2}}}^3\dr^j F_{1,0}}_{\partial_1\Omega,a-{1\over 2},\ell}^2.
\end{align*}
The bound for $\|\dr^j\partial_\theta^\ell u_{\theta,F_{1,0}}\|_{\Omega,a}$ uses Lemma \ref{A coefficient estimates}, whose squared brace is $|\lambda|^2\{\theta_b,\<\lambda\>^{-1}\}^4$. In the term with $|\lambda|^{2\ell}$ we then take $\<\lambda\>^{-1}$ in two of these four factors, which absorbs $|\lambda_1|^2$, and in the term with $\theta_b^{-1}$ we take $\theta_b$ in one of them and use $|\lambda_1|^2\lesssim 1+|\lambda|^2$. The other bounds use Lemma \ref{A coefficient estimates} in the same way:
\begin{align*}
&\|(r\grad)^ju_{\theta,F_{1,0}}\|_{\Omega,a}^2\\
&\lesssim_j\int_{\Re\lambda=a-1}\int_0^{\theta_b}|(\lambda,\partial_\theta)^j\hat u_{\theta,F_{1,0}}|^2\d\theta\d\Im\lambda\\
&=\int_{\Re\lambda=a-1} \abs{\lambda_1\over4\lambda}^2 \int_0^{\theta_b}|(\lambda,\partial_\theta)^jA_1|^2\d\theta|\hat F_{1,0}|^2\d\Im\lambda\\
&\lesssim\int_{\Re\lambda=a-1}\{1,|\lambda_{\{-1,0,1\}}|^{-2}\} {\{\theta_b,|\lambda_{\{-1,0,1\}}|^{-1}\} \sum_{j'=0}^{j}|\lambda|^{2j'} +\theta_b^{-1}\sum_{j'=0}^{j-1}|\lambda|^{2j'} \over\Gamma(\lambda,\theta_b)^2}|\hat F_{1,0}|^2\d\Im\lambda\\
&\lesssim\Gamma_{a,\theta_b}^{-2}\Gamma_{a-2,\theta_b}^{-2}\brac{\norm{\set{\theta_b^{1\over 2},\<\dr\>^{-{1\over 2}}}F_{1,0}}_{\partial_1\Omega,a-{1\over 2},j^-}^2+\theta_b^{-1}\|F_{1,0}\|_{\partial_1\Omega,a-{1\over 2},(j^-)^-}^2},
\end{align*}
\begin{align*}
\|\dr^j\partial_\theta^\ell u_{r,F_{1,0}}\|_{\Omega,a}^2
&=\int_{\Re\lambda=a-1}\int_0^{\theta_b}|\lambda^j\partial_\theta^\ell\hat u_{r,F_{1,0}}|^2\d\theta\d\Im\lambda\\
&=\int_{\Re\lambda=a-1}\abs{1\over4\lambda}^2 \int_0^{\theta_b}|\lambda^j\partial_\theta^\ell\partial_\theta A_1|^2\d\theta|\hat F_{1,0}|^2\d\Im\lambda\\
&\lesssim\int_{\Re\lambda=a-1} |\lambda|^{2j}\{\theta_b^{1-\ell},|\lambda_{\{-1,0,1\}}|^{-1}\}^2 {\{\theta_b,|\lambda_{\{-1,0,1\}}|^{-1}\}|\lambda|^{2\ell} \over\Gamma(\lambda,\theta_b)^2}|\hat F_{1,0}|^2\d\Im\lambda\\
&\lesssim\Gamma_{a,\theta_b}^{-2}\Gamma_{a-2,\theta_b}^{-2}\norm{\set{\theta_b^{1\over 2},\<\dr\>^{-{1\over 2}}}\set{\theta_b^{{1\over 2}(1-\ell)},\<\dr\>^{-{1\over 2}}}^2\dr^j F_{1,0}}_{\partial_1\Omega,a-{1\over 2},\ell}^2,
\end{align*}
\begin{align*}
&\|(r\grad)^ju_{r,F_{1,0}}\|_{\Omega,a}^2\\
&\lesssim_j\int_{\Re\lambda=a-1}\int_0^{\theta_b}|(\lambda,\partial_\theta)^j\hat u_{r,F_{1,0}}|^2\d\theta\d\Im\lambda\\
&=\int_{\Re\lambda=a-1}\abs{1\over4\lambda}^2 \int_0^{\theta_b}|(\lambda,\partial_\theta)^j\partial_\theta A_1|^2\d\theta|\hat F_{1,0}|^2\d\Im\lambda\\
&\lesssim\int_{\Re\lambda=a-1} \{|\lambda|,1,|\lambda_{\{-1,0,1\}}|^{-1}\}^2 {\{\theta_b,|\lambda_{\{-1,0,1\}}|^{-1}\}\sum_{j'=0}^{j}|\lambda|^{2j'} +\theta_b^{-1}\sum_{j'=-1}^{j-1}|\lambda|^{2j'}\over\Gamma(\lambda,\theta_b)^2}|\hat F_{1,0}|^2\d\Im\lambda\\
&\lesssim\Gamma_{a,\theta_b}^{-2}\Gamma_{a-2,\theta_b}^{-2}\brac{\norm{\set{\theta_b^{1\over 2},\<\dr\>^{-{1\over 2}}}F_{1,0}}_{\partial_1\Omega,a-{1\over 2},j^-}^2+\theta_b^{-1}\|F_{1,0}\|_{\partial_1\Omega,a-{1\over 2},(j^-)^-}^2}.
\end{align*}
With $j+1$ in place of $j$, and since $(j+1)^-=j$, the bounds for $(r\grad)^ju_{\theta,F_{1,0}}$ and $(r\grad)^ju_{r,F_{1,0}}$ are those of the statement for $(r\grad)^{j+1}\mb u_{F_{1,0}}$. Further,
\begin{align*}
\|(r\grad)^ju_{\theta,F_{1,1}}\|_{\Omega,a}^2
&\lesssim_j\int_{\Re\lambda=a-1}\int_0^{\theta_b}|(\lambda,\partial_\theta)^j\hat u_{\theta,F_{1,1}}|^2\d\theta\d\Im\lambda\\
&=\int_{\Re\lambda=a-1}\int_0^{\theta_b}\abs{(\lambda,\partial_\theta)^j\brac{\lambda A_1-A_{-1}\over 2}}^2\d\theta|\hat F_{1,1}|^2\d\Im\lambda\\
&\lesssim\int_{\Re\lambda=a-1}{\{\theta_b,\<\lambda\>^{-1}\}\sum_{j'=0}^j|\lambda|^{2j'}+\theta_b^{-1}\sum_{j'=0}^{j-1}|\lambda|^{2j'} \over\Gamma(\lambda,\theta_b)^2}|\hat F_{1,1}|^2\d\Im\lambda\\
&\lesssim\Gamma_{a,\theta_b}^{-2}\Gamma_{a-2,\theta_b}^{-2}\brac{\norm{\set{\theta_b^{1\over 2},\<\dr\>^{-{1\over 2}}}F_{1,1}}_{\partial_1\Omega,a-{1\over 2},j}^2+\theta_b^{-1}\|F_{1,1}\|_{\partial_1\Omega,a-{1\over 2},j^-}^2}\\
\|(r\grad)^ju_{r,F_{1,1}}\|_{\Omega,a}^2
&\lesssim_j\int_{\Re\lambda=a-1}\int_0^{\theta_b}|(\lambda,\partial_\theta)^j\hat u_{r,F_{1,1}}|^2\d\theta\d\Im\lambda\\
&=\int_{\Re\lambda=a-1}\abs{1\over\lambda_1}^2\int_0^{\theta_b}\abs{(\lambda,\partial_\theta)^j\brac{\lambda\partial_\theta A_1-\partial_\theta A_{-1}\over 2}}^2\d\theta|\hat F_{1,1}|^2\d\Im\lambda\\
&\lesssim\int_{\Re\lambda=a-1}\abs{\lambda\over\lambda_1}^2{\{\theta_b,\<\lambda\>^{-1}\}\sum_{j'=0}^{j}|\lambda|^{2j'}+\theta_b^{-1}\sum_{j'=-1}^{j-1}|\lambda|^{2j'} \over\Gamma(\lambda,\theta_b)^2}|\hat F_{1,1}|^2\d\Im\lambda\\
&\lesssim\Gamma_{a,\theta_b}^{-2}\Gamma_{a-2,\theta_b}^{-2}\brac{\norm{\set{\theta_b^{1\over 2},\<\dr\>^{-{1\over 2}}}\dr_1^{-1}\dr F_{1,1}}_{\partial_1\Omega,a-{1\over 2},j}^2+\theta_b^{-1}\|\dr_1^{-1}F_{1,1}\|_{\partial_1\Omega,a-{1\over 2},j}^2}.
\end{align*}
Using Lemma \ref{cosh bound} on the line $\Re\lambda=a-2$, where it gives $|\sin(\lambda\theta_b)|\gtrsim\Gamma_{a-2,\theta_b}|\lambda\theta_b|$ when $|\Im\lambda\theta_b|<1$ and $|\sin(\lambda\theta_b)|\gtrsim e^{|\Im\lambda|\theta_b}$ otherwise, we can also estimate $p_\bullet$, with $|\lambda_{\{0,1,2\}}|$ a choice among $|\lambda|$, $|\lambda_1|$ and $|\lambda_2|$, which are $\sim_a\<\lambda\>$ on this line:
\begin{align*}
  \|(r\grad)^jp_{F_{1,0}}\|_{\Omega,a-1}^2
  &\lesssim_j \int_{\Re\lambda=a-2}\int_0^{\theta_b} \abs{(\lambda,\partial_\theta)^j {\cos(\lambda\theta)\over\sin(\lambda\theta_b)}\hat F_{1,0}(\lambda_1)}^2\d\theta\d\Im\lambda \\
  &\lesssim\int_{\Re\lambda=a-2}\mb 1_{\{|\Im\lambda\theta_b|\geq 1\}}\int_0^{\theta_b}|\lambda|^{2j}| \hat F_{1,0}(\lambda_1)|^2e^{2|\Im\lambda|(\theta-\theta_b)}\d\theta\d\Im\lambda\\
  &\quad+\Gamma_{a-2,\theta_b}^{-2}\int_{\Re\lambda=a-2}\mb 1_{\{|\Im\lambda\theta_b|< 1\}}\int_0^{\theta_b} {|\lambda|^{2(j-1)}\over\theta_b^2}|\hat F_{1,0}(\lambda_1)|^2\d\theta\d\Im\lambda\\
  &\lesssim\int_{\Re\lambda=a-2}\{\theta_b,|\lambda_{\{0,1,2\}}|^{-1}\}|\lambda|^{2j}| \hat F_{1,0}(\lambda_1)|^2\d\Im\lambda\\
  &\quad+{1\over\theta_b\Gamma_{a-2,\theta_b}^{2}}\int_{\Re\lambda=a-2} |\lambda|^{2(j-1)}|\hat F_{1,0}(\lambda_1)|^2\d\Im\lambda\\
  &\lesssim\norm{\set{\theta_b^{1\over 2},\<\dr\>^{-{1\over 2}}} \dr_{-1}^j F_{1,0}}_{\partial_1\Omega,a-{1\over 2}}^2+{1\over\theta_b\Gamma_{a-2,\theta_b}^{2}}\|\dr_{-1}^{j-1}F_{1,0}\|_{\partial_1\Omega,a-{1\over 2}}^2\\
  \|(r\grad)^jp_{F_{1,1}}\|_{\Omega,a-1}^2
  &\lesssim\norm{\set{\theta_b^{1\over 2},\<\dr\>^{-{1\over 2}}} \dr_{-1}^j\dr F_{1,1}}_{\partial_1\Omega,a-{1\over 2}}^2+{1\over\theta_b\Gamma_{a-2,\theta_b}^{2}}\|\dr_{-1}^{j-1}\dr F_{1,1}\|_{\partial_1\Omega,a-{1\over 2}}^2.
\end{align*}
The estimates for $\mb u_{F_{0,\bullet}}$ and $p_{F_{0,\bullet}}$ follow by the reflection $\theta\mapsto\theta_b-\theta$. For the uniqueness, the difference of two solutions in the class of the statement solves the homogeneous system and the solution formulae above determine it.
\end{proof}

The implicit constants in Proposition \ref{Mellin system positive time} depend on the distance of $a$ from $0$, $1$ and $2$, and every application below has $a\geq 1-\kappa>0$.

\subsubsection{Trading time derivatives for radial scaling}\label{trading time derivatives}

Using Proposition \ref{Mellin system positive time}, we can obtain control of $(\phi,\mb u,p)$ in a higher spatial scaling provided we have control over the time derivative of the solution.

The next proposition has a parameter $\tilde\kappa\in[0,\kappa)$, which we call the shift. Theorem \ref{linear solution in time} controls $\grad\mb u$ and $p$ at the two weights $-\kappa$ and $0$, and hence at every weight $\tilde\kappa-\kappa$ between them. The trade raises this weight by one in \eqref{time for radial trade estimate 1} and by two in \eqref{time for radial trade estimate 2}. At $\tilde\kappa=0$ the weights reached are those of the corner ladder of Proposition \ref{corner ladder}, and for $\tilde\kappa>0$ they are higher by $\tilde\kappa$. The remainders of the corner expansions then vanish at the contact point, at almost every time, to the orders $r^{2+\tilde\kappa-\kappa}$ for the velocity and the profile and $r^{1+\tilde\kappa-\kappa}$ for the pressure. The shift cannot reach $\kappa$: the velocity would then be at the weights $1$ and $2$ excluded in Proposition \ref{Mellin system positive time}, whose Mellin lines $\Re\lambda=0$ and $\Re\lambda=1$ pass through poles of the Mellin representation, as detailed in the proof. So the exponents of these orders stay below $2$ and $1$. For every shift the positive-time estimate controls less at $t=0$ than $\|\phi_0^{\smash[t]{[2]}}\|_{\partial_1\Omega,5/2-\kappa,l+7/2}$, by $1/2-\tilde\kappa$ in the weight; see Remark \ref{linear theory closing remark}.

\begin{proposition}[Trading time derivatives for radial scaling]\label{time for radial trade}
Take the solution $(\phi,\mb u,p)$ with respect to data $(\mb F,\mb f,\phi_0)$ according to Theorem \ref{linear solution in time} with $l$ replaced by $l+2$. Let $\tilde\kappa\in[0,\kappa)$, with $\kappa\in(1/4,1/2)$ as in that theorem, and suppose that $(3+\tilde\kappa-\kappa)\theta_b<\pi$. If in addition we have control over 
\begin{align*}
&\nn\partial_t\mb u\nn_{[0,T],\Omega:(\varsigma,1-\kappa)(0,l+1,2)}+\nn r\grad\partial_t\mb u\nn_{[0,T],\Omega:(\varsigma,1)(0,l+1,1)}\\
&\quad+\nn\partial_tp\nn_{[0,T],\Omega:(\varsigma,-\kappa)(0,l+1,1)}+\nn\partial_tp\nn_{[0,T],\Omega:(\varsigma,0)(0,l+1,1)}\\
&\quad+\nn\mb F\nn_{[0,T],\Omega:(\varsigma,\tilde\kappa-\kappa)(0,l+2)}+\nn\mb f\nn_{[0,T],\partial_1\Omega:(\varsigma,\tilde\kappa-\kappa+{1\over 2})(0,l+{5\over 2})}+\nn\partial_tf_\theta\nn_{[0,T],\partial_1\Omega:(\varsigma,\tilde\kappa-\kappa-{1\over 2})(0,l+{3\over 2})},
\end{align*}
then we have
\begin{align}
&(\sin\theta_b)^2\nn\phi\nn_{[0,T],\partial_1\Omega:(\varsigma,{3\over 2}+\tilde\kappa-\kappa)(0,l+{7\over 2})}+\nn\mb u^{[1]}\nn_{[0,T],\Omega:(\varsigma,2+\tilde\kappa-\kappa)(0,l+2,2)}+\nn p\nn_{[0,T],\Omega:(\varsigma,1+\tilde\kappa-\kappa)(0,l+2,1)}\nonumber\\
&\lesssim_{\theta_b,\kappa,\tilde\kappa}\nn\mb F\nn_{[0,T],\Omega:(\varsigma,\tilde\kappa-\kappa)(0,l+2)}+\nn\mb f\nn_{[0,T],\partial_1\Omega:(\varsigma,\tilde\kappa-\kappa+{1\over 2})(0,l+{5\over 2})}+\nn\partial_tf_\theta\nn_{[0,T],\partial_1\Omega:(\varsigma,\tilde\kappa-\kappa-{1\over 2})(0,l+{3\over 2})}\nonumber\\
&\quad+\nn\mb u\nn_{[0,T],\Omega:(\varsigma,1+\tilde\kappa-\kappa)(1,l+1,2)}+\nn p\nn_{[0,T],\Omega:(\varsigma,\tilde\kappa-\kappa)(1,l+1,1)},\label{time for radial trade estimate 1}
\end{align}
where
\begin{align*}
\mb u^{[1]} &= \mb u-\mb u(\mb 0).
\end{align*}
If we further have control over 
\begin{align*}
&\nn\partial_t\mb u^{[1]}\nn_{[0,T],\Omega:(\varsigma,2+\tilde\kappa-\kappa)(0,l+1,2)}+\nn\partial_tp\nn_{[0,T],\Omega:(\varsigma,1+\tilde\kappa-\kappa)(0,l+1,1)}\\
&\quad+\nn\mb F\nn_{[0,T],\Omega:(\varsigma,1+\tilde\kappa-\kappa)(0,l+2)}+\nn\mb f\nn_{[0,T],\partial_1\Omega:(\varsigma,1+\tilde\kappa-\kappa+{1\over 2})(0,l+{5\over 2})}+\nn\partial_tf_\theta\nn_{[0,T],\partial_1\Omega:(\varsigma,1+\tilde\kappa-\kappa-{1\over 2})(0,l+{3\over 2})},
\end{align*}
then we have
\begin{align}
&(\sin\theta_b)^2\nn\phi^{[2]}\nn_{[0,T],\partial_1\Omega:(\varsigma,{5\over 2}+\tilde\kappa-\kappa)(0,l+{7\over 2})}+\nn\mb u^{[2]}\nn_{[0,T],\Omega:(\varsigma,3+\tilde\kappa-\kappa)(0,l+2,2)}+\nn p^{[2]}\nn_{[0,T],\Omega:(\varsigma,2+\tilde\kappa-\kappa)(0,l+2,1)}\nonumber\\
&\lesssim_{\theta_b,\kappa,\tilde\kappa}\nn\mb F\nn_{[0,T],\Omega:(\varsigma,1+\tilde\kappa-\kappa)(0,l+2)}+\nn\mb f\nn_{[0,T],\partial_1\Omega:(\varsigma,1+\tilde\kappa-\kappa+{1\over 2})(0,l+{5\over 2})}\nonumber\\
&\quad+\nn\partial_tf_\theta\nn_{[0,T],\partial_1\Omega:(\varsigma,1+\tilde\kappa-\kappa-{1\over 2})(0,l+{3\over 2})}+\nn\mb u^{[1]}\nn_{[0,T],\Omega:(\varsigma,2+\tilde\kappa-\kappa)(1,l+1,2)}+\nn p\nn_{[0,T],\Omega:(\varsigma,1+\tilde\kappa-\kappa)(1,l+1,1)},\label{time for radial trade estimate 2}
\end{align}
where
\begin{align*}
\mb u^{[2]} &= \mb u-\mb u(\mb 0)-{u_1(\mb 0)\over 2}\mb P_0(r,\theta-\theta_b)\\
p^{[2]} &= p+c_3+u_1(\mb 0){\ln r\over\theta_b}\\
\phi^{[2]} &= \phi-cr-u_1(\mb 0){r\ln r\over\theta_b(\sin\theta_b)^2}\\
\mb P_0(r,\theta)&=r\brac{\brac{{\theta\over\theta_b} -{\sin(2\theta)\over\sin(2\theta_b)}}\mb e_{\theta} -\brac{{1\over2\theta_b}-{\cos(2\theta)\over\sin(2\theta_b)}}\mb e_{r}}\\
\mb P_1(r,\theta)&={r\over\theta_b}\brac{\theta\mb e_\theta-{1\over 2}\mb e_r}\\
c_3&=(\sin\theta_b)^2c+\brac{{1\over 2\theta_b}+{1\over\sin(2\theta_b)}}u_1(\mb 0),
\end{align*}
and $c(t)$ and $c_3(t)$ are functions of $t$ alone, given in the proof. In $\mb P_0(r,\theta-\theta_b)$ the argument $\theta-\theta_b$ is inserted in the coefficients only, $\mb e_\theta$ and $\mb e_r$ being taken at $\theta$.
\end{proposition}
\begin{proof}
By Theorem \ref{linear solution in time}, $(\phi,\mb u,p)$ solves \eqref{system: time side}, in particular,
\begin{subequations}\label{u p phi equations}
\begin{alignat}{3}
\lpc\mb u-\grad p&=\mb F\qquad&&\text{on}\quad\Omega\\
\grad\cdot\mb u&=0\qquad&&\text{on}\quad\Omega\\
u_\theta&=0\qquad&&\text{on}\quad\partial_0\Omega\\
u_r-r^{-1}\partial_\theta u_r&=0\qquad&&\text{on}\quad\partial_0\Omega\\
\mb e_r\cdot \mathbb T[\mb u,p]\cdot\mb e_\theta &=f_r\qquad&&\text{on}\quad\partial_1\Omega\\
\mb e_\theta\cdot \mathbb T[\mb u,p]\cdot\mb e_\theta-(\sin\theta_b)^2(\partial_r\phi)&=f_\theta\qquad&&\text{on}\quad\partial_1\Omega\\
\partial_{t}\phi-(\sin\theta_b)^{-2}\partial_r u_\theta &=0\qquad&&\text{on}\quad\partial_1\Omega,
\end{alignat}
\end{subequations}
which we can rewrite as
\begin{subequations}\label{u p phi rewritten}
\begin{alignat}{3}
r^2(\lpc\mb u-\grad p)&=r^2\mb F\qquad&&\text{on}\quad\Omega\\
\grad\cdot\mb u&=0\qquad&&\text{on}\quad\Omega\\
\partial_\theta u_r+\dr_{-1}u_\theta&=ru_r\qquad&&\text{on}\quad\partial_0\Omega\\
u_\theta&=0\qquad&&\text{on}\quad\partial_0\Omega\\
\partial_\theta u_r+\dr_{-1}u_\theta &=rf_r\qquad&&\text{on}\quad\partial_1\Omega\\
u_\theta&=\dr_{-1}^{-1}\dr^{-1}\brac{r^2\mb e_\theta\cdot \mathbb T[\partial_t\mb u,\partial_t p]\cdot\mb e_\theta-r^2\partial_tf_\theta}\nonumber\\
&=-\dr_{-1}^{-1}\dr^{-1}\brac{2r\dr\partial_tu_r+r^2\partial_tp+r^2\partial_tf_\theta}\qquad&&\text{on}\quad\partial_1\Omega.\label{traded condition}
\end{alignat}
\end{subequations}
On $\partial_1\Omega$, the profile evolution equation and the time derivative of the surface tension condition of \eqref{u p phi equations} give
\begin{align*}
\dr u_\theta=(\sin\theta_b)^2r\partial_t\phi,\qquad\qquad(\sin\theta_b)^2\dr\partial_t\phi=r(\mb e_\theta\cdot\mathbb T[\partial_t\mb u,\partial_tp]\cdot\mb e_\theta-\partial_tf_\theta).
\end{align*}
Since $r\dr^{-1}=\dr_{-1}^{-1}r$, this gives the first line of \eqref{traded condition}, and the second follows from
\begin{align}
r\mb e_\theta\cdot\mathbb T[\mb u,p]\cdot\mb e_\theta=-2\dr u_r-rp,\label{normal stress identity}
\end{align}
which holds because $\grad\cdot\mb u=0$. We write $\mathcal H:=2r\dr\partial_tu_r+r^2\partial_tp+r^2\partial_tf_\theta$ on $\partial_1\Omega$ for the bracket in \eqref{traded condition}.

For $M\in\{-\tilde\kappa,0,1,2\}$, let
\begin{align*}
\mb F_\Omega^{[M]} &= r^2\mb F\\
F_{0,0}^{[M]} & = r(u_r-\delta_{M2}u_r(0))|_{\theta=0}\\
F_{0,1}^{[M]} & = 0\\
F_{1,0}^{[M]} & = rf_r\\
F_{1,1}^{[M]} & = -\dr_{-1}^{-1}\dr^{-1}\mathcal H,
\end{align*}
or in Mellin variables
\begin{align}
\hat{\mb F}^{[M]}_\Omega(\lambda,\theta) &= \hat{\mb F}(\lambda_{-2},\theta)\nonumber\\
\hat F_{0,0}^{[M]}(\lambda) & = \begin{cases} \hat u_r(\lambda_{-1},0), & M=-\tilde\kappa,0,1\\ \lambda_{-1}^{-1}(\dr u_r)^\wedge(\lambda_{-1},0), & M=2 \end{cases}\nonumber\\
\hat F_{0,1}^{[M]}(\lambda) &= 0\nonumber\\
\hat F_{1,0}^{[M]}(\lambda) 
& =  \hat f_r(\lambda_{-1})\nonumber\\
\hat F_{1,1}^{[M]}(\lambda) &=-\lambda_{-1}^{-1}\lambda^{-1}\hat{\mathcal H}(\lambda).\label{Mellin data positive time}
\end{align}
Then by Proposition \ref{Mellin system positive time}, whose admissibility condition at the largest weight $a=3+\tilde\kappa-\kappa$ used below is the hypothesis $(3+\tilde\kappa-\kappa)\theta_b<\pi$, we have $(\mb u^{\smash[t]{[M]}}_{F_\bullet},p^{\smash[t]{[M]}}_{F_\bullet})$ with the following bounds, in which $\Gamma(M,\kappa,\theta_b):=\Gamma_{M+1-\kappa,\theta_b}\Gamma_{M-1-\kappa,\theta_b}$ is $\Gamma(\lambda,\theta_b)$ of Lemma \ref{A coefficient estimates} on the line $\Re\lambda=M-\kappa$:
\begin{subequations}\label{time for radial trade estimate master}
\begin{align}
  &\|\mb u^{[M]}_{\mb F_\Omega}\|_{\Omega,M+1+\tilde\kappa-\kappa,2+l,2} +\|p^{[M]}_{\mb F_\Omega}\|_{\Omega,M+\tilde\kappa-\kappa,2+l,1}\nonumber\\ 
  &\hspace{3cm}\lesssim_{\theta_b,\kappa,\tilde\kappa}\|r^2\mb F\|_{\Omega,M+1+\tilde\kappa-\kappa,2+l}
  \lesssim_l\|r\mb F\|_{\Omega,M+\tilde\kappa-\kappa,2+l}\\
  &\Gamma(M,\kappa-\tilde\kappa,\theta_b)\|\mb u^{[M]}_{F_{0,0}}\|_{\Omega,M+1+\tilde\kappa-\kappa,2+l,2} + \Gamma_{M-1+\tilde\kappa-\kappa,\theta_b}\|p^{[M]}_{F_{0,0}}\|_{\Omega,M+\tilde\kappa-\kappa,2+l,1}\nonumber\\ 
  &\hspace{3cm}\lesssim_{\theta_b,\kappa,\tilde\kappa}\|u_r-\delta_{M2} u_r(0)\|_{\partial_0\Omega,M-{1\over 2}+\tilde\kappa-\kappa,{5\over 2}+l}\\
  &\Gamma(M,\kappa-\tilde\kappa,\theta_b)\|\mb u^{[M]}_{F_{1,0}}\|_{\Omega,M+1+\tilde\kappa-\kappa,2+l,2} + \Gamma_{M-1+\tilde\kappa-\kappa,\theta_b}\|p^{[M]}_{F_{1,0}}\|_{\Omega,M+\tilde\kappa-\kappa,2+l,1}\nonumber\\
  &\hspace{3cm}\lesssim_{\theta_b,\kappa,\tilde\kappa}\|f_r\|_{\partial_1\Omega,M-{1\over 2}+\tilde\kappa-\kappa,{5\over 2}+l}\\
  &\Gamma(M,\kappa-\tilde\kappa,\theta_b)\|\mb u^{[M]}_{F_{1,1}}\|_{\Omega,M+1+\tilde\kappa-\kappa,2+l,2} + \Gamma_{M-1+\tilde\kappa-\kappa,\theta_b}\|p^{[M]}_{F_{1,1}}\|_{\Omega,M+\tilde\kappa-\kappa,2+l,1}\nonumber\\
  &\hspace{3cm}\lesssim_{\theta_b,\kappa,\tilde\kappa}(\|\dr\partial_tu_r\|_{\partial_1\Omega,M-{1\over 2}+\tilde\kappa-\kappa,{3\over 2}+l}+\|(\partial_tp,\partial_tf_\theta)\|_{\partial_1\Omega,M-{3\over 2}+\tilde\kappa-\kappa,{3\over 2}+l}).
\end{align}
\end{subequations}
The constants degenerate as $\tilde\kappa\uparrow\kappa$: the weight $M+1+\tilde\kappa-\kappa$ then approaches the values $1$ and $2$ excluded in Proposition \ref{Mellin system positive time}, and on the lines $\Re\lambda=\tilde\kappa-\kappa$ and $\Re\lambda=1+\tilde\kappa-\kappa$ the factor $\lambda_{-1}^{-1}\lambda^{-1}$ of $\hat F^{\smash[t]{[M]}}_{1,1}$ is $\lesssim(\kappa-\tilde\kappa)^{-1}$. In the bound for the surface tension channel $F_{1,1}$, named as in the proof of Proposition \ref{third rung}, $\|\partial_tp\|_{\partial_1\Omega,M-3/2+\tilde\kappa-\kappa,l+3/2}\lesssim_{\theta_b}\|\partial_tp\|_{\Omega,M-1+\tilde\kappa-\kappa,l+1,1}$ by Lemma \ref{trace}. The norm $\|\partial_tp\|_{\Omega,M-1+\tilde\kappa-\kappa,l+1,1}$ is assumed for $M=2$, and for $M=1$ it follows from the assumed $\nn\partial_tp\nn_{[0,T],\Omega:(\varsigma,-\kappa)(0,l+1,1)}$ and $\nn\partial_tp\nn_{[0,T],\Omega:(\varsigma,0)(0,l+1,1)}$, since $-\kappa\leq\tilde\kappa-\kappa<0$. This is the only use of $\partial_tp$.

On the two lowest rungs $M\in\{-\tilde\kappa,0\}$ only the finiteness of the RHS is used, since these rungs serve only to identify $(\mb u,p)$ with $(\mb u^{[0]},p^{[0]})$. For these $M$, $\|r\mb F\|_{\Omega,M+\tilde\kappa-\kappa,2+l}$ and $\|f_r\|_{\partial_1\Omega,M-1/2+\tilde\kappa-\kappa,5/2+l}$ are finite, since $\mathfrak D_{0,l+2}[\mb F,\mb f;\phi_0]$ of Theorem \ref{linear solution in time} is finite. For the surface tension channel $F_{1,1}$ the finiteness is read off from the combined form $F_{1,1}=u_\theta|_{\partial_1\Omega}$ rather than from the separated one, and for the slip channel $F_{0,0}$ from $F_{0,0}=ru_r|_{\theta=0}=\partial_\theta u_r|_{\theta=0}$, by the Navier slip condition of \eqref{u p phi equations}. By Lemma \ref{trace}, $\|u_\theta\|_{\partial_1\Omega,M+1/2+\tilde\kappa-\kappa,l+7/2}$ and $\|\partial_\theta u_r\|_{\partial_0\Omega,M+1/2+\tilde\kappa-\kappa,l+5/2}$ are finite, since $\|\mb u\|_{\Omega,M+1+\tilde\kappa-\kappa,l+2,2}$ is finite by Theorem \ref{linear solution in time} and, at $M=0$, by Hardy's inequality at infinity, Lemma \ref{limit at infinity}.

For $M=1$, by Lemma \ref{trace},
\begin{align*}
\|u_r\|_{\partial_0\Omega,{1\over 2}+\tilde\kappa-\kappa,{5\over 2}+l}&\lesssim_{\theta_b}\|\mb u\|_{\Omega,1+\tilde\kappa-\kappa,l+2,2},\\
\|\dr\partial_tu_r\|_{\partial_1\Omega,{1\over 2}+\tilde\kappa-\kappa,{3\over 2}+l}&\lesssim_{\theta_b}\|r\grad\partial_t\mb u\|_{\Omega,1+\tilde\kappa-\kappa,l+1,1},
\end{align*}
where $\|\mb u\|_{\Omega,1+\tilde\kappa-\kappa,l+2,2}$ is finite as shown for $M=0$, and $\|r\grad\partial_t\mb u\|_{\Omega,1+\tilde\kappa-\kappa,l+1,1}$ by the hypotheses. The other norms on the RHSs of \eqref{time for radial trade estimate master} are finite by the hypotheses.

The estimates \eqref{time for radial trade estimate master} already give us the part of the estimates \eqref{time for radial trade estimate 1} and \eqref{time for radial trade estimate 2} without the $\phi$ term, the control of which we will show later. The solutions for each $M$ could be different due to the presence of poles, so we first find the relationship between $\mb u^{[M]}$ and $\mb u$. Let
\begin{align*}
  \mb u^{[M]} & :=\mb u^{[M]}_{\mb F_\Omega}+\mb u^{[M]}_{F_{0,0}}+\mb u^{[M]}_{F_{1,0}}+\mb u^{[M]}_{F_{1,1}}\\
  p^{[M]} & :=p^{[M]}_{\mb F_\Omega}+p^{[M]}_{F_{0,0}}+p^{[M]}_{F_{1,0}}+p^{[M]}_{F_{1,1}}.
\end{align*}
At $M=-\tilde\kappa$ the weights $M+1+\tilde\kappa-\kappa=1-\kappa$ and $M+\tilde\kappa-\kappa=-\kappa$ in \eqref{time for radial trade estimate master} are those at which Theorem \ref{linear solution in time} bounds $\mb u$ and $p$. Since $(\mb u,p)$ solves \eqref{system: Mellin positive time} with the data $F^{\smash[t]{[-\tilde\kappa]}}_\bullet$, by \eqref{u p phi rewritten}, the uniqueness in Proposition \ref{Mellin system positive time} at $a=1-\kappa$ gives $\mb u=\mb u^{[-\tilde\kappa]}$ and $p=p^{[-\tilde\kappa]}$.

We now locate the poles of the Mellin transforms. By \eqref{u theta formulae positive time} and \eqref{u r formulae positive time}, the factors of $\hat{\mb u}_{F_\bullet}$ multiplying the data have poles only at $\lambda=-1$, from the factor $\lambda_1^{-1}$ of $\hat u_{r,F_{1,1}}$, and where $|\Re\lambda|\geq\pi/\theta_b-1$; the factors $\lambda^{-1}$ of $\hat{\mb u}_{F_{1,0}}$ and $\hat{\mb u}_{F_{0,0}}$ are removable, since $A_1(0,\ph)=0$, and $A_{\pm 1}$ are regular at $\lambda=\pm 1$. By \eqref{p formulae positive time}, those of $\hat p_{F_\bullet}$ have poles only at $\lambda=0$ and where $|\Re\lambda|\geq\pi/\theta_b$. Since $\pi/\theta_b>3+\tilde\kappa-\kappa$, all other poles in the strips below come from the data, namely from the factor $\lambda_{-1}^{-1}\lambda^{-1}$ of $\hat F^{\trimtop{0.9pt}{[M]}}_{1,1}$ at $\lambda=0$ and $\lambda=1$, and from the factor $\lambda_{-1}^{-1}$ of $\hat F^{\smash[t]{[2]}}_{0,0}$ at $\lambda=1$.

Since $-1<-\kappa$ and $\tilde\kappa-\kappa<0$, the factors of $\hat{\mb u}_{F_\bullet}$ multiplying the data have no poles in the strip $\Re\lambda\in[-\kappa,\tilde\kappa-\kappa]$. The data are holomorphic there by Lemma \ref{block matching} with $m=0$, since $\|r\mb F\|_{\Omega,M+\tilde\kappa-\kappa}$, $\|f_r\|_{\partial_1\Omega,M-1/2+\tilde\kappa-\kappa}$, $\|u_\theta\|_{\partial_1\Omega,M+1/2+\tilde\kappa-\kappa}$ and $\|\partial_\theta u_r\|_{\partial_0\Omega,M+1/2+\tilde\kappa-\kappa}$ are finite for $M\in\{-\tilde\kappa,0\}$. Hence $\mb u^{\trimtop{1.5pt}{[-\tilde\kappa]}}_\bullet=\mb u^{\trimtop{1.5pt}{[0]}}_\bullet$. By \eqref{residue rule}, \eqref{u theta formulae positive time}, \eqref{u r formulae positive time} and \eqref{Mellin data positive time}, where $\lambda A_1(\lambda,\theta)$ vanishes at $\lambda=0$ and $A_{-1}(0,\theta)=-2\sin\theta/\sin\theta_b$, we have
\begin{align*}
  (\mb u^{[1]}_{F_{1,1}}-\mb u^{[0]}_{F_{1,1}})(r,\theta)&= {1\over\sqrt{2\pi}}\brac{\int_{\Re\lambda=1+\tilde\kappa-\kappa}-\int_{\Re\lambda=\tilde\kappa-\kappa}} r^{\lambda}\hat{\mb u}_{F_{1,1}}(\lambda,\theta)\d\Im\lambda\\
  &=\sqrt{2\pi}\Res_{\lambda=0}(r^{\lambda}\hat{\mb u}_{F_{1,1}})
  =\sqrt{2\pi}\lim_{\lambda\to 0}\brac{\lambda r^{\lambda}\hat{\mb u}_{F_{1,1}}}\\
  &=\sqrt{2\pi}{\hat{\mathcal H}(0)\over\sin\theta_b} \brac{\sin\theta\mb e_{\theta}-\cos\theta\mb e_{r}}
  =-\sqrt{2\pi}{\hat{\mathcal H}(0)\over\sin\theta_b}\mb e_1.
\end{align*}
The three terms of $\mathcal H$ are controlled separately only on the lines of the rungs $M\geq 1$, and not at $\lambda=0$; but $\mathcal H=-\dr\dr_{-1}u_\theta$ on $\partial_1\Omega$, by \eqref{traded condition}, so that on the line $\Re\lambda=\tilde\kappa-\kappa$ of the rung $M=0$ the Mellin transform $\hat{\mathcal H}=-\lambda\lambda_{-1}\hat u_\theta(\ph,\theta_b)$ is controlled through the combined form, and $\hat{\mathcal H}$ is holomorphic on the strip between the two lines; $\hat{\mathcal H}(0)$ is its value there.
By the location of the poles above, and since the data of the other channels are holomorphic on the strip $\Re\lambda\in[\tilde\kappa-\kappa,1+\tilde\kappa-\kappa]$ by Lemma \ref{block matching} with $m=0$, $\|r\mb F\|_{\Omega,M+\tilde\kappa-\kappa}$, $\|f_r\|_{\partial_1\Omega,M-1/2+\tilde\kappa-\kappa}$ and $\|u_r\|_{\partial_0\Omega,M-1/2+\tilde\kappa-\kappa}$ being finite for $M\in\{0,1\}$, $\hat{\mb u}_{F_{\{(0,0),(1,0),\Omega\}}}$ has no poles in this strip, so $\mb u^{\trimtop{0.9pt}{[0]}}_{F_{\{(0,0),(1,0),\Omega\}}}=\mb u^{\trimtop{0.9pt}{[1]}}_{F_{\{(0,0),(1,0),\Omega\}}}$. Therefore
\[\mb u^{[1]}-\mb u^{[0]}=c_0(t)\mb e_1,\]
for some $c_0(t)$. Since $\|\mb u^{[1]}\|_{\Omega,2+\tilde\kappa-\kappa,2+l,2}<\infty$, by Lemma \ref{continuity at the corner} we have $\mb u^{[1]}(\mb 0)=\mb 0$. So we have $c_0(t)\mb e_1=-\mb u^{[0]}(\mb 0)$. Therefore $\mb u^{[1]}=\mb u-\mb u(\mb 0)$.

By \eqref{residue rule}, \eqref{u theta formulae positive time}, \eqref{u r formulae positive time} and \eqref{Mellin data positive time}, where $(A_1-A_{-1})(1,\theta)=2\theta/\theta_b$, we have
\begin{align*}
  (\mb u^{[2]}_{F_{1,1}}-\mb u^{[1]}_{F_{1,1}})(r,\theta)&= {1\over\sqrt{2\pi}}\brac{\int_{\Re\lambda=2+\tilde\kappa-\kappa}-\int_{\Re\lambda=1+\tilde\kappa-\kappa}} r^{\lambda}\hat{\mb u}_{F_{1,1}}(\lambda,\theta)\d\Im\lambda\\
  &=\sqrt{2\pi}\Res_{\lambda=1}(r^{\lambda}\hat{\mb u}_{F_{1,1}})
  =\sqrt{2\pi}\lim_{\lambda\to 1}\brac{\lambda_{-1} r^{\lambda}\hat{\mb u}_{F_{1,1}}}\\
  &=-\sqrt{2\pi}{\hat{\mathcal H}(1)\over\theta_b} r\brac{\theta\mb e_{\theta} -{1\over 2}\mb e_{r}}\\
  &=-\sqrt{2\pi}\hat{\mathcal H}(1)\mb P_1(r,\theta).
\end{align*}
Note that $\hat F^{\smash[t]{[M]}}_{1,0}(\lambda)=\hat f_r(\lambda_{-1})$ is holomorphic on the strip $\Re\lambda\in[1+\tilde\kappa-\kappa,2+\tilde\kappa-\kappa]$ crossed here, since $\hat f_r$ is holomorphic on $\Re\lambda\in[\tilde\kappa-\kappa,1+\tilde\kappa-\kappa]$ by Lemma \ref{block matching} with $m=0$, $\|f_r\|_{\partial_1\Omega,1/2+\tilde\kappa-\kappa}$ and $\|f_r\|_{\partial_1\Omega,3/2+\tilde\kappa-\kappa}$ being finite by the hypotheses; so $\mb u^{\smash[t]{[2]}}_{F_{1,0}}=\mb u^{\smash[t]{[1]}}_{F_{1,0}}$. Likewise $\hat{\mb F}^{\smash[t]{[M]}}_\Omega(\lambda,\theta)=\hat{\mb F}(\lambda_{-2},\theta)$ is holomorphic on that strip by Lemma \ref{block matching} with $m=0$, since $\|\mb F\|_{\Omega,\tilde\kappa-\kappa}$ and $\|\mb F\|_{\Omega,1+\tilde\kappa-\kappa}$ are finite by the hypotheses, and the factors of $\hat{\mb u}_{\mb F_\Omega}$ multiplying the data have no poles there, as shown above; so $\mb u^{\smash[t]{[2]}}_{\mb F_\Omega}=\mb u^{\smash[t]{[1]}}_{\mb F_\Omega}$.

By \eqref{Mellin data positive time}, $\hat F^{\smash[t]{[2]}}_{0,0}(\lambda)=\lambda_{-1}^{-1}(\dr u_r)^\wedge(\lambda_{-1},0)$, which by Lemma \ref{Mellin properties}(i), (ii) is the Mellin transform of $F^{\smash[t]{[2]}}_{0,0}=r(u_r-u_r(0))=r\dr^{-1}\dr u_r$ on $\partial_0\Omega$; the constant $u_r(0)$, which has no Mellin transform, has disappeared. Here $(\dr u_r)^\wedge(\ph,0)$ is holomorphic on the strip $\Re\lambda\in[\tilde\kappa-\kappa,1+\tilde\kappa-\kappa]$ by Lemma \ref{block matching} with $m=0$, since $\|\dr u_r\|_{\partial_0\Omega,1/2+\tilde\kappa-\kappa}$ and $\|\dr u_r\|_{\partial_0\Omega,3/2+\tilde\kappa-\kappa}=\|\dr u^{\smash[t]{[1]}}_r\|_{\partial_0\Omega,3/2+\tilde\kappa-\kappa}$ are finite by Lemma \ref{trace}, $\|\mb u\|_{\Omega,1+\tilde\kappa-\kappa,l+2,2}$ and $\|\mb u^{[1]}\|_{\Omega,2+\tilde\kappa-\kappa,l+2,2}$ being finite. By \eqref{residue rule}, \eqref{u theta formulae positive time} and \eqref{u r formulae positive time}, where $A_1(1,\theta)=\theta/\theta_b-\sin(2\theta)/\sin(2\theta_b)$, we have
\begin{align*}
  &(\mb u^{[2]}_{F_{0,0}}-\mb u^{[1]}_{F_{0,0}})(r,\theta)\\
  &={1\over\sqrt{2\pi}}\brac{\int_{\Re\lambda=2+\tilde\kappa-\kappa}-\int_{\Re\lambda=1+\tilde\kappa-\kappa}} r^{\lambda}\hat{\mb u}_{F_{0,0}}(\lambda,\theta)\d\Im\lambda\\
  &=\sqrt{2\pi}\Res_{\lambda=1}(r^{\lambda}\hat{\mb u}_{F_{0,0}})
  =\sqrt{2\pi}\lim_{\lambda\to 1}\brac{\lambda_{-1} r^{\lambda}\hat{\mb u}_{F_{0,0}}}\\
  &=-{\sqrt{2\pi}\over 2}(\dr u_r)^\wedge(0,0) r\brac{\brac{{\theta_b-\theta\over\theta_b} -{\sin(2(\theta_b-\theta))\over\sin(2\theta_b)}}\mb e_{\theta} +\brac{{1\over2\theta_b}-{\cos(2(\theta_b-\theta))\over\sin(2\theta_b)}}\mb e_{r}}\\
  &={\sqrt{2\pi}\over 2}(\dr u_r)^\wedge(0,0)\mb P_0(r,\theta-\theta_b).
\end{align*}
Hence, with $c_{0,0}:={\sqrt{2\pi}\over 2}(\dr u_r)^\wedge(0,0)$, where $(\dr u_r)^\wedge(0,0)$ is the value at $\lambda=0$ of the Mellin transform of $\dr u_r(\ph,0)$, and $c_1:=-\sqrt{2\pi}\hat{\mathcal H}(1)$, the residues above give
\begin{align}
\mb u^{[2]}-\mb u^{[1]}=c_{0,0}\mb P_0(r,\theta-\theta_b) +c_{1}\mb P_1(r,\theta).\label{u2 minus u1}
\end{align}
Since $\|\grad\mb u^{\smash[t]{[2]}}\|_{\Omega,2+\tilde\kappa-\kappa,1,1}\lesssim\|\mb u^{\smash[t]{[2]}}\|_{\Omega,3+\tilde\kappa-\kappa,2+l,2}<\infty$, Lemma \ref{continuity at the corner} with $a=1+\tilde\kappa-\kappa$ gives $\grad\mb u^{\smash[t]{[2]}}(\mb 0)=\mb 0$. As $\grad\mb P_0$ and $\grad\mb P_1$ are homogeneous of degree $0$, \eqref{u2 minus u1} shows that $\grad\mb u^{[1]}$ has a limit along each ray from the corner. On $\partial_0\Omega$ the $\mb e_\theta$-components of $\mb u^{[1]}$, $\mb P_0(r,\theta-\theta_b)$ and $\mb P_1$ vanish, so that there $r^{-1}\partial_\theta$ of their $\mb e_r$-components is $\mb e_r\cdot(\mb e_\theta\cdot\grad)$ of the fields, by $r^{-1}\partial_\theta v_r=\mb e_r\cdot(\mb e_\theta\cdot\grad)\mb v+r^{-1}v_\theta$. Since moreover $\partial_\theta P_{1,r}=0$ and $P_{0,\theta}(r,0)=0$, \eqref{u2 minus u1} gives
\begin{align}
\lim_{r\to 0}r^{-1}\partial_\theta u^{[1]}_r(r,0)&=-c_{0,0}r^{-1}\partial_\theta P_{0,r}(r,-\theta_b)=-2c_{0,0},\label{corner limit slip}\\
\lim_{r\to 0}\partial_ru^{[1]}_\theta(r,\theta_b)&=-c_1\partial_rP_{1,\theta}(r,\theta_b)=-c_1.\label{corner limit angle}
\end{align}
By the Navier slip condition of \eqref{u p phi equations}, and since $\mb u(\mb 0)\cdot\mb e_2=0$ by $u_\theta=0$ on $\partial_0\Omega$, the LHS of \eqref{corner limit slip} is $\lim_{r\to 0}u_r(r,0)=u_1(\mb 0)$. By Theorem \ref{linear solution in time}, $\partial_ru_\theta(0,\theta_b)=0$ holds weakly, as a natural boundary condition, so it holds classically once $\partial_ru_\theta(\ph,\theta_b)$ has a limit at the corner (Remark \ref{corner condition and uniqueness}). By \eqref{corner limit angle} it has one, and the LHS of \eqref{corner limit angle} is therefore $0$.
Hence $c_{0,0}=-{1\over 2}u_1(\mb 0)$ and $c_1=0$, and \eqref{u2 minus u1} becomes
\begin{align}
\mb u^{[2]}=\mb u-\mb u(\mb 0)-{u_1(\mb 0)\over 2}\mb P_0(r,\theta-\theta_b).\label{u2 and u relationship}
\end{align}

Now we will find the relationship between $p^{[M]}$ and $p$. Inserting \eqref{Mellin data positive time} in \eqref{p formulae positive time}, where $\hat{\mb F}^{\smash[t]{[M]}}_\Omega(\lambda_1,\theta)=\hat{\mb F}(\lambda_{-1},\theta)$, so that the standalone term $\hat F_{\Omega,r}(\lambda_1,\theta)$ of $\hat p_{\mb F_\Omega}$ becomes the function $\hat F_r(\lambda_{-1},\theta)$ of $\theta\in[0,\theta_b]$, we have
\begin{align*}
  \hat p_{\mb F_\Omega}^{[M]} 
  &=-{1\over\lambda}\brac{\hat F_r(\lambda_{-1},\theta)+\sum_{k=1}^\infty {\Lambda_k(\lambda\hat F_{\theta,k}(\lambda_{-1})+\Lambda_k\hat F_{r,k}(\lambda_{-1}))(f_k'/\Lambda_k) \over\lambda^2-\Lambda_k^2}}\\
  \hat p_{F_{1,1}}^{[M]}
  &={2\cos(\lambda\theta)\over\lambda\sin(\lambda\theta_b)}\hat{\mathcal H}(\lambda_1)\\
  \hat p_{F_{1,0}}^{[M]}
  &={\cos(\lambda\theta)\over\sin(\lambda\theta_b)}\hat f_r(\lambda)\\
  \hat p_{F_{0,0}}^{[M]}
  &=-{\cos(\lambda(\theta_b-\theta))\over\lambda\sin(\lambda\theta_b)} (\dr u_r)^\wedge(\lambda,0),
\end{align*}
where $\hat F_{\theta,k}$ and $\hat F_{r,k}$ are the coefficients of $\hat F_\theta$ and $\hat F_r$ as in the proof of Proposition \ref{Mellin system positive time}; for $M\leq 1$ we have $\lambda^{-1}(\dr u_r)^\wedge(\lambda,0)=\hat u_r(\lambda,0)$. By the location of the poles above, and since the factor $\lambda_1$ of $\hat p_{F_{1,1}}$ in \eqref{p formulae positive time} cancels the factor $\lambda_1^{-1}$ of $\hat F^{\smash[t]{[M]}}_{1,1}(\lambda_1)$, these have no poles in the strip $\Re\lambda\in[-1-\kappa,\tilde\kappa-\kappa]$, so that $p^{\smash[t]{[1]}}=p^{\smash[t]{[0]}}=p^{\smash[t]{[-\tilde\kappa]}}$, and in the strip $\Re\lambda\in[\tilde\kappa-\kappa,1+\tilde\kappa-\kappa]$ their only pole is at $\lambda=0$. By \eqref{residue rule} and \eqref{mean of F r}, the channel $\mb F_\Omega$ gives
\begin{align}
  (p^{[2]}_{\mb F_\Omega}-p^{[1]}_{\mb F_\Omega})(r,\theta)&= {1\over\sqrt{2\pi}}\brac{\int_{\Re\lambda=1+\tilde\kappa-\kappa}-\int_{\Re\lambda=\tilde\kappa-\kappa}} r^{\lambda}\hat p_{\mb F_\Omega}(\lambda,\theta)\d\Im\lambda\nonumber\\
  &=\sqrt{2\pi}\Res_{\lambda=0}(r^{\lambda}\hat p_{\mb F_\Omega})
  =\sqrt{2\pi}\lim_{\lambda\to 0}\brac{\lambda r^{\lambda}\hat p_{\mb F_\Omega}}\nonumber\\
  &=-\sqrt{2\pi}\brac{\hat F_r(-1,\theta)-\sum_{k=1}^\infty\hat F_{r,k}(-1){f_k'\over\Lambda_k}}
  =-{\sqrt{2\pi}\over\theta_b}\int_0^{\theta_b}\hat F_r(-1,\xi)\d\xi.\label{pressure residue Omega}
\end{align}
The function $\lambda\cos(\lambda\theta)/\sin(\lambda\theta_b)$ is even in $\lambda$, so that it is $\theta_b^{-1}+O(\lambda^2)$, and \eqref{residue rule} with $m=2$ gives
\begin{align}
  (p^{[2]}_{F_{1,1}}-p^{[1]}_{F_{1,1}})(r,\theta)&= {1\over\sqrt{2\pi}}\brac{\int_{\Re\lambda=1+\tilde\kappa-\kappa}-\int_{\Re\lambda=\tilde\kappa-\kappa}} r^{\lambda}\hat p_{F_{1,1}}(\lambda,\theta)\d\Im\lambda\nonumber\\
  &=\sqrt{2\pi}\Res_{\lambda=0}(r^{\lambda}\hat p_{F_{1,1}})
  =\sqrt{2\pi}\lim_{\lambda\to 0}\partial_\lambda\brac{\lambda^2 r^{\lambda}\hat p_{F_{1,1}}}\nonumber\\
  &=2\sqrt{2\pi}\lim_{\lambda\to 0}\partial_\lambda\brac{\lambda e^{\lambda\ln r}{\cos(\lambda\theta)\over\sin(\lambda\theta_b)}\hat{\mathcal H}(\lambda_1)}\nonumber\\
  &=2\sqrt{2\pi}\brac{{\ln r\over\theta_b}\hat{\mathcal H}(1)+{1\over\theta_b}\eva{\partial_\lambda\hat{\mathcal H}}_{\lambda=1}}.\label{pressure residue 11}
\end{align}
The channel $F_{1,0}$ contributes a simple pole, $\hat p_{F_{1,0}}$ having the factor $1/\sin(\lambda\theta_b)$ but no $1/\lambda$:
\begin{align}
  (p^{[2]}_{F_{1,0}}-p^{[1]}_{F_{1,0}})(r,\theta)&=\sqrt{2\pi}\Res_{\lambda=0}\brac{r^\lambda{\cos(\lambda\theta)\over\sin(\lambda\theta_b)}\hat f_r(\lambda)}={\sqrt{2\pi}\over\theta_b}\hat f_r(0).\label{pressure residue 10}
\end{align}
The channel $F_{0,0}$ contributes a double pole, as $F_{1,1}$ does, and by the same parity argument
\begin{align}
  (p^{[2]}_{F_{0,0}}-p^{[1]}_{F_{0,0}})(r,\theta)&=\sqrt{2\pi}\lim_{\lambda\to 0}\partial_\lambda\brac{\lambda^2 r^{\lambda}\hat p_{F_{0,0}}}\nonumber\\
  &=-\sqrt{2\pi}\lim_{\lambda\to 0}\partial_\lambda\brac{e^{\lambda\ln r}{\lambda\cos(\lambda(\theta_b-\theta))\over\sin(\lambda\theta_b)}(\dr u_r)^\wedge(\lambda,0)}\nonumber\\
  &=-{\sqrt{2\pi}\over\theta_b}\brac{(\dr u_r)^\wedge(0,0)\ln r+\eva{\partial_\lambda(\dr u_r)^\wedge(\lambda,0)}_{\lambda=0}}.\label{pressure residue 00}
\end{align}
Adding \eqref{pressure residue Omega}--\eqref{pressure residue 00}, with $c_{0,0}$ and $c_1$ as in \eqref{u2 minus u1}, we get
\begin{align}
  (p^{[2]}-p^{[1]})(r,\theta)&=c_3+c_2{\ln r\over\theta_b},\qquad\qquad c_2=-2(c_{0,0}+c_1),\label{p2 and p relationship}
\end{align}
where $c_3=c_3(t)$, independent of $r$ and $\theta$, is the sum of the terms without $\ln r$ in \eqref{pressure residue Omega}--\eqref{pressure residue 00}. Inserting $c_{0,0}=-{1\over 2}u_1(\mb 0)$ and $c_1=0$, found from \eqref{corner limit slip} and \eqref{corner limit angle}, we get $c_2=u_1(\mb 0)$: the coefficient of $\theta_b^{-1}\ln r$ in $p^{\smash[t]{[2]}}-p$ is the contact-point velocity, and it is supplied by the slip channel $F_{0,0}$, the surface tension channel $F_{1,1}$ contributing nothing to it under the contact angle condition $\partial_ru_\theta(0,\theta_b)=0$. Therefore, as $p=p^{[1]}$,
\begin{align}
  p^{[2]}&=p+c_3+u_1(\mb 0){\ln r\over\theta_b}.\label{p2 formula}
\end{align}

Now we will obtain the estimates for $\phi$. By \eqref{u p phi equations} and \eqref{normal stress identity}, we have
\begin{align*}
\dr\phi  &= (\sin\theta_b)^{-2}r\brac{\mb e_\theta\cdot\mathbb T[\mb u, p]\cdot\mb e_\theta -f_\theta}\\
&= (\sin\theta_b)^{-2}\brac{-2\dr u_r -rp-rf_\theta} \qquad\qquad\text{on}\quad\partial_1\Omega.
\end{align*}
On $\partial_1\Omega$, by $\mb u^{\smash[t]{[1]}}=\mb u-\mb u(\mb 0)$, \eqref{u2 and u relationship} and \eqref{p2 formula}, we have
\begin{align*}
\dr u_r &= \dr(u_r-u_r(0)) =\dr u_r^{[1]}\\
\dr u_r &= \dr u_r^{[2]}-{u_1(\mb 0)\over 2}\brac{{1\over 2\theta_b}-{1\over\sin(2\theta_b)}}r
= \dr u_r^{[2]}- \dr(c'(t)r)\\
rp & = rp^{[2]} -\brac{c_3(t)-{u_1(\mb 0)\over\theta_b}}r -{u_1(\mb 0)\over\theta_b}\dr(r\ln r),
\end{align*}
where $c'(t):=u_1(\mb 0)(1/(4\theta_b)-1/(2\sin(2\theta_b)))$. Collecting the terms proportional to $r$, we get that $\phi^{\smash[t]{[2]}}$ of the statement, with
\begin{align}\label{c formula}
c(t)=(\sin\theta_b)^{-2}\brac{2c'(t)+c_3(t)-{u_1(\mb 0)\over\theta_b}},
\end{align}
satisfies
\begin{align*}
\dr\phi^{[2]}=(\sin\theta_b)^{-2}(-2\dr u^{[2]}_r-rp^{[2]}-rf_\theta)\qquad\qquad\text{on}\quad\partial_1\Omega.
\end{align*}
By the definition of $c'(t)$, \eqref{c formula} is the relation between $c$ and $c_3$ of the statement. At each boundary weight $a$ we have $\|\phi\|_{\partial_1\Omega,a,l+7/2}^2=\|\phi\|_{\partial_1\Omega,a}^2+\|\dr\phi\|_{\partial_1\Omega,a,l+5/2}^2$, and $\|\phi\|_{\partial_1\Omega,a}\leq\sup|\lambda|^{-1}\|\dr\phi\|_{\partial_1\Omega,a}$ by \eqref{Plancherel's identity}, the supremum being over the line $\Re\lambda=a-1/2$. For $a>1/2$ this uses that $\phi$, and with it $\phi^{\smash[t]{[2]}}$, vanishes at the corner, by the profile evolution equation, since $\partial_ru_\theta(0,\theta_b)=0$ and $\phi_0(0)=0$. So
\begin{align*}
\|\phi\|_{\partial_1\Omega,a,l+{7\over 2}}\leq\brac{1+\sup_{\Re\lambda=a-{1\over 2}}|\lambda|^{-1}}\|\dr\phi\|_{\partial_1\Omega,a,l+{5\over 2}},
\end{align*}
where the supremum is $(\kappa-\tilde\kappa)^{-1}$ at $a=1/2+\tilde\kappa-\kappa$, and is below $2$ and $2/3$ at $a=3/2+\tilde\kappa-\kappa$ and $a=5/2+\tilde\kappa-\kappa$.
With this and Lemma \ref{trace}, the identities above give
\begin{align*}
(\sin\theta_b)^2\|\phi\|_{\partial_1\Omega,{1\over 2}+\tilde\kappa-\kappa,{7\over 2}+l}
&\lesssim_{\kappa-\tilde\kappa}\|-2\dr u_r -rp-rf_\theta\|_{\partial_1\Omega,{1\over 2}+\tilde\kappa-\kappa,{5\over 2}+l}\\
&\lesssim_l\|\mb u\|_{\Omega,1+\tilde\kappa-\kappa,l+2,2} +\|p\|_{\Omega,\tilde\kappa-\kappa,l+2,1} +\|f_\theta\|_{\partial_1\Omega,-{1\over 2}+\tilde\kappa-\kappa,{5\over 2}+l}\\
(\sin\theta_b)^2\|\phi\|_{\partial_1\Omega,{3\over 2}+\tilde\kappa-\kappa,{7\over 2}+l}
&\lesssim\|-2\dr u_r -rp-rf_\theta\|_{\partial_1\Omega,{3\over 2}+\tilde\kappa-\kappa,{5\over 2}+l}\\
&\lesssim_l\|\mb u^{[1]}\|_{\Omega,2+\tilde\kappa-\kappa,l+2,2} +\|p\|_{\Omega,1+\tilde\kappa-\kappa,l+2,1} +\|f_\theta\|_{\partial_1\Omega,{1\over 2}+\tilde\kappa-\kappa,{5\over 2}+l}\\
(\sin\theta_b)^2\|\phi^{[2]}\|_{\partial_1\Omega,{5\over 2}+\tilde\kappa-\kappa,{7\over 2}+l}
&\lesssim_l\|\mb u^{[2]}\|_{\Omega,3+\tilde\kappa-\kappa,l+2,2} +\|p^{[2]}\|_{\Omega,2+\tilde\kappa-\kappa,l+2,1} +\|f_\theta\|_{\partial_1\Omega,{3\over 2}+\tilde\kappa-\kappa,{5\over 2}+l}.\qedhere
\end{align*}
\end{proof}

\begin{remark}[Comparison with the corner ladder]\label{radial trade and corner ladder}
With $\smash[b]{\mb P_{\mb u_{F_\bullet},\dr_{-1}}}$ as in \eqref{data step relations}, we have
\begin{align*}
\mb P_0(r,\theta-\theta_b)=\mb P_{\mb u_{F_{0,0}},\dr_{-1}}+{1\over 2\theta_b}\mb P_{\mb u_{F_{1,1}},\dr_{-1}},\qquad\qquad\mb P_{\mb u_{F_{1,1}},\dr_{-1}}=2\theta_b\mb P_1.
\end{align*}
Under \eqref{C1 in the linear theory} the two degree one fields of \eqref{data step relations} therefore combine into $\mb u^{\smash[t]{[2]}}=\mb u-\mb u(\mb 0)-{1\over 2}u_1(\mb 0)\mb P_0(r,\theta-\theta_b)$, which is \eqref{u2 and u relationship}; and the coefficient $(\sin\theta_b)^2\smash[b]{c_{\phi,\dr_{-1}^2}}=u_1(\mb 0)/\theta_b$ of $\ln r$ in $p^{\smash[t]{[2]}}-p$ in Proposition \ref{corner ladder} agrees with that in Proposition \ref{time for radial trade}.
\end{remark}

\subsubsection{The full regularity theorem}\label{full regularity theorem}

We obtain two time derivatives of the solution of Theorem \ref{linear solution in time} by applying that theorem to the time derivatives of the data, with the initial profiles $\dot\phi_0$ and $\ddot\phi_0$ of the initial jet of Theorem \ref{initial jet}, and trade them for radial scaling by Proposition \ref{time for radial trade}.

\begin{theorem}[Solution of the linear system, full regularity]\label{linear solution full regularity}
Take the setting of Theorem \ref{linear solution in time}. Let $\tilde\kappa\in[0,\kappa)$ and suppose that $(3+\tilde\kappa-\kappa)\theta_b<\pi$. Suppose that $(\phi_0,\dot\phi_0,\ddot\phi_0)$ is the initial jet of Theorem \ref{initial jet} for the data $(\mb F(0),\mb f(0),\partial_t\mb F(0),\partial_t\mb f(0))$. That is, as in that theorem, $(\mb u_0,p_0)$ is the stationary solution with data $(\phi_0,\mb F(0),\mb f(0))$ satisfying the compatibility condition \eqref{C1 in the linear theory}, $(\dot{\mb u}_0,\dot p_0)$ is the stationary solution with data $(\dot\phi_0,\partial_t\mb F(0),\partial_t\mb f(0))$, and on $\partial_1\Omega$ we have $\dot\phi_0=(\sin\theta_b)^{-2}\partial_ru_{0,\theta}$ and $\ddot\phi_0=(\sin\theta_b)^{-2}\partial_r\dot u_{0,\theta}$. Suppose further that
\begin{align}\label{improved theorem C2}
\|\ddot\phi_0\|_{\partial_1\Omega,{3\over 2}-\kappa,l+{3\over 2}}<\infty,
\end{align}
and that $\mb F$ and $\mb f$ have the regularity that makes the RHS of \eqref{Improved solution estimate} finite. Then $(\mb u,p,\partial_t\mb u,\partial_tp)$ takes the initial value $(\mb u_0,p_0,\dot{\mb u}_0,\dot p_0)$, and
\begin{align}\label{Improved solution estimate}
&(\sin\theta_b)^2\nn\phi\nn_{[0,T],\partial_1\Omega:(\varsigma,{1\over 2}-\kappa)(2,l+{3\over 2})}+\nn\mb u\nn_{[0,T],\Omega:(\varsigma,1-\kappa)(2,l,2)}+\nn p\nn_{[0,T],\Omega:(\varsigma,-\kappa)(2,l,1)}+\nonumber\\
&(\sin\theta_b)^2\nn\phi\nn_{[0,T],\partial_1\Omega:(\varsigma,{3\over 2}+\tilde\kappa-\kappa)(1,l+{5\over 2})}+\nn\mb u^{[1]}\nn_{[0,T],\Omega:(\varsigma,2+\tilde\kappa-\kappa)(1,l+1,2)}+\nn p\nn_{[0,T],\Omega:(\varsigma,1+\tilde\kappa-\kappa)(1,l+1,1)}\nonumber\\
&\lesssim_{\theta_b,\kappa,\tilde\kappa,l}\brac{1+\varsigma^{-1}}\Bigg(\sum_{k=0}^1\brac{\nn\mb F\nn_{[0,T],\Omega:(\varsigma,-\kappa k-1)(2,l)}+\nn\mb f\nn_{[0,T],\partial_1\Omega:(\varsigma,-\kappa k-{1\over 2})(2,l+{1\over 2})}}\nonumber\\
&\quad+\nn\mb F\nn_{[0,T],\Omega:(\varsigma,\varpi-1)(2)}+\nn\mb f\nn_{[0,T],\partial_1\Omega:(\varsigma,\varpi-{1\over 2})(2,{1\over 2})}\nonumber\\
&\quad+\nn\mb F\nn_{[0,T],\Omega:(\varsigma,\tilde\kappa-\kappa)(1,l+1)}+\nn\mb f\nn_{[0,T],\partial_1\Omega:(\varsigma,\tilde\kappa-\kappa+{1\over 2})(1,l+{3\over 2})}\nonumber\\
&\quad+(\sin\theta_b)^2\brac{\|\phi_0\|_{\partial_1\Omega,{1\over 2}-\kappa,l+{7\over 2}}+\|\dot\phi_0\|_{\partial_1\Omega,{1\over 2}-\kappa,l+{5\over 2}}+\|\ddot\phi_0\|_{\partial_1\Omega,{1\over 2}-\kappa,l+{3\over 2}}}\Bigg),
\end{align}
where $\mb u^{\smash[t]{[1]}}$ is as defined in Proposition \ref{time for radial trade}. If in addition the norms $\nn\mb F\nn_{[0,T],\Omega:(\varsigma,\tilde\kappa-\kappa+1)(0,l+2)}$ and $\nn\mb f\nn_{[0,T],\partial_1\Omega:(\varsigma,\tilde\kappa-\kappa+3/2)(0,l+5/2)}$ are finite, then \eqref{Improved solution estimate} holds with
\begin{align}\label{Improved solution M2 norms}
(\sin\theta_b)^2\nn\phi^{[2]}\nn_{[0,T],\partial_1\Omega:(\varsigma,{5\over 2}+\tilde\kappa-\kappa)(0,l+{7\over 2})}+\nn\mb u^{[2]}\nn_{[0,T],\Omega:(\varsigma,3+\tilde\kappa-\kappa)(0,l+2,2)}+\nn p^{[2]}\nn_{[0,T],\Omega:(\varsigma,2+\tilde\kappa-\kappa)(0,l+2,1)}
\end{align}
added to its LHS and $1+\varsigma^{-1}$ times these two norms added to its RHS, with $\phi^{\smash[t]{[2]}}$, $\mb u^{\smash[t]{[2]}}$ and $p^{\smash[t]{[2]}}$ as defined in Proposition \ref{time for radial trade}.
\end{theorem}
\begin{proof}
Step 1 solves the problem and the differentiated problem, and Step 2 shows, by Lemma \ref{time derivative of the solution}, that the solution of the second is the time derivative of the solution of the first. Step 3 does the same one time derivative higher, and Step 4 trades the time derivatives for radial scaling by Proposition \ref{time for radial trade}.

\emph{Step 1: the triples $(\phi,\mb u,p)$ and $(\dot\phi,\dot{\mb u},\dot p)$.} Applying Theorem \ref{linear solution in time} to $(\phi_0,\mb F,\mb f)$ we get $(\mb u,p)$ such that
\begin{align}\label{system: time side undifferentiated}
(\phi,\mb u,p)\text{ solves }\mathfrak W_t[\mb F,\mb f;\phi_0]\text{ a.e., with }\phi(0,t)=0\text{ weakly},
\end{align}
with regularity
\begin{align}
\mathfrak E_{0,l+2}[\phi,\mb u,p]\lesssim_{\theta_b}\mathfrak D_{0,l+2}[\mb F,\mb f;\phi_0],\label{time solution estimate undifferentiated}
\end{align}
where $\mathfrak W_t$ is the system defined after System \ref{reference wedge system}, and $\mathfrak E_{n,j}$ and $\mathfrak D_{n,j}$ are the two sides of \eqref{time solution estimate} with the counts $(0,l)$ replaced by $(n,j)$, as in Theorem \ref{linear solution in time}. That theorem also applies to $(\dot\phi_0,\partial_t\mb F,\partial_t\mb f)$, since $\dot\phi_0(0)=0$ by \eqref{C1 in the linear theory} and Proposition \ref{corner ladder}, and since $\|\dot\phi_0\|_{\partial_1\Omega,1/2-\kappa,l+5/2}$ is supplied by Theorem \ref{initial jet}. It gives $(\dot{\mb u},\dot p)$ such that
\begin{align}\label{system: time side differentiated}
(\dot\phi,\dot{\mb u},\dot p)\text{ solves }\mathfrak W_t[\partial_t\mb F,\partial_t\mb f;\dot\phi_0]\text{ a.e., with }\dot\phi(0,t)=0\text{ weakly},
\end{align}
with regularity
\begin{align}
\mathfrak E_{0,l+1}[\dot\phi,\dot{\mb u},\dot p]\lesssim_{\theta_b}\mathfrak D_{0,l+1}[\partial_t\mb F,\partial_t\mb f;\dot\phi_0].\label{time solution estimate differentiated}
\end{align}
By Theorem \ref{initial jet} there exist $(\mb u_0,p_0)$ such that
\begin{align}\label{system: stationary initial}
(\mb u_0,p_0)\text{ solves }\eqref{system: stationary linear}\text{ with data }(\phi_0,\mb F(0),\mb f(0)),
\end{align}
with regularity
\begin{align*}
&\sum_{M=0}^2\brac{\|\dot\phi_0\|_{\partial_1\Omega,M-{1\over 2}-\kappa,l+{5\over 2}}+\|\mb u_0^{[M]}\|_{\Omega,M+1-\kappa,l+2,2}+\|p_0^{[M]}\|_{\Omega,M-\kappa,l+2,1}}\\
&\quad+\sum_{M=0}^1\brac{\|\ddot\phi_0\|_{\partial_1\Omega,M-{1\over 2}-\kappa,l+{3\over 2}}+\|\dot{\mb u}_0^{[M]}\|_{\Omega,M+1-\kappa,l+1,2}+\|\dot p_0^{[M]}\|_{\Omega,M-\kappa,l+1,1}}\\
&\lesssim_{\theta_b,\kappa}\sum_{M=0}^2\brac{\|\mb F(0)\|_{\Omega,M-1-\kappa,l+2}+\|\mb f(0)\|_{\partial_1\Omega,M-{1\over 2}-\kappa,l+{5\over 2}}+(\sin\theta_b)^2\|\phi_0^{[M]}\|_{\partial_1\Omega,M+{1\over 2}-\kappa,l+{7\over 2}}}\\
&\quad+\sum_{M=0}^1\brac{\|\partial_t\mb F(0)\|_{\Omega,M-1-\kappa,l+1}+\|\partial_t\mb f(0)\|_{\partial_1\Omega,M-{1\over 2}-\kappa,l+{3\over 2}}}.
\end{align*}

\emph{Step 2: the first time derivative.} Lemma \ref{time derivative of the solution} with $l+1$ in place of $l$ applies to the two triples of Step 1. Its hypotheses on $(\mb u_0,p_0)$ are those of Lemma \ref{stationary as variational} at the count $l+1$. The pair $(\mb u_0,p_0)$ is the solution of Proposition \ref{stationary system} at the count $l+2$, so it lies in the class of \eqref{data estimate} with the Navier slip term. It has the corner value and the rung $\|\mb u_0^{\smash[t]{[1]}}\|_{\Omega,2-\kappa,l+2,2}$ of Theorem \ref{initial jet}. The contact angle condition \eqref{stationary corner condition} holds by \eqref{C1 in the linear theory} and Proposition \ref{corner ladder}. The norms of the data are finite by the hypotheses and Theorem \ref{initial jet}. So $(\partial_t\phi,\partial_t\mb u,\partial_tp)=(\dot\phi,\dot{\mb u},\dot p)$ and $(\mb u,p)(0)=(\mb u_0,p_0)$. Therefore adding \eqref{time solution estimate undifferentiated} and \eqref{time solution estimate differentiated} we get
\begin{align}\label{time solution estimate first derivative}
\mathfrak E_{1,l+1}[\phi,\mb u,p]
&\lesssim_{\theta_b}\brac{1+\varsigma^{-1}}\Bigg(\sum_{k=0}^1\brac{\nn\mb F\nn_{[0,T],\Omega:(\varsigma,-\kappa k-1)(1,l+1)}+\nn\mb f\nn_{[0,T],\partial_1\Omega:(\varsigma,-\kappa k-{1\over 2})(1,l+{3\over 2})}}\nonumber\\
&\quad+\nn\mb F\nn_{[0,T],\Omega:(\varsigma,\varpi-1)(2)}+\nn\mb f\nn_{[0,T],\partial_1\Omega:(\varsigma,\varpi-{1\over 2})(2,{1\over 2})}\nonumber\\
&\quad+(\sin\theta_b)^2\brac{\|\phi_0\|_{\partial_1\Omega,{1\over 2}-\kappa,l+{7\over 2}}+\|\dot\phi_0\|_{\partial_1\Omega,{1\over 2}-\kappa,l+{5\over 2}}}\Bigg).
\end{align}

\emph{Step 3: the second time derivative.} We do the same with the second time derivative. Theorem \ref{linear solution in time} applies at the count $l$ to $(\ddot\phi_0,\partial_t^2\mb F,\partial_t^2\mb f)$, since $\ddot\phi_0(0)=0$ by \eqref{improved theorem C2} and Lemma \ref{embeddings}(i), and since $\|\ddot\phi_0\|_{\partial_1\Omega,1/2-\kappa,l+3/2}$ is supplied by Theorem \ref{initial jet}. Let $(\ddot\phi,\ddot{\mb u},\ddot p)$ be its solution. Lemma \ref{time derivative of the solution}, with $(\partial_t\mb F,\partial_t\mb f,\dot\phi_0)$ and $(\dot{\mb u}_0,\dot p_0)$ in place of $(\mb F,\mb f,\phi_0)$ and $(\mb u_0,p_0)$, applies to $(\dot\phi,\dot{\mb u},\dot p)$ and $(\ddot\phi,\ddot{\mb u},\ddot p)$, the derived profile being $\ddot\phi_0$. Its hypotheses on $(\dot{\mb u}_0,\dot p_0)$ are those of Lemma \ref{stationary as variational} at the count $l$. The pair is the solution of Proposition \ref{stationary system} at the count $l+1$, so it lies in the class of \eqref{data estimate} with the Navier slip term. It has the corner value and the rung $\|\dot{\mb u}_0^{\trimtop{0.9pt}{[1]}}\|_{\Omega,2-\kappa,l+1,2}$ of Theorem \ref{initial jet}. Its contact angle condition \eqref{stationary corner condition} is the pointwise vanishing of $\ddot\phi_0$ at the contact point just noted, and its profile datum $s\dot\phi_0-\ddot\phi_0$ has the finite norms required, by Theorem \ref{initial jet} and \eqref{improved theorem C2}. The norms of $\partial_t\mb F(0)$ and $\partial_t\mb f(0)$ are finite by the hypotheses. We get $(\partial_t\mb u,\partial_tp)(0)=(\dot{\mb u},\dot p)(0)=(\dot{\mb u}_0,\dot p_0)$ and $(\partial_t^2\phi,\partial_t^2\mb u,\partial_t^2p)=(\ddot\phi,\ddot{\mb u},\ddot p)$. Therefore adding \eqref{time solution estimate first derivative} and $\mathfrak E_{0,l}[\ddot\phi,\ddot{\mb u},\ddot p]\lesssim_{\theta_b}\mathfrak D_{0,l}[\partial_t^2\mb F,\partial_t^2\mb f;\ddot\phi_0]$ we get
\begin{align}\label{time solution estimate second derivative}
\mathfrak E_{2,l}[\phi,\mb u,p]\nonumber
&\lesssim_{\theta_b}\brac{1+\varsigma^{-1}}\Bigg(\sum_{k=0}^1\brac{\nn\mb F\nn_{[0,T],\Omega:(\varsigma,-\kappa k-1)(2,l)}+\nn\mb f\nn_{[0,T],\partial_1\Omega:(\varsigma,-\kappa k-{1\over 2})(2,l+{1\over 2})}}\nonumber\\
&\quad+\nn\mb F\nn_{[0,T],\Omega:(\varsigma,\varpi-1)(2)}+\nn\mb f\nn_{[0,T],\partial_1\Omega:(\varsigma,\varpi-{1\over 2})(2,{1\over 2})}\nonumber\\
&\quad+(\sin\theta_b)^2\brac{\|\phi_0\|_{\partial_1\Omega,{1\over 2}-\kappa,l+{7\over 2}}+\|\dot\phi_0\|_{\partial_1\Omega,{1\over 2}-\kappa,l+{5\over 2}}+\|\ddot\phi_0\|_{\partial_1\Omega,{1\over 2}-\kappa,l+{3\over 2}}}\Bigg).
\end{align}

\emph{Step 4: trading the time derivatives for radial scaling.} The RHSs of \eqref{time for radial trade estimate 1}, applied to $(\phi,\mb u,p)$ and, with $l-1$ in place of $l$, to $(\partial_t\phi,\partial_t\mb u,\partial_tp)$, contain the norms on the left below. Since $-\kappa\leq\tilde\kappa-\kappa<0$, we have, using Lemma \ref{limit at infinity} for the terms of $\mb u$ and $\partial_t\mb u$ without derivatives,
\begin{align*}
\nn\mb u\nn_{[0,T],\Omega:(\varsigma,1+\tilde\kappa-\kappa)(1,l+1,2)}+\nn p\nn_{[0,T],\Omega:(\varsigma,\tilde\kappa-\kappa)(1,l+1,1)}&\lesssim_{\kappa,\tilde\kappa}\mathfrak E_{1,l+1}[\phi,\mb u,p],\\
\nn\partial_t\mb u\nn_{[0,T],\Omega:(\varsigma,1+\tilde\kappa-\kappa)(1,l,2)}+\nn\partial_tp\nn_{[0,T],\Omega:(\varsigma,\tilde\kappa-\kappa)(1,l,1)}&\lesssim_{\kappa,\tilde\kappa}\mathfrak E_{2,l}[\phi,\mb u,p],\\
\nn\partial_tf_\theta\nn_{[0,T],\partial_1\Omega:(\varsigma,\tilde\kappa-\kappa-{1\over 2})(1,l+{1\over 2})}&\lesssim_{\kappa,l}\sum_{k=0}^1\nn\mb f\nn_{[0,T],\partial_1\Omega:(\varsigma,-\kappa k-{1\over 2})(2,l+{1\over 2})}.
\end{align*}
The hypotheses of Proposition \ref{time for radial trade} on $\partial_t\mb u$ and $\partial_tp$ are contained in $\mathfrak E_{1,l+1}[\phi,\mb u,p]$ and $\mathfrak E_{2,l}[\phi,\mb u,p]$, and those on $\mb F$ and $\mb f$ are bounded by the RHS of \eqref{Improved solution estimate}, using the last line above for $\partial_tf_\theta$. So \eqref{time for radial trade estimate 1}, under the hypothesis $(3+\tilde\kappa-\kappa)\theta_b<\pi$ of that proposition, together with \eqref{time solution estimate first derivative} and \eqref{time solution estimate second derivative}, gives \eqref{Improved solution estimate}. When the two further norms of the statement are finite, \eqref{time for radial trade estimate 2}, applied to $(\phi,\mb u,p)$, gives the norms \eqref{Improved solution M2 norms}, its hypotheses on $\partial_t\mb u^{[1]}$ and $\partial_tp$ being contained in $\nn\mb u^{[1]}\nn_{[0,T],\Omega:(\varsigma,2+\tilde\kappa-\kappa)(1,l+1,2)}$ and $\nn p\nn_{[0,T],\Omega:(\varsigma,1+\tilde\kappa-\kappa)(1,l+1,1)}$ on the left of \eqref{Improved solution estimate}. The constants of Proposition \ref{time for radial trade} are stationary, so the factor $1+\varsigma^{-1}$ of the three applications of Theorem \ref{linear solution in time} is the only dependence on $\varsigma$.
\end{proof}

The implicit constant in \eqref{Improved solution estimate} degenerates as $\tilde\kappa\uparrow\kappa$, as explained in the proof of Proposition \ref{time for radial trade}.

\begin{remark}\label{linear theory closing remark}
We have left $\|\dot\phi_0\|_{\partial_1\Omega,1/2-\kappa,l+5/2}+\|\ddot\phi_0\|_{\partial_1\Omega,1/2-\kappa,l+3/2}$ on the RHS of \eqref{Improved solution estimate}, rather than converting them back into norms of $\phi_0$ using the estimate of Theorem \ref{initial jet}. Doing so would also bring norms such as $\|\partial_t\mb F(0)\|_{\Omega,-\kappa,l+1}$ onto the RHS, and these cannot be absorbed into $\nn\mb F\nn$. Indeed, counting $\partial_t$ as $\partial_r$, a norm $\nn u\nn_{[0,T],O:(\varsigma,a)(1,j)}$ is formally $\nn u\nn_{[0,T],O:(\varsigma,a+1/2)(1/2,j+1/2)}$, which at $t=0$ controls at best $\|\ph\|_{O,a+1/2,j+1/2}$. On the RHS, $\partial_t^2\mb F$ occurs at the count $l$ only in $\nn\mb F\nn_{[0,T],\Omega:(\varsigma,-\kappa k-1)(2,l)}$, $k=0,1$. At the higher weight, $k=0$, the norm $\nn\partial_t\mb F\nn_{[0,T],\Omega:(\varsigma,-1)(1,l)}$ contained in $\nn\mb F\nn_{[0,T],\Omega:(\varsigma,-1)(2,l)}$ controls at best $\|\ph\|_{\Omega,-1/2,l+1/2}$, which falls short of $\|\partial_t\mb F(0)\|_{\Omega,-\kappa,l+1}$ by $1/2$ in the count and, since $\kappa<1/2$, by $1/2-\kappa$ in the weight. Likewise, converting $\|\ddot\phi_0\|_{\partial_1\Omega,1/2-\kappa,l+3/2}$ brings $\|\phi_0^{\smash[t]{[2]}}\|_{\partial_1\Omega,5/2-\kappa,l+7/2}$, whereas on the LHS the norm of $\phi$ with a time derivative at the highest weight, $\nn\phi\nn_{[0,T],\partial_1\Omega:(\varsigma,3/2+\tilde\kappa-\kappa)(1,l+5/2)}$, controls at best $\|\ph\|_{\partial_1\Omega,2+\tilde\kappa-\kappa,l+3}$, which falls short of $\|\phi_0^{\smash[t]{[2]}}\|_{\partial_1\Omega,5/2-\kappa,l+7/2}$ by $1/2$ in the count and, since $\tilde\kappa<\kappa<1/2$, by $1/2-\tilde\kappa$ in the weight. So the initial data is not in the trace space of the positive-time estimate. Radial smoothing supplies the missing count, as in Step 1 of the proof of Theorem \ref{nonlinear positive-time existence}, but not the missing weight, since \eqref{smoothing gain} gains counts at a fixed weight. Consequently, the nonlinear problem is solved in two steps: the initial jet is constructed from $\phi_0$ by a nonlinear problem at the initial time, and then prescribed in a nonlinear problem for positive time. This order is also forced by the hypothesis of Theorem \ref{linear solution full regularity} that the jet is compatible with the values of the forcing at $t=0$, which for the nonlinear problem depend on the jet itself; see Remark \ref{positive-time closing remark}.
\end{remark}

\section{Existence for the full nonlinear system}\label{nonlinear existence}

We can now solve System \ref{reference wedge system} in full by a fixed point argument in the norms supplied by the linear theory of the last two sections. As explained in Remark \ref{linear theory closing remark}, the initial data of our problem are not in the trace space of the positive-time estimate, and the jet must be compatible with the values of the forcing at $t=0$, which depend on the jet. So the argument splits into two nonlinear problems. We first construct the jet from $\phi_0$ by solving the nonlinear problem at the initial time, in Subsection \ref{initial-time construction}, and then solve for positive time with this jet prescribed, in Subsection \ref{positive-time construction}. What makes the fixed point work is that the nonlinearities $\mb N$ and $\mb N_\partial$ are at least quadratic, with at least one factor in every summand coming from the flattening geometry $\bs\nu$, which is of the size $\mathcal D$ of the data. So for $\mathcal D$ small the map that sends a pair $(\mb v,q)$ to the solution of the linear system with forcing $(\mb N,\mb N_\partial)$ evaluated at that pair is a contraction. The weights need care: the linear theory is a ladder of estimates indexed by $M\in\{0,1,2\}$, of which the top rung is only available after the explicit corner fields of \eqref{data step relations} have been subtracted. Lemma \ref{decomposition estimates} and Remark \ref{metric remark} keep track of those subtractions along the iteration.

\subsection{Construction of the initial-time solution}\label{initial-time construction}

In this subsection we produce, from an initial profile $\phi_0$ alone, the stationary pairs and derived profiles $(\mb u,p,\dot{\mb u},\dot p,\dot\phi,\ddot\phi)$ at $t=0$, of size controlled by that of $\phi_0$, and we identify which condition on $\phi_0$ is needed for the contact angle condition to hold also after one time differentiation. The output is the compatible triple $(\phi_0,\dot\phi,\ddot\phi)$ that Theorem \ref{linear solution full regularity} takes as input.

Throughout this subsection we fix $\kappa\in(0,1/2]$ and an integer $l\geq 1$, and we work under the standing angle condition
\begin{align}\label{standing angle condition}
2-\kappa<\Lambda^*(\theta_b)\qquad\text{ and }\qquad\brac{2-{\kappa\over 2}}\theta_b<{\pi\over 2},
\end{align}
where $\Lambda^*(\theta_b)$ is as in Definition \ref{zero-free strip}; the two conditions are explained in Remark \ref{angle remark}. We also assume that $\theta_b$ is small enough for Proposition \ref{stationary system} to apply at the derivative parameters $l+2$ and $l+1$, which is the unquantified smallness that the linear theory of Section \ref{linear existence} requires.
We are given an initial profile $\phi_0$ with $\phi_0(0)=0$ and $\lim_{r\to\infty}\phi_0=0$, and we write $\phi_0^{\smash[t]{[M]}}$, $\smash[b]{c_{\phi_0,\dr_{-1}}}$ and $\smash[b]{c_{\phi_0,\dr_{-1}^2}}$ for the decomposition \eqref{data step relations} applied to $\phi_0$. Every object below lives at the initial time $t=0$ and we suppress the time argument throughout, writing $(\mb u,p)$ for $(\mb u(0),p(0))$, $\dot\phi$ for $\partial_t\phi|_{t=0}$, and so on. All implicit constants are allowed to depend on $\theta_b,\kappa$ and $l$, so that $\lesssim_{\theta_b}$ stands for $\lesssim_{\theta_b,\kappa,l}$; in particular all positive and negative powers of $\sin\theta_b$ and of $\theta_b$, and all resonance factors $|\cos(a\theta_b)|^{-1}$, $\Gamma_{a}$, $\Gamma_{a,\theta_b}$ and $(\Gamma^H_{a,\theta_b})^{-1}$ arising at the admissible weights $a$, are absorbed into $\lesssim_{\theta_b}$.

\begin{definition}[Data, profile and solution norms]\label{data profile and solution norms}
We measure the size of a profile $\psi$, decomposed as in \eqref{data step relations}, by
\begin{align}\label{data norm}
\mathcal D(\psi):=\sum_{M=0}^2(\sin\theta_b)^2\|\psi^{[M]}\|_{\partial_1\Omega,M+{1\over 2}-\kappa,l+{7\over 2}},\qquad\mathcal D:=\mathcal D(\phi_0),
\end{align}
the size of a profile $\dot\psi$ at the weights of the derived profile by
\begin{align}\label{profile norm}
\dot{\mathcal D}(\dot\psi):=\sum_{M=0}^2(\sin\theta_b)^2\|\dot\psi\|_{\partial_1\Omega,M-{1\over 2}-\kappa,l+{5\over 2}},
\end{align}
and the size of a pair $(\mb v,q)$, respectively $(\dot{\mb v},\dot q)$, by
\begin{subequations}\label{solution norms}
\begin{align}
\mathcal E(\mb v,q)&:=\sum_{M=0}^2\brac{\|\mb v^{[M]}\|_{\Omega,M+1-\kappa,l+2,2}+\|q^{[M]}\|_{\Omega,M-\kappa,l+2,1}},\\
\dot{\mathcal E}(\dot{\mb v},\dot q)&:=\sum_{M=0}^1\brac{\|\dot{\mb v}^{[M]}\|_{\Omega,M+1-\kappa,l+1,2}+\|\dot q^{[M]}\|_{\Omega,M-\kappa,l+1,1}},
\end{align}
\end{subequations}
where the decompositions are those of \eqref{data step relations}, the profile occurring there being frozen at $\phi_0$, so that the two constants $\smash[b]{c_{\phi,\dr_{-1}}}$ and $\smash[b]{c_{\phi,\dr_{-1}^2}}$ are frozen at $\smash[b]{c_{\phi_0,\dr_{-1}}}$ and $\smash[b]{c_{\phi_0,\dr_{-1}^2}}$. On pairs we use the distance
\begin{align}\label{contraction metric}
d\brac{(\mb v,q),(\mb v',q')}:=\sum_{M=0}^2\brac{\|\mb v^{[M]}-\mb v'^{[M]}\|_{\Omega,M+1-\kappa,l+2,2}+\|q^{[M]}-q'^{[M]}\|_{\Omega,M-\kappa,l+2,1}}.
\end{align}
\end{definition}

\begin{remark}[The distance of the iteration]\label{metric remark}
The terms of \eqref{data step relations} carrying the frozen constants do not depend on the pair, so they cancel in $d$: the differences $\mb v^{[M]}-\mb v'^{[M]}$ and $q^{[M]}-q'^{[M]}$ are the pieces of $(\mb v-\mb v',q-q')$ computed with profile datum $0$, which is the profile datum carried by the difference of two solutions of \eqref{system: stationary linear}. Thus $d$ is the translation-invariant metric attached to the norm obtained from $\mathcal E$ by setting the frozen constants to $0$, and it is with respect to $d$, and not with respect to $\mathcal E$, that the contraction argument of Theorem \ref{nonlinear initial data} is run. This care is necessary because $\mathcal E$ is not a norm: the maps $\mb v\mapsto\mb v^{[2]}$ and $q\mapsto q^{[2]}$ are affine rather than linear, and unless $\smash[b]{c_{\phi_0,\dr_{-1}^2}}=0$ we have $\mathcal E(\mb 0,0)=\infty$, since the degree one field $\smash[b]{\mb P_{\mb u_{F_{1,1}},\dr_{-1}}}$ is of size $r$ and hence has infinite norm $\|\ph\|_{\Omega,3-\kappa}$, the divergence being at the origin. No such care is needed at the dotted level, where only $M\leq 1$ occurs: there $\dot{\mb v}^{[0]}=\dot{\mb v}$, $\dot{\mb v}^{[1]}=\dot{\mb v}-\dot{\mb v}(\mb 0)$ and $\dot q^{[0]}=\dot q^{[1]}=\dot q$ are linear in the pair, and $\dot{\mathcal E}$ is a norm.
\end{remark}

\begin{remark}[On the angle condition]\label{angle remark}
Two angle conditions are used below, and \eqref{standing angle condition} collects them.

The first, $2-\kappa<\Lambda^*(\theta_b)$, is the hypothesis of the linear ladder of Proposition \ref{corner ladder}, which Step 1 of the proof of Theorem \ref{nonlinear initial data} climbs up to $M=2$. Under it the Green's function of Lemma \ref{Green's function initial time} is controlled on the closed strip $-\kappa\leq\Re\lambda\leq 2-\kappa$ that the ladder crosses, and the factor $(\Gamma^H_{2-\kappa,\theta_b})^{-1}$ in \eqref{um data reg master} blows up as $2-\kappa$ approaches $\Lambda^*(\theta_b)$. By Remark \ref{H numerics remark} this condition holds for $\theta_b$ up to about $73^\circ$--$75^\circ$ when $\kappa\in(1/4,1/2)$.

The second, $(2-\kappa/2)\theta_b<\pi/2$, is the admissibility of the extension estimates of Lemma \ref{mixed harmonic extension} at every weight used below. All our applications of that lemma occur at its weight parameter $a+1$ with $a$ in the range \eqref{margin ranges}, so that $a\leq 2-\kappa/2$. This condition is the more restrictive: it limits the angle to $\theta_b<\pi/(4-\kappa)$, that is, to about $48^\circ$--$51^\circ$.
\end{remark}

\begin{lemma}[Decomposition estimates]\label{decomposition estimates}
With $\kappa$ and $l$ as fixed at the beginning of Subsection \ref{initial-time construction}, let $c\in\R$ and let $(\mb v,q)$ be a pair, with $\mb v^{[M]}$ given by \eqref{data step relations} with $c$ in place of $\smash[b]{c_{\phi,\dr_{-1}^2}}$.
\begin{enumerate}
\item Let $\Omega_{1,2}=\Omega\cap\{1<r<2\}$. Every linear combination of finitely many fixed functions in $L^2(\Omega_{1,2})$ that are linearly independent there has coefficients bounded in absolute value by a constant, depending only on these functions, times the $L^2(\Omega_{1,2})$ norm of the combination. Every function $g$ on $\Omega$ satisfies $\|g\|_{L^2(\Omega_{1,2})}\lesssim_a\|g\|_{\Omega,a}$ for every $a\in\R$. The same holds on $\partial_1\Omega\cap\{1<r<2\}$, with $\|\ph\|_{\partial_1\Omega,a}$ in place of $\|\ph\|_{\Omega,a}$.
\item Whenever the RHSs are finite, we have
\begin{align}
|\mb v(\mb 0)|&\lesssim_{\theta_b}\|\mb v\|_{\Omega,1-\kappa}+\|\mb v^{[1]}\|_{\Omega,2-\kappa},\label{corner coefficients 1}\\
|v_1(\mb 0)|+(\sin\theta_b)^2|c|&\lesssim_{\theta_b}\|\mb v^{[1]}\|_{\Omega,2-\kappa}+\|\mb v^{[2]}\|_{\Omega,3-\kappa};\label{corner coefficients 2}
\end{align}
in particular, for $c=\smash[b]{c_{\phi_0,\dr_{-1}^2}}$, both LHSs are $\lesssim_{\theta_b}\mathcal E(\mb v,q)$. Similarly, for the profile,
\begin{align}\label{profile coefficients}
(\sin\theta_b)^2\brac{|c_{\phi_0,\dr_{-1}}|+|c_{\phi_0,\dr_{-1}^2}|}\lesssim_{\theta_b}\mathcal D.
\end{align}
\item If $v_\theta|_{\partial_0\Omega}=0$ and $\|\mb v^{\smash[t]{[1]}}\|_{\Omega,2-\kappa,1,1}<\infty$, then $\mb v$ is continuous at $\mb 0$ and $v_y(\mb 0)=0$, that is, $\mb v(\mb 0)=v_1(\mb 0)\mb e_1$.
\end{enumerate}
\end{lemma}
\begin{proof}
(i) The Gram matrix of the fixed functions in $L^2(\Omega_{1,2})$ is invertible, which gives the bound on the coefficients. On $\Omega_{1,2}$ all the weights $r^{-a}$ are comparable to $1$, which gives the second bound. The same argument applies on $\partial_1\Omega\cap\{1<r<2\}$.

(ii) The constant $\mb v(\mb 0)=\mb v-\mb v^{[1]}$ is a combination of $\mb e_1$ and $\mb e_2$, so (i) gives \eqref{corner coefficients 1}. For \eqref{corner coefficients 2} we apply (i) to
\begin{align*}
\mb v^{[1]}-\mb v^{[2]}={v_1(\mb 0)\over 2}\mb P_{\mb u_{F_{0,0}},\dr_{-1}}+{(\sin\theta_b)^2c\over 4}\mb P_{\mb u_{F_{1,1}},\dr_{-1}}:
\end{align*}
the two degree one fields are linearly independent in $L^2(\Omega_{1,2})$, since their $\mb e_r$ components $r\cos(2(\theta_b-\theta))/\sin(2\theta_b)$ and $-r$ are not proportional as functions of $\theta$. Finally (i) on $\partial_1\Omega\cap\{1<r<2\}$, applied to $\phi_0-\phi_0^{[2]}=\smash[b]{c_{\phi_0,\dr_{-1}}}r+\smash[b]{c_{\phi_0,\dr_{-1}^2}}r\ln r$, gives \eqref{profile coefficients}.

(iii) By Lemma \ref{continuity at the corner} at $a=1-\kappa$, the function $\mb v^{[1]}=\mb v-\mb v(\mb 0)$ extends continuously to $\overline\Omega$ and vanishes at $\mb 0$, so $\mb v$ is continuous at $\mb 0$. Since $\mb e_\theta=\mb e_2$ on $\partial_0\Omega$, the component $v_y=v_\theta$ vanishes there, and hence at $\mb 0$.
\end{proof}

\subsubsection{Estimates for the flattening geometry}

Here and below $\psi$ and $\dot\psi$ denote arbitrary profiles, that is, functions on $\partial_1\Omega$, whereas $\phi_0$ is the initial profile, which carries the data norm $\mathcal D$. Recall from Lemma \ref{Piola transform} that
\begin{align}\label{flattening geometry}
\bs\nu[\psi]=-(\sin\theta_b)\grad\E\dr^{-1}(r\psi)=\grad\E\chi,\qquad\text{ where }\qquad\chi=-(\sin\theta_b)\dr^{-1}(r\psi),
\end{align}
and set
\begin{align*}
\mu_1:=(1+\nu_1)^{-1}-1=-{\nu_1\over 1+\nu_1}.
\end{align*}
The factor $\sin\theta_b$ is the one converting $\partial_y$ into $\partial_r$ along $\partial_1\Omega$, where $y=r\sin\theta_b$. Along the ray $\partial_1\Omega$ the unit vector $\mb e_r=(\cos\theta_b,\sin\theta_b)$ is constant and $\partial_r$ is tangential, hence commutes with restriction to the ray; since $\E\chi|_{\partial_1\Omega}=\chi$, $\dr\chi=-(\sin\theta_b)r\psi$ and $\cot\theta_b=b^{-1}$, we obtain
\begin{align}\label{profile trace}
\eva{(b^{-1}\nu_1+\nu_2)}_{\partial_1\Omega}=\eva{(\sin\theta_b)^{-1}\nu_r}_{\partial_1\Omega}=(\sin\theta_b)^{-1}\eva{\mb e_r\cdot\grad\E\chi}_{\partial_1\Omega}=(\sin\theta_b)^{-1}\partial_r\chi&=-\psi,\nonumber\\
\eva{r^{-1}\dr\brac{b^{-1}\nu_1+\nu_2}}_{\partial_1\Omega}=\eva{(\sin\theta_b)^{-1}r^{-1}\dr\nu_r}_{\partial_1\Omega}&=-\partial_r\psi,
\end{align}
the second line following from the first by applying $\partial_r$.
Thus $\psi$ and $\partial_r\psi$ are traces of components of $\bs\nu$ and of $r^{-1}\dr\bs\nu$, so that the boundary nonlinearity is the trace of a bulk expression, \eqref{schematic boundary}, estimated in Case 4 of the proof of Proposition \ref{nonlinear estimates}.

\begin{lemma}[Geometry estimates]\label{geometry estimates}
With $\kappa$ and $l$ as fixed at the beginning of Subsection \ref{initial-time construction}, let $\psi$ be a function on $\partial_1\Omega$ and $\bs\nu=\bs\nu[\psi]$ as in \eqref{flattening geometry}.
\begin{enumerate}
\item Suppose $|a\theta_b|<\pi/2$. Then for every $m\in\N_0$,
\begin{align}
\|\bs\nu\|_{\Omega,a,0,m}&\lesssim_{\theta_b,m}\|\psi\|_{\partial_1\Omega,a-{1\over 2},\max\set{m-{1\over 2},0}},\label{geometry bulk}\\
\|r^{1-a}(r\grad)^{m}\bs\nu\|_{L^\infty(\Omega)}&\lesssim_{\theta_b,m}\|\psi\|_{\partial_1\Omega,a-{1\over 2},m+{3\over 2}}.\label{geometry sup}
\end{align}
\item Under \eqref{standing angle condition}, for $\psi=\phi_0$ and
\begin{align}\label{margin ranges}
a\in\sbrac{1-\kappa,2-{\kappa\over 2}},
\end{align}
the hypothesis of (i) holds, and the RHSs of \eqref{geometry bulk} and \eqref{geometry sup} are $\lesssim_{\theta_b}\mathcal D$ for $m\leq l+4$ and $m\leq l+2$ respectively.
\item Under \eqref{standing angle condition}, for $\psi=\dot\psi$ with $\dot{\mathcal D}=\dot{\mathcal D}(\dot\psi)<\infty$ and $a\in[1-\kappa,2-\kappa]$, the RHSs of \eqref{geometry bulk} and \eqref{geometry sup} are $\lesssim_{\theta_b}\dot{\mathcal D}$ for $m\leq l+3$ and $m\leq l+1$ respectively.
\item There is $\delta_0=\delta_0(\theta_b,\kappa,l)>0$ such that $\mathcal D\leq\delta_0$ implies
\begin{align}\label{geometry smallness}
\|\bs\nu[\phi_0]\|_{L^\infty(\Omega)}+\|\phi_0\|_{L^\infty(\partial_1\Omega)}\leq{1\over 2}.
\end{align}
\end{enumerate}
\end{lemma}
\begin{proof}
(i) Since $\bs\nu=-(\sin\theta_b)\grad\E\dr^{-1}(r\psi)$ by \eqref{flattening geometry}, the estimate of Lemma \ref{mixed harmonic extension} for $r\grad\E\phi$, applied at the weight $a+1$ and with $\phi=\dr^{-1}(r\psi)$, so that $\dr\phi=r\psi$, gives
\begin{align*}
\|\bs\nu\|_{\Omega,a,0,m}&\lesssim_m|\sin\theta_b|\|r\grad\E\dr^{-1}(r\psi)\|_{\Omega,a+1,0,m}\lesssim_m{|\sin\theta_b|\over|\cos(a\theta_b)|}\|r\psi\|_{\partial_1\Omega,a+{1\over 2},\max\{m-{1\over 2},0\}}\\
&\lesssim_{\theta_b,m}\|\psi\|_{\partial_1\Omega,a-{1\over 2},\max\{m-{1\over 2},0\}}.
\end{align*}
This proves \eqref{geometry bulk}. Combining \eqref{geometry bulk} with Lemma \ref{embeddings}(ii) we then get
\begin{align*}
\|r^{1-a}(r\grad)^{m}\bs\nu\|_{L^\infty(\Omega)}\lesssim_{\theta_b}\|(r\grad)^{m}\bs\nu\|_{\Omega,a,1,1}&\lesssim_{\theta_b}\|\psi\|_{\partial_1\Omega,a-{1\over 2},m+{3\over 2}}.
\end{align*}

(ii) The hypothesis of (i) holds on \eqref{margin ranges}, since $|a|\theta_b\leq(2-\kappa/2)\theta_b<\pi/2$ by \eqref{standing angle condition}. For $a\in[1-\kappa,2-\kappa]$ the claim follows from $\|\phi_0\|_{\partial_1\Omega,1/2-\kappa,l+7/2}$ and $\|\phi_0\|_{\partial_1\Omega,3/2-\kappa,l+7/2}$ in $\mathcal D$, the counts being $\max\{m-1/2,0\}\leq l+7/2$ and $m+3/2\leq l+7/2$, which are the two stated bounds on $m$. For $a\in(2-\kappa,2-\kappa/2]$ it suffices to show that $\|\phi_0\|_{\partial_1\Omega,a-1/2,k}\lesssim_{\theta_b}\mathcal D$ for $k\leq l+7/2$. We split at $r=1$. On $\{r>1\}$ the norm $\|\phi_0\|_{\partial_1\Omega,1/2-\kappa,l+7/2}$ dominates, since $a-1/2\geq 1/2-\kappa$. On $\{r<1\}$ we write
\[\phi_0=\phi_0^{[2]}+c_{\phi_0,\dr_{-1}}r+c_{\phi_0,\dr_{-1}^2}r\ln r\]
as in \eqref{data step relations}: the norm $\|\phi_0^{\smash[t]{[2]}}\|_{\partial_1\Omega,5/2-\kappa,l+7/2}$ in $\mathcal D$ dominates $\|\phi_0^{\smash[t]{[2]}}\|_{\partial_1\Omega,a-1/2,k}$ there, since $a-1/2\leq 5/2-\kappa$, while $\dr^{j'}r=r$ and $\dr^{j'}(r\ln r)=r\ln r+j'r$ satisfy
\begin{align*}
\int_0^1r^{1-2a}\brac{r^2+r^2\ln^2r}\d r\leq\int_0^1r^{\kappa-1}\brac{1+\ln^2r}\d r<\infty,
\end{align*}
since $3-2a\geq\kappa-1$, their coefficients being $\lesssim_{\theta_b}\mathcal D$ by \eqref{profile coefficients}. This is the only point at which the weight of $\|\phi_0\|_{\partial_1\Omega,3/2-\kappa,l+7/2}$ is crossed, and it is crossed by at most $\kappa/2$.

(iii) For $a\in[1-\kappa,2-\kappa]$ the claim follows from $\|\dot\psi\|_{\partial_1\Omega,1/2-\kappa,l+5/2}$ and $\|\dot\psi\|_{\partial_1\Omega,3/2-\kappa,l+5/2}$ in $\dot{\mathcal D}$, the counts being $\max\{m-1/2,0\}\leq l+5/2$ and $m+3/2\leq l+5/2$.

(iv) By \eqref{geometry sup} at $a=1$ and $m=0$ and by Lemma \ref{embeddings}(i) at $a=0$,
\begin{align*}
\|\bs\nu[\phi_0]\|_{L^\infty(\Omega)}+\|\phi_0\|_{L^\infty(\partial_1\Omega)}\lesssim_{\theta_b}\|\phi_0\|_{\partial_1\Omega,{1\over 2},{3\over 2}}\lesssim_{\theta_b}\mathcal D,
\end{align*}
so it suffices to take $\delta_0$ small enough.
\end{proof}

\begin{lemma}[Composition]\label{composition of the geometry}
Let $\bs\nu=\bs\nu[\psi]$ with $\|\bs\nu\|_{L^\infty(\Omega)}\leq 1/2$. We write $\mathfrak G$ for any function $\mathfrak g\circ\bs\nu$ with $\mathfrak g\in C^\infty([-1/2,1/2]^2)$, possibly different at each occurrence, and $\mathfrak G_0$ for one with in addition $\mathfrak g(\mb 0)=0$.
\begin{enumerate}
\item We have $\|\mathfrak G\|_{L^\infty(\Omega)}\lesssim_{\mathfrak g}1$ and
\begin{align}\label{composition}
\mathfrak G_0=\nu_1\mathfrak G+\nu_2\mathfrak G,\qquad\qquad (r\grad)\mathfrak G=\brac{(r\grad)\nu_1}\mathfrak G+\brac{(r\grad)\nu_2}\mathfrak G.
\end{align}
\item For $\psi=\phi_0$ with $\mathcal D\leq 1$ and $a$ in the range \eqref{margin ranges},
\begin{align}\label{composition bounds}
\|\mathfrak G_0\|_{\Omega,a,0,m}+\|r^{1-a}(r\grad)^{m'}\mathfrak G_0\|_{L^\infty(\Omega)}\lesssim_{\theta_b,\mathfrak g}\mathcal D\qquad\text{ for }\qquad m\leq l+4,\quad m'\leq l+2.
\end{align}
\item The linearisation of a factor $\mathfrak G=\mathfrak g\circ\bs\nu$ at $\bs\nu=\bs\nu[\phi_0]$ in the direction $\bs\nu[\dot\psi]$ is $\nu_1[\dot\psi](\partial_1\mathfrak g)(\bs\nu)+\nu_2[\dot\psi](\partial_2\mathfrak g)(\bs\nu)$. If $\mathcal D\leq 1$ and $a\in[1-\kappa,2-\kappa]$, it obeys the bounds \eqref{composition bounds} with $\mathcal D$ replaced by $\dot{\mathcal D}=\dot{\mathcal D}(\dot\psi)$, for $m\leq l+3$ and $m'\leq l+1$.
\item Let $\bs\nu'=\bs\nu[\psi']$ be a second geometry with $\|\bs\nu'\|_{L^\infty(\Omega)}\leq 1/2$, and write $\tilde{\mathfrak G}$ for any function $\tilde{\mathfrak g}\circ(\bs\nu,\bs\nu')$ with $\tilde{\mathfrak g}\in C^\infty([-1/2,1/2]^4)$, possibly different at each occurrence. Then $\mathfrak g(\bs\nu)-\mathfrak g(\bs\nu')=(\nu_1-\nu_1')\tilde{\mathfrak G}+(\nu_2-\nu_2')\tilde{\mathfrak G}$, and (i) holds for the factors $\tilde{\mathfrak G}$, with $(r\grad)\tilde{\mathfrak G}$ a sum of products of $\tilde{\mathfrak G}$ with $(r\grad)\nu_1$, $(r\grad)\nu_2$, $(r\grad)\nu_1'$ and $(r\grad)\nu_2'$.
\end{enumerate}
\end{lemma}
\begin{proof}
(i) The hypothesis gives $|\nu_1|,|\nu_2|\leq|\bs\nu|\leq 1/2$, so that $\bs\nu$ takes values in the square, on which $\mathfrak g$ and all its derivatives are bounded; this is the sup norm bound. Since $\mathfrak g(\mb 0)=0$ and the square is convex, Hadamard's lemma gives $\mathfrak g(\mb t)=t_1\mathfrak g_1(\mb t)+t_2\mathfrak g_2(\mb t)$ with $\mathfrak g_1,\mathfrak g_2$ smooth on the square; this is the first identity of \eqref{composition}. The second is the chain rule, $\partial_1\mathfrak g$ and $\partial_2\mathfrak g$ being again smooth on the square.

(ii) Iterating the second identity of \eqref{composition} and expanding the first by the Leibniz rule, we find that $(r\grad)^{m}\mathfrak G_0$ is a finite sum of terms $\mathfrak Gf_1\cdots f_J$, where $J\geq 1$ and each $f_{j'}$ is a component of $(r\grad)^{m_{j'}}\bs\nu$, with $m_1+\dots+m_J=m$. We order the factors so that $m_1$ is the largest count. Then $m_{j'}\leq m/2$ for $j'\geq 2$, since $m_1+m_{j'}\leq m$. We place $f_1$ in the weighted norm and the other factors in the sup norm \eqref{geometry sup} at $a=1$, where the weight is trivial:
\begin{align*}
\|\mathfrak Gf_1\cdots f_J\|_{\Omega,a}\leq\|\mathfrak G\|_{L^\infty(\Omega)}\|(r\grad)^{m_1}\bs\nu\|_{\Omega,a}\prod_{j'=2}^J\|(r\grad)^{m_{j'}}\bs\nu\|_{L^\infty(\Omega)}\lesssim_{\theta_b,\mathfrak g}\mathcal D^J\leq\mathcal D,
\end{align*}
by \eqref{geometry bulk} with the count $m_1\leq m\leq l+4$, by \eqref{geometry sup} with the counts $m_{j'}\leq(l+4)/2\leq l+2$, both within the ranges of Lemma \ref{geometry estimates}(ii), and since $\mathcal D\leq 1$. The sup bound in \eqref{composition bounds} follows in the same way, with $f_1$ in the weighted sup norm \eqref{geometry sup} at the weight $a$, and with all counts at most $m'\leq l+2$.

(iii) The formula is the chain rule. Applying $(r\grad)^m$ to it gives a finite sum of terms of the form in (ii) with exactly one factor built from $\dot\psi$. They are placed as in (ii); the counts asked of the factor built from $\dot\psi$, at most $m\leq l+3$ in the weighted norm and at most $(l+3)/2\leq l+1$ in the sup norm, and at most $m'\leq l+1$ for the sup bound, are those of Lemma \ref{geometry estimates}(iii) at the weights $a\in[1-\kappa,2-\kappa]$. No smallness of $\dot{\mathcal D}$ is involved, the factor $\bs\nu[\dot\psi]$ entering linearly.

(iv) Hadamard's lemma along the segment from $\bs\nu'$ to $\bs\nu$, which stays in the square by convexity, gives the formula for $\mathfrak g(\bs\nu)-\mathfrak g(\bs\nu')$, with $\tilde{\mathfrak g}$ smooth on $[-1/2,1/2]^4$. The rest is the proof of (i), with the square replaced by $[-1/2,1/2]^4$.
\end{proof}

\begin{lemma}[Weighted products]\label{weighted products}
Let $j\in\N_0$ and $a_1,a_2\in\R$. Then
\begin{align}\label{bulk product}
\|fg\|_{\Omega,a_1+a_2,j}\lesssim_{\theta_b,j}\|f\|_{\Omega,a_1+1,j,2}\|g\|_{\Omega,a_2,j}.
\end{align}
\end{lemma}
\begin{proof}
We expand $\dr^{j'}(fg)$ by the Leibniz rule and bound each term by
\begin{align*}
\|r^{-a_1-a_2}(\dr^{j_1}f)(\dr^{j_2}g)\|_{L^2(\Omega)}\leq\|r^{-a_1}\dr^{j_1}f\|_{L^\infty(\Omega)}\|\dr^{j_2}g\|_{\Omega,a_2}.
\end{align*}
By Lemma \ref{embeddings}(ii),
\begin{align*}
\|r^{-a_1}\dr^{j_1}f\|_{L^\infty(\Omega)}\lesssim\theta_b^{-{1\over 2}}\|\dr^{j_1}f\|_{\Omega,a_1+1,1,1}\leq\theta_b^{-{1\over 2}}\|f\|_{\Omega,a_1+1,j+1,1}\lesssim\theta_b^{-{1\over 2}}\|f\|_{\Omega,a_1+1,j,2},
\end{align*}
which is why the second full derivative is needed on the right of \eqref{bulk product}.
\end{proof}

\subsubsection{Estimates of the nonlinearities}

We first record the schematic structure of the nonlinearities of System \ref{reference wedge system} at $(\phi_0,\mb v,q)$, the pair $(\mb v,q)$ in place of $(\mb u,p)$. By Lemma \ref{Piola transform} we have $\grad_\phi=-(1+\mu_1)\bs\nu\partial_x$ and $\mb u_\phi=-(1+\mu_1)\bs\nu^\perp v_y$. Recall the classes $\mathfrak G$ and $\mathfrak G_0$ of Lemma \ref{composition of the geometry}, of which $1+\mu_1$ and $\mu_1$ are respectively examples. Then
\begin{align}\label{schematic bulk}
\mb N[\phi_0,\mb v,q]=\sum\mathfrak G\cdot\brac{\grad^{k_1}\bs\nu}\brac{\grad^{k_2}\bs\nu}\cdot\grad^{j}Z,\qquad\qquad Z\in\set{v_y,\ \mb v,\ q},
\end{align}
where $\mb v$ stands for any of its components, each summand carries one or two factors of $\bs\nu$ type, the factor $\grad^{k_2}\bs\nu$ being absent in the first case, and the total number of derivatives satisfies $k_1+k_2+j\leq 2$ if $Z$ is a velocity component and $k_1+k_2+j\leq 1$ if $Z=q$. Indeed, we expand the three groups of $\mb N$ by the Leibniz rule. A derivative falling on a factor $\mathfrak G$ produces a component of $\grad\bs\nu$ times a factor $\mathfrak G$, by the second identity of \eqref{composition}. A factor of $\bs\nu$ type without derivative is absorbed into $\mathfrak G$, since $\nu_1\mathfrak g(\bs\nu)$ and $\nu_2\mathfrak g(\bs\nu)$ are again smooth functions of $\bs\nu$; but when no factor of $\bs\nu$ type is differentiated, one of them is kept. Since at most two derivatives fall on the geometry, one or two factors of $\bs\nu$ type remain, and the groups contribute the following summands:
\begin{align*}
\begin{array}{c|c|c}
\text{group of }\mb N&Z&\text{derivatives}\\\hline
-(\grad+\grad_\phi)\cdot(\grad+\grad_\phi)\mb u_\phi&v_y&k_1+k_2+j=2\\
-(\grad\cdot\grad_\phi+\grad_\phi\cdot\grad+\grad_\phi\cdot\grad_\phi)\mb u&\text{a component of }\mb v&k_1+k_2+j=2,\ j\geq 1\\
\grad_\phi p&q&j=1,\ k_1=0
\end{array}
\end{align*}
In particular two derivatives fall on the geometry only when $Z=v_y$. 

For the boundary nonlinearity we use \eqref{profile trace}, which turns every factor built from the profile into a factor of geometry type. In \eqref{profile trace}, $\phi_0$ and $\partial_r\phi_0$ are exhibited as traces of $-(b^{-1}\nu_1+\nu_2)$ and of $-r^{-1}\dr(b^{-1}\nu_1+\nu_2)$. So the composite data factors of System \ref{reference wedge system}, namely the surface tension difference $G(\phi_0)$ and the normal $\mb n_\phi=\phi_0\mb B(\phi_0)$ with $G,\mb B$ smooth and $G(0)=0$, are traces of factors of type $\mathfrak G_0$. Hence, with $\mb e$ denoting any of the bounded degree zero fields $\mb e_1,\mb e_2,\mb e_r,\mb e_\theta$, $\mb N_\partial$ is the trace of a bulk expression of the same shape as \eqref{schematic bulk},
\begin{align}\label{schematic boundary}
\mb N_\partial[\phi_0,\mb v,q]=\eva{\brac{\underbrace{\mathfrak G_0\cdot\brac{r^{-1}\dr\nu_1}\mb e+\mathfrak G_0\cdot\brac{r^{-1}\dr\nu_2}\mb e}_{\text{pure geometry}}+\sum\mathfrak G\cdot\brac{\grad^{k}\bs\nu}\cdot\grad^{j}Z}}_{\partial_1\Omega},\qquad k+j\leq 1,
\end{align}
where each summand of the second sum carries at least one factor of $\bs\nu$ type, and where $j=0$ whenever $Z=q$, the pressure entering $\mb N_\partial$ only through the undifferentiated product $p\mb n_\phi$. The groups of $\mb N_\partial$ contribute the following summands, with $Z=v_y$ in those with $k=1$, the geometry being differentiated only through $\mb u_\phi$:
\begin{align*}
\begin{array}{c|c|c}
\text{group of }\mb N_\partial&Z&\text{derivatives}\\\hline
(\sin\theta_b)^{-1}\brac{(1+(b^{-1}-\phi_0)^2)^{-{3\over 2}}-(\sin\theta_b)^3}(\partial_r\phi_0)\mb e_\theta&\text{---}&\text{pure geometry}\\
(\sin\theta_b)^{-1}(1+(b^{-1}-\phi_0)^2)^{-{3\over 2}}(\partial_r\phi_0)\mb n_\phi&\text{---}&\text{pure geometry}\\
-\brac{(\grad+\grad_\phi)(\mb u+\mb u_\phi)+((\grad+\grad_\phi)(\mb u+\mb u_\phi))^T}\mb n_\phi&\text{a component of }\mb v&k+j=1\\
-\brac{\grad_\phi(\mb u+\mb u_\phi)+\grad\mb u_\phi+(\grad_\phi(\mb u+\mb u_\phi))^T+(\grad\mb u_\phi)^T}\mb e_\theta&\text{a component of }\mb v&k+j=1\\
p\mb n_\phi&q&k=j=0
\end{array}
\end{align*}

We call the factors $\grad^jZ$ solution factors and the other factors, all built from $\bs\nu$, geometry factors. The first sum of \eqref{schematic boundary} collects the two summands of $\mb N_\partial$ carrying no solution factor, the pure geometry terms of the first two rows; both are quadratic in the geometry, the second factor of $\bs\nu$ type sitting inside $\mathfrak G_0$. By \eqref{schematic bulk}--\eqref{schematic boundary}, $\mb N$ is linear and $\mb N_\partial$ is affine in $(\mb v,q)$ for fixed $\phi_0$. The formal time linearisations $D\mb N$ and $D\mb N_\partial$ at $(\phi_0,\mb v,q)$ in an arbitrary direction $(\dot\psi,\dot{\mb v},\dot q)$ have the same form, with exactly one factor in each summand replaced by its linearisation, namely $\bs\nu\mapsto\bs\nu[\dot\psi]$ and $Z\mapsto\dot Z$. By the chain rule, as in the second identity of \eqref{composition}, the linearisation of a factor $\mathfrak G$ is again a $\mathfrak G$ times a component of $\bs\nu[\dot\psi]$. After \eqref{profile trace} the profile direction $\dot\psi$ enters $D\mb N_\partial$ only through $\bs\nu[\dot\psi]$, as it enters $D\mb N$.

\begin{proposition}[Nonlinear estimates]\label{nonlinear estimates}
With $\kappa$ and $l$ as fixed at the beginning of Subsection \ref{initial-time construction}, assume \eqref{standing angle condition} and $\mathcal D\leq\delta_0$, with $\delta_0$ as in Lemma \ref{geometry estimates}(iv). Let $(\mb v,q)$ satisfy $v_\theta|_{\partial_0\Omega}=0$ with $\mathcal E=\mathcal E(\mb v,q)<\infty$, and let $(\dot{\mb v},\dot q,\dot\psi)$ be such that $\dot{\mathcal E}=\dot{\mathcal E}(\dot{\mb v},\dot q)<\infty$ and $\dot{\mathcal D}=\dot{\mathcal D}(\dot\psi)<\infty$. Then
\begin{align}
\sum_{M=0}^2\|\mb N[\phi_0,\mb v,q]\|_{\Omega,M-1-\kappa,l+2}+\sum_{M=0}^2\|\mb N_\partial[\phi_0,\mb v,q]\|_{\partial_1\Omega,M-{1\over 2}-\kappa,l+{5\over 2}}&\lesssim_{\theta_b}\mathcal D\brac{\mathcal D+\mathcal E},\label{N estimate}\\
\sum_{M=0}^1\|D\mb N\|_{\Omega,M-1-\kappa,l+1}+\sum_{M=0}^1\|D\mb N_\partial\|_{\partial_1\Omega,M-{1\over 2}-\kappa,l+{3\over 2}}&\lesssim_{\theta_b}\mathcal D\dot{\mathcal E}+\brac{\mathcal D+\mathcal E}\dot{\mathcal D}.\label{DN estimate}
\end{align}
Moreover, in the norms on the left of \eqref{N estimate} and \eqref{DN estimate}, every summand of \eqref{schematic bulk}--\eqref{schematic boundary} and of their linearisations is bounded by a product of norms of its factors, each of which is bounded by $\mathcal D$, $\mathcal E$, $\dot{\mathcal D}$ or $\dot{\mathcal E}$. If one factor is replaced by the difference of two factors of the same kind, the same bound holds with the norm of that factor evaluated at the difference. In particular, for $(\mb v,q)$ and $(\mb v',q')$ as above,
\begin{align}
&\sum_{M=0}^2\|\mb N[\phi_0,\mb v,q]-\mb N[\phi_0,\mb v',q']\|_{\Omega,M-1-\kappa,l+2}\nonumber\\
&\qquad+\sum_{M=0}^2\|\mb N_\partial[\phi_0,\mb v,q]-\mb N_\partial[\phi_0,\mb v',q']\|_{\partial_1\Omega,M-{1\over 2}-\kappa,l+{5\over 2}}
\lesssim_{\theta_b}\mathcal Dd\brac{(\mb v,q),(\mb v',q')}.\label{N difference}
\end{align}
\end{proposition}
\begin{proof}
We estimate the summands of \eqref{schematic bulk}. We call the norms on the left of \eqref{N estimate} and \eqref{DN estimate} the targets, and their weights $M-1-\kappa$ and $M-1/2-\kappa$ the target weights. The factors $\mathfrak G$ are suppressed from the notation: they are bounded by Lemma \ref{composition of the geometry}(i), and the Leibniz expansion stays within the schematic form \eqref{schematic bulk} by the second identity of \eqref{composition}. Likewise a second factor of $\bs\nu$ type carrying no derivative is bounded by \eqref{geometry smallness} and suppressed, and one carrying derivatives is bounded by \eqref{geometry sup}, at the cost of one further factor $\mathcal D\leq\delta_0$. We spell out four representative cases; after the Leibniz expansion all remaining summands are of one of these forms and are estimated by the identical pattern.

\emph{Case 1: all derivatives on the solution.} Consider $\bs\nu\grad^2\mb v$ at level $M$ and split $\mb v$ by \eqref{data step relations}. For $\mb v^{[M]}$, \eqref{bulk product} with $(a_1,a_2)=(0,M-1-\kappa)$ gives
\begin{align*}
\|\bs\nu\grad^2\mb v^{[M]}\|_{\Omega,M-1-\kappa,l+2}
\lesssim_{\theta_b}\|\bs\nu\|_{\Omega,1,l+2,2}\|\grad^2\mb v^{[M]}\|_{\Omega,M-1-\kappa,l+2}
\lesssim_{\theta_b}\mathcal D\mathcal E,
\end{align*}
where the geometry factor is estimated by \eqref{geometry bulk} at $a=1$,
\begin{align}\label{basic geometry bound}
\|\bs\nu\|_{\Omega,1,l+2,2}\lesssim\|\bs\nu\|_{\Omega,1,0,l+4}\lesssim_{\theta_b}\|\phi_0\|_{\partial_1\Omega,{1\over 2},l+{7\over 2}}\lesssim_{\theta_b}\mathcal D.
\end{align}
For $M\geq 1$ the pieces subtracted in $\mb v^{\smash[t]{[M]}}$ have to be estimated as well: the constant $\mb v(\mb 0)$ is annihilated by $\grad^2$, while the degree one fields, present for $M=2$, satisfy $|\dr^{j'}\grad^2\mb P|\lesssim_{\theta_b} r^{-1}$, so that by \eqref{corner coefficients 2} and \eqref{geometry bulk} at $a=M-\kappa\in[1-\kappa,2-\kappa]$,
\begin{align*}
\brac{|v_1(\mb 0)|+(\sin\theta_b)^2|c_{\phi_0,\dr_{-1}^2}|}\|\bs\nu\grad^2\mb P\|_{\Omega,M-1-\kappa,l+2}
\lesssim_{\theta_b}\mathcal E\|\bs\nu\|_{\Omega,M-\kappa,l+2}
\lesssim_{\theta_b}\mathcal D\mathcal E.
\end{align*}

\emph{Case 2: all derivatives on the geometry.} Consider $(\grad^2\bs\nu)v_y$ at $M=2$, so that the target weight is $1-\kappa$. By Lemma \ref{decomposition estimates}(iii), $v_y(\mb 0)=0$, and therefore
\begin{align}\label{v_y splitting}
v_y=v_y^{\deg 1}+v_y^{[2]},\qquad
v_y^{\deg 1}:={v_1(\mb 0)\over 2}\mb P_{\mb u_{F_{0,0}},\dr_{-1}}\cdot\mb e_2+{(\sin\theta_b)^2c_{\phi_0,\dr_{-1}^2}\over 4}\mb P_{\mb u_{F_{1,1}},\dr_{-1}}\cdot\mb e_2,
\end{align}
with $|\dr^{j'}v_y^{\deg 1}|\lesssim_{\theta_b}r(|v_1(\mb 0)|+(\sin\theta_b)^2|\smash[b]{c_{\phi_0,\dr_{-1}^2}}|)\lesssim_{\theta_b}r\mathcal E$ by \eqref{corner coefficients 2}. For the degree one part we therefore have
\begin{align*}
\|(\dr^{j'}\grad^2\bs\nu)v_y^{\deg 1}\|_{\Omega,1-\kappa}
\lesssim_{\theta_b}\mathcal E\|\grad^2\bs\nu\|_{\Omega,-\kappa,l+2}
\lesssim_{\theta_b}\mathcal E\|\bs\nu\|_{\Omega,2-\kappa,0,l+4}
\lesssim_{\theta_b}\mathcal E\mathcal D,
\end{align*}
by \eqref{geometry bulk} at $a=2-\kappa$ with count $l+4$, i.e. by $\|\phi_0\|_{\partial_1\Omega,3/2-\kappa,l+7/2}$. This is the term that uses the full count of this norm in $\mathcal D$.

For the remainder piece we expand $\dr^{l+2}$ by the Leibniz rule and put the solution factor in the weighted sup norm by Lemma \ref{embeddings}(ii) with its parameter equal to $2-\kappa$; the weights split exactly as $-(1-\kappa)-2=-1-(2-\kappa)$, so that, by \eqref{basic geometry bound},
\begin{align*}
\|(\dr^{j_1}\grad^2\bs\nu)(\dr^{j_2}v_y^{[2]})\|_{\Omega,1-\kappa}
&\leq\|r^{-(2-\kappa)}\dr^{j_2}v_y^{[2]}\|_{L^\infty(\Omega)}\|\dr^{j_1}\grad^2\bs\nu\|_{\Omega,-1}\\
&\lesssim_{\theta_b}\|\mb v^{[2]}\|_{\Omega,3-\kappa,j_2,2}\|\bs\nu\|_{\Omega,1,0,l+4}
\lesssim_{\theta_b}\mathcal E\mathcal D.
\end{align*}

The analogous terms at $M\in\{0,1\}$ are estimated in the same way with strictly easier weights. At $M=1$, where the target weight is $-\kappa$, we have $v_y=v_y^{\smash[t]{[1]}}$, since $v_y(\mb 0)=0$, and the solution factor in the weighted sup norm with its parameter equal to $1-\kappa$, together with the same geometry factor, gives
\begin{align*}
\|(\dr^{j_1}\grad^2\bs\nu)(\dr^{j_2}v_y^{[1]})\|_{\Omega,-\kappa}
\leq\|r^{-(1-\kappa)}\dr^{j_2}v_y^{[1]}\|_{L^\infty(\Omega)}\|\dr^{j_1}\grad^2\bs\nu\|_{\Omega,-1}
\lesssim_{\theta_b}\mathcal E\mathcal D,
\end{align*}
the weights again splitting exactly, $-(-\kappa)-2=-1-(1-\kappa)$. At $M=0$ no splitting of $\mb v$ is needed and the geometry is always placed in the sup norm at $a=1$ or in weight $1$ norms. The terms $(\grad\bs\nu)(\grad v_y)$ are estimated in the same way, each factor carrying one full derivative fewer than in Cases 1 and 2, except when all $l+2$ radial derivatives fall on $\grad v_y$, where Lemma \ref{embeddings}(ii) would ask the count $(l+3,2)$, beyond that of $\mathcal E$. In that case we put the geometry in the weighted sup norm instead, at $a=1$; at $M=2$ the weights again split exactly, $-(1-\kappa)-1=-(2-\kappa)$, and
\begin{align*}
\|(\grad\bs\nu)(\dr^{l+2}\grad v_y^{[2]})\|_{\Omega,1-\kappa}
\leq\|r\grad\bs\nu\|_{L^\infty(\Omega)}\|\dr^{l+2}\grad v_y^{[2]}\|_{\Omega,2-\kappa}
\lesssim_{\theta_b}\mathcal D\mathcal E,
\end{align*}
by \eqref{geometry sup} at $a=1$ and Lemma \ref{geometry estimates}(ii).

\emph{Case 3: the pressure terms.} Here $\grad_\phi q=-(1+\mu_1)\bs\nu\partial_xq$, and we split $q$ by \eqref{data step relations}. The constants are annihilated by $\partial_x$, while the logarithm in $q-q^{[2]}$ satisfies $|\dr^{j'}\partial_x(2\ln r+1)|\lesssim r^{-1}$ and therefore contributes, at $M=2$,
\begin{align*}
(\sin\theta_b)^2|c_{\phi_0,\dr_{-1}^2}|\|\bs\nu r^{-1}\|_{\Omega,1-\kappa,l+2}
\lesssim\mathcal D\|\bs\nu\|_{\Omega,2-\kappa,l+2}\lesssim_{\theta_b}\mathcal D^2,
\end{align*}
by \eqref{profile coefficients}. The regular part is handled as in Case 1, by \eqref{bulk product} and \eqref{basic geometry bound},
\begin{align*}
\|\bs\nu\partial_xq^{[M]}\|_{\Omega,M-1-\kappa,l+2}\lesssim_{\theta_b}\|\bs\nu\|_{\Omega,1,l+2,2}\|\grad q^{[M]}\|_{\Omega,M-1-\kappa,l+2}\lesssim_{\theta_b}\mathcal D\mathcal E.
\end{align*}

\emph{Case 4: the boundary terms.} By \eqref{schematic boundary}, every summand of $\mb N_\partial$ is $H|_{\partial_1\Omega}$, where $H$ is a bulk expression. Since $\|H|_{\partial_1\Omega}\|_{\partial_1\Omega,M-1/2-\kappa,l+5/2}\lesssim_{\theta_b}\|H\|_{\Omega,M-\kappa,l+2,1}$ by Lemma \ref{trace}, it suffices to prove that
\begin{align}\label{boundary bulk target}
\|H\|_{\Omega,M-\kappa,l+2,1}\lesssim_{\theta_b}\mathcal D\brac{\mathcal D+\mathcal E}.
\end{align}

For the summands of \eqref{schematic boundary} containing a solution factor we have
\begin{align*}
\|H\|_{\Omega,M-\kappa,l+2,1}\lesssim\|r^{-1}H\|_{\Omega,M-1-\kappa,l+2}+\|\grad H\|_{\Omega,M-1-\kappa,l+2}.
\end{align*}
By \eqref{schematic boundary}, $\grad H$ is a sum of summands of the form \eqref{schematic bulk}, in which two derivatives fall on the geometry only when $Z=v_y$, as Case 2 requires. The function $r^{-1}H$ is of the same form with $r^{-1}$ in place of one $\grad$. So the terms on the right are of the form of those of Cases 1--3, in the same norm, and are estimated in the same way, except for the pressure.

As noted after \eqref{schematic boundary}, the pressure enters only through the summand $H=\bs\nu q$, from $p\mb n_\phi$, in which it carries no derivative, whereas in Case 3 it carries one. For its regular part $q^{[M]}$ this difference is immaterial: the placements of $\grad\mb v^{[M]}$ in the summands $\bs\nu\grad\mb v$ use only $\|\grad\mb v^{[M]}\|_{\Omega,M-\kappa,l+2,1}$, so they apply to $q^{[M]}$, since $\|q^{[M]}\|_{\Omega,M-\kappa,l+2,1}\leq\mathcal E$. It matters for the terms subtracted in $q^{[2]}$ at $M=2$, which Case 3 differentiates: by \eqref{data step relations}, $q-q^{[2]}$ is a constant plus a multiple of $\ln r$, with coefficients $\lesssim_{\theta_b}\mathcal D+\mathcal E$ by \eqref{corner coefficients 2} and \eqref{profile coefficients}. Since $1+|\ln r|\lesssim_\kappa r^{-\kappa/2}+r^{\kappa/2}$, while $\dr\ln r=1$ and $r\grad\ln r=\mb e_r$ are bounded, we place the geometry at the weights $2-\kappa\pm\kappa/2$, which lie in \eqref{margin ranges}, and obtain by \eqref{geometry bulk} and Lemma \ref{geometry estimates}(ii)
\begin{align*}
\|\bs\nu(q-q^{[2]})\|_{\Omega,2-\kappa,l+2,1}\lesssim_{\theta_b}\brac{\mathcal D+\mathcal E}(\|\bs\nu\|_{\Omega,2-{\kappa\over 2},l+2,1}+\|\bs\nu\|_{\Omega,2-{3\kappa\over 2},l+2,1})\lesssim_{\theta_b}\mathcal D\brac{\mathcal D+\mathcal E}.
\end{align*}

There remain the two pure geometry terms of \eqref{schematic boundary}; we write $\Xi$ for either of them, $\mathfrak G_0\cdot(r^{-1}\dr\nu_1)\mb e$ or $\mathfrak G_0\cdot(r^{-1}\dr\nu_2)\mb e$. They carry no solution factor, and for them \eqref{boundary bulk target} must give $\lesssim_{\theta_b}\mathcal D^2$. Expand $\dr^{j'}(r\grad)^{m'}\Xi$ with $j'\leq l+2$ and $m'\leq 1$ by \eqref{composition}. Every resulting term is then of the form $r^{-1}AB$, where $B$ is the factor carrying the most derivatives and $A$ is the product of the remaining ones, of which there is at least one, since $\mathfrak G_0$ vanishes when $\bs\nu=\mb 0$; all of these factors are derivatives of components of $\bs\nu$, of order at most $l+4$ in total, multiplied by a bounded $\mathfrak G$. The factors making up $A$ carry at most half of the available derivatives each, hence at most $l+2$, so that they may be placed in the sup norms \eqref{geometry sup}, and we obtain
\begin{align*}
\|r^{-1}AB\|_{\Omega,M-\kappa}
\leq\|r^{1-a}A\|_{L^\infty(\Omega)}\|B\|_{\Omega,M+2-\kappa-a}
\lesssim_{\theta_b}\mathcal D^2,
\end{align*}
where we choose $a=1$ when $M\leq 1$ and $a=2-\kappa/2$ when $M=2$. Both weights lie in \eqref{margin ranges}: the sup weight is $a\in\{1,2-\kappa/2\}$, while the bulk weight is $M+2-\kappa-a\in\{1-\kappa,2-\kappa\}$ when $M\leq 1$ and $2-\kappa/2$ when $M=2$; at $M=2$ the sup weight $a$ and the bulk weight $4-\kappa-a$ both lie in \eqref{margin ranges} only for $a=2-\kappa/2$. At $M=2$ the choice $a=1$ would instead ask for \eqref{geometry bulk} at the weight $3-\kappa$, where its RHS is $\|\phi_0\|_{\partial_1\Omega,5/2-\kappa,l+7/2}$, which is infinite unless the block of $\phi_0$ at $d=1$ vanishes, since $\|r\|_{\partial_1\Omega,5/2-\kappa}$ and $\|r\ln r\|_{\partial_1\Omega,5/2-\kappa}$ are infinite. This proves \eqref{boundary bulk target} for the pure geometry terms, and with it \eqref{N estimate}.

For \eqref{DN estimate}, each summand of $D\mb N$ and $D\mb N_\partial$ is a summand of \eqref{schematic bulk}--\eqref{schematic boundary} with one factor replaced by its linearisation. If the replaced factor is a solution factor, the estimates above apply with $(\dot{\mb v},\dot q)$ in place of $(\mb v,q)$, with one derivative fewer, $l+1$ in place of $l+2$, and with $M\leq 1$, producing the terms $\mathcal D\dot{\mathcal E}$. No side condition is asked of $\dot{\mb v}$: the splitting \eqref{v_y splitting} is needed only at $M=2$, and at $M=1$ the constant $\dot v_y(\mb 0)$, which need not vanish, is paired with the geometry at the weight $2-\kappa$,
\begin{align*}
\|(\dr^{l+1}\grad^2\bs\nu)\dot v_y(\mb 0)\|_{\Omega,-\kappa}\leq|\dot{\mb v}(\mb 0)|\|\grad^2\bs\nu\|_{\Omega,-\kappa,l+1}\lesssim_{\theta_b}\dot{\mathcal E}\|\bs\nu\|_{\Omega,2-\kappa,0,l+3}\lesssim_{\theta_b}\mathcal D\dot{\mathcal E},
\end{align*}
by \eqref{corner coefficients 1} and Lemma \ref{geometry estimates}(ii). If the replaced factor is a geometry factor, which by \eqref{profile trace} now includes all the factors built from the profile, we use Lemma \ref{geometry estimates}(iii) for $\bs\nu[\dot\psi]$; for instance, in the analogue of Case 2 at $M=1$, where the target weight is $-\kappa$ and $v_y=v_y^{\smash[t]{[1]}}$, the solution factor is again in the weighted sup norm with its parameter equal to $1-\kappa$, and for $j_1+j_2\leq l+1$
\begin{align*}
\|(\dr^{j_1}\grad^2\bs\nu[\dot\psi])(\dr^{j_2}v_y^{[1]})\|_{\Omega,-\kappa}
&\leq\|r^{-(1-\kappa)}\dr^{j_2}v_y^{[1]}\|_{L^\infty(\Omega)}\|\dr^{j_1}\grad^2\bs\nu[\dot\psi]\|_{\Omega,-1}\\
&\lesssim_{\theta_b}\mathcal E\|\bs\nu[\dot\psi]\|_{\Omega,1,0,l+3}
\lesssim_{\theta_b}\mathcal E\|\dot\psi\|_{\partial_1\Omega,{1\over 2},l+{5\over 2}}
\lesssim_{\theta_b}\mathcal E\dot{\mathcal D},
\end{align*}
the count $l+5/2$ of $\dot\psi$ being exactly the budget of $\dot{\mathcal D}$; the other pairings ask for at most this count. The pure geometry terms of \eqref{schematic boundary}, being quadratic in $\bs\nu$, linearise into terms in each of which exactly one factor is of $\bs\nu[\dot\psi]$ type, and the splitting used for $\Xi$ applies unchanged, the reduced target count $l+3/2$ asking for at most $l+3$ derivatives of $\bs\nu[\dot\psi]$, which is its budget in Lemma \ref{geometry estimates}(iii). This produces the terms $(\mathcal D+\mathcal E)\dot{\mathcal D}$. Every placement above bounds a summand by a product of norms of its factors, which gives the multilinear clause of the proposition. Since $\mb N$ is linear and $\mb N_\partial$ is affine in the pair, the pure geometry terms of \eqref{schematic boundary} cancel in the differences of \eqref{N difference}, and by Remark \ref{metric remark} the proof above applies to the differences with the frozen constants set to $0$, which gives \eqref{N difference}.
\end{proof}

\begin{remark}[The origin of the second angle condition]\label{sharpness of the second angle condition}
The mixed harmonic extension $\E$ carries the resonant mode $r^{\pi/(2\theta_b)}$ of the Dirichlet--Neumann wedge problem, with a coefficient proportional to $\hat\phi_0(\pi/(2\theta_b)-1)$; consequently $\bs\nu[\phi_0]=\grad\E\chi$ behaves like $r^{\pi/(2\theta_b)-1}$ at the corner for generic data. The summand $(\grad^2\bs\nu)u_y$ of $\lpc\mb u_\phi$ in $\mb N$, in which $u_y$ vanishes linearly at the corner, behaves there like $r^{\pi/(2\theta_b)-2}$. Since the measure on $\Omega$ is $r\d r\d\theta$, this summand lies in the $M=2$ forcing class $\|\ph\|_{\Omega,1-\kappa}$ if and only if
\begin{align*}
-2(1-\kappa)+2\brac{{\pi\over 2\theta_b}-2}+1>-1,\qquad\text{ i.e. }\qquad {\pi\over 2\theta_b}>2-\kappa.
\end{align*}
The pure geometry term $\Xi$ of $\mb N_\partial$ in Case 4 of the proof of Proposition \ref{nonlinear estimates}, quadratic in the geometry with one derivative, behaves like $r^{\pi/\theta_b-3}$. Its bulk target $\|\Xi\|_{\Omega,2-\kappa}$ there is finite if and only if
\begin{align*}
-2(2-\kappa)+2\brac{{\pi\over\theta_b}-3}+1>-1,\qquad\text{ i.e. }\qquad {\pi\over 2\theta_b}>2-{\kappa\over 2},
\end{align*}
which is the second condition of \eqref{standing angle condition}. For $(2-\kappa/2)\theta_b\geq\pi/2$ the bulk expression $\Xi$ therefore leaves the class $\|\ph\|_{\Omega,2-\kappa}$, and \eqref{boundary bulk target} fails for $H=\Xi$ at $M=2$. On $\partial_1\Omega$, however, the pure geometry terms are functions of $\phi_0$ and $\partial_r\phi_0$, which the mode does not enter, its gradient being normal to $\partial_1\Omega$; and the summand $(\grad^2\bs\nu)u_y$ needs only the condition of the first display. So the second condition of \eqref{standing angle condition} comes from bounding the trace through $\Xi$ in Case 4. What the second condition of \eqref{standing angle condition} records is where the explicit part of the argument stops, at the resonance $r^{\pi/(2\theta_b)}$ of the flattening, and not the range of angles actually covered: the linear theory of Section \ref{linear existence} asks only for an unquantified smallness, see Remark \ref{smallness remark}.
\end{remark}

\subsubsection{The contact point and the derived profiles}

\begin{lemma}[Corner values]\label{corner values}
Let $(\mb u,p)$ be a pair with $\mathcal E=\mathcal E(\mb u,p)<\infty$ and let $\dot\phi$ be the profile derived from it as in Lemma \ref{Piola transform}, that is, $\dot\phi:=(\sin\theta_b)^{-2}\partial_ru_\theta(\ph,\theta_b)$ on $\partial_1\Omega$. Then $\partial_ru_\theta(r,\theta_b)$ has a limit $\partial_ru_\theta(0,\theta_b)$ as $r\to 0$, and
\begin{align}\label{corner identity}
\partial_ru_\theta(r,\theta_b)=\partial_ru_\theta(0,\theta_b)+\partial_ru^{[2]}_\theta(r,\theta_b),\quad\text{ where }\quad\partial_ru_\theta(0,\theta_b)=-{u_1(\mb 0)\over 2}+{\theta_b(\sin\theta_b)^2\over 2}c_{\phi_0,\dr_{-1}^2},
\end{align}
with $|\partial_r\mb u^{\smash[t]{[2]}}|\lesssim_{\theta_b}\mathcal Er^{1-\kappa}$, and
\begin{subequations}
\begin{align}\label{dotphi bound}
(\sin\theta_b)^2\|\dot\phi\|_{\partial_1\Omega,M-{1\over 2}-\kappa,l+{5\over 2}}\lesssim_{\theta_b}\mathcal E\qquad\text{ for }\qquad M\in\{0,1\}.
\end{align}
If moreover $\partial_ru_\theta(0,\theta_b)=0$, that is, if $\dot\phi$ vanishes at the contact point, then \eqref{dotphi bound} holds also for $M=2$, and $\dot{\mathcal D}(\dot\phi)\lesssim_{\theta_b}\mathcal E$; conversely, if $\partial_ru_\theta(0,\theta_b)\neq 0$, then the LHS of \eqref{dotphi bound} is infinite for $M=2$. The same estimates hold at the dotted level: if $\dot{\mathcal E}(\dot{\mb u},\dot p)<\infty$ and $\ddot\phi:=(\sin\theta_b)^{-2}\partial_r\dot u_\theta(\ph,\theta_b)$, then
\begin{align}\label{ddotphi bound}
(\sin\theta_b)^2\|\ddot\phi\|_{\partial_1\Omega,M-{1\over 2}-\kappa,l+{3\over 2}}\lesssim_{\theta_b}\dot{\mathcal E}\qquad\text{ for }\qquad M\in\{0,1\}.
\end{align}
\end{subequations}
\end{lemma}
\begin{proof}
By \eqref{data step relations},
\begin{align}\label{u expansion}
\mb u=\mb u(\mb 0)+{u_1(\mb 0)\over 2}\mb P_{\mb u_{F_{0,0}},\dr_{-1}}+{(\sin\theta_b)^2c_{\phi_0,\dr_{-1}^2}\over 4}\mb P_{\mb u_{F_{1,1}},\dr_{-1}}+\mb u^{[2]}.
\end{align}
On $\theta=\theta_b$ the $\mb e_\theta$ components of the two degree one fields are $-r$ and $2\theta_br$, while the $\mb e_\theta$ component of the constant $\mb u(\mb 0)$ does not depend on $r$; applying $\partial_r$ gives $\partial_ru_\theta(r,\theta_b)=-u_1(\mb 0)/2+\theta_b(\sin\theta_b)^2\smash[b]{c_{\phi_0,\dr_{-1}^2}}/2+\partial_ru^{\trimtop{0.9pt}{[2]}}_\theta(r,\theta_b)$. The bound on $\partial_r\mb u^{[2]}$ follows from Lemma \ref{embeddings}(ii) at $a=1-\kappa$ applied to $\grad\mb u^{[2]}$,
\begin{align*}
\|r^{-(1-\kappa)}\grad\mb u^{[2]}\|_{L^\infty(\Omega)}\lesssim_{\theta_b}\|\grad\mb u^{[2]}\|_{\Omega,2-\kappa,1,1}\lesssim\|\mb u^{[2]}\|_{\Omega,3-\kappa,1,2}\leq\mathcal E,
\end{align*}
and it shows that $\partial_ru^{\smash[t]{[2]}}_\theta(r,\theta_b)\to 0$ as $r\to 0$, which gives the limit $\partial_ru_\theta(0,\theta_b)$ and \eqref{corner identity}.

For \eqref{dotphi bound} we have, for every $M$,
\begin{align}\label{dotphi weight shift}
(\sin\theta_b)^2\|\dot\phi\|_{\partial_1\Omega,M-{1\over 2}-\kappa,l+{5\over 2}}
=\|\partial_ru_\theta\|_{\partial_1\Omega,M-{1\over 2}-\kappa,l+{5\over 2}}
\lesssim\|\dr u_\theta\|_{\partial_1\Omega,M+{1\over 2}-\kappa,l+{5\over 2}}.
\end{align}
Since $\dr$ annihilates constants we have $\dr\mb u=\dr\mb u^{\smash[t]{[1]}}$, so that $\dr u_\theta=\dr u^{\smash[t]{[M]}}_\theta$ for $M\in\{0,1\}$, whereas for $M=2$ the two degree one fields remain and, since $\dr\mb P=\mb P$,
\begin{align}\label{dr u expansion}
\dr u_\theta(r,\theta_b)=r\partial_ru_\theta(0,\theta_b)+\dr u^{[2]}_\theta(r,\theta_b).
\end{align}
For every $M\leq 2$ we have $\|\dr u^{\smash[t]{[M]}}_\theta\|_{\partial_1\Omega,M+1/2-\kappa,l+5/2}\lesssim_{\theta_b}\|\mb u^{\smash[t]{[M]}}\|_{\Omega,M+1-\kappa,l+2,2}\leq\mathcal E$ by Lemma \ref{trace}. For $M\in\{0,1\}$ this is \eqref{dotphi bound}. For $M=2$ the same bound controls $\dr u^{\trimtop{1.2pt}{[2]}}_\theta(r,\theta_b)$, whereas
\begin{align*}
\|r\|_{\partial_1\Omega,{5\over 2}-\kappa}^2=\int_0^\infty r^{-3+2\kappa}\d r=\infty,
\end{align*}
the divergence being at the origin; so $\|\dot\phi\|_{\partial_1\Omega,3/2-\kappa,l+5/2}$ is finite if and only if $\partial_ru_\theta(0,\theta_b)=0$, in which case \eqref{dotphi bound} holds for $M=2$ as well. The dotted statements are identical, with the budgets $(l+1,2)$ in place of $(l+2,2)$ and with $M\leq 1$.
\end{proof}

\subsubsection{The second-order contact angle condition}

This subsection supplies what Theorem \ref{jet corner structure} needs and the theory so far does not provide: the decomposition of the \emph{dotted} pair at the level $M=2$. As explained in Remark \ref{third rung remark}, it costs one further rung of the ladder for the undotted pair, the profile datum of the dotted problem being $\dot\phi$, and it brings the coefficient $\smash[b]{c_{\dot\phi,\dr_{-1}^3}}$ of a logarithm into the extended second compatibility condition. Throughout we strengthen the standing angle smallness to
\begin{align}\label{third rung angle condition}
3-\kappa<\Lambda^*(\theta_b)\qquad\text{ and }\qquad(3-\kappa)\theta_b<{\pi\over 2},
\end{align}
which implies \eqref{standing angle condition}. The first is the admissibility condition of Proposition \ref{third rung}, holding for $\theta_b$ up to about $45^\circ$--$50^\circ$ when $\kappa\in(1/4,1/2)$ by Remark \ref{H numerics remark}; the second is the admissibility of the harmonic extension of Lemma \ref{mixed harmonic extension} at the geometry weight $3-\kappa$, at which the corner remainder of the flattening has to be read in Proposition \ref{nonlinear corner structure} below, and it also covers the geometry weights of Proposition \ref{positive-time regularity upgrade}. The second is the more restrictive, limiting the angle to $\theta_b<\pi/(6-2\kappa)$, that is, to about $33^\circ$--$36^\circ$. We measure a profile $\psi$ by
\begin{align}\label{data norm +}
\mathcal D_+(\psi):=\mathcal D(\psi)+(\sin\theta_b)^2\brac{\sum_{k=1}^3|c_{\psi,\dr_{-2}^k}|+\|\psi^{[3]}\|_{\partial_1\Omega,{7\over 2}-\kappa,l+{7\over 2}}},\qquad\mathcal D_+:=\mathcal D_+(\phi_0),
\end{align}
where $\psi^{\smash[t]{[3]}}$ is $\psi$ minus its blocks \eqref{higher corner coefficients} at $(d,m)=(1,2)$ and $(2,3)$; so $\mathcal D_+(\psi)<\infty$ asks of $\psi$ a corner expansion of one order more than $\mathcal D(\psi)<\infty$, with second-order terms $r^2$, $r^2\ln r$ and $r^2\ln^2r$ whose coefficients are otherwise free.

The corner structure of the nonlinearities, Proposition \ref{nonlinear corner structure} below, rests on splitting every factor of \eqref{schematic bulk}--\eqref{schematic boundary} into an explicit corner part and a remainder. The solution factors carry the corner blocks of $(\mb u,p)$ through \eqref{data step relations}, and the geometry factors carry those of $\phi_0$ through the following lemma, whose general form in (i) also splits the linearised geometry in Proposition \ref{nonlinear corner structure}(ii) and the geometry at positive time in \eqref{upgrade corner splitting}.

\begin{lemma}[Corner splitting of the geometry]\label{corner splitting of the geometry}
Assume \eqref{third rung angle condition} and $\mathcal D_+<\infty$, and fix $\iota_1\in(0,\kappa)$, depending only on $\theta_b$ and $\kappa$, with $(3-\kappa+\iota_1)\theta_b<\pi/2$.
\begin{enumerate}
\item We have
\begin{align}\label{geometry corner}
\bs\nu[\phi_0]=\bs\nu_{\rm c}+\bs\nu^{[2]},\qquad
\bs\nu_{\rm c}:=c_{\phi_0,\dr_{-1}}\bs\nu[r]+c_{\phi_0,\dr_{-1}^2}\bs\nu[r\ln r]\in\mathcal Q_1^2,\qquad
\bs\nu^{[2]}:=\bs\nu[\phi_0^{[2]}],
\end{align}
where
\begin{align}\label{block geometry}
\bs\nu[r]&:=-{\sin\theta_b\over\cos(2\theta_b)}(x\mb e_1-y\mb e_2),\nonumber\\
\bs\nu[r\ln r]&:=(\ln r+\theta_b\tan(2\theta_b))\bs\nu[r]+{\sin\theta_b\over\cos(2\theta_b)}\theta(y\mb e_1+x\mb e_2).
\end{align}
More generally, let $\psi$ and $\psi^{[2]}$ be profiles that differ by the block $c_1r+c_2r\ln r+c_3r\ln^2r/2$, with constants $c_1,c_2,c_3$, and with $\|\psi\|_{\partial_1\Omega,1/2-\kappa}$ and $\|\psi^{[2]}\|_{\partial_1\Omega,5/2-\kappa-\iota}$ finite for some $\iota\in[0,1-\kappa)$. Then $\bs\nu[\psi]-\bs\nu[\psi^{\smash[t]{[2]}}]\in\mathcal Q_1^3$, and if $c_3=0$, then $\bs\nu[\psi]-\bs\nu[\psi^{\smash[t]{[2]}}]=c_1\bs\nu[r]+c_2\bs\nu[r\ln r]$.
\item For every $\iota\in[-\iota_1,1-\kappa)$,
\begin{align}
\|r^{-(2-\kappa-\iota)}(r\grad)^{m'}\bs\nu^{[2]}\|_{L^\infty(\Omega)}&\lesssim_{\theta_b,\iota}\mathcal D_+\qquad\text{ for }\qquad m'\leq l+2,\label{geometry remainder sup}\\
\|\bs\nu^{[2]}\|_{\Omega,3-\kappa-\iota,0,m}&\lesssim_{\theta_b,\iota}\mathcal D_+\qquad\text{ for }\qquad m\leq l+4.\label{geometry remainder bulk}
\end{align}
\item If $\smash[b]{c_{\phi_0,\dr_{-1}}}$ and $\smash[b]{c_{\phi_0,\dr_{-1}^2}}$ do not both vanish, then $\|\bs\nu^{\smash[t]{[2]}}\|_{\Omega\cap\{r>1\},a}=\infty$ for every $a\leq 2$.
\end{enumerate}
\end{lemma}
\begin{proof}
Such $\iota_1$ exists by the second condition of \eqref{third rung angle condition}.

(i) We prove the general statement by computing $\bs\nu$ with the Mellin transform and moving the line of integration across the pole that the block produces at $\lambda=1$. The case $\psi=\phi_0$, $\psi^{\smash[t]{[2]}}=\phi_0^{\smash[t]{[2]}}$, $\iota=0$ of the general statement is \eqref{geometry corner}: by \eqref{ladder phi2} the two profiles differ by the block $\smash[b]{c_{\phi_0,\dr_{-1}}}r+\smash[b]{c_{\phi_0,\dr_{-1}^2}}r\ln r$, and the norms $\|\phi_0\|_{\partial_1\Omega,1/2-\kappa}$ and $\|\phi_0^{[2]}\|_{\partial_1\Omega,5/2-\kappa}$ are finite since $\mathcal D<\infty$.

Let $0<a\theta_b<\pi/2$, and let $g$ be a profile with $\|g\|_{\partial_1\Omega,a-1/2}<\infty$, so that its Mellin transform is computed on the line $\Re\lambda=a-1$ by \eqref{weight and line}. With $\chi:=-(\sin\theta_b)\dr^{-1}(rg)$, as in \eqref{flattening geometry}, we have $\bs\nu[g]=\grad\E(\chi-\chi(0))$ by Lemma \ref{mixed harmonic extension}, applied at the weight $a+1$. Since $\dr\chi=-(\sin\theta_b)rg$ and $\grad=r^{-1}(\mb e_r\dr+\mb e_\theta\partial_\theta)$, Lemma \ref{Mellin properties}(i)--(ii) and the solution $\hat q_{f_1}$ of \eqref{mixed Mellin system} give
\begin{align*}
(\E(\chi-\chi(0)))^\wedge(\lambda,\theta)=-{\sin\theta_b\cos(\lambda\theta)\over\lambda\cos(\lambda\theta_b)}\hat g(\lambda_{-1}),\qquad
(\grad u)^\wedge(\lambda,\theta)=\lambda_1\hat u(\lambda_1,\theta)\mb e_r+\partial_\theta\hat u(\lambda_1,\theta)\mb e_\theta,
\end{align*}
the first on the line $\Re\lambda=a$, and the second on the line $\Re\lambda=a-1$ for a function $u$ whose Mellin transform is computed on the line $\Re\lambda=a$; here $\lambda_{\pm 1}=\lambda\pm 1$ as in Definition \ref{scaled radial derivatives}. Composing the two, and using $\cos(\lambda_1\theta)\mb e_r-\sin(\lambda_1\theta)\mb e_\theta=\cos(\lambda\theta)\mb e_1-\sin(\lambda\theta)\mb e_2$, we obtain on the line $\Re\lambda=a-1$
\begin{align*}
(\bs\nu[g])^\wedge(\lambda,\theta)=\mb m(\lambda,\theta)\hat g(\lambda),\qquad
\mb m(\lambda,\theta):=-{\sin\theta_b\over\cos(\lambda_1\theta_b)}(\cos(\lambda\theta)\mb e_1-\sin(\lambda\theta)\mb e_2).
\end{align*}

We apply this identity to $\psi$ at $a=1-\kappa$ and to $\psi^{\smash[t]{[2]}}$ at $a=3-\kappa-\iota$, both values of $a$ satisfying $0<a\theta_b<\pi/2$ by \eqref{third rung angle condition}. By Lemma \ref{block matching} with $(a_0,d,a_1)=(-\kappa,1,2-\kappa-\iota)$, the Mellin transforms of $\psi$ on the line $\Re\lambda=-\kappa$ and of $\psi^{[2]}$ on the line $\Re\lambda=2-\kappa-\iota$ are the restrictions of one function $\hat\psi$, meromorphic on the strip $-\kappa\leq\Re\lambda\leq 2-\kappa-\iota$ with its only pole at $\lambda=1$, and square integrable on the parts $|\Im\lambda|\geq 1$ of the lines of the strip, uniformly in the line. By \eqref{principal part}, whose coefficients are the $c_k$ by the same lemma, as $\lambda\to 1$,
\begin{align*}
\hat\psi(\lambda)=-{1\over\sqrt{2\pi}}\brac{{c_1\over\lambda-1}+{c_2\over(\lambda-1)^2}+{c_3\over(\lambda-1)^3}}+O(1).
\end{align*}

On the strip, $\Re(\lambda_1\theta_b)\in[(1-\kappa)\theta_b,(3-\kappa)\theta_b]\subset(0,\pi/2)$, so $\mb m$ is holomorphic there; it is also bounded there, by Lemma \ref{cosh bound}, since $\cos(\Re\lambda_1\theta_b)\geq\cos((3-\kappa)\theta_b)>0$. Hence the two horizontal sides of the rectangular contour at $\Im\lambda=\pm Y$ drop out along a sequence $Y\to\infty$, and the residue theorem, as in \eqref{residue rule}, gives
\begin{align*}
\bs\nu[\psi]-\bs\nu[\psi^{[2]}]&=-\sqrt{2\pi}\Res_{\lambda=1}\brac{r^\lambda\mb m(\lambda,\theta)\hat\psi(\lambda)}\\
&=c_1r\mb m(1,\theta)+c_2(r\ln r\mb m(1,\theta)+r\partial_\lambda\mb m(1,\theta))\\
&\quad+{c_3\over 2}(r\ln^2r\mb m(1,\theta)+2r\ln r\partial_\lambda\mb m(1,\theta)+r\partial_\lambda^2\mb m(1,\theta)),
\end{align*}
the second equality by the expansion of $\hat\psi$ above, since $r^\lambda\mb m(\lambda,\theta)$ is holomorphic at $\lambda=1$ and $\partial_\lambda r^\lambda=r^\lambda\ln r$. The RHS lies in $\mathcal Q_1^3$, since $\mb m(1,\ph)$, $\partial_\lambda\mb m(1,\ph)$ and $\partial_\lambda^2\mb m(1,\ph)$ are smooth on $[0,\theta_b]$. If $c_3=0$, the RHS is $c_1\bs\nu[r]+c_2\bs\nu[r\ln r]$: by \eqref{block geometry} and $(x,y)=r(\cos\theta,\sin\theta)$ we have $r\mb m(1,\theta)=\bs\nu[r]$, and $r\ln r\mb m(1,\theta)+r\partial_\lambda\mb m(1,\theta)=\bs\nu[r\ln r]$ since
\begin{align*}
\partial_\lambda\mb m(1,\theta)=\theta_b\tan(2\theta_b)\mb m(1,\theta)+{\sin\theta_b\over\cos(2\theta_b)}\theta(\sin\theta\mb e_1+\cos\theta\mb e_2).
\end{align*}

(ii) We apply \eqref{geometry sup} and \eqref{geometry bulk} at the weight $3-\kappa-\iota$, which is admissible since $(3-\kappa-\iota)\theta_b\leq(3-\kappa+\iota_1)\theta_b<\pi/2$, to the profile $\phi_0^{\smash[t]{[2]}}$. It remains to show that $\|\phi_0^{\smash[t]{[2]}}\|_{\partial_1\Omega,5/2-\kappa-\iota,l+7/2}\lesssim_{\theta_b,\iota}\mathcal D_+$. We split at $r=1$. On $\{r<1\}$, $\phi_0^{\trimtop{1.5pt}{[2]}}$ is $\phi_0^{\trimtop{1.5pt}{[3]}}$ plus the block of $\phi_0$ at $d=2$; this block is integrable against the weight, since $\iota>-\kappa$, and the norm of $\phi_0^{\smash[t]{[3]}}$ in $\mathcal D_+$ dominates. On $\{r>1\}$ the norm $\|\phi_0\|_{\partial_1\Omega,1/2-\kappa,l+7/2}$ in $\mathcal D$ dominates, and the subtracted block $\smash[b]{c_{\phi_0,\dr_{-1}}}r+\smash[b]{c_{\phi_0,\dr_{-1}^2}}r\ln r$ has an integrable tail, since $\iota<1-\kappa$. The coefficients of the blocks are bounded by $\mathcal D_+$ at $d=2$ and by \eqref{profile coefficients} at $d=1$.

(iii) Let $a\leq 2$ and $a':=\max\{a,1-\kappa\}$, so that $r^{-2a'}\leq\min\{r^{-2a},r^{-2(1-\kappa)}\}$ on $\{r>1\}$ and $3-2a'\geq-1$. By \eqref{block geometry}, in which $x\mb e_1-y\mb e_2$ and $y\mb e_1+x\mb e_2$ are orthogonal and of length $r$,
\begin{align*}
\|\bs\nu_{\rm c}\|_{\Omega\cap\{r>1\},a'}^2\geq{\theta_b\sin^2\theta_b\over\cos^2(2\theta_b)}\int_1^\infty|c_{\phi_0,\dr_{-1}}+c_{\phi_0,\dr_{-1}^2}(\ln r+\theta_b\tan(2\theta_b))|^2r^{3-2a'}\d r=\infty,
\end{align*}
since the first factor of the integrand is a positive constant if $\smash[b]{c_{\phi_0,\dr_{-1}^2}}=0$\vspace{1pt}, and tends to infinity as $r\to\infty$ otherwise. By (i), $\bs\nu^{[2]}=\bs\nu[\phi_0]-\bs\nu_{\rm c}$, and by Lemma \ref{geometry estimates}(ii), $\|\bs\nu[\phi_0]\|_{\Omega,1-\kappa}\lesssim_{\theta_b}\mathcal D$; hence
\begin{align*}
\|\bs\nu^{[2]}\|_{\Omega\cap\{r>1\},a}&\geq\|\bs\nu^{[2]}\|_{\Omega\cap\{r>1\},a'}\geq\|\bs\nu_{\rm c}\|_{\Omega\cap\{r>1\},a'}-\|\bs\nu[\phi_0]\|_{\Omega,1-\kappa}=\infty.\qedhere
\end{align*}
\end{proof}

The splitting $\mb u^{\smash[t]{[2]}}=\mb u^{\smash[t]{[3]}}_{\rm c}+\mb u^{\smash[t]{[3]}}$, $p^{\smash[t]{[2]}}=p^{\smash[t]{[3]}}_{\rm c}+p^{\smash[t]{[3]}}$ in the hypothesis of Proposition \ref{nonlinear corner structure} below is unique when it exists, by Lemma \ref{corner blocks under limits}(i). Neither $(\mb u,p)$ nor $(\dot{\mb u},\dot p)$ is assumed to solve the Stokes system: only $u_\theta|_{\partial_0\Omega}=\dot u_\theta|_{\partial_0\Omega}=0$, the splittings and the finiteness of $\mathcal E$, $\mathfrak X$, $\dot{\mathcal E}$ and $\dot{\mathfrak X}$ enter, which is what the iterations in the proof of Theorem \ref{jet corner structure} supply. Part (ii) of Proposition \ref{nonlinear corner structure} gives the holomorphy of $(r^2D\mb N)^\wedge$ and $(rD\mb N_\partial)^\wedge$ across $\Re\lambda=1$, the line that the dotted second rung of Proposition \ref{third rung}(i) crosses.

\begin{proposition}[Corner structure of the nonlinearities]\label{nonlinear corner structure}
With $\kappa$ and $l$ as fixed at the beginning of Subsection \ref{initial-time construction}, assume \eqref{third rung angle condition}, $\mathcal D_+<\infty$ and $\mathcal D\leq\delta_0$, with $\delta_0$ as in Lemma \ref{geometry estimates}(iv). Let $(\mb u,p)$ be a pair with $u_\theta|_{\partial_0\Omega}=0$ and $\mathcal E=\mathcal E(\mb u,p)<\infty$, whose pieces split as $\mb u^{[2]}=\mb u^{\trimtop{0.5pt}{[3]}}_{\rm c}+\mb u^{[3]}$ and $p^{[2]}=p^{\trimtop{0.5pt}{[3]}}_{\rm c}+p^{[3]}$, with $\mb u^{\smash[t]{[3]}}_{\rm c}\in\mathcal Q_2^3$, $p^{\smash[t]{[3]}}_{\rm c}\in\mathcal Q_1^3$ and
\begin{align*}
\mathfrak X:=\|\mb u^{[3]}\|_{\Omega,4-\kappa,2,2}+\Gamma_{2-\kappa}^{-1}\|p^{[3]}\|_{\Omega,3-\kappa,2,1}+\|\mb u^{[3]}_{\rm c}\|_{\mathcal Q_2^3}+\|p^{[3]}_{\rm c}\|_{\mathcal Q_1^3}<\infty.
\end{align*}
Then the following hold.
\begin{enumerate}
\item We have $\mb N[\phi_0,\mb u,p]=\mb F^{\smash[t]{[3]}}_{\rm c}+\mb F^{\smash[t]{[3]}}$ and $\mb N_\partial[\phi_0,\mb u,p]=\mb f^{\smash[t]{[3]}}_{\rm c}+\mb f^{\smash[t]{[3]}}$ as in \eqref{third rung data classes}, where the coefficients of $\mb F^{\trimtop{0.5pt}{[3]}}_{\rm c}\in\mathcal Q_0^2$ and $\mb f^{\trimtop{0.5pt}{[3]}}_{\rm c}\in\mathcal Q_1^3$ in \eqref{corner classes} are quadratic forms, with coefficients depending only on $\theta_b$ and $\theta$, in the first-order corner data $(\smash[b]{c_{\phi_0,\dr_{-1}}},\smash[b]{c_{\phi_0,\dr_{-1}^2}},u_1(\mb 0))$ of the profile and of the pair, and
\begin{align}\label{N third rung estimate}
\mathcal R_2\brac{\mb N[\phi_0,\mb u,p],\mb N_\partial[\phi_0,\mb u,p]}\lesssim_{\theta_b}\mathcal D_+\brac{\mathcal D_++\mathcal E}+\mathcal D\mathfrak X,
\end{align}
with $\mathcal R_2(\mb F,\mb f)$ as in \eqref{third rung data norm}.
\item Let $\dot\phi$ be a profile and $(\dot{\mb u},\dot p)$ a pair with $\dot u_\theta|_{\partial_0\Omega}=0$ and $\dot{\mathcal E}=\dot{\mathcal E}(\dot{\mb u},\dot p)<\infty$, let
\begin{align}\label{derived data norm plus}
\dot{\mathcal D}_+:=\dot{\mathcal D}(\dot\phi)+(\sin\theta_b)^2\brac{\sum_{k=2}^3|c_{\dot\phi,\dr_{-1}^k}|+\|\dot\phi^{[2]}\|_{\partial_1\Omega,{5\over 2}-\kappa,{5\over 2}}}<\infty,
\end{align}
and let $\dot{\mb u}^{[2]}$ and $\dot p^{[2]}$, given by \eqref{log decomposition} with $(\dot{\mb u},\dot p,\dot\phi)$ in place of $(\mb u,p,\phi)$, satisfy $\dot{\mathfrak X}:=\|\dot{\mb u}^{[2]}\|_{\Omega,3-\kappa,1,2}+\Gamma_{1-\kappa}^{-1}\|\dot p^{[2]}\|_{\Omega,2-\kappa,1,1}<\infty$. Then the linearisations $D\mb N$ and $D\mb N_\partial$ at $(\phi_0,\mb u,p)$ in the direction $(\dot\phi,\dot{\mb u},\dot p)$ have Mellin transforms $(r^2D\mb N)^\wedge(\lambda,\theta)$ and $(rD\mb N_\partial)^\wedge(\lambda)$ that are holomorphic across $\Re\lambda=1$, and
\begin{align}\label{DN second rung estimate}
\|D\mb N\|_{\Omega,1-\kappa,1}+\|D\mb N_\partial\|_{\partial_1\Omega,{3\over 2}-\kappa,{3\over 2}}
\lesssim_{\theta_b}\mathcal D\dot{\mathfrak X}+\mathcal D_+\dot{\mathcal E}+\brac{\mathcal D_++\mathcal E+\mathfrak X}\dot{\mathcal D}_+.
\end{align}
\item For $j\leq l+2$, statement (i) holds with $\mathcal R_j$ in place of $\mathcal R_2$ and with $\|\mb u^{\smash[t]{[3]}}\|_{\Omega,4-\kappa,j,2}$ and $\|p^{\smash[t]{[3]}}\|_{\Omega,3-\kappa,j,1}$ in $\mathfrak X$. For $j\leq l+1$, statement (ii) holds with $\|D\mb N\|_{\Omega,1-\kappa,j}+\|D\mb N_\partial\|_{\partial_1\Omega,3/2-\kappa,j+1/2}$ on the left of \eqref{DN second rung estimate}, with $\|\dot{\mb u}^{[2]}\|_{\Omega,3-\kappa,j,2}$ and $\|\dot p^{[2]}\|_{\Omega,2-\kappa,j,1}$ in $\dot{\mathfrak X}$, and with $\|\dot\phi^{[2]}\|_{\partial_1\Omega,5/2-\kappa,j+3/2}$ in place of $\|\dot\phi^{[2]}\|_{\partial_1\Omega,5/2-\kappa,5/2}$ in $\dot{\mathcal D}_+$. Moreover, let $\tilde\kappa\in(0,\kappa)$. In the setting of (i) at $j=l+2$,
\begin{align}\label{shifted N estimate}
\|\mb N\|_{\Omega,1+\tilde\kappa-\kappa,l+2}+\|\mb N_\partial\|_{\partial_1\Omega,{3\over 2}+\tilde\kappa-\kappa,l+{5\over 2}}
\lesssim_{\theta_b,\tilde\kappa}\mathcal D_+\brac{\mathcal D_++\mathcal E}+\mathcal D\mathfrak X,
\end{align}
and in the setting of (ii) at $j\leq l+1$,
\begin{align}\label{shifted DN estimate}
\|D\mb N\|_{\Omega,1+\tilde\kappa-\kappa,j}+\|D\mb N_\partial\|_{\partial_1\Omega,{3\over 2}+\tilde\kappa-\kappa,j+{1\over 2}}
\lesssim_{\theta_b,\tilde\kappa}\mathcal D\dot{\mathfrak X}+\mathcal D_+\dot{\mathcal E}+\brac{\mathcal D_++\mathcal E+\mathfrak X}\dot{\mathcal D}_+,
\end{align}
the arguments of the nonlinearities being those of (i) and (ii).
\end{enumerate}
\end{proposition}
\begin{proof}
(i) We write $\mb N=\mb N[\phi_0,\mb u,p]$ and $\mb N_\partial=\mb N_\partial[\phi_0,\mb u,p]$. We split every factor of the summands of \eqref{schematic bulk}--\eqref{schematic boundary} into a corner part and a remainder, collect the products of corner parts of degree zero in the bulk and of degree one on the boundary into $(\mb F^{\smash[t]{[3]}}_{\rm c},\mb f^{\smash[t]{[3]}}_{\rm c})$ (Step 1), and bound $\mb F^{\smash[t]{[3]}}:=\mb N-\mb F^{\smash[t]{[3]}}_{\rm c}$, and through Lemma \ref{trace} the bulk expression whose trace is $\mb f^{[3]}:=\mb N_\partial-\mb f^{\trimtop{0.9pt}{[3]}}_{\rm c}$, on $\{r>1\}$ (Step 2) and on $\{r<1\}$ (Step 3).

\emph{Step 1: the corner parts.} By Lemma \ref{corner splitting of the geometry}(i), by \eqref{u expansion} with $\mb u(\mb 0)=u_1(\mb 0)\mb e_1$ from Lemma \ref{decomposition estimates}(iii), by \eqref{ladder p2}, by the splitting of the hypothesis, and by Taylor's formula with integral remainder for a factor $\mathfrak G=\mathfrak g\circ\bs\nu$, the factors split as
\begin{align}
\bs\nu&=\bs\nu_{\rm c}+\bs\nu^{[2]},\nonumber\\*
\mb u&=u_1(\mb 0)\mb e_1+\mb u^{\deg 1}+\mb u^{[3]}_{\rm c}+\mb u^{[3]},\nonumber\\*
p&=(p-p^{[2]})+p^{[3]}_{\rm c}+p^{[3]},\label{factor splitting}\\*
\mathfrak G&=\mathfrak g(\mb 0)+(D\mathfrak g)(\mb 0)\bs\nu_{\rm c}+(D\mathfrak g)(\mb 0)\bs\nu^{[2]}+\nu_1^2\mathfrak G+\nu_1\nu_2\mathfrak G+\nu_2^2\mathfrak G,\nonumber
\end{align}
where the factors $\mathfrak G$ on the right are as in Lemma \ref{composition of the geometry}, $(D\mathfrak g)(\mb 0)\bs\nu_{\rm c}:=(\partial_1\mathfrak g)(\mb 0)\nu_{{\rm c},1}+(\partial_2\mathfrak g)(\mb 0)\nu_{{\rm c},2}$ with $\nu_{{\rm c},1},\nu_{{\rm c},2}$ the components of $\bs\nu_{\rm c}$, likewise for $\bs\nu^{[2]}$, and
\begin{align*}
\mb u^{\deg 1}&:={u_1(\mb 0)\over 2}\mb P_{\mb u_{F_{0,0}},\dr_{-1}}+{(\sin\theta_b)^2c_{\phi_0,\dr_{-1}^2}\over 4}\mb P_{\mb u_{F_{1,1}},\dr_{-1}},\\
p-p^{[2]}&=-{u_1(\mb 0)\over\sin(2\theta_b)}-(\sin\theta_b)^2c_{\phi_0,\dr_{-1}}-{(\sin\theta_b)^2c_{\phi_0,\dr_{-1}^2}\over 2}(2\ln r+1),
\end{align*}
so that $u^{\deg 1}_y$ is as in \eqref{v_y splitting}. The corner parts are $\bs\nu_{\rm c}$, $u_1(\mb 0)\mb e_1+\mb u^{\deg 1}+\mb u^{\smash[t]{[3]}}_{\rm c}$, $(p-p^{\smash[t]{[2]}})+p^{\smash[t]{[3]}}_{\rm c}$ and $\mathfrak g(\mb 0)+(D\mathfrak g)(\mb 0)\bs\nu_{\rm c}$, and the remainders are $\bs\nu^{\smash[t]{[2]}}$, $\mb u^{\smash[t]{[3]}}$, $p^{\smash[t]{[3]}}$ and $(D\mathfrak g)(\mb 0)\bs\nu^{\smash[t]{[2]}}+\nu_1^2\mathfrak G+\nu_1\nu_2\mathfrak G+\nu_2^2\mathfrak G$. In the classes \eqref{corner classes},
\begin{align*}
\bs\nu_{\rm c},\ (D\mathfrak g)(\mb 0)\bs\nu_{\rm c}\in\mathcal Q_1^2,\qquad
\mb u^{\deg 1}\in\mathcal Q_1^1,\qquad
\mb u^{[3]}_{\rm c}\in\mathcal Q_2^3,\qquad
p-p^{[2]}\in\mathcal Q_0^2,\qquad
p^{[3]}_{\rm c}\in\mathcal Q_1^3.
\end{align*}
We call $d$ the degree of $r^d\sum_k(\ln r)^k\mb q_k(\theta)$, also when $\mb q_k\in L^2([0,\theta_b])$. A derivative lowers the degree by one and in products the degrees add:
\begin{align*}
\grad\mathcal Q_d^m\subset\mathcal Q_{d-1}^m\quad\text{for smooth coefficients},\qquad\mathcal Q_d^m\mathcal Q_{d'}^{m'}\subset\mathcal Q_{d+d'}^{m+m'-1}.
\end{align*}

In the summands of \eqref{schematic bulk}--\eqref{schematic boundary} with one factor of $\bs\nu$ type we replace every factor by the term of lowest degree of its corner part,
\begin{align*}
\mathfrak G\mapsto\mathfrak g(\mb 0),\qquad
\mathfrak G_0\mapsto(D\mathfrak g)(\mb 0)\bs\nu_{\rm c},\qquad
\bs\nu\mapsto\bs\nu_{\rm c},\qquad
u_y\mapsto u^{\deg 1}_y,\qquad
\mb u\mapsto\mb u^{\deg 1},\qquad
p\mapsto p-p^{[2]},
\end{align*}
where $\mathfrak G_0$ is the factor of the pure geometry terms, for which $\mathfrak g(\mb 0)=0$, and where $u_1(\mb 0)\mb e_1$ is omitted: by the tables after \eqref{schematic bulk} and \eqref{schematic boundary}, a velocity factor $\grad^jZ$ with $j=0$ has $Z=u_y$, while $\grad\mb e_1=\mb 0$ and $\mb e_1\cdot\mb e_2=0$. This gives the corner part of $(\mb N,\mb N_\partial)$,
\begin{align}
\mb F^{[3]}_{\rm c}&:=\sum\mathfrak g(\mb 0)\cdot(\grad^{k_1}\bs\nu_{\rm c})\cdot\grad^jZ^{\deg 1}+\sum\mathfrak g(\mb 0)\cdot\bs\nu_{\rm c}\partial_x(p-p^{[2]}),\label{nonlinear corner part}\\
\mb f^{[3]}_{\rm c}&:=\bigg((D\mathfrak g)(\mb 0)\bs\nu_{\rm c}\cdot(r^{-1}\dr\nu_{{\rm c},1})\mb e+(D\mathfrak g)(\mb 0)\bs\nu_{\rm c}\cdot(r^{-1}\dr\nu_{{\rm c},2})\mb e\nonumber\\
&\qquad\quad+\sum\mathfrak g(\mb 0)\cdot(\grad^k\bs\nu_{\rm c})\cdot\grad^jZ^{\deg 1}+\sum\mathfrak g(\mb 0)\cdot\bs\nu_{\rm c}(p-p^{[2]})\bigg)\bigg|_{\partial_1\Omega},\nonumber
\end{align}
where the sums with $Z^{\deg 1}$ run over the summands of \eqref{schematic bulk} with one factor of $\bs\nu$ type and of the second sum of \eqref{schematic boundary} with $Z\in\{u_y,\mb u\}$, $Z^{\deg 1}\in\{u^{\deg 1}_y,\mb u^{\deg 1}\}$ being the image of $Z$, and the sums with $p-p^{[2]}$ run over those with $Z=p$. By the same tables, $k_1+j=2$ and $k+j=1$ in the sums with $Z^{\deg 1}$. Since $\grad^k\bs\nu_{\rm c}\in\mathcal Q^2_{1-k}$, $\grad^jZ^{\deg 1}\in\mathcal Q^1_{1-j}$, $\partial_x(p-p^{[2]})\in\mathcal Q^1_{-1}$ and $r^{-1}\dr\bs\nu_{\rm c}\in\mathcal Q^2_0$, the terms of \eqref{nonlinear corner part} satisfy
\begin{align*}
(\grad^{k_1}\bs\nu_{\rm c})\grad^jZ^{\deg 1},\ \bs\nu_{\rm c}\partial_x(p-p^{[2]})\in\mathcal Q^2_0,\qquad
(\grad^k\bs\nu_{\rm c})\grad^jZ^{\deg 1},\ \bs\nu_{\rm c}(p-p^{[2]}),\ \bs\nu_{\rm c}(r^{-1}\dr\bs\nu_{\rm c})\in\mathcal Q^3_1,
\end{align*}
so that $\mb F^{\smash[t]{[3]}}_{\rm c}\in\mathcal Q_0^2$ and $\mb f^{\smash[t]{[3]}}_{\rm c}\in\mathcal Q_1^3$, the latter with $\theta$-independent coefficients, as \eqref{third rung data classes} requires. Only $\bs\nu_{\rm c}$, its derivatives and $p-p^{[2]}$ carry $\ln r$, each with a coefficient proportional to $\smash[b]{c_{\phi_0,\dr_{-1}^2}}$ by \eqref{block geometry} and the display of $p-p^{[2]}$, whereas $\mb u^{\deg 1}$ carries none. So the terms $r\ln^2r$ of $\mb f^{\smash[t]{[3]}}_{\rm c}$ come from $\bs\nu_{\rm c}(r^{-1}\dr\bs\nu_{\rm c})$ and $\bs\nu_{\rm c}(p-p^{[2]})$, with coefficients $\smash[b]{c_{\phi_0,\dr_{-1}^2}^2}$ times constants depending only on $\theta_b$.

Every other product of corner parts has a higher degree. The summands of \eqref{schematic bulk} with two factors of $\bs\nu$ type give $(\grad^{k_1}\bs\nu_{\rm c})(\grad^{k_2}\bs\nu_{\rm c})\grad^jZ^{\deg 1}\in\mathcal Q^3_{3-k_1-k_2-j}=\mathcal Q^3_1$, and replacing $\mathfrak g(\mb 0)$ by $(D\mathfrak g)(\mb 0)\bs\nu_{\rm c}\in\mathcal Q_1^2$, $\mb u^{\deg 1}\in\mathcal Q_1^1$ by $\mb u^{\trimtop{0.3pt}{[3]}}_{\rm c}\in\mathcal Q_2^3$, or $p-p^{[2]}\in\mathcal Q_0^2$ by $p^{\trimtop{0.3pt}{[3]}}_{\rm c}\in\mathcal Q_1^3$ raises the degree by one. So the products of corner parts not in \eqref{nonlinear corner part} have degree at least one in the bulk and at least two on the boundary, and $(\mb F^{\trimtop{0.5pt}{[3]}}_{\rm c},\mb f^{\trimtop{0.5pt}{[3]}}_{\rm c})$ contains neither $\mb u^{\trimtop{0.5pt}{[3]}}_{\rm c}$ nor $p^{\trimtop{0.5pt}{[3]}}_{\rm c}$. The coefficients $\smash[b]{c_{\phi_0,\dr_{-2}^k}}$, $k\in\{1,2,3\}$, of $\phi_0$ enter only through $\bs\nu^{\smash[t]{[2]}}=\bs\nu[\phi_0^{\smash[t]{[2]}}]$, so $(\mb F^{\smash[t]{[3]}}_{\rm c},\mb f^{\smash[t]{[3]}}_{\rm c})$ does not depend on them.

By \eqref{geometry corner}, $\bs\nu_{\rm c}$ is linear in $(\smash[b]{c_{\phi_0,\dr_{-1}}},\smash[b]{c_{\phi_0,\dr_{-1}^2}})$, and $\mb u^{\deg 1}$ and $p-p^{[2]}$ are linear in $(\smash[b]{c_{\phi_0,\dr_{-1}}},\smash[b]{c_{\phi_0,\dr_{-1}^2}},u_1(\mb 0))$. Every term of \eqref{nonlinear corner part} is a constant times the product of a derivative of $\bs\nu_{\rm c}$ with a derivative of $\bs\nu_{\rm c}$, $\mb u^{\deg 1}$ or $p-p^{\smash[t]{[2]}}$. So the coefficients of $\mb F^{\smash[t]{[3]}}_{\rm c}$ and $\mb f^{\smash[t]{[3]}}_{\rm c}$ are quadratic forms in $(\smash[b]{c_{\phi_0,\dr_{-1}}},\smash[b]{c_{\phi_0,\dr_{-1}^2}},u_1(\mb 0))$, and by \eqref{profile coefficients} and \eqref{corner coefficients 2}
\begin{align*}
\|\mb F^{[3]}_{\rm c}\|_{\mathcal Q_0^2}+\|\mb f^{[3]}_{\rm c}\|_{\mathcal Q_1^3}\lesssim_{\theta_b}(|c_{\phi_0,\dr_{-1}}|+|c_{\phi_0,\dr_{-1}^2}|)(|c_{\phi_0,\dr_{-1}}|+|c_{\phi_0,\dr_{-1}^2}|+|u_1(\mb 0)|)\lesssim_{\theta_b}\mathcal D(\mathcal D+\mathcal E).
\end{align*}

We write $\mb N_\partial=H|_{\partial_1\Omega}$ and $\mb f^{\smash[t]{[3]}}_{\rm c}=H_{\rm c}|_{\partial_1\Omega}$, with $H$ and $H_{\rm c}$ the bulk expressions in \eqref{schematic boundary} and \eqref{nonlinear corner part}. Since $\|\mb f^{[3]}\|_{\partial_1\Omega,5/2-\kappa,5/2}\lesssim_{\theta_b}\|H-H_{\rm c}\|_{\Omega,3-\kappa,2,1}$ by Lemma \ref{trace}, it remains to prove
\begin{align}\label{remainder bulk target}
\|\mb F^{[3]}\|_{\Omega,2-\kappa,2}+\|H-H_{\rm c}\|_{\Omega,3-\kappa,2,1}\lesssim_{\theta_b}\mathcal D_+\brac{\mathcal D_++\mathcal E}+\mathcal D\mathfrak X.
\end{align}

\emph{Step 2: the remainders on $\{r>1\}$.} Since $r^{-(2-\kappa)}\leq r^{-(1-\kappa)}$ and $r^{-(3-\kappa)}\leq r^{-(2-\kappa)}$ on $\{r>1\}$, \eqref{N estimate} and \eqref{boundary bulk target} at $M=2$ give
\begin{align*}
&\|\mb F^{[3]}\|_{\Omega\cap\{r>1\},2-\kappa,2}+\|H-H_{\rm c}\|_{\Omega\cap\{r>1\},3-\kappa,2,1}\\*
&\qquad\leq\|\mb N\|_{\Omega,1-\kappa,l+2}+\|H\|_{\Omega,2-\kappa,l+2,1}+\|\mb F^{[3]}_{\rm c}\|_{\Omega\cap\{r>1\},2-\kappa,2}+\|H_{\rm c}\|_{\Omega\cap\{r>1\},3-\kappa,2,1}\\*
&\qquad\lesssim_{\theta_b}\mathcal D(\mathcal D+\mathcal E),
\end{align*}
where $\mb F^{\smash[t]{[3]}}_{\rm c}$ and $H_{\rm c}$, of degree zero and one, are bounded on $\{r>1\}$ by their coefficients, which are bounded as in Step 1, since $\int_1^\infty r^{-3+2\kappa}(1+\ln r)^4\d r<\infty$.

\emph{Step 3: the remainders on $\{r<1\}$.} We insert \eqref{factor splitting} into the summands of \eqref{schematic bulk}--\eqref{schematic boundary}, but we do not split the geometry factors of the products with $\mb u^{[3]}_{\rm c}$, $\mb u^{[3]}$, $p^{[3]}_{\rm c}$ or $p^{[3]}$: these products are linear in $\mathfrak X$ and must carry the factor $\mathcal D$ of the term $\mathcal D\mathfrak X$ of \eqref{N third rung estimate}, whereas $\bs\nu^{[2]}$ is only bounded by $\mathcal D_+$. For the summand $\bs\nu\grad^2\mb u$ of Case 1 of the proof of Proposition \ref{nonlinear estimates}, with $\mathfrak G=1$, for instance, $\grad^2(u_1(\mb 0)\mb e_1)=\mb 0$ and
\begin{align*}
\bs\nu\grad^2\mb u=\bs\nu_{\rm c}\grad^2\mb u^{\deg 1}+\bs\nu^{[2]}\grad^2\mb u^{\deg 1}+\bs\nu\grad^2\mb u^{[3]}_{\rm c}+\bs\nu\grad^2\mb u^{[3]},
\end{align*}
where $\bs\nu_{\rm c}\grad^2\mb u^{\deg 1}$ is a term of $\mb F^{\smash[t]{[3]}}_{\rm c}$, $\bs\nu^{\smash[t]{[2]}}\grad^2\mb u^{\deg 1}$ is bounded by the first placement below, $\bs\nu\grad^2\mb u^{\smash[t]{[3]}}_{\rm c}$ with the products with $\mb u^{\smash[t]{[3]}}_{\rm c}$, and $\bs\nu\grad^2\mb u^{\smash[t]{[3]}}$ by the second placement.

Besides $(\mb F^{\trimtop{0.3pt}{[3]}}_{\rm c},\mb f^{\trimtop{0.3pt}{[3]}}_{\rm c})$ the expansion gives three kinds of products, which we bound in turn: products of corner parts of degree at least one in the bulk and two on the boundary; products with $\mb u^{\trimtop{0.5pt}{[3]}}_{\rm c}$ or $p^{\trimtop{0.5pt}{[3]}}_{\rm c}$; and products with $\mb u^{[3]}$, $p^{[3]}$ or one of the remainders $\bs\nu^{[2]}$, $(D\mathfrak g)(\mb 0)\bs\nu^{[2]}$ and $\nu_1^2\mathfrak G+\nu_1\nu_2\mathfrak G+\nu_2^2\mathfrak G$ of the geometry. The remainder $\nu_1^2\mathfrak G+\nu_1\nu_2\mathfrak G+\nu_2^2\mathfrak G$ is split further: by the Leibniz rule and the second identity of \eqref{composition},
\begin{align*}
(r\grad)^m\brac{\nu_1^2\mathfrak G+\nu_1\nu_2\mathfrak G+\nu_2^2\mathfrak G}=\sum\mathfrak Gf_1\cdots f_J,\qquad J\geq 2,
\end{align*}
where each of $f_1,\dots,f_J$ is a component of some $(r\grad)^{m'}\bs\nu=(r\grad)^{m'}\bs\nu_{\rm c}+(r\grad)^{m'}\bs\nu^{[2]}$.

The targets $\|\mb F^{[3]}\|_{\Omega,2-\kappa,2}$ and $\|H-H_{\rm c}\|_{\Omega,3-\kappa,2,1}$ and the norms of $\bs\nu^{[2]}$ are taken on $\{r<1\}$, where $\|f\|_{\Omega\cap\{r<1\},a}\leq\|f\|_{\Omega\cap\{r<1\},a'}$ for $a\leq a'$. For $\bs\nu^{[2]}$ this is necessary: by Lemma \ref{corner splitting of the geometry}(iii), $\|\bs\nu^{[2]}\|_{\Omega\cap\{r>1\},a}=\infty$ for every $a\leq 2$ unless $\smash[b]{c_{\phi_0,\dr_{-1}}}=\smash[b]{c_{\phi_0,\dr_{-1}^2}}=0$.

\emph{Products of corner parts.} A product $P$ of corner parts of degree $d$, with $(\ln r)^k$ the highest power of the logarithm in it, satisfies
\begin{align}\label{pointwise corner bound}
|\dr^{j'}P|+|\dr^{j'}r\grad P|\lesssim_{\theta_b}\mathcal D(\mathcal D+\mathcal E)r^d(1+|\ln r|)^k\qquad\text{on }\{r<1\},
\end{align}
its coefficients being bounded as in Step 1. So $\|P\|_{\Omega\cap\{r<1\},2-\kappa,2}$ for $d\geq 1$ and $\|P\|_{\Omega\cap\{r<1\},3-\kappa,2,1}$ for $d\geq 2$ are $\lesssim_{\theta_b}\mathcal D(\mathcal D+\mathcal E)$, since
\begin{align*}
\|r(1+|\ln r|)^k\|_{\Omega\cap\{r<1\},2-\kappa}^2=\|r^2(1+|\ln r|)^k\|_{\Omega\cap\{r<1\},3-\kappa}^2=\theta_b\int_0^1r^{-1+2\kappa}(1+|\ln r|)^{2k}\d r<\infty.
\end{align*}
The same bound, with the derivatives already expanded, holds for the products of a term $\mathfrak Gf_1\cdots f_J$ of the display above in which $f_1,\dots,f_J$ are corner parts with derivatives of the other corner parts of the summand, $\mathfrak G$ being bounded by Lemma \ref{composition of the geometry}(i); for them $d\geq 2$.

\emph{Products with $\mb u^{\smash[t]{[3]}}_{\rm c}$ or $p^{\smash[t]{[3]}}_{\rm c}$.} Here the geometry is placed in the sup norm \eqref{geometry sup} at $a=2-\kappa/2$, where it is $\lesssim_{\theta_b}\mathcal D$ by Lemma \ref{geometry estimates}(ii); for instance, for $j_1+j_2\leq 2$,
\begin{align*}
\|(\dr^{j_1}\bs\nu)\dr^{j_2}\grad^2\mb u^{[3]}_{\rm c}\|_{\Omega\cap\{r<1\},2-\kappa}
&\leq\|r^{{\kappa\over 2}-1}\dr^{j_1}\bs\nu\|_{L^\infty(\Omega)}\|r^{1-{\kappa\over 2}}\dr^{j_2}\grad^2\mb u^{[3]}_{\rm c}\|_{\Omega\cap\{r<1\},2-\kappa}\\
&\lesssim_{\theta_b}\mathcal D\|\mb u^{[3]}_{\rm c}\|_{\mathcal Q_2^3},
\end{align*}
since $\int_0^1r^{-1+\kappa}(\ln r)^{2k}\d r<\infty$. When the geometry carries $(r\grad)^m\bs\nu$ with $m>l+2$, beyond the counts $m\leq l+2$ for which Lemma \ref{geometry estimates}(ii) bounds \eqref{geometry sup} by $\mathcal D$, the factor $\grad^j\mb u^{\smash[t]{[3]}}_{\rm c}$ has $j\leq 1$, and $|\dr^{j_2}\grad^j\mb u^{\smash[t]{[3]}}_{\rm c}|\lesssim_{\theta_b}r^{2-j}(1+|\ln r|)^2\mathfrak X$, since $H^2([0,\theta_b])\subset C^1([0,\theta_b])$. Then the geometry is placed in \eqref{geometry bulk} at $a=2-\kappa/2$: with $k_1+j=2$,
\begin{align*}
\|(\dr^{j_1}\grad^{k_1}\bs\nu)\dr^{j_2}\grad^j\mb u^{[3]}_{\rm c}\|_{\Omega\cap\{r<1\},2-\kappa}\lesssim_{\theta_b}\mathfrak X\|\bs\nu\|_{\Omega,2-{\kappa\over 2},0,j_1+k_1}\lesssim_{\theta_b}\mathcal D\mathfrak X,
\end{align*}
and likewise for $p^{\smash[t]{[3]}}_{\rm c}$, with $|\dr^{j_2}\grad^jp^{\smash[t]{[3]}}_{\rm c}|\lesssim_{\theta_b}r^{1-j}(1+|\ln r|)^2\mathfrak X$ for $j\leq 1$.

\emph{The remaining products.} These contain $\mb u^{[3]}$, $p^{[3]}$ or a remainder of the geometry. On them we run the proof of Proposition \ref{nonlinear estimates} with the targets $\|\ph\|_{\Omega\cap\{r<1\},2-\kappa,2}$ and $\|\ph\|_{\Omega\cap\{r<1\},3-\kappa,2,1}$ in place of $\|\ph\|_{\Omega,1-\kappa,l+2}$ and $\|\ph\|_{\Omega,2-\kappa,l+2,1}$. In every placement one factor is a remainder and carries one more power of $r$ at the corner, and the corner factors are placed by \eqref{corner product bound} below: a corner factor of degree $d$ raises the weight by $d$, less a margin for the logarithm, since $r^\iota|\ln r|^k$ is bounded on $\{r<1\}$ for $\iota>0$,
\begin{align}\label{corner product bound}
\|r^d(\ln r)^kf\|_{\Omega\cap\{r<1\},a+d-\iota}\lesssim_{\iota,k}\|f\|_{\Omega\cap\{r<1\},a},
\end{align}
with $\iota=0$ allowed when $k=0$. Likewise, by \eqref{block geometry}, $|\dr^{j'}(r\grad)^{m'}\bs\nu_{\rm c}|\lesssim_{\theta_b,\iota}(|\smash[b]{c_{\phi_0,\dr_{-1}}}|+|\smash[b]{c_{\phi_0,\dr_{-1}^2}}|)r^{1-\iota}$ on $\{r<1\}$. We fix $\iota\in(0,\iota_1)$. Apart from the products of two remainders in the pure geometry terms, only the factor carrying the extra power of $r$ is placed above the weights of the proof of Proposition \ref{nonlinear estimates}. Every other remainder of the geometry stays where that proof places $\bs\nu$, at a weight $a\leq 2-\kappa/2$, where its norms on $\{r<1\}$ are $\lesssim_{\theta_b}\mathcal D$, since $\bs\nu^{[2]}=\bs\nu-\bs\nu_{\rm c}$, by Lemma \ref{geometry estimates}(ii), \eqref{profile coefficients} and the bound on $\bs\nu_{\rm c}$ above with $\kappa/4$ in place of $\iota$. Three placements are new.

First, where the proof of Proposition \ref{nonlinear estimates} places $\bs\nu$ in $\|\ph\|_{\Omega,a,0,m}$ or in $\|r^{1-a}(r\grad)^{m'}\ph\|_{L^\infty(\Omega)}$ with $a\in[1-\kappa,2-\kappa]$, we place $\bs\nu^{\smash[t]{[2]}}$ one weight higher, using, for $m\leq l+4$ and $m'\leq l+2$,
\begin{align*}
\|\bs\nu^{[2]}\|_{\Omega\cap\{r<1\},a+1,0,m}&\leq\|\bs\nu^{[2]}\|_{\Omega,3-\kappa,0,m}\lesssim_{\theta_b}\mathcal D_+,\\
\|r^{-a}(r\grad)^{m'}\bs\nu^{[2]}\|_{L^\infty(\Omega\cap\{r<1\})}&\leq\|r^{-(2-\kappa)}(r\grad)^{m'}\bs\nu^{[2]}\|_{L^\infty(\Omega)}\lesssim_{\theta_b}\mathcal D_+,
\end{align*}
by \eqref{geometry remainder bulk} and \eqref{geometry remainder sup} at $\iota=0$. No loss of weight can be afforded here, and smallness is of no help: the other factor may be $\grad^2\mb u^{\deg 1}$, whose fields $\mb P\in\{\smash[b]{\mb P_{\mb u_{F_{0,0}},\dr_{-1}}},\smash[b]{\mb P_{\mb u_{F_{1,1}},\dr_{-1}}}\}$ have the small coefficients $u_1(\mb 0)/2$ and $(\sin\theta_b)^2\smash[b]{c_{\phi_0,\dr_{-1}^2}}/4$ but carry no extra power of $r$ at the corner, $|\dr^{j_2}\grad^2\mb P|\lesssim_{\theta_b}r^{-1}$. So, for $j_1+j_2\leq 2$,
\begin{align*}
\|(\dr^{j_1}\bs\nu^{[2]})\dr^{j_2}\grad^2\mb P\|_{\Omega\cap\{r<1\},2-\kappa}\lesssim_{\theta_b}\|\dr^{j_1}\bs\nu^{[2]}\|_{\Omega\cap\{r<1\},3-\kappa}\lesssim\|\bs\nu^{[2]}\|_{\Omega,3-\kappa,0,2}\lesssim_{\theta_b}\mathcal D_+,
\end{align*}
and $\|\bs\nu^{\smash[t]{[2]}}\grad^2\mb P\|_{\Omega\cap\{r<1\},2-\kappa}$ is finite exactly when $\|\bs\nu^{\smash[t]{[2]}}\|_{\Omega\cap\{r<1\},3-\kappa}$ is. This finiteness is what the second condition of \eqref{third rung angle condition} is for, and that condition cannot be weakened by the argument, since for $(3-\kappa)\theta_b\geq\pi/2$ the geometry carries, for generic data, the term $r^{\pi/(2\theta_b)-1}$ of Remark \ref{sharpness of the second angle condition}, whose product with $\grad^2\mb P$ lies outside the third-rung classes.

Second, $\mb u^{[3]}$ and $p^{[3]}$ are placed in their norms $\|\mb u^{[3]}\|_{\Omega,4-\kappa,2,2}$ and $\|p^{[3]}\|_{\Omega,3-\kappa,2,1}$ in $\mathfrak X$, where the proof of Proposition \ref{nonlinear estimates} uses $\|\mb u^{[2]}\|_{\Omega,3-\kappa,l+2,2}$ and $\|p^{[2]}\|_{\Omega,2-\kappa,l+2,1}$. In Case 2, for instance, where both radial derivatives of $\|\ph\|_{\Omega,2-\kappa,2}$ fall on $u^{[3]}_y$,
\begin{align*}
\|(\grad^2\bs\nu)\dr^2u^{[3]}_y\|_{\Omega\cap\{r<1\},2-\kappa}
\leq\|r^{-(3-\kappa)}\dr^2u^{[3]}_y\|_{L^\infty(\Omega)}\|\grad^2\bs\nu\|_{\Omega,-1}
\lesssim_{\theta_b}\mathcal D\|\mb u^{[3]}\|_{\Omega,4-\kappa,2,2}.
\end{align*}
These terms are linear in $\mathfrak X$ and carry the factor $\mathcal D$, which gives the term $\mathcal D\mathfrak X$ of \eqref{N third rung estimate}.

Third, in the pure geometry terms $\mathfrak G_0\cdot(r^{-1}\dr\nu_1)\mb e$ and $\mathfrak G_0\cdot(r^{-1}\dr\nu_2)\mb e$ we split $\mathfrak G_0$ by \eqref{factor splitting}, where $\mathfrak g(\mb 0)=0$, and $\bs\nu$ by Lemma \ref{corner splitting of the geometry}(i). Apart from the terms of $H_{\rm c}$ and the products with $\nu_1^2\mathfrak G+\nu_1\nu_2\mathfrak G+\nu_2^2\mathfrak G$, this gives the products $r^{-1}AB$ with $A\in\{(D\mathfrak g)(\mb 0)\bs\nu_{\rm c},(D\mathfrak g)(\mb 0)\bs\nu^{[2]}\}$ and $B$ a component of $\dr\bs\nu_{\rm c}$ or of $\dr\bs\nu^{[2]}$, not both corner parts, and their derivatives. In a mixed product the corner factor is placed in the bound on $|\dr^{j'}(r\grad)^{m'}\bs\nu_{\rm c}|$ above, where its logarithm costs the margin $\iota$, and the remainder factor in $\|\ph\|_{\Omega,3-\kappa+\iota,0,m}$, available by \eqref{geometry remainder bulk} with $-\iota$ in place of $\iota$, since $\iota<\iota_1$; with $A=(D\mathfrak g)(\mb 0)\bs\nu_{\rm c}$ and $B$ a component of $\dr\bs\nu^{[2]}$, for instance,
\begin{align*}
\|r^{-1}AB\|_{\Omega\cap\{r<1\},3-\kappa}\leq\|r^{-(1-\iota)}A\|_{L^\infty(\Omega\cap\{r<1\})}\|B\|_{\Omega,3-\kappa+\iota}\lesssim_{\theta_b}\mathcal D\mathcal D_+,
\end{align*}
by \eqref{profile coefficients} and \eqref{geometry remainder bulk}. In a product of two remainders the factor with fewer derivatives is placed in \eqref{geometry remainder sup} at $\iota=0$; for instance,
\begin{align*}
\|r^{-1}AB\|_{\Omega\cap\{r<1\},3-\kappa}\leq\|r^{-(2-\kappa)}A\|_{L^\infty(\Omega)}\|B\|_{\Omega\cap\{r<1\},2}\lesssim_{\theta_b}\mathcal D_+^2,
\end{align*}
since $\|B\|_{\Omega\cap\{r<1\},2}\leq\|B\|_{\Omega,3-\kappa}\lesssim_{\theta_b}\mathcal D_+$ by \eqref{geometry remainder bulk}.

(ii) We write $D\mb N_\partial=DH|_{\partial_1\Omega}$, with $DH$ the linearisation of $H$. Since $\|D\mb N_\partial\|_{\partial_1\Omega,3/2-\kappa,3/2}\lesssim_{\theta_b}\|DH\|_{\Omega,2-\kappa,1,1}$ by Lemma \ref{trace}, \eqref{DN second rung estimate} follows from
\begin{align}\label{DH target}
\|D\mb N\|_{\Omega,1-\kappa,1}+\|DH\|_{\Omega,2-\kappa,1,1}\lesssim_{\theta_b}\mathcal D\dot{\mathfrak X}+\mathcal D_+\dot{\mathcal E}+\brac{\mathcal D_++\mathcal E+\mathfrak X}\dot{\mathcal D}_+.
\end{align}
Every summand of $D\mb N$ and $DH$ is a summand of \eqref{schematic bulk}--\eqref{schematic boundary} with exactly one factor replaced by its linearisation: $\bs\nu$ by $\bs\nu[\dot\phi]$, $Z\in\{u_y,\mb u,p\}$ by $\dot Z\in\{\dot u_y,\dot{\mb u},\dot p\}$, or $\mathfrak G=\mathfrak g\circ\bs\nu$ by $\nu_1[\dot\phi](\partial_1\mathfrak g)(\bs\nu)+\nu_2[\dot\phi](\partial_2\mathfrak g)(\bs\nu)$, by Lemma \ref{composition of the geometry}(iii). The linearised factors split as in (i),
\begin{align*}
\dot{\mb u}&=\dot u_1(\mb 0)\mb e_1+(\dot{\mb u}^{[1]}-\dot{\mb u}^{[2]})+\dot{\mb u}^{[2]},&\dot{\mb u}^{[1]}-\dot{\mb u}^{[2]}&\in\mathcal Q_1^2,\\
\dot p&=(\dot p-\dot p^{[2]})+\dot p^{[2]},&\dot p-\dot p^{[2]}&\in\mathcal Q_0^3,\\
\bs\nu[\dot\phi]&=(\bs\nu[\dot\phi]-\bs\nu[\dot\phi^{[2]}])+\bs\nu[\dot\phi^{[2]}],&\bs\nu[\dot\phi]-\bs\nu[\dot\phi^{[2]}]&\in\mathcal Q_1^3,
\end{align*}
by \eqref{log decomposition} with $(\dot{\mb u},\dot p,\dot\phi)$ in place of $(\mb u,p,\phi)$, where $\dot{\mb u}(\mb 0)=\dot u_1(\mb 0)\mb e_1$ by Lemma \ref{decomposition estimates}(iii) and $\dot{\mb u}^{[1]}=\dot{\mb u}-\dot{\mb u}(\mb 0)$, and by Lemma \ref{corner splitting of the geometry}(i) at $\iota=0$, applied to $\dot\phi$ and to $\dot\phi^{[2]}$ given by \eqref{triple subtracted profile} with $\dot\phi$ in place of $\phi$, whose norms $\|\dot\phi\|_{\partial_1\Omega,1/2-\kappa}$ and $\|\dot\phi^{[2]}\|_{\partial_1\Omega,5/2-\kappa}$ are finite since $\dot{\mathcal D}_+<\infty$.
Every summand of $DH$ has a factor vanishing at the corner: $\bs\nu$ or $\bs\nu[\dot\phi]$ without derivative, $u_y$ or $\dot u_y$, by Lemma \ref{decomposition estimates}(iii), or $\mathfrak G_0$ or its linearisation, which carries a component of $\bs\nu[\dot\phi]$ by Lemma \ref{composition of the geometry}(iii). Hence, by the count of Step 1 of (i), the products of corner parts in $DH$ have degree at least one and those in $D\mb N$ degree at least zero, so that they are integrable on $\{r<1\}$ at the weights $2-\kappa$ and $1-\kappa$ of \eqref{DH target}, as in Step 3 of (i). Their coefficients are $\lesssim_{\theta_b}\mathcal D\dot{\mathcal E}+(\mathcal D+\mathcal E+\mathfrak X)\dot{\mathcal D}_+$. Indeed, those of the corner parts of the linearised factors are $\dot u_1(\mb 0)$, with $|\dot u_1(\mb 0)|\lesssim_{\theta_b}\dot{\mathcal E}$ by \eqref{corner coefficients 1}, and $\smash[b]{c_{\dot\phi,\dr_{-1}^k}}$, $k\in\{1,2,3\}$, which are $\lesssim_{\theta_b}\dot{\mathcal D}_+$, for $k=1$ by Lemma \ref{decomposition estimates}(i) applied to $\dot\phi-\dot\phi^{[2]}$; those of the other corner parts are bounded as in Step 1 and by $\mathfrak X$.

We prove \eqref{DH target} as \eqref{DN estimate}, at $M=2$. If the linearised factor is a solution factor, the placements of (i) apply with $(\dot{\mb u},\dot p)$ in place of $(\mb u,p)$, since $\dot u_y(\mb 0)=0$ by Lemma \ref{decomposition estimates}(iii). Now $\dot{\mb u}^{[2]}$ and $\dot p^{[2]}$ are placed in their norms $\|\dot{\mb u}^{[2]}\|_{\Omega,3-\kappa,1,2}$ and $\|\dot p^{[2]}\|_{\Omega,2-\kappa,1,1}$ in $\dot{\mathfrak X}$ and produce the terms $\mathcal D\dot{\mathfrak X}$; in the analogue of Case 2, for $j_1+j_2\leq 1$,
\begin{align*}
\|(\dr^{j_1}\grad^2\bs\nu)(\dr^{j_2}\dot u^{[2]}_y)\|_{\Omega,1-\kappa}
\leq\|r^{-(2-\kappa)}\dr^{j_2}\dot u^{[2]}_y\|_{L^\infty(\Omega)}\|\dr^{j_1}\grad^2\bs\nu\|_{\Omega,-1}
\lesssim_{\theta_b}\mathcal D\|\dot{\mb u}^{[2]}\|_{\Omega,3-\kappa,1,2}.
\end{align*}
The coefficients of $\dot{\mb u}^{\smash[t]{[1]}}-\dot{\mb u}^{\smash[t]{[2]}}$ and of $\dot p-\dot p^{\smash[t]{[2]}}$, bounded as above, and the first placement of (i) produce the terms $\mathcal D_+(\dot{\mathcal E}+\dot{\mathcal D}_+)$.

If the linearised factor is a geometry factor, $\bs\nu[\dot\phi]$ is placed by Lemma \ref{geometry estimates}(i) at weights $a\leq 2-\iota$, where the RHSs of \eqref{geometry bulk} and \eqref{geometry sup} are bounded by $\|\dot\phi\|_{\partial_1\Omega,a-1/2,5/2}$ with $a-1/2\leq 3/2-\iota$. We split $\dot\phi$ at $r=1$. On $\{r<1\}$, $\dot\phi-\dot\phi^{[2]}=\smash[b]{c_{\dot\phi,\dr_{-1}}}r+\smash[b]{c_{\dot\phi,\dr_{-1}^2}}r\ln r+\smash[b]{c_{\dot\phi,\dr_{-1}^3}}r\ln^2r/2$ by \eqref{triple subtracted profile}, and $\int_0^1r^{-1+2\iota}(1+|\ln r|)^4\d r<\infty$, while the norm $\|\dot\phi^{[2]}\|_{\partial_1\Omega,5/2-\kappa,5/2}$ in $\dot{\mathcal D}_+$ dominates; on $\{r>1\}$ the norm $\|\dot\phi\|_{\partial_1\Omega,-1/2-\kappa,l+5/2}$ in $\dot{\mathcal D}$ dominates. The other factors are placed as in (i), the geometry factors other than $\bs\nu[\dot\phi]$ at weights $a\leq 2-\kappa/2$, and these terms are bounded by $(\mathcal D_++\mathcal E+\mathfrak X)\dot{\mathcal D}_+$; in the analogue of Case 2 with $u^{\deg 1}_y$ of \eqref{v_y splitting}, for instance, for $j'\leq 1$,
\begin{align*}
\|(\dr^{j'}\grad^2\bs\nu[\dot\phi])u^{\deg 1}_y\|_{\Omega,1-\kappa}
\lesssim_{\theta_b}\mathcal E\|\grad^2\bs\nu[\dot\phi]\|_{\Omega,-\kappa,1}
\lesssim_{\theta_b}\mathcal E\|\bs\nu[\dot\phi]\|_{\Omega,2-\kappa,0,3}
\lesssim_{\theta_b}\mathcal E\dot{\mathcal D}_+.
\end{align*}

The holomorphy across $\Re\lambda=1$ follows from \eqref{DN estimate} at $M=1$ and \eqref{DN second rung estimate}. By them $r^2D\mb N$ has finite norms $\|\ph\|_{\Omega,2-\kappa}$ and $\|\ph\|_{\Omega,3-\kappa}$, and $rD\mb N_\partial$ finite norms $\|\ph\|_{\partial_1\Omega,3/2-\kappa}$ and $\|\ph\|_{\partial_1\Omega,5/2-\kappa}$. So by Lemma \ref{block matching} with $m=0$, applied at each $\theta$ in the bulk, their Mellin transforms on the lines $\Re\lambda=1-\kappa$ and $\Re\lambda=2-\kappa$ are the restrictions of one function holomorphic on the strip between them.

(iii) The proofs of (i) and (ii) hold at the count $j$ in place of $2$, respectively $1$, unchanged. There the geometry factors enter through $\|\bs\nu\|_{\Omega,a,0,j+2}$, $\|\bs\nu^{[2]}\|_{\Omega,3-\kappa,0,j+2}$, $\|\bs\nu^{[2]}\|_{\Omega,3-\kappa+\iota,0,j+2}$ and $\|\bs\nu[\dot\phi]\|_{\Omega,a,0,j+2}$. Lemma \ref{geometry estimates}(ii) bounds $\|\bs\nu\|_{\Omega,a,0,j+2}$, and \eqref{geometry remainder bulk} the norms of $\bs\nu^{[2]}$, for $j+2\leq l+4$. Lemma \ref{geometry estimates}(iii) bounds $\|\bs\nu[\dot\phi]\|_{\Omega,a,0,j+2}$ at the weights $a\leq 2-\kappa$ for $j+2\leq l+3$, as in the proof of \eqref{DN estimate}; at the weights $a\in(2-\kappa,2-\iota]$ of (ii), $\|\dot\phi^{[2]}\|_{\partial_1\Omega,5/2-\kappa,j+3/2}$ in $\dot{\mathcal D}_+$ takes the place of $\|\dot\phi^{[2]}\|_{\partial_1\Omega,5/2-\kappa,5/2}$ on $\{r<1\}$. Every other placement holds at every count.

For the shifted estimates, Lemma \ref{trace} bounds the LHSs of \eqref{shifted N estimate} and \eqref{shifted DN estimate} by a constant times $\|\mb N\|_{\Omega,1+\tilde\kappa-\kappa,l+2}+\|H\|_{\Omega,2+\tilde\kappa-\kappa,l+2,1}$ and $\|D\mb N\|_{\Omega,1+\tilde\kappa-\kappa,j}+\|DH\|_{\Omega,2+\tilde\kappa-\kappa,j,1}$. On $\{r>1\}$ these are dominated by $\|\mb N\|_{\Omega,1-\kappa,l+2}+\|H\|_{\Omega,2-\kappa,l+2,1}\lesssim_{\theta_b}\mathcal D(\mathcal D+\mathcal E)$, by \eqref{N estimate} and \eqref{boundary bulk target} at $M=2$, and by $\|D\mb N\|_{\Omega,1-\kappa,j}+\|DH\|_{\Omega,2-\kappa,j,1}$, which \eqref{DH target} at the count $j$ bounds by the RHS of \eqref{shifted DN estimate}.

On $\{r<1\}$ we have $\|\mb F^{[3]}\|_{\Omega\cap\{r<1\},1+\tilde\kappa-\kappa,l+2}\leq\|\mb F^{[3]}\|_{\Omega,2-\kappa,l+2}$ and $\|H-H_{\rm c}\|_{\Omega\cap\{r<1\},2+\tilde\kappa-\kappa,l+2,1}\leq\|H-H_{\rm c}\|_{\Omega,3-\kappa,l+2,1}$, which \eqref{remainder bulk target} at the count $l+2$ bounds. The corner parts $\mb F^{\trimtop{0.3pt}{[3]}}_{\rm c}$ and $H_{\rm c}$, of degree zero and one, satisfy \eqref{pointwise corner bound}, and they are integrable against the shifted weights because $\tilde\kappa<\kappa$:
\begin{align*}
\|(1+|\ln r|)^k\|_{\Omega\cap\{r<1\},1+\tilde\kappa-\kappa}^2=\|r(1+|\ln r|)^k\|_{\Omega\cap\{r<1\},2+\tilde\kappa-\kappa}^2=\theta_b\int_0^1r^{2(\kappa-\tilde\kappa)-1}(1+|\ln r|)^{2k}\d r<\infty.
\end{align*}
This proves \eqref{shifted N estimate}.

For \eqref{shifted DN estimate} on $\{r<1\}$ we use the placements of the proof of \eqref{DN estimate} at $M=2$, shifted by $\tilde\kappa$, with $(\dot{\mb u},\dot p)$ split by \eqref{log decomposition}. Outside the pure geometry terms the shift is carried by one geometry factor, by $\bs\nu[\dot\phi]$ when it is present, while the other geometry factors stay at their weights in that proof. Against $\mb u^{[2]}$, $p^{[2]}$, $\dot{\mb u}^{[2]}$ and $\dot p^{[2]}$ this factor is placed at the weight $1+\tilde\kappa$ instead of $1$, within the ranges of Lemma \ref{geometry estimates}(ii)--(iii), since $1+\tilde\kappa<1+\kappa\leq 2-\kappa$. Against the terms subtracted in \eqref{data step relations} and \eqref{log decomposition} it is placed at $2-\kappa+\tilde\kappa$ instead of $2-\kappa$, plus a margin $(\kappa-\tilde\kappa)/2$ for a logarithm as in \eqref{corner product bound}, and in the pure geometry terms both factors are placed at $2-(\kappa-\tilde\kappa)/2$ instead of $2-\kappa/2$. These weights lie in $[1,2)$, and for $a\in[1,2)$, $m\leq j+2$ and $m'\leq j$,
\begin{align*}
\|\bs\nu\|_{\Omega\cap\{r<1\},a,0,m}+\|r^{1-a}(r\grad)^{m'}\bs\nu\|_{L^\infty(\Omega\cap\{r<1\})}&\lesssim_{\theta_b,a}\mathcal D_+,\\
\|\bs\nu[\dot\phi]\|_{\Omega\cap\{r<1\},a,0,m}+\|r^{1-a}(r\grad)^{m'}\bs\nu[\dot\phi]\|_{L^\infty(\Omega\cap\{r<1\})}&\lesssim_{\theta_b,a}\dot{\mathcal D}_+.
\end{align*}
Indeed, $\bs\nu=\bs\nu_{\rm c}+\bs\nu^{[2]}$, where $|\dr^{j'}(r\grad)^{m'}\bs\nu_{\rm c}|\lesssim_{\theta_b}\mathcal Dr(1+|\ln r|)$ on $\{r<1\}$ by \eqref{block geometry} and \eqref{profile coefficients}, and $r^{2-a}(1+|\ln r|)$ is bounded on $\{r<1\}$ and $\int_0^1r^{3-2a}(1+|\ln r|)^2\d r<\infty$, while on $\{r<1\}$ the norms of $\bs\nu^{\smash[t]{[2]}}$ at the weight $a$ are dominated by those at $3-\kappa$, which \eqref{geometry remainder sup}--\eqref{geometry remainder bulk} at $\iota=0$ bound. The bound on $\bs\nu[\dot\phi]$ is that of (ii), at the count $j$ and with $2-a$ in place of $\iota$. In the analogues of Case 2 displayed in (ii), for instance, for $j_1+j_2\leq j$ and $j'\leq j$,
\begin{align*}
\|(\dr^{j_1}\grad^2\bs\nu)(\dr^{j_2}\dot u^{[2]}_y)\|_{\Omega\cap\{r<1\},1+\tilde\kappa-\kappa}
&\leq\|r^{-(2-\kappa)}\dr^{j_2}\dot u^{[2]}_y\|_{L^\infty(\Omega)}\|\dr^{j_1}\grad^2\bs\nu\|_{\Omega,\tilde\kappa-1}
\lesssim_{\theta_b}\mathcal D\|\dot{\mb u}^{[2]}\|_{\Omega,3-\kappa,j,2},\\
\|(\dr^{j'}\grad^2\bs\nu[\dot\phi])u^{\deg 1}_y\|_{\Omega\cap\{r<1\},1+\tilde\kappa-\kappa}
&\lesssim_{\theta_b}\mathcal E\|\bs\nu[\dot\phi]\|_{\Omega\cap\{r<1\},2+\tilde\kappa-\kappa,0,j+2}
\lesssim_{\theta_b,\tilde\kappa}\mathcal E\dot{\mathcal D}_+,
\end{align*}
by Lemma \ref{geometry estimates}(ii) at $a=1+\tilde\kappa$ and by the bound above at $a=2+\tilde\kappa-\kappa$.
\end{proof}

\begin{lemma}[Second-order corner values]\label{second-order corner values}
Let $(\dot{\mb u},\dot p)$ be a dotted pair with $\dot{\mathfrak X}<\infty$, with $\dot{\mathfrak X}$ as in Proposition \ref{nonlinear corner structure}(ii), for the decomposition of Proposition \ref{third rung}(i),
\begin{align}\label{dotted log decomposition}
\dot{\mb u}=\dot{\mb u}(\mb 0)+{\dot u_1(\mb 0)\over 2}\mb P_{\mb u_{F_{0,0}},\dr_{-1}}+{(\sin\theta_b)^2c_{\dot\phi,\dr_{-1}^2}\over 4}\mb P_{\mb u_{F_{1,1}},\dr_{-1}}+{(\sin\theta_b)^2c_{\dot\phi,\dr_{-1}^3}\over 4}\mb P_{\mb u_{F_{1,1}},\dr_{-1}^2}+\dot{\mb u}^{[2]}.
\end{align}
Then
\begin{align}\label{second-order corner identity}
\partial_r\dot u_\theta(r,\theta_b)&={\theta_b(\sin\theta_b)^2c_{\dot\phi,\dr_{-1}^3}\over 2}\ln r\nonumber\\
&\quad+\brac{-{\dot u_1(\mb 0)\over 2}+{\theta_b(\sin\theta_b)^2c_{\dot\phi,\dr_{-1}^2}\over 2}+{(\sin\theta_b)^2c_{\dot\phi,\dr_{-1}^3}\over 4}\brac{2\theta_b+\mb Q_\star(\theta_b)\cdot\mb e_\theta(\theta_b)}}+\partial_r\dot u^{[2]}_\theta(r,\theta_b),
\end{align}
where $\mb Q_\star$ is as in \eqref{log field} and $|\partial_r\dot{\mb u}^{[2]}|\lesssim_{\theta_b}\dot{\mathfrak X}r^{1-\kappa}$. The function $\ddot\phi=(\sin\theta_b)^{-2}\partial_r\dot u_\theta(\ph,\theta_b)$ extends continuously to the contact point with value $0$ if and only if
\begin{align}\label{C2 explicit}
c_{\dot\phi,\dr_{-1}^3}=0\qquad\text{ and }\qquad c_{\dot\phi,\dr_{-1}^2}={\dot u_1(\mb 0)\over\theta_b(\sin\theta_b)^2}.\tag{C2$'$}
\end{align}
In that case $(\sin\theta_b)^2\|\ddot\phi\|_{\partial_1\Omega,3/2-\kappa,3/2}\lesssim_{\theta_b}\dot{\mathfrak X}$, and if $\dot{\mathfrak X}$ is formed at the counts $(j,2)$ in place of $(1,2)$, the same bound holds at the count $j+{1\over 2}$ in place of ${3\over 2}$.
\end{lemma}
\begin{proof}
On $\theta=\theta_b$ the $\mb e_\theta$ components of $\smash[b]{\mb P_{\mb u_{F_{0,0}},\dr_{-1}}}$ and $\smash[b]{\mb P_{\mb u_{F_{1,1}},\dr_{-1}}}$ are $-r$ and $2\theta_br$, as in the proof of Lemma \ref{corner values}, and that of \eqref{log field} is $2\theta_br\ln r+(\mb Q_\star(\theta_b)\cdot\mb e_\theta(\theta_b))r$; applying $\partial_r$ to \eqref{dotted log decomposition} gives \eqref{second-order corner identity}, the constant $\dot{\mb u}(\mb 0)$ contributing nothing. The bound on $\partial_r\dot{\mb u}^{[2]}$ is Lemma \ref{embeddings}(ii) at $a=1-\kappa$ applied to $\grad\dot{\mb u}^{[2]}$, as in Lemma \ref{corner values}, using the counts $(1,2)$ of $\dot{\mathfrak X}$. The function $\ln r$ has no limit and no weighted norm at the origin, so the continuous extension by $0$ forces first the logarithmic and then the constant coefficient of \eqref{second-order corner identity} to vanish, and conversely; once the logarithmic coefficient vanishes, the vanishing of the constant one is the second equation of \eqref{C2 explicit}. The remaining statements follow as in Lemma \ref{corner values}, $\|\ddot\phi\|_{\partial_1\Omega,3/2-\kappa,3/2}$ being obtained from $\|\dr\dot u^{\smash[t]{[2]}}_\theta\|_{\partial_1\Omega,5/2-\kappa,3/2}$ by \eqref{dotphi weight shift} and Lemma \ref{trace}, at the count $j+{1\over 2}$ when $\dot{\mathfrak X}$ carries the counts $(j,2)$.
\end{proof}

The second equation of \eqref{C2 explicit} is the dotted \eqref{C1 in the linear theory}.

\subsubsection{The initial-time theorem}

The next theorem constructs from $\phi_0$ the pairs $(\mb u,p)$ and $(\dot{\mb u},\dot p)$ and the profiles $\dot\phi$ and $\ddot\phi$ derived from them; the triple $(\phi_0,\dot\phi,\ddot\phi)$ is the initial jet of the nonlinear problem. The condition \eqref{C2} of the theorem is the analogue, one time derivative up, of what \eqref{C1} is by Lemma \ref{corner values}: the finiteness of $\|\dot\phi\|_{\partial_1\Omega,3/2-\kappa,l+5/2}$. The positive-time existence theorem of Subsection \ref{positive-time construction} uses all three parts of the theorem.

\begin{theorem}[The initial jet of the nonlinear problem]\label{nonlinear initial data}
Let $\kappa\in(0,1/2]$, let $l\geq 1$ be an integer, assume \eqref{standing angle condition}, and let $\theta_b$ be small enough for Proposition \ref{stationary system} to apply at the derivative parameters $l+2$ and $l+1$. There exists $\delta>0$, depending only on $\theta_b,\kappa$ and $l$, such that if $\mathcal D\leq\delta$ then the following hold.
\begin{enumerate}
\item There exists a unique pair $(\mb u,p)$ with $\mathcal E(\mb u,p)<\infty$ satisfying a.e. the stationary subsystem of System \ref{reference wedge system}:
\begin{subequations}\label{system: stationary nonlinear}
\begin{alignat}{3}
\lpc\mb u-\grad p&=\mb N[\phi_0,\mb u,p]\qquad&&\text{on}\quad\Omega\\
\grad\cdot\mb u&=0\qquad&&\text{on}\quad\Omega\\
u_\theta=0,\qquad u_r-r^{-1}\partial_\theta u_r&=0\qquad&&\text{on}\quad\partial_0\Omega\\
\mathbb T[\mb u,p]\mb e_\theta-(\sin\theta_b)^2(\partial_r\phi_0)\mb e_\theta&=\mb N_\partial[\phi_0,\mb u,p]\qquad&&\text{on}\quad\partial_1\Omega.
\end{alignat}
\end{subequations}
Setting $\dot\phi:=(\sin\theta_b)^{-2}\partial_ru_\theta|_{\partial_1\Omega}$, we have $\mathcal E(\mb u,p)+\sum_{M=0}^1(\sin\theta_b)^2\|\dot\phi\|_{\partial_1\Omega,M-1/2-\kappa,l+5/2}\lesssim_{\theta_b}\mathcal D$.
\item Assume in addition that the pair constructed in (i) satisfies the first compatibility condition, that is, \eqref{C1 in the linear theory} for the profile $\phi_0$:
\begin{align}\label{C1}
c_{\phi_0,\dr_{-1}^2}={u_1(\mb 0)\over\theta_b(\sin\theta_b)^2}.\tag{C1}
\end{align}
Then $\dot\phi$ vanishes at the contact point, $\dot{\mathcal D}(\dot\phi)\lesssim_{\theta_b}\mathcal D$, and there exists a unique pair $(\dot{\mb u},\dot p)$ with $\dot{\mathcal E}(\dot{\mb u},\dot p)<\infty$ satisfying a.e.\ the time-differentiated system
\begin{subequations}\label{system: dotted}
\begin{alignat}{3}
\lpc\dot{\mb u}-\grad\dot p&=D\mb N[\phi_0,\mb u,p](\dot\phi,\dot{\mb u},\dot p)\qquad&&\text{on}\quad\Omega\\
\grad\cdot\dot{\mb u}&=0\qquad&&\text{on}\quad\Omega\\
\dot u_\theta=0,\qquad \dot u_r-r^{-1}\partial_\theta\dot u_r&=0\qquad&&\text{on}\quad\partial_0\Omega\\
\mathbb T[\dot{\mb u},\dot p]\mb e_\theta-(\sin\theta_b)^2(\partial_r\dot\phi)\mb e_\theta&=D\mb N_\partial[\phi_0,\mb u,p](\dot\phi,\dot{\mb u},\dot p)\qquad&&\text{on}\quad\partial_1\Omega.
\end{alignat}
\end{subequations}
Setting $\ddot\phi:=(\sin\theta_b)^{-2}\partial_r\dot u_\theta|_{\partial_1\Omega}$, the tuple $(\mb u,p,\dot{\mb u},\dot p,\dot\phi,\ddot\phi)$ satisfies
\begin{align}
&\sum_{M=0}^2\bigg((\sin\theta_b)^2\|\dot\phi\|_{\partial_1\Omega,M-{1\over 2}-\kappa,l+{5\over 2}}+\|\mb u^{[M]}\|_{\Omega,M+1-\kappa,l+2,2}+\|p^{[M]}\|_{\Omega,M-\kappa,l+2,1}\bigg)\nonumber\\
&\quad+\sum_{M=0}^1\bigg((\sin\theta_b)^2\|\ddot\phi\|_{\partial_1\Omega,M-{1\over 2}-\kappa,l+{3\over 2}}+\|\dot{\mb u}^{[M]}\|_{\Omega,M+1-\kappa,l+1,2}+\|\dot p^{[M]}\|_{\Omega,M-\kappa,l+1,1}\bigg)\nonumber\\
&\lesssim_{\theta_b}\mathcal D.\label{main estimate}
\end{align}
\item Assume in addition that the derived second profile satisfies the second compatibility condition, that is, \eqref{improved theorem C2} of Theorem \ref{linear solution full regularity} for the profile $\ddot\phi$:
\begin{align}\label{C2}
\|\ddot\phi\|_{\partial_1\Omega,{3\over 2}-\kappa,l+{3\over 2}}<\infty.\tag{C2}
\end{align}
Then $\ddot\phi$ vanishes at the contact point. If moreover $\kappa\in(1/4,1/2)$, then the triple $(\phi_0,\dot\phi,\ddot\phi)$ satisfies the hypotheses of Theorem \ref{linear solution full regularity}.
\end{enumerate}
\end{theorem}
\begin{proof}
Step 1 constructs $(\mb u,p)$ and Step 3 constructs $(\dot{\mb u},\dot p)$, both as fixed points of contractions; Step 2 derives the profile $\dot\phi$ and shows that \eqref{C1} makes it an admissible datum for Step 3, and Step 4 assembles parts (i) and (ii). A short paragraph proves part (iii).

\emph{Step 1: the pair $(\mb u,p)$.} Let $C_{\rm lin}$ denote the constant of Propositions \ref{stationary system} and \ref{corner ladder}, applied with their derivative parameter equal to $l+2$ for $(\mb u,p)$ and to $l+1$ for $(\dot{\mb u},\dot p)$, together with the trace constants of Lemma \ref{trace}; and let $C_{\rm nl}$ denote the constant of Proposition \ref{nonlinear estimates}. We take $\delta\leq\delta_0$ small enough that
\begin{align}\label{delta conditions}
C_{\rm lin}C_{\rm nl}\delta\leq{1\over 2}\qquad\text{ and }\qquad C_{\rm nl}\delta(1+2C_{\rm lin})\leq{1\over 2},
\end{align}
and small enough that Lemma \ref{regularised flattening} applies to $\chi'=-\phi_0$, using $\|\phi_0\|_{\partial_1\Omega,1/2,3/2}\lesssim_{\theta_b}\mathcal D$; then $\mb N[\phi_0,\ph,\ph]$ and $\mb N_\partial[\phi_0,\ph,\ph]$ are well defined and \eqref{geometry smallness} holds. Let $\rho:=2C_{\rm lin}\mathcal D$ and
\begin{align*}
B:=\set{(\mb v,q):\mathcal E(\mb v,q)\leq\rho,\ \grad\cdot\mb v=0,\ v_\theta|_{\partial_0\Omega}=0},
\end{align*}
and define $\mathcal T_1(\mb v,q):=(\mb u,p)$ to be the solution of \eqref{system: stationary linear} with data
\begin{align*}
(\phi,\mb F,\mb f)=\brac{\phi_0,\mb N[\phi_0,\mb v,q],\mb N_\partial[\phi_0,\mb v,q]}
\end{align*}
given by Proposition \ref{stationary system}, upgraded to the full ladder by Proposition \ref{corner ladder} for $M=1,2$, whose angle hypothesis $2-\kappa<\Lambda^*(\theta_b)$ is the first condition of \eqref{standing angle condition}. All the data norms required are available, by \eqref{N estimate} and $\mathcal D$, whose counts $l+2$, $l+5/2$ and $l+7/2$ are the required ones. Proposition \ref{corner ladder} is applied inductively in $M$: on the right of \eqref{um data reg master}, $\|u_r^{\trimtop{0.9pt}{[M-1]}}\|_{\partial_0\Omega,M-1/2-\kappa,l+5/2}\lesssim_{\theta_b}\|\mb u^{[M-1]}\|_{\Omega,M-\kappa,l+2,2}$ by Lemma \ref{trace}, a norm already bounded at the rung $M-1$. The set $B$ is not empty, since the solution of \eqref{system: stationary linear} with data $(\phi_0,\mb 0,\mb 0)$ has $\mathcal E\leq C_{\rm lin}\mathcal D\leq\rho$ by \eqref{data estimate} and \eqref{um data reg master}, and so lies in $B$.

By \eqref{N estimate} and \eqref{data estimate} we then have, for $(\mb v,q)\in B$,
\begin{align*}
\mathcal E\brac{\mathcal T_1(\mb v,q)}\leq C_{\rm lin}\brac{\mathcal D+C_{\rm nl}\mathcal D(\mathcal D+\rho)}=C_{\rm lin}\mathcal D\brac{1+C_{\rm nl}\mathcal D(1+2C_{\rm lin})}\leq{3\over 2}C_{\rm lin}\mathcal D\leq\rho
\end{align*}
by \eqref{delta conditions}, and hence $\mathcal T_1(B)\subseteq B$, since the outputs of Proposition \ref{stationary system} satisfy $\grad\cdot\mb u=0$ and $u_\theta|_{\partial_0\Omega}=0$. By \eqref{N difference}, the uniqueness statement of Proposition \ref{stationary system} and Remark \ref{metric remark}, the difference of two outputs is the solution with profile datum $0$ and forcing the difference of the forcings, so that
\begin{align*}
d\brac{\mathcal T_1(\mb v,q),\mathcal T_1(\mb v',q')}\leq C_{\rm lin}C_{\rm nl}\mathcal Dd\brac{(\mb v,q),(\mb v',q')}\leq{1\over 2}d\brac{(\mb v,q),(\mb v',q')}
\end{align*}
for $(\mb v,q),(\mb v',q')\in B$.

The set $B$ is complete for $d$. Weighted $L^2$ spaces are complete, so the decomposed pieces of a $d$ Cauchy sequence converge. Their limits are again the pieces of a single pair, because by Lemma \ref{decomposition estimates}(ii) the corner coefficients converge as well. The conditions $\grad\cdot\mb v=0$ and $v_\theta|_{\partial_0\Omega}=0$ pass to the limit, and the limit satisfies $\mathcal E\leq\rho$ by Fatou's lemma. The Banach fixed point theorem therefore yields a unique fixed point $(\mb u,p)\in B$. By construction it solves \eqref{system: stationary nonlinear} and satisfies $\mathcal E\leq 2C_{\rm lin}\mathcal D$; the bounds on $\dot\phi$ for $M\in\{0,1\}$ are \eqref{dotphi bound}.

The fixed point is the only solution of \eqref{system: stationary nonlinear} with $\mathcal E<\infty$, whatever its size. Indeed, let $(\mb u',p')$ be such a solution. It lies in the class of \eqref{data estimate} at the derivative parameter $l+2$, whose norms are finite by the rungs $M=0$ and $M=1$ of $\mathcal E$; for the unweighted term this reads, for every pair $(\mb v,q)$ and every count $j$, since $\grad\mb v=\grad\mb v^{[1]}$,
\begin{align}\label{unweighted by weighted}
\|r\grad\mb v\|_{\Omega,1,j,1}\sim\|\grad\mb v\|_{\Omega,0,j,1}\leq\|\grad\mb v\|_{\Omega,-\kappa,j,1}+\|\grad\mb v^{[1]}\|_{\Omega,1-\kappa,j,1}\lesssim\|\mb v\|_{\Omega,1-\kappa,j,2}+\|\mb v^{[1]}\|_{\Omega,2-\kappa,j,2},
\end{align}
whose RHS is finite at $j=l+2$ when $\mathcal E(\mb v,q)<\infty$. The map $\mathcal T_1$ is defined on every pair with finite $\mathcal E$ satisfying $\grad\cdot\mb v=0$ and $v_\theta|_{\partial_0\Omega}=0$, the data norms being finite by \eqref{N estimate}. The contraction estimate holds for any two such pairs, since \eqref{N difference} does not involve the radius $\rho$. By \eqref{system: stationary nonlinear}, $(\mb u',p')$ solves \eqref{system: stationary linear} with profile datum $\phi_0$ and forcing $(\mb N,\mb N_\partial)[\phi_0,\mb u',p']$, and so does $\mathcal T_1(\mb u',p')$, by the definition of $\mathcal T_1$. So, by the uniqueness statement of Proposition \ref{stationary system}, $\mathcal T_1(\mb u',p')=(\mb u',p')$. Hence $d((\mb u,p),(\mb u',p'))=d(\mathcal T_1(\mb u,p),\mathcal T_1(\mb u',p'))\leq{1\over 2}d((\mb u,p),(\mb u',p'))$, and $(\mb u',p')=(\mb u,p)$.

\emph{Step 2: compatibility and the derived profile.} Under \eqref{C1}, \eqref{corner identity} gives $\partial_ru_\theta(0,\theta_b)=0$, that is, $\dot\phi$ vanishes at the contact point. So Lemma \ref{corner values}, together with $\mathcal E\leq 2C_{\rm lin}\mathcal D$ from Step 1, gives
\begin{align}\label{P bound}
\dot{\mathcal D}(\dot\phi)\lesssim_{\theta_b}\mathcal E\lesssim_{\theta_b}\mathcal D.
\end{align}
In particular the norm $\|\dot\phi\|_{\partial_1\Omega,1/2,l+5/2}$, required by Proposition \ref{stationary system} at the derivative parameter $l+1$, is bounded by $\|\dot\phi\|_{\partial_1\Omega,1/2-\kappa,l+5/2}$ and $\|\dot\phi\|_{\partial_1\Omega,3/2-\kappa,l+5/2}$ in $\dot{\mathcal D}$.

\emph{Step 3: the pair $(\dot{\mb u},\dot p)$.} Let $\rho_2:=3C_{\rm lin}\dot{\mathcal D}(\dot\phi)$, let $B_2:=\{(\dot{\mb v},\dot q):\dot{\mathcal E}(\dot{\mb v},\dot q)\leq\rho_2\}$ and define $\mathcal T_2(\dot{\mb v},\dot q):=(\dot{\mb u},\dot p)$ to be the solution of \eqref{system: stationary linear} at the derivative parameter $l+1$ with data
\begin{align*}
\brac{\dot\phi,\ D\mb N[\phi_0,\mb u,p](\dot\phi,\dot{\mb v},\dot q),\ D\mb N_\partial[\phi_0,\mb u,p](\dot\phi,\dot{\mb v},\dot q)},
\end{align*}
upgraded by Proposition \ref{corner ladder} for $M=1$, the data norms being available as in Step 1 with \eqref{DN estimate} and $\dot{\mathcal D}$ in place of \eqref{N estimate} and $\mathcal D$. By \eqref{DN estimate} with $\mathcal E\leq 2C_{\rm lin}\mathcal D$, for $(\dot{\mb v},\dot q)\in B_2$,
\begin{align*}
\dot{\mathcal E}\brac{\mathcal T_2(\dot{\mb v},\dot q)}
\leq C_{\rm lin}(\dot{\mathcal D}+C_{\rm nl}\mathcal D\rho_2+C_{\rm nl}\mathcal D(1+2C_{\rm lin})\dot{\mathcal D})
\leq{3\over 2}C_{\rm lin}\dot{\mathcal D}+{1\over 2}\rho_2=\rho_2
\end{align*}
by \eqref{delta conditions}, where $\dot{\mathcal D}=\dot{\mathcal D}(\dot\phi)$. Since $D\mb N$ and $D\mb N_\partial$ are affine in $(\dot{\mb v},\dot q)$, the difference of two outputs is the solution with data $(0,D\mb N,D\mb N_\partial)$ evaluated at the difference and at the profile direction $0$, so that \eqref{DN estimate} with $\dot{\mathcal D}=0$ gives the contraction constant $C_{\rm lin}C_{\rm nl}\mathcal D\leq 1/2$. As $\dot{\mathcal E}$ is a norm by Remark \ref{metric remark}, $B_2$ is a complete metric space and we obtain a unique fixed point $(\dot{\mb u},\dot p)\in B_2$ solving \eqref{system: dotted}, with $\dot{\mathcal E}\leq 3C_{\rm lin}\dot{\mathcal D}\lesssim_{\theta_b}\mathcal D$ by \eqref{P bound}. As in Step 1, a solution of \eqref{system: dotted} with $\dot{\mathcal E}<\infty$ lies in the class of \eqref{data estimate} at the derivative parameter $l+1$, by \eqref{unweighted by weighted} with $j=l+1$, and the contraction constant $C_{\rm lin}C_{\rm nl}\mathcal D$ does not depend on $\rho_2$; so the fixed point is the only solution of \eqref{system: dotted} with $\dot{\mathcal E}<\infty$.

\emph{Step 4: assembly.} Define $\ddot\phi:=(\sin\theta_b)^{-2}\partial_r\dot u_\theta|_{\partial_1\Omega}$; its bounds for $M\in\{0,1\}$ are \eqref{ddotphi bound}. Summing the bounds of Steps 1--3 gives \eqref{main estimate}. All the equations of the tuple hold: \eqref{system: stationary nonlinear} and \eqref{system: dotted} by the two fixed points, the profile evolution equation and its time derivative by the definitions of $\dot\phi$ and $\ddot\phi$, the contact angle condition by the hypothesis $\phi_0(0)=0$, and its first time derivative by Step 2.

\emph{Part (iii).} Under \eqref{C2}, Lemma \ref{embeddings}(i) at $a=1-\kappa$ gives $|\ddot\phi(r)|\lesssim_{\theta_b}r^{1-\kappa}\|\ddot\phi\|_{\partial_1\Omega,3/2-\kappa,l+3/2}$, hence the vanishing at the contact point. For $\kappa\in(1/4,1/2)$, the range of Theorem \ref{linear solution in time}, the hypotheses of Theorem \ref{linear solution full regularity} are met: $\dot\phi$ and $\ddot\phi$ are the profiles derived from the two pairs, $(\mb u,p)$ satisfies \eqref{C1 in the linear theory}, which is \eqref{C1}, and $\|\ddot\phi\|_{\partial_1\Omega,3/2-\kappa,l+3/2}$ is finite, as the theorem asks.
\end{proof}

By part (iii) of Theorem \ref{nonlinear initial data}, the twice time-differentiated contact angle condition holds pointwise under \eqref{C2}.

Under a second-order corner expansion of the data and the strengthened angle condition \eqref{third rung angle condition}, the next theorem carries the corner decompositions of the pairs of Theorem \ref{nonlinear initial data} one order further, and it resolves \eqref{C2} into the two explicit equations \eqref{C2 explicit}. It is used for the rung $M=2$ along the flow.

\begin{theorem}[Corner structure of the initial jet]\label{jet corner structure}
Let $\phi_0$ satisfy the hypotheses of part (ii) of Theorem \ref{nonlinear initial data}, and let $(\mb u,p)$, $(\dot{\mb u},\dot p)$, $\dot\phi$ and $\ddot\phi$ be as in that part. Assume in addition the strengthened angle condition \eqref{third rung angle condition} and that the extended data norm \eqref{data norm +} is finite, $\mathcal D_+<\infty$. After shrinking the constant $\delta$ of that theorem by a factor depending only on $\theta_b,\kappa$ and $l$, we have $\mb u^{[2]}=\mb u^{\trimtop{0.4pt}{[3]}}_{\rm c}+\mb u^{[3]}$ and $p^{[2]}=p^{\trimtop{0.4pt}{[3]}}_{\rm c}+p^{[3]}$ as in Proposition \ref{third rung}(ii), and $\dot{\mb u}$ has the decomposition \eqref{dotted log decomposition}. With $\mathfrak X$ at the count $j=l+2$ and $\dot{\mathfrak X}$ at the count $j=l+1$, as in Proposition \ref{nonlinear corner structure}(iii), and with $\|\dot\phi^{[2]}\|_{\partial_1\Omega,5/2-\kappa,l+5/2}$ in place of $\|\dot\phi^{[2]}\|_{\partial_1\Omega,5/2-\kappa,5/2}$ in $\dot{\mathcal D}_+$ of \eqref{derived data norm plus},
\begin{align}\label{full counts estimate}
\mathfrak X+\dot{\mathfrak X}+\dot{\mathcal D}_+\lesssim_{\theta_b}\mathcal D_+(1+\mathcal D_+)^3.
\end{align}
Moreover, the second compatibility condition \eqref{C2} holds if and only if the extended second compatibility condition \eqref{C2 explicit} of Lemma \ref{second-order corner values} holds, that is, if and only if $\smash[b]{c_{\dot\phi,\dr_{-1}^3}}=0$ and $\smash[b]{c_{\dot\phi,\dr_{-1}^2}}=\dot u_1(\mb 0)/(\theta_b(\sin\theta_b)^2)$; in that case $(\sin\theta_b)^2\|\ddot\phi\|_{\partial_1\Omega,3/2-\kappa,l+3/2}\lesssim_{\theta_b}\mathcal D_+(1+\mathcal D_+)^3$.
\end{theorem}
\begin{proof}
We use the notation of the proof of Theorem \ref{nonlinear initial data}, in which $(\mb u,p)$ and $(\dot{\mb u},\dot p)$ are the fixed points of the contractions $\mathcal T_1$ of $B$ and $\mathcal T_2$ of $B_2$. Step 1 carries the third-rung structure along the iteration of $\mathcal T_1$, Step 2 gives $\dot\phi$ its block at $d=1$, Step 3 carries the decomposition \eqref{dotted log decomposition} along the iteration of $\mathcal T_2$, and Step 4 proves the equivalence of \eqref{C2} and \eqref{C2 explicit}.

\emph{Step 1: the third rung of $(\mb u,p)$.} We give $(\mb u,p)$ its third-rung structure at the count $j=l+2$, by propagating the structure along the iteration of $\mathcal T_1$ and passing to the limit with Lemma \ref{corner blocks under limits}(i). Let $(\mb v^{(0)},q^{(0)})$ be the solution of \eqref{system: stationary linear} with data $(\phi_0,\mb 0,\mb 0)$, which lies in $B$ by Step 1 of the proof of Theorem \ref{nonlinear initial data}, and let $(\mb v^{(n+1)},q^{(n+1)}):=\mathcal T_1(\mb v^{(n)},q^{(n)})$, so that $(\mb v^{(n)},q^{(n)})\to(\mb u,p)$ in the metric $d$. Every iterate is a solution of the linear system \eqref{system: stationary linear} with the $M\leq 2$ conclusions of Proposition \ref{corner ladder}. So Proposition \ref{third rung}(ii) at $j=l+2$ applies to it as soon as its data admit the decompositions \eqref{third rung data classes}, the trace hypothesis being available at the count $l+{5\over 2}$ from the counts $(l+2,2)$ of $\mathcal E$ by Lemma \ref{trace}. For $n=0$ the data admit them, since the forcing vanishes and the profile has the blocks by the definition of $\mathcal D_+$, with $\mathcal R_{l+2}\lesssim_{\theta_b}\mathcal D_+$: the norm $\|\dr\phi_0^{\smash[t]{[3]}}\|_{\partial_1\Omega,7/2-\kappa,l+5/2}$ in $\mathcal R_{l+2}$ is bounded by $\|\phi_0^{\smash[t]{[3]}}\|_{\partial_1\Omega,7/2-\kappa,l+7/2}$ in $\mathcal D_+$. Hence $\mathfrak X^{(0)}\lesssim_{\theta_b}\mathcal D_++\mathcal E^{(0)}\lesssim_{\theta_b}\mathcal D_+$. If the $n$-th iterate carries the structure with $\mathfrak X^{(n)}<\infty$, then Proposition \ref{nonlinear corner structure}(i), in the form of part (iii) at the count $j=l+2$, applies to it: its hypotheses ask of the pair only the splitting, $v^{(n)}_\theta|_{\partial_0\Omega}=0$, which is \eqref{stationary no penetration}, and the finiteness of $\mathcal E$ and $\mathfrak X$. So the data of the next iterate admit the decompositions, with $\mathcal R_{l+2}(\mb N,\mb N_\partial)\lesssim_{\theta_b}\mathcal D_+(\mathcal D_++\mathcal E^{(n)})+\mathcal D\mathfrak X^{(n)}$; the remaining terms of $\mathcal R_{l+2}$, those of the profile, are $\lesssim_{\theta_b}\mathcal D_+$ as for $n=0$. Hence \eqref{third rung estimate} gives
\begin{align*}
\mathfrak X^{(n+1)}\leq C_3\brac{\mathcal D_++\mathcal D_+\brac{\mathcal D_++\mathcal E^{(n)}}+\mathcal D\mathfrak X^{(n)}+\mathcal E^{(n+1)}}\leq A+{1\over 2}\mathfrak X^{(n)},
\end{align*}
with $C_3=C_3(\theta_b,\kappa,l)$ and $A\lesssim_{\theta_b}\mathcal D_+(1+\mathcal D_+)$, after shrinking $\delta$ so that $C_3\delta\leq{1\over 2}$, using $\mathcal E^{(n)}\leq 2C_{\rm lin}\mathcal D$ and $\mathcal D\leq\mathcal D_+$; here the trace term and the corner coefficients on the right of \eqref{third rung estimate} are bounded by $\mathcal E^{(n+1)}$ and $\mathcal D$ through Lemmas \ref{trace} and \ref{decomposition estimates}(ii). Hence $\sup_n\mathfrak X^{(n)}\leq 2A+\mathfrak X^{(0)}\lesssim_{\theta_b}\mathcal D_+(1+\mathcal D_+)$.

The convergence in the metric $d$ carries the decomposed pieces in the weighted spaces and the corner coefficients, by Lemma \ref{decomposition estimates}(ii), so the pieces $\mb v^{(n)[2]}$ and $q^{(n)[2]}$ converge to $\mb u^{[2]}$ and $p^{[2]}$ in the sense of distributions. Lemma \ref{corner blocks under limits}(i), applied with the classes $\mathcal Q_2^3$ at the weight $4-\kappa>3$ and $\mathcal Q_1^3$ at $3-\kappa>2$, therefore gives $(\mb u,p)$ the third-rung decompositions, with $\mathfrak X\lesssim_{\theta_b}\mathcal D_+(1+\mathcal D_+)$.

\emph{Step 2: the block of $\dot\phi$.} By Proposition \ref{third rung}(iii) at $j=l+2$, whose additional hypothesis \eqref{C1 in the linear theory} is \eqref{C1}, the derived profile $\dot\phi$ carries the block \eqref{higher corner coefficients} at $(d,m)=(1,3)$, with $\|\dot\phi^{[2]}\|_{\partial_1\Omega,5/2-\kappa,l+5/2}$ bounded by \eqref{derived profile third rung}. The RHS of \eqref{derived profile third rung} is that of \eqref{third rung estimate}, which for $(\mb u,p)$ is $\lesssim_{\theta_b}\mathcal D_+(1+\mathcal D_+)$ as in Step 1, by \eqref{N third rung estimate} at the count $l+2$ and the bound on $\mathfrak X$ of Step 1. Together with \eqref{P bound} this gives $\dot{\mathcal D}_+\lesssim_{\theta_b}\mathcal D_+(1+\mathcal D_+)$.

\emph{Step 3: the decomposition of $(\dot{\mb u},\dot p)$.} We propagate the decomposition of Proposition \ref{third rung}(i) at $j=l+1$ along the iteration of $\mathcal T_2$, started at the solution $(\dot{\mb v}^{(0)},\dot q^{(0)})$ of \eqref{system: stationary linear} with data $(\dot\phi,\mb 0,\mb 0)$. This solution lies in $B_2$, since $\dot{\mathcal E}\leq C_{\rm lin}\dot{\mathcal D}\leq\rho_2$, and it carries the decomposition with $\dot{\mathfrak X}^{(0)}\lesssim_{\theta_b}\dot{\mathcal D}_++\dot{\mathcal E}^{(0)}\lesssim_{\theta_b}\mathcal D_+(1+\mathcal D_+)$, its profile datum having the block of Step 2. If the $n$-th iterate carries the decomposition with $\dot{\mathfrak X}^{(n)}<\infty$, then by Proposition \ref{nonlinear corner structure}(ii) and (iii) at $j=l+1$ the Mellin transforms of the forcing of the next iterate are holomorphic across $\Re\lambda=1$, with $\|D\mb N\|_{\Omega,1-\kappa,l+1}+\|D\mb N_\partial\|_{\partial_1\Omega,3/2-\kappa,l+3/2}$ bounded by \eqref{DN second rung estimate} with $\dot{\mathfrak X}^{(n)}$ in place of $\dot{\mathfrak X}$. The hypotheses of part (iii) of that proposition are the finiteness of $\dot{\mathcal E}$, which $B_2$ supplies, and of $\|\dot\phi^{[2]}\|_{\partial_1\Omega,5/2-\kappa,l+5/2}$, which Step 2 supplies; the splitting of $(\mb u,p)$ that the proposition asks for is that of Step 1. Proposition \ref{third rung}(i) at $j=l+1$ then gives, as in Step 1,
\begin{align*}
\dot{\mathfrak X}^{(n+1)}\lesssim_{\theta_b}\mathcal D\dot{\mathfrak X}^{(n)}+\mathcal D_+\dot{\mathcal E}^{(n)}+\brac{\mathcal D_++\mathcal E+\mathfrak X}\dot{\mathcal D}_++\dot{\mathcal E}^{(n+1)}+\dot{\mathcal D}_+\leq A_2+{1\over 2}\dot{\mathfrak X}^{(n)},
\end{align*}
with $A_2\lesssim_{\theta_b}\mathcal D_+(1+\mathcal D_+)^3$, using $\dot{\mathcal E}^{(n)}\leq\rho_2\lesssim_{\theta_b}\mathcal D$ and the bounds on $\mathfrak X$ and $\dot{\mathcal D}_+$ of Steps 1 and 2, after a further shrinking of $\delta$, so that $\sup_n\dot{\mathfrak X}^{(n)}\lesssim_{\theta_b}\mathcal D_+(1+\mathcal D_+)^3$. Lemma \ref{corner blocks under limits}(i), now with the constant and the three degree one fields of \eqref{dotted log decomposition} as the block class at the weight $3-\kappa>2$, gives the fixed point $(\dot{\mb u},\dot p)$ of $\mathcal T_2$ the decomposition \eqref{dotted log decomposition}, with $\dot{\mathfrak X}\lesssim_{\theta_b}\mathcal D_+(1+\mathcal D_+)^3$. Here the coefficient $\dot u_1^{(n)}(\mb 0)$ converges by Lemma \ref{decomposition estimates}(ii), and the profile coefficients are those of the fixed datum $\dot\phi$. Together with Steps 1 and 2 this proves \eqref{full counts estimate}.

\emph{Step 4: the conditions \eqref{C2} and \eqref{C2 explicit}.} Lemma \ref{second-order corner values} applies to $(\dot{\mb u},\dot p)$ by Step 3, with $\dot{\mathfrak X}$ formed at the counts $(l+1,2)$. If \eqref{C2} holds, then $\ddot\phi$ extends continuously by $0$ to the contact point by Theorem \ref{nonlinear initial data}(iii), so \eqref{C2 explicit} holds by that lemma. Conversely, under \eqref{C2 explicit} the lemma gives $(\sin\theta_b)^2\|\ddot\phi\|_{\partial_1\Omega,3/2-\kappa,l+3/2}\lesssim_{\theta_b}\dot{\mathfrak X}\lesssim_{\theta_b}\mathcal D_+(1+\mathcal D_+)^3$, so \eqref{C2} holds.
\end{proof}

The fixed points of Theorem \ref{nonlinear initial data} depend Lipschitz continuously on the initial profile. To state this, for two profiles $\phi_0,\phi_0'$ and two pairs we write $\mathcal E^{\phi_0-\phi_0'}(\mb u-\mb u',p-p')$ for the norm \eqref{solution norms} of the difference, formed with the decomposition constant $\smash[b]{c_{\phi_0-\phi_0',\dr_{-1}^2}}=\smash[b]{c_{\phi_0,\dr_{-1}^2}}-\smash[b]{c_{\phi_0',\dr_{-1}^2}}$ of the difference profile, the decomposition \eqref{data step relations} being linear in the profile. By Remark \ref{metric remark} this is the distance $d$ of \eqref{contraction metric} when $\phi_0=\phi_0'$, and in general it is the norm attached to the difference of the two affine decompositions: $\mb u^{[2]}-\mb u'^{[2]}$ is exactly the piece $(\mb u-\mb u')^{[2]}$ of the difference formed with the constant $\smash[b]{c_{\phi_0-\phi_0',\dr_{-1}^2}}$.

\begin{proposition}[Lipschitz dependence at the initial time]\label{Lipschitz initial}
Let $\phi_0$ and $\phi_0'$ satisfy the hypotheses of part (i) of Theorem \ref{nonlinear initial data}, and let $(\mb u,p)$ and $(\mb u',p')$ be the pairs of part (i). After shrinking the constant $\delta$ of that theorem by a factor depending only on $\theta_b,\kappa$ and $l$, we have
\begin{align}\label{Lipschitz initial estimate}
\mathcal E^{\phi_0-\phi_0'}(\mb u-\mb u',p-p')+|\mb u(\mb 0)-\mb u'(\mb 0)|+|u_1(\mb 0)-u_1'(\mb 0)|\lesssim_{\theta_b}\mathcal D(\phi_0-\phi_0'),
\end{align}
and the nonlinearities satisfy
\begin{align}\label{Lipschitz forcing estimate}
&\sum_{M=0}^2\brac{\|\mb N[\phi_0,\mb u,p]-\mb N[\phi_0',\mb u',p']\|_{\Omega,M-1-\kappa,l+2}+\|\mb N_\partial[\phi_0,\mb u,p]-\mb N_\partial[\phi_0',\mb u',p']\|_{\partial_1\Omega,M-{1\over 2}-\kappa,l+{5\over 2}}}\nonumber\\
&\qquad\lesssim_{\theta_b}\brac{\mathcal D(\phi_0)+\mathcal D(\phi_0')}\mathcal D(\phi_0-\phi_0').
\end{align}
If both profiles satisfy \eqref{C1}, then for the derived profiles and the dotted pairs of part (ii) also
\begin{align}\label{Lipschitz dotted estimate}
\dot{\mathcal D}(\dot\phi-\dot\phi')+\dot{\mathcal E}(\dot{\mb u}-\dot{\mb u}',\dot p-\dot p')+\sum_{M=0}^1(\sin\theta_b)^2\|\ddot\phi-\ddot\phi'\|_{\partial_1\Omega,M-{1\over 2}-\kappa,l+{3\over 2}}\lesssim_{\theta_b}\mathcal D(\phi_0-\phi_0'),
\end{align}
and, writing $(D\mb N,D\mb N_\partial)$ for the linearisations of the nonlinearities at $(\phi_0,\mb u,p)$ in the direction $(\dot\phi,\dot{\mb u},\dot p)$, and $(D\mb N',D\mb N_\partial')$ for those at $(\phi_0',\mb u',p')$ in the direction $(\dot\phi',\dot{\mb u}',\dot p')$,
\begin{align}\label{Lipschitz dotted forcing estimate}
\sum_{M=0}^1\brac{\|D\mb N-D\mb N'\|_{\Omega,M-1-\kappa,l+1}+\|D\mb N_\partial-D\mb N_\partial'\|_{\partial_1\Omega,M-{1\over 2}-\kappa,l+{3\over 2}}}\lesssim_{\theta_b}\brac{\mathcal D(\phi_0)+\mathcal D(\phi_0')}\mathcal D(\phi_0-\phi_0').
\end{align}
In particular the map $\phi_0\mapsto u_1(\mb 0)$, and the map $\phi_0\mapsto\dot u_1(\mb 0)$ on the set cut out by \eqref{C1}, are Lipschitz on $\{\mathcal D\leq\delta\}$ with constants depending only on $\theta_b,\kappa$ and $l$.
\end{proposition}
\begin{proof}
The difference $(\mb u-\mb u',p-p')$ solves the stationary system \eqref{system: stationary linear} with profile datum $\phi_0-\phi_0'$ and forcing
\begin{align*}
\brac{\mb N[\phi_0,\mb u,p]-\mb N[\phi_0,\mb u',p']}+\brac{\mb N[\phi_0,\mb u',p']-\mb N[\phi_0',\mb u',p']},
\end{align*}
and likewise for $\mb N_\partial$. The linear theory is linear: the solution of Proposition \ref{stationary system} with these data is the difference of the two solutions, and Proposition \ref{corner ladder} gives it the ladder \eqref{um data reg master} with the decomposition \eqref{data step relations} formed with the coefficients of the profile datum $\phi_0-\phi_0'$, which are the differences of the coefficients. Hence, with $C_{\rm lin}$ as in the proof of Theorem \ref{nonlinear initial data},
\begin{align*}
\mathcal E^{\phi_0-\phi_0'}(\mb u-\mb u',p-p')&\leq C_{\rm lin}\mathcal D(\phi_0-\phi_0')\\
&\quad+C_{\rm lin}\sum_{M=0}^2\Big(\|\text{forcing}\|_{\Omega,M-1-\kappa,l+2}+\|\text{boundary forcing}\|_{\partial_1\Omega,M-{1\over 2}-\kappa,l+{5\over 2}}\Big).
\end{align*}
The first difference is linear in the pair difference, and the proof of \eqref{N estimate} applies to it with the pair $(\mb u-\mb u',p-p')$ decomposed with the constant $\smash[b]{c_{\phi_0-\phi_0',\dr_{-1}^2}}$ and with $\mathcal E$ formed accordingly: the decomposition constant enters that proof only through \eqref{corner coefficients 2}, which Lemma \ref{decomposition estimates}(ii) provides for an arbitrary value of the constant, and through \eqref{profile coefficients} in Cases 3 and 4, which for the difference profile reads $(\sin\theta_b)^2|\smash[b]{c_{\phi_0-\phi_0',\dr_{-1}^2}}|\lesssim_{\theta_b}\mathcal D(\phi_0-\phi_0')$. This gives the bound $C_{\rm nl}\mathcal D(\phi_0)(\mathcal E^{\phi_0-\phi_0'}(\mb u-\mb u',p-p')+\mathcal D(\phi_0-\phi_0'))$ for the first difference, the second term collecting the decomposition constants of the difference pair, which enter through \eqref{profile coefficients} in Cases 3 and 4 of that proof.

In the second difference the pair is fixed and the geometry changes: $\bs\nu[\phi_0]-\bs\nu[\phi_0']=\bs\nu[\phi_0-\phi_0']$ by the linearity of the extension. Both geometries lie in the square of Lemma \ref{composition of the geometry} by \eqref{geometry smallness}, so by part (iv) of that lemma the difference of two composition factors is $(\nu_1[\phi_0-\phi_0'])\tilde{\mathfrak G}+(\nu_2[\phi_0-\phi_0'])\tilde{\mathfrak G}$, with bounded factors $\tilde{\mathfrak G}$ whose derivatives produce factors of both geometries. Every summand of the difference is therefore a summand of the schematic form \eqref{schematic bulk}--\eqref{schematic boundary} in which exactly one factor of $\bs\nu$ type is replaced by the difference geometry $\bs\nu[\phi_0-\phi_0']$, including the pure geometry terms of \eqref{schematic boundary}, which are quadratic in the geometry and whose difference splits as $A(B-B')+(A-A')B'$; so the multilinear clause of Proposition \ref{nonlinear estimates} applies. The placements of its proof apply verbatim, the norms of the difference geometry being supplied by Lemma \ref{geometry estimates}(ii) for the profile $\phi_0-\phi_0'$, whose proof is linear in the profile and produces $\mathcal D(\phi_0-\phi_0')$, including on the whole range \eqref{margin ranges}; the difference factor carrying no derivative, which the proof of \eqref{N estimate} suppresses by \eqref{geometry smallness}, is bounded instead by $\|\bs\nu[\phi_0-\phi_0']\|_{L^\infty(\Omega)}\lesssim_{\theta_b}\mathcal D(\phi_0-\phi_0')$ by \eqref{geometry sup} at $a=1$. This gives the bound $C_{\rm nl}\mathcal D(\phi_0-\phi_0')(\mathcal D(\phi_0)+\mathcal D(\phi_0')+\mathcal E(\mb u',p'))$ for the second difference, with $\mathcal E(\mb u',p')\leq 2C_{\rm lin}\mathcal D(\phi_0')\leq 2C_{\rm lin}\delta$.

After shrinking $\delta$ so that $C_{\rm lin}C_{\rm nl}\delta\leq{1\over 2}$, which is \eqref{delta conditions}, the part of the first bound carrying $\mathcal E^{\phi_0-\phi_0'}$ is absorbed into the LHS and \eqref{Lipschitz initial estimate} follows, the corner values being controlled by Lemma \ref{decomposition estimates}(ii) applied to the difference pair. Inserting \eqref{Lipschitz initial estimate} and the bound on $\mathcal E(\mb u',p')$ into the bounds for the two differences gives \eqref{Lipschitz forcing estimate}.

For \eqref{Lipschitz dotted estimate}, the derived profiles satisfy $\dot\phi-\dot\phi'=(\sin\theta_b)^{-2}\partial_r(u_\theta-u_\theta')$, and Lemma \ref{corner values} applied to the difference pair, decomposed with the constant $\smash[b]{c_{\phi_0-\phi_0',\dr_{-1}^2}}$, gives $\dot{\mathcal D}(\dot\phi-\dot\phi')\lesssim_{\theta_b}\mathcal E^{\phi_0-\phi_0'}(\mb u-\mb u',p-p')$: the value $\partial_r(u_\theta-u_\theta')(0,\theta_b)$ in \eqref{corner identity} for the difference vanishes, since $\partial_ru_\theta(0,\theta_b)$ and $\partial_ru_\theta'(0,\theta_b)$ vanish under \eqref{C1}. The difference of the dotted pairs solves the dotted system with profile datum $\dot\phi-\dot\phi'$ and forcing $(D\mb N-D\mb N',D\mb N_\partial-D\mb N_\partial')$. Since the linearisation is linear in the direction,
\begin{align*}
D\mb N-D\mb N'=D\mb N[\phi_0,\mb u,p](\dot\phi-\dot\phi',\dot{\mb u}-\dot{\mb u}',\dot p-\dot p')+\brac{D\mb N[\phi_0,\mb u,p]-D\mb N[\phi_0',\mb u',p']}(\dot\phi',\dot{\mb u}',\dot p'),
\end{align*}
and likewise for $D\mb N_\partial-D\mb N_\partial'$. The linearisation in the difference of the directions is bounded through \eqref{DN estimate} by $C_{\rm nl}(\mathcal D\dot{\mathcal E}(\dot{\mb u}-\dot{\mb u}',\dot p-\dot p')+(\mathcal D+\mathcal E)\dot{\mathcal D}(\dot\phi-\dot\phi'))$. By the multilinear clause, the difference of the linearisations at the two base points is a sum of terms in which one factor of the base point $(\phi_0,\mb u,p)$ is replaced by its difference from $(\phi_0',\mb u',p')$. As in the first part of the proof, this sum is bounded by $C_{\rm nl}(\mathcal D(\phi_0-\phi_0')+\mathcal E^{\phi_0-\phi_0'}(\mb u-\mb u',p-p'))(\dot{\mathcal E}'+\dot{\mathcal D}')$, with $\dot{\mathcal E}'+\dot{\mathcal D}'\lesssim_{\theta_b}\mathcal D(\phi_0')$. Absorbing $C_{\rm lin}C_{\rm nl}\mathcal D\dot{\mathcal E}(\dot{\mb u}-\dot{\mb u}',\dot p-\dot p')$ as in Step 3 of the proof of Theorem \ref{nonlinear initial data} gives the bound on $\dot{\mathcal E}$, and the bounds on $\ddot\phi-\ddot\phi'$ are \eqref{ddotphi bound} for the difference. Inserting \eqref{Lipschitz initial estimate}, \eqref{Lipschitz dotted estimate} and $\mathcal E(\mb u,p)\leq 2C_{\rm lin}\mathcal D(\phi_0)$ into the bounds for the two kinds of terms gives \eqref{Lipschitz dotted forcing estimate}. The last sentence of the statement is \eqref{corner coefficients 2} for the difference pair and \eqref{corner coefficients 1} for the dotted difference pair.
\end{proof}

\begin{remark}[The role of the compatibility conditions]\label{initial data closing remark}
Without \eqref{C1}, part (i) of Theorem \ref{nonlinear initial data} still holds in full, but the construction necessarily stops there. Indeed, by \eqref{corner identity} the derived profile $\dot\phi$ then has a non-zero contact-point value. So the differentiated contact angle condition fails, and $\dot\phi$ is not an admissible datum for Proposition \ref{stationary system}, $\|\dot\phi\|_{\partial_1\Omega,1/2}$ being infinite. Correspondingly, the term $(\sin\theta_b)^2\|\dot\phi\|_{\partial_1\Omega,3/2-\kappa,l+5/2}$ of \eqref{main estimate} is infinite by Lemma \ref{corner values}. Thus, given part (i), the condition \eqref{C1} is necessary as well as sufficient for part (ii), and necessary for part (iii) and for Theorem \ref{jet corner structure}. Since $u_1(\mb 0)$ is determined by $\phi_0$ through the fixed point of Step 1 of the proof of Theorem \ref{nonlinear initial data}, this confirms that the compatibility condition \eqref{C1 in the linear theory} is a constraint on $\phi_0$, and not a local one at the contact point. The second compatibility condition \eqref{C2} plays no role in \eqref{main estimate}: it is the admissibility of $\ddot\phi$ as a profile datum for a further level of the construction. Under the hypotheses of Theorem \ref{jet corner structure}, Lemma \ref{second-order corner values} resolves \eqref{C2} into the two equations \eqref{C2 explicit}. The first states that the coefficient of the second-order logarithm of the derived profile at the contact point vanishes, and the second is the dotted \eqref{C1}.

The set of data cut out by the compatibility conditions is described in Theorem \ref{compatible set} below: near $\phi_0=0$ the data satisfying \eqref{C1} alone form a Lipschitz graph of codimension one, and those satisfying \eqref{C1} and \eqref{C2 explicit} one of codimension three. The mechanism is an implicit function argument whose inputs are Proposition \ref{Lipschitz initial} and Lemmas \ref{corner value response} and \ref{nonlinear corner value response}.
\end{remark}

\subsubsection{The compatible set}

We now make the description of Remark \ref{initial data closing remark} precise. The condition \eqref{C1} prescribes the coefficient $\smash[b]{c_{\phi_0,\dr_{-1}^2}}$ in terms of the corner value $u_1(\mb 0)$, which depends on $\phi_0$ through the fixed point of Theorem \ref{nonlinear initial data}. For the linear system \eqref{system: stationary linear} the corner value is a linear functional of the data, and the next lemma computes its part $\Lambda$ due to the profile datum. The kernel of $\Lambda$ is used below to choose coordinates on the space of data.

\begin{lemma}[Linear response of the corner value]\label{corner value response}
Let $(\mb u,p)$ be the solution of \eqref{system: stationary linear} with data $(\phi_0,\mb F,\mb f)$ given by Propositions \ref{stationary system} and \ref{corner ladder} at the rungs $M\leq 1$, with $\phi_0(0)=0$. Write $\Lambda(\phi_0)$ for the value $u_1(\mb 0)$ of the solution with data $(\phi_0,\mb 0,\mb 0)$, and $\ell_0(\mb F,\mb f)$ for that of the solution with data $(0,\mb F,\mb f)$.
\begin{enumerate}
\item We have $\mb u(\mb 0)=u_1(\mb 0)\mb e_1$ and
\begin{align}\label{corner value response formula}
u_1(\mb 0)=\Lambda(\phi_0)+\ell_0(\mb F,\mb f).
\end{align}
\item The functionals $\Lambda$ and $\ell_0$ are linear, and whenever the RHSs are finite,
\begin{align*}
|\Lambda(\phi_0)|&\lesssim_{\theta_b}(\sin\theta_b)^2\brac{\|\phi_0\|_{\partial_1\Omega,{1\over 2}-\kappa,l+{5\over 2}}+\|\phi_0\|_{\partial_1\Omega,{3\over 2}-\kappa,l+{5\over 2}}},\\
|\ell_0(\mb F,\mb f)|&\lesssim_{\theta_b}\sum_{M=0}^1\brac{\|\mb F\|_{\Omega,M-1-\kappa,l+1}+\|\mb f\|_{\partial_1\Omega,M-{1\over 2}-\kappa,l+{3\over 2}}}.
\end{align*}
In particular $|\Lambda(\phi_0)|\lesssim_{\theta_b}\mathcal D(\phi_0)$, and $|\Lambda(\dot\psi)|\lesssim_{\theta_b}\dot{\mathcal D}(\dot\psi)$ for every profile $\dot\psi$ with $\dot{\mathcal D}(\dot\psi)<\infty$.
\item If $\mb F=\mb 0$ and $\mb f=\mb 0$, then
\begin{align}\label{Lambda formula}
\Lambda(\phi_0)={1\over\sin(2\theta_b)+2\theta_b}\brac{-(\sin\theta_b)^3\int_0^\infty\phi_0(r){\d r\over r}+\int_0^\infty u_r(r,0)\ln r\,\d r},
\end{align}
both integrals converging absolutely, and
\begin{align}\label{slip term bound}
\abs{\int_0^\infty u_r(r,0)\ln r\,\d r}\lesssim_{\theta_b}\sum_{k=0}^1(\sin\theta_b)^2\|\phi_0\|_{\partial_1\Omega,{1\over 2}-\kappa k,l+{5\over 2}}.
\end{align}
\item $\Lambda$ does not vanish identically on $C_c^\infty((0,\infty))$.
\end{enumerate}
\end{lemma}
\begin{proof}
Parts (i) and (ii) follow from the linear theory. For (iii) we bound the integral of $u_r(\ph,0)\ln r$ and then read the corner value off the residue at $\lambda=0$ in the proof of Proposition \ref{corner ladder}. Part (iv) follows from (iii).

\emph{Part (i).} We have $\mb u(\mb 0)=u_1(\mb 0)\mb e_1$ by Proposition \ref{corner ladder}. By the uniqueness in Proposition \ref{stationary system} the solution depends linearly on the data, which gives \eqref{corner value response formula} and the linearity of $\Lambda$ and $\ell_0$.

\emph{Part (ii).} By \eqref{corner coefficients 1}, by \eqref{data estimate} at the derivative parameter $l+1$, and by \eqref{um data reg master} at $M=1$ and $j=1$, whose term $\|u_r\|_{\partial_0\Omega,1/2-\kappa,3/2}$ is bounded by $\|\mb u\|_{\Omega,1-\kappa,l+1,2}$ by Lemma \ref{trace},
\begin{align*}
|\mb u(\mb 0)|&\lesssim_{\theta_b}\|\mb u\|_{\Omega,1-\kappa,l+1,2}+\|\mb u^{[1]}\|_{\Omega,2-\kappa,1,2}\\
&\lesssim_{\theta_b}\sum_{k=0}^1\brac{\|\mb F\|_{\Omega,-1-\kappa k,l+1}+\|\mb f\|_{\partial_1\Omega,-{1\over 2}-\kappa k,l+{3\over 2}}+(\sin\theta_b)^2\|\phi_0\|_{\partial_1\Omega,{1\over 2}-\kappa k,l+{5\over 2}}}\\
&\qquad+\|\mb F\|_{\Omega,-\kappa,1}+\|\mb f\|_{\partial_1\Omega,{1\over 2}-\kappa,{3\over 2}}+(\sin\theta_b)^2\|\phi_0\|_{\partial_1\Omega,{3\over 2}-\kappa,{5\over 2}}.
\end{align*}
The RHS is bounded by the RHSs in (ii). A finite RHS in the bound on $\Lambda$ gives $\phi_0(0)=0$ by Lemma \ref{embeddings}(i) at $a=1-\kappa$, so that $\Lambda(\phi_0)$ is defined. The norms in this bound are dominated by summands of $\mathcal D$, where $\phi_0^{\trimtop{0.5pt}{[0]}}=\phi_0^{\trimtop{0.5pt}{[1]}}=\phi_0$, and they are summands of $\dot{\mathcal D}$; this gives the two particular bounds.

\emph{Part (iii): the two integrals.} Let $\mb F=\mb 0$ and $\mb f=\mb 0$. By the Cauchy--Schwarz inequality on $\{r<1\}$ and on $\{r>1\}$,
\begin{align*}
\int_0^\infty|u_r(r,0)\ln r|\d r\leq\|\ln r\|_{L^2((0,1))}\|u_r\|_{\partial_0\Omega,0}+\|r^{-{1\over 2}-\kappa}\ln r\|_{L^2((1,\infty))}\|u_r\|_{\partial_0\Omega,-{1\over 2}-\kappa}.
\end{align*}
The norm $\|u_r\|_{\partial_0\Omega,0}$ is part of the Navier slip term of Proposition \ref{stationary system}. By \eqref{stationary Navier slip}, $\|u_r\|_{\partial_0\Omega,-1/2-\kappa}=\|\partial_\theta u_r\|_{\partial_0\Omega,1/2-\kappa}$, and this norm is $\lesssim_{\theta_b}\|\mb u\|_{\Omega,1-\kappa,0,2}$ by Lemma \ref{trace}. By \eqref{data estimate} at the derivative parameter $l+1$ both norms are bounded by the RHS of \eqref{slip term bound}, which proves \eqref{slip term bound} and the absolute convergence of the second integral in \eqref{Lambda formula}. The same two norms show that the Mellin transform of the slip channel $F_{0,0}=ru_r$ of \eqref{M01 data def}, which is
\begin{align*}
\hat F_{0,0}(\lambda)={1\over\sqrt{2\pi}}\int_0^\infty r^{-\lambda}u_r(r,0)\d r
\end{align*}
by Lemma \ref{Mellin properties}(i), converges absolutely and is holomorphic on the strip $-\kappa<\Re\lambda<1/2$, and the second integral in \eqref{Lambda formula} is $-\sqrt{2\pi}\hat F_{0,0}'(0)$. Likewise the Mellin transform $\hat\phi_0$ converges absolutely and is holomorphic on the strip $-\kappa<\Re\lambda<1-\kappa$ between the lines \eqref{weight and line} of the weights $1/2-\kappa$ and $3/2-\kappa$ in (ii), and $\sqrt{2\pi}\hat\phi_0(0)$ is the first integral in \eqref{Lambda formula}.

\emph{Part (iii): the corner value.} In the proof of Proposition \ref{corner ladder}, \eqref{residue rule} at $d=0$ gives $\mb u^{[1]}-\mb u^{[0]}=\sqrt{2\pi}\Res_{\lambda=0}(r^\lambda\hat{\mb u})$, where $\hat{\mb u}$ is given by the representation of Proposition \ref{Mellin system initial time} and has at most a double pole at $\lambda=0$. With $m=2$ the RHS is $\sqrt{2\pi}\ln r\lim_{\lambda\to 0}\lambda^2\hat{\mb u}+\sqrt{2\pi}\lim_{\lambda\to 0}\partial_\lambda(\lambda^2\hat{\mb u})$, and $\mb u^{[1]}-\mb u^{[0]}=-u_1(\mb 0)\mb e_1$ by \eqref{data step relations} and (i). Comparing the coefficients of $\ln r$ and of $1$ in the $\mb e_\theta$ component, where $\mb e_1\cdot\mb e_\theta=-\sin\theta$, we get
\begin{align}\label{corner value residue}
\lim_{\lambda\to 0}\lambda^2\hat u_\theta(\lambda,\theta)=0,\qquad\qquad\sqrt{2\pi}\lim_{\lambda\to 0}\partial_\lambda\brac{\lambda^2\hat u_\theta(\lambda,\theta)}=u_1(\mb 0)\sin\theta.
\end{align}
For the data $(\phi_0,\mb 0,\mb 0)$ the only nonzero channels of \eqref{M01 data def} are $F_{0,0}$ and $F_{1,1}=(\sin\theta_b)^2\dr\phi_0$, whose Mellin transform is $(\sin\theta_b)^2\lambda\hat\phi_0(\lambda)$ by Lemma \ref{Mellin properties}(ii). So the representation reduces to
\begin{align*}
\lambda^2\hat u_\theta(\lambda,\theta)=-(\sin\theta_b)^2\lambda^2G_\lambda(\theta,\theta_b)\lambda_1\lambda_{-1}\lambda\hat\phi_0(\lambda)+\partial_\xi(\lambda^2G_\lambda)(\theta,0)\lambda_1\hat F_{0,0}(\lambda).
\end{align*}
We insert $\lim_{\lambda\to 0}\lambda^2G_\lambda(\theta,\xi)=-\sin\theta\sin\xi/(\sin(2\theta_b)+2\theta_b)$ and $\lim_{\lambda\to 0}\partial_\lambda(\lambda^2G_\lambda(\theta,\xi))=0$ from Lemma \ref{Green's function initial time}, both for all $\theta,\xi$ by the symmetry of $G_\lambda$, and $\lambda_1\lambda_{-1}=-1$ at $\lambda=0$, where $\lambda_{\pm 1}=\lambda\pm 1$ as in Definition \ref{scaled radial derivatives}. Since $\lambda\hat\phi_0(\lambda)$ vanishes at $\lambda=0$ with derivative $\hat\phi_0(0)$, we get
\begin{align*}
\lim_{\lambda\to 0}\lambda^2\hat u_\theta(\lambda,\theta)&=-{\sin\theta\over\sin(2\theta_b)+2\theta_b}\hat F_{0,0}(0),\\
\lim_{\lambda\to 0}\partial_\lambda\brac{\lambda^2\hat u_\theta(\lambda,\theta)}&=-{\sin\theta\over\sin(2\theta_b)+2\theta_b}\brac{(\sin\theta_b)^3\hat\phi_0(0)+\hat F_{0,0}(0)+\hat F_{0,0}'(0)}.
\end{align*}
Comparing with \eqref{corner value residue} gives $\hat F_{0,0}(0)=0$, and then $u_1(\mb 0)=-\sqrt{2\pi}((\sin\theta_b)^3\hat\phi_0(0)+\hat F_{0,0}'(0))/(\sin(2\theta_b)+2\theta_b)$, which is \eqref{Lambda formula}.

\emph{Part (iv).} If $\Lambda$ vanished on $C_c^\infty((0,\infty))$, then by \eqref{Lambda formula} and \eqref{slip term bound} we would have
\begin{align*}
(\sin\theta_b)^3\abs{\int_0^\infty\phi{\d r\over r}}\lesssim_{\theta_b}\sum_{k=0}^1(\sin\theta_b)^2\|\phi\|_{\partial_1\Omega,{1\over 2}-\kappa k,l+{5\over 2}}\qquad\text{ for }\qquad\phi\in C_c^\infty((0,\infty)).
\end{align*}
Let $\zeta_n$ be a smooth cut-off with $0\leq\zeta_n\leq 1$, equal to $1$ on $[e^{-n},{1\over 2}]$ and supported in $[e^{-n-1},1]$, with $|\dr^{j'}\zeta_n|\lesssim_{j'}1$ for $j'\in\N_0$, and let $\phi_n:=n^{-1/2}\zeta_n$. Then
\begin{align*}
\int_0^1|\dr^{j'}\phi_n|^2{\d r\over r}\lesssim_{j'}{n+1\over n},\qquad\qquad\int_0^\infty\phi_n{\d r\over r}\geq{n-\ln 2\over\sqrt n}.
\end{align*}
The bound on $\int_0^1|\dr^{j'}\phi_n|^2\d r/r$ bounds the RHS above uniformly in $n$: $\|\ph\|_{\partial_1\Omega,1/2}^2=\int_0^\infty|\ph|^2\d r/r$, and $\|\ph\|_{\partial_1\Omega,1/2-\kappa}$ is dominated by $\|\ph\|_{\partial_1\Omega,1/2}$ on $\{r<1\}$. The lower bound on $\int_0^\infty\phi_n\d r/r$ tends to infinity, a contradiction.
\end{proof}

For the pair $(\mb u,p)$ of Theorem \ref{nonlinear initial data}(i) the forcing $(\mb F,\mb f)$ in \eqref{corner value response formula} is $(\mb N,\mb N_\partial)$, which is at least quadratic. So the corner value is $\Lambda$ of the profile up to a term that is Lipschitz in the profile with a constant of order $\delta$, and a change of the profile in the kernel of $\Lambda$ changes the corner value only at that order. This is the smallness behind the contractions in the proof of Theorem \ref{compatible set} below.

\begin{lemma}[Nonlinear response of the corner value]\label{nonlinear corner value response}
Let $\phi_0$, $\phi_0'$ and the pairs $(\mb u,p)$, $(\mb u',p')$ be as in Proposition \ref{Lipschitz initial}.
\begin{enumerate}
\item We have $|u_1(\mb 0)-u_1'(\mb 0)-\Lambda(\phi_0-\phi_0')|\lesssim_{\theta_b}\delta\mathcal D(\phi_0-\phi_0')$.
\item If both profiles satisfy \eqref{C1}, then for the dotted pairs and the derived profiles of part (ii) of Theorem \ref{nonlinear initial data} also $|\dot u_1(\mb 0)-\dot u_1'(\mb 0)-\Lambda(\dot\phi-\dot\phi')|\lesssim_{\theta_b}\delta\mathcal D(\phi_0-\phi_0')$.
\end{enumerate}
\end{lemma}
\begin{proof}
(i) By \eqref{system: stationary nonlinear} the pair $(\mb u,p)$ solves \eqref{system: stationary linear} with profile datum $\phi_0$ and forcing $(\mb N,\mb N_\partial)[\phi_0,\mb u,p]$, and likewise $(\mb u',p')$. So by \eqref{corner value response formula} the LHS of (i) is $|\ell_0|$ of the difference of the two forcings. The norms in the bound on $\ell_0$ of Lemma \ref{corner value response}(ii) are dominated by the LHS of \eqref{Lipschitz forcing estimate}, and the RHS of \eqref{Lipschitz forcing estimate} is at most $2\delta\mathcal D(\phi_0-\phi_0')$ since $\mathcal D(\phi_0),\mathcal D(\phi_0')\leq\delta$. This proves (i).

(ii) The dotted pair solves \eqref{system: stationary linear} with profile datum $\dot\phi$ and with the forcing $(D\mb N,D\mb N_\partial)[\phi_0,\mb u,p]$ in the direction $(\dot\phi,\dot{\mb u},\dot p)$, by \eqref{system: dotted}, and likewise the primed one. As in (i), the LHS of (ii) is $|\ell_0|$ of the difference of the two forcings. The norms in the bound on $\ell_0$ of Lemma \ref{corner value response}(ii) are those on the left of \eqref{Lipschitz dotted forcing estimate}, and (ii) follows as (i).
\end{proof}

We now fix coordinates on the space of data. The condition \eqref{C1} involves the coefficient of $r\ln r$ of the profile, and by Lemma \ref{profile channel} the condition \eqref{C2 explicit} involves those of $r^2\ln r$ and $r^2\ln^2r$. We take representatives of these three blocks in the kernel of $\Lambda$. Let $\zeta\in C^\infty([0,\infty))$ be a cut-off equal to $1$ on $[0,{1\over 2}]$ and to $0$ on $[1,\infty)$, let $\psi_\Lambda\in C_c^\infty((0,\infty))$ with $\Lambda(\psi_\Lambda)=1$, which exists by Lemma \ref{corner value response}(iv), and for $(d,k)\in\{(1,2),(2,2),(2,3)\}$ and a scale $\omega\in(0,1]$ put
\begin{align}\label{block representatives}
\eta_{d,k}:=\zeta(\omega^{1-d}r){r^d(\ln r)^{k-1}\over(k-1)!}-\Lambda\brac{\zeta(\omega^{1-d}r){r^d(\ln r)^{k-1}\over(k-1)!}}\psi_\Lambda(r),
\end{align}
and $\eta_1:=\eta_{1,2}$, $\eta_2:=\eta_{2,2}$, $\eta_3:=\eta_{2,3}$. The blocks at $d=2$ are thus cut off at the scale $\omega$, which is fixed below so that \eqref{nondegeneracy condition} holds. On the spaces
\begin{align*}
\mathcal X_{\rm dat}:=\set{\phi_0:\mathcal D(\phi_0)<\infty,\ \phi_0(0)=0,\ \lim_{r\to\infty}\phi_0=0},\qquad\qquad\mathcal X_{\rm dat}^+:=\set{\phi_0\in\mathcal X_{\rm dat}:\mathcal D_+(\phi_0)<\infty},
\end{align*}
normed by $\mathcal D$ and by $\mathcal D_+$, we use the block coordinates
\begin{align*}
b_1(\phi_0):=c_{\phi_0,\dr_{-1}^2},\qquad b_2(\phi_0):=c_{\phi_0,\dr_{-2}^2},\qquad b_3(\phi_0):=c_{\phi_0,\dr_{-2}^3},
\end{align*}
$b_1$ on both spaces and $b_2,b_3$ on $\mathcal X_{\rm dat}^+$, and we let $\mathcal Z:=\{b_1=0\}\subset\mathcal X_{\rm dat}$ and $\mathcal Z^+:=\{b_1=b_2=b_3=0\}\subset\mathcal X_{\rm dat}^+$. Finally, for $\eta\in\mathcal X_{\rm dat}$ we write $\mb u_{\rm lin}[\eta]$ for the velocity of the solution of \eqref{system: stationary linear} with data $(\eta,\mb 0,\mb 0)$, and $\dot\phi_{\rm lin}[\eta]:=(\sin\theta_b)^{-2}\partial_r(\mb e_\theta\cdot\mb u_{\rm lin}[\eta])(\ph,\theta_b)$ for its derived profile.

\begin{lemma}[Block coordinates on the data spaces]\label{block coordinates}
With the representatives \eqref{block representatives} and the block coordinates above, the following hold.
\begin{enumerate}
\item $\mathcal X_{\rm dat}$ and $\mathcal X_{\rm dat}^+$ are Banach spaces, $b_1$ is a bounded linear functional on $\mathcal X_{\rm dat}$, and $b_1,b_2,b_3$ are bounded linear functionals on $\mathcal X_{\rm dat}^+$.
\item For $m,m'\in\{1,2,3\}$ with $m\neq m'$ we have $\Lambda(\eta_m)=0$, $\mathcal D_+(\eta_m)<\infty$, $b_m(\eta_m)=1$, $b_m(\eta_{m'})=0$ and $\smash[b]{c_{\eta_m,\dr_{-1}}}=0$.
\item $\mathcal Z$ and $\mathcal Z^+$ are closed subspaces of codimension one and three, with $\mathcal X_{\rm dat}=\mathcal Z\oplus\R\eta_1$ and $\mathcal X_{\rm dat}^+=\mathcal Z^+\oplus\R\eta_1\oplus\R\eta_2\oplus\R\eta_3$: every $\phi_0\in\mathcal X_{\rm dat}$ is $\psi+b_1(\phi_0)\eta_1$ with $\psi\in\mathcal Z$, and every $\phi_0\in\mathcal X_{\rm dat}^+$ is $\psi+\sum_{m=1}^3b_m(\phi_0)\eta_m$ with $\psi\in\mathcal Z^+$.
\item For $m\in\{2,3\}$ the derived profile $\dot\phi_{\rm lin}[\eta_m]$ vanishes at the contact point, and $\dot{\mathcal D}(\dot\phi_{\rm lin}[\eta_m])\lesssim_{\theta_b}\mathcal D(\eta_m)$.
\item $\mathcal D(\eta_2)\to 0$ as $\omega\to 0$.
\end{enumerate}
\end{lemma}
\begin{proof}
(i) The functionals are linear, and they are bounded by \eqref{profile coefficients} for $b_1$ and by $\mathcal D_+$ for $b_2,b_3$. Let $(\phi_0^{(n)})$ be a Cauchy sequence for $\mathcal D$. Its pieces $\phi_0^{(n)[M]}$ converge in the weighted spaces of $\mathcal D$, and its corner coefficients converge by \eqref{profile coefficients} applied to the differences. The limits are the pieces of the limit profile by Lemma \ref{block matching}, and $\phi_0(0)=0$ and $\lim_{r\to\infty}\phi_0=0$ pass to the limit by Lemma \ref{embeddings}(i). The same argument applies to $\mathcal D_+$, with \eqref{data norm +} and the pieces $\phi_0^{(n)[3]}$.

(ii) Let $(d,k)\in\{(1,2),(2,2),(2,3)\}$. We have $\Lambda(\eta_{d,k})=0$ by \eqref{block representatives} and the linearity of $\Lambda$. By Lemma \ref{truncated powers} with $n=k-1$, the Mellin transform of $\mb 1_{(0,1)}r^d(\ln r)^{k-1}/(k-1)!$ is $-(2\pi)^{-1/2}(\lambda-d)^{-k}$ for $\Re\lambda<d$, while $(\zeta(\omega^{1-d}r)-\mb 1_{(0,1)})r^d(\ln r)^{k-1}$ and $\psi_\Lambda$ have compact support in $(0,\infty)$ and therefore entire Mellin transforms. So the Mellin transform of $\eta_{d,k}$ continues from $\Re\lambda<d$ to a function meromorphic on $\C$ whose only pole is at $\lambda=d$, with principal part $-(2\pi)^{-1/2}(\lambda-d)^{-k}$, and \eqref{principal part} gives the coefficients. Finally, $\eta_{d,k}$ is smooth on $(0,\infty)$ with compact support in $[0,\infty)$, and it coincides with its block near $r=0$. Each subtracted block has a finite tail at infinity at the weight at which it is subtracted, ${5\over 2}-\kappa$ in $\mathcal D$ and ${7\over 2}-\kappa$ in $\mathcal D_+$; so $\mathcal D_+(\eta_{d,k})<\infty$.

(iii) Let $\phi_0\in\mathcal X_{\rm dat}^+$. By (ii), the profile $\psi:=\phi_0-\sum_{m=1}^3b_m(\phi_0)\eta_m$ has $b_{m'}(\psi)=0$ for every $m'$, so that $\psi\in\mathcal Z^+$; and applying $b_{m'}$ to $\phi_0=\psi+\sum_{m=1}^3z_m\eta_m$ with $\psi\in\mathcal Z^+$ gives $z_{m'}=b_{m'}(\phi_0)$. The same holds in $\mathcal X_{\rm dat}$ with $b_1$ and $\eta_1$ alone. The subspaces are closed by (i).

(iv) Let $m\in\{2,3\}$. The size $\mathcal E$ of the solution with data $(\eta_m,\mb 0,\mb 0)$, formed with $\eta_m$ in place of $\phi_0$, is $\lesssim_{\theta_b}\mathcal D(\eta_m)$ by Propositions \ref{stationary system} and \ref{corner ladder} and Lemma \ref{trace}, as in Step 1 of the proof of Theorem \ref{nonlinear initial data}, and $\mathcal D(\eta_m)<\infty$ by (ii). By Lemma \ref{corner value response}(i) and part (ii), $\mb u_{\rm lin}[\eta_m](\mb 0)=\Lambda(\eta_m)\mb e_1=\mb 0$, and $\smash[b]{c_{\eta_m,\dr_{-1}^2}}=b_1(\eta_m)=0$ by part (ii). So $\partial_ru_\theta(0,\theta_b)$ in \eqref{corner identity}, formed for this solution with $\eta_m$ in place of $\phi_0$, vanishes, and Lemma \ref{corner values} gives the claim.

(v) By (ii) and \eqref{data step relations}, $\eta_2^{\smash[t]{[M]}}=\eta_2$ for $M\leq 2$. Hence, by \eqref{block representatives} and the bound on $\Lambda$ in Lemma \ref{corner value response}(ii),
\begin{align*}
\mathcal D(\eta_2)\lesssim_{\theta_b}\brac{1+\mathcal D(\psi_\Lambda)}\sum_{M=0}^2\|\zeta(\omega^{-1}r)r^2\ln r\|_{\partial_1\Omega,M+{1\over 2}-\kappa,l+{7\over 2}}.
\end{align*}
By the change of variables $r\mapsto\omega r$ in Definition \ref{Mellin transform}, the Mellin transform of $\zeta(\omega^{-1}r)r^2\ln r$ is $\omega^{2-\lambda}$ times that of $\zeta(r)r^2\ln(\omega r)$. On the line $\Re\lambda=a-1/2$ of \eqref{weight and line} the factor $\omega^{2-\lambda}$ has the modulus $\omega^{5/2-a}$, so that, for $a\leq 5/2-\kappa$,
\begin{align*}
\|\zeta(\omega^{-1}r)r^2\ln r\|_{\partial_1\Omega,a,l+{7\over 2}}&=\omega^{{5\over 2}-a}\|\zeta(r)r^2\ln(\omega r)\|_{\partial_1\Omega,a,l+{7\over 2}}\\
&\leq\omega^\kappa(1+|\ln\omega|)\brac{\|\zeta(r)r^2\ln r\|_{\partial_1\Omega,a,l+{7\over 2}}+\|\zeta(r)r^2\|_{\partial_1\Omega,a,l+{7\over 2}}}.
\end{align*}
The norms on the right are finite since $a<5/2$, and $\omega^\kappa(1+|\ln\omega|)\to 0$ as $\omega\to 0$.
\end{proof}

By Lemma \ref{corner value response}(ii) and Lemma \ref{block coordinates}(iv) we may define, for $m\in\{2,3\}$,
\begin{align}\label{a_i definition}
a_m:={\Lambda(\dot\phi_{\rm lin}[\eta_m])\over\theta_b(\sin\theta_b)^2}.
\end{align}
These two numbers depend on $\theta_b$ and on the representatives \eqref{block representatives}. In Step 3 of the proof of Theorem \ref{compatible set} below, $\theta_b(\sin\theta_b)^2a_m$ is the coefficient of the coordinate along $\eta_m$ in the linear part of $\dot u_1(\mb 0)$. In Step 2 there, the number $\gamma(\theta_b)-a_2$ must not vanish, with $\gamma$ as in \eqref{gamma formula}. By Lemma \ref{corner value response}(ii) and Lemma \ref{block coordinates}(iv), $|a_2|\lesssim_{\theta_b}\mathcal D(\eta_2)$, which tends to $0$ as $\omega\to 0$ by Lemma \ref{block coordinates}(v). Since $\theta_b<\pi/(4-\kappa)<\pi/3$ by \eqref{standing angle condition}, we have $\gamma(\theta_b)\neq 0$ by Lemma \ref{profile channel}. We therefore fix $\omega\in(0,1]$ such that
\begin{align}\label{nondegeneracy condition}
a_2\neq\gamma(\theta_b).
\end{align}

In the coordinates of Lemma \ref{block coordinates}(iii), the next theorem solves \eqref{C1} for $b_1$, and \eqref{C1} and \eqref{C2 explicit} for $(b_1,b_2,b_3)$, as Lipschitz functions of the component of the datum in $\mathcal Z$, respectively $\mathcal Z^+$.

\begin{theorem}[The compatible set]\label{compatible set}
Let $\kappa$, $l$, $\theta_b$ and $\delta$ be as in Theorem \ref{nonlinear initial data}, with $\delta$ shrunk as in Proposition \ref{Lipschitz initial} and further by a factor depending only on $\theta_b$, $\kappa$, $l$ and the representatives \eqref{block representatives}.
\begin{enumerate}
\item There are $\delta'\in(0,\delta]$, $C\geq 1$ and a map $g_1:\{\psi\in\mathcal Z:\mathcal D(\psi)\leq\delta'\}\to\R$, Lipschitz with respect to $\mathcal D$, with $g_1(0)=0$ and $|g_1(\psi)|\leq C\mathcal D(\psi)$, such that for $\psi$ in that ball and $|z|\leq C\delta'$ the profile $\psi+z\eta_1$ satisfies \eqref{C1} if and only if $z=g_1(\psi)$.
\item Assume in addition \eqref{third rung angle condition}. There are $\delta''\in(0,\delta']$, $C'\geq 1$ and a map $(g_2,g_3):B'':=\{\psi\in\mathcal Z^+:\mathcal D(\psi)\leq\delta''\}\to\R^2$, Lipschitz with respect to $\mathcal D$, with $(g_2,g_3)(0)=0$, such that for $\psi\in B''$ and $|z_1|,|z_2|,|z_3|\leq C'\delta''$ the profile $\psi+\sum_{m=1}^3z_m\eta_m$ satisfies \eqref{C1} and \eqref{C2 explicit} if and only if $z_2=g_2(\psi)$, $z_3=g_3(\psi)$ and $z_1=g_1(\psi+z_2\eta_2+z_3\eta_3)$.
\end{enumerate}
\end{theorem}
\begin{proof}
Both parts are proved by the contraction mapping principle in the block coordinates, with contraction constants of order $\delta$ by Lemma \ref{nonlinear corner value response}. Part (ii) solves \eqref{C1} for $z_1$ by part (i) (Step 1), writes \eqref{C2 explicit} as a fixed point equation for $(z_2,z_3)$ whose linear part is triangular by Lemma \ref{profile channel} (Step 2), bounds the remainder (Step 3) and concludes (Step 4).

Throughout, $\mb u[\phi_0]$ denotes the velocity of the pair of Theorem \ref{nonlinear initial data}(i) for the datum $\phi_0$, and all constants depend only on $\theta_b$, $\kappa$, $l$ and the representatives \eqref{block representatives}.

\emph{Part (i).} By \eqref{corner coefficients 2} and Theorem \ref{nonlinear initial data}(i) there is a constant $C_{\rm cv}$ with $|u_1[\phi_0](\mb 0)|\leq C_{\rm cv}\theta_b(\sin\theta_b)^2\mathcal D(\phi_0)$ for every $\phi_0\in\mathcal X_{\rm dat}$ with $\mathcal D(\phi_0)\leq\delta$. Put $C:=2C_{\rm cv}+1$, and let $\delta'$ be so small that $\delta'(1+C\mathcal D(\eta_1))\leq\delta$. For $\psi\in\mathcal Z$ with $\mathcal D(\psi)\leq\delta'$ and $|z|\leq C\delta'$ we then have $\mathcal D(\psi+z\eta_1)\leq\delta$, and $b_1(\psi+z\eta_1)=z$ by Lemma \ref{block coordinates}(ii); so \eqref{C1} for $\psi+z\eta_1$ reads $z=\Phi(\psi,z)$, with
\begin{align*}
\Phi(\psi,z):={u_1[\psi+z\eta_1](\mb 0)\over\theta_b(\sin\theta_b)^2}.
\end{align*}
For two such $\psi,\psi'$ and $z,z'$, Lemma \ref{nonlinear corner value response}(i), with $\Lambda(\eta_1)=0$, and Proposition \ref{Lipschitz initial} give
\begin{align*}
|\Phi(\psi,z)-\Phi(\psi,z')|\leq L|z-z'|\quad\text{ with }\quad L\lesssim\delta\mathcal D(\eta_1),\qquad\qquad|\Phi(\psi,z)-\Phi(\psi',z)|\lesssim\mathcal D(\psi-\psi');
\end{align*}
if $\Lambda(\psi-\psi')=0$, the bound on $\Phi(\psi,z)-\Phi(\psi',z)$ holds with the factor $\delta$ on the right, by Lemma \ref{nonlinear corner value response}(i). The implicit constant in $L\lesssim\delta\mathcal D(\eta_1)$ does not depend on $\delta'$, so the shrinking of $\delta$ in the statement gives $L\leq{1\over 2}$.

Since $|\Phi(\psi,0)|\leq C_{\rm cv}\mathcal D(\psi)$, for $|z|\leq C\delta'$ we get $|\Phi(\psi,z)|\leq C_{\rm cv}\delta'+{1\over 2}C\delta'\leq C\delta'$. Hence $\Phi(\psi,\ph)$ is a contraction of $[-C\delta',C\delta']$, and its unique fixed point $g_1(\psi)$ gives the claim. We have $g_1(0)=0$ since $u_1[0](\mb 0)=0$, and $|g_1(\psi)|\leq(1-L)^{-1}|\Phi(\psi,0)|\leq C\mathcal D(\psi)$. Finally, by the two bounds on $\Phi$ above,
\begin{align}\label{g1 Lipschitz}
|g_1(\psi)-g_1(\psi')|\leq{|\Phi(\psi,g_1(\psi))-\Phi(\psi',g_1(\psi))|\over 1-L}\lesssim\begin{cases}\mathcal D(\psi-\psi'),&\\\delta\mathcal D(\psi-\psi')&\text{if }\Lambda(\psi-\psi')=0.\end{cases}
\end{align}

\emph{Part (ii), Step 1: the condition \eqref{C1}.} Let $\psi\in B''$ and $\mb z:=(z_2,z_3)$ with $|\mb z|:=\max\{|z_2|,|z_3|\}\leq C'\delta''$. Here $C'\geq 2C$ is fixed in Step 4 in terms of the constants of Steps 1--3, which do not depend on $\delta$, $\delta'$, $\delta''$ or $C'$, and $\delta''$ is small in terms of $C'$. Then $\psi+z_2\eta_2+z_3\eta_3$ lies in $\mathcal Z$, by Lemma \ref{block coordinates}(ii), with $\mathcal D\leq\delta'$, and $C'\delta''\leq C\delta'$. By part (i), for $|z_1|\leq C'\delta''$ the profile $\psi+\sum_{m=1}^3z_m\eta_m$ therefore satisfies \eqref{C1} if and only if $z_1=h(\mb z)$, where
\begin{align*}
h(\mb z):=g_1(\psi+z_2\eta_2+z_3\eta_3).
\end{align*}
Since the arguments of $g_1$ at $\mb z$ and $\mb z'$ differ by $(z_2-z_2')\eta_2+(z_3-z_3')\eta_3\in\ker\Lambda$, by Lemma \ref{block coordinates}(ii), \eqref{g1 Lipschitz} gives
\begin{align}\label{h bounds}
|h(\mb z)-h(\mb z')|\lesssim\delta|\mb z-\mb z'|,\qquad\qquad|h(0)|=|g_1(\psi)|\leq C\mathcal D(\psi).
\end{align}
It remains to solve \eqref{C2 explicit} for $\mb z$ along the profiles $\phi_0:=\psi+h(\mb z)\eta_1+z_2\eta_2+z_3\eta_3$. They satisfy \eqref{C1} by part (i). Since $|h(\mb z)|\leq C\mathcal D(\psi+z_2\eta_2+z_3\eta_3)$, they have $\mathcal D(\phi_0)\leq\delta$ by the choice of $\delta'$, and $\mathcal D_+(\phi_0)<\infty$ by $\psi\in\mathcal Z^+$ and Lemma \ref{block coordinates}(ii). So they satisfy the hypotheses of Theorem \ref{jet corner structure}, with $\delta$ shrunk as there, and we write $(\mb u,p)$, $(\dot{\mb u},\dot p)$ and $\dot\phi$ for their pairs and derived profile; the dotted pair carries the decomposition \eqref{dotted log decomposition}. By Lemma \ref{block coordinates}(ii), $b_2(\phi_0)=z_2$ and $b_3(\phi_0)=z_3$.

\emph{Step 2: \eqref{C2 explicit} as a fixed point equation.} By Lemma \ref{profile channel}, with $\gamma=\gamma(\theta_b)\neq 0$ and with the constant $t_{23}$ and the functionals $\ell_2,\ell_3$ of \eqref{triangular profile channel},
\begin{align*}
c_{\dot\phi,\dr_{-1}^3}=\gamma z_3+\ell_3\brac{\mb F^{[3]}_{\rm c},\mb f^{[3]}_{\rm c}},\qquad
c_{\dot\phi,\dr_{-1}^2}=\gamma z_2+t_{23}z_3+\ell_2\brac{\mb F^{[3]}_{\rm c},\mb f^{[3]}_{\rm c}},
\end{align*}
where $(\mb F^{\smash[t]{[3]}}_{\rm c},\mb f^{\smash[t]{[3]}}_{\rm c})$ is the corner part of $(\mb N,\mb N_\partial)[\phi_0,\mb u,p]$ of Proposition \ref{nonlinear corner structure}(i). By that proposition, its coefficients are quadratic forms in the first-order corner data $(\smash[b]{c_{\phi_0,\dr_{-1}}},\smash[b]{c_{\phi_0,\dr_{-1}^2}},u_1(\mb 0))$. So \eqref{C2 explicit} reads $A\mb z+R(\mb z)=0$, with $\mb z$ as a column vector, the equation for $\smash[b]{c_{\dot\phi,\dr_{-1}^2}}$ giving the first row and that for $\smash[b]{c_{\dot\phi,\dr_{-1}^3}}$ the second, where
\begin{align*}
A:=\mat{\gamma-a_2&t_{23}-a_3\\0&\gamma},\qquad\qquad
\begin{aligned}
R_1(\mb z)&:=\ell_2\brac{\mb F^{[3]}_{\rm c},\mb f^{[3]}_{\rm c}}-{\dot u_1(\mb 0)\over\theta_b(\sin\theta_b)^2}+a_2z_2+a_3z_3,\\
R_2(\mb z)&:=\ell_3\brac{\mb F^{[3]}_{\rm c},\mb f^{[3]}_{\rm c}}.
\end{aligned}
\end{align*}
The determinant of $A$ is $\gamma(\gamma-a_2)$, nonzero by $\gamma\neq 0$ and \eqref{nondegeneracy condition}. Hence \eqref{C2 explicit} is the fixed point equation $\mb z=-A^{-1}R(\mb z)$.

\emph{Step 3: the remainder.} We show that there are constants $C_R$ and $C_L$ with, for $|\mb z|,|\mb z'|\leq C'\delta''$,
\begin{align}\label{R bounds}
|R(0)|\leq C_R\mathcal D(\psi),\qquad\qquad|R(\mb z)-R(\mb z')|\leq C_L\delta|\mb z-\mb z'|.
\end{align}
At $\mb z=0$ the profile is $\psi+h(0)\eta_1$, whose norm $\mathcal D$ is $\lesssim\mathcal D(\psi)$ by \eqref{h bounds}. Its first-order corner data are $\lesssim\mathcal D(\psi)$ by \eqref{profile coefficients}, \eqref{corner coefficients 2} and Theorem \ref{nonlinear initial data}(i), and $|\dot u_1(\mb 0)|\lesssim\mathcal D(\psi)$ by \eqref{corner coefficients 1} and \eqref{main estimate}. This gives the bound on $R(0)$.

For the Lipschitz bound let $\phi_0$ and $\phi_0'$ be the profiles at $\mb z$ and at $\mb z'$, with their pairs and derived profiles, so that
\begin{align*}
\phi_0-\phi_0'=(h(\mb z)-h(\mb z'))\eta_1+(z_2-z_2')\eta_2+(z_3-z_3')\eta_3\in\ker\Lambda,\qquad\mathcal D(\phi_0-\phi_0')\lesssim|\mb z-\mb z'|.
\end{align*}
The first-order corner data change by $\lesssim\delta|\mb z-\mb z'|$. Indeed $\smash[b]{c_{\phi_0,\dr_{-1}}}$ does not change, by Lemma \ref{block coordinates}(ii); $\smash[b]{c_{\phi_0,\dr_{-1}^2}}=b_1(\phi_0)=h(\mb z)$ by Lemma \ref{block coordinates}(ii), which changes by $\lesssim\delta|\mb z-\mb z'|$ by \eqref{h bounds}; and $u_1(\mb 0)$ changes by $\lesssim\delta|\mb z-\mb z'|$ by Lemma \ref{nonlinear corner value response}(i). Since the quadratic forms of Step 2 are evaluated at arguments $\lesssim\delta$, we get, for $m\in\{2,3\}$ and with the primed corner part formed at $\phi_0'$,
\begin{align*}
\abs{\ell_m\brac{\mb F^{[3]}_{\rm c},\mb f^{[3]}_{\rm c}}-\ell_m\brac{\mb F'^{[3]}_{\rm c},\mb f'^{[3]}_{\rm c}}}\lesssim|h(\mb z)-h(\mb z')|+|u_1(\mb 0)-u_1'(\mb 0)|\lesssim\delta|\mb z-\mb z'|.
\end{align*}
This bounds the change of $R_2$ and of the term with $\ell_2$ in $R_1$.

For $\dot u_1(\mb 0)$, Lemma \ref{nonlinear corner value response}(ii) gives
\begin{align*}
\dot u_1(\mb 0)-\dot u_1'(\mb 0)=\Lambda(\dot\phi-\dot\phi')+O(\delta|\mb z-\mb z'|).
\end{align*}
To compute $\Lambda(\dot\phi-\dot\phi')$ we split the difference of the velocities by the linearity of the solution operator of Proposition \ref{stationary system},
\begin{align*}
\mb u-\mb u'=(h(\mb z)-h(\mb z'))\mb u_{\rm lin}[\eta_1]+\sum_{m=2}^3(z_m-z_m')\mb u_{\rm lin}[\eta_m]+\mb w,
\end{align*}
where $\mb w$ is the velocity of the solution with profile datum $0$ and with the difference of the forcings $(\mb N,\mb N_\partial)[\phi_0,\mb u,p]$ and $(\mb N,\mb N_\partial)[\phi_0',\mb u',p']$. The size $\mathcal E$ of this solution, formed with $0$ in place of $\phi_0$, is $\lesssim\delta|\mb z-\mb z'|$ by Propositions \ref{stationary system} and \ref{corner ladder} and by \eqref{Lipschitz forcing estimate}. Correspondingly
\begin{align*}
\dot\phi-\dot\phi'=\sum_{m=2}^3(z_m-z_m')\dot\phi_{\rm lin}[\eta_m]+\varrho,
\end{align*}
where $\varrho$ is the derived profile of $(h(\mb z)-h(\mb z'))\mb u_{\rm lin}[\eta_1]+\mb w$.

Both $\dot\phi$ and $\dot\phi'$ vanish at the contact point, by Theorem \ref{nonlinear initial data}(ii), and so do $\dot\phi_{\rm lin}[\eta_2]$ and $\dot\phi_{\rm lin}[\eta_3]$, by Lemma \ref{block coordinates}(iv); hence so does $\varrho$. Consider the solution with velocity $(h(\mb z)-h(\mb z'))\mb u_{\rm lin}[\eta_1]+\mb w$ and profile datum $(h(\mb z)-h(\mb z'))\eta_1$. Its size $\mathcal E$, formed with this profile datum in place of $\phi_0$, is $\lesssim|h(\mb z)-h(\mb z')|+\delta|\mb z-\mb z'|\lesssim\delta|\mb z-\mb z'|$, by \eqref{h bounds}. So Lemma \ref{corner values} gives $\dot{\mathcal D}(\varrho)\lesssim\delta|\mb z-\mb z'|$, and Lemma \ref{corner value response}(ii) gives $|\Lambda(\varrho)|\lesssim\dot{\mathcal D}(\varrho)$. With \eqref{a_i definition} we obtain
\begin{align*}
\dot u_1(\mb 0)-\dot u_1'(\mb 0)=\theta_b(\sin\theta_b)^2\brac{a_2(z_2-z_2')+a_3(z_3-z_3')}+O(\delta|\mb z-\mb z'|),
\end{align*}
whose linear part is cancelled by the term $a_2z_2+a_3z_3$ of $R_1$. So $R_1$ changes by $\lesssim\delta|\mb z-\mb z'|$, which proves \eqref{R bounds}.

\emph{Step 4: conclusion.} Let $C_A$ be the operator norm of $A^{-1}$ for $|\ph|$, and fix $C':=\max\{2C,2C_AC_R\}$. The constants $C_A$, $C_L$ and that of \eqref{h bounds} do not depend on $\delta$, $\delta'$, $\delta''$ or $C'$, so the shrinking of $\delta$ in the statement gives $C_AC_L\delta\leq{1\over 2}$ and $|h(\mb z)-h(\mb z')|\leq{1\over 2}|\mb z-\mb z'|$. Then the map $\mb z\mapsto-A^{-1}R(\mb z)$ is a contraction of $\{|\mb z|\leq C'\delta''\}$ with constant ${1\over 2}$. It maps this set into itself, since
\begin{align*}
|A^{-1}R(\mb z)|\leq C_A\brac{C_R\mathcal D(\psi)+C_L\delta|\mb z|}\leq{1\over 2}C'\delta''+{1\over 2}C'\delta''.
\end{align*}
Its fixed point is $(g_2(\psi),g_3(\psi))$, and $|h(g_2(\psi),g_3(\psi))|\leq C\mathcal D(\psi)+{1\over 2}C'\delta''\leq C'\delta''$ by \eqref{h bounds}, so that $z_1=h(\mb z)$ lies in the range of the statement. With Step 1 this proves the equivalence. For $\psi=0$ all corner data vanish at $\mb z=0$, so that $R(0)=0$ and $(g_2,g_3)(0)=0$.

For the Lipschitz bound let $\psi,\psi'\in B''$, and let $R'$ be formed with $\psi'$ in place of $\psi$. As in part (i),
\begin{align*}
|(g_2,g_3)(\psi)-(g_2,g_3)(\psi')|\leq 2C_A|R(g_2(\psi),g_3(\psi))-R'(g_2(\psi),g_3(\psi))|.
\end{align*}
At fixed $\mb z$ the two profiles differ by $\psi-\psi'+(g_1(\psi+z_2\eta_2+z_3\eta_3)-g_1(\psi'+z_2\eta_2+z_3\eta_3))\eta_1$, whose norm $\mathcal D$ is $\lesssim\mathcal D(\psi-\psi')$ by \eqref{g1 Lipschitz}. So $R$ and $R'$ differ there by $\lesssim\mathcal D(\psi-\psi')$, by \eqref{profile coefficients}, Proposition \ref{Lipschitz initial}, and \eqref{corner coefficients 1} with \eqref{Lipschitz dotted estimate}, the terms $a_2z_2+a_3z_3$ being equal.
\end{proof}

By Theorem \ref{compatible set}, in the splittings of Lemma \ref{block coordinates}(iii) the graphs described in Remark \ref{initial data closing remark} are the sets
\begin{align*}
&\set{\psi+g_1(\psi)\eta_1:\psi\in\mathcal Z},\\
&\set{\psi+g_1(\psi+g_2(\psi)\eta_2+g_3(\psi)\eta_3)\eta_1+g_2(\psi)\eta_2+g_3(\psi)\eta_3:\psi\in\mathcal Z^+},
\end{align*}
with $\psi$ in the domains of $g_1$ and of $(g_2,g_3)$; both contain $\phi_0=0$.

\begin{remark}[Examples of compatible data]\label{compatible data examples}
Besides $\phi_0=0$, the graphs contain a point for every small smooth profile with compact support away from the contact point. Let $\psi\in C_c^\infty((0,\infty))$. Its Mellin transform is entire, so that $\psi$ has no corner blocks. Its norms $\mathcal D(\psi)$ and $\mathcal D_+(\psi)$ are finite, so that $\psi\in\mathcal Z^+$. If $\mathcal D(\psi)\leq\delta'$, which holds for a small multiple of any such function, then the point
\begin{align*}
\phi_0=\psi+g_1(\psi)\eta_1
\end{align*}
of the first graph satisfies \eqref{C1}. It is nonzero if $\psi\neq 0$, and distinct functions $\psi$ give distinct points, since $\psi$ is the component of $\phi_0$ in $\mathcal Z$ by Lemma \ref{block coordinates}(iii). Near the contact point the correction $g_1(\psi)\eta_1$ is $g_1(\psi)r\ln r$, by \eqref{block representatives}. Since $b_1(\psi)=0$ and $b_1(\eta_1)=1$, the condition \eqref{C1} for $\phi_0$ reads $g_1(\psi)=u_1(\mb 0)/(\theta_b(\sin\theta_b)^2)$, where $u_1(\mb 0)$ is the corner value of the pair of Theorem \ref{nonlinear initial data}(i) for $\phi_0$.

To first order the coefficient $g_1(\psi)$ is explicit. By \eqref{corner value response formula} and $\Lambda(\eta_1)=0$ we have $u_1(\mb 0)=\Lambda(\psi)+\ell_0(\mb N,\mb N_\partial)$, with the nonlinearities taken at $(\phi_0,\mb u,p)$, and $|\ell_0(\mb N,\mb N_\partial)|\lesssim_{\theta_b}\mathcal D(\phi_0)^2$ by Lemma \ref{corner value response}(ii), \eqref{N estimate} and Theorem \ref{nonlinear initial data}(i). With $C$ as in Theorem \ref{compatible set}(i), this gives
\begin{align}\label{first-order coefficient}
\abs{g_1(\psi)-{\Lambda(\psi)\over\theta_b(\sin\theta_b)^2}}\lesssim_{\theta_b}\mathcal D(\phi_0)^2\leq\brac{1+C\mathcal D(\eta_1)}^2\mathcal D(\psi)^2,
\end{align}
where $\Lambda(\psi)$ is given by \eqref{Lambda formula} with $\mb u_{\rm lin}[\psi]$ in place of $\mb u$.

Likewise, under the hypotheses of Theorem \ref{compatible set}(ii), a function $\psi\in C_c^\infty((0,\infty))$ with $\mathcal D(\psi)\leq\delta''$ lies in $B''$, and the point of the second graph with this $\psi$ satisfies \eqref{C1} and \eqref{C2 explicit}. It differs from $\psi$ by a linear combination, with coefficients of order $\mathcal D(\psi)$, of $\eta_1$, $\eta_2$ and $\eta_3$, which near the contact point are $r\ln r$, $r^2\ln r$ and $r^2\ln^2r/2$.

A profile with compact support away from the contact point has $b_1=0$, so that it satisfies \eqref{C1} exactly when $u_1(\mb 0)=0$, that is, when it leaves the contact point initially at rest. Such profiles are obtained from every small $\psi\in C_c^\infty((0,\infty))$ by a correction along the profile $\psi_\Lambda$ of \eqref{block representatives}. Let $\mathcal D(\psi)$ be small, and let $z_\Lambda\in\R$ with $|z_\Lambda|\leq|\Lambda(\psi)|+\mathcal D(\psi)$. Then $\psi+z_\Lambda\psi_\Lambda$ lies in $\mathcal Z$, with $\mathcal D(\psi+z_\Lambda\psi_\Lambda)\lesssim\mathcal D(\psi)$ by the bound $|\Lambda(\psi)|\lesssim_{\theta_b}\mathcal D(\psi)$ of Lemma \ref{corner value response}(ii), so that by Theorem \ref{compatible set}(i) it satisfies \eqref{C1} if and only if $g_1(\psi+z_\Lambda\psi_\Lambda)=0$. By \eqref{first-order coefficient} with $\psi+z_\Lambda\psi_\Lambda$ in place of $\psi$, and since $\Lambda(\psi_\Lambda)=1$,
\begin{align*}
|\theta_b(\sin\theta_b)^2g_1(\psi+z_\Lambda\psi_\Lambda)-\Lambda(\psi)-z_\Lambda|\leq C_\Lambda\mathcal D(\psi)^2,
\end{align*}
with $C_\Lambda$ depending only on $\theta_b$, $\kappa$, $l$ and the representatives \eqref{block representatives}. If $\psi\neq 0$ and $2C_\Lambda\mathcal D(\psi)\leq 1$, the numbers $z_\Lambda^\pm:=-\Lambda(\psi)\pm 2C_\Lambda\mathcal D(\psi)^2$ lie in the range of $z_\Lambda$, and $g_1(\psi+z_\Lambda^\pm\psi_\Lambda)$ has the sign $\pm$. Since $g_1$ is Lipschitz, by Theorem \ref{compatible set}(i), the intermediate value theorem gives a zero $z_\Lambda\in[z_\Lambda^-,z_\Lambda^+]$ of $g_1(\psi+z_\Lambda\psi_\Lambda)$, and for this $z_\Lambda$ the profile $\psi+z_\Lambda\psi_\Lambda$ satisfies \eqref{C1}. For $\psi\in\ker\Lambda$ the correction $z_\Lambda\psi_\Lambda$ is therefore of second order.
\end{remark}


\subsection{Construction of the positive-time solution}\label{positive-time construction}

With the initial jet of Theorem \ref{nonlinear initial data} constructed, it remains to solve System \ref{reference wedge system} on a time interval $[0,T]$. The scheme is again a contraction, now in the space-time norms $\nn\ph\nn$ of Theorem \ref{linear solution full regularity}, whose hypotheses on $(\phi_0,\dot\phi_0,\ddot\phi_0)$ are what Theorem \ref{nonlinear initial data} provides; the nonlinear estimates are the space-time analogues of Proposition \ref{nonlinear estimates}. The smallness needed to close the argument now comes also from the shortness of the time interval, in addition to the smallness of the data that Lemma \ref{regularised flattening} requires in any case.

Throughout this subsection we fix $\kappa\in(1/4,1/2)$, which is the range in which the linear theory of Sections \ref{linear existence} and \ref{corner expansions} operates, and an integer $l\geq 1$. We also fix the auxiliary exponent $\varpi$ of the linear theory as
\begin{align}\label{varpi choice}
\varpi:=\min\set{\kappa,1-2\kappa}\in(0,\kappa],
\end{align}
the reason for this choice being explained in the proof of Proposition \ref{space-time nonlinear estimates}. We work under the standing angle condition \eqref{standing angle condition} and assume that $\theta_b$ is small enough for the linear theory to apply. The hypothesis $(3+\tilde\kappa-\kappa)\theta_b<\pi$ of Theorem \ref{linear solution full regularity} then holds for every $\tilde\kappa\in[0,\kappa)$, by $(2-\kappa/2)\theta_b<\pi/2$ in \eqref{standing angle condition}, since $3+\tilde\kappa-\kappa<3<2(2-\kappa/2)$. The strengthened angle condition \eqref{third rung angle condition} and the condition $\mathcal D_+<\infty$ of Theorem \ref{jet corner structure} are needed only in Proposition \ref{positive-time regularity upgrade}, where they are stated, and in the results that use its conclusion \eqref{upgrade estimate}.

The exponential weight $\varsigma>0$ of Definition \ref{time norms} is fixed, and the time interval is $[0,T]$ with $T>0$; Corollary \ref{polynomial lifespan} chooses $\varsigma$ in terms of $T$. We write
\begin{align}\label{Theta definition}
\Theta_T:=e^{\varsigma T}T^{1\over 2},
\end{align}
which is the constant of the embedding \eqref{sup in time 1}, and through which alone the length of the time interval will enter the estimates. All implicit constants are allowed to depend on $\theta_b$, $\kappa$, $\tilde\kappa$ and $l$, with the same absorption conventions as in the last subsection, but never on $T$ or on the data, and not on $\varsigma$ unless this is displayed or stated. The dependence on $\varsigma$ is tracked where it matters for the lifespan, in Theorem \ref{nonlinear positive-time existence} and Corollary \ref{polynomial lifespan}; in the nonlinear estimates below $\varsigma$ enters only through $\Theta_T$.

We are given an initial profile $\phi_0$, with $\phi_0(0)=0$ and $\lim_{r\to\infty}\phi_0=0$ as in the last subsection, whose size we measure by the data norm $\mathcal D$ of \eqref{data norm}, at the same parameter $l$ as the initial-time theory. We assume $\mathcal D\leq\delta$, where $\delta$ is at most the constant of Theorem \ref{nonlinear initial data} and will be shrunk further below, and that $\phi_0$ satisfies the two compatibility conditions \eqref{C1} and \eqref{C2} of that theorem, that is, its hypotheses up to part (iii). We write
\begin{align*}
(\mb u_0,p_0,\dot{\mb u}_0,\dot p_0,\dot\phi_0,\ddot\phi_0)
\end{align*}
for the tuple produced by Theorem \ref{nonlinear initial data}, which that theorem denotes $(\mb u,p,\dot{\mb u},\dot p,\dot\phi,\ddot\phi)$; the unsubscripted letters $(\phi,\mb u,p)$ are from now on reserved for the evolving solution. By \eqref{main estimate} the tuple is bounded by $\mathcal D$. We abbreviate the values of the nonlinearities on this tuple by
\begin{align}\label{initial forcing values}
(\mb N_0,\mb N_{\partial,0}):=(\mb N,\mb N_\partial)[\phi_0,\mb u_0,p_0],\qquad
(\dot{\mb N}_0,\dot{\mb N}_{\partial,0}):=(D\mb N,D\mb N_\partial)[\phi_0,\mb u_0,p_0](\dot\phi_0,\dot{\mb u}_0,\dot p_0).
\end{align}
By construction, $(\mb u_0,p_0)$ is the solution of the stationary system \eqref{system: stationary linear} with data $(\phi_0,\mb N_0,\mb N_{\partial,0})$ and satisfies the compatibility condition \eqref{C1 in the linear theory}, which is \eqref{C1}; $(\dot{\mb u}_0,\dot p_0)$ is the solution with data $(\dot\phi_0,\dot{\mb N}_0,\dot{\mb N}_{\partial,0})$; $\dot\phi_0$ and $\ddot\phi_0$ are the profiles derived from the two pairs; and $\ddot\phi_0$ satisfies \eqref{C2}. The triple $(\phi_0,\dot\phi_0,\ddot\phi_0)$ is therefore the initial jet of Theorem \ref{initial jet}, with the hypothesis \eqref{improved theorem C2} of Theorem \ref{linear solution full regularity}, for the forcing values \eqref{initial forcing values}.

\subsubsection{Space-time norms and the iteration space}

The iteration runs in the norms of the rungs $M=0$ and $M=1$ of Theorem \ref{linear solution full regularity} at $\tilde\kappa=0$; the rung $M=2$ is recovered for the solution in Proposition \ref{positive-time regularity upgrade}.

\begin{definition}[Space-time norms and the iteration space]\label{space-time norms}
With the time norms of Definition \ref{time norms}, for a triple $(\psi,\mb v,q)$ on $[0,T]$ we define
\begin{align}
\mathcal E_T(\psi,\mb v,q)&:=(\sin\theta_b)^2\nn\psi\nn_{[0,T],\partial_1\Omega:(\varsigma,{1\over 2}-\kappa)(2,l+{3\over 2})}+\nn\mb v\nn_{[0,T],\Omega:(\varsigma,1-\kappa)(2,l,2)}+\nn q\nn_{[0,T],\Omega:(\varsigma,-\kappa)(2,l,1)}\nonumber\\
&\quad+(\sin\theta_b)^2\nn\psi\nn_{[0,T],\partial_1\Omega:(\varsigma,{3\over 2}-\kappa)(1,l+{5\over 2})}+\nn\mb v^{[1]}\nn_{[0,T],\Omega:(\varsigma,2-\kappa)(1,l+1,2)}\nonumber\\
&\quad+\nn q\nn_{[0,T],\Omega:(\varsigma,1-\kappa)(1,l+1,1)},\label{ET definition}
\end{align}
where $\mb v^{\smash[t]{[1]}}=\mb v-\mb v(\mb 0,t)$ as in Proposition \ref{time for radial trade}, and for a pair of forcings $(\mb F,\mb f)$ on $[0,T]$
\begin{align}
\mathcal N_T(\mb F,\mb f)&:=\sum_{k=0}^1\brac{\nn\mb F\nn_{[0,T],\Omega:(\varsigma,-\kappa k-1)(2,l)}+\nn\mb f\nn_{[0,T],\partial_1\Omega:(\varsigma,-\kappa k-{1\over 2})(2,l+{1\over 2})}}\nonumber\\
&\quad+\nn\mb F\nn_{[0,T],\Omega:(\varsigma,\varpi-1)(2)}+\nn\mb f\nn_{[0,T],\partial_1\Omega:(\varsigma,\varpi-{1\over 2})(2,{1\over 2})}\nonumber\\
&\quad+\nn\mb F\nn_{[0,T],\Omega:(\varsigma,-\kappa)(1,l+1)}+\nn\mb f\nn_{[0,T],\partial_1\Omega:(\varsigma,-\kappa+{1\over 2})(1,l+{3\over 2})}.\label{NT definition}
\end{align}
Given a radius $\rho\geq 0$, the iteration space is
\begin{align}
B_T:=\Big\{(\psi,\mb v,q):\ &\mathcal E_T(\psi,\mb v,q)\leq\rho,\quad\grad\cdot\mb v=0,\quad v_\theta|_{\partial_0\Omega}=0,\nonumber\\
&(\psi,\partial_t\psi,\mb v,\partial_t\mb v,q,\partial_tq)(0)=(\phi_0,\dot\phi_0,\mb u_0,\dot{\mb u}_0,p_0,\dot p_0)\Big\}.\label{iteration ball}
\end{align}
On triples we use the distance
\begin{align}\label{positive time metric}
d_T\brac{(\psi,\mb v,q),(\psi',\mb v',q')}:=\mathcal E_T(\psi-\psi',\mb v-\mb v',q-q').
\end{align}
\end{definition}

The norms $\mathcal E_T$ and $\mathcal N_T$ are the LHS and the forcing norms on the right of \eqref{Improved solution estimate} at $\tilde\kappa=0$, whose data terms are bounded by the data norm $\mathcal D$ of \eqref{data norm} and by \eqref{main estimate} at $M=1$. Hence, whenever Theorem \ref{linear solution full regularity} applies with the initial jet $(\phi_0,\dot\phi_0,\ddot\phi_0)$, its solution satisfies
\begin{align}\label{linear estimate renamed}
\mathcal E_T(\phi,\mb u,p)\leq C_{{\rm lin},\varsigma}\brac{\mathcal N_T(\mb F,\mb f)+\mathcal D},
\end{align}
with a constant $C_{{\rm lin},\varsigma}\geq 1$ of the form $C(\theta_b,\kappa,l)(1+\varsigma^{-1})$, which we fix for the rest of this subsection.

The values at $t=0$ in $B_T$ are well defined by Lemma \ref{supremum in time}, since each of $\psi,\partial_t\psi,\mb v,\partial_t\mb v,q,\partial_tq$ occurs in $\mathcal E_T$ together with its time derivative. Unlike the size $\mathcal E$ of the initial-time theory, see Remark \ref{metric remark}, $\mathcal E_T$ is a norm and $d_T$ its metric: the only subtraction in $\mathcal E_T$ is the corner value $\mb v(\mb 0,t)$, a linear functional of the triple.

\subsubsection{Space-time estimates of the nonlinearities}

The estimates of this subsection are the space-time counterparts of Lemma \ref{geometry estimates} and Proposition \ref{nonlinear estimates}. The initial data and the factor $\Theta_T$ enter them only through \eqref{sup in time 1}. For a triple satisfying the initial conditions of $B_T$ we write throughout
\begin{align}\label{K_T definition}
\mathcal K_T:=\mathcal D+\Theta_T\mathcal E_T(\psi,\mb v,q),
\end{align}
which is the bound that \eqref{sup in time 1} gives for suprema in time along the iteration, the values at $t=0$ being norms of $(\phi_0,\dot\phi_0,\mb u_0,\dot{\mb u}_0,p_0,\dot p_0)$, bounded by $\mathcal D$ through \eqref{main estimate} and Lemma \ref{geometry estimates}(ii)--(iii).

\begin{lemma}[Space-time geometry estimates]\label{space-time geometry estimates}
Let $(\psi,\mb v,q)$ be a triple with $\mathcal E_T=\mathcal E_T(\psi,\mb v,q)<\infty$ satisfying the initial conditions of $B_T$, and let $\bs\nu(t)=\bs\nu[\psi(t)]$.
\begin{enumerate}
\item (Profile norms and corner value) We have
\begin{align}\label{spacetime coefficients}
\|\mb v(\mb 0,\ph)\|_{L^2_\varsigma[0,T]}\lesssim_{\theta_b}\mathcal E_T,\qquad\qquad
\sup_{t\in[0,T]}|\mb v(\mb 0,t)|\lesssim_{\theta_b}\mathcal K_T,
\end{align}
where $\|\ph\|_{L^2_\varsigma[0,T]}$ is the $L^2[0,T]$ norm against the weight $e^{-\varsigma t}$. Moreover,
\begin{align}\label{profile spacetime ladder}
\nn\psi\nn_{[0,T],\partial_1\Omega:(\varsigma,a-{1\over 2})(0,k)}\lesssim_{\theta_b}\mathcal E_T,\qquad\qquad
\sup_{t\in[0,T]}\|\psi(t)\|_{\partial_1\Omega,a-{1\over 2},k'}\lesssim_{\theta_b}\mathcal K_T,
\end{align}
for the weights $a$ and the counts $k$ and $k'$ in the first row of the following table, and the same bounds hold for $\partial_t\psi$ and $\partial_t^2\psi$ in the second and third rows:
\begin{align*}
\begin{array}{c|c|c}
&\text{space-time norm}&\text{supremum in time}\\\hline
\psi&a\in[1-\kappa,2-\kappa],\ k\leq l+{7\over 2}&a\in[1-\kappa,2-\kappa],\ k'\leq l+{5\over 2}\\
\partial_t\psi&a\in[1-\kappa,2-\kappa],\ k\leq l+{5\over 2}&a=1-\kappa,\ k'\leq l+{3\over 2}\\
\partial_t^2\psi&a=1-\kappa,\ k\leq l+{3\over 2}&\text{---}
\end{array}
\end{align*}
\item (Geometry) We have
\begin{align}
\nn\bs\nu\nn_{[0,T],\Omega:(\varsigma,a)(0,0,m)}+\norm{\|r^{1-a}(r\grad)^{m'}\bs\nu\|_{L^\infty(\Omega)}}_{L^2_\varsigma[0,T]}&\lesssim_{\theta_b}\mathcal E_T,\label{space time geometry L2}\\
\sup_{t\in[0,T]}\|\bs\nu(t)\|_{\Omega,a,0,m}+\sup_{t\in[0,T]}\|r^{1-a}(r\grad)^{m'}\bs\nu(t)\|_{L^\infty(\Omega)}&\lesssim_{\theta_b}\mathcal K_T,\label{space time geometry sup}
\end{align}
for the weights $a$ and the counts $m$ and $m'$ in the first row of the following table, and the same bounds hold for $\bs\nu[\partial_t\psi]$ and $\bs\nu[\partial_t^2\psi]$ in the second and third rows:
\begin{align*}
\begin{array}{c|c|c}
&\eqref{space time geometry L2}&\eqref{space time geometry sup}\\\hline
\bs\nu&a\in[1-\kappa,2-\kappa],\ m\leq l+4,\ m'\leq l+2&a\in[1-\kappa,2-\kappa],\ m\leq l+3,\ m'\leq l+1\\
\bs\nu[\partial_t\psi]&a\in[1-\kappa,2-\kappa],\ m\leq l+3,\ m'\leq l+1&a=1-\kappa,\ m\leq l+2,\ m'\leq l\\
\bs\nu[\partial_t^2\psi]&a=1-\kappa,\ m\leq l+2,\ m'\leq l&\text{---}
\end{array}
\end{align*}
\item (Smallness) There is $\delta_1=\delta_1(\theta_b,\kappa,l)>0$, not depending on $\varsigma$, such that $\mathcal K_T\leq\delta_1$ implies
\begin{align}\label{space time geometry smallness}
\sup_{t\in[0,T]}\brac{\|\bs\nu(t)\|_{L^\infty(\Omega)}+\|\psi(t)\|_{L^\infty(\partial_1\Omega)}+\|\chi'(t)\|_{\partial_1\Omega,{1\over 2},{3\over 2}}}\leq\min\set{{1\over 2},\epsilon_{\bs\Psi}(\theta_b)},
\end{align}
with $\epsilon_{\bs\Psi}(\theta_b)$ the constant of Lemma \ref{regularised flattening}. Moreover, if $\mathcal K_T\leq\delta_1$, then at every $t\in[0,T]$ the map $\bs\Psi_t$ is a diffeomorphism $\Omega\to\Omega_t$, the transformation of Lemma \ref{Piola transform} and the nonlinearities of System \ref{reference wedge system} are defined, and the factors $\mathfrak G$ obey \eqref{composition} and \eqref{composition bounds} with $\mathcal D$ replaced by $\mathcal K_T$ at the counts of \eqref{space time geometry sup} and by $\mathcal E_T$ at the counts of \eqref{space time geometry L2}.
\end{enumerate}
\end{lemma}
\begin{proof}
(i) Lemma \ref{decomposition estimates}(ii) holds at each fixed $t$, and the first bound of \eqref{spacetime coefficients} follows by integration in time. The second follows by \eqref{sup in time 1} from the first bound for $\partial_t\mb v$, whose norms are likewise contained in $\mathcal E_T$, the initial value being $|\mb u_0(\mb 0)|\lesssim_{\theta_b}\mathcal D$ by \eqref{corner coefficients 1} and \eqref{main estimate}. The space-time bound of \eqref{profile spacetime ladder} follows from $\nn\psi\nn_{[0,T],\partial_1\Omega:(\varsigma,1/2-\kappa)(2,l+3/2)}$ and $\nn\psi\nn_{[0,T],\partial_1\Omega:(\varsigma,3/2-\kappa)(1,l+5/2)}$ in $\mathcal E_T$. The supremum bound follows by \eqref{sup in time 1}, the initial value being $\|\phi_0\|_{\partial_1\Omega,a-1/2,k'}\lesssim_{\theta_b}\mathcal D$ by Lemma \ref{geometry estimates}(ii), and the derivative term $\Theta_T\nn\partial_t\psi\nn_{[0,T],\partial_1\Omega:(\varsigma,a-1/2)(0,k')}$ being $\lesssim_{\theta_b}\Theta_T\mathcal E_T$ by the space-time bound for $\partial_t\psi$. The bounds for $\partial_t\psi$ and $\partial_t^2\psi$ follow in the same way, with the initial value $\dot\phi_0$ bounded by \eqref{main estimate}; $\partial_t^2\psi$ occurs in $\mathcal E_T$ only in $\nn\psi\nn_{[0,T],\partial_1\Omega:(\varsigma,1/2-\kappa)(2,l+3/2)}$, which restricts the supremum of $\partial_t\psi$ to $a=1-\kappa$.

(ii) Given (i), the bounds follow by applying \eqref{geometry bulk} and \eqref{geometry sup} at each fixed $t$, as in Lemma \ref{geometry estimates}: the boundary weight is $a-{1\over 2}$, and in each row of the table the counts $\max\{m-{1\over 2},0\}$ and $m'+{3\over 2}$ are at most the count $k$ of (i) for \eqref{space time geometry L2} and at most $k'$ for \eqref{space time geometry sup}.

(iii) By \eqref{space time geometry sup} at $a=1$, $m'=0$, by Lemma \ref{embeddings}(i) at $a=0$, and by the supremum bound of \eqref{profile spacetime ladder} at $a=1$, $k'={3\over 2}$, recalling $\dr\chi=-(\sin\theta_b)r\psi$, the LHS of \eqref{space time geometry smallness} is $\lesssim_{\theta_b}\mathcal K_T$ with a constant that does not involve $\varsigma$, which enters (i) and (ii) only through $\Theta_T$ inside $\mathcal K_T$, and we take $\delta_1\leq 1$ small enough. Then $\bs\Psi_t$ is a diffeomorphism by Lemma \ref{regularised flattening}, so that Lemma \ref{Piola transform} applies, and $\|\bs\nu(t)\|_{L^\infty(\Omega)}\leq 1/2$, so that Lemma \ref{composition of the geometry}(i) applies at each $t$. For the bounds on the factors $\mathfrak G$, by \eqref{composition} and the Leibniz rule $(r\grad)^m\mathfrak G_0$ is a finite sum of terms $\mathfrak Gf_1\cdots f_J$, where each $f_{j'}$ is a component of $(r\grad)^{m_{j'}}\bs\nu$, $m_1+\dots+m_J=m$ and $m_1\geq m_{j'}$, so that $m_{j'}\leq m/2$ for $j'\geq 2$. For $m\leq l+4$ we place $f_1$ in \eqref{space time geometry L2} and the others, with $m_{j'}\leq\lfloor(l+4)/2\rfloor\leq l+1$ since $l\geq 1$, in the suprema in time of \eqref{space time geometry sup}:
\begin{align*}
\nn\mathfrak Gf_1\cdots f_J\nn_{[0,T],\Omega:(\varsigma,a)}&\leq\sup_{t\in[0,T]}\|\mathfrak G(t)\|_{L^\infty(\Omega)}\nn(r\grad)^{m_1}\bs\nu\nn_{[0,T],\Omega:(\varsigma,a)}\prod_{j'=2}^J\sup_{t\in[0,T]}\|(r\grad)^{m_{j'}}\bs\nu(t)\|_{L^\infty(\Omega)}\\
&\lesssim_{\theta_b,\mathfrak g}\mathcal E_T\mathcal K_T^{J-1}\leq\mathcal E_T,
\end{align*}
using $\mathcal K_T\leq\delta_1\leq 1$. The other bounds follow in the same way.
\end{proof}

\begin{proposition}[Space-time nonlinear estimates]\label{space-time nonlinear estimates}
Let $(\psi,\mb v,q)$ and $(\psi',\mb v',q')$ be triples satisfying $\grad\cdot\mb v=\grad\cdot\mb v'=0$, $v_\theta|_{\partial_0\Omega}=v'_\theta|_{\partial_0\Omega}=0$ and the initial conditions of $B_T$, with $\mathcal E_T=\mathcal E_T(\psi,\mb v,q)$ and $\mathcal E_T'=\mathcal E_T(\psi',\mb v',q')$ finite, and assume $\mathcal K_T\leq\delta_1$, with $\delta_1$ as in Lemma \ref{space-time geometry estimates}(iii). Then $(\mb N,\mb N_\partial)[\psi,\mb v,q]$, formed with the time-dependent profile $\psi(t)$, satisfies
\begin{align}
\mathcal N_T\brac{\mb N[\psi,\mb v,q],\mb N_\partial[\psi,\mb v,q]}&\lesssim_{\theta_b}\brac{\mathcal D+\Theta_T\mathcal E_T}\mathcal E_T,\label{ST N estimate}\\
\intertext{and, if also $\mathcal K_T':=\mathcal D+\Theta_T\mathcal E_T'\leq\delta_1$,}
\mathcal N_T\brac{\mb N[\psi,\mb v,q]-\mb N[\psi',\mb v',q'],\ \mb N_\partial[\psi,\mb v,q]-\mb N_\partial[\psi',\mb v',q']}&\lesssim_{\theta_b}\brac{\mathcal D+\Theta_T\brac{\mathcal E_T+\mathcal E_T'}}d_T,\label{ST N difference}
\end{align}
where $d_T=d_T((\psi,\mb v,q),(\psi',\mb v',q'))$. The bound \eqref{ST N estimate} also holds for $\nn H\nn_{[0,T],\Omega:(\varsigma,1-\kappa)(0,l+2,1)}$, where $H$ is the bulk expression of \eqref{schematic boundary} whose trace is $\mb N_\partial[\psi,\mb v,q]$.
\end{proposition}
\begin{proof}
Step 1 reduces $\mathcal N_T$ to targets of the form of those of Proposition \ref{nonlinear estimates}, carrying up to two time derivatives, and Step 2 places the time derivatives. Steps 3 and 4 estimate the bulk and the boundary targets at the rungs $M\leq 1$, Step 5 the corner targets, and Step 6 the difference \eqref{ST N difference}.

By Lemma \ref{space-time geometry estimates}(iii) the nonlinearities are defined along $(\psi,\mb v,q)$, and along $(\psi',\mb v',q')$ when $\mathcal K_T'\leq\delta_1$. In this proof we abbreviate
\begin{align*}
\nn u\nn_{O,a,j,m}:=\nn u\nn_{[0,T],O:(\varsigma,a)(0,j,m)}=\norm{\|u(t)\|_{O,a,j,m}}_{L^2_\varsigma[0,T]},
\end{align*}
where $O$ is $\Omega$, $\partial_1\Omega$ or $\Omega\cap\{r<1\}$ and vanishing counts are omitted as for the spatial norms.

\emph{Step 1: the targets.} By Definition \ref{time norms}, $\mathcal N_T(\mb N,\mb N_\partial)$ is at most the sum, over $a\in\{-1-\kappa,-1\}$ with $n\leq 2$ and $a=-\kappa$ with $n\leq 1$, of the norms $\nn\partial_t^n\mb N\nn_{\Omega,a,l+2-n}$ and $\nn\partial_t^n\mb N_\partial\nn_{\partial_1\Omega,a+1/2,l+5/2-n}$, and of the same norms at $a=\varpi-1$ with the counts $2-n$ and $5/2-n$, $n\leq 2$. By Lemma \ref{trace}, applied at each $t$,
\begin{align*}
\nn\partial_t^n\mb N_\partial\nn_{\partial_1\Omega,a+{1\over 2},l+{5\over 2}-n}\lesssim_{\theta_b}\nn\partial_t^nH\nn_{\Omega,a+1,l+2-n,1},
\end{align*}
and likewise with the count $5/2-n$ on the left and $(2-n,1)$ on the right. For $n\leq 1$ every weight now lies between the $M=0$ and $M=1$ weights, $-1-\kappa$ and $-\kappa$ for $\mb N$, $-\kappa$ and $1-\kappa$ for $H$. For $n=2$ the $M=1$ weight is absent from $\mathcal N_T$, and we split at $r=1$: for $a\in\{-1,\varpi-1\}$,
\begin{align*}
\nn\partial_t^2\mb N\nn_{\Omega,a,l}&\leq\nn\partial_t^2\mb N\nn_{\Omega,-1-\kappa,l}+\nn\partial_t^2\mb N\nn_{\Omega\cap\{r<1\},\varpi-1,l},\\
\nn\partial_t^2H\nn_{\Omega,a+1,l,1}&\leq\nn\partial_t^2H\nn_{\Omega,-\kappa,l,1}+\nn\partial_t^2H\nn_{\Omega\cap\{r<1\},\varpi,l,1}.
\end{align*}
Hence
\begin{align}\label{space-time targets}
\mathcal N_T\brac{\mb N,\mb N_\partial}\lesssim_{\theta_b}{}&\sum_{M=0}^1\sum_{n=0}^{2-M}\brac{\nn\partial_t^n\mb N\nn_{\Omega,M-1-\kappa,l+2-n}+\nn\partial_t^nH\nn_{\Omega,M-\kappa,l+2-n,1}}\nonumber\\
&+\nn\partial_t^2\mb N\nn_{\Omega\cap\{r<1\},\varpi-1,l}+\nn\partial_t^2H\nn_{\Omega\cap\{r<1\},\varpi,l,1}.
\end{align}
The first line collects the targets of \eqref{N estimate} and \eqref{boundary bulk target} at $M\leq 1$, with the counts lowered by the number $n$ of time derivatives. The second line collects the corner targets, which lie $\varpi+\kappa$ above the $M=0$ targets; the choice \eqref{varpi choice} of $\varpi$ enters only in their estimate. No weight beyond the $M=1$ weights occurs, so the decompositions $\mb v^{[2]}$ and $q^{[2]}$ of the stationary proof are not needed.

\emph{Step 2: the placement in time.} By the linearity of $\psi\mapsto\bs\nu[\psi]$ and the chain rule, as in the second identity of \eqref{composition},
\begin{align}\label{time chain rule}
\partial_t\bs\nu=\bs\nu[\partial_t\psi],\qquad
\partial_t\mathfrak G=\brac{\nu_1[\partial_t\psi]}\mathfrak G+\brac{\nu_2[\partial_t\psi]}\mathfrak G,\qquad
\partial_t\mb N[\psi,\mb v,q]=D\mb N[\psi,\mb v,q](\partial_t\psi,\partial_t\mb v,\partial_tq),
\end{align}
and likewise for $\mb N_\partial$, for $H$ and for the second time derivatives. Hence $\partial_t^n\mb N$ and $\partial_t^nH$ are sums of summands of \eqref{schematic bulk}--\eqref{schematic boundary} whose factors carry in addition $n$ time derivatives in total, either all on one factor or, when $n=2$, one on each of two factors; a time derivative falling on a factor $\mathfrak G$ produces a component of $\bs\nu[\partial_t\psi]$. After the Leibniz expansion, the proof of Proposition \ref{nonlinear estimates} at $M\leq 1$, at the count lowered by $n$, bounds each resulting product at each fixed $t$ by one factor in a weighted $L^2$ norm in space times the others in weighted $L^\infty$ norms in space. We integrate in time by
\begin{align*}
\nn XY\nn_{O,a}\leq\sup_{t\in[0,T]}\|X(t)\|_{\mathcal X}\norm{\|Y(t)\|_{\mathcal Y}}_{L^2_\varsigma[0,T]},
\end{align*}
valid whenever $\|X(t)Y(t)\|_{O,a}\leq\|X(t)\|_{\mathcal X}\|Y(t)\|_{\mathcal Y}$ at every $t$ for norms $\|\ph\|_{\mathcal X}$ and $\|\ph\|_{\mathcal Y}$ in space, where $Y$ is one of the factors, which we call distinguished, and $X$ is the product of the others. Thus the distinguished factor is measured in a space-time norm, $L^2_\varsigma$ in time, and the others in the supremum in time, each factor being in $L^2$ or in $L^\infty$ in space. As a rule the distinguished factor carries the most time derivatives and, among those, the most spatial derivatives, with one exception: if two factors carry one time derivative each and one of them is a geometry factor, a geometry factor among these two is distinguished. Which factor is placed in $L^\infty$ in space is decided by the counts, independently of this choice. The distinguished factor contributes $\mathcal E_T$, by \eqref{ET definition} and \eqref{space time geometry L2}, and the others contribute $\mathcal K_T$: the geometry by \eqref{space time geometry sup}, and the solution by
\begin{align}\label{solution sup}
&\sup_{t\in[0,T]}\brac{\|\mb v(t)\|_{\Omega,1-\kappa,l+1,2}+\|\mb v^{[1]}(t)\|_{\Omega,2-\kappa,l+1,2}+\|q(t)\|_{\Omega,-\kappa,l+1,1}+\|q(t)\|_{\Omega,1-\kappa,l+1,1}}\nonumber\\
&+\sup_{t\in[0,T]}\brac{\|\partial_t\mb v(t)\|_{\Omega,1-\kappa,l,2}+\|\partial_tq(t)\|_{\Omega,-\kappa,l,1}}\lesssim_{\theta_b}\mathcal K_T,
\end{align}
which follows from \eqref{sup in time 1}, the values at $t=0$ being norms of $(\mb u_0,p_0,\dot{\mb u}_0,\dot p_0)$, bounded by \eqref{main estimate}, and the norms of the next time derivative being contained in $\mathcal E_T$. The factors $\mathfrak G$, and any factor of $\bs\nu$ type beyond the first carrying no derivative in space or in time, are bounded by \eqref{space time geometry smallness} and suppressed, as in the stationary proof.

The counts are within budget. A factor carrying $n'$ time derivatives is available at its count in the stationary proof lowered by $n'$ in the space-time norms, and lowered by $n'+1$ in the supremum in time, while the target lowers the total count by $n$. A factor in the supremum in time carries $n'<n$ time derivatives, unless $n=0$, when it carries at most half of the spatial derivatives, which fits the supremum in time for $l\geq 1$, as in Lemma \ref{space-time geometry estimates}(iii).

The weights at which a factor carrying $n'$ time derivatives is available are listed below, from Lemma \ref{space-time geometry estimates}(ii), $\mathcal E_T$ and \eqref{solution sup}. For a solution factor, $M=0$ and $M=1$ refer to its weights on the two rungs of $\mathcal E_T$, $1-\kappa$ for $\mb v$ and $2-\kappa$ for $\mb v^{[1]}$, respectively $-\kappa$ and $1-\kappa$ for $q$, every weight in between being then available:
\begin{align}\label{available weights}
\begin{array}{c|cc|cc}
&\multicolumn{2}{c|}{\text{solution factor}}&\multicolumn{2}{c}{\text{geometry factor}}\\
n'&\text{space-time norm}&\text{supremum in time}&\text{space-time norm}&\text{supremum in time}\\\hline
0&M\in\{0,1\}&M\in\{0,1\}&a\in[1-\kappa,2-\kappa]&a\in[1-\kappa,2-\kappa]\\
1&M\in\{0,1\}&M=0&a\in[1-\kappa,2-\kappa]&a=1-\kappa\\
2&M=0&\text{---}&a=1-\kappa&\text{---}
\end{array}
\end{align}
The $M=1$ weight applies to every solution factor that occurs: the derivatives of $\mb v$ are those of $\mb v^{\smash[t]{[1]}}$, the only undifferentiated velocity factor is $v_y=v_y^{\smash[t]{[1]}}$, since $v_y(\mb 0,t)=0$ by Lemma \ref{decomposition estimates}(iii), and $q^{\smash[t]{[1]}}=q$.

\emph{Step 3: the bulk targets.} We estimate the terms $\nn\partial_t^n\mb N\nn_{\Omega,M-1-\kappa,l+2-n}$ of \eqref{space-time targets} on the summands $\bs\nu\grad^2\mb v$ of Case 1 and $(\grad^2\bs\nu)v_y$ of Case 2 of the stationary proof. By the Leibniz rule and \eqref{time chain rule} they produce the terms
\begin{align*}
(\dr^{j_1}\bs\nu[\partial_t^{n_1}\psi])(\dr^{j_2}\grad^2\partial_t^{n_2}\mb v),\qquad
(\dr^{j_1}\grad^2\bs\nu[\partial_t^{n_1}\psi])(\dr^{j_2}\partial_t^{n_2}v_y),\qquad
j_1+j_2\leq l+2-n,\quad n_1+n_2=n.
\end{align*}
In the terms from Case 1 the geometry factor is placed in $L^\infty$ in space and in those from Case 2 the solution factor, through Lemma \ref{embeddings}(ii), as in the stationary proof; the other factor is in $L^2$ in space. The pressure terms of Case 3 are placed as those from Case 1, with $\grad q$ and the weights of $q$ in place of $\grad^2\mb v$ and those of $\mb v$. The summands $(\grad\bs\nu)\grad\mb v$, $(\grad\bs\nu)\grad v_y$ and $(\grad\bs\nu)(\grad\bs\nu)v_y$, and the terms with a further geometry factor produced by a time derivative on a factor $\mathfrak G$, are placed likewise, the factor in $L^\infty$ in space being decided by the counts as in the stationary proof; the one new case among them is treated at the end of this step.

When $n\geq 1$ and all time derivatives fall on the solution factor, it is distinguished, and so it is when $n=0$ and it carries the most spatial derivatives. It is then in the space-time norm at its weight $M+1-\kappa$, which \eqref{available weights} supplies for every $n\leq 2-M$, and the geometry factor, carrying no time derivative, is in the supremum in time at $a=1$:
\begin{align*}
\nn(\dr^{j_1}\bs\nu)(\dr^{j_2}\grad^2\partial_t^n\mb v)\nn_{\Omega,M-1-\kappa}
&\leq\sup_{t\in[0,T]}\|\dr^{j_1}\bs\nu(t)\|_{L^\infty(\Omega)}\nn\grad^2\partial_t^n\mb v\nn_{\Omega,M-1-\kappa,j_2},\\
\nn(\dr^{j_1}\grad^2\bs\nu)(\dr^{j_2}\partial_t^nv_y)\nn_{\Omega,M-1-\kappa}
&\leq\sup_{t\in[0,T]}\|\grad^2\bs\nu(t)\|_{\Omega,-1,j_1}\norm{\|r^{\kappa-M}\dr^{j_2}\partial_t^nv_y\|_{L^\infty(\Omega)}}_{L^2_\varsigma[0,T]}.
\end{align*}
The suprema in time are bounded by \eqref{space time geometry sup} at $a=1$, and the space-time norms by the norm of $\mb v^{\smash[t]{[M]}}$ in $\mathcal E_T$. So both RHSs are $\lesssim_{\theta_b}\mathcal K_T\mathcal E_T$, with no loss of weight, also when $n=2$.

When $n\geq 1$ and all time derivatives fall on the geometry factor, it is distinguished, and so it is when $n=0$ and it carries the most spatial derivatives. For $n\leq 1$ it is in the space-time norm at $a=1$, and the solution factor, carrying no time derivative, is in the supremum in time at its weight $M+1-\kappa$:
\begin{align*}
\nn(\dr^{j_1}\bs\nu[\partial_t^n\psi])(\dr^{j_2}\grad^2\mb v)\nn_{\Omega,M-1-\kappa}
&\leq\norm{\|\dr^{j_1}\bs\nu[\partial_t^n\psi]\|_{L^\infty(\Omega)}}_{L^2_\varsigma[0,T]}\sup_{t\in[0,T]}\|\grad^2\mb v(t)\|_{\Omega,M-1-\kappa,j_2},\\
\nn(\dr^{j_1}\grad^2\bs\nu[\partial_t^n\psi])(\dr^{j_2}v_y)\nn_{\Omega,M-1-\kappa}
&\leq\sup_{t\in[0,T]}\|r^{\kappa-M}\dr^{j_2}v_y(t)\|_{L^\infty(\Omega)}\nn\grad^2\bs\nu[\partial_t^n\psi]\nn_{\Omega,-1,j_1},
\end{align*}
the space-time norms being bounded by \eqref{space time geometry L2} at $a=1$ and the suprema in time by \eqref{solution sup} at the weight $M+1-\kappa$.

Now let $n=2$, hence $M=0$, with both time derivatives on the geometry factor. By \eqref{available weights} the geometry factor is then available at $a=1-\kappa$ only, a loss of $\kappa$ against $a=1$. The solution factor carries no time derivative and is in the supremum in time, where every weight between its two rungs is available, and raising its weight from $1-\kappa$ to $1$ recovers the loss:
\begin{align*}
\nn(\dr^{j_1}\bs\nu[\partial_t^2\psi])(\dr^{j_2}\grad^2\mb v)\nn_{\Omega,-1-\kappa}
&\leq\norm{\|r^\kappa\dr^{j_1}\bs\nu[\partial_t^2\psi]\|_{L^\infty(\Omega)}}_{L^2_\varsigma[0,T]}\sup_{t\in[0,T]}\|\grad^2\mb v(t)\|_{\Omega,-1,j_2},\\
\nn(\dr^{j_1}\grad^2\bs\nu[\partial_t^2\psi])(\dr^{j_2}v_y)\nn_{\Omega,-1-\kappa}
&\leq\sup_{t\in[0,T]}\|\dr^{j_2}v_y(t)\|_{L^\infty(\Omega)}\nn\grad^2\bs\nu[\partial_t^2\psi]\nn_{\Omega,-1-\kappa,j_1},
\end{align*}
the space-time norms being bounded by \eqref{space time geometry L2} for $\bs\nu[\partial_t^2\psi]$ at $a=1-\kappa$ and the suprema in time by \eqref{solution sup} at the weight $1$.

Now let $n=2$ with one time derivative on each factor. The geometry factor is distinguished by the rule of Step 2, at $a=1$, and the solution factor is in the supremum in time at its $M=0$ weight $1-\kappa$, which \eqref{available weights} supplies for one time derivative:
\begin{align*}
\nn(\dr^{j_1}\bs\nu[\partial_t\psi])(\dr^{j_2}\grad^2\partial_t\mb v)\nn_{\Omega,-1-\kappa}
&\leq\norm{\|\dr^{j_1}\bs\nu[\partial_t\psi]\|_{L^\infty(\Omega)}}_{L^2_\varsigma[0,T]}\sup_{t\in[0,T]}\|\grad^2\partial_t\mb v(t)\|_{\Omega,-1-\kappa,j_2},\\
\nn(\dr^{j_1}\grad^2\bs\nu[\partial_t\psi])(\dr^{j_2}\partial_tv_y)\nn_{\Omega,-1-\kappa}
&\leq\sup_{t\in[0,T]}\|r^\kappa\dr^{j_2}\partial_tv_y(t)\|_{L^\infty(\Omega)}\nn\grad^2\bs\nu[\partial_t\psi]\nn_{\Omega,-1,j_1},
\end{align*}
with no loss of weight, the space-time norms being bounded by \eqref{space time geometry L2} for $\bs\nu[\partial_t\psi]$ and the suprema in time by \eqref{solution sup} for $\partial_t\mb v$.

A loss occurs again when the two time derivatives fall on two geometry factors, as in the summand $(\grad\bs\nu)(\grad\bs\nu)v_y$ of $\lpc\mb u_\phi$, or on a geometry factor and a factor $\mathfrak G$. The one carrying fewer spatial derivatives is then in the supremum in time, where $\bs\nu[\partial_t\psi]$ is available at $a=1-\kappa$ only, and the solution factor is again raised to the weight $1$: for $j_3\leq j_1$,
\begin{align*}
&\nn(\dr^{j_1}\grad\bs\nu[\partial_t\psi])(\dr^{j_3}\grad\bs\nu[\partial_t\psi])(\dr^{j_2}v_y)\nn_{\Omega,-1-\kappa}\\
&\qquad\leq\sup_{t\in[0,T]}\|r^{1+\kappa}\dr^{j_3}\grad\bs\nu[\partial_t\psi(t)]\|_{L^\infty(\Omega)}\sup_{t\in[0,T]}\|\dr^{j_2}v_y(t)\|_{L^\infty(\Omega)}\nn\grad\bs\nu[\partial_t\psi]\nn_{\Omega,0,j_1},
\end{align*}
the supremum in time of the geometry factor being bounded by \eqref{space time geometry sup} for $\bs\nu[\partial_t\psi]$ at $a=1-\kappa$, that of $v_y$ by \eqref{solution sup} at the weight $1$, and the space-time norm by \eqref{space time geometry L2} for $\bs\nu[\partial_t\psi]$ at $a=1$.

\emph{Step 4: the boundary targets.} We estimate the terms $\nn\partial_t^nH\nn_{\Omega,M-\kappa,l+2-n,1}$ of \eqref{space-time targets}. The summands of $H$ containing a solution factor are placed in time as in Step 3 and in space as in Case 4 of the stationary proof at $M\leq 1$, the full derivative of the target falling on either factor. For instance, in the summand $(\grad\bs\nu)v_y$ with both time derivatives on the geometry factor, the full derivative on the solution factor and $j_2<j_1$, the loss of $\kappa$ is recovered as in Step 3:
\begin{align*}
\nn(\dr^{j_1}\grad\bs\nu[\partial_t^2\psi])(\dr^{j_2}(r\grad)v_y)\nn_{\Omega,-\kappa}
\leq\sup_{t\in[0,T]}\|\dr^{j_2}(r\grad)v_y(t)\|_{L^\infty(\Omega)}\nn\grad\bs\nu[\partial_t^2\psi]\nn_{\Omega,-\kappa,j_1},
\end{align*}
the supremum in time being bounded by \eqref{solution sup} at the weight $1$, through Lemma \ref{embeddings}(ii), and the space-time norm by \eqref{space time geometry L2} for $\bs\nu[\partial_t^2\psi]$ at $a=1-\kappa$.

In the pure geometry terms $r^{-1}AB$ of Case 4, which come from the summands of $H$ carrying no solution factor, $B$ is the factor carrying the most derivatives and $A$ the product of the others, all of them derivatives of components of $\bs\nu$ times bounded factors $\mathfrak G$. The stationary proof places $A$ in $L^\infty$ in space at $a=1$ and $B$ in $L^2$ in space at $a=M+1-\kappa$. When $A$ carries no time derivative, $B$ is distinguished, at $a=M+1-\kappa$, which \eqref{available weights} supplies for every $n\leq 2-M$:
\begin{align*}
\nn r^{-1}AB\nn_{\Omega,M-\kappa}\leq\sup_{t\in[0,T]}\|A(t)\|_{L^\infty(\Omega)}\nn B\nn_{\Omega,M+1-\kappa}.
\end{align*}
When $n=1$ and the time derivative falls on a factor of $A$, this factor is distinguished, at $a=1$, and $B$ is in the supremum in time:
\begin{align*}
\nn r^{-1}AB\nn_{\Omega,M-\kappa}\leq\norm{\|A\|_{L^\infty(\Omega)}}_{L^2_\varsigma[0,T]}\sup_{t\in[0,T]}\|B(t)\|_{\Omega,M+1-\kappa}.
\end{align*}
When $n=2$, hence $M=0$, and a time derivative falls on $A$, a factor of $A$ is available at $a=1-\kappa$ only: the one carrying both time derivatives, which is distinguished, or one carrying a single time derivative in the supremum in time. The weight of $B$ is then raised by $\kappa$, to $a=1$, which \eqref{available weights} allows, since $B$ carries no time derivative, or carries one and is distinguished:
\begin{align*}
\nn r^{-1}AB\nn_{\Omega,-\kappa}&\leq\norm{\|r^\kappa A\|_{L^\infty(\Omega)}}_{L^2_\varsigma[0,T]}\sup_{t\in[0,T]}\|B(t)\|_{\Omega,1}&&\text{if $B$ carries no time derivative},\\
\nn r^{-1}AB\nn_{\Omega,-\kappa}&\leq\sup_{t\in[0,T]}\|r^\kappa A(t)\|_{L^\infty(\Omega)}\nn B\nn_{\Omega,1}&&\text{if $B$ carries one}.
\end{align*}
When a time derivative falls on $A$, $A$ is bounded factor by factor: the distinguished factor, when it belongs to $A$, is in the space-time norm and the others in the supremum in time, and for $n=2$ the factor of $A$ at $a=1-\kappa$ carries $r^\kappa$. The factors of $A$, and $B$, are bounded by \eqref{space time geometry L2} and \eqref{space time geometry sup}, so that all the RHSs are $\lesssim_{\theta_b}\mathcal K_T\mathcal E_T$.

\emph{Step 5: the corner targets.} It remains to estimate the second line of \eqref{space-time targets}. Its weights $\varpi-1$ and $\varpi$ exceed by $\varpi+\kappa$ the weights $-1-\kappa$ and $-\kappa$ of the targets of the first line with $M=0$ and $n=2$, which have the same counts. On $\{r<1\}$ a weighted norm is non-decreasing in the weight, since $r^{-a}\leq r^{-a'}$ there for $a\leq a'$. As $\varpi+\kappa\leq 1-\kappa$ by \eqref{varpi choice}, this gives
\begin{align}\label{corner target reduction}
\nn\partial_t^2\mb N\nn_{\Omega\cap\{r<1\},\varpi-1,l}+\nn\partial_t^2H\nn_{\Omega\cap\{r<1\},\varpi,l,1}
\leq\nn\partial_t^2\mb N\nn_{\Omega\cap\{r<1\},-2\kappa,l}+\nn\partial_t^2H\nn_{\Omega\cap\{r<1\},1-2\kappa,l,1},
\end{align}
the weights on the right being those of the $M=0$ targets raised by $1-\kappa$.

To reach these weights we raise the weights of the factors above those of the proof of Proposition \ref{nonlinear estimates} at $M=0$, where the target weight splits exactly between the factors of every product. We call the weight of a factor in that proof its stationary weight. It is the $M=0$ weight for a solution factor, $1-\kappa$ for $\mb v$ and $-\kappa$ for $q$, and $a=1$ for a geometry factor, with one exception: in a pure geometry term $r^{-1}AB$ of Step 4, all of whose factors are geometry factors, $B$ has the stationary weight $a=1-\kappa$. The gain of a factor is its weight in the placement below minus its stationary weight. Raising the weight of one factor raises by the same amount the weight at which H\"older's inequality bounds the product, so a term is bounded at the weights on the right of \eqref{corner target reduction} once the gains of its factors add up to at least $1-\kappa$.

By \eqref{available weights} a geometry factor is available at the weights $a\in[1-\kappa,2-\kappa]$, and only at $a=1-\kappa$ when it carries two time derivatives, or one in the supremum in time. As for a solution factor, we call $a=1-\kappa$ and $a=2-\kappa$ its $M=0$ and $M=1$ weights. The gains at these weights are
\begin{align*}
\begin{array}{c|ccc|cc}
&\multicolumn{3}{c|}{\text{weight}}&\multicolumn{2}{c}{\text{gain}}\\
&\text{stationary}&M=0&M=1&M=0&M=1\\\hline
\mb v&1-\kappa&1-\kappa&2-\kappa&0&1\\
q&-\kappa&-\kappa&1-\kappa&0&1\\
B&1-\kappa&1-\kappa&2-\kappa&0&1\\
\text{other geometry factor}&1&1-\kappa&2-\kappa&-\kappa&1-\kappa
\end{array}
\end{align*}
So only a geometry factor other than $B$ can lose, by $\kappa$ at its $M=0$ weight, because its stationary weight is $a=1$ and not $a=1-\kappa$.

A term of $\partial_t^2\mb N$ or $\partial_t^2H$ contains exactly one solution factor or, if it is a pure geometry term $r^{-1}AB$, the factor $B$. In both cases it contains at least one geometry factor other than $B$ not suppressed in Step 2. A time derivative falling on a factor $\mathfrak G$ produces a further such factor, a component of $\bs\nu[\partial_t\psi]$, by \eqref{time chain rule}, so that the time derivatives fall only on the solution factor, on $B$ or on geometry factors other than $B$. A geometry factor carrying no time derivative is in the supremum in time, at $a=1$ or at its $M=1$ weight, and does not lose. In space the factors are placed as in Steps 3 and 4, the weights used are available by \eqref{available weights}, and the counts are within budget as in Step 2.

Consider first a term with a solution factor. By Step 2 the two time derivatives fall on one factor or on two, which gives the four cases below. We display the terms of the summands $(\grad^2\bs\nu)v_y$ and $(\grad\bs\nu)(\grad\bs\nu)v_y$ of $\lpc\mb u_\phi$, with $j_1+j_2\leq l$, respectively $j_1+j_2+j_3\leq l$.

Let first both time derivatives fall on the solution factor, as in the terms $(\dr^{j_1}\grad^2\bs\nu)(\dr^{j_2}\partial_t^2v_y)$. The solution factor is distinguished, at its $M=0$ weight, and gains nothing. The geometry factors carry no time derivative and are in the supremum in time at their $M=1$ weights, each gaining $1-\kappa$:
\begin{align}\label{corner target solution}
\nn(\dr^{j_1}\grad^2\bs\nu)(\dr^{j_2}\partial_t^2v_y)\nn_{\Omega\cap\{r<1\},-2\kappa}
\leq\sup_{t\in[0,T]}\|\grad^2\bs\nu(t)\|_{\Omega,-\kappa,j_1}\norm{\|r^\kappa\dr^{j_2}\partial_t^2v_y\|_{L^\infty(\Omega)}}_{L^2_\varsigma[0,T]},
\end{align}
where the supremum in time is bounded by \eqref{space time geometry sup} at $a=2-\kappa$, and the space-time norm by the norm of $\mb v$ in $\mathcal E_T$.

Now let both time derivatives fall on a geometry factor, as in the terms $(\dr^{j_1}\grad^2\bs\nu[\partial_t^2\psi])(\dr^{j_2}v_y)$. This geometry factor is distinguished, at its $M=0$ weight, and loses $\kappa$. The solution factor carries no time derivative and is in the supremum in time at its $M=1$ weight, gaining $1$:
\begin{align}\label{corner target geometry}
\nn(\dr^{j_1}\grad^2\bs\nu[\partial_t^2\psi])(\dr^{j_2}v_y)\nn_{\Omega\cap\{r<1\},-2\kappa}
\leq\sup_{t\in[0,T]}\|r^{-(1-\kappa)}\dr^{j_2}v_y(t)\|_{L^\infty(\Omega)}\nn\grad^2\bs\nu[\partial_t^2\psi]\nn_{\Omega,-1-\kappa,j_1},
\end{align}
where the supremum in time is bounded by \eqref{solution sup} at the weight $2-\kappa$ of $\mb v^{\smash[t]{[1]}}$, recalling $v_y=v_y^{\smash[t]{[1]}}$, and the space-time norm by \eqref{space time geometry L2} for $\bs\nu[\partial_t^2\psi]$ at $a=1-\kappa$.

Now let one time derivative fall on the solution factor and one on a geometry factor, as in the terms $(\dr^{j_1}\grad^2\bs\nu[\partial_t\psi])(\dr^{j_2}\partial_tv_y)$. By the rule of Step 2 the geometry factor is distinguished, at its $M=1$ weight, and gains $1-\kappa$. The solution factor is in the supremum in time at its $M=0$ weight and gains nothing:
\begin{align}\label{corner target one each}
\nn(\dr^{j_1}\grad^2\bs\nu[\partial_t\psi])(\dr^{j_2}\partial_tv_y)\nn_{\Omega\cap\{r<1\},-2\kappa}
\leq\sup_{t\in[0,T]}\|r^\kappa\dr^{j_2}\partial_tv_y(t)\|_{L^\infty(\Omega)}\nn\grad^2\bs\nu[\partial_t\psi]\nn_{\Omega,-\kappa,j_1},
\end{align}
where the supremum in time is bounded by \eqref{solution sup} at the weight $1-\kappa$ of $\partial_t\mb v$, and the space-time norm by \eqref{space time geometry L2} for $\bs\nu[\partial_t\psi]$ at $a=2-\kappa$.

Finally let the two time derivatives fall on two geometry factors, as in the terms of $\partial_t^2((\grad\bs\nu)(\grad\bs\nu)v_y)$ with one time derivative on each factor $\grad\bs\nu$. The one carrying more spatial derivatives is distinguished, at $a=1$, and gains nothing; the other is in the supremum in time at its $M=0$ weight and loses $\kappa$. The solution factor carries no time derivative and is in the supremum in time at its $M=1$ weight, gaining $1$: for $j_3\leq j_1$,
\begin{align}
&\nn(\dr^{j_1}\grad\bs\nu[\partial_t\psi])(\dr^{j_3}\grad\bs\nu[\partial_t\psi])(\dr^{j_2}v_y)\nn_{\Omega\cap\{r<1\},-2\kappa}\nonumber\\
&\qquad\leq\sup_{t\in[0,T]}\|r^{1+\kappa}\dr^{j_3}\grad\bs\nu[\partial_t\psi(t)]\|_{L^\infty(\Omega)}\sup_{t\in[0,T]}\|r^{-(1-\kappa)}\dr^{j_2}v_y(t)\|_{L^\infty(\Omega)}\nn\grad\bs\nu[\partial_t\psi]\nn_{\Omega,0,j_1},\label{corner target two geometry}
\end{align}
where the suprema in time are bounded by \eqref{space time geometry sup} for $\bs\nu[\partial_t\psi]$ at $a=1-\kappa$ and by \eqref{solution sup} at the weight $2-\kappa$ of $\mb v^{[1]}$, and the space-time norm by \eqref{space time geometry L2} for $\bs\nu[\partial_t\psi]$ at $a=1$.

In a pure geometry term $r^{-1}AB$ the factor $B$ takes the place of the solution factor, and the four cases are treated in the same way. The only difference arises when one time derivative falls on $B$ and one on a factor of $A$, as in $r^{-1}(\nu_1[\partial_t\psi])\dr^{l+1}\nu_2[\partial_t\psi]$. Then $B$, carrying the most spatial derivatives, is distinguished by the rule of Step 2, at its $M=1$ weight, and gains $1$, while the factor of $A$ is in the supremum in time at its $M=0$ weight and loses $\kappa$:
\begin{align}\label{corner target pure geometry}
\nn r^{-1}AB\nn_{\Omega\cap\{r<1\},1-2\kappa}\leq\sup_{t\in[0,T]}\|r^\kappa A(t)\|_{L^\infty(\Omega)}\nn B\nn_{\Omega,2-\kappa},
\end{align}
where in the supremum in time the factor of $A$ carrying the time derivative is bounded by \eqref{space time geometry sup} for $\bs\nu[\partial_t\psi]$ at $a=1-\kappa$ and the other factors of $A$ by \eqref{space time geometry sup} at $a=1$, and $B$ by \eqref{space time geometry L2} for $\bs\nu[\partial_t\psi]$ at $a=2-\kappa$.

In every case the gains add up to at least $1-\kappa$, and to exactly $1-\kappa$ in \eqref{corner target solution}--\eqref{corner target pure geometry}. By \eqref{corner target reduction} the corner targets are therefore $\lesssim_{\theta_b}\mathcal K_T\mathcal E_T$. Since the gain $1-\kappa$ is attained, these placements reach the corner targets only if $\varpi+\kappa\leq 1-\kappa$. With $\varpi=\kappa$ this would require $\kappa\leq 1/3$, which is why $\varpi$ is chosen by \eqref{varpi choice} rather than equal to $\kappa$, the value for which, in the linear theory, the weight $-\varpi$ is itself an energy weight; the linear theory permits any $\varpi\in(0,\kappa]$. This proves \eqref{ST N estimate}, and the bound on $H$ in the statement, $\nn H\nn_{\Omega,1-\kappa,l+2,1}$ being one of the targets of \eqref{space-time targets}.

\emph{Step 6: the difference.} The summands of \eqref{schematic bulk}--\eqref{schematic boundary} are multilinear in their factors, so the difference of the two nonlinearities, and of their time derivatives, is a sum of terms in each of which exactly one factor $\mathfrak d$ is a difference: a component of $\bs\nu[\psi-\psi']$, possibly with time derivatives, a solution difference, or a factor $\mathfrak G-\mathfrak G'$, which is $(\nu_1-\nu_1')\tilde{\mathfrak G}+(\nu_2-\nu_2')\tilde{\mathfrak G}$ by Lemma \ref{composition of the geometry}(iv). As in the multilinear clause of Proposition \ref{nonlinear estimates}, we place every term as in Steps 2 to 5. If $\mathfrak d$ is distinguished, it contributes $d_T$, by the space-time bounds of Lemma \ref{space-time geometry estimates}, which do not involve the initial values, applied to the difference triple, $\mathcal E_T$ being a norm, and the other factors contribute $\mathcal K_T$ or $\mathcal K_T'$. If $\mathfrak d$ is placed in the supremum in time, its value at $t=0$ vanishes, since the two triples satisfy the initial conditions of $B_T$. So \eqref{sup in time 1} gives
\begin{align*}
\sup_{t\in[0,T]}\|\mathfrak d(t)\|\leq\Theta_T\nn\partial_t\mathfrak d\nn\lesssim_{\theta_b}\Theta_Td_T,
\end{align*}
where $\|\ph\|$ is the spatial norm through which the supremum in time is bounded in Steps 3 to 5 and $\nn\ph\nn$ the corresponding space-time norm. The distinguished factor then contributes $\mathcal E_T$ or $\mathcal E_T'$. Both cases are bounded by the RHS of \eqref{ST N difference}. Unlike in \eqref{N difference}, the pure geometry terms of \eqref{schematic boundary} do not cancel here, the two profiles being different; their difference is placed in the same way.
\end{proof}

\subsubsection{The positive-time theorem}

Each step of the iteration solves \eqref{system: time side undifferentiated} with the profile datum $\phi_0$ by Theorem \ref{linear solution full regularity}, always with the initial jet $(\phi_0,\dot\phi_0,\ddot\phi_0)$. The next lemma shows that this theorem applies whenever the forcing has the initial values \eqref{initial forcing values}, and that its solution then satisfies the initial conditions of $B_T$. The initial values of a forcing pair with $\mathcal N_T(\mb F,\mb f)<\infty$ are well defined by Lemma \ref{supremum in time}.

\begin{lemma}[The initial jet along the iteration]\label{iteration jet}
Let $(\mb F,\mb f)$ be a forcing pair with $\mathcal N_T(\mb F,\mb f)<\infty$ whose initial values are those of \eqref{initial forcing values},
\begin{align}\label{forcing jet}
(\mb F,\mb f)(0)=(\mb N_0,\mb N_{\partial,0}),\qquad\qquad
(\partial_t\mb F,\partial_t\mb f)(0)=(\dot{\mb N}_0,\dot{\mb N}_{\partial,0}).
\end{align}
Then the following hold.
\begin{enumerate}
\item Theorem \ref{linear solution full regularity} applies at $\tilde\kappa=0$ to the data $(\phi_0,\mb F,\mb f)$ with the initial jet $(\phi_0,\dot\phi_0,\ddot\phi_0)$.
\item The solution $(\phi,\mb u,p)$ of \eqref{system: time side undifferentiated} given by Theorem \ref{linear solution full regularity} satisfies \eqref{linear estimate renamed} and the initial conditions of $B_T$,
\begin{align*}
(\phi,\partial_t\phi,\mb u,\partial_t\mb u,p,\partial_tp)(0)=(\phi_0,\dot\phi_0,\mb u_0,\dot{\mb u}_0,p_0,\dot p_0).
\end{align*}
\end{enumerate}
\end{lemma}
\begin{proof}
(i) By \eqref{forcing jet}, the data at $t=0$ to which Theorem \ref{linear solution full regularity} refers are
\begin{align*}
(\phi_0,\mb F(0),\mb f(0))=(\phi_0,\mb N_0,\mb N_{\partial,0}),\qquad\qquad
(\dot\phi_0,\partial_t\mb F(0),\partial_t\mb f(0))=(\dot\phi_0,\dot{\mb N}_0,\dot{\mb N}_{\partial,0}).
\end{align*}
By the identification made after \eqref{initial forcing values}, the pairs $(\mb u_0,p_0)$ and $(\dot{\mb u}_0,\dot p_0)$ solve the stationary system \eqref{system: stationary linear} with these data, the first satisfying \eqref{C1 in the linear theory}; $\dot\phi_0$ and $\ddot\phi_0$ are the profiles derived from them; and $\ddot\phi_0$ satisfies \eqref{improved theorem C2}. These are the hypotheses of Theorem \ref{linear solution full regularity} on the initial jet, once the two pairs are known to be the stationary solutions of Proposition \ref{stationary system}. By the uniqueness statement of that proposition it suffices that they lie in the class of \eqref{data estimate}. For $(\mb u_0,p_0)$ the norms of \eqref{data estimate} are finite by the rungs $M=0$ and $M=1$ of $\mathcal E(\mb u_0,p_0)\lesssim_{\theta_b}\mathcal D$, the unweighted term by \eqref{unweighted by weighted} with $j=l+2$.
The dotted pair is treated in the same way, with $j=l+1$ and $\dot{\mathcal E}(\dot{\mb u}_0,\dot p_0)\lesssim_{\theta_b}\mathcal D$. The hypothesis of Theorem \ref{linear solution full regularity} on the forcing at $\tilde\kappa=0$ is $\mathcal N_T(\mb F,\mb f)<\infty$.

(ii) Theorem \ref{linear solution full regularity} gives \eqref{Improved solution estimate}, that is, \eqref{linear estimate renamed}, and the initial values $(\mb u_0,p_0,\dot{\mb u}_0,\dot p_0)$ of $(\mb u,p,\partial_t\mb u,\partial_tp)$. For the profile, $\phi(\ph,0)=\phi_0$ is part of \eqref{system: time side undifferentiated}, and the profile evolution equation gives $\partial_t\phi(0)=(\sin\theta_b)^{-2}\partial_ru_\theta(0)|_{\partial_1\Omega}=(\sin\theta_b)^{-2}\partial_ru_{0,\theta}|_{\partial_1\Omega}=\dot\phi_0$, the last equality being the definition of $\dot\phi_0$ in Theorem \ref{nonlinear initial data}.
\end{proof}

\begin{theorem}[Existence of the positive-time solution]\label{nonlinear positive-time existence}
Let $\kappa\in(1/4,1/2)$, let $l\geq 1$ be an integer, let $\varsigma>0$, assume the standing angle condition \eqref{standing angle condition}, and let $\theta_b$ be small enough for the linear theory of Sections \ref{linear existence} and \ref{corner expansions} to apply. There exist $\delta,C>0$, depending only on $\theta_b,\kappa,l$ and $\varsigma$, and of the form
\begin{align}\label{delta tau form}
\delta\geq c\brac{1+\varsigma^{-1}}^{-4},\qquad\qquad C\leq c^{-1}\brac{1+\varsigma^{-1}},
\end{align}
for a constant $c=c(\theta_b,\kappa,l)\in(0,1]$, such that the following holds. Suppose the initial profile $\phi_0$, with $\phi_0(0)=0$ and $\lim_{r\to\infty}\phi_0=0$, satisfies the two compatibility conditions \eqref{C1} and \eqref{C2} of Theorem \ref{nonlinear initial data}, and suppose $T>0$ is such that
\begin{align}\label{positive time smallness}
\brac{1+\Theta_T}\mathcal D\leq\delta,\qquad\qquad\Theta_T=e^{\varsigma T}T^{1\over 2}.
\end{align}
Then there exists a triple $(\phi,\mb u,p)$ on $[0,T]$ satisfying a.e.\ System \ref{reference wedge system} in full,
\begin{subequations}\label{system: positive time}
\begin{alignat}{3}
\lpc\mb u-\grad p&=\mb N[\phi,\mb u,p]\qquad&&\text{on}\quad\Omega\\
\grad\cdot\mb u&=0\qquad&&\text{on}\quad\Omega\\
u_\theta&=0\qquad&&\text{on}\quad\partial_0\Omega\\
u_r-r^{-1}\partial_\theta u_r&=0\qquad&&\text{on}\quad\partial_0\Omega\\
\mathbb T[\mb u,p]\mb e_\theta-(\sin\theta_b)^2(\partial_r\phi)\mb e_\theta&=\mb N_\partial[\phi,\mb u,p]\qquad&&\text{on}\quad\partial_1\Omega\\
\partial_t\phi-(\sin\theta_b)^{-2}\partial_ru_\theta&=0\qquad&&\text{on}\quad\partial_1\Omega\\
\lim_{r\to\infty}\phi(r,t)&=0\qquad&&\text{for all}\quad t\\
\phi(\ph,0)&=\phi_0
\intertext{and satisfying weakly}
\phi(0,t)&=0\qquad&&\text{for all}\quad t,
\end{alignat}
\end{subequations}
with the regularity and the initial values
\begin{align}\label{positive time main estimate}
\mathcal E_T(\phi,\mb u,p)\leq C\mathcal D,\qquad\qquad
(\partial_t\phi,\partial_t^2\phi,\mb u,\partial_t\mb u,p,\partial_tp)(0)=(\dot\phi_0,\ddot\phi_0,\mb u_0,\dot{\mb u}_0,p_0,\dot p_0),
\end{align}
and it is the only triple in the iteration space $B_T$ of \eqref{iteration ball} with radius $\rho=C\mathcal D$ that satisfies \eqref{system: positive time}.
\end{theorem}
\begin{proof}
Step 1 constructs an element of the ball $B_T$ of \eqref{iteration ball} and fixes the constants. Step 2 shows that the iteration map sends $B_T$ into itself, Step 3 that it is a contraction on the complete metric space $(B_T,d_T)$, and Step 4 that its fixed point solves \eqref{system: positive time} and is unique.

\emph{Step 1: the ball is not empty.} The pair $(\mb N_0+t\dot{\mb N}_0,\mb N_{\partial,0}+t\dot{\mb N}_{\partial,0})$, which has the initial values \eqref{forcing jet}, need not have finite $\mathcal N_T$, since the norms in $\mathcal N_T$ without time derivative ask one radial derivative more of $(\dot{\mb N}_0,\dot{\mb N}_{\partial,0})$ than \eqref{DN estimate} provides; we supply it by smoothing radially at a scale vanishing with $t$. Let $\mathfrak s(t):=te^{-t}$, so that $\mathfrak s(0)=0$, $\dot{\mathfrak s}(0)=1$, and $\mathfrak s$, $\dot{\mathfrak s}$ and $\ddot{\mathfrak s}$ are bounded on $[0,\infty)$ with $\mathfrak s\leq 1$, and let $S_{\mathfrak s(t)}$ be the radial smoothing of Lemma \ref{radial smoothing} at the scale $\mathfrak s(t)$. The smoothed flat forcing is
\begin{align}\label{flat forcing}
(\mb F^\flat,\mb f^\flat):=\brac{\mb N_0+\mathfrak s S_{\mathfrak s}\dot{\mb N}_0,\ \mb N_{\partial,0}+\mathfrak s S_{\mathfrak s}\dot{\mb N}_{\partial,0}},
\qquad\quad
\mathcal N_T(\mb F^\flat,\mb f^\flat)\lesssim_{\theta_b}\brac{1+\varsigma^{-1}}\mathcal D^2,
\end{align}
the estimate holding uniformly in $T$.

The pair satisfies \eqref{forcing jet}. Its part $(\mb N_0,\mb N_{\partial,0})$ does not depend on time, and $\partial_tS_{\mathfrak s}=2\dot{\mathfrak s}\mathfrak s\dr^2S_{\mathfrak s}$ for $t>0$ by Lemma \ref{radial smoothing}(iii) and the chain rule. Hence the time derivatives of the pair are those of its part $\mathfrak sS_{\mathfrak s}(\dot{\mb N}_0,\dot{\mb N}_{\partial,0})$, with
\begin{align}
\partial_t(\mathfrak sS_{\mathfrak s})&=\dot{\mathfrak s}S_{\mathfrak s}+2\dot{\mathfrak s}\mathfrak s^2\dr^2S_{\mathfrak s},\nonumber\\*
\partial_t^2(\mathfrak sS_{\mathfrak s})&=\ddot{\mathfrak s}S_{\mathfrak s}+6\dot{\mathfrak s}^2\mathfrak s\dr^2S_{\mathfrak s}+2\ddot{\mathfrak s}\mathfrak s^2\dr^2S_{\mathfrak s}+4\dot{\mathfrak s}^2\mathfrak s^3\dr^4S_{\mathfrak s}.\label{smoothing time derivatives}
\end{align}
As $t\to 0$, $S_{\mathfrak s}$ applied to a fixed function tends to that function, and $\mathfrak s^2\dr^2S_{\mathfrak s}=(\mathfrak s\dr)^2S_{\mathfrak s}$ applied to it tends to $0$, by the continuity stated in Lemma \ref{radial smoothing}(iii). So at $t=0$, where $\mathfrak s$ vanishes, $\dot{\mathfrak s}=1$ and $S_0$ is the identity, the pair takes the values $(\mb N_0,\mb N_{\partial,0})$ and its time derivative the values $(\dot{\mb N}_0,\dot{\mb N}_{\partial,0})$.

For the estimate of \eqref{flat forcing} we first check the counts. The norms in $\mathcal N_T$ ask for $\partial_t^n(\mb F^\flat,\mb f^\flat)$, $n\leq 2$, at the counts $l+2-n$ and $l+{5\over 2}-n$. The part $(\mb N_0,\mb N_{\partial,0})$ occurs only for $n=0$, at the counts $l+2$ and $l+{5\over 2}$ supplied by \eqref{N estimate}. The rest is $\mathfrak sS_{\mathfrak s}$ and its time derivatives \eqref{smoothing time derivatives} applied to $(\dot{\mb N}_0,\dot{\mb N}_{\partial,0})$, which \eqref{DN estimate} supplies only at the counts $l+1$ and $l+{3\over 2}$. Each operator in $\partial_t^n(\mathfrak sS_{\mathfrak s})$ therefore asks of $S_{\mathfrak s}\dot{\mb N}_0$ the count $l+2-n$ plus the number of its radial derivatives $\dr$. We call the amount $k$ by which the count asked exceeds the count $l+1$ available the excess of the operator; for $\dot{\mb N}_{\partial,0}$ both counts are raised by $1/2$ and the excess is the same. The power of $\mathfrak s$ in front of every operator is at least its excess:
\begin{align*}
\begin{array}{c|c|c|c}
n&\text{operators}&\text{excess }k&\text{power of }\mathfrak s\\\hline
0&\mathfrak sS_{\mathfrak s}&1&1\\
1&S_{\mathfrak s},\ \mathfrak s^2\dr^2S_{\mathfrak s}&0,\ 2&0,\ 2\\
2&S_{\mathfrak s},\ \mathfrak s\dr^2S_{\mathfrak s},\ \mathfrak s^2\dr^2S_{\mathfrak s},\ \mathfrak s^3\dr^4S_{\mathfrak s}&-1,\ 1,\ 1,\ 3&0,\ 1,\ 2,\ 3
\end{array}
\end{align*}
When the excess $k$ is positive, \eqref{smoothing gain} for $S_\sigma$ itself, at $\sigma=\mathfrak s$, costs a factor $\mathfrak s^{-k}$, which the power of $\mathfrak s$ in front pays, since $\mathfrak s\leq 1$. For the operator $\mathfrak s^3\dr^4S_{\mathfrak s}$, for instance, the same bound with $k=3$ gives, at each weight $a$ of $\mathcal N_T$,
\begin{align*}
\|\mathfrak s^3\dr^4S_{\mathfrak s}\dot{\mb N}_0\|_{\Omega,a,l}\leq\mathfrak s^3\|S_{\mathfrak s}\dot{\mb N}_0\|_{\Omega,a,l+4}\lesssim\|\dot{\mb N}_0\|_{\Omega,a,l+1}.
\end{align*}
A non-positive excess costs nothing, by the same bound with $k=0$.

So the norms in $\mathcal N_T$ of the pair factor into $L^2_\varsigma[0,\infty)$ norms of bounded functions of time against stationary norms of the initial forcing values \eqref{initial forcing values}. The functions of time are the constant $1$, in front of the time-independent part $(\mb N_0,\mb N_{\partial,0})$, whose $L^2_\varsigma[0,\infty)$ norm is $(2\varsigma)^{-1/2}\leq 1+\varsigma^{-1}$, and products of $\mathfrak s$, $\dot{\mathfrak s}$, $\ddot{\mathfrak s}$ and of powers of $\mathfrak s$, which remain bounded once the smoothing costs have been paid as above, so that their $L^2_\varsigma[0,\infty)$ norms are $\lesssim(2\varsigma)^{-1/2}\leq 1+\varsigma^{-1}$. The stationary norms are $\lesssim_{\theta_b}\mathcal D^2$ by \eqref{N estimate}, \eqref{DN estimate} and \eqref{main estimate}, which proves \eqref{flat forcing}.

We now fix the constants. Let $C_{{\rm lin},\varsigma}\geq 1$ be the constant of \eqref{linear estimate renamed}, and let $C_{{\rm nl},\varsigma}\geq 1$ dominate the constants of Lemma \ref{space-time geometry estimates}, of Proposition \ref{space-time nonlinear estimates} and of \eqref{flat forcing}, so that $C_{{\rm lin},\varsigma}+C_{{\rm nl},\varsigma}\lesssim_{\theta_b}1+\varsigma^{-1}$. We set $\rho:=C_{\rm ball}\mathcal D$ and $C:=C_{\rm ball}$, with $C_{\rm ball}:=2\max\{C_{{\rm lin},\varsigma},C_{{\rm nl},\varsigma}\}$, and we take $\delta$ at most the constant of Theorem \ref{nonlinear initial data} and small enough that
\begin{align}\label{positive time delta conditions}
(1+C_{\rm ball})\delta\leq\delta_1\qquad\text{ and }\qquad 2C_{{\rm lin},\varsigma}C_{{\rm nl},\varsigma}(1+C_{\rm ball})^2\delta\leq 1,
\end{align}
with $\delta_1$ as in Lemma \ref{space-time geometry estimates}(iii). Since $\delta_1$ and the constant of Theorem \ref{nonlinear initial data} do not depend on $\varsigma$, the choices $C=C_{\rm ball}\lesssim_{\theta_b}1+\varsigma^{-1}$ and $\delta\gtrsim_{\theta_b}\min\{\delta_1C_{\rm ball}^{-1},(C_{{\rm lin},\varsigma}C_{{\rm nl},\varsigma}C_{\rm ball}^2)^{-1}\}\gtrsim_{\theta_b}(1+\varsigma^{-1})^{-4}$ are of the form \eqref{delta tau form}. For every triple of the ball $B_T$ of \eqref{iteration ball} with radius $\rho$ we then have, by \eqref{positive time smallness},
\begin{align*}
\mathcal K_T\leq\mathcal D+\Theta_T\rho\leq(1+C_{\rm ball})(1+\Theta_T)\mathcal D\leq(1+C_{\rm ball})\delta\leq\delta_1,
\end{align*}
so that Lemma \ref{space-time geometry estimates}(iii) and Proposition \ref{space-time nonlinear estimates} are available on all of $B_T$.

Let $w^{(0)}$ be the solution of \eqref{system: time side undifferentiated} with data $(\phi_0,\mb F^\flat,\mb f^\flat)$ given by Lemma \ref{iteration jet}; by \eqref{linear estimate renamed} and \eqref{flat forcing},
\begin{align*}
\mathcal E_T(w^{(0)})\leq C_{{\rm lin},\varsigma}\brac{C_{{\rm nl},\varsigma}\mathcal D^2+\mathcal D}\leq{3\over 2}C_{{\rm lin},\varsigma}\mathcal D\leq\rho
\end{align*}
by \eqref{positive time delta conditions}, so that $w^{(0)}\in B_T$.

\emph{Step 2: the map.} For $(\psi,\mb v,q)\in B_T$ let $\mathcal T(\psi,\mb v,q):=(\phi,\mb u,p)$ be the solution of \eqref{system: time side undifferentiated} with data
\begin{align*}
(\phi_0,\mb F,\mb f)=\brac{\phi_0,\ \mb N[\psi,\mb v,q],\ \mb N_\partial[\psi,\mb v,q]}
\end{align*}
given by Lemma \ref{iteration jet}. This is well defined: the forcing has $\mathcal N_T<\infty$ by \eqref{ST N estimate}, and its initial values are \eqref{forcing jet}, since the initial conditions of $B_T$ give $(\mb N[\psi,\mb v,q])(0)=\mb N[\phi_0,\mb u_0,p_0]=\mb N_0$ and, by the chain rule \eqref{time chain rule}, $(\partial_t\mb N[\psi,\mb v,q])(0)=D\mb N[\phi_0,\mb u_0,p_0](\dot\phi_0,\dot{\mb u}_0,\dot p_0)=\dot{\mb N}_0$, and likewise for $\mb N_\partial$. By the same lemma the output satisfies the initial conditions, and it satisfies the conditions $\grad\cdot\mb v=0$ and $v_\theta|_{\partial_0\Omega}=0$ of $B_T$ as part of the linear system. By \eqref{linear estimate renamed}, \eqref{ST N estimate} and \eqref{positive time delta conditions},
\begin{align*}
\mathcal E_T\brac{\mathcal T(\psi,\mb v,q)}&\leq C_{{\rm lin},\varsigma}\brac{\mathcal D+C_{{\rm nl},\varsigma}\brac{\mathcal D+\Theta_T\rho}\rho}
\leq C_{{\rm lin},\varsigma}\mathcal D\brac{1+C_{{\rm nl},\varsigma}C_{\rm ball}(1+C_{\rm ball})(1+\Theta_T)\mathcal D}\\
&\leq{3\over 2}C_{{\rm lin},\varsigma}\mathcal D\leq\rho
\end{align*}
by \eqref{positive time delta conditions}, so that $\mathcal T(B_T)\subseteq B_T$.

\emph{Step 3: contraction and completeness.} For two elements of $B_T$, the difference of the two outputs solves \eqref{system: time side undifferentiated} with profile datum $0$ and forcing the difference of the forcings, so by \eqref{linear estimate renamed} with $\mathcal D$ replaced by $0$, \eqref{ST N difference} and \eqref{positive time delta conditions},
\begin{align*}
d_T\brac{\mathcal T(\psi,\mb v,q),\mathcal T(\psi',\mb v',q')}
\leq C_{{\rm lin},\varsigma}C_{{\rm nl},\varsigma}\brac{\mathcal D+2\Theta_T\rho}d_T
\leq C_{{\rm lin},\varsigma}C_{{\rm nl},\varsigma}(1+C_{\rm ball})^2(1+\Theta_T)\mathcal Dd_T\leq{1\over 2}d_T.
\end{align*}
The set $B_T$ is complete for $d_T$. The space-time weighted $L^2$ spaces underlying $\mathcal E_T$ are complete, so the pieces of a $d_T$ Cauchy sequence converge. By \eqref{spacetime coefficients} applied to differences, the corner values $\mb v(\mb 0,\ph)$ converge in $L^2$ of time as well, so the limits are the pieces of a single triple. The conditions $\grad\cdot\mb v=0$ and $v_\theta|_{\partial_0\Omega}=0$ pass to the limit. So do the initial conditions, because, by \eqref{sup in time 3} applied to differences, each of $\psi,\partial_t\psi,\mb v,\partial_t\mb v,q,\partial_tq$ converges in $C([0,T])$ with values in the corresponding weighted space. Finally, the limit satisfies $\mathcal E_T\leq\rho$ by Fatou's lemma. The Banach fixed point theorem therefore yields a unique fixed point $(\phi,\mb u,p)\in B_T$.

\emph{Step 4: the solved system and uniqueness.} Being a fixed point, $(\phi,\mb u,p)$ is the solution of \eqref{system: time side undifferentiated} with forcing $(\mb N,\mb N_\partial)[\phi,\mb u,p]$. This is the statement that $(\phi,\mb u,p)$ satisfies \eqref{system: positive time} a.e., with the contact angle condition holding weakly. Moreover, \eqref{positive time main estimate} holds with $C=C_{\rm ball}$: the estimate is membership in $B_T$, the initial values of $(\partial_t\phi,\mb u,\partial_t\mb u,p,\partial_tp)$ are given by Lemma \ref{iteration jet}(ii), and differentiating the profile evolution equation in time gives $\partial_t^2\phi(0)=(\sin\theta_b)^{-2}\partial_r\partial_tu_\theta(0)|_{\partial_1\Omega}=(\sin\theta_b)^{-2}\partial_r\dot u_{0,\theta}=\ddot\phi_0$. The trace at $t=0$ is meaningful since $\partial_t\mb u$ carries a further time derivative in $\mathcal E_T$. For the uniqueness, let $(\phi',\mb u',p')$ be a triple in $B_T$ that satisfies \eqref{system: positive time}. Then both $(\phi',\mb u',p')$ and $\mathcal T(\phi',\mb u',p')$ solve \eqref{system: time side undifferentiated} with the same data $(\phi_0,\mb N[\phi',\mb u',p'],\mb N_\partial[\phi',\mb u',p'])$ and make the norms \eqref{uniqueness class} finite, these being controlled along $B_T$ by the rungs $M=0$ and $M=1$ of $\mathcal E_T$; so they coincide, by Lemma \ref{linear evolution uniqueness}. Thus $(\phi',\mb u',p')$ is a fixed point of $\mathcal T$ in $B_T$, and by Step 3 it equals $(\phi,\mb u,p)$.
\end{proof}

Reversing the flattening turns the solution of Theorem \ref{nonlinear positive-time existence} into a solution of the problem in the moving domain.

\begin{corollary}[Return to the moving domain]\label{moving domain}
Let $(\phi,\mb u,p)$ be the solution of Theorem \ref{nonlinear positive-time existence} on $[0,T]$. Then for every $t\in[0,T]$ the map $\bs\Psi_t$ of Lemma \ref{regularised flattening} is a diffeomorphism $\Omega\to\Omega_t$, and reversing the transformation of Lemma \ref{Piola transform} turns $(\phi,\mb u,p)$ into a solution on $[0,T]$ of the moving domain problem of System \ref{dimensionless system} in the sense of \eqref{system: positive time}, the contact angle condition holding weakly. Its free boundary is $X=b^{-1}y+\chi$ with $\chi'(y,t)=-\phi(r,t)$, $y=r\sin\theta_b$, and $\partial_t\chi=-(\sin\theta_b)^{-1}u_\theta|_{\partial_1\Omega}$, the constant of integration being fixed at $t=0$ by the initial position of the contact point.
\end{corollary}
\begin{proof}
By \eqref{positive time main estimate} with $C=C_{\rm ball}$, \eqref{positive time smallness} and \eqref{positive time delta conditions}, the solution has $\mathcal K_T\leq\mathcal D+\Theta_TC_{\rm ball}\mathcal D\leq(1+C_{\rm ball})\delta\leq\delta_1$, with $\mathcal K_T$ as in \eqref{K_T definition}. So by Lemma \ref{space-time geometry estimates}(iii) the map $\bs\Psi_t$ is a diffeomorphism $\Omega\to\Omega_t$ for every $t\in[0,T]$, and reversing the transformation of Lemma \ref{Piola transform} turns the equations of \eqref{system: positive time} other than the profile evolution equation into those of System \ref{dimensionless system}, the conditions on $\phi$ at $r=0$, as $r\to\infty$ and at $t=0$ becoming those on $X$ since $\chi'=-\phi$. The two conditions on $\chi$ in the statement are consistent: by \eqref{r and y on the free boundary} and the profile evolution equation of \eqref{system: positive time}, the derivative in $y$ of the second is $\partial_t\chi'=-(\sin\theta_b)^{-2}\partial_ru_\theta=-\partial_t\phi$. The equalities of \eqref{transformed u_theta} other than the last do not use the profile evolution equation; with $\partial_t\chi=-(\sin\theta_b)^{-1}u_\theta|_{\partial_1\Omega}$ and $\sqrt{1+b^2}\sin\theta_b=b$ they give $\partial_tX=\partial_t\chi=u_x-u_yX'$ for the velocity on $\Omega_t$, which is the profile evolution equation of System \ref{dimensionless system}.
\end{proof}

The smallness condition \eqref{positive time smallness} contains a lower bound on the lifespan. For a fixed weight $\varsigma$ it is only logarithmic in $\mathcal D^{-1}$, because $\Theta_T$ grows exponentially; but the weight is ours to choose, and the form \eqref{delta tau form} of the constants lets us choose it in terms of $T$.

\begin{corollary}[Polynomial lifespan]\label{polynomial lifespan}
Let $\theta_b$, $\kappa$, $l$ and $\phi_0$ be as in Theorem \ref{nonlinear positive-time existence}. There is a constant $c_*=c_*(\theta_b,\kappa,l)\in(0,1]$ such that, if $0<\mathcal D\leq c_*$, the conclusions of Theorem \ref{nonlinear positive-time existence} hold on $[0,T_*]$ for the weight $\varsigma=1/T_*$, where
\begin{align}\label{polynomial lifespan formula}
T_*:=\brac{c_*\over\mathcal D}^{2\over 9}\geq 1.
\end{align}
That is, there is a triple $(\phi,\mb u,p)$ solving \eqref{system: positive time} on $[0,T_*]$ with the initial values \eqref{positive time main estimate}, unique in the iteration space $B_{T_*}$ of \eqref{iteration ball} for that weight and the radius $\rho=C\mathcal D$ of the theorem, and with
\begin{align}\label{polynomial lifespan estimate}
\mathcal E_{T_*}(\phi,\mb u,p)\leq c_*^{-1}\brac{1+T_*}\mathcal D,
\end{align}
the norms of $\mathcal E_{T_*}$ being formed with $\varsigma=1/T_*$. The same bound holds, up to the factor $e$, for the unweighted norms obtained by setting $\varsigma=0$ in $\mathcal E_{T_*}$.
\end{corollary}
\begin{proof}
Apply Theorem \ref{nonlinear positive-time existence} with $\varsigma:=1/T$ for a $T\geq 1$ to be chosen. Then $\Theta_T=eT^{1/2}$, and by \eqref{delta tau form} the constant of the theorem satisfies $\delta\geq c(1+T)^{-4}\geq c2^{-4}T^{-4}$, so \eqref{positive time smallness} holds as soon as $(1+e)T^{1/2}\mathcal D\leq c2^{-4}T^{-4}$, that is, as soon as $T^{9/2}\leq c'/\mathcal D$ with $c':=c/(16(1+e))$. With $c_*:=c'$ this is $T\leq T_*$, and $T_*\geq 1$ is the condition $\mathcal D\leq c_*$; the estimate \eqref{polynomial lifespan estimate} is \eqref{positive time main estimate} with $C\leq c^{-1}(1+\varsigma^{-1})=c^{-1}(1+T_*)\leq c_*^{-1}(1+T_*)$. The bound for the unweighted norms follows since $e^{-1}\leq e^{-t/T_*}\leq 1$ on $[0,T_*]$.
\end{proof}

For $\mathcal D=0$ the zero triple solves \eqref{system: positive time} on every interval and is the solution of Theorem \ref{nonlinear positive-time existence} for every $T$, so Corollary \ref{polynomial lifespan} holds trivially with $T_*:=\infty$.

\begin{remark}[On the lifespan]\label{lifespan remark}
Two comments on Corollary \ref{polynomial lifespan}.
\begin{enumerate}
\item The exponent ${2\over 9}$ is what the form \eqref{delta tau form} of the constants gives, and is not claimed to be optimal.
\item The exponential weight is a device of the $L^2$-in-time framework and not a feature of the problem. The linear evolution is dissipative: for zero data the energy identity \eqref{truncated identity} of Lemma \ref{linear evolution uniqueness} reads
\begin{align*}
{1\over 2}(\sin\theta_b)^4\partial_t\|\phi\|_{L^2(\partial_1\Omega)}^2=-2\|\S(\mb u)\|_{L^2(\Omega)}^2-\|u_r\|_{L^2(\partial_0\Omega)}^2\leq 0,
\end{align*}
so that the profile does not grow in $L^2$ at all, and the weight $e^{-\varsigma t}$ is needed only because the framework asks for square integrability on $[0,\infty)$ of a solution that need not decay. Posed on the unbounded wedge, the linear problem has no spectral gap, and the expected large-time behaviour of the profile is algebraic decay rather than exponential; a global small-data theorem would need decay estimates for the linear semigroup in the weighted spaces, which the present framework does not provide.
\end{enumerate}
\end{remark}

The solution also depends Lipschitz continuously on the initial profile, so that the problem is locally well posed in the sense of Hadamard on the set of compatible data.

\begin{proposition}[Lipschitz dependence at positive times]\label{Lipschitz positive}
Let $\phi_0$ and $\phi_0'$ both satisfy the hypotheses of Theorem \ref{nonlinear positive-time existence}, with the same $T$, and let $(\phi,\mb u,p)$ and $(\phi',\mb u',p')$ be the two solutions. After shrinking the constant $\delta$ of that theorem by a factor depending only on $\theta_b,\kappa$ and $l$, we have
\begin{align}\label{Lipschitz positive estimate}
\mathcal E_T(\phi-\phi',\mb u-\mb u',p-p')\leq C'\mathcal D(\phi_0-\phi_0'),
\end{align}
with $C'=C'(\theta_b,\kappa,l,\varsigma)\lesssim_{\theta_b}1+\varsigma^{-1}$, where $\mathcal D(\phi_0-\phi_0')$ is the data norm of the difference profile as in Proposition \ref{Lipschitz initial}. In particular the corner values $\mb u(\mb 0,t)$ and the profiles are Lipschitz functions of the datum uniformly in $t\in[0,T]$.
\end{proposition}
\begin{proof}
Write $w=(\phi,\mb u,p)$, $w'=(\phi',\mb u',p')$ and $\bar w:=w-w'$. By the two systems \eqref{system: positive time}, $\bar w$ solves the linear system \eqref{system: time side undifferentiated} with profile datum $\phi_0-\phi_0'$, forcing $(\mb N,\mb N_\partial)[w]-(\mb N,\mb N_\partial)[w']$, and the initial values equal to the difference of the two tuples of Theorem \ref{nonlinear initial data}. This difference is the tuple of Theorem \ref{initial jet} for the difference of the forcing values \eqref{initial forcing values}, by the linearity of Proposition \ref{stationary system} and of the derived profiles. Its entry $\ddot\phi_0-\ddot\phi_0'$ satisfies \eqref{improved theorem C2}, since both $\ddot\phi_0$ and $\ddot\phi_0'$ satisfy \eqref{C2}. So Theorem \ref{linear solution full regularity} applies to $\bar w$, and \eqref{Improved solution estimate} at $\tilde\kappa=0$, with its data terms now bounded by Proposition \ref{Lipschitz initial}, gives
\begin{align*}
\mathcal E_T(\bar w)\leq C_{{\rm lin},\varsigma}\brac{\mathcal N_T\brac{(\mb N,\mb N_\partial)[w]-(\mb N,\mb N_\partial)[w']}+\mathcal D(\phi_0-\phi_0')}.
\end{align*}

The forcing difference is estimated as in \eqref{ST N difference}, with one change. There the two triples shared their initial values, so that every difference factor placed in a supremum in time had initial value $0$ and was bounded through \eqref{sup in time 1} by $\Theta_Td_T(w,w')$. Here the difference factors have the initial values of the difference of the two tuples, bounded by Proposition \ref{Lipschitz initial}, and \eqref{sup in time 1} bounds them by $\mathcal D(\phi_0-\phi_0')+\Theta_T\mathcal E_T(\bar w)$ instead. The geometry differences $\bs\nu[\phi-\phi']$ are handled by Lemma \ref{space-time geometry estimates} applied to the triple $\bar w$: the proof of that lemma is linear in the triple, and the initial values of $\bar w$ are those of the difference profile. Every other placement of the proof of \eqref{ST N difference} is unchanged, so
\begin{align*}
\mathcal N_T\brac{(\mb N,\mb N_\partial)[w]-(\mb N,\mb N_\partial)[w']}\lesssim_{\theta_b}\brac{\mathcal K_T+\mathcal K_T'}\mathcal E_T(\bar w)+\brac{\mathcal E_T(w)+\mathcal E_T(w')}\mathcal D(\phi_0-\phi_0'),
\end{align*}
where $\mathcal K_T=\mathcal D(\phi_0)+\Theta_T\mathcal E_T(w)$ and $\mathcal K_T'$ likewise. The term with $\mathcal D(\phi_0-\phi_0')$ comes from the initial values of the difference factors in the supremum, paired with $w$ or $w'$ in $L^2$ in time. The quantities $\mathcal K_T$, $\mathcal K_T'$, $\mathcal E_T(w)$ and $\mathcal E_T(w')$ are all at most $(1+C_{\rm ball})\delta$ by \eqref{positive time smallness} and \eqref{positive time main estimate}. We shrink $\delta$ so that the factor $\lesssim_{\theta_b}C_{{\rm lin},\varsigma}(\mathcal K_T+\mathcal K_T')$ multiplying $\mathcal E_T(\bar w)$ is at most ${1\over 2}$. This is a condition of the form $C(\theta_b)C_{{\rm lin},\varsigma}(1+C_{\rm ball})\delta\leq{1\over 2}$, satisfied for $\delta\lesssim_{\theta_b}(1+\varsigma^{-1})^{-2}$, so it preserves the form \eqref{delta tau form}. We obtain $\mathcal E_T(\bar w)\lesssim_{\theta_b}C_{{\rm lin},\varsigma}\mathcal D(\phi_0-\phi_0')$, which is \eqref{Lipschitz positive estimate} with $C'\lesssim_{\theta_b}C_{{\rm lin},\varsigma}\lesssim_{\theta_b}1+\varsigma^{-1}$. The last claim follows from \eqref{spacetime coefficients} and \eqref{profile spacetime ladder} for $\bar w$, whose initial values are bounded by Proposition \ref{Lipschitz initial}.
\end{proof}

It remains to recover, for the solution just constructed, the rung $M=2$ of Theorem \ref{linear solution full regularity}, which the iteration did not use. We recover it afterwards, by replaying the iteration of the proof, rather than including it in the iteration from the start, for two reasons. First, the existence theorem does not need it: the estimates of the nonlinearities at the rungs $M\leq 1$ close among themselves, and without the rung $M=2$ the construction holds under the standing angle condition \eqref{standing angle condition}. The rung $M=2$, on the other hand, asks for the hypotheses of Theorem \ref{jet corner structure}: at the start of the iteration, through the third rung of the initial data, and along the iteration, through geometry weights beyond those that \eqref{standing angle condition} admits. Second, the rung $M=2$ carries no time derivative, so that its remainders and the coefficients $c(t)$ and $c_3(t)$ of Proposition \ref{time for radial trade} have no bound on their suprema in time. In the estimates of the nonlinearities they can therefore only be placed in $L^2$ in time, against suprema supplied by the rungs $M\leq 1$, and they enter linearly and with a small factor. So the rung $M=2$ is carried along an iteration whose convergence is already secured at $M\leq 1$, and it is then bounded by $\mathcal D$ alone, although the start of the iteration has $M=2$ norms that are finite but not small. It is here that the shift $\tilde\kappa$ and the corner decompositions at the level $M=2$ enter.

\begin{proposition}[Improved regularity of the positive-time solution]\label{positive-time regularity upgrade}
Let $(\phi,\mb u,p)$ be the solution of Theorem \ref{nonlinear positive-time existence}, assume in addition the strengthened angle condition \eqref{third rung angle condition} and the finiteness of the extended data norm $\mathcal D_+$ of \eqref{data norm +}, and fix $\tilde\kappa\in(\kappa/2,\kappa)$. After shrinking the constant $\delta$ of Theorem \ref{nonlinear positive-time existence} by a factor depending only on $\theta_b,\kappa,\tilde\kappa,l$ and $\varsigma$, in a way that preserves the form \eqref{delta tau form}, we have
\begin{align}
&(\sin\theta_b)^2\nn\phi\nn_{[0,T],\partial_1\Omega:(\varsigma,{3\over 2}+\tilde\kappa-\kappa)(1,l+{5\over 2})}+\nn\mb u^{[1]}\nn_{[0,T],\Omega:(\varsigma,2+\tilde\kappa-\kappa)(1,l+1,2)}+\nn p\nn_{[0,T],\Omega:(\varsigma,1+\tilde\kappa-\kappa)(1,l+1,1)}\nonumber\\
&+(\sin\theta_b)^2\nn\phi^{[2]}\nn_{[0,T],\partial_1\Omega:(\varsigma,{5\over 2}+\tilde\kappa-\kappa)(0,l+{7\over 2})}+\nn\mb u^{[2]}\nn_{[0,T],\Omega:(\varsigma,3+\tilde\kappa-\kappa)(0,l+2,2)}+\nn p^{[2]}\nn_{[0,T],\Omega:(\varsigma,2+\tilde\kappa-\kappa)(0,l+2,1)}\nonumber\\
&\lesssim_{\theta_b}\mathcal D,\label{upgrade estimate}
\end{align}
with the decompositions of Proposition \ref{time for radial trade}. In particular the solution carries the full regularity of Theorem \ref{linear solution full regularity}, \eqref{Improved solution estimate} together with the norms \eqref{Improved solution M2 norms}, at this $\tilde\kappa$, for the forcing $(\mb N,\mb N_\partial)[\phi,\mb u,p]$.
\end{proposition}
\begin{proof}
Write $\mathcal E_T^{\smash[t]{[1]}}$ for the sum of the three norms of $\phi$, $\mb u^{\smash[t]{[1]}}$ and $p$ in \eqref{upgrade estimate}, and $\mathcal E_T^{\smash[t]{[2]}}$ for the sum of the three norms of $\phi^{[2]}$, $\mb u^{[2]}$ and $p^{[2]}$; and, for a triple $w=(\psi,\mb v,q)$ satisfying the initial conditions of $B_T$, write $\mathcal E_T^{\trimtop{0.5pt}{[1]}}(w)$ and $\mathcal E_T^{\trimtop{0.5pt}{[2]}}(w)$ for the same norms formed on $w$, with the decompositions of Proposition \ref{time for radial trade}. Their coefficients are the corner value $\mb v(\mb 0,t)$, which also fixes the coefficient $v_1(\mb 0,t)/(\theta_b(\sin\theta_b)^2)$ of $r\ln r$ in $\psi$, and the functions $c(t)$ and $c_3(t)$ of that proposition. We write $\mathcal K_T^{\smash[t]{[1]}}:=\mathcal D+\Theta_T(\mathcal E_T+\mathcal E_T^{\smash[t]{[1]}})$, the counterpart of $\mathcal K_T$ for suprema in time of the norms in $\mathcal E_T^{\smash[t]{[1]}}$, and we use the abbreviation $\nn u\nn_{O,a,j,m}=\nn u\nn_{[0,T],O:(\varsigma,a)(0,j,m)}$ of the proof of Proposition \ref{space-time nonlinear estimates}. The two sums are recovered along the iteration of Theorem \ref{nonlinear positive-time existence}, with the solution operator $\mathcal T$ of its proof: Step 1 bounds $\mathcal E_T^{\smash[t]{[1]}}(\mathcal Tw)$ for $w\in B_T$, Step 2 shows that the start has finite rungs, Step 3 estimates the nonlinearities at the rung $M=2$, and Step 4 proves $\mathcal E_T^{\smash[t]{[2]}}(\mathcal Tw)\leq A+{1\over 2}\mathcal E_T^{\smash[t]{[2]}}(w)$ along the iterates and passes to the limit.

\emph{Step 1: the shifted rung $M=1$.} Let $w\in B_T$. We apply \eqref{time for radial trade estimate 1} at the exponent $\tilde\kappa$ to $\mathcal Tw$ and to $\partial_t\mathcal Tw$. Its hypotheses on the solution are $\lesssim_{\theta_b}C_{{\rm lin},\varsigma}(\mathcal N_T((\mb N,\mb N_\partial)[w])+\mathcal D)$, as in the proof of Theorem \ref{linear solution full regularity}. Its forcing hypotheses, evaluated on $(\mb N,\mb N_\partial)[w]$, are $\nn\partial_tf_\theta\nn_{[0,T],\partial_1\Omega:(\varsigma,\tilde\kappa-\kappa-1/2)(1,l+1/2)}$, bounded by \eqref{ST N estimate}, and the $M=1$ targets of \eqref{space-time targets} raised by $\tilde\kappa$, for which we claim
\begin{align}\label{shifted M1 targets}
\sum_{n=0}^1\brac{\nn\partial_t^n\mb N[w]\nn_{\Omega,\tilde\kappa-\kappa,l+2-n}+\nn\partial_t^n\mb N_\partial[w]\nn_{\partial_1\Omega,\tilde\kappa-\kappa+{1\over 2},l+{5\over 2}-n}}\lesssim_{\theta_b}\mathcal K_T\mathcal E_T(w).
\end{align}
This follows from the placements of Steps 2--4 of the proof of Proposition \ref{space-time nonlinear estimates} at $M=1$, with the shift placed on one geometry factor. Those placements put every geometry factor at $a=1$, except the factor $B$ of a pure geometry term, which sits at $a=2-\kappa$ and cannot be raised; we raise the factor $A$ of a pure geometry term and the geometry factor otherwise. The raised weight $1+\tilde\kappa$ lies in $[1-\kappa,2-\kappa]$, since $\tilde\kappa<\kappa<1/2$, and it is available by \eqref{available weights}, since for $n\leq 1$ a geometry factor carrying a time derivative is distinguished. For instance, for the summand $(\grad^2\bs\nu)v_y$ with all the radial derivatives on the solution factor,
\begin{align*}
\nn(\grad^2\bs\nu)\dr^{l+2-n}\partial_t^nv_y\nn_{\Omega,\tilde\kappa-\kappa}
\leq\sup_{t\in[0,T]}\|\grad^2\bs\nu(t)\|_{\Omega,\tilde\kappa-1}\norm{\|r^{\kappa-1}\dr^{l+2-n}\partial_t^nv_y\|_{L^\infty(\Omega)}}_{L^2_\varsigma[0,T]}
\lesssim_{\theta_b}\mathcal K_T\mathcal E_T(w),
\end{align*}
with the geometry at $a=1+\tilde\kappa$ in \eqref{space time geometry sup}, and the solution by the norm of $\mb v^{\smash[t]{[1]}}$ in $\mathcal E_T$, through Lemma \ref{embeddings}(ii), as in Step 3 of that proof. By \eqref{time for radial trade estimate 1}, \eqref{ST N estimate} and \eqref{shifted M1 targets}, since $\mathcal K_T\leq\delta_1$ and $\mathcal E_T(w)\leq\rho$ on $B_T$,
\begin{align}\label{shifted M1 bound}
\mathcal E_T^{[1]}(\mathcal Tw)\lesssim_{\theta_b}C_{{\rm lin},\varsigma}\mathcal D\qquad\text{ for every }w\in B_T,
\end{align}
and in particular for the solution, which is its own image under $\mathcal T$.

For a triple of $B_T$ with finite $\mathcal E_T^{\smash[t]{[1]}}$, the norms of $\partial_t\psi$, $\partial_t\mb v^{\smash[t]{[1]}}$ and $\partial_tq$ in $\mathcal E_T^{\smash[t]{[1]}}$ enlarge the ranges of Lemma \ref{space-time geometry estimates}. The space-time bounds \eqref{space time geometry L2} hold for $a\in[1-\kappa,2+\tilde\kappa-\kappa]$, which is admissible by \eqref{third rung angle condition}, with $\mathcal E_T+\mathcal E_T^{\smash[t]{[1]}}$ in place of $\mathcal E_T$, and \eqref{sup in time 1} gives
\begin{align}\label{shifted sup}
\sup_{t\in[0,T]}\brac{\|\psi(t)\|_{\partial_1\Omega,{3\over 2}+\tilde\kappa-\kappa,l+{5\over 2}}+\|\mb v^{[1]}(t)\|_{\Omega,2+\tilde\kappa-\kappa,l+1,2}+\|q(t)\|_{\Omega,1+\tilde\kappa-\kappa,l+1,1}}\lesssim_{\theta_b}\mathcal K_T^{[1]},
\end{align}
the initial values being bounded by $\mathcal D$ and \eqref{main estimate} through the splitting at $r=1$. Hence \eqref{space time geometry sup} holds for $a\in[1-\kappa,2+\tilde\kappa-\kappa]$ with $\mathcal K_T^{\smash[t]{[1]}}$ in place of $\mathcal K_T$.

\emph{Step 2: the start.} Let $w^{(0)}\in B_T$ be the start of Step 1 of the proof of Theorem \ref{nonlinear positive-time existence}, the solution with the smoothed flat forcing $(\mb F^\flat,\mb f^\flat)$ of \eqref{flat forcing}, and $w^{(n+1)}:=\mathcal T(w^{(n)})$, so that $w^{(n)}$ converges to the solution in the metric $d_T$. The rungs of $w^{(0)}$ are finite, though not necessarily small. Indeed, \eqref{time for radial trade estimate 1} and \eqref{time for radial trade estimate 2} bound them by hypotheses on the solution, controlled as in Step 1 and by $\mathcal E_T^{\smash[t]{[1]}}(w^{(0)})$, and by bounded functions of time times $\tilde\kappa$-shifted stationary norms of $(\mb N_0,\mb N_{\partial,0})$ and, at the counts lowered by one, of $(\dot{\mb N}_0,\dot{\mb N}_{\partial,0})$, the missing count being paid through \eqref{smoothing gain} as in Step 1 of that proof. At the $M\leq 1$ levels these norms are bounded by the unshifted $M\leq 1$ and $M=2$ norms, supplied by \eqref{N estimate}, \eqref{DN estimate} and \eqref{DN second rung estimate} at the counts $(l+1,l+{3\over 2})$. At the level $M=2$ they are \eqref{shifted N estimate} at $j=l+2$ and \eqref{shifted DN estimate} at $j=l+1$, whose hypotheses $\mathfrak X<\infty$ and $\dot{\mathfrak X}<\infty$ at these counts are given by Theorem \ref{jet corner structure}, with $\delta$ shrunk as there; this is where the hypotheses of that theorem enter.

\emph{Step 3: the forcing at the rung $M=2$.} We show that every triple $w\in B_T$ carrying the decompositions of Proposition \ref{time for radial trade}, with $\mathcal E_T^{\smash[t]{[1]}}(w)$ and $\mathcal E_T^{\smash[t]{[2]}}(w)$ finite, satisfies, with $(\mb F,\mb f)=(\mb N,\mb N_\partial)[w]$,
\begin{align}\label{M2 forcing estimate}
\nn\mb F\nn_{[0,T],\Omega:(\varsigma,1+\tilde\kappa-\kappa)(0,l+2)}+\nn\mb f\nn_{[0,T],\partial_1\Omega:(\varsigma,{3\over 2}+\tilde\kappa-\kappa)(0,l+{5\over 2})}
\lesssim_{\theta_b}\mathcal K_T^{[1]}\brac{\mathcal E_T+\mathcal E_T^{[1]}+\mathcal E_T^{[2]}},
\end{align}
all quantities on the right being those of $w$, which we write without arguments. For the coefficients, Lemma \ref{decomposition estimates}(i), applied at each $t$ to $\psi-\psi^{[2]}$ and $q-q^{[2]}$, and \eqref{spacetime coefficients} give
\begin{align}\label{upgrade coefficients}
\|\mb v(\mb 0,\ph)\|_{L^2_\varsigma[0,T]}+\|c\|_{L^2_\varsigma[0,T]}+\|c_3\|_{L^2_\varsigma[0,T]}\lesssim_{\theta_b}\mathcal E_T+\mathcal E_T^{[2]},\qquad\qquad
\sup_{t\in[0,T]}|\mb v(\mb 0,t)|\lesssim_{\theta_b}\mathcal K_T.
\end{align}
No bound on the suprema of $c$ and $c_3$ is available, since it would ask for the time derivative at the rung $M=2$.

By Lemma \ref{trace}, applied at each $t$, the LHS of \eqref{M2 forcing estimate} is at most a constant times $\nn\mb N\nn_{\Omega,1+\tilde\kappa-\kappa,l+2}+\nn H\nn_{\Omega,2+\tilde\kappa-\kappa,l+2,1}$, with $H$ as in Proposition \ref{space-time nonlinear estimates}. On $\{r>1\}$ these norms are dominated by $\nn\mb N\nn_{\Omega,-\kappa,l+2}$ and $\nn H\nn_{\Omega,1-\kappa,l+2,1}$, which are $\lesssim_{\theta_b}\mathcal K_T\mathcal E_T$ by \eqref{ST N estimate} and the bound on $H$ in Proposition \ref{space-time nonlinear estimates}. On $\{r<1\}$ the $M=2$ placements of the proof of Proposition \ref{nonlinear estimates}, shifted by $\tilde\kappa$, have to be read in time. No time derivative occurs, so the only constraint is that each product has exactly one factor in the space-time norm, the others being in the supremum in time; but several quantities are available in the space-time norm only, as we see by splitting the factors at the level $M=2$:
\begin{align}\label{upgrade corner splitting}
\bs\nu&=\bs\nu^{[2]}+c(t)\bs\nu[r]+{v_1(\mb 0,t)\over\theta_b(\sin\theta_b)^2}\bs\nu[r\ln r],\nonumber\\
\mb v&=\mb v(\mb 0,t)+{v_1(\mb 0,t)\over 2}\mb P_0(r,\theta-\theta_b)+\mb v^{[2]},\qquad\qquad
q=-c_3(t)-{v_1(\mb 0,t)\over\theta_b}\ln r+q^{[2]}.
\end{align}
Here $\bs\nu^{\smash[t]{[2]}}=\bs\nu[\psi^{\smash[t]{[2]}}]$, and $\bs\nu[r]$ and $\bs\nu[r\ln r]$ are the explicit fields of \eqref{block geometry}, so that $|\dr^j(r\grad)^m\bs\nu[r]|+|\dr^j(r\grad)^m\mb P_0|\lesssim_{\theta_b,j,m}r$ and $|\dr^j(r\grad)^m\bs\nu[r\ln r]|\lesssim_{\theta_b,j,m}r(1+|\ln r|)$. In the supremum in time we have, by Step 1, the unsplit solution factors at their weights in $\mathcal E_T^{\smash[t]{[1]}}$, the unsplit geometry at the weights $a\leq 2+\tilde\kappa-\kappa$, and $v_1(\mb 0,t)$. In the space-time norm only we have the coefficients $c$ and $c_3$, the remainders $\mb v^{\smash[t]{[2]}}$ and $q^{\smash[t]{[2]}}$ at the weights $3+\tilde\kappa-\kappa$ and $2+\tilde\kappa-\kappa$ of $\mathcal E_T^{\smash[t]{[2]}}$, and the remainder $\bs\nu^{\smash[t]{[2]}}$, at the weights up to $2+\tilde\kappa-\iota_2$ for a fixed $\iota_2\in[\kappa/2,\tilde\kappa)$: for $m'\leq l+2$,
\begin{align}\label{upgrade geometry remainder}
\nn\bs\nu^{[2]}\nn_{\Omega,2+\tilde\kappa-\iota_2,0,l+4}+\norm{\|r^{-(1+\tilde\kappa-\iota_2)}(r\grad)^{m'}\bs\nu^{[2]}\|_{L^\infty(\Omega)}}_{L^2_\varsigma[0,T]}
\lesssim_{\theta_b}\nn\psi^{[2]}\nn_{\partial_1\Omega,{3\over 2}+\tilde\kappa-\iota_2,l+{7\over 2}},
\end{align}
by \eqref{geometry bulk}--\eqref{geometry sup}, admissible since $2+\tilde\kappa-\iota_2\leq 2+\kappa/2<3-\kappa$ and $(3-\kappa)\theta_b<\pi/2$ by \eqref{third rung angle condition}. The RHS is $\lesssim_{\theta_b}\mathcal E_T+\mathcal E_T^{[2]}$ by the splitting at $r=1$: on $\{r<1\}$ the norm $\nn\psi^{[2]}\nn_{\partial_1\Omega,5/2+\tilde\kappa-\kappa,l+7/2}$ in $\mathcal E_T^{[2]}$ dominates; on $\{r>1\}$ the norm $\nn\psi\nn_{\partial_1\Omega,1/2-\kappa,l+7/2}$, contained in $\mathcal E_T$, dominates for $\psi$, while the subtracted terms have integrable tails, since ${3\over 2}+\tilde\kappa-\iota_2>{3\over 2}$, and coefficients bounded by \eqref{upgrade coefficients}. Since $\nn\psi^{\smash[t]{[2]}}\nn_{\partial_1\Omega,3/2+\tilde\kappa-\iota_2,l+7/2}$ and $\nn\psi\nn_{\partial_1\Omega,1/2-\kappa,l+7/2}$ are finite, the first line of \eqref{upgrade corner splitting} holds at almost every $t$ by Lemma \ref{corner splitting of the geometry}(i), applied to the block $c(t)r+v_1(\mb 0,t)r\ln r/(\theta_b(\sin\theta_b)^2)$, which has no term in $r\ln^2r$, at the boundary weight ${3\over 2}+\tilde\kappa-\iota_2\in({3\over 2},{5\over 2}-\kappa)$. The factor in the space-time norm is chosen as follows.

\emph{Terms with a solution factor.} In a term with a solution factor, other than those with all the radial derivatives on the geometry, the solution factor is split by \eqref{upgrade corner splitting}. Its remainder is distinguished, and the geometry is in the supremum in time at $a=1$, in $L^\infty$ or in $L^2$ in space according to the counts. The weights then add up to the target; for instance, for $j_1+j_2\leq l+2$ with $j_1\leq l+1$,
\begin{align*}
\nn(\dr^{j_1}\grad^2\bs\nu)(\dr^{j_2}v_y^{[2]})\nn_{\Omega\cap\{r<1\},1+\tilde\kappa-\kappa}&\leq\sup_{t\in[0,T]}\|\dr^{j_1}\grad^2\bs\nu(t)\|_{\Omega,-1}\norm{\|r^{-(2+\tilde\kappa-\kappa)}\dr^{j_2}v_y^{[2]}\|_{L^\infty(\Omega)}}_{L^2_\varsigma[0,T]}\\
&\lesssim_{\theta_b}\mathcal K_T\mathcal E_T^{[2]},
\end{align*}
by \eqref{space time geometry sup} at $a=1$ and Lemma \ref{embeddings}(ii). The corner part of the solution factor is paired with the geometry in $L^2$ in space at $a=2+\tilde\kappa-\kappa$, which is in the space-time norm when the coefficient is $v_1(\mb 0,t)$ and in the supremum in time when it is $c_3$. The coefficient $c_3$ occurs only in the undifferentiated pressure of $p\mb n_\phi$, and so does the term $v_1(\mb 0,t)(\ln r)\bs\nu$. The logarithm in this term costs a margin $\iota\in(0,\kappa-\iota_2]$: $v_1(\mb 0,t)$ is placed in the supremum in time, and $\bs\nu$ is split with its remainder distinguished at $a=2+\tilde\kappa-\kappa+\iota$.

\emph{Terms with all the radial derivatives on the geometry.} The placement above respects the counts unless the geometry carries all the radial derivatives, which would put it in the supremum in time beyond the counts of \eqref{space time geometry sup}: this happens in $(\dr^{l+2}\bs\nu)\grad^2\mb v$, $(\dr^{l+2}\bs\nu)\grad q$ and $(\dr^{l+2}\grad^2\bs\nu)v_y$ and, in $H$, in $(\grad\bs\nu)v_y$ with all the derivatives on the geometry. In these terms the geometry is split instead. Its remainder is distinguished at $a=2$, with the unsplit solution factor in the supremum in time at its weight in $\mathcal E_T^{\smash[t]{[1]}}$, and the weights again add up to the target. The piece $c(t)\bs\nu[r]$, which by the pointwise bound is of the size of a geometry factor at $a=2$, has $c$ distinguished and the solution factor again in the supremum in time. The piece with $\bs\nu[r\ln r]$ has its coefficient in the supremum in time and the solution factor split, the remainder of the latter being distinguished and absorbing the logarithm. For $(\dr^{l+2}\grad^2\bs\nu)v_y$,
\begin{align*}
\nn(\dr^{l+2}\grad^2\bs\nu^{[2]})v_y\nn_{\Omega\cap\{r<1\},1+\tilde\kappa-\kappa}&\leq\sup_{t\in[0,T]}\|r^{-(1+\tilde\kappa-\kappa)}v_y(t)\|_{L^\infty(\Omega)}\nn\grad^2\bs\nu^{[2]}\nn_{\Omega\cap\{r<1\},0,l+2},\\
\nn c(\dr^{l+2}\grad^2\bs\nu[r])v_y\nn_{\Omega\cap\{r<1\},1+\tilde\kappa-\kappa}&\lesssim_{\theta_b}\|c\|_{L^2_\varsigma[0,T]}\sup_{t\in[0,T]}\|v_y(t)\|_{\Omega,2+\tilde\kappa-\kappa}.
\end{align*}
Here $\nn\grad^2\bs\nu^{\smash[t]{[2]}}\nn_{\Omega\cap\{r<1\},0,l+2}$ is bounded by \eqref{upgrade geometry remainder}, since $\tilde\kappa-\iota_2>0$, and the suprema in time are $\lesssim_{\theta_b}\mathcal K_T^{\smash[t]{[1]}}$ by \eqref{shifted sup}, $v_y$ being $v_y^{[1]}$, the first through Lemma \ref{embeddings}(ii).

\emph{Pure geometry terms.} In a pure geometry term $r^{-1}AB$ of Case 4, $B$ carries the most derivatives and is split, $B^{[2]}$, $B[r]$ and $B[r\ln r]$ denoting the same derivative of $\bs\nu^{[2]}$, $\bs\nu[r]$ and $\bs\nu[r\ln r]$. The remainder is placed as in the stationary proof at $M=2$, with the margin $\iota_2$ in place of $\kappa/2$:
\begin{align*}
\nn r^{-1}AB^{[2]}\nn_{\Omega\cap\{r<1\},2+\tilde\kappa-\kappa}\leq\sup_{t\in[0,T]}\|r^{-(1-\kappa+\iota_2)}A(t)\|_{L^\infty(\Omega)}\nn B^{[2]}\nn_{\Omega,2+\tilde\kappa-\iota_2},
\end{align*}
with $A$ in the supremum in time at $a=2-\kappa+\iota_2\leq 2+\tilde\kappa-\kappa$. The piece $c(t)B[r]$ has $c$ distinguished, with $A$ in $L^2$ in space at $a=2+\tilde\kappa-\kappa$ and in the supremum in time. The piece with $B[r\ln r]$ has its coefficient in the supremum in time; one factor of $A$ is split in turn, with its remainder distinguished at $a=2+\tilde\kappa-\kappa+\iota$.

\emph{Products of corner parts.} What remains are products of corner parts, with further geometry factors in the supremum in time at $a=1$. They are explicit functions of degree at least $0$ in $\mb N$ and at least $1$ in $H$, up to powers of $\ln r$, times two coefficients, one of which is $v_1(\mb 0,t)$ or the coefficient of $r\ln r$ and is placed in the supremum in time. They are integrable against the target weights on $\{r<1\}$ because $\tilde\kappa<\kappa$: the measure on $\Omega$ being $r\d r\d\theta$, for $\mb N$
\begin{align*}
\int_0^1r^{-2(1+\tilde\kappa-\kappa)}(1+|\ln r|)^{2k}r\,\d r=\int_0^1r^{2(\kappa-\tilde\kappa)-1}(1+|\ln r|)^{2k}\d r<\infty,
\end{align*}
and for $H$ the weight and the degree are both one higher.
Every term is thus a supremum in time, $\lesssim_{\theta_b}\mathcal K_T^{[1]}$, times a space-time norm, $\lesssim_{\theta_b}\mathcal E_T+\mathcal E_T^{[1]}+\mathcal E_T^{[2]}$ by Step 1, \eqref{upgrade coefficients} and \eqref{upgrade geometry remainder}. This proves \eqref{M2 forcing estimate}.

\emph{Step 4: the rung $M=2$ along the iteration.} The forcing hypotheses of \eqref{time for radial trade estimate 2} are the norms on the left of \eqref{M2 forcing estimate}, bounded in Step 3, and $\nn\partial_tf_\theta\nn_{\partial_1\Omega,1/2+\tilde\kappa-\kappa,l+3/2}$, bounded in Step 1. Inserting \eqref{M2 forcing estimate} for $w^{(n)}$ into \eqref{time for radial trade estimate 2} for $w^{(n+1)}$, whose remaining hypotheses are the norms of $\partial_t\mb v^{[1]}$ and $\partial_tq$ in $\mathcal E_T^{\trimtop{1.0pt}{[1]}}(w^{(n+1)})$, bounded by \eqref{shifted M1 bound}, gives
\begin{align*}
\mathcal E_T^{[2]}(w^{(n+1)})\lesssim_{\theta_b}C_{{\rm lin},\varsigma}\mathcal D+\mathcal K_T^{[1]}(w^{(n)})\brac{\mathcal E_T(w^{(n)})+\mathcal E_T^{[1]}(w^{(n)})+\mathcal E_T^{[2]}(w^{(n)})}.
\end{align*}
By Step 2 and induction on $n$, every iterate carries the decompositions with finite rungs. For $n\geq 1$, \eqref{shifted M1 bound}, $\mathcal E_T(w^{(n)})\leq\rho$ and \eqref{positive time smallness} give $\mathcal K_T^{\smash[t]{[1]}}(w^{(n)})\lesssim_{\theta_b}C_{{\rm lin},\varsigma}(1+C_{\rm ball})\delta$. So after shrinking $\delta$ so that $C(\theta_b,\kappa,\tilde\kappa,l)C_{{\rm lin},\varsigma}(1+C_{\rm ball})\delta\leq{1\over 2}$, which holds for $\delta\lesssim_{\theta_b}(1+\varsigma^{-1})^{-2}$ and therefore preserves the form \eqref{delta tau form}, the RHS is at most $A+{1\over 2}\mathcal E_T^{\smash[t]{[2]}}(w^{(n)})$ with $A\lesssim_{\theta_b}\mathcal D$, and $\limsup_n\mathcal E_T^{\smash[t]{[2]}}(w^{(n)})\leq 2A$, whatever the size of the start.

 In the limit, $w^{(n)}$ converges to the solution in $d_T$, hence in the sense of distributions, while the coefficients of the decompositions, in $L^2_\varsigma[0,T]$, and their remainders at the rung $M=2$ are bounded uniformly in $n$, by \eqref{upgrade coefficients} and the bound just proved. So Lemma \ref{corner blocks under limits}(iii), with the corner blocks of Proposition \ref{time for radial trade} as the block class, at the weights $3+\tilde\kappa-\kappa>2$ and $2+\tilde\kappa-\kappa>1$ in the bulk and ${5\over 2}+\tilde\kappa-\kappa>{3\over 2}$ on the boundary, gives the solution the decompositions and the bound $\mathcal E_T^{\smash[t]{[2]}}\leq\liminf_n\mathcal E_T^{\smash[t]{[2]}}(w^{(n)})\lesssim_{\theta_b}\mathcal D$. Its uniqueness clause identifies the limiting decompositions with those of Proposition \ref{time for radial trade} for the solution. This proves \eqref{upgrade estimate}.
\end{proof}

By their proofs, the shrinkings of $\delta$ in Propositions \ref{Lipschitz positive} and \ref{positive-time regularity upgrade} preserve the form \eqref{delta tau form}. So Corollary \ref{polynomial lifespan}, with a smaller $c_*$, covers these propositions as well, on $[0,T_*]$ with the weight $\varsigma=1/T_*$; the constant of \eqref{upgrade estimate} is not tracked in $\varsigma$.

The regularity just established improves, through the profile evolution equation, the corner behaviour of the second time derivative of the profile at almost every time --- not only at $t=0$, where the compatibility conditions were assumed, but along the whole flow. The improvement stops strictly before the critical weight ${1\over 2}$, a threshold that is structural, see Remark \ref{positive-time closing remark}.

\begin{proposition}[The second time derivative at the contact point]\label{second time derivative at the corner}
Let $(\phi,\mb u,p)$ be the solution of Theorem \ref{nonlinear positive-time existence} with the regularity \eqref{upgrade estimate}, and let $\tilde\kappa\in(\kappa/2,\kappa)$ be as there. Then
\begin{align}\label{second derivative at the corner estimate}
(\sin\theta_b)^2\nn\partial_t^2\phi\nn_{[0,T],\partial_1\Omega:(\varsigma,{1\over 2}+\tilde\kappa-\kappa)(0,l+{1\over 2})}\lesssim_{\theta_b}\mathcal D,
\end{align}
and at a.e.\ $t\in[0,T]$
\begin{align*}
|\partial_t^2\phi(r,t)|\lesssim_{\theta_b}r^{\tilde\kappa-\kappa}\|\partial_t^2\phi(t)\|_{\partial_1\Omega,{1\over 2}+\tilde\kappa-\kappa,{3\over 2}}.
\end{align*}
\end{proposition}
\begin{proof}
Differentiating the profile evolution equation of \eqref{system: positive time} in time gives $\partial_t^2\phi=(\sin\theta_b)^{-2}\partial_r\partial_tu_\theta$ on $\partial_1\Omega$ a.e., both sides lying in the class of $\mathcal E_T$, as in Step 4 of the proof of Theorem \ref{nonlinear positive-time existence}. Since $\mb e_\theta$ is independent of $r$ and $\dr$ annihilates the constant $\partial_t\mb u(\mb 0,t)$, we have $\partial_r\partial_tu_\theta=\mb e_\theta\cdot r^{-1}\dr\partial_t\mb u^{[1]}$, so that, by Lemma \ref{trace}, at a.e.\ time
\begin{align*}
(\sin\theta_b)^2\|\partial_t^2\phi(t)\|_{\partial_1\Omega,{1\over 2}+\tilde\kappa-\kappa,l+{1\over 2}}
\lesssim_{\theta_b}\|\partial_t\mb u^{[1]}(t)\|_{\Omega,2+\tilde\kappa-\kappa,l+1,2},
\end{align*}
whose square integrates in time to $\nn\partial_t\mb u^{\smash[t]{[1]}}\nn_{[0,T],\Omega:(\varsigma,2+\tilde\kappa-\kappa)(0,l+1,2)}^2$, bounded by \eqref{upgrade estimate}. This proves \eqref{second derivative at the corner estimate}, and the pointwise bound is Lemma \ref{embeddings}(i) at $a=\tilde\kappa-\kappa$.
\end{proof}

The contact angle condition, differentiated once in time, holds pointwise along the solution at almost every time. In the decomposition of $\phi(t)$ in Proposition \ref{time for radial trade}, this condition fixes the coefficient $u_1(\mb 0,t)/(\theta_b(\sin\theta_b)^2)$ of $r\ln r$, and so it propagates the first compatibility condition.

\begin{corollary}[The first compatibility condition along the flow]\label{C1 along the flow}
Let $(\phi,\mb u,p)$ be the solution of Theorem \ref{nonlinear positive-time existence}. For a.e.\ $t\in[0,T]$ the derived profile $\partial_t\phi(t)$ extends continuously to the contact point by $0$. Under the hypotheses of Proposition \ref{positive-time regularity upgrade}, for a.e.\ $t\in[0,T]$ the profile $\phi(t)$ satisfies \eqref{C1 in the linear theory}, that is, the first compatibility condition \eqref{C1} formed at the time $t$:
\begin{align}\label{C1 along the flow identity}
c_{\phi(t),\dr_{-1}^2}={u_1(\mb 0,t)\over\theta_b(\sin\theta_b)^2}.
\end{align}
Under these hypotheses the data norm $\mathcal D(\phi(t))$ is also finite for a.e.\ $t\in[0,T]$, and so is $\mathcal E(\mb u(t),p(t))$, formed with $\phi(t)$ in place of $\phi_0$.
\end{corollary}
\begin{proof}
By Lemma \ref{embeddings}(i) at $a=1-\kappa$, $|\partial_t\phi(r,t)|\lesssim_{\theta_b}r^{1-\kappa}\|\partial_t\phi(t)\|_{\partial_1\Omega,3/2-\kappa,l+5/2}$, which is finite for a.e.\ $t$ since $\mathcal E_T<\infty$; this is the first claim. By \eqref{upgrade estimate}, for a.e.\ $t$ the remainder $\phi^{\smash[t]{[2]}}(t)=\phi(t)-c(t)r-u_1(\mb 0,t)r\ln r/(\theta_b(\sin\theta_b)^2)$ of Proposition \ref{time for radial trade} has finite norm $\|\phi^{[2]}(t)\|_{\partial_1\Omega,5/2+\tilde\kappa-\kappa,l+7/2}$, and $\|\phi(t)\|_{\partial_1\Omega,3/2-\kappa}$ is finite since $\mathcal E_T<\infty$. At such $t$, Lemma \ref{block matching} with $(a_0,d,a_1)=(1-\kappa,1,2+\tilde\kappa-\kappa)$, applied to $f_0=\phi(t)$ and $f_1=\phi^{\smash[t]{[2]}}(t)$, gives $\smash[b]{c_{\phi(t),\dr_{-1}}}=c(t)$ and \eqref{C1 along the flow identity}. So the decomposition \eqref{data step relations} of $\phi(t)$ is that of Proposition \ref{time for radial trade}. By \eqref{C1 along the flow identity} and Remark \ref{radial trade and corner ladder}, the same holds for $(\mb u,p)(t)$, decomposed as in \eqref{data step relations} with $\phi(t)$ in place of $\phi_0$. For $p^{[2]}$ this also uses the relation between $c$ and $c_3$ of Proposition \ref{time for radial trade}. Hence $\mathcal D(\phi(t))$ and $\mathcal E(\mb u(t),p(t))$ are finite, by the norms in $\mathcal E_T$ and by \eqref{upgrade estimate}, together with the splitting at $r=1$.
\end{proof}

\begin{remark}[On the positive-time theorem]\label{positive-time closing remark}
Three comments on Theorem \ref{nonlinear positive-time existence}.
\begin{enumerate}
\item (Initial values) For a triple with finite $\mathcal E_T$ and for the tuple of Theorem \ref{nonlinear initial data},
\begin{align*}
\|\partial_t\psi(0)\|_{\partial_1\Omega,{1\over 2}-\kappa,l+2}+\|\partial_t\mb v(0)\|_{\Omega,1-\kappa,l+{1\over 2},2}&\lesssim\mathcal E_T,\\
\|\dot\phi_0\|_{\partial_1\Omega,{1\over 2}-\kappa,l+{5\over 2}}+\|\dot{\mb u}_0\|_{\Omega,1-\kappa,l+1,2}&\lesssim_{\theta_b}\mathcal D,
\end{align*}
the first bound by \eqref{sup in time 4}, with a constant depending on $T$ and $\varsigma$, and the second by \eqref{main estimate}. Step 1 of the proof uses the corresponding extra half radial derivative of the forcing values \eqref{initial forcing values} to construct an element of $B_T$ by radial smoothing. Conversely, the initial data is not in the trace space of the positive-time estimate, by Remark \ref{linear theory closing remark}, so the iteration does not control the initial values of its iterates. Moreover, the initial jet must be compatible with the values of the forcing at $t=0$, which depend on the initial jet. Indeed, \eqref{C1} ties $\smash[b]{c_{\phi_0,\dr_{-1}^2}}$ to the corner value $u_1(\mb 0)$ of the stationary solution with these values, so iterates with their own initial values would satisfy \eqref{C1} only in the limit. For these two reasons Theorem \ref{nonlinear initial data} computes the initial values, with norms bounded by $\mathcal D$ through \eqref{main estimate}, and $B_T$ prescribes them; Theorem \ref{nonlinear positive-time existence} takes exactly what Theorem \ref{nonlinear initial data} gives.
\item (Restart at a later time) Assume in addition the hypotheses of Proposition \ref{positive-time regularity upgrade}, and consider a restart of the construction at a later time $t_0$. By this proposition the solution is as regular as the linear theory allows. By Corollary \ref{C1 along the flow}, for a.e.\ $t_0$, the profile $\phi(t_0)$ has finite data norm $\mathcal D$ at the same parameter $l$ and satisfies the first compatibility condition, and $\mathcal E(\mb u(t_0),p(t_0))$, formed with $\phi(t_0)$ in place of $\phi_0$, is finite. The slice $(\mb u,p)(t_0)$ solves \eqref{system: stationary nonlinear} with $\phi(t_0)$ in place of $\phi_0$, by \eqref{system: positive time}. So, wherever $\mathcal D(\phi(t_0))$ is small, the slice is the pair that Theorem \ref{nonlinear initial data}(i) constructs from $\phi(t_0)$, by the uniqueness stated there, and parts (i) and (ii) of that theorem apply to $\phi(t_0)$. A hypothesis at the parameter $l+1$, by contrast, would not be reproduced by the class. What a restart would lack even there is the hypothesis of Theorem \ref{nonlinear initial data}(iii), the second compatibility condition for $\phi(t_0)$.
\item (The second compatibility condition along the flow) Assume in addition the hypotheses of Proposition \ref{positive-time regularity upgrade}. Whether \eqref{C2 explicit} holds along the flow, and with it whether the solution can be restarted, is left open: of the two compatibility conditions, both hypotheses at the initial time, only the first is a theorem at almost every later time. The contact angle condition holds classically at every time, by Theorem \ref{the contact point}(iii), and once differentiated in time it holds pointwise at a.e.\ time, by Corollary \ref{C1 along the flow}; the pattern stops at the second derivative, since the trace at the contact point need not commute with two time derivatives. At $t=0$ the corner value of $\partial_t^2\phi$ vanishes if and only if \eqref{C2 explicit} holds, by Lemma \ref{second-order corner values}. The obstruction is the second-order logarithm that the quadratic nonlinearity creates out of the first-order corner data.

Along the flow, Proposition \ref{second time derivative at the corner} gives $\|\partial_t^2\phi(t)\|_{\partial_1\Omega,1/2+\tilde\kappa-\kappa,l+1/2}<\infty$ at a.e.\ time, which still admits a bounded nonzero corner value. The finiteness of $\|\partial_t^2\phi(t)\|_{\partial_1\Omega,1/2}$ is the weak form of the twice time-differentiated contact angle condition, and hence of \eqref{C2 explicit} along the flow. The bound of Proposition \ref{second time derivative at the corner} cannot reach the critical weight ${1\over 2}$. The threshold is structural: the boundary weight ${1\over 2}$ corresponds to the Mellin line $\Re\lambda=0$, the line through the double pole $\lambda=0$ of the system, and every shift of the linear theory is capped by the first pole it meets, which is the constraint $\tilde\kappa<\kappa$ of Proposition \ref{time for radial trade}. Nor does a cut-off in time circumvent the threshold. A cut-off makes the initial data of the differentiated problem compatible to every order, but at the price of a boundary forcing that the unweighted block of \eqref{time solution estimate} asks, for the twice-differentiated problem, at the same boundary weight ${1\over 2}$.

Crossing $\Re\lambda=0$ would ask for a dotted second rung carried along the flow, that is, the time-dependent form of the corner-block subtraction of Propositions \ref{third rung} and \ref{nonlinear corner structure}. The line would be crossed not by an estimate but by the vanishing of the coefficients of the subtracted blocks, forced by the once-differentiated contact angle condition, as for the first compatibility condition in Corollary \ref{C1 along the flow}. We do not carry this out.
\end{enumerate}
\end{remark}

\subsubsection{The contact point}

We collect what the solution does at the contact point. Put together, the statements below describe the flow near a moving contact line with Navier slip and a fixed angle; their physical content is discussed in Remark \ref{contact point physics}.

\begin{theorem}[The contact point]\label{the contact point}
Let $(\phi,\mb u,p)$ be the solution of Theorem \ref{nonlinear positive-time existence} and let
\begin{align}\label{contact-point velocity}
U(t):=u_1(\mb 0,t)
\end{align}
be the horizontal velocity at the contact point. Then parts (i) and (iii) below hold, and, if the solution has the regularity \eqref{upgrade estimate} for some $\tilde\kappa\in(\kappa/2,\kappa)$, so does part (ii).
\begin{enumerate}
\item (Velocity) $U$ is continuous on $[0,T]$, with $\sup_{[0,T]}|U|\lesssim_{\theta_b}(1+(1+\varsigma^{-1})\Theta_T)\mathcal D$, and $\mb u(\mb 0,t)=U(t)\mb e_1$ for every $t$.
\item (Expansions) For almost every $t\in[0,T]$, with the field $\mb P_0$ of Proposition \ref{time for radial trade} and functions $c(t)$, $c_3(t)$ of time alone,
\begin{align}
\mb u(r,\theta,t)&=U(t)\mb e_1+{U(t)\over 2}\mb P_0(r,\theta-\theta_b)+O\brac{r^{2+\tilde\kappa-\kappa}},\label{contact point u expansion}\\
p(r,\theta,t)&=-c_3(t)-{U(t)\over\theta_b}\ln r+O\brac{r^{1+\tilde\kappa-\kappa}},\label{contact point p expansion}\\
\phi(r,t)&=c(t)r+{U(t)\over\theta_b(\sin\theta_b)^2}r\ln r+O\brac{r^{2+\tilde\kappa-\kappa}},\label{contact point phi expansion}
\end{align}
as $r\to 0$, uniformly in $\theta$, with the constants in the $O$ terms in $L^2$ of time against $e^{-\varsigma t}$ and bounded by $C(\theta_b,\kappa,\tilde\kappa,l,\varsigma)\mathcal D$. The functions $c$ and $c_3$ are related by
\begin{align*}
c_3(t)=(\sin\theta_b)^2c(t)+\brac{{1\over 2\theta_b}+{1\over\sin(2\theta_b)}}U(t).
\end{align*}
At $t=0$ the expansions hold with the coefficients of Theorem \ref{nonlinear initial data}.
\item (Contact angle and position) $\phi(0,t)=0$ for every $t\in[0,T]$ and $|\phi(r,t)|\lesssim_{\theta_b}r^{1-\kappa}(1+(1+\varsigma^{-1})\Theta_T)\mathcal D$. In the moving domain of Corollary \ref{moving domain} the contact point moves with the velocity $U$,
\begin{align}\label{contact point position}
X(0,t)=X(0,0)+\int_0^tU(t')\d t',
\end{align}
the fluid velocity at the contact point of the moving domain being $U(t)\mb e_1$ as well.
\end{enumerate}
\end{theorem}
\begin{proof}
(i) By \eqref{spacetime coefficients}, $\mb u(\mb 0,t)$ is bounded by $\mathcal K_T=\mathcal D+\Theta_T\mathcal E_T$, and $\mathcal K_T\lesssim_{\theta_b}(1+(1+\varsigma^{-1})\Theta_T)\mathcal D$ by \eqref{positive time main estimate} and \eqref{delta tau form}. The value $\mb u(\mb 0,t)$ is continuous in $t$ because, by Lemma \ref{supremum in time} applied to the $M\leq 1$ rungs of $\mathcal E_T$, $\mb u$ and $\mb u^{[1]}$ have representatives continuous in time with values in the weighted spaces at the weights $1-\kappa$ and $2-\kappa$, and \eqref{corner coefficients 1} applied to $\mb u(t)-\mb u(t')$ bounds $|\mb u(\mb 0,t)-\mb u(\mb 0,t')|$. The second component vanishes by Lemma \ref{decomposition estimates}(iii).

(ii) The decompositions are those of Proposition \ref{time for radial trade}, carried by the solution by Proposition \ref{positive-time regularity upgrade}, with $\mb u^{[2]}$, $p^{[2]}$ and $\phi^{[2]}$ in the weighted spaces of \eqref{upgrade estimate} at almost every time; the remainders are bounded pointwise by Lemma \ref{embeddings}(ii) at $a=2+\tilde\kappa-\kappa$ and $a=1+\tilde\kappa-\kappa$ for $\mb u^{[2]}$ and $p^{[2]}$, and by Lemma \ref{embeddings}(i) at $a=2+\tilde\kappa-\kappa$ for $\phi^{[2]}$, the norms on the right being summands of \eqref{upgrade estimate}. The absence of a logarithm in \eqref{contact point u expansion} is the continuity of $\mb u$ at $\mb 0$ stated in Proposition \ref{corner ladder}, applied at almost every time, and the coefficients of the logarithms in \eqref{contact point p expansion} and \eqref{contact point phi expansion} are those of Proposition \ref{time for radial trade}, read off \eqref{data step relations} under \eqref{C1 in the linear theory}, which holds at almost every time by Corollary \ref{C1 along the flow}; at $t=0$ they are those of Theorem \ref{nonlinear initial data}, and \eqref{C1} is the statement that the coefficient of $r\ln r$ in $\phi_0$ equals $U(0)/(\theta_b(\sin\theta_b)^2)$. The relation between $c$ and $c_3$ is that of Proposition \ref{time for radial trade}.

(iii) By \eqref{sup in time 1}, with the initial value bounded by $\mathcal D$ and $\nn\partial_t\phi\nn_{[0,T],\partial_1\Omega:(\varsigma,3/2-\kappa)(0,l+5/2)}$ by $\mathcal E_T$, we have $\sup_{t\in[0,T]}\|\phi(t)\|_{\partial_1\Omega,3/2-\kappa,l+5/2}\lesssim_{\theta_b}\mathcal K_T$, which is bounded as in (i), and Lemma \ref{embeddings}(i) at $a=1-\kappa$ gives the bound on $|\phi(r,t)|$, hence $\phi(0,t)=0$, at every $t$. For \eqref{contact point position}, in Corollary \ref{moving domain} the free boundary is $X=b^{-1}y+\chi$ with $\partial_t\chi=-(\sin\theta_b)^{-1}u_\theta|_{\partial_1\Omega}$, and at $y=0$ this is $\partial_tX(0,t)=-(\sin\theta_b)^{-1}\mb u(\mb 0,t)\cdot\mb e_\theta(\theta_b)=-(\sin\theta_b)^{-1}U(t)\mb e_1\cdot\mb e_\theta(\theta_b)=U(t)$. That the flattening $\bs\Psi_t$ of Lemma \ref{regularised flattening} has the Jacobian $I$ at the corner follows from $|\bs\nu[\phi(t)]|\lesssim_{\theta_b}r^{1-\kappa}$, which is \eqref{space time geometry sup} at $a=2-\kappa$ and $m'=0$; hence the Piola transform of Lemma \ref{Piola transform} does not change the velocity at the contact point.
\end{proof}

By Theorem \ref{the contact point}(iii), the solution of Corollary \ref{moving domain} satisfies the contact angle condition $X'(0,t)=b^{-1}$ of System \ref{dimensionless system} classically for all $t\in[0,T]$. Together with Corollary \ref{polynomial lifespan} this proves part (i) of Theorem \ref{main theorem}.

\begin{remark}[The physics of the contact point]\label{contact point physics}
Three comments on Theorem \ref{the contact point}, in the dimensional variables of the first system of Section \ref{setup}, related to the dimensionless ones by measuring lengths in units of the slip length $\beta$, velocities in units of $\gamma/\mu$ and pressures in units of $\gamma/\beta$.
\begin{enumerate}
\item The logarithm of the pressure is the slip-regularised form of the Huh--Scriven singularity of the no-slip moving contact line \cite{HuhScriven}. With Navier slip the velocity is continuous at the contact point and the force on the substrate is finite. What remains of the singularity is the logarithm $-\mu U(\beta\theta_b)^{-1}\ln(r/\beta)$ of the pressure, with coefficient the viscous stress $\mu U/\beta$ at the slip scale divided by the angle. In molecular-dynamics simulations of the moving contact line, the normal stress along the substrate, averaged over the adjacent layer of fluid, has a jump across the contact line region that is likewise proportional to the speed of the contact line \cite{RenE}. The theorem proves this asymptotics for the free-boundary problem with surface tension, and not only for the wedge with prescribed boundaries; the price is the two compatibility conditions, without which the expansion \eqref{contact point phi expansion} cannot be propagated.
\item The first compatibility condition is a normal stress balance made explicit. On $\partial_1\Omega$ the $\mb e_\theta$ component of the surface tension condition reads $-2\partial_ru_r-p=(\sin\theta_b)^2\partial_r\phi+f_\theta$, and inserting the logarithm of the pressure gives $\partial_r\phi=U\ln r/(\theta_b(\sin\theta_b)^2)+O(1)$, whose integral is the $r\ln r$ term of \eqref{contact point phi expansion}. In the variables of the moving domain, where $\phi=b^{-1}-X'$, the slope of the free boundary near the contact point is
\begin{align*}
X'(y,t)=b^{-1}-c_\beta(t)y-{\mathrm{Ca}(t)\over\theta_b(\sin\theta_b)^3}{y\over\beta}\ln{y\over\beta}+O\brac{\brac{y\over\beta}^{2+\tilde\kappa-\kappa}},\qquad\mathrm{Ca}(t):={\mu U(t)\over\gamma},
\end{align*}
with $c_\beta(t)=(\beta\sin\theta_b)^{-1}c(t)-\mathrm{Ca}(t)\ln(\sin\theta_b)/(\beta\theta_b(\sin\theta_b)^3)$ a function of time alone, $y=r\sin\theta_b$ the height above the substrate, and $\mathrm{Ca}$ the capillary number of the contact line. This is the counterpart, at scales $y\ll\beta$, of the Cox--Voinov relation \cite{Voinov,Cox}, which describes the apparent angle at scales $y\gg\beta$ by a logarithm of $y/\beta$ in the cube of the angle itself. Here the logarithm sits in the slope, multiplied by $y$. So the angle is preserved at the contact point and the correction to it vanishes there, which is the fixed contact angle of System \ref{dimensionless system}.

The sign is the expected one. In the orientation of Section \ref{setup} the fluid occupies $\{x\geq X\}$, so $U>0$ is a receding contact line. Then $\phi<0$ sufficiently close to the contact point, that is $X'>b^{-1}$ as $y\to 0$, since $\ln(y/\beta)<0$ there and the logarithm dominates the linear term. So the free surface $x=X(y)$ meets the substrate at the angle $\theta_b$ but bends towards the substrate immediately above it. The apparent angle at scales $y\lesssim\beta$ is therefore lowered, as it is for a receding contact line in the Cox--Voinov description at larger scales, and it is raised for an advancing one.
\item The contact-line velocity $U(t)$ is continuous in time, and the position \eqref{contact point position} of the contact point is a $C^1$ function of time. The velocity is not, however, a free parameter: $U(0)$ is determined by the initial profile through the fixed point of Theorem \ref{nonlinear initial data}, and the first compatibility condition then fixes the coefficient of $y\ln y$ in the initial slope in terms of it. This is the sense in which the class of initial free boundaries satisfying the first compatibility condition is of codimension one, Theorem \ref{compatible set}.
\end{enumerate}
\end{remark}

\section*{Acknowledgements}

Manuel V. Gnann and King Ming Lam are supported by the Dutch Research Council (NWO) through the grant VI.Vidi.223.019, and King Ming Lam additionally through the grant OCENW.M20.194 awarded to Manuel V. Gnann.

The first drafts of the mostly routine nonlinear estimates of Section \ref{nonlinear existence} and of the third rung residue calculus of Section \ref{third rung subsection} were generated with Claude, from the linear theory and the first and second rung residue calculus written by the authors, and were then refined and adjusted by them. The rest of the mathematics was written without the use of AI, with the exception of the proof of Korn's inequality; Claude and ChatGPT were used in the preparation of the final version to refine the language, to catch and correct typos and slips, and as an editing tool for systematic updates to notation, labelling and referencing. All AI-generated output was reviewed and confirmed by the authors before inclusion.

\appendix

\section{Useful identities}\label{useful identities}

We collect here the identities used throughout the paper: the dictionary between the polar and the Cartesian descriptions of a vector field on the wedge, and the two elementary Mellin transforms behind the corner expansions of Section \ref{corner expansions}.

\subsection{Polar coordinates}\label{polar coordinates}

The reference wedge $\Omega$ of \eqref{fixed wedge reference domain} is a cone, so that the natural coordinates on it are the polar ones $(r,\theta)$ and the natural frame is $\mb e_r,\mb e_\theta$; the boundary conditions of System \ref{dimensionless system} and the corner expansions of Section \ref{corner expansions}, on the other hand, are read in the Cartesian components. The identities below translate between the two, in increasing order of the number of derivatives they involve.

The components are related by $u_1=u_r\cos\theta-u_\theta\sin\theta$ and $u_2=u_r\sin\theta+u_\theta\cos\theta$, and inversely by $u_r=u_1\cos\theta+u_2\sin\theta$ and $u_\theta=-u_1\sin\theta+u_2\cos\theta$; the frame depends on $\theta$ alone, with $\partial_\theta\mb e_r=\mb e_\theta$ and $\partial_\theta\mb e_\theta=-\mb e_r$, and both fields are homogeneous of degree $0$ in $r$, so that $\dr\mb e_r=\dr\mb e_\theta=\mb 0$.

The gradient and its rotation by a right angle are
\begin{align*}
\grad&=\mb e_r\partial_r+\mb e_\theta r^{-1}\partial_\theta\\
\grad^\perp&=(-\partial_2,\partial_1)=\mb e_\theta\partial_r-\mb e_rr^{-1}\partial_\theta,
\end{align*}
so that $\partial_r=\mb e_r\cdot\grad=\cos\theta\partial_1+\sin\theta\partial_2$ and, the frame being $\theta$ dependent,
\begin{align*}
\grad\cdot\mb u&=(\mb e_r\partial_r+\mb e_\theta r^{-1}\partial_\theta)\cdot(\mb e_ru_r+\mb e_\theta u_\theta)=\partial_ru_r+r^{-1}\partial_\theta u_\theta+r^{-1}u_r\\
\grad^\perp\cdot\mb u&=(\mb e_\theta\partial_r-\mb e_rr^{-1}\partial_\theta)\cdot(\mb e_ru_r+\mb e_\theta u_\theta)=\partial_ru_\theta-r^{-1}\partial_\theta u_r+r^{-1}u_\theta.
\end{align*}
As a tensor,
\begin{align}\label{grad polar}
\grad\mb u&=(\mb e_1\partial_1+\mb e_2\partial_2)\otimes(u_1\mb e_1+u_2\mb e_2)=\sum_{i,j\in\{1,2\}}(\partial_iu_j)\mb e_i\otimes\mb e_j\nonumber\\
&=(\mb e_r\partial_r+\mb e_\theta r^{-1}\partial_\theta)\otimes(u_r\mb e_r+u_\theta\mb e_\theta)\nonumber\\
&=\mb e_r\otimes(\mb e_r\partial_ru_r+\mb e_\theta\partial_ru_\theta)+r^{-1}\mb e_\theta\otimes\brac{\mb e_r(\partial_\theta u_r-u_\theta)+\mb e_\theta(\partial_\theta u_\theta+u_r)}.
\end{align}
The entries of $r\grad\mb u$ in this frame are $\dr u_r$, $\dr u_\theta$, $\partial_\theta u_r-u_\theta$ and $\partial_\theta u_\theta+u_r$; these are the four functions transformed in the proof of Lemma \ref{Korn on the wedge}. The symmetric gradient $\S(\mb u)={1\over 2}((\grad\mb u)+(\grad\mb u)^T)$ therefore has the two components
\begin{subequations}\label{stress compo}
\begin{align}
\mb e_\theta\cdot\S(\mb u)\cdot\mb e_\theta&=r^{-1}(\partial_\theta u_\theta+u_r)=-\partial_ru_r+\grad\cdot\mb u\\
\mb e_r\cdot\S(\mb u)\cdot\mb e_\theta&={1\over 2}(\partial_r u_\theta+r^{-1}(\partial_\theta u_r-u_\theta))=\partial_ru_\theta-{1\over 2}\grad^\perp\cdot\mb u
\end{align}
\end{subequations}
the first expression in each being the one read off from \eqref{grad polar} and the second the one used on $\partial_1\Omega$, where the first of the two reduces further to $-\partial_ru_r$ because $\mb u$ is divergence-free. The shear $\partial_ru_\theta$, in which the profile evolution equation and the contact angle condition \eqref{contact angle without profile} are expressed, is in particular a combination of the Cartesian first derivatives with coefficients bounded on $\Omega$,
\begin{align}\label{radial shear cartesian}
\partial_ru_\theta&=r^{-2}\dr_{-1}(ru_\theta)=r^{-2}(x\partial_1+y\partial_2-1)(u_2x-u_1y)\nonumber\\
&=r^{-2}(x^2\partial_1u_2-xy\partial_1u_1-y^2\partial_2u_1+xy\partial_2u_2),
\end{align}
so that its value at the contact point along $\partial_1\Omega$ is read off from $\grad\mb u(\mb 0)$; this is used in the proof of Proposition \ref{corner ladder}.

At second order the double rotation is the Laplacian up to the divergence,
\begin{align}
\grad^\perp\grad^\perp\cdot\mb u&=\grad^\perp(-\partial_2u_1+\partial_1u_2)
=\mat{\partial_2^2u_1-\partial_1\partial_2u_2\\-\partial_1\partial_2u_1+\partial_1^2u_2}=\mat{\partial_2^2u_1+\partial_1^2u_1\\\partial_2^2u_2+\partial_1^2u_2}-\grad\grad\cdot\mb u\nonumber\\
&=\lpc\mb u-\grad\grad\cdot\mb u,\label{lpc form}
\end{align}
so that $\grad^\perp\grad^\perp\cdot\mb u=\lpc\mb u$ for divergence-free $\mb u$; this is the identity behind the choice of the test field $\smash[b]{\mb V^\theta_{a,j}}$ of Definition \ref{test functions}. Finally, the angular derivative and the gradient fail to commute by exactly one rotation, $[\partial_\theta,\grad]=\grad^\perp$.

\subsection{Powers and logarithms on the Mellin side}\label{Mellin identities}

The substitution $r=e^x$ identifies the Mellin transform of Subsection \ref{norms and Mellin transform} with the Fourier transform, which is where its elementary properties, and its $L^2$ theory, come from.

The corner blocks \eqref{higher corner coefficients} of Section \ref{corner expansions} are sums of terms $r^d(\ln r)^n$, and the integral defining the Mellin transform of such a term converges for no $\lambda$. The proof of Lemma \ref{block matching} below uses instead the Mellin transforms of its restrictions to $(0,1)$ and to $[1,\infty)$.

\begin{lemma}[Mellin transforms of truncated powers]\label{truncated powers}
Let $d\in\R$ and $n\in\N_0$. The Mellin transform of $\mb 1_{(0,1)}r^d(\ln r)^n$ for $\Re\lambda<d$ and that of $-\mb 1_{[1,\infty)}r^d(\ln r)^n$ for $\Re\lambda>d$ are both equal to
\begin{align*}
-{n!\over\sqrt{2\pi}}{1\over(\lambda-d)^{n+1}}.
\end{align*}
\end{lemma}
\begin{proof}
For $n=0$ the two Mellin transforms are $(2\pi)^{-1/2}\int_0^1r^{d-\lambda-1}\d r$ and $-(2\pi)^{-1/2}\int_1^\infty r^{d-\lambda-1}\d r$, both equal to $-(2\pi)^{-1/2}(\lambda-d)^{-1}$. Applying $(-\partial_\lambda)^n$ under the integral sign, which multiplies the integrands by $(\ln r)^n$, gives the general case.
\end{proof}

In Section \ref{corner expansions} a datum with a finite norm at one weight is written as a block plus a remainder with a finite norm at a higher weight, and the line of integration of the Mellin transform of the solution is moved from the line of the lower weight to that of the higher one. Moving the line requires the Mellin transform of the datum on the first line and that of the remainder on the second to be the traces of one function meromorphic between the two lines, and the next lemma shows that they are. For short, Sections \ref{corner expansions} and \ref{nonlinear existence} call such a function meromorphic, or holomorphic, on the closed strip between the two lines, and refer to its traces as its restrictions to the two lines.

\begin{lemma}[Matching a corner block across a Mellin line]\label{block matching}
Let $a_0<d<a_1$, and let $f_0$ and $f_1$, square integrable on $(0,\infty)$ against $r^{-2a_0-1}\d r$ and $r^{-2a_1-1}\d r$ respectively, differ by a block,
\begin{align*}
f_0-f_1=\sum_{k=1}^mc_k{r^d(\ln r)^{k-1}\over(k-1)!},
\end{align*}
with constants $c_k$. Then the Mellin transforms of $f_0$ on the line $\Re\lambda=a_0$ and of $f_1$ on the line $\Re\lambda=a_1$ are the traces, as limits in $L^2$ from the lines $\Re\lambda=a$ between them, of one function meromorphic on the strip $a_0<\Re\lambda<a_1$, with its only pole at $\lambda=d$, of order at most $m$, and square integrable on the parts $|\Im\lambda|\geq 1$ of these lines, uniformly in $a$; in the notation \eqref{higher corner coefficients}, $\smash[b]{c_{f_0,\dr_{-d}^k}}=c_k$ for $k\leq m$.
\end{lemma}
\begin{proof}
The function equal to $f_1$ on $(0,1)$ and to $f_0$ on $[1,\infty)$ is square integrable against $r^{-2a-1}\d r$ uniformly in $a\in[a_0,a_1]$, since $r^{-2a}\leq r^{-2a_1}$ for $r<1$ and $r^{-2a}\leq r^{-2a_0}$ for $r\geq 1$. So its Mellin transform is holomorphic on the strip $a_0<\Re\lambda<a_1$, by the Cauchy--Schwarz inequality, and square integrable on the lines $\Re\lambda=a$ uniformly in $a\in[a_0,a_1]$, by \eqref{Plancherel's identity}. By Lemma \ref{Mellin properties}(i), \eqref{Plancherel's identity} and dominated convergence, the restrictions of this Mellin transform to the lines $\Re\lambda=a$ converge in $L^2$, as $a\to a_0$ and as $a\to a_1$, to the Mellin transforms of the function on the lines $\Re\lambda=a_0$ and $\Re\lambda=a_1$. Adding $\mb 1_{(0,1)}(f_0-f_1)$ to this function gives $f_0$, and adding $-\mb 1_{[1,\infty)}(f_0-f_1)$ gives $f_1$. By Lemma \ref{truncated powers} with $n=k-1$, the Mellin transforms of these two terms, for $\Re\lambda<d$ and for $\Re\lambda>d$ respectively, are both
\begin{align*}
-{1\over\sqrt{2\pi}}\sum_{k=1}^m{c_k\over(\lambda-d)^k},
\end{align*}
which is square integrable on the parts $|\Im\lambda|\geq 1$ of the lines uniformly and, on the lines $\Re\lambda=a\neq d$, depends continuously on $a$ in $L^2$. Adding it to the Mellin transform of the function above gives the function sought, and by \eqref{principal part}, $\smash[b]{c_{f_0,\dr_{-d}^k}}=c_k$.
\end{proof}

\section{Useful inequalities}\label{useful inequalities}

\subsection{Interpolation in the weight}\label{interpolation in the weight}

The norms of Definition \ref{weighted Sobolev norms} are $L^2$ norms against a power of $r$, so they interpolate in the weight.

\begin{proposition}[H\"older's inequality]\label{Holder's inequality}
Let $(S,\mu)$ be a measure space and let $r$ be a positive measurable function on it. If $a\leq b\leq c$ and $a<c$, then
\begin{align*}
\|r^{b}u\|_{L^2(S,\mu)}\leq\|r^{a}u\|_{L^2(S,\mu)}^{{c-b\over c-a}}\|r^{c}u\|_{L^2(S,\mu)}^{{b-a\over c-a}}.
\end{align*}
And when $a\leq b\leq(a+c)/2$ we have
\begin{align*}
\|r^{2b}uv\|_{L^1(S,\mu)}
&\leq\|r^{a}u\|_{L^2(S,\mu)}\|r^{a}v\|_{L^2(S,\mu)}^{{a+c-2b\over c-a}}\|r^{c}v\|_{L^2(S,\mu)}^{2{b-a\over c-a}}
\end{align*}
\end{proposition}
\begin{proof}
Both follow from H\"older's inequality at the exponents $1/\gamma$ and $1/(1-\gamma)$, where $b=a\gamma+c(1-\gamma)$; the second needs a further splitting of the $u,v$ factor, admissible because $b\leq(a+c)/2$ makes $\gamma\geq{1\over 2}$.
\end{proof}

\subsection{Hardy's inequality and the behaviour at infinity}\label{Hardy and infinity}

Hardy's inequality trades a weight for a derivative at either end of the ray, and is the one dimensional engine behind the decay statements that follow. We state it with the constant $1/k$, which is sharp, and with the existence of the boundary value at the relevant end as part of the conclusion rather than as a hypothesis, since that is the form in which it is applied.

\begin{proposition}[Hardy's inequality]\label{Hardy's inequality}
Let $k>0$ and let $u$ be locally absolutely continuous on $(0,\infty)$. If the RHS of the first inequality below is finite, then $u(\infty):=\lim_{r\to\infty}u(r)$ exists and
\begin{align*}
\|r^{-{1\over 2}+k}(u-u(\infty))\|_{L^2([0,\infty))}\leq{1\over k}\|r^{{1\over 2}+k}u'\|_{L^2([0,\infty))};
\end{align*}
if the RHS of the second inequality is finite, then $u(0):=\lim_{r\to 0}u(r)$ exists and
\begin{align*}
\|r^{-{1\over 2}-k}(u-u(0))\|_{L^2([0,\infty))}\leq{1\over k}\|r^{{1\over 2}-k}u'\|_{L^2([0,\infty))}.
\end{align*}
\end{proposition}
\begin{proof}
For the first inequality, the Cauchy--Schwarz inequality gives, for $0<r<R$,
\begin{align*}
|u(R)-u(r)|\leq\int_r^R|u'|\d s
&\leq\brac{\int_r^\infty s^{-1-2k}\d s}^{1\over 2}\brac{\int_r^\infty s^{1+2k}|u'|^2\d s}^{1\over 2}\\
&={r^{-k}\over\sqrt{2k}}\brac{\int_r^\infty s^{1+2k}|u'|^2\d s}^{1\over 2},
\end{align*}
so that $u(\infty)$ exists and $v:=u-u(\infty)$ satisfies $r^{2k}|v(r)|^2\to 0$ as $r\to\infty$. For $0<r_0<r_1$, integrating by parts and dropping the boundary term at $r_0$, which has a favourable sign,
\begin{align*}
\int_{r_0}^{r_1}r^{2k-1}|v|^2\d r&={r_1^{2k}|v(r_1)|^2-r_0^{2k}|v(r_0)|^2\over 2k}-{1\over k}\int_{r_0}^{r_1}r^{2k}\Re(\bar vv')\d r\\
&\leq{r_1^{2k}|v(r_1)|^2\over 2k}+{1\over k}\brac{\int_{r_0}^{r_1}r^{2k-1}|v|^2\d r}^{1\over 2}\brac{\int_0^\infty r^{2k+1}|u'|^2\d r}^{1\over 2}.
\end{align*}
The LHS is finite, $v$ being continuous on $[r_0,r_1]$, so this is a quadratic inequality for its square root, which shows that it stays bounded as $r_1\to\infty$; letting $r_1\to\infty$ the boundary term vanishes, and the inequality on $[r_0,\infty)$ follows after dividing by the square root of the LHS. Letting $r_0\to 0$ by monotone convergence gives the claim. The second inequality is the same argument with the roles of $0$ and $\infty$ exchanged, the favourable boundary term now being the one at $r_1$.
\end{proof}

On the wedge it is applied ray by ray. The form in which the decay at infinity is available for the elements of $H_s^l$, whose norm controls $\grad\mb u$ but not $\mb u$ itself, is the following: the gradient bound pins $\mb u$ down at infinity up to a constant, and the constant vanishes by the boundary term on $\partial_0\Omega$. It is what makes Hardy's inequality applicable at the weight $1-\kappa$ throughout Section \ref{linear existence}.

\begin{lemma}[Limits at infinity and Hardy's inequality]\label{limit at infinity}
Let $k>0$ and suppose $\|\grad\mb u\|_{\Omega,-k}<\infty$. Then there is $\mb c\in\C^2$ such that $\mb u(r,\theta)\to\mb c$ as $r\to\infty$ for a.e.\ $\theta\in[0,\theta_b]$, and
\begin{align*}
\|\mb u-\mb c\|_{\Omega,1-k}\leq k^{-1}\|\grad\mb u\|_{\Omega,-k}.
\end{align*}
If moreover $\|\mb u\|_{\partial_0\Omega,0}<\infty$, then $\mb c=\mb 0$. In particular every $\mb u\in H_s^l$ satisfies $\|\mb u\|_{\Omega,1-\kappa}\leq\kappa^{-1}\|\grad\mb u\|_{\Omega,-\kappa}$.
\end{lemma}
\begin{proof}
By Fubini's theorem, for a.e.\ $\theta$ we have
\begin{align*}
\int_0^\infty r^{1+2k}|\partial_r\mb u(r,\theta)|^2\d r\leq\int_0^\infty r^{1+2k}|\grad\mb u(r,\theta)|^2\d r<\infty,
\end{align*}
so by the first inequality of Proposition \ref{Hardy's inequality} the limit $\mb c(\theta):=\mb u(\infty,\theta)$ exists and
\begin{align*}
\int_0^\infty r^{2k-1}|\mb u(r,\theta)-\mb c(\theta)|^2\d r\leq k^{-2}\int_0^\infty r^{1+2k}|\grad\mb u(r,\theta)|^2\d r.
\end{align*}
The limit does not depend on $\theta$: for a.e.\ $r$ the function $\mb u(r,\ph)$ is absolutely continuous on $[0,\theta_b]$ with $|\partial_\theta\mb u|\leq r|\grad\mb u|$, so that for a.e.\ $\theta_1,\theta_2$
\begin{align*}
\int_R^{2R}|\mb u(r,\theta_1)-\mb u(r,\theta_2)|^2{\d r\over r}
\leq\theta_b\int_R^{2R}\int_0^{\theta_b}r|\grad\mb u|^2\d\theta\d r
\leq\theta_bR^{-2k}\|\grad\mb u\|_{\Omega,-k}^2\to 0\qquad(R\to\infty),
\end{align*}
while the LHS tends to $(\ln 2)|\mb c(\theta_1)-\mb c(\theta_2)|^2$. Integrating the Hardy bound over $\theta$ gives the displayed inequality. If $\|\mb u\|_{\partial_0\Omega,0}<\infty$, the same estimate with $\theta_2=0$ shows $\int_R^{2R}|\mb u(r,0)-\mb c|^2r^{-1}\d r\to 0$, whereas $\int_R^{2R}|\mb u(r,0)|^2r^{-1}\d r\leq R^{-1}\|\mb u\|_{\partial_0\Omega,0}^2\to 0$; hence $(\ln 2)|\mb c|^2=0$. The last statement is the case $k=\kappa$, the norm of $H_s^l$ containing $\|\grad\mb u\|_{\Omega,-\kappa}$ and $\|u_r\|_{\partial_0\Omega,0}=\|\mb u\|_{\partial_0\Omega,0}$, using $u_\theta|_{\partial_0\Omega}=0$.
\end{proof}

\subsection{Continuity at the corner and the embeddings}\label{corner continuity and embeddings}

At the other end of the ray a weight $a>0$ to spare buys pointwise information at the corner. The mechanism of both lemmas below is the fundamental theorem of calculus in $r$ along a fixed ray. The increment along the ray is estimated by the Cauchy--Schwarz inequality against $r^{a-1/2}$, which is square integrable at the origin exactly because $a>0$, and by the trace estimate of Lemma \ref{trace}. The trace estimate is run on the longer of the two intervals into which $\theta_0$ divides $[0,\theta_b]$, so that the ray may be an interior one without loss in $\theta_b$. The value at the corner is then forced to vanish by the finiteness of the norm at hand.

\begin{lemma}[Continuity at the corner]\label{continuity at the corner}
Suppose $a>0$. If $\|u\|_{\Omega,1+a,1,1}<\infty$, then $u$ extends to a continuous function on $\overline\Omega$, and $u(\mb 0)=0$.
\end{lemma}
\begin{proof}
We have
\begin{align*}
|u(r_1,\theta_0)-u(r_0,\theta_0)|
&=\abs{\int_{r_0}^{r_1}\partial_r u(r,\theta_0)\d r}
\leq\|r^{{1\over 2}-a}\partial_r u(r,\theta_0)\|_{L^2_r([r_0,r_1])}\|r^{a-{1\over 2}}\|_{L^2_r([r_0,r_1])}\\
&=\brac{r_1^{2a}-r_0^{2a}\over 2a}^{1\over 2}\|r^{{1\over 2}-a}\partial_r u(r,\theta_0)\|_{L^2_r([r_0,r_1])}\\
&\leq\brac{r_1^{2a}-r_0^{2a}\over 2a}^{1\over 2}\brac{\brac{2\over\theta_b}^{1\over 2}\|\partial_ru\|_{\Omega,a}+\sqrt 2\|\partial_ru\|_{\Omega,a}^{1\over 2}\|r\grad\partial_ru\|_{\Omega,a}^{1\over 2}}\\
&\lesssim(r_1^{2a}-r_0^{2a})^{1\over 2}\|u\|_{\Omega,1+a,1,1}
\end{align*}
by the proof of Lemma \ref{trace}, run on the longer of the two intervals $[0,\theta_0]$ and $[\theta_0,\theta_b]$, whose length is at least $\theta_b/2$, in place of $[0,\theta_b]$, and
\begin{align*}
&|u(r_1,\theta_1)-u(r_1,\theta_0)|\\
&=\abs{\int_{\theta_0}^{\theta_1}\partial_\theta u(r_1,\theta)\d\theta}
\leq\|r_1^{-{1\over 2}-a}\partial_\theta u(r_1,\theta)\|_{L^2_\theta([\theta_0,\theta_1])}\|r_1^{a+{1\over 2}}\|_{L^2_\theta([\theta_0,\theta_1])}\\
&\leq r_1^{a+{1\over 2}}|\theta_1-\theta_0|^{1\over 2}\brac{\int_{\theta_0}^{\theta_1}r_1^{-1-2a}|\partial_\theta u(r_1,\theta)|^2\d\theta}^{1\over 2}\\
&\leq r_1^{a+{1\over 2}}|\theta_1-\theta_0|^{1\over 2}\brac{-\int_{\theta_0}^{\theta_1}\int_{r_1}^\infty\partial_r(r^{-1-2a}|\partial_\theta u(r,\theta)|^2)\d r\d\theta}^{1\over 2}\\
&\leq r_1^{a+{1\over 2}}|\theta_1-\theta_0|^{1\over 2}\brac{\int_{\theta_0}^{\theta_1}\int_{r_1}^\infty\brac{(1+2a)r^{-2-2a}|\partial_\theta u(r,\theta)|^2+2r^{-1-2a}|\partial_\theta u(r,\theta)||\partial_r\partial_\theta u(r,\theta)|}\d r\d\theta}^{1\over 2}\\
&\leq r_1^{a+{1\over 2}}|\theta_1-\theta_0|^{1\over 2}r_1^{-{1\over 2}}\brac{\int_{\theta_0}^{\theta_1}\int_{r_1}^\infty\brac{(1+2a)r^{-1-2a}|\partial_\theta u|^2+2r^{-{1\over 2}-2a}r^{1\over 2}|\partial_\theta u||\partial_r\partial_\theta u|}\d r\d\theta}^{1\over 2}\\
&\leq r_1^{a}|\theta_1-\theta_0|^{1\over 2}\brac{(1+2a)^{1\over 2}\|r^{-1}\partial_\theta u\|_{\Omega,a}+\sqrt 2\|r^{-1}\partial_\theta u\|_{\Omega,a}^{1\over 2}\|\partial_r(r\grad u)\|_{\Omega,a}^{1\over 2}}\\
&\lesssim r_1^{a}|\theta_1-\theta_0|^{1\over 2}\|u\|_{\Omega,1+a,1,1}
\end{align*}
This proves that $u$ extends continuously to $\overline\Omega$. Now we must have $u(\mb 0)=0$, otherwise $\|u\|_{\Omega,1+a}=\infty$.
\end{proof}

The same two moves, carried out with the weight inside the derivative, give the embeddings into the weighted supremum norms, on the ray and on the wedge. It is in this form that the nonlinear estimates of Section \ref{nonlinear existence} take the factor carrying the fewer derivatives out of $L^2$.

\begin{lemma}[Embeddings]\label{embeddings}
Let $a\in\R$. We have
\begin{enumerate}
\item $\|r^{-a}u\|_{L^\infty([0,\infty))}\lesssim_a\|u\|_{[0,\infty),a+1/2,1}$.
\item $\|r^{-a}u\|_{L^\infty(\Omega)}\lesssim_a\theta_b^{-1/2}\|u\|_{\Omega,a+1,1,1}$.
\end{enumerate}
\end{lemma}
\begin{proof}
We have
\begin{enumerate}
  \item Note that $\liminf_{r\to 0}|r^{-a}u(r)|=0=\liminf_{r\to \infty}|r^{-a}u(r)|$, otherwise $\|u\|_{[0,\infty),a+1/2}=\infty$. We have
      \begin{align*}
|r_1^{-a}u(r_1)|^2-|r_0^{-a}u(r_0)|^2
&=\int_{r_0}^{r_1}\partial_r(r^{-2a}u(r)\bar u(r))\d r\\
&\leq 2\int_{r_0}^{r_1}r^{-2a}|\partial_ru(r)||u(r)|\d r
-2a\int_{r_0}^{r_1}r^{-2a-1}u(r)\bar u(r)\d r\\
&\lesssim_a\|u\|_{[0,\infty),a+{1\over 2},1}^2.
\end{align*}
Sending $r_0$ to 0 or $\infty$ along a sequence realising the lower limit gives us the desired result.
  \item This is (i) on the ray $\theta=\theta_0$, the one dimensional norm along the ray being supplied by the trace estimate of Lemma \ref{trace}, run on the longer of the two intervals $[0,\theta_0]$ and $[\theta_0,\theta_b]$ as in the proof of Lemma \ref{continuity at the corner}:
  \begin{align*}
  \|u(\ph,\theta_0)\|_{[0,\infty),a+{1\over 2},1}^2\lesssim\theta_b^{-1}\|u\|_{\Omega,1+a,1}^2+\|u\|_{\Omega,1+a,1}\|\grad u\|_{\Omega,a,1}\lesssim\theta_b^{-1}\|u\|_{\Omega,a+1,1,1}^2,
  \end{align*}
  the last step because $\theta_b<\pi/2$. (A nonzero limit along a single ray would not by itself make the bulk norm $\|u\|_{\Omega,1+a}$ infinite; it is the finiteness of the one dimensional norm along the ray that excludes it.)\qedhere
\end{enumerate}
\end{proof}

\section{Resonance factors}\label{resonance factors}

The solution formulae on the Mellin side are quotients, and on a Mellin line $\Re\lambda=c$ we bound their denominators from below. The constants in these lower bounds degenerate as the line approaches a zero of the denominator that is not carried by an explicit factor of the bound; these constants, or their reciprocals, are the resonance factors of this paper. There are four.
\begin{itemize}
\item $|\cos(c\theta_b)|^{-1}$, for $|c|\theta_b<\pi/2$. On the line $\Re\lambda=c$ Lemma \ref{cosh bound} gives $|\cos(\lambda\theta_b)|\gtrsim\cos(c\theta_b)e^{|\Im\lambda|\theta_b}$. The denominator $\cos(\lambda\theta_b)$ is that of the mixed problem of Lemma \ref{mixed wedge problem}, on which the harmonic extension of Lemma \ref{mixed harmonic extension} and the Helmholtz projection of Definition \ref{the two projections} are built. The factor blows up as $|c|\theta_b\uparrow\pi/2$, that is, as the line approaches a zero $\pm\pi/(2\theta_b)$ of $\cos(\lambda\theta_b)$; for the harmonic extension the zero $\pi/(2\theta_b)$ gives the resonance $r^{\pi/(2\theta_b)}$ of the flattening.
\item $\Gamma_{c,\theta}=\min\{1,\pi-|c\theta|\}$ of \eqref{sine resonance factor}, for $|c\theta|<\pi$, the distance of $c\theta$ from the zeros $\pm\pi$ of the sine, capped at $1$. Where $|\Im\lambda\theta|<C$ on the line $\Re\lambda=c$, Lemma \ref{cosh bound} gives $|\sin(\lambda\theta)|\gtrsim_C\Gamma_{c,\theta}|\lambda\theta|$, the zero $\lambda=0$ being carried by the factor $|\lambda\theta|$. The denominator $\sin(\lambda\theta_b)$ occurs in the Neumann problem of Lemma \ref{Neumann wedge problem}, and the denominators $\sin(\lambda\theta_b)$ and $\sin(\lambda_{\pm 1}\theta_b)$ in the solution formulae of Subsection \ref{positive-time regularity}, which carry them in place of $H$; there $\Gamma(\lambda,\theta_b)$ and $\Gamma(M,\kappa,\theta_b)$ are products of two factors $\Gamma_{c,\theta_b}$.
\item $\Gamma_c=\<c\>/|c|$ of Proposition \ref{Mellin system initial time}, for real $c\neq 0$, the supremum of $\<\lambda\>/|\lambda|$ over the line $\Re\lambda=c$, so that $|\lambda|\geq\Gamma_c^{-1}\<\lambda\>$ there. It blows up as the line approaches $\lambda=0$. In that proposition $\Gamma_{a-1}$, $\Gamma_a$ and $\Gamma_{a-2}$ bound the negative powers of $\lambda$, $\lambda_1$ and $\lambda_{-1}$ that the solution formulae carry on the line $\Re\lambda=a-1$, the pressure formula being evaluated at $\lambda_{-1}$.
\item $\Gamma^H_{c,\theta_b}$ of \eqref{Gamma H definition}, for $|c|<\Lambda^*(\theta_b)$. The characteristic function $H$ of \eqref{H definition} is, up to the factor $4\lambda$, the denominator of the Green's function of Lemma \ref{Green's function initial time}, and on the line $\Re\lambda=c$ we have $|H(\lambda)|\geq\Gamma^H_{c,\theta_b}\min\{|\lambda|\theta_b,1\}e^{2|\Im\lambda|\theta_b}$ by \eqref{H lower bound}. The factor degenerates as $|c|\uparrow\Lambda^*(\theta_b)$, that is, as the line approaches the edge of the zero-free strip of Definition \ref{zero-free strip}.
\end{itemize}
Without subscripts, $\Gamma$ is Euler's Gamma function in the proof of Theorem \ref{linear solution in time}, and it abbreviates $\Gamma(\lambda,\theta_b)$ in the proof of Lemma \ref{A coefficient estimates}.

\end{document}